\documentclass{article}
\pdfoutput=1	%
\usepackage{fullpage}
\usepackage[breaklinks,hidelinks]{hyperref}
\hypersetup{
  pdftitle={Internal Reasoning and Consistency in Reflective Grounded Arithmetic},
  pdfauthor={Bryan Ford}}

\usepackage{amsmath}
\usepackage{amssymb}
\usepackage{amsthm}

\newtheorem{thm}{Theorem}[section]
\newtheorem{cor}[thm]{Corollary}	%

\usepackage{times}
\usepackage{xspace}
\usepackage{cite}		%
\usepackage{bussproofs}
\usepackage{stmaryrd}
\usepackage[noabbrev,capitalize]{cleveref}
\usepackage{thmtools}
\usepackage{xstring}
\usepackage{graphicx}
\usepackage{subcaption}
\usepackage{nicematrix}		%
\usepackage{longtable}
\usepackage{pifont}

\usepackage{isabelle}
\usepackage{isabellesym}
\usepackage{pdfpages}
\usepackage{ragged2e}
\newcommand{\SNIPTEXT}[1]{#1} %
\newcommand{\SNIP}[2]{\expandafter\newcommand\csname snippet--#1\endcsname{#2}}
\IfFileExists{snips.tex}{\SNIP{ConRGA-Paper-Statements:theory:ConRGA-Paper-Statements-0}{%
\isakeywordONE{theory}\isamarkupfalse%
\ ConRGA{\isacharunderscore}{\kern0pt}Paper{\isacharunderscore}{\kern0pt}Statements\isanewline%
}%
\SNIP{ConRGA-Paper-Statements:theory:ConRGA-Paper-Statements-1}{%
\ \ \isakeywordTWO{imports}\isanewline%
}%
\SNIP{ConRGA-Paper-Statements:theory:ConRGA-Paper-Statements-2}{%
\ \ \ \ {\isachardoublequoteopen}RGA{\isacharunderscore}{\kern0pt}Kit{\isachardot}{\kern0pt}Entailment{\isacharunderscore}{\kern0pt}Kit{\isachardoublequoteclose}\ {\isachardoublequoteopen}RGA{\isacharunderscore}{\kern0pt}Kit{\isachardot}{\kern0pt}Code{\isacharunderscore}{\kern0pt}Index{\isachardoublequoteclose}\ {\isachardoublequoteopen}RGA{\isacharunderscore}{\kern0pt}Kit{\isachardot}{\kern0pt}Term{\isacharunderscore}{\kern0pt}Coders{\isachardoublequoteclose}\isanewline%
}%
\SNIP{ConRGA-Paper-Statements:theory:ConRGA-Paper-Statements-3}{%
\ \ \ \ {\isachardoublequoteopen}RGA{\isacharunderscore}{\kern0pt}Decider{\isachardot}{\kern0pt}Decider{\isacharunderscore}{\kern0pt}Truth{\isachardoublequoteclose}\ {\isachardoublequoteopen}RGA{\isacharunderscore}{\kern0pt}Checker{\isachardot}{\kern0pt}Chk{\isacharunderscore}{\kern0pt}Faithful{\isachardoublequoteclose}\ {\isachardoublequoteopen}RGA{\isacharunderscore}{\kern0pt}Checker{\isachardot}{\kern0pt}Chk{\isacharunderscore}{\kern0pt}Completeness{\isachardoublequoteclose}\isanewline%
}%
\SNIP{ConRGA-Paper-Statements:theory:ConRGA-Paper-Statements-4}{%
\ \ \ \ {\isachardoublequoteopen}RGA{\isacharunderscore}{\kern0pt}QF{\isachardot}{\kern0pt}QF{\isacharunderscore}{\kern0pt}Last{\isacharunderscore}{\kern0pt}Arms{\isachardoublequoteclose}\ {\isachardoublequoteopen}RGA{\isacharunderscore}{\kern0pt}Goedel{\isachardot}{\kern0pt}Goedel{\isacharunderscore}{\kern0pt}II{\isachardoublequoteclose}\ {\isachardoublequoteopen}RGA{\isacharunderscore}{\kern0pt}Goedel{\isachardot}{\kern0pt}Goedel{\isacharunderscore}{\kern0pt}II{\isacharunderscore}{\kern0pt}PR{\isachardoublequoteclose}\ {\isachardoublequoteopen}RGA{\isacharunderscore}{\kern0pt}Goedel{\isachardot}{\kern0pt}Goedel{\isacharunderscore}{\kern0pt}Cut{\isachardoublequoteclose}\isanewline%
}%
\SNIP{ConRGA-Paper-Statements:theory:ConRGA-Paper-Statements-5}{%
\ \ \ \ {\isachardoublequoteopen}RGA{\isacharunderscore}{\kern0pt}Goedel{\isachardot}{\kern0pt}Entailment{\isacharunderscore}{\kern0pt}Goedel{\isachardoublequoteclose}\ {\isachardoublequoteopen}RGA{\isacharunderscore}{\kern0pt}Goedel{\isachardot}{\kern0pt}Goedel{\isacharunderscore}{\kern0pt}Conditions{\isachardoublequoteclose}\ {\isachardoublequoteopen}ConRGA{\isachardot}{\kern0pt}Wfd{\isacharunderscore}{\kern0pt}General{\isachardoublequoteclose}\isanewline%
}%
\SNIP{ConRGA-Paper-Statements:theory:ConRGA-Paper-Statements-6}{%
\ \ \ \ {\isachardoublequoteopen}ConRGA{\isachardot}{\kern0pt}Guarded{\isacharunderscore}{\kern0pt}Residual{\isachardoublequoteclose}\ {\isachardoublequoteopen}RGA{\isacharunderscore}{\kern0pt}Liar{\isacharunderscore}{\kern0pt}Div{\isachardot}{\kern0pt}Entailment{\isacharunderscore}{\kern0pt}Truth{\isachardoublequoteclose}\ {\isachardoublequoteopen}RGA{\isacharunderscore}{\kern0pt}Liar{\isacharunderscore}{\kern0pt}Div{\isachardot}{\kern0pt}Entailment{\isacharunderscore}{\kern0pt}Liar{\isachardoublequoteclose}\isanewline%
}%
\SNIP{ConRGA-Paper-Statements:theory:ConRGA-Paper-Statements-7}{%
\ \ \ \ {\isachardoublequoteopen}RGA{\isacharunderscore}{\kern0pt}Liar{\isacharunderscore}{\kern0pt}Div{\isachardot}{\kern0pt}Truth{\isacharunderscore}{\kern0pt}Limits{\isachardoublequoteclose}\ {\isachardoublequoteopen}ConOGA{\isachardot}{\kern0pt}ConOGA{\isacharunderscore}{\kern0pt}Level{\isachardoublequoteclose}\isanewline%
}%
\SNIP{ConRGA-Paper-Statements:theory:ConRGA-Paper-Statements-8}{%
\ \ \ \ {\isachardoublequoteopen}RGA{\isacharunderscore}{\kern0pt}Reflect{\isachardot}{\kern0pt}Reflect{\isacharunderscore}{\kern0pt}QF{\isachardoublequoteclose}\ {\isachardoublequoteopen}RGA{\isacharunderscore}{\kern0pt}Poly{\isachardot}{\kern0pt}Poly{\isacharunderscore}{\kern0pt}Reflect{\isachardoublequoteclose}\ {\isachardoublequoteopen}RGA{\isacharunderscore}{\kern0pt}Tin{\isachardot}{\kern0pt}Tin{\isacharunderscore}{\kern0pt}Univ{\isadigit{2}}{\isachardoublequoteclose}\isanewline%
}%
\SNIP{ConRGA-Paper-Statements:theory:ConRGA-Paper-Statements-9}{%
\isakeywordTWO{begin}%
\endisatagtheory
{\isafoldtheory}%
\isadelimtheory
\endisadelimtheory
\SNIPTEXT{\begin{isamarkuptext}%
A document-only theory: it proves nothing and defines nothing.  For each headline theorem
  of the ConRGA paper and for five core definitions it has one \isa{thm} command naming the
  fact, followed by the fact's statement as Isabelle prints it (a \isa{thm} antiquotation in a
  \isa{text{\isacharunderscore}{\kern0pt}raw} block).  The theories that state the results are the imports above.

  The session's document output is cut into named snippets by the cartouches tool
  (\url{https://github.com/astahfrom/cartouches}) into \isa{doc{\isacharslash}{\kern0pt}conrga{\isacharslash}{\kern0pt}snips{\isachardot}{\kern0pt}tex}.
  The snippet of the fact \isa{NAME} is \isa{ConRGA{\isacharminus}{\kern0pt}Paper{\isacharminus}{\kern0pt}Statements{\isacharcolon}{\kern0pt}thm{\isacharcolon}{\kern0pt}NAME}, with every
  underscore of \isa{NAME} written as a hyphen; its line 0 is the command
  \isa{thm\ NAME} and its line 1 is the printed statement.  So, in the paper,
  \isa{{\isacharbackslash}{\kern0pt}Snippet{\isacharbraceleft}{\kern0pt}ConRGA{\isacharminus}{\kern0pt}Paper{\isacharminus}{\kern0pt}Statements{\isacharcolon}{\kern0pt}thm{\isacharcolon}{\kern0pt}NAME{\isacharbraceright}{\kern0pt}} shows both lines and
  \isa{{\isacharbackslash}{\kern0pt}GetSnip{\isacharbraceleft}{\kern0pt}ConRGA{\isacharminus}{\kern0pt}Paper{\isacharminus}{\kern0pt}Statements{\isacharcolon}{\kern0pt}thm{\isacharcolon}{\kern0pt}NAME{\isacharminus}{\kern0pt}{\isadigit{1}}{\isacharbraceright}{\kern0pt}} the statement alone.  The
  regeneration recipe is in \isa{doc{\isacharslash}{\kern0pt}conrga{\isacharslash}{\kern0pt}briefs{\isacharslash}{\kern0pt}ADDENDUM{\isacharminus}{\kern0pt}T{\isadigit{1}}{\isachardot}{\kern0pt}md}.

  Rules of this theory, so that the snippets stay well formed: every statement is printed inline
  (a display-mode antiquotation would put its line breaks into separate snippets and split a
  LaTeX group across them), and no prose stands between two facts (it would be appended to the
  preceding snippet); prose names commands with \isa{text} antiquotations only (the tool reads a
  command antiquotation as a command and would start a snippet there).  Schematic variables are printed without question marks (the session
  option \isa{show{\isacharunderscore}{\kern0pt}question{\isacharunderscore}{\kern0pt}marks}).

  One notation is added for printing only.  The biconditional \isa{a\ \isactrlbold {\isasymlongleftrightarrow}\ b} that the
  statements use is folded by the printer, inside \isa{rga}, into the abbreviation
  \isa{gprop{\isachardot}{\kern0pt}trmIff{\isacharprime}{\kern0pt}}: the trunk's \isa{trmIff{\isacharprime}{\kern0pt}} (\isa{Grounded{\isacharunderscore}{\kern0pt}Logic}) as
  re-exported by \isa{sublocale\ gprop{\isacharcolon}{\kern0pt}\ gd{\isacharunderscore}{\kern0pt}prop{\isacharunderscore}{\kern0pt}minimal} (\isa{RGA{\isacharunderscore}{\kern0pt}Sound}), an
  instantiation that drops the infix syntax.  It would print as a prefix constant with its
  qualified name.  The \isa{notation} below, the first command inside the context, gives it
  back the infix syntax, so the statements print as their sources write them.  Other
  abbreviations that the printer does not fold back appear expanded: \isa{\isactrlbold {\isasymforall}b}
  prints as \isa{\isactrlbold A\ \isactrlbold {\isasymcdot}\ {\isacharparenleft}{\kern0pt}\isactrlbold {\isasymlambda}b{\isacharparenright}{\kern0pt}}, and the tests \isa{andT},
  \isa{sgnT} and \isa{iszT} of the internal kit print through \isa{gapSub}, and OGA's
  entailment \isa{{\isasymGamma}\ {\isasymturnstile}\isactrlsub O\ c} (the abbreviation \isa{oent} of \isa{ConOGA{\isacharunderscore}{\kern0pt}Seam}) prints as
  \isa{prf{\isacharunderscore}{\kern0pt}derives{\isachardot}{\kern0pt}entails\ oga{\isachardot}{\kern0pt}derives\ {\isasymGamma}\ c}.%
\end{isamarkuptext}}\isamarkuptrue%
}%
\SNIP{ConRGA-Paper-Statements:context:rga-0}{%
\isakeywordONE{context}\isamarkupfalse%
\ rga%
}%
\SNIP{ConRGA-Paper-Statements:notation:gprop-0}{%
\isakeywordONE{notation}\isamarkupfalse%
\ gprop{\isachardot}{\kern0pt}trmIff{\isacharprime}{\kern0pt}\ {\isacharparenleft}{\kern0pt}\isakeywordTWO{infixr}\ {\isachardoublequoteopen}\isactrlbold {\isasymlongleftrightarrow}{\isachardoublequoteclose}\ {\isadigit{1}}{\isadigit{2}}{\isadigit{5}}{\isacharparenright}{\kern0pt}%
}%
\SNIP{ConRGA-Paper-Statements:thm:deduction-theorem-exact-0}{%
\isakeywordONE{thm}\isamarkupfalse%
\ deduction{\isacharunderscore}{\kern0pt}theorem{\isacharunderscore}{\kern0pt}exact%
\isanewline%
}%
\SNIP{ConRGA-Paper-Statements:thm:deduction-theorem-exact-1}{%
\isa{\isaconst{wft}\ \isavar{A}\ {\isasymLongrightarrow}\ {\isacharparenleft}{\kern0pt}{\isasymforall}\isabound{B}{\isachardot}{\kern0pt}\ \isavar{A}\ {\isasymtriangleright}\ \isavar{{\isasymGamma}}\ {\isasymturnstile}\ \isabound{B}\ {\isasymlongrightarrow}\ \isavar{{\isasymGamma}}\ {\isasymturnstile}\ \isavar{A}\ \isactrlbold {\isasymlongrightarrow}\ \isabound{B}{\isacharparenright}{\kern0pt}\ {\isacharequal}{\kern0pt}\ {\isacharparenleft}{\kern0pt}\isavar{{\isasymGamma}}\ {\isasymturnstile}\ \isavar{A}\ \isactrlbold {\isasymor}\ \isactrlbold {\isasymnot}\isavar{A}{\isacharparenright}{\kern0pt}}%
}%
\SNIP{ConRGA-Paper-Statements:thm:sigma-deduction-theorem-0}{%
\isakeywordONE{thm}\isamarkupfalse%
\ sigma{\isacharunderscore}{\kern0pt}deduction{\isacharunderscore}{\kern0pt}theorem%
\isanewline%
}%
\SNIP{ConRGA-Paper-Statements:thm:sigma-deduction-theorem-1}{%
\isa{\isactrlbold v\isactrlsub {\isadigit{0}}\ \isactrlbold N\ {\isasymtriangleright}\ {\isasymupharpoonleft}\isavar{{\isasymGamma}}\ {\isasymturnstile}\ \isavar{D}\ \isactrlbold B\ {\isasymLongrightarrow}\ {\isacharparenleft}{\kern0pt}\isactrlbold {\isasymnot}\isactrlbold A\ \isactrlbold {\isasymcdot}\ {\isacharparenleft}{\kern0pt}\isactrlbold {\isasymlambda}{\isacharparenleft}{\kern0pt}\isactrlbold {\isasymnot}\isavar{D}{\isacharparenright}{\kern0pt}{\isacharparenright}{\kern0pt}\ {\isasymtriangleright}\ \isavar{{\isasymGamma}}\ {\isasymturnstile}\ \isavar{C}{\isacharparenright}{\kern0pt}\ {\isacharequal}{\kern0pt}\ {\isacharparenleft}{\kern0pt}\isavar{{\isasymGamma}}\ {\isasymturnstile}\ \isactrlbold A\ \isactrlbold {\isasymcdot}\ {\isacharparenleft}{\kern0pt}\isactrlbold {\isasymlambda}{\isacharparenleft}{\kern0pt}\isavar{D}\ \isactrlbold {\isasymlongrightarrow}\ {\isasymupharpoonleft}\isavar{C}{\isacharparenright}{\kern0pt}{\isacharparenright}{\kern0pt}{\isacharparenright}{\kern0pt}}%
}%
\SNIP{ConRGA-Paper-Statements:thm:rga-consistent-via-wtrue-0}{%
\isakeywordONE{thm}\isamarkupfalse%
\ rga{\isacharunderscore}{\kern0pt}consistent{\isacharunderscore}{\kern0pt}via{\isacharunderscore}{\kern0pt}wtrue%
\isanewline%
}%
\SNIP{ConRGA-Paper-Statements:thm:rga-consistent-via-wtrue-1}{%
\isa{\isavar{P}\ {\isasymP}{\isasymlbrakk}{\isasymemptyset}\ {\isasymturnstile}\ \isactrlbold {\isadigit{0}}\ \isactrlbold {\isacharequal}{\kern0pt}\ \isactrlbold n{\isadigit{1}}{\isasymrbrakk}\ {\isasymLongrightarrow}\ \isaconst{False}}%
}%
\SNIP{ConRGA-Paper-Statements:thm:rga-complete-int2-0}{%
\isakeywordONE{thm}\isamarkupfalse%
\ rga{\isacharunderscore}{\kern0pt}complete{\isacharunderscore}{\kern0pt}int{\isadigit{2}}%
\isanewline%
}%
\SNIP{ConRGA-Paper-Statements:thm:rga-complete-int2-1}{%
\isa{{\isasymemptyset}\ {\isasymturnstile}\ \isactrlbold A\ \isactrlbold {\isasymcdot}\ {\isacharparenleft}{\kern0pt}\isactrlbold {\isasymlambda}\isactrlbold A\ \isactrlbold {\isasymcdot}\ {\isacharparenleft}{\kern0pt}\isactrlbold {\isasymlambda}{\isacharparenleft}{\kern0pt}\isaconst{gapSub}\ {\isacharparenleft}{\kern0pt}\isaconst{gapSub}\ {\isacharparenleft}{\kern0pt}\isactrlbold n\isaconst{Suc}\ {\isadigit{0}}{\isacharparenright}{\kern0pt}\ {\isacharparenleft}{\kern0pt}\isaconst{gapSub}\ {\isacharparenleft}{\kern0pt}\isactrlbold n\isaconst{Suc}\ {\isadigit{0}}{\isacharparenright}{\kern0pt}\ {\isacharparenleft}{\kern0pt}\isaconst{gapSub}\ {\isacharparenleft}{\kern0pt}\isaconst{gapSub}\ {\isacharparenleft}{\kern0pt}\isactrlbold n\isaconst{Suc}\ {\isadigit{0}}{\isacharparenright}{\kern0pt}\ {\isacharparenleft}{\kern0pt}\isaconst{gapSub}\ {\isacharparenleft}{\kern0pt}\isactrlbold n\isaconst{Suc}\ {\isadigit{0}}{\isacharparenright}{\kern0pt}\ {\isacharparenleft}{\kern0pt}\isaconst{iWft}\ \isactrlbold v\isactrlsub {\isadigit{0}}{\isacharparenright}{\kern0pt}{\isacharparenright}{\kern0pt}{\isacharparenright}{\kern0pt}\ {\isacharparenleft}{\kern0pt}\isaconst{gapSub}\ {\isacharparenleft}{\kern0pt}\isactrlbold n\isaconst{Suc}\ {\isadigit{0}}{\isacharparenright}{\kern0pt}\ {\isacharparenleft}{\kern0pt}\isaconst{iIsCode}\ \isactrlbold v\isactrlsub {\isadigit{0}}{\isacharparenright}{\kern0pt}{\isacharparenright}{\kern0pt}{\isacharparenright}{\kern0pt}{\isacharparenright}{\kern0pt}{\isacharparenright}{\kern0pt}\ {\isacharparenleft}{\kern0pt}\isaconst{gapSub}\ {\isacharparenleft}{\kern0pt}\isactrlbold n\isaconst{Suc}\ {\isadigit{0}}{\isacharparenright}{\kern0pt}\ {\isacharparenleft}{\kern0pt}\isaconst{iPrDecTrueAt}\ \isactrlbold v\isactrlsub {\isadigit{0}}\ {\isacharparenleft}{\kern0pt}\isaconst{trmVar}\ {\isadigit{1}}{\isacharparenright}{\kern0pt}{\isacharparenright}{\kern0pt}{\isacharparenright}{\kern0pt}\ \isactrlbold {\isasymlongrightarrow}\ \isaconst{iHasCert}\ {\isacharparenleft}{\kern0pt}\isaconst{cpPair}\ {\isacharparenleft}{\kern0pt}\isactrlbold n{\isasymlbrace}{\isasymemptyset}{\isasymrbrace}{\isacharparenright}{\kern0pt}\ \isactrlbold v\isactrlsub {\isadigit{0}}{\isacharparenright}{\kern0pt}{\isacharparenright}{\kern0pt}{\isacharparenright}{\kern0pt}{\isacharparenright}{\kern0pt}}%
}%
\SNIP{ConRGA-Paper-Statements:thm:complete-ent-0}{%
\isakeywordONE{thm}\isamarkupfalse%
\ complete{\isacharunderscore}{\kern0pt}ent%
\isanewline%
}%
\SNIP{ConRGA-Paper-Statements:thm:complete-ent-1}{%
\isa{\isaconst{iWtrue}\ \isactrlbold v\isactrlsub {\isadigit{0}}\ {\isasymtriangleright}\ \isaconst{iIsCode}\ \isactrlbold v\isactrlsub {\isadigit{0}}\ \isactrlbold {\isasymnoteq}\ \isactrlbold {\isadigit{0}}\ {\isasymtriangleright}\ \isaconst{iWft}\ \isactrlbold v\isactrlsub {\isadigit{0}}\ \isactrlbold {\isasymnoteq}\ \isactrlbold {\isadigit{0}}\ {\isasymtriangleright}\ \isactrlbold v\isactrlsub {\isadigit{0}}\ \isactrlbold N\ {\isasymtriangleright}\ {\isasymemptyset}\ {\isasymturnstile}\ \isaconst{iHasCert}\ {\isacharparenleft}{\kern0pt}\isaconst{cpPair}\ {\isacharparenleft}{\kern0pt}\isactrlbold n{\isasymlbrace}{\isasymemptyset}{\isasymrbrace}{\isacharparenright}{\kern0pt}\ \isactrlbold v\isactrlsub {\isadigit{0}}{\isacharparenright}{\kern0pt}}%
}%
\SNIP{ConRGA-Paper-Statements:thm:establishment-transparency-0}{%
\isakeywordONE{thm}\isamarkupfalse%
\ establishment{\isacharunderscore}{\kern0pt}transparency%
\isanewline%
}%
\SNIP{ConRGA-Paper-Statements:thm:establishment-transparency-1}{%
\isa{\isaconst{wft}\ \isavar{t}\ {\isasymLongrightarrow}\ {\isacharparenleft}{\kern0pt}{\isasymemptyset}\ {\isasymturnstile}\ \isavar{t}{\isacharparenright}{\kern0pt}\ {\isacharequal}{\kern0pt}\ {\isacharparenleft}{\kern0pt}{\isasymemptyset}\ {\isasymturnstile}\ \isaconst{iWtrue}\ {\isacharparenleft}{\kern0pt}\isactrlbold n{\isasymlbrace}\isavar{t}{\isasymrbrace}{\isacharparenright}{\kern0pt}{\isacharparenright}{\kern0pt}}%
}%
\SNIP{ConRGA-Paper-Statements:thm:tbicond-liar-explosive-diag-0}{%
\isakeywordONE{thm}\isamarkupfalse%
\ tbicond{\isacharunderscore}{\kern0pt}liar{\isacharunderscore}{\kern0pt}explosive{\isacharunderscore}{\kern0pt}diag%
\isanewline%
}%
\SNIP{ConRGA-Paper-Statements:thm:tbicond-liar-explosive-diag-1}{%
\isa{{\isasymlbrakk}\isaconst{wft}\ \isavar{T}{\isacharsemicolon}{\kern0pt}\ \isaconst{Syntax{\isacharunderscore}{\kern0pt}Term{\isachardot}{\kern0pt}closed}\ \isavar{T}{\isasymrbrakk}\ {\isasymLongrightarrow}\ \isavar{T}\ \isactrlbold {\isasymcdot}\ {\isacharparenleft}{\kern0pt}\isactrlbold n{\isasymlbrace}\isactrlbold {\isasymnot}\isavar{T}\ \isactrlbold {\isasymcdot}\ \isaconst{iDiagSub}\ {\isacharparenleft}{\kern0pt}\isactrlbold n{\isasymlbrace}\isactrlbold {\isasymnot}\isavar{T}\ \isactrlbold {\isasymcdot}\ \isaconst{iDiagSub}\ \isactrlbold v\isactrlsub {\isadigit{0}}\ \isactrlbold v\isactrlsub {\isadigit{0}}{\isasymrbrace}{\isacharparenright}{\kern0pt}\ {\isacharparenleft}{\kern0pt}\isactrlbold n{\isasymlbrace}\isactrlbold {\isasymnot}\isavar{T}\ \isactrlbold {\isasymcdot}\ \isaconst{iDiagSub}\ \isactrlbold v\isactrlsub {\isadigit{0}}\ \isactrlbold v\isactrlsub {\isadigit{0}}{\isasymrbrace}{\isacharparenright}{\kern0pt}{\isasymrbrace}{\isacharparenright}{\kern0pt}\ \isactrlbold {\isasymlongleftrightarrow}\ \isactrlbold {\isasymnot}\isavar{T}\ \isactrlbold {\isasymcdot}\ \isaconst{iDiagSub}\ {\isacharparenleft}{\kern0pt}\isactrlbold n{\isasymlbrace}\isactrlbold {\isasymnot}\isavar{T}\ \isactrlbold {\isasymcdot}\ \isaconst{iDiagSub}\ \isactrlbold v\isactrlsub {\isadigit{0}}\ \isactrlbold v\isactrlsub {\isadigit{0}}{\isasymrbrace}{\isacharparenright}{\kern0pt}\ {\isacharparenleft}{\kern0pt}\isactrlbold n{\isasymlbrace}\isactrlbold {\isasymnot}\isavar{T}\ \isactrlbold {\isasymcdot}\ \isaconst{iDiagSub}\ \isactrlbold v\isactrlsub {\isadigit{0}}\ \isactrlbold v\isactrlsub {\isadigit{0}}{\isasymrbrace}{\isacharparenright}{\kern0pt}\ {\isasymtriangleright}\ \isavar{{\isasymGamma}}\ {\isasymturnstile}\ \isactrlbold {\isadigit{0}}\ \isactrlbold {\isacharequal}{\kern0pt}\ \isactrlbold n{\isadigit{1}}}%
}%
\SNIP{ConRGA-Paper-Statements:thm:ordinary-liar-neither-0}{%
\isakeywordONE{thm}\isamarkupfalse%
\ ordinary{\isacharunderscore}{\kern0pt}liar{\isacharunderscore}{\kern0pt}neither%
\isanewline%
}%
\SNIP{ConRGA-Paper-Statements:thm:ordinary-liar-neither-1}{%
\isa{{\isasymemptyset}\ {\isasymturnstile}\ \isactrlbold {\isasymnot}\isaconst{iWtrue}\ {\isacharparenleft}{\kern0pt}\isactrlbold n{\isasymlbrace}\isaconst{trmLiar}{\isasymrbrace}{\isacharparenright}{\kern0pt}\ {\isasymand}\ {\isasymemptyset}\ {\isasymturnstile}\ \isactrlbold {\isasymnot}\isaconst{iWtrue}\ {\isacharparenleft}{\kern0pt}\isactrlbold n{\isasymlbrace}\isactrlbold {\isasymnot}\isaconst{trmLiar}{\isasymrbrace}{\isacharparenright}{\kern0pt}\ {\isasymand}\ {\isasymnot}\ {\isasymemptyset}\ {\isasymturnstile}\ \isaconst{trmLiar}\ {\isasymand}\ {\isasymnot}\ {\isasymemptyset}\ {\isasymturnstile}\ \isactrlbold {\isasymnot}\isaconst{trmLiar}}%
}%
\SNIP{ConRGA-Paper-Statements:thm:refuted-unestablished-trichotomy-0}{%
\isakeywordONE{thm}\isamarkupfalse%
\ refuted{\isacharunderscore}{\kern0pt}unestablished{\isacharunderscore}{\kern0pt}trichotomy%
\isanewline%
}%
\SNIP{ConRGA-Paper-Statements:thm:refuted-unestablished-trichotomy-1}{%
\isa{{\isacharparenleft}{\kern0pt}{\isasymforall}\isabound{t}{\isachardot}{\kern0pt}\ \isaconst{afree}\ \isabound{t}\ {\isasymlongrightarrow}\ {\isasymemptyset}\ {\isasymturnstile}\ \isactrlbold {\isasymnot}\isabound{t}\ {\isasymlongrightarrow}\ {\isasymemptyset}\ {\isasymturnstile}\ \isactrlbold {\isasymnot}\isaconst{iWtrue}\ {\isacharparenleft}{\kern0pt}\isactrlbold n{\isasymlbrace}\isabound{t}{\isasymrbrace}{\isacharparenright}{\kern0pt}{\isacharparenright}{\kern0pt}\ {\isasymand}\ {\isacharparenleft}{\kern0pt}{\isasymforall}\isabound{t}{\isachardot}{\kern0pt}\ {\isasymemptyset}\ {\isasymturnstile}\ \isactrlbold {\isasymnot}\isabound{t}\ {\isasymlongrightarrow}\ \isaconst{ConRGA}\ {\isasymtriangleright}\ {\isasymemptyset}\ {\isasymturnstile}\ \isactrlbold {\isasymnot}\isaconst{iWtrue}\ {\isacharparenleft}{\kern0pt}\isactrlbold n{\isasymlbrace}\isabound{t}{\isasymrbrace}{\isacharparenright}{\kern0pt}{\isacharparenright}{\kern0pt}\ {\isasymand}\ {\isacharparenleft}{\kern0pt}{\isasymforall}\isabound{u}{\isachardot}{\kern0pt}\ \isaconst{wft}\ \isabound{u}\ {\isasymlongrightarrow}\ {\isasymemptyset}\ {\isasymturnstile}\ \isactrlbold {\isasymnot}\isaconst{iWtrue}\ {\isacharparenleft}{\kern0pt}\isactrlbold n{\isasymlbrace}\isactrlbold A\ \isactrlbold {\isasymcdot}\ \isabound{u}{\isasymrbrace}{\isacharparenright}{\kern0pt}\ {\isasymlongrightarrow}\ {\isasymemptyset}\ {\isasymturnstile}\ \isaconst{ConRGA}{\isacharparenright}{\kern0pt}}%
}%
\SNIP{ConRGA-Paper-Statements:thm:con-tschema-decided-0}{%
\isakeywordONE{thm}\isamarkupfalse%
\ con{\isacharunderscore}{\kern0pt}tschema{\isacharunderscore}{\kern0pt}decided%
\isanewline%
}%
\SNIP{ConRGA-Paper-Statements:thm:con-tschema-decided-1}{%
\isa{{\isasymlbrakk}\isaconst{wft}\ \isavar{t}{\isacharsemicolon}{\kern0pt}\ {\isasymemptyset}\ {\isasymturnstile}\ \isavar{t}\ {\isasymor}\ {\isasymemptyset}\ {\isasymturnstile}\ \isactrlbold {\isasymnot}\isavar{t}{\isasymrbrakk}\ {\isasymLongrightarrow}\ \isaconst{ConRGA}\ {\isasymtriangleright}\ \isavar{{\isasymGamma}}\ {\isasymturnstile}\ \isaconst{iWtrue}\ {\isacharparenleft}{\kern0pt}\isactrlbold n{\isasymlbrace}\isavar{t}{\isasymrbrace}{\isacharparenright}{\kern0pt}\ \isactrlbold {\isasymlongleftrightarrow}\ \isavar{t}}%
}%
\SNIP{ConRGA-Paper-Statements:thm:goedelG-under-ConRGA-0}{%
\isakeywordONE{thm}\isamarkupfalse%
\ goedelG{\isacharunderscore}{\kern0pt}under{\isacharunderscore}{\kern0pt}ConRGA%
\isanewline%
}%
\SNIP{ConRGA-Paper-Statements:thm:goedelG-under-ConRGA-1}{%
\isa{\isaconst{ConRGA}\ {\isasymtriangleright}\ {\isasymemptyset}\ {\isasymturnstile}\ \isaconst{goedelG}}%
}%
\SNIP{ConRGA-Paper-Statements:thm:rga-goedel1-0}{%
\isakeywordONE{thm}\isamarkupfalse%
\ rga{\isacharunderscore}{\kern0pt}goedel{\isadigit{1}}%
\isanewline%
}%
\SNIP{ConRGA-Paper-Statements:thm:rga-goedel1-1}{%
\isa{{\isasymnot}\ {\isasymemptyset}\ {\isasymturnstile}\ \isaconst{goedelG}}%
}%
\SNIP{ConRGA-Paper-Statements:thm:rga-goedel1-unrefutable-0}{%
\isakeywordONE{thm}\isamarkupfalse%
\ rga{\isacharunderscore}{\kern0pt}goedel{\isadigit{1}}{\isacharunderscore}{\kern0pt}unrefutable%
\isanewline%
}%
\SNIP{ConRGA-Paper-Statements:thm:rga-goedel1-unrefutable-1}{%
\isa{{\isasymnot}\ {\isasymemptyset}\ {\isasymturnstile}\ \isactrlbold {\isasymnot}\isaconst{goedelG}}%
}%
\SNIP{ConRGA-Paper-Statements:thm:rga-goedel2-0}{%
\isakeywordONE{thm}\isamarkupfalse%
\ rga{\isacharunderscore}{\kern0pt}goedel{\isadigit{2}}%
\isanewline%
}%
\SNIP{ConRGA-Paper-Statements:thm:rga-goedel2-1}{%
\isa{{\isasymnot}\ {\isasymemptyset}\ {\isasymturnstile}\ \isaconst{ConRGA}}%
}%
\SNIP{ConRGA-Paper-Statements:thm:rga-goedel2-internal-contra-0}{%
\isakeywordONE{thm}\isamarkupfalse%
\ rga{\isacharunderscore}{\kern0pt}goedel{\isadigit{2}}{\isacharunderscore}{\kern0pt}internal{\isacharunderscore}{\kern0pt}contra%
\isanewline%
}%
\SNIP{ConRGA-Paper-Statements:thm:rga-goedel2-internal-contra-1}{%
\isa{\isactrlbold {\isasymnot}\isactrlbold A\ \isactrlbold {\isasymcdot}\ {\isacharparenleft}{\kern0pt}\isactrlbold {\isasymlambda}{\isacharparenleft}{\kern0pt}\isactrlbold {\isasymnot}\isaconst{iChk}\ \isactrlbold v\isactrlsub {\isadigit{0}}\ {\isacharparenleft}{\kern0pt}{\isasymupharpoonleft}{\isacharparenleft}{\kern0pt}\isactrlbold n{\isasymlbrace}\isaconst{ConRGA}{\isasymrbrace}{\isacharparenright}{\kern0pt}{\isacharparenright}{\kern0pt}{\isacharparenright}{\kern0pt}{\isacharparenright}{\kern0pt}\ {\isasymtriangleright}\ {\isasymemptyset}\ {\isasymturnstile}\ \isactrlbold {\isasymnot}\isaconst{ConRGA}}%
}%
\SNIP{ConRGA-Paper-Statements:thm:rga-goedel2-internal-0}{%
\isakeywordONE{thm}\isamarkupfalse%
\ rga{\isacharunderscore}{\kern0pt}goedel{\isadigit{2}}{\isacharunderscore}{\kern0pt}internal%
\isanewline%
}%
\SNIP{ConRGA-Paper-Statements:thm:rga-goedel2-internal-1}{%
\isa{\isaconst{ConRGA}\ {\isasymtriangleright}\ {\isasymemptyset}\ {\isasymturnstile}\ \isactrlbold {\isasymnot}\isactrlbold {\isasymnot}\isactrlbold A\ \isactrlbold {\isasymcdot}\ {\isacharparenleft}{\kern0pt}\isactrlbold {\isasymlambda}{\isacharparenleft}{\kern0pt}\isactrlbold {\isasymnot}\isaconst{iChk}\ \isactrlbold v\isactrlsub {\isadigit{0}}\ {\isacharparenleft}{\kern0pt}{\isasymupharpoonleft}{\isacharparenleft}{\kern0pt}\isactrlbold n{\isasymlbrace}\isaconst{ConRGA}{\isasymrbrace}{\isacharparenright}{\kern0pt}{\isacharparenright}{\kern0pt}{\isacharparenright}{\kern0pt}{\isacharparenright}{\kern0pt}}%
}%
\SNIP{ConRGA-Paper-Statements:thm:rga-goedel2-sentence-0}{%
\isakeywordONE{thm}\isamarkupfalse%
\ rga{\isacharunderscore}{\kern0pt}goedel{\isadigit{2}}{\isacharunderscore}{\kern0pt}sentence%
\isanewline%
}%
\SNIP{ConRGA-Paper-Statements:thm:rga-goedel2-sentence-1}{%
\isa{{\isasymemptyset}\ {\isasymturnstile}\ \isactrlbold A\ \isactrlbold {\isasymcdot}\ {\isacharparenleft}{\kern0pt}\isactrlbold {\isasymlambda}{\isacharparenleft}{\kern0pt}\isaconst{iChk}\ \isactrlbold v\isactrlsub {\isadigit{0}}\ {\isacharparenleft}{\kern0pt}\isactrlbold n{\isasymlbrace}\isaconst{ConRGA}{\isasymrbrace}{\isacharparenright}{\kern0pt}\ \isactrlbold {\isasymlongrightarrow}\ \isactrlbold {\isasymnot}\isaconst{ConRGA}{\isacharparenright}{\kern0pt}{\isacharparenright}{\kern0pt}}%
}%
\SNIP{ConRGA-Paper-Statements:thm:goedel2T-sentence-0}{%
\isakeywordONE{thm}\isamarkupfalse%
\ goedel{\isadigit{2}}T{\isacharunderscore}{\kern0pt}sentence%
\isanewline%
}%
\SNIP{ConRGA-Paper-Statements:thm:goedel2T-sentence-1}{%
\isa{\isavar{Q}\ {\isasymP}{\isasymlbrakk}\isaconst{ConRGA}\ {\isasymtriangleright}\ {\isasymemptyset}\ {\isasymturnstile}\ \isaconst{goedelG}{\isasymrbrakk}\ {\isasymLongrightarrow}\ {\isasymemptyset}\ {\isasymturnstile}\ \isactrlbold A\ \isactrlbold {\isasymcdot}\ {\isacharparenleft}{\kern0pt}\isactrlbold {\isasymlambda}{\isacharparenleft}{\kern0pt}\isaconst{iChk}\ \isactrlbold v\isactrlsub {\isadigit{0}}\ {\isacharparenleft}{\kern0pt}\isactrlbold n{\isasymlbrace}\isaconst{ConRGA}{\isasymrbrace}{\isacharparenright}{\kern0pt}\ \isactrlbold {\isasymlongrightarrow}\ \isaconst{iChk}\ {\isacharparenleft}{\kern0pt}\isaconst{goedel{\isadigit{2}}T}\ \isavar{Q}\ \isactrlbold {\isasymcdot}\ \isactrlbold v\isactrlsub {\isadigit{0}}{\isacharparenright}{\kern0pt}\ \isaconst{cEq{\isadigit{0}}{\isadigit{1}}}{\isacharparenright}{\kern0pt}{\isacharparenright}{\kern0pt}}%
}%
\SNIP{ConRGA-Paper-Statements:thm:rga-goedel2-sentence-witness-0}{%
\isakeywordONE{thm}\isamarkupfalse%
\ rga{\isacharunderscore}{\kern0pt}goedel{\isadigit{2}}{\isacharunderscore}{\kern0pt}sentence{\isacharunderscore}{\kern0pt}witness%
\isanewline%
}%
\SNIP{ConRGA-Paper-Statements:thm:rga-goedel2-sentence-witness-1}{%
\isa{{\isasymexists}\isabound{T}{\isachardot}{\kern0pt}\ \isaconst{Syntax{\isacharunderscore}{\kern0pt}Term{\isachardot}{\kern0pt}closed}\ \isabound{T}\ {\isasymand}\ {\isasymemptyset}\ {\isasymturnstile}\ \isactrlbold A\ \isactrlbold {\isasymcdot}\ {\isacharparenleft}{\kern0pt}\isactrlbold {\isasymlambda}{\isacharparenleft}{\kern0pt}\isaconst{iChk}\ \isactrlbold v\isactrlsub {\isadigit{0}}\ {\isacharparenleft}{\kern0pt}\isactrlbold n{\isasymlbrace}\isaconst{ConRGA}{\isasymrbrace}{\isacharparenright}{\kern0pt}\ \isactrlbold {\isasymlongrightarrow}\ \isaconst{iChk}\ {\isacharparenleft}{\kern0pt}\isabound{T}\ \isactrlbold {\isasymcdot}\ \isactrlbold v\isactrlsub {\isadigit{0}}{\isacharparenright}{\kern0pt}\ \isaconst{cEq{\isadigit{0}}{\isadigit{1}}}{\isacharparenright}{\kern0pt}{\isacharparenright}{\kern0pt}}%
}%
\SNIP{ConRGA-Paper-Statements:thm:bs-goedel2-0}{%
\isakeywordONE{thm}\isamarkupfalse%
\ bs{\isacharunderscore}{\kern0pt}goedel{\isadigit{2}}%
\isanewline%
}%
\SNIP{ConRGA-Paper-Statements:thm:bs-goedel2-1}{%
\isa{\isactrlbold {\isasymnot}\isaconst{goedelG}\ {\isasymtriangleright}\ {\isasymemptyset}\ {\isasymturnstile}\ \isactrlbold {\isasymnot}\isactrlbold A\ \isactrlbold {\isasymcdot}\ {\isacharparenleft}{\kern0pt}\isactrlbold {\isasymlambda}{\isacharparenleft}{\kern0pt}\isactrlbold {\isasymnot}\isaconst{iChk}\ \isactrlbold v\isactrlsub {\isadigit{0}}\ {\isacharparenleft}{\kern0pt}{\isasymupharpoonleft}{\isacharparenleft}{\kern0pt}\isactrlbold n{\isasymlbrace}\isactrlbold {\isasymnot}\isactrlbold {\isadigit{0}}\ \isactrlbold {\isacharequal}{\kern0pt}\ \isactrlbold {\isadigit{0}}{\isasymrbrace}{\isacharparenright}{\kern0pt}{\isacharparenright}{\kern0pt}{\isacharparenright}{\kern0pt}{\isacharparenright}{\kern0pt}\ {\isasymand}\ \isactrlbold {\isasymnot}\isactrlbold A\ \isactrlbold {\isasymcdot}\ {\isacharparenleft}{\kern0pt}\isactrlbold {\isasymlambda}{\isacharparenleft}{\kern0pt}\isactrlbold {\isasymnot}\isaconst{iChk}\ \isactrlbold v\isactrlsub {\isadigit{0}}\ {\isacharparenleft}{\kern0pt}{\isasymupharpoonleft}{\isacharparenleft}{\kern0pt}\isactrlbold n{\isasymlbrace}\isactrlbold {\isasymnot}\isactrlbold {\isasymnot}\isactrlbold A\ \isactrlbold {\isasymcdot}\ {\isacharparenleft}{\kern0pt}\isactrlbold {\isasymlambda}{\isacharparenleft}{\kern0pt}\isactrlbold {\isasymnot}\isaconst{iChk}\ \isactrlbold v\isactrlsub {\isadigit{0}}\ {\isacharparenleft}{\kern0pt}{\isasymupharpoonleft}{\isacharparenleft}{\kern0pt}\isactrlbold n{\isasymlbrace}\isactrlbold {\isasymnot}\isactrlbold {\isadigit{0}}\ \isactrlbold {\isacharequal}{\kern0pt}\ \isactrlbold {\isadigit{0}}{\isasymrbrace}{\isacharparenright}{\kern0pt}{\isacharparenright}{\kern0pt}{\isacharparenright}{\kern0pt}{\isacharparenright}{\kern0pt}{\isasymrbrace}{\isacharparenright}{\kern0pt}{\isacharparenright}{\kern0pt}{\isacharparenright}{\kern0pt}{\isacharparenright}{\kern0pt}\ {\isasymtriangleright}\ {\isasymemptyset}\ {\isasymturnstile}\ \isactrlbold {\isasymnot}\isactrlbold A\ \isactrlbold {\isasymcdot}\ {\isacharparenleft}{\kern0pt}\isactrlbold {\isasymlambda}{\isacharparenleft}{\kern0pt}\isactrlbold {\isasymnot}\isaconst{iChk}\ \isactrlbold v\isactrlsub {\isadigit{0}}\ {\isacharparenleft}{\kern0pt}{\isasymupharpoonleft}{\isacharparenleft}{\kern0pt}\isactrlbold n{\isasymlbrace}\isactrlbold {\isasymnot}\isactrlbold {\isadigit{0}}\ \isactrlbold {\isacharequal}{\kern0pt}\ \isactrlbold {\isadigit{0}}{\isasymrbrace}{\isacharparenright}{\kern0pt}{\isacharparenright}{\kern0pt}{\isacharparenright}{\kern0pt}{\isacharparenright}{\kern0pt}\ {\isasymand}\ \isactrlbold {\isasymnot}\isactrlbold A\ \isactrlbold {\isasymcdot}\ {\isacharparenleft}{\kern0pt}\isactrlbold {\isasymlambda}{\isacharparenleft}{\kern0pt}\isactrlbold {\isasymnot}\isaconst{iChk}\ \isactrlbold v\isactrlsub {\isadigit{0}}\ {\isacharparenleft}{\kern0pt}{\isasymupharpoonleft}{\isacharparenleft}{\kern0pt}\isactrlbold n{\isasymlbrace}\isactrlbold {\isasymnot}\isactrlbold {\isadigit{0}}\ \isactrlbold {\isacharequal}{\kern0pt}\ \isactrlbold {\isadigit{0}}{\isasymrbrace}{\isacharparenright}{\kern0pt}{\isacharparenright}{\kern0pt}{\isacharparenright}{\kern0pt}{\isacharparenright}{\kern0pt}\ {\isasymtriangleright}\ {\isasymemptyset}\ {\isasymturnstile}\ \isactrlbold {\isasymnot}\isaconst{goedelG}}%
}%
\SNIP{ConRGA-Paper-Statements:thm:rga-con-qf-int-0}{%
\isakeywordONE{thm}\isamarkupfalse%
\ rga{\isacharunderscore}{\kern0pt}con{\isacharunderscore}{\kern0pt}qf{\isacharunderscore}{\kern0pt}int%
\isanewline%
}%
\SNIP{ConRGA-Paper-Statements:thm:rga-con-qf-int-1}{%
\isa{{\isasymemptyset}\ {\isasymturnstile}\ \isactrlbold A\ \isactrlbold {\isasymcdot}\ {\isacharparenleft}{\kern0pt}\isactrlbold {\isasymlambda}{\isacharparenleft}{\kern0pt}\isactrlbold {\isasymnot}\isaconst{iChkQf}\ \isactrlbold v\isactrlsub {\isadigit{0}}\ \isaconst{cEq{\isadigit{0}}{\isadigit{1}}}{\isacharparenright}{\kern0pt}{\isacharparenright}{\kern0pt}}%
}%
\SNIP{ConRGA-Paper-Statements:thm:rga-con-frag-int-0}{%
\isakeywordONE{thm}\isamarkupfalse%
\ rga{\isacharunderscore}{\kern0pt}con{\isacharunderscore}{\kern0pt}frag{\isacharunderscore}{\kern0pt}int%
\isanewline%
}%
\SNIP{ConRGA-Paper-Statements:thm:rga-con-frag-int-1}{%
\isa{{\isasymemptyset}\ {\isasymturnstile}\ \isactrlbold A\ \isactrlbold {\isasymcdot}\ {\isacharparenleft}{\kern0pt}\isactrlbold {\isasymlambda}{\isacharparenleft}{\kern0pt}\isaconst{iStat}\ \isactrlbold v\isactrlsub {\isadigit{0}}\ \isactrlbold {\isasymand}\ \isaconst{iFragV}\ \isactrlbold v\isactrlsub {\isadigit{0}}\ \isactrlbold {\isasymlongrightarrow}\ \isactrlbold {\isasymnot}\isaconst{iChk}\ \isactrlbold v\isactrlsub {\isadigit{0}}\ \isaconst{cEq{\isadigit{0}}{\isadigit{1}}}{\isacharparenright}{\kern0pt}{\isacharparenright}{\kern0pt}}%
}%
\SNIP{ConRGA-Paper-Statements:thm:rga-con-lvi-under-lf-all-0}{%
\isakeywordONE{thm}\isamarkupfalse%
\ rga{\isacharunderscore}{\kern0pt}con{\isacharunderscore}{\kern0pt}lvi{\isacharunderscore}{\kern0pt}under{\isacharunderscore}{\kern0pt}lf{\isacharunderscore}{\kern0pt}all%
\isanewline%
}%
\SNIP{ConRGA-Paper-Statements:thm:rga-con-lvi-under-lf-all-1}{%
\isa{\isaconst{iLFI}\ {\isasymtriangleright}\ {\isasymemptyset}\ {\isasymturnstile}\ \isactrlbold A\ \isactrlbold {\isasymcdot}\ {\isacharparenleft}{\kern0pt}\isactrlbold {\isasymlambda}{\isacharparenleft}{\kern0pt}\isaconst{iLvClassI}\ \isactrlbold v\isactrlsub {\isadigit{0}}\ \isactrlbold {\isasymlongrightarrow}\ \isactrlbold {\isasymnot}\isaconst{iChk}\ \isactrlbold v\isactrlsub {\isadigit{0}}\ \isaconst{cEq{\isadigit{0}}{\isadigit{1}}}{\isacharparenright}{\kern0pt}{\isacharparenright}{\kern0pt}}%
}%
\SNIP{ConRGA-Paper-Statements:thm:rga-con-lvi-oga-0}{%
\isakeywordONE{thm}\isamarkupfalse%
\ rga{\isacharunderscore}{\kern0pt}con{\isacharunderscore}{\kern0pt}lvi{\isacharunderscore}{\kern0pt}oga%
\isanewline%
}%
\SNIP{ConRGA-Paper-Statements:thm:rga-con-lvi-oga-1}{%
\isa{\isaconst{prf{\isacharunderscore}{\kern0pt}derives{\isachardot}{\kern0pt}entails}\ \isaconst{oga{\isachardot}{\kern0pt}derives}\ {\isasymemptyset}\ {\isacharparenleft}{\kern0pt}\isactrlbold A\ \isactrlbold {\isasymcdot}\ {\isacharparenleft}{\kern0pt}\isactrlbold {\isasymlambda}{\isacharparenleft}{\kern0pt}\isaconst{iLvClassI}\ \isactrlbold v\isactrlsub {\isadigit{0}}\ \isactrlbold {\isasymlongrightarrow}\ \isactrlbold {\isasymnot}\isaconst{iChk}\ \isactrlbold v\isactrlsub {\isadigit{0}}\ \isaconst{cEq{\isadigit{0}}{\isadigit{1}}}{\isacharparenright}{\kern0pt}{\isacharparenright}{\kern0pt}{\isacharparenright}{\kern0pt}}%
}%
\SNIP{ConRGA-Paper-Statements:thm:rga-con-under-res-0}{%
\isakeywordONE{thm}\isamarkupfalse%
\ rga{\isacharunderscore}{\kern0pt}con{\isacharunderscore}{\kern0pt}under{\isacharunderscore}{\kern0pt}res%
\isanewline%
}%
\SNIP{ConRGA-Paper-Statements:thm:rga-con-under-res-1}{%
\isa{\isaconst{iRes}\ \isactrlbold v\isactrlsub {\isadigit{0}}\ {\isasymtriangleright}\ \isactrlbold v\isactrlsub {\isadigit{0}}\ \isactrlbold N\ {\isasymtriangleright}\ {\isasymemptyset}\ {\isasymturnstile}\ \isactrlbold {\isasymnot}\isaconst{iChk}\ \isactrlbold v\isactrlsub {\isadigit{0}}\ \isaconst{cEq{\isadigit{0}}{\isadigit{1}}}}%
}%
\SNIP{ConRGA-Paper-Statements:thm:ConRGA-def-0}{%
\isakeywordONE{thm}\isamarkupfalse%
\ ConRGA{\isacharunderscore}{\kern0pt}def%
\isanewline%
}%
\SNIP{ConRGA-Paper-Statements:thm:ConRGA-def-1}{%
\isa{\isaconst{ConRGA}\ {\isasymequiv}\ \isactrlbold A\ \isactrlbold {\isasymcdot}\ {\isacharparenleft}{\kern0pt}\isactrlbold {\isasymlambda}{\isacharparenleft}{\kern0pt}\isactrlbold {\isasymnot}\isaconst{iChk}\ \isactrlbold v\isactrlsub {\isadigit{0}}\ \isaconst{cEq{\isadigit{0}}{\isadigit{1}}}{\isacharparenright}{\kern0pt}{\isacharparenright}{\kern0pt}}%
}%
\SNIP{ConRGA-Paper-Statements:thm:iChk-def-0}{%
\isakeywordONE{thm}\isamarkupfalse%
\ iChk{\isacharunderscore}{\kern0pt}def%
\isanewline%
}%
\SNIP{ConRGA-Paper-Statements:thm:iChk-def-1}{%
\isa{\isaconst{iChk}\ \isavar{p}\ \isavar{c}\ {\isasymequiv}\ {\isacharparenleft}{\kern0pt}\isaconst{cmp}\ \isaconst{validChkIdx}\ \isactrlbold {\isasymcdot}\ \isavar{p}\ \isactrlbold {\isasymcdot}\ \isavar{p}\ \isactrlbold {\isacharequal}{\kern0pt}\ \isactrlbold n\isaconst{Suc}\ {\isadigit{0}}\ \isactrlbold {\isasymand}\ \isavar{p}\ \isactrlbold {\isasymnoteq}\ \isactrlbold {\isadigit{0}}{\isacharparenright}{\kern0pt}\ \isactrlbold {\isasymand}\ \isaconst{cmp}\ \isaconst{length{\isacharunderscore}{\kern0pt}encIdx}\ \isactrlbold {\isasymcdot}\ \isaconst{cpFst}\ {\isacharparenleft}{\kern0pt}\isaconst{cpFst}\ {\isacharparenleft}{\kern0pt}\isaconst{trmPredR}\ \isavar{p}{\isacharparenright}{\kern0pt}{\isacharparenright}{\kern0pt}\ \isactrlbold {\isasymcdot}\ \isaconst{cpFst}\ {\isacharparenleft}{\kern0pt}\isaconst{cpFst}\ {\isacharparenleft}{\kern0pt}\isaconst{trmPredR}\ \isavar{p}{\isacharparenright}{\kern0pt}{\isacharparenright}{\kern0pt}\ \isactrlbold {\isacharequal}{\kern0pt}\ \isactrlbold {\isadigit{0}}\ \isactrlbold {\isasymand}\ \isaconst{cpSnd}\ {\isacharparenleft}{\kern0pt}\isaconst{cpFst}\ {\isacharparenleft}{\kern0pt}\isaconst{trmPredR}\ \isavar{p}{\isacharparenright}{\kern0pt}{\isacharparenright}{\kern0pt}\ \isactrlbold {\isacharequal}{\kern0pt}\ \isavar{c}}%
}%
\SNIP{ConRGA-Paper-Statements:thm:iWtrue-def-0}{%
\isakeywordONE{thm}\isamarkupfalse%
\ iWtrue{\isacharunderscore}{\kern0pt}def%
\isanewline%
}%
\SNIP{ConRGA-Paper-Statements:thm:iWtrue-def-1}{%
\isa{\isaconst{iWtrue}\ \isavar{c}\ {\isasymequiv}\ \isactrlbold {\isasymnot}\isactrlbold A\ \isactrlbold {\isasymcdot}\ {\isacharparenleft}{\kern0pt}\isactrlbold {\isasymlambda}{\isacharparenleft}{\kern0pt}\isactrlbold {\isasymnot}\isaconst{iPrDecTrueAt}\ {\isacharparenleft}{\kern0pt}{\isasymupharpoonleft}\isavar{c}{\isacharparenright}{\kern0pt}\ \isactrlbold v\isactrlsub {\isadigit{0}}{\isacharparenright}{\kern0pt}{\isacharparenright}{\kern0pt}}%
}%
\SNIP{ConRGA-Paper-Statements:thm:iPrDecTrueAt-def-0}{%
\isakeywordONE{thm}\isamarkupfalse%
\ iPrDecTrueAt{\isacharunderscore}{\kern0pt}def%
\isanewline%
}%
\SNIP{ConRGA-Paper-Statements:thm:iPrDecTrueAt-def-1}{%
\isa{\isaconst{iPrDecTrueAt}\ \isavar{c}\ \isavar{s}\ {\isasymequiv}\ \isaconst{gapSub}\ {\isacharparenleft}{\kern0pt}\isaconst{gapSub}\ {\isacharparenleft}{\kern0pt}\isactrlbold n\isaconst{Suc}\ {\isadigit{0}}{\isacharparenright}{\kern0pt}\ {\isacharparenleft}{\kern0pt}\isaconst{gapSub}\ {\isacharparenleft}{\kern0pt}\isaconst{iPrDecV}\ \isavar{s}\ \isavar{c}{\isacharparenright}{\kern0pt}\ {\isacharparenleft}{\kern0pt}\isactrlbold n{\isadigit{2}}{\isacharparenright}{\kern0pt}{\isacharparenright}{\kern0pt}{\isacharparenright}{\kern0pt}\ {\isacharparenleft}{\kern0pt}\isaconst{gapSub}\ {\isacharparenleft}{\kern0pt}\isactrlbold n{\isadigit{2}}{\isacharparenright}{\kern0pt}\ {\isacharparenleft}{\kern0pt}\isaconst{iPrDecV}\ \isavar{s}\ \isavar{c}{\isacharparenright}{\kern0pt}{\isacharparenright}{\kern0pt}}%
}%
\SNIP{ConRGA-Paper-Statements:thm:cEq01-def-0}{%
\isakeywordONE{thm}\isamarkupfalse%
\ cEq{\isadigit{0}}{\isadigit{1}}{\isacharunderscore}{\kern0pt}def%
\isanewline%
}%
\SNIP{ConRGA-Paper-Statements:thm:cEq01-def-1}{%
\isa{\isaconst{cEq{\isadigit{0}}{\isadigit{1}}}\ {\isasymequiv}\ \isactrlbold n{\isasymlbrace}\isactrlbold {\isadigit{0}}\ \isactrlbold {\isacharequal}{\kern0pt}\ \isactrlbold n\isaconst{Suc}\ {\isadigit{0}}{\isasymrbrace}}%
}%
\SNIP{ConRGA-Paper-Statements:thm:qfm-tin-0}{%
\isakeywordONE{thm}\isamarkupfalse%
\ qfm{\isacharunderscore}{\kern0pt}tin%
\isanewline%
}%
\SNIP{ConRGA-Paper-Statements:thm:qfm-tin-1}{%
\isa{{\isasymlbrakk}\isaconst{qfm}\ \isavar{p}{\isacharsemicolon}{\kern0pt}\ \isavar{{\isasymGamma}}\ {\isasymturnstile}\ \isavar{v}\ \isactrlbold N{\isacharsemicolon}{\kern0pt}\ \isavar{{\isasymGamma}}\ {\isasymturnstile}\ {\isacharbrackleft}{\kern0pt}{\isadigit{0}}{\isasymZpinj}\isavar{v}{\isacharbrackright}{\kern0pt}\isavar{p}{\isasymrbrakk}\ {\isasymLongrightarrow}\ \isavar{{\isasymGamma}}\ {\isasymturnstile}\ \isaconst{iWtrue}\ {\isacharparenleft}{\kern0pt}\isaconst{iSubM}\ {\isacharparenleft}{\kern0pt}\isaconst{updSubM}\ \isactrlbold {\isadigit{0}}\ {\isacharparenleft}{\kern0pt}\isactrlbold n{\isadigit{0}}{\isacharparenright}{\kern0pt}\ {\isacharparenleft}{\kern0pt}\isaconst{cmp}\ \isaconst{num{\isacharunderscore}{\kern0pt}codeIdx}\ \isactrlbold {\isasymcdot}\ \isactrlbold {\isadigit{0}}\ \isactrlbold {\isasymcdot}\ \isavar{v}{\isacharparenright}{\kern0pt}{\isacharparenright}{\kern0pt}\ {\isacharparenleft}{\kern0pt}\isactrlbold n{\isasymlbrace}\isavar{p}{\isasymrbrace}{\isacharparenright}{\kern0pt}{\isacharparenright}{\kern0pt}}%
}%
\SNIP{ConRGA-Paper-Statements:thm:qfm-tout-0}{%
\isakeywordONE{thm}\isamarkupfalse%
\ qfm{\isacharunderscore}{\kern0pt}tout%
\isanewline%
}%
\SNIP{ConRGA-Paper-Statements:thm:qfm-tout-1}{%
\isa{{\isasymlbrakk}\isaconst{qfm}\ \isavar{p}{\isacharsemicolon}{\kern0pt}\ \isavar{{\isasymGamma}}\ {\isasymturnstile}\ \isavar{v}\ \isactrlbold N{\isacharsemicolon}{\kern0pt}\ \isavar{{\isasymGamma}}\ {\isasymturnstile}\ \isaconst{iWtrue}\ {\isacharparenleft}{\kern0pt}\isaconst{iSubM}\ {\isacharparenleft}{\kern0pt}\isaconst{updSubM}\ \isactrlbold {\isadigit{0}}\ {\isacharparenleft}{\kern0pt}\isactrlbold n{\isadigit{0}}{\isacharparenright}{\kern0pt}\ {\isacharparenleft}{\kern0pt}\isaconst{cmp}\ \isaconst{num{\isacharunderscore}{\kern0pt}codeIdx}\ \isactrlbold {\isasymcdot}\ \isactrlbold {\isadigit{0}}\ \isactrlbold {\isasymcdot}\ \isavar{v}{\isacharparenright}{\kern0pt}{\isacharparenright}{\kern0pt}\ {\isacharparenleft}{\kern0pt}\isactrlbold n{\isasymlbrace}\isavar{p}{\isasymrbrace}{\isacharparenright}{\kern0pt}{\isacharparenright}{\kern0pt}{\isasymrbrakk}\ {\isasymLongrightarrow}\ \isavar{{\isasymGamma}}\ {\isasymturnstile}\ {\isacharbrackleft}{\kern0pt}{\isadigit{0}}{\isasymZpinj}\isavar{v}{\isacharbrackright}{\kern0pt}\isavar{p}}%
}%
\SNIP{ConRGA-Paper-Statements:thm:ptm-tout-eq-0}{%
\isakeywordONE{thm}\isamarkupfalse%
\ ptm{\isacharunderscore}{\kern0pt}tout{\isacharunderscore}{\kern0pt}eq%
\isanewline%
}%
\SNIP{ConRGA-Paper-Statements:thm:ptm-tout-eq-1}{%
\isa{{\isasymlbrakk}\isaconst{ptm}\ \isavar{a}{\isacharsemicolon}{\kern0pt}\ \isaconst{ptm}\ \isavar{b}{\isacharsemicolon}{\kern0pt}\ \isavar{{\isasymGamma}}\ {\isasymturnstile}\ \isavar{v}\ \isactrlbold N{\isacharsemicolon}{\kern0pt}\ \isavar{{\isasymGamma}}\ {\isasymturnstile}\ \isaconst{iWtrue}\ {\isacharparenleft}{\kern0pt}\isaconst{iSubM}\ {\isacharparenleft}{\kern0pt}\isaconst{updSubM}\ \isactrlbold {\isadigit{0}}\ {\isacharparenleft}{\kern0pt}\isactrlbold n{\isadigit{0}}{\isacharparenright}{\kern0pt}\ {\isacharparenleft}{\kern0pt}\isaconst{cmp}\ \isaconst{num{\isacharunderscore}{\kern0pt}codeIdx}\ \isactrlbold {\isasymcdot}\ \isactrlbold {\isadigit{0}}\ \isactrlbold {\isasymcdot}\ \isavar{v}{\isacharparenright}{\kern0pt}{\isacharparenright}{\kern0pt}\ {\isacharparenleft}{\kern0pt}\isactrlbold n{\isasymlbrace}\isavar{a}\ \isactrlbold {\isacharequal}{\kern0pt}\ \isavar{b}{\isasymrbrace}{\isacharparenright}{\kern0pt}{\isacharparenright}{\kern0pt}{\isasymrbrakk}\ {\isasymLongrightarrow}\ \isavar{{\isasymGamma}}\ {\isasymturnstile}\ {\isacharbrackleft}{\kern0pt}{\isadigit{0}}{\isasymZpinj}\isavar{v}{\isacharbrackright}{\kern0pt}\isavar{a}\ \isactrlbold {\isacharequal}{\kern0pt}\ {\isacharbrackleft}{\kern0pt}{\isadigit{0}}{\isasymZpinj}\isavar{v}{\isacharbrackright}{\kern0pt}\isavar{b}}%
}%
\SNIP{ConRGA-Paper-Statements:thm:iPrDecVal-app-open-0}{%
\isakeywordONE{thm}\isamarkupfalse%
\ iPrDecVal{\isacharunderscore}{\kern0pt}app{\isacharunderscore}{\kern0pt}open%
\isanewline%
}%
\SNIP{ConRGA-Paper-Statements:thm:iPrDecVal-app-open-1}{%
\isa{{\isasymlbrakk}\isavar{{\isasymGamma}}\ {\isasymturnstile}\ \isavar{c{\isadigit{1}}{\isachardot}{\kern0pt}{\isadigit{0}}}\ \isactrlbold N{\isacharsemicolon}{\kern0pt}\ \isavar{{\isasymGamma}}\ {\isasymturnstile}\ \isavar{c{\isadigit{2}}{\isachardot}{\kern0pt}{\isadigit{0}}}\ \isactrlbold N{\isacharsemicolon}{\kern0pt}\ \isavar{{\isasymGamma}}\ {\isasymturnstile}\ \isaconst{iWft}\ \isavar{c{\isadigit{1}}{\isachardot}{\kern0pt}{\isadigit{0}}}\ \isactrlbold {\isasymnoteq}\ \isactrlbold {\isadigit{0}}{\isacharsemicolon}{\kern0pt}\ \isavar{{\isasymGamma}}\ {\isasymturnstile}\ \isaconst{iWft}\ \isavar{c{\isadigit{2}}{\isachardot}{\kern0pt}{\isadigit{0}}}\ \isactrlbold {\isasymnoteq}\ \isactrlbold {\isadigit{0}}{\isacharsemicolon}{\kern0pt}\ \isavar{{\isasymGamma}}\ {\isasymturnstile}\ \isaconst{iIsCode}\ \isavar{c{\isadigit{1}}{\isachardot}{\kern0pt}{\isadigit{0}}}\ \isactrlbold {\isasymnoteq}\ \isactrlbold {\isadigit{0}}{\isacharsemicolon}{\kern0pt}\ \isavar{{\isasymGamma}}\ {\isasymturnstile}\ \isaconst{iIsCode}\ \isavar{c{\isadigit{2}}{\isachardot}{\kern0pt}{\isadigit{0}}}\ \isactrlbold {\isasymnoteq}\ \isactrlbold {\isadigit{0}}{\isacharsemicolon}{\kern0pt}\ \isavar{{\isasymGamma}}\ {\isasymturnstile}\ \isaconst{iQf}\ \isavar{c{\isadigit{1}}{\isachardot}{\kern0pt}{\isadigit{0}}}\ \isactrlbold {\isasymnoteq}\ \isactrlbold {\isadigit{0}}{\isacharsemicolon}{\kern0pt}\ \isavar{{\isasymGamma}}\ {\isasymturnstile}\ \isaconst{iQf}\ \isavar{c{\isadigit{2}}{\isachardot}{\kern0pt}{\isadigit{0}}}\ \isactrlbold {\isasymnoteq}\ \isactrlbold {\isadigit{0}}{\isacharsemicolon}{\kern0pt}\ \isavar{{\isasymGamma}}\ {\isasymturnstile}\ \isavar{x}\ \isactrlbold N{\isacharsemicolon}{\kern0pt}\ \isavar{{\isasymGamma}}\ {\isasymturnstile}\ \isavar{y}\ \isactrlbold N{\isacharsemicolon}{\kern0pt}\ \isavar{{\isasymGamma}}\ {\isasymturnstile}\ \isaconst{iPrDecVal}\ \isavar{c{\isadigit{1}}{\isachardot}{\kern0pt}{\isadigit{0}}}\ \isavar{x}{\isacharsemicolon}{\kern0pt}\ \isavar{{\isasymGamma}}\ {\isasymturnstile}\ \isaconst{iPrDecVal}\ \isavar{c{\isadigit{2}}{\isachardot}{\kern0pt}{\isadigit{0}}}\ \isavar{y}{\isasymrbrakk}\ {\isasymLongrightarrow}\ \isavar{{\isasymGamma}}\ {\isasymturnstile}\ \isaconst{iPrDecVal}\ {\isacharparenleft}{\kern0pt}\isaconst{mkAppS}\ {\isacharparenleft}{\kern0pt}\isaconst{mkAppS}\ {\isacharparenleft}{\kern0pt}\isactrlbold n{\isasymlbrace}\isaconst{cmp}\ \isavar{f}{\isasymrbrace}{\isacharparenright}{\kern0pt}\ \isavar{c{\isadigit{1}}{\isachardot}{\kern0pt}{\isadigit{0}}}{\isacharparenright}{\kern0pt}\ \isavar{c{\isadigit{2}}{\isachardot}{\kern0pt}{\isadigit{0}}}{\isacharparenright}{\kern0pt}\ {\isacharparenleft}{\kern0pt}\isaconst{cmp}\ \isavar{f}\ \isactrlbold {\isasymcdot}\ \isavar{x}\ \isactrlbold {\isasymcdot}\ \isavar{y}{\isacharparenright}{\kern0pt}}%
}%
\SNIP{ConRGA-Paper-Statements:thm:iPrDecVal-app-img-0}{%
\isakeywordONE{thm}\isamarkupfalse%
\ iPrDecVal{\isacharunderscore}{\kern0pt}app{\isacharunderscore}{\kern0pt}img%
\isanewline%
}%
\SNIP{ConRGA-Paper-Statements:thm:iPrDecVal-app-img-1}{%
\isa{{\isasymlbrakk}\isavar{{\isasymGamma}}\ {\isasymturnstile}\ \isavar{v}\ \isactrlbold N{\isacharsemicolon}{\kern0pt}\ \isavar{{\isasymGamma}}\ {\isasymturnstile}\ \isavar{a}\ \isactrlbold N{\isacharsemicolon}{\kern0pt}\ \isavar{{\isasymGamma}}\ {\isasymturnstile}\ \isavar{b}\ \isactrlbold N{\isacharsemicolon}{\kern0pt}\ \isavar{{\isasymGamma}}\ {\isasymturnstile}\ \isaconst{iWft}\ \isavar{a}\ \isactrlbold {\isasymnoteq}\ \isactrlbold {\isadigit{0}}{\isacharsemicolon}{\kern0pt}\ \isavar{{\isasymGamma}}\ {\isasymturnstile}\ \isaconst{iWft}\ \isavar{b}\ \isactrlbold {\isasymnoteq}\ \isactrlbold {\isadigit{0}}{\isacharsemicolon}{\kern0pt}\ \isavar{{\isasymGamma}}\ {\isasymturnstile}\ \isaconst{iIsCode}\ \isavar{a}\ \isactrlbold {\isasymnoteq}\ \isactrlbold {\isadigit{0}}{\isacharsemicolon}{\kern0pt}\ \isavar{{\isasymGamma}}\ {\isasymturnstile}\ \isaconst{iIsCode}\ \isavar{b}\ \isactrlbold {\isasymnoteq}\ \isactrlbold {\isadigit{0}}{\isacharsemicolon}{\kern0pt}\ \isavar{{\isasymGamma}}\ {\isasymturnstile}\ \isaconst{iQf}\ \isavar{a}\ \isactrlbold {\isasymnoteq}\ \isactrlbold {\isadigit{0}}{\isacharsemicolon}{\kern0pt}\ \isavar{{\isasymGamma}}\ {\isasymturnstile}\ \isaconst{iQf}\ \isavar{b}\ \isactrlbold {\isasymnoteq}\ \isactrlbold {\isadigit{0}}{\isacharsemicolon}{\kern0pt}\ \isavar{{\isasymGamma}}\ {\isasymturnstile}\ \isavar{x}\ \isactrlbold N{\isacharsemicolon}{\kern0pt}\ \isavar{{\isasymGamma}}\ {\isasymturnstile}\ \isavar{y}\ \isactrlbold N{\isacharsemicolon}{\kern0pt}\ \isavar{{\isasymGamma}}\ {\isasymturnstile}\ \isaconst{iPrDecVal}\ {\isacharparenleft}{\kern0pt}\isaconst{iSubM}\ {\isacharparenleft}{\kern0pt}\isaconst{updSubM}\ \isactrlbold {\isadigit{0}}\ {\isacharparenleft}{\kern0pt}\isactrlbold n{\isadigit{0}}{\isacharparenright}{\kern0pt}\ {\isacharparenleft}{\kern0pt}\isaconst{cmp}\ \isaconst{num{\isacharunderscore}{\kern0pt}codeIdx}\ \isactrlbold {\isasymcdot}\ \isactrlbold {\isadigit{0}}\ \isactrlbold {\isasymcdot}\ \isavar{v}{\isacharparenright}{\kern0pt}{\isacharparenright}{\kern0pt}\ \isavar{a}{\isacharparenright}{\kern0pt}\ \isavar{x}{\isacharsemicolon}{\kern0pt}\ \isavar{{\isasymGamma}}\ {\isasymturnstile}\ \isaconst{iPrDecVal}\ {\isacharparenleft}{\kern0pt}\isaconst{iSubM}\ {\isacharparenleft}{\kern0pt}\isaconst{updSubM}\ \isactrlbold {\isadigit{0}}\ {\isacharparenleft}{\kern0pt}\isactrlbold n{\isadigit{0}}{\isacharparenright}{\kern0pt}\ {\isacharparenleft}{\kern0pt}\isaconst{cmp}\ \isaconst{num{\isacharunderscore}{\kern0pt}codeIdx}\ \isactrlbold {\isasymcdot}\ \isactrlbold {\isadigit{0}}\ \isactrlbold {\isasymcdot}\ \isavar{v}{\isacharparenright}{\kern0pt}{\isacharparenright}{\kern0pt}\ \isavar{b}{\isacharparenright}{\kern0pt}\ \isavar{y}{\isasymrbrakk}\ {\isasymLongrightarrow}\ \isavar{{\isasymGamma}}\ {\isasymturnstile}\ \isaconst{iPrDecVal}\ {\isacharparenleft}{\kern0pt}\isaconst{iSubM}\ {\isacharparenleft}{\kern0pt}\isaconst{updSubM}\ \isactrlbold {\isadigit{0}}\ {\isacharparenleft}{\kern0pt}\isactrlbold n{\isadigit{0}}{\isacharparenright}{\kern0pt}\ {\isacharparenleft}{\kern0pt}\isaconst{cmp}\ \isaconst{num{\isacharunderscore}{\kern0pt}codeIdx}\ \isactrlbold {\isasymcdot}\ \isactrlbold {\isadigit{0}}\ \isactrlbold {\isasymcdot}\ \isavar{v}{\isacharparenright}{\kern0pt}{\isacharparenright}{\kern0pt}\ {\isacharparenleft}{\kern0pt}\isaconst{mkAppS}\ {\isacharparenleft}{\kern0pt}\isaconst{mkAppS}\ {\isacharparenleft}{\kern0pt}\isactrlbold n{\isasymlbrace}\isaconst{cmp}\ \isavar{f}{\isasymrbrace}{\isacharparenright}{\kern0pt}\ \isavar{a}{\isacharparenright}{\kern0pt}\ \isavar{b}{\isacharparenright}{\kern0pt}{\isacharparenright}{\kern0pt}\ {\isacharparenleft}{\kern0pt}\isaconst{cmp}\ \isavar{f}\ \isactrlbold {\isasymcdot}\ \isavar{x}\ \isactrlbold {\isasymcdot}\ \isavar{y}{\isacharparenright}{\kern0pt}}%
}%
\SNIP{ConRGA-Paper-Statements:thm:iLift-zero-code-0}{%
\isakeywordONE{thm}\isamarkupfalse%
\ iLift{\isacharunderscore}{\kern0pt}zero{\isacharunderscore}{\kern0pt}code%
\isanewline%
}%
\SNIP{ConRGA-Paper-Statements:thm:iLift-zero-code-1}{%
\isa{{\isasymlbrakk}\isavar{{\isasymGamma}}\ {\isasymturnstile}\ \isavar{c}\ \isactrlbold N{\isacharsemicolon}{\kern0pt}\ \isavar{{\isasymGamma}}\ {\isasymturnstile}\ \isaconst{iWft}\ \isavar{c}\ \isactrlbold {\isasymnoteq}\ \isactrlbold {\isadigit{0}}{\isacharsemicolon}{\kern0pt}\ \isavar{{\isasymGamma}}\ {\isasymturnstile}\ \isaconst{iIsCode}\ \isavar{c}\ \isactrlbold {\isasymnoteq}\ \isactrlbold {\isadigit{0}}{\isasymrbrakk}\ {\isasymLongrightarrow}\ \isavar{{\isasymGamma}}\ {\isasymturnstile}\ \isaconst{iLift}\ \isactrlbold {\isadigit{0}}\ {\isacharparenleft}{\kern0pt}\isaconst{cpPair}\ \isactrlbold {\isadigit{0}}\ \isavar{c}{\isacharparenright}{\kern0pt}\ \isactrlbold {\isacharequal}{\kern0pt}\ \isavar{c}}%
}%
\SNIP{ConRGA-Paper-Statements:thm:iSubM-updSubM-hit-code-0}{%
\isakeywordONE{thm}\isamarkupfalse%
\ iSubM{\isacharunderscore}{\kern0pt}updSubM{\isacharunderscore}{\kern0pt}hit{\isacharunderscore}{\kern0pt}code%
\isanewline%
}%
\SNIP{ConRGA-Paper-Statements:thm:iSubM-updSubM-hit-code-1}{%
\isa{{\isasymlbrakk}\isavar{{\isasymGamma}}\ {\isasymturnstile}\ \isavar{z}\ \isactrlbold N{\isacharsemicolon}{\kern0pt}\ \isavar{{\isasymGamma}}\ {\isasymturnstile}\ \isavar{c}\ \isactrlbold N{\isacharsemicolon}{\kern0pt}\ \isavar{{\isasymGamma}}\ {\isasymturnstile}\ \isaconst{iWft}\ \isavar{c}\ \isactrlbold {\isasymnoteq}\ \isactrlbold {\isadigit{0}}{\isacharsemicolon}{\kern0pt}\ \isavar{{\isasymGamma}}\ {\isasymturnstile}\ \isaconst{iIsCode}\ \isavar{c}\ \isactrlbold {\isasymnoteq}\ \isactrlbold {\isadigit{0}}{\isasymrbrakk}\ {\isasymLongrightarrow}\ \isavar{{\isasymGamma}}\ {\isasymturnstile}\ \isaconst{iSubM}\ {\isacharparenleft}{\kern0pt}\isaconst{updSubM}\ \isavar{z}\ {\isacharparenleft}{\kern0pt}\isactrlbold n\isavar{i}{\isacharparenright}{\kern0pt}\ \isavar{c}{\isacharparenright}{\kern0pt}\ {\isacharparenleft}{\kern0pt}\isactrlbold n{\isasymlbrace}\isaconst{trmVar}\ \isavar{i}{\isasymrbrace}{\isacharparenright}{\kern0pt}\ \isactrlbold {\isacharequal}{\kern0pt}\ \isavar{c}}%
}%
\SNIP{ConRGA-Paper-Statements:thm:recT-unfold-0}{%
\isakeywordONE{thm}\isamarkupfalse%
\ recT{\isacharunderscore}{\kern0pt}unfold%
\isanewline%
}%
\SNIP{ConRGA-Paper-Statements:thm:recT-unfold-1}{%
\isa{{\isasymlbrakk}\isavar{{\isasymGamma}}\ {\isasymturnstile}\ \isavar{m}\ \isactrlbold N{\isacharsemicolon}{\kern0pt}\ \isavar{{\isasymGamma}}\ {\isasymturnstile}\ \isavar{n}\ \isactrlbold N{\isacharsemicolon}{\kern0pt}\ \isavar{{\isasymGamma}}\ {\isasymturnstile}\ \isavar{n}\ \isactrlbold {\isasymnoteq}\ \isactrlbold {\isadigit{0}}{\isacharsemicolon}{\kern0pt}\ \isavar{{\isasymGamma}}\ {\isasymturnstile}\ \isaconst{gapSub}\ {\isacharparenleft}{\kern0pt}\isactrlbold S\ \isaconst{cmp}\ \isavar{S{\isadigit{1}}{\isachardot}{\kern0pt}{\isadigit{0}}}\ \isactrlbold {\isasymcdot}\ \isavar{m}\ \isactrlbold {\isasymcdot}\ \isavar{n}{\isacharparenright}{\kern0pt}\ \isavar{n}\ \isactrlbold {\isacharequal}{\kern0pt}\ \isactrlbold {\isadigit{0}}{\isacharsemicolon}{\kern0pt}\ \isavar{{\isasymGamma}}\ {\isasymturnstile}\ \isaconst{gapSub}\ {\isacharparenleft}{\kern0pt}\isactrlbold S\ \isaconst{cmp}\ \isavar{S{\isadigit{2}}{\isachardot}{\kern0pt}{\isadigit{0}}}\ \isactrlbold {\isasymcdot}\ \isavar{m}\ \isactrlbold {\isasymcdot}\ \isavar{n}{\isacharparenright}{\kern0pt}\ \isavar{n}\ \isactrlbold {\isacharequal}{\kern0pt}\ \isactrlbold {\isadigit{0}}{\isasymrbrakk}\ {\isasymLongrightarrow}\ \isavar{{\isasymGamma}}\ {\isasymturnstile}\ \isaconst{lstHd}\ {\isacharparenleft}{\kern0pt}\isactrlbold R\ \isactrlbold {\isasymcdot}\ {\isacharparenleft}{\kern0pt}\isaconst{cmp}\ {\isacharparenleft}{\kern0pt}\isaconst{prSub}\ \isaconst{lstConsIdx}\ \isavar{B}\ {\isacharparenleft}{\kern0pt}\isaconst{prInit}\ \isaconst{rInit{\isadigit{0}}}{\isacharparenright}{\kern0pt}{\isacharparenright}{\kern0pt}\ \isactrlbold {\isasymcdot}\ \isavar{m}\ \isactrlbold {\isasymcdot}\ \isactrlbold {\isadigit{0}}{\isacharparenright}{\kern0pt}\ \isactrlbold {\isasymcdot}\ \isaconst{prRecH}\ {\isacharparenleft}{\kern0pt}\isaconst{prSub}\ \isaconst{lstConsIdx}\ {\isacharparenleft}{\kern0pt}\isaconst{recHIdx}\ \isavar{S{\isadigit{1}}{\isachardot}{\kern0pt}{\isadigit{0}}}\ \isavar{S{\isadigit{2}}{\isachardot}{\kern0pt}{\isadigit{0}}}\ \isavar{H}{\isacharparenright}{\kern0pt}\ {\isacharparenleft}{\kern0pt}\isaconst{prInit}\ \isaconst{rInitR}{\isacharparenright}{\kern0pt}{\isacharparenright}{\kern0pt}\ \isaconst{covStepIdx}\ \isavar{m}\ \isactrlbold {\isasymcdot}\ \isavar{n}{\isacharparenright}{\kern0pt}\ \isactrlbold {\isacharequal}{\kern0pt}\ \isaconst{cmp}\ \isavar{H}\ \isactrlbold {\isasymcdot}\ \isaconst{cpPair}\ \isavar{m}\ \isavar{n}\ \isactrlbold {\isasymcdot}\ \isaconst{cpPair}\ {\isacharparenleft}{\kern0pt}\isaconst{lstHd}\ {\isacharparenleft}{\kern0pt}\isactrlbold R\ \isactrlbold {\isasymcdot}\ {\isacharparenleft}{\kern0pt}\isaconst{cmp}\ {\isacharparenleft}{\kern0pt}\isaconst{prSub}\ \isaconst{lstConsIdx}\ \isavar{B}\ {\isacharparenleft}{\kern0pt}\isaconst{prInit}\ \isaconst{rInit{\isadigit{0}}}{\isacharparenright}{\kern0pt}{\isacharparenright}{\kern0pt}\ \isactrlbold {\isasymcdot}\ \isavar{m}\ \isactrlbold {\isasymcdot}\ \isactrlbold {\isadigit{0}}{\isacharparenright}{\kern0pt}\ \isactrlbold {\isasymcdot}\ \isaconst{prRecH}\ {\isacharparenleft}{\kern0pt}\isaconst{prSub}\ \isaconst{lstConsIdx}\ {\isacharparenleft}{\kern0pt}\isaconst{recHIdx}\ \isavar{S{\isadigit{1}}{\isachardot}{\kern0pt}{\isadigit{0}}}\ \isavar{S{\isadigit{2}}{\isachardot}{\kern0pt}{\isadigit{0}}}\ \isavar{H}{\isacharparenright}{\kern0pt}\ {\isacharparenleft}{\kern0pt}\isaconst{prInit}\ \isaconst{rInitR}{\isacharparenright}{\kern0pt}{\isacharparenright}{\kern0pt}\ \isaconst{covStepIdx}\ \isavar{m}\ \isactrlbold {\isasymcdot}\ {\isacharparenleft}{\kern0pt}\isaconst{cmp}\ \isavar{S{\isadigit{1}}{\isachardot}{\kern0pt}{\isadigit{0}}}\ \isactrlbold {\isasymcdot}\ \isavar{m}\ \isactrlbold {\isasymcdot}\ \isavar{n}{\isacharparenright}{\kern0pt}{\isacharparenright}{\kern0pt}{\isacharparenright}{\kern0pt}\ {\isacharparenleft}{\kern0pt}\isaconst{lstHd}\ {\isacharparenleft}{\kern0pt}\isactrlbold R\ \isactrlbold {\isasymcdot}\ {\isacharparenleft}{\kern0pt}\isaconst{cmp}\ {\isacharparenleft}{\kern0pt}\isaconst{prSub}\ \isaconst{lstConsIdx}\ \isavar{B}\ {\isacharparenleft}{\kern0pt}\isaconst{prInit}\ \isaconst{rInit{\isadigit{0}}}{\isacharparenright}{\kern0pt}{\isacharparenright}{\kern0pt}\ \isactrlbold {\isasymcdot}\ \isavar{m}\ \isactrlbold {\isasymcdot}\ \isactrlbold {\isadigit{0}}{\isacharparenright}{\kern0pt}\ \isactrlbold {\isasymcdot}\ \isaconst{prRecH}\ {\isacharparenleft}{\kern0pt}\isaconst{prSub}\ \isaconst{lstConsIdx}\ {\isacharparenleft}{\kern0pt}\isaconst{recHIdx}\ \isavar{S{\isadigit{1}}{\isachardot}{\kern0pt}{\isadigit{0}}}\ \isavar{S{\isadigit{2}}{\isachardot}{\kern0pt}{\isadigit{0}}}\ \isavar{H}{\isacharparenright}{\kern0pt}\ {\isacharparenleft}{\kern0pt}\isaconst{prInit}\ \isaconst{rInitR}{\isacharparenright}{\kern0pt}{\isacharparenright}{\kern0pt}\ \isaconst{covStepIdx}\ \isavar{m}\ \isactrlbold {\isasymcdot}\ {\isacharparenleft}{\kern0pt}\isaconst{cmp}\ \isavar{S{\isadigit{2}}{\isachardot}{\kern0pt}{\isadigit{0}}}\ \isactrlbold {\isasymcdot}\ \isavar{m}\ \isactrlbold {\isasymcdot}\ \isavar{n}{\isacharparenright}{\kern0pt}{\isacharparenright}{\kern0pt}{\isacharparenright}{\kern0pt}}%
}%
\SNIP{ConRGA-Paper-Statements:thm:ev-addP-0}{%
\isakeywordONE{thm}\isamarkupfalse%
\ ev{\isacharunderscore}{\kern0pt}addP%
\isanewline%
}%
\SNIP{ConRGA-Paper-Statements:thm:ev-addP-1}{%
\isa{{\isasymlbrakk}\isavar{{\isasymGamma}}\ {\isasymturnstile}\ \isavar{v}\ \isactrlbold N{\isacharsemicolon}{\kern0pt}\ \isavar{{\isasymGamma}}\ {\isasymturnstile}\ \isavar{p}\ \isactrlbold N{\isacharsemicolon}{\kern0pt}\ \isavar{{\isasymGamma}}\ {\isasymturnstile}\ \isavar{q}\ \isactrlbold N{\isasymrbrakk}\ {\isasymLongrightarrow}\ \isavar{{\isasymGamma}}\ {\isasymturnstile}\ \isaconst{cmp}\ \isaconst{evIdx}\ \isactrlbold {\isasymcdot}\ \isavar{v}\ \isactrlbold {\isasymcdot}\ {\isacharparenleft}{\kern0pt}\isaconst{cmp}\ \isaconst{plAddIdx}\ \isactrlbold {\isasymcdot}\ \isavar{p}\ \isactrlbold {\isasymcdot}\ \isavar{q}{\isacharparenright}{\kern0pt}\ \isactrlbold {\isacharequal}{\kern0pt}\ \isaconst{gapAdd}\ {\isacharparenleft}{\kern0pt}\isaconst{cmp}\ \isaconst{evIdx}\ \isactrlbold {\isasymcdot}\ \isavar{v}\ \isactrlbold {\isasymcdot}\ \isavar{p}{\isacharparenright}{\kern0pt}\ {\isacharparenleft}{\kern0pt}\isaconst{cmp}\ \isaconst{evIdx}\ \isactrlbold {\isasymcdot}\ \isavar{v}\ \isactrlbold {\isasymcdot}\ \isavar{q}{\isacharparenright}{\kern0pt}}%
}%
\SNIP{ConRGA-Paper-Statements:thm:ev-mulP-0}{%
\isakeywordONE{thm}\isamarkupfalse%
\ ev{\isacharunderscore}{\kern0pt}mulP%
\isanewline%
}%
\SNIP{ConRGA-Paper-Statements:thm:ev-mulP-1}{%
\isa{{\isasymlbrakk}\isavar{{\isasymGamma}}\ {\isasymturnstile}\ \isavar{v}\ \isactrlbold N{\isacharsemicolon}{\kern0pt}\ \isavar{{\isasymGamma}}\ {\isasymturnstile}\ \isavar{p}\ \isactrlbold N{\isacharsemicolon}{\kern0pt}\ \isavar{{\isasymGamma}}\ {\isasymturnstile}\ \isavar{q}\ \isactrlbold N{\isasymrbrakk}\ {\isasymLongrightarrow}\ \isavar{{\isasymGamma}}\ {\isasymturnstile}\ \isaconst{cmp}\ \isaconst{evIdx}\ \isactrlbold {\isasymcdot}\ \isavar{v}\ \isactrlbold {\isasymcdot}\ {\isacharparenleft}{\kern0pt}\isaconst{cmp}\ \isaconst{plMulIdx}\ \isactrlbold {\isasymcdot}\ \isavar{p}\ \isactrlbold {\isasymcdot}\ \isavar{q}{\isacharparenright}{\kern0pt}\ \isactrlbold {\isacharequal}{\kern0pt}\ \isaconst{cmp}\ \isaconst{mulIdx}\ \isactrlbold {\isasymcdot}\ {\isacharparenleft}{\kern0pt}\isaconst{cmp}\ \isaconst{evIdx}\ \isactrlbold {\isasymcdot}\ \isavar{v}\ \isactrlbold {\isasymcdot}\ \isavar{p}{\isacharparenright}{\kern0pt}\ \isactrlbold {\isasymcdot}\ {\isacharparenleft}{\kern0pt}\isaconst{cmp}\ \isaconst{evIdx}\ \isactrlbold {\isasymcdot}\ \isavar{v}\ \isactrlbold {\isasymcdot}\ \isavar{q}{\isacharparenright}{\kern0pt}}%
}%
\SNIP{ConRGA-Paper-Statements:thm:norm-correct-0}{%
\isakeywordONE{thm}\isamarkupfalse%
\ norm{\isacharunderscore}{\kern0pt}correct%
\isanewline%
}%
\SNIP{ConRGA-Paper-Statements:thm:norm-correct-1}{%
\isa{{\isasymlbrakk}\isavar{{\isasymGamma}}\ {\isasymturnstile}\ \isavar{v}\ \isactrlbold N{\isacharsemicolon}{\kern0pt}\ \isavar{{\isasymGamma}}\ {\isasymturnstile}\ \isavar{x}\ \isactrlbold N{\isacharsemicolon}{\kern0pt}\ \isavar{{\isasymGamma}}\ {\isasymturnstile}\ \isaconst{iWft}\ \isavar{x}\ \isactrlbold {\isasymnoteq}\ \isactrlbold {\isadigit{0}}{\isacharsemicolon}{\kern0pt}\ \isavar{{\isasymGamma}}\ {\isasymturnstile}\ \isaconst{iIsCode}\ \isavar{x}\ \isactrlbold {\isasymnoteq}\ \isactrlbold {\isadigit{0}}{\isacharsemicolon}{\kern0pt}\ \isavar{{\isasymGamma}}\ {\isasymturnstile}\ \isaconst{iQf}\ \isavar{x}\ \isactrlbold {\isasymnoteq}\ \isactrlbold {\isadigit{0}}{\isacharsemicolon}{\kern0pt}\ \isavar{{\isasymGamma}}\ {\isasymturnstile}\ \isaconst{cpFst}\ {\isacharparenleft}{\kern0pt}\isaconst{cmp}\ \isaconst{plNormIdx}\ \isactrlbold {\isasymcdot}\ \isactrlbold {\isadigit{0}}\ \isactrlbold {\isasymcdot}\ \isavar{x}{\isacharparenright}{\kern0pt}\ \isactrlbold {\isacharequal}{\kern0pt}\ \isactrlbold {\isadigit{0}}{\isasymrbrakk}\ {\isasymLongrightarrow}\ \isavar{{\isasymGamma}}\ {\isasymturnstile}\ \isaconst{iPrDecVal}\ {\isacharparenleft}{\kern0pt}\isaconst{iSubM}\ {\isacharparenleft}{\kern0pt}\isaconst{updSubM}\ \isactrlbold {\isadigit{0}}\ {\isacharparenleft}{\kern0pt}\isactrlbold n{\isadigit{0}}{\isacharparenright}{\kern0pt}\ {\isacharparenleft}{\kern0pt}\isaconst{cmp}\ \isaconst{num{\isacharunderscore}{\kern0pt}codeIdx}\ \isactrlbold {\isasymcdot}\ \isactrlbold {\isadigit{0}}\ \isactrlbold {\isasymcdot}\ \isavar{v}{\isacharparenright}{\kern0pt}{\isacharparenright}{\kern0pt}\ \isavar{x}{\isacharparenright}{\kern0pt}\ {\isacharparenleft}{\kern0pt}\isaconst{cmp}\ \isaconst{evIdx}\ \isactrlbold {\isasymcdot}\ \isavar{v}\ \isactrlbold {\isasymcdot}\ \isaconst{cpSnd}\ {\isacharparenleft}{\kern0pt}\isaconst{cmp}\ \isaconst{plNormIdx}\ \isactrlbold {\isasymcdot}\ \isactrlbold {\isadigit{0}}\ \isactrlbold {\isasymcdot}\ \isavar{x}{\isacharparenright}{\kern0pt}{\isacharparenright}{\kern0pt}}%
}%
\SNIP{ConRGA-Paper-Statements:thm:nmCor-allI-0}{%
\isakeywordONE{thm}\isamarkupfalse%
\ nmCor{\isacharunderscore}{\kern0pt}allI%
\isanewline%
}%
\SNIP{ConRGA-Paper-Statements:thm:nmCor-allI-1}{%
\isa{{\isasymemptyset}\ {\isasymturnstile}\ \isactrlbold A\ \isactrlbold {\isasymcdot}\ \isaconst{nmCorMot}}%
}%
\SNIP{ConRGA-Paper-Statements:thm:poly-reflect-sentence-0}{%
\isakeywordONE{thm}\isamarkupfalse%
\ poly{\isacharunderscore}{\kern0pt}reflect{\isacharunderscore}{\kern0pt}sentence%
\isanewline%
}%
\SNIP{ConRGA-Paper-Statements:thm:poly-reflect-sentence-1}{%
\isa{{\isasymemptyset}\ {\isasymturnstile}\ \isaconst{polyReflS}}%
}%
\SNIP{ConRGA-Paper-Statements:thm:poly-reflect-0}{%
\isakeywordONE{thm}\isamarkupfalse%
\ poly{\isacharunderscore}{\kern0pt}reflect%
\isanewline%
}%
\SNIP{ConRGA-Paper-Statements:thm:poly-reflect-1}{%
\isa{{\isasymlbrakk}\isaconst{ptmA}\ \isavar{a}{\isacharsemicolon}{\kern0pt}\ \isaconst{ptmA}\ \isavar{b}{\isacharsemicolon}{\kern0pt}\ \isaconst{trimH}\ {\isacharparenleft}{\kern0pt}\isaconst{normH}\ \isavar{a}{\isacharparenright}{\kern0pt}\ {\isacharequal}{\kern0pt}\ \isaconst{trimH}\ {\isacharparenleft}{\kern0pt}\isaconst{normH}\ \isavar{b}{\isacharparenright}{\kern0pt}{\isacharsemicolon}{\kern0pt}\ \isavar{{\isasymGamma}}\ {\isasymturnstile}\ \isavar{v}\ \isactrlbold N{\isasymrbrakk}\ {\isasymLongrightarrow}\ \isavar{{\isasymGamma}}\ {\isasymturnstile}\ {\isacharbrackleft}{\kern0pt}{\isadigit{0}}{\isasymZpinj}\isavar{v}{\isacharbrackright}{\kern0pt}\isavar{a}\ \isactrlbold {\isacharequal}{\kern0pt}\ {\isacharbrackleft}{\kern0pt}{\isadigit{0}}{\isasymZpinj}\isavar{v}{\isacharbrackright}{\kern0pt}\isavar{b}}%
}%
\SNIP{ConRGA-Paper-Statements:thm:rfl-binom-by-reflection-0}{%
\isakeywordONE{thm}\isamarkupfalse%
\ rfl{\isacharunderscore}{\kern0pt}binom{\isacharunderscore}{\kern0pt}by{\isacharunderscore}{\kern0pt}reflection%
\isanewline%
}%
\SNIP{ConRGA-Paper-Statements:thm:rfl-binom-by-reflection-1}{%
\isa{\isavar{{\isasymGamma}}\ {\isasymturnstile}\ \isavar{v}\ \isactrlbold N\ {\isasymLongrightarrow}\ \isavar{{\isasymGamma}}\ {\isasymturnstile}\ \isaconst{rflL}\ \isavar{v}\ \isactrlbold {\isacharequal}{\kern0pt}\ \isaconst{rflR}\ \isavar{v}}%
}%
\SNIP{ConRGA-Paper-Statements:thm:rfl-cube-by-reflection-0}{%
\isakeywordONE{thm}\isamarkupfalse%
\ rfl{\isacharunderscore}{\kern0pt}cube{\isacharunderscore}{\kern0pt}by{\isacharunderscore}{\kern0pt}reflection%
\isanewline%
}%
\SNIP{ConRGA-Paper-Statements:thm:rfl-cube-by-reflection-1}{%
\isa{\isavar{{\isasymGamma}}\ {\isasymturnstile}\ \isavar{v}\ \isactrlbold N\ {\isasymLongrightarrow}\ \isavar{{\isasymGamma}}\ {\isasymturnstile}\ \isaconst{plCubeL}\ \isavar{v}\ \isactrlbold {\isacharequal}{\kern0pt}\ \isaconst{plCubeR}\ \isavar{v}}%
}%
\SNIP{ConRGA-Paper-Statements:thm:rfl-prod-by-reflection-0}{%
\isakeywordONE{thm}\isamarkupfalse%
\ rfl{\isacharunderscore}{\kern0pt}prod{\isacharunderscore}{\kern0pt}by{\isacharunderscore}{\kern0pt}reflection%
\isanewline%
}%
\SNIP{ConRGA-Paper-Statements:thm:rfl-prod-by-reflection-1}{%
\isa{\isavar{{\isasymGamma}}\ {\isasymturnstile}\ \isavar{v}\ \isactrlbold N\ {\isasymLongrightarrow}\ \isavar{{\isasymGamma}}\ {\isasymturnstile}\ \isaconst{plProdL}\ \isavar{v}\ \isactrlbold {\isacharequal}{\kern0pt}\ \isaconst{plProdR}\ \isavar{v}}%
}%
\SNIP{ConRGA-Paper-Statements:thm:pr2fi-normIdx-ptmA-0}{%
\isakeywordONE{thm}\isamarkupfalse%
\ pr{\isadigit{2}}fi{\isacharunderscore}{\kern0pt}normIdx{\isacharunderscore}{\kern0pt}ptmA%
\isanewline%
}%
\SNIP{ConRGA-Paper-Statements:thm:pr2fi-normIdx-ptmA-1}{%
\isa{\isaconst{ptmA}\ \isavar{t}\ {\isasymLongrightarrow}\ \isaconst{pr{\isadigit{2}}fi}\ \isaconst{plNormIdx}\ {\isadigit{0}}\ {\isasymlbrace}\isavar{t}{\isasymrbrace}\ {\isacharequal}{\kern0pt}\ {\isasymlangle}{\isadigit{0}}{\isacharcomma}{\kern0pt}{\isasymlbrace}\isaconst{normH}\ \isavar{t}{\isasymrbrace}{\isasymrangle}}%
}%
\SNIP{ConRGA-Paper-Statements:thm:pr2fi-polyEqIdx-0}{%
\isakeywordONE{thm}\isamarkupfalse%
\ pr{\isadigit{2}}fi{\isacharunderscore}{\kern0pt}polyEqIdx%
\isanewline%
}%
\SNIP{ConRGA-Paper-Statements:thm:pr2fi-polyEqIdx-1}{%
\isa{{\isasymlbrakk}\isaconst{ptmA}\ \isavar{a}{\isacharsemicolon}{\kern0pt}\ \isaconst{ptmA}\ \isavar{b}{\isasymrbrakk}\ {\isasymLongrightarrow}\ \isaconst{pr{\isadigit{2}}fi}\ \isaconst{polyEqIdx}\ {\isasymlbrace}\isavar{a}{\isasymrbrace}\ {\isasymlbrace}\isavar{b}{\isasymrbrace}\ {\isacharequal}{\kern0pt}\ {\isasymlbrace}\isaconst{trimH}\ {\isacharparenleft}{\kern0pt}\isaconst{normH}\ \isavar{a}{\isacharparenright}{\kern0pt}\ {\isacharequal}{\kern0pt}\ \isaconst{trimH}\ {\isacharparenleft}{\kern0pt}\isaconst{normH}\ \isavar{b}{\isacharparenright}{\kern0pt}{\isasymrbrace}}%
}%
\SNIP{ConRGA-Paper-Statements:thm:pr-dec-eval-0}{%
\isakeywordONE{thm}\isamarkupfalse%
\ pr{\isacharunderscore}{\kern0pt}dec{\isacharunderscore}{\kern0pt}eval%
\isanewline%
}%
\SNIP{ConRGA-Paper-Statements:thm:pr-dec-eval-1}{%
\isa{\isaconst{prDecEvalOK}\ \isavar{f}}%
}%
\SNIP{ConRGA-Paper-Statements:thm:tru-univ-of-cert-0}{%
\isakeywordONE{thm}\isamarkupfalse%
\ tru{\isacharunderscore}{\kern0pt}univ{\isacharunderscore}{\kern0pt}of{\isacharunderscore}{\kern0pt}cert%
\isanewline%
}%
\SNIP{ConRGA-Paper-Statements:thm:tru-univ-of-cert-1}{%
\isa{{\isasymlbrakk}\isavar{{\isasymGamma}}\ {\isasymturnstile}\ \isavar{ca}\ \isactrlbold N{\isacharsemicolon}{\kern0pt}\ \isavar{{\isasymGamma}}\ {\isasymturnstile}\ \isavar{p}\ \isactrlbold N{\isacharsemicolon}{\kern0pt}\ \isavar{{\isasymGamma}}\ {\isasymturnstile}\ \isaconst{cmp}\ \isaconst{provesChkIdx}\ \isactrlbold {\isasymcdot}\ \isaconst{openJ}\ \isavar{ca}\ \isactrlbold {\isasymcdot}\ \isavar{p}\ \isactrlbold {\isasymnoteq}\ \isactrlbold {\isadigit{0}}{\isasymrbrakk}\ {\isasymLongrightarrow}\ \isavar{{\isasymGamma}}\ {\isasymturnstile}\ \isaconst{iWtrue}\ {\isacharparenleft}{\kern0pt}\isaconst{mkAppS}\ \isaconst{mkAS}\ \isavar{ca}{\isacharparenright}{\kern0pt}}%
}%
\SNIP{ConRGA-Paper-Statements:thm:tru-univ-cert-ent-0}{%
\isakeywordONE{thm}\isamarkupfalse%
\ tru{\isacharunderscore}{\kern0pt}univ{\isacharunderscore}{\kern0pt}cert{\isacharunderscore}{\kern0pt}ent%
\isanewline%
}%
\SNIP{ConRGA-Paper-Statements:thm:tru-univ-cert-ent-1}{%
\isa{\isaconst{iHasCert}\ {\isacharparenleft}{\kern0pt}\isaconst{openJ}\ \isactrlbold v\isactrlsub {\isadigit{0}}{\isacharparenright}{\kern0pt}\ {\isasymtriangleright}\ \isactrlbold v\isactrlsub {\isadigit{0}}\ \isactrlbold N\ {\isasymtriangleright}\ {\isasymemptyset}\ {\isasymturnstile}\ \isaconst{iWtrue}\ {\isacharparenleft}{\kern0pt}\isaconst{mkAppS}\ \isaconst{mkAS}\ \isactrlbold v\isactrlsub {\isadigit{0}}{\isacharparenright}{\kern0pt}}%
}%
\SNIP{ConRGA-Paper-Statements:thm:tru-univ-cert-int-0}{%
\isakeywordONE{thm}\isamarkupfalse%
\ tru{\isacharunderscore}{\kern0pt}univ{\isacharunderscore}{\kern0pt}cert{\isacharunderscore}{\kern0pt}int%
\isanewline%
}%
\SNIP{ConRGA-Paper-Statements:thm:tru-univ-cert-int-1}{%
\isa{{\isasymemptyset}\ {\isasymturnstile}\ \isactrlbold A\ \isactrlbold {\isasymcdot}\ {\isacharparenleft}{\kern0pt}\isactrlbold {\isasymlambda}\isactrlbold A\ \isactrlbold {\isasymcdot}\ {\isacharparenleft}{\kern0pt}\isactrlbold {\isasymlambda}{\isacharparenleft}{\kern0pt}\isaconst{cmp}\ \isaconst{provesChkIdx}\ \isactrlbold {\isasymcdot}\ \isaconst{openJ}\ {\isacharparenleft}{\kern0pt}\isaconst{trmVar}\ {\isadigit{1}}{\isacharparenright}{\kern0pt}\ \isactrlbold {\isasymcdot}\ \isactrlbold v\isactrlsub {\isadigit{0}}\ \isactrlbold {\isasymnoteq}\ \isactrlbold {\isadigit{0}}\ \isactrlbold {\isasymlongrightarrow}\ \isaconst{iWtrue}\ {\isacharparenleft}{\kern0pt}\isaconst{mkAppS}\ \isaconst{mkAS}\ {\isacharparenleft}{\kern0pt}\isaconst{trmVar}\ {\isadigit{1}}{\isacharparenright}{\kern0pt}{\isacharparenright}{\kern0pt}{\isacharparenright}{\kern0pt}{\isacharparenright}{\kern0pt}{\isacharparenright}{\kern0pt}}%
}%
\SNIP{ConRGA-Paper-Statements:thm:tru-univ-cert-of-dec-0}{%
\isakeywordONE{thm}\isamarkupfalse%
\ tru{\isacharunderscore}{\kern0pt}univ{\isacharunderscore}{\kern0pt}cert{\isacharunderscore}{\kern0pt}of{\isacharunderscore}{\kern0pt}dec%
\isanewline%
}%
\SNIP{ConRGA-Paper-Statements:thm:tru-univ-cert-of-dec-1}{%
\isa{{\isasymlbrakk}\isavar{{\isasymGamma}}\ {\isasymturnstile}\ \isavar{ca}\ \isactrlbold N{\isacharsemicolon}{\kern0pt}\ \isavar{{\isasymGamma}}\ {\isasymturnstile}\ \isavar{s}\ \isactrlbold N{\isacharsemicolon}{\kern0pt}\ \isavar{{\isasymGamma}}\ {\isasymturnstile}\ \isaconst{iPrDecV}\ \isavar{s}\ {\isacharparenleft}{\kern0pt}\isaconst{mkAppS}\ \isaconst{mkAS}\ \isavar{ca}{\isacharparenright}{\kern0pt}\ \isactrlbold {\isacharequal}{\kern0pt}\ \isactrlbold n{\isadigit{2}}{\isasymrbrakk}\ {\isasymLongrightarrow}\ \isavar{{\isasymGamma}}\ {\isasymturnstile}\ \isaconst{iHasCert}\ {\isacharparenleft}{\kern0pt}\isaconst{openJ}\ \isavar{ca}{\isacharparenright}{\kern0pt}}%
}%
\SNIP{ConRGA-Paper-Statements:thm:tru-univ-cert-conv-0}{%
\isakeywordONE{thm}\isamarkupfalse%
\ tru{\isacharunderscore}{\kern0pt}univ{\isacharunderscore}{\kern0pt}cert{\isacharunderscore}{\kern0pt}conv%
\isanewline%
}%
\SNIP{ConRGA-Paper-Statements:thm:tru-univ-cert-conv-1}{%
\isa{{\isasymlbrakk}\isavar{{\isasymGamma}}\ {\isasymturnstile}\ \isavar{ca}\ \isactrlbold N{\isacharsemicolon}{\kern0pt}\ \isavar{{\isasymGamma}}\ {\isasymturnstile}\ \isaconst{iWtrue}\ {\isacharparenleft}{\kern0pt}\isaconst{mkAppS}\ \isaconst{mkAS}\ \isavar{ca}{\isacharparenright}{\kern0pt}{\isasymrbrakk}\ {\isasymLongrightarrow}\ \isavar{{\isasymGamma}}\ {\isasymturnstile}\ \isaconst{iHasCert}\ {\isacharparenleft}{\kern0pt}\isaconst{openJ}\ \isavar{ca}{\isacharparenright}{\kern0pt}}%
}%
\SNIP{ConRGA-Paper-Statements:thm:tru-univ-cert-conv-int-0}{%
\isakeywordONE{thm}\isamarkupfalse%
\ tru{\isacharunderscore}{\kern0pt}univ{\isacharunderscore}{\kern0pt}cert{\isacharunderscore}{\kern0pt}conv{\isacharunderscore}{\kern0pt}int%
\isanewline%
}%
\SNIP{ConRGA-Paper-Statements:thm:tru-univ-cert-conv-int-1}{%
\isa{{\isasymemptyset}\ {\isasymturnstile}\ \isactrlbold A\ \isactrlbold {\isasymcdot}\ {\isacharparenleft}{\kern0pt}\isactrlbold {\isasymlambda}\isactrlbold A\ \isactrlbold {\isasymcdot}\ {\isacharparenleft}{\kern0pt}\isactrlbold {\isasymlambda}{\isacharparenleft}{\kern0pt}\isaconst{iPrDecTrueAt}\ {\isacharparenleft}{\kern0pt}\isaconst{mkAppS}\ \isaconst{mkAS}\ {\isacharparenleft}{\kern0pt}\isaconst{trmVar}\ {\isadigit{1}}{\isacharparenright}{\kern0pt}{\isacharparenright}{\kern0pt}\ \isactrlbold v\isactrlsub {\isadigit{0}}\ \isactrlbold {\isasymnoteq}\ \isactrlbold {\isadigit{0}}\ \isactrlbold {\isasymlongrightarrow}\ \isaconst{iHasCert}\ {\isacharparenleft}{\kern0pt}\isaconst{openJ}\ {\isacharparenleft}{\kern0pt}\isaconst{trmVar}\ {\isadigit{1}}{\isacharparenright}{\kern0pt}{\isacharparenright}{\kern0pt}{\isacharparenright}{\kern0pt}{\isacharparenright}{\kern0pt}{\isacharparenright}{\kern0pt}}%
}%
\SNIP{ConRGA-Paper-Statements:thm:all-omega-ctx-rule-fails-0}{%
\isakeywordONE{thm}\isamarkupfalse%
\ all{\isacharunderscore}{\kern0pt}omega{\isacharunderscore}{\kern0pt}ctx{\isacharunderscore}{\kern0pt}rule{\isacharunderscore}{\kern0pt}fails%
\isanewline%
}%
\SNIP{ConRGA-Paper-Statements:thm:all-omega-ctx-rule-fails-1}{%
\isa{{\isasymnot}\ {\isacharparenleft}{\kern0pt}{\isasymforall}\isabound{{\isasymGamma}}{\isachardot}{\kern0pt}\ \isabound{{\isasymGamma}}\ {\isasymturnstile}\ \isactrlbold A\ \isactrlbold {\isasymcdot}\ {\isacharparenleft}{\kern0pt}\isactrlbold {\isasymlambda}\isaconst{iWtrue}\ {\isacharparenleft}{\kern0pt}\isaconst{mkAppS}\ {\isacharparenleft}{\kern0pt}\isactrlbold n{\isasymlbrace}\isactrlbold {\isasymlambda}{\isacharparenleft}{\kern0pt}\isactrlbold {\isasymnot}\isaconst{iPrDecTrueAt}\ \isaconst{goedelWS}\ \isactrlbold v\isactrlsub {\isadigit{0}}{\isacharparenright}{\kern0pt}{\isasymrbrace}{\isacharparenright}{\kern0pt}\ {\isacharparenleft}{\kern0pt}\isaconst{iNumT}\ \isactrlbold {\isasymcdot}\ \isactrlbold v\isactrlsub {\isadigit{0}}{\isacharparenright}{\kern0pt}{\isacharparenright}{\kern0pt}{\isacharparenright}{\kern0pt}\ {\isasymlongrightarrow}\ \isabound{{\isasymGamma}}\ {\isasymturnstile}\ \isaconst{iWtrue}\ {\isacharparenleft}{\kern0pt}\isaconst{mkAppS}\ \isaconst{mkAS}\ {\isacharparenleft}{\kern0pt}\isactrlbold n{\isasymlbrace}\isactrlbold {\isasymlambda}{\isacharparenleft}{\kern0pt}\isactrlbold {\isasymnot}\isaconst{iPrDecTrueAt}\ \isaconst{goedelWS}\ \isactrlbold v\isactrlsub {\isadigit{0}}{\isacharparenright}{\kern0pt}{\isasymrbrace}{\isacharparenright}{\kern0pt}{\isacharparenright}{\kern0pt}{\isacharparenright}{\kern0pt}}%
}%
\SNIP{ConRGA-Paper-Statements:thm:tru-univ-reflect-0}{%
\isakeywordONE{thm}\isamarkupfalse%
\ tru{\isacharunderscore}{\kern0pt}univ{\isacharunderscore}{\kern0pt}reflect%
\isanewline%
}%
\SNIP{ConRGA-Paper-Statements:thm:tru-univ-reflect-1}{%
\isa{{\isasymlbrakk}\isaconst{wft}\ \isavar{u}{\isacharsemicolon}{\kern0pt}\ \isaconst{Syntax{\isacharunderscore}{\kern0pt}Term{\isachardot}{\kern0pt}closed}\ \isavar{u}{\isasymrbrakk}\ {\isasymLongrightarrow}\ \isaconst{iChk}\ \isactrlbold v\isactrlsub {\isadigit{0}}\ {\isacharparenleft}{\kern0pt}\isactrlbold n{\isasymlbrace}\isactrlbold A\ \isactrlbold {\isasymcdot}\ \isavar{u}{\isasymrbrace}{\isacharparenright}{\kern0pt}\ {\isasymtriangleright}\ \isactrlbold v\isactrlsub {\isadigit{0}}\ \isactrlbold N\ {\isasymtriangleright}\ {\isasymemptyset}\ {\isasymturnstile}\ \isaconst{iWtrue}\ {\isacharparenleft}{\kern0pt}\isactrlbold n{\isasymlbrace}\isactrlbold A\ \isactrlbold {\isasymcdot}\ \isavar{u}{\isasymrbrace}{\isacharparenright}{\kern0pt}}%
}%
\SNIP{ConRGA-Paper-Statements:thm:iLiftS-numC-0}{%
\isakeywordONE{thm}\isamarkupfalse%
\ iLiftS{\isacharunderscore}{\kern0pt}numC%
\isanewline%
}%
\SNIP{ConRGA-Paper-Statements:thm:iLiftS-numC-1}{%
\isa{\isavar{{\isasymGamma}}\ {\isasymturnstile}\ \isavar{k}\ \isactrlbold N\ {\isasymLongrightarrow}\ \isavar{{\isasymGamma}}\ {\isasymturnstile}\ \isaconst{cmp}\ \isaconst{liftPHdIdx}\ \isactrlbold {\isasymcdot}\ {\isacharparenleft}{\kern0pt}\isactrlbold n{\isadigit{1}}{\isacharparenright}{\kern0pt}\ \isactrlbold {\isasymcdot}\ \isaconst{cpPair}\ \isactrlbold {\isadigit{0}}\ {\isacharparenleft}{\kern0pt}\isaconst{cmp}\ \isaconst{num{\isacharunderscore}{\kern0pt}codeIdx}\ \isactrlbold {\isasymcdot}\ \isactrlbold {\isadigit{0}}\ \isactrlbold {\isasymcdot}\ \isavar{k}{\isacharparenright}{\kern0pt}\ \isactrlbold {\isacharequal}{\kern0pt}\ \isaconst{cmp}\ \isaconst{num{\isacharunderscore}{\kern0pt}codeIdx}\ \isactrlbold {\isasymcdot}\ \isactrlbold {\isadigit{0}}\ \isactrlbold {\isasymcdot}\ \isavar{k}}%
}%
\SNIP{ConRGA-Paper-Statements:thm:uinCert-proves-0}{%
\isakeywordONE{thm}\isamarkupfalse%
\ uinCert{\isacharunderscore}{\kern0pt}proves%
\isanewline%
}%
\SNIP{ConRGA-Paper-Statements:thm:uinCert-proves-1}{%
\isa{{\isasymlbrakk}\isaconst{wft}\ \isavar{G}{\isacharsemicolon}{\kern0pt}\ \isaconst{Syntax{\isacharunderscore}{\kern0pt}Term{\isachardot}{\kern0pt}closed}\ \isavar{G}{\isacharsemicolon}{\kern0pt}\ \isavar{Q}\ {\isasymP}{\isasymlbrakk}{\isasymemptyset}\ {\isasymturnstile}\ \isactrlbold A\ \isactrlbold {\isasymcdot}\ {\isacharparenleft}{\kern0pt}\isactrlbold {\isasymlambda}\isactrlbold A\ \isactrlbold {\isasymcdot}\ {\isacharparenleft}{\kern0pt}\isavar{G}\ \isactrlbold {\isasymcdot}\ \isactrlbold v\isactrlsub {\isadigit{0}}{\isacharparenright}{\kern0pt}{\isacharparenright}{\kern0pt}{\isasymrbrakk}{\isacharsemicolon}{\kern0pt}\ \isavar{{\isasymGamma}}\ {\isasymturnstile}\ \isavar{x}\ \isactrlbold N{\isasymrbrakk}\ {\isasymLongrightarrow}\ \isavar{{\isasymGamma}}\ {\isasymturnstile}\ \isaconst{cmp}\ \isaconst{provesChkIdx}\ \isactrlbold {\isasymcdot}\ \isaconst{openJ}\ {\isacharparenleft}{\kern0pt}\isaconst{mkAppS}\ {\isacharparenleft}{\kern0pt}\isactrlbold n{\isasymlbrace}\isavar{G}{\isasymrbrace}{\isacharparenright}{\kern0pt}\ {\isacharparenleft}{\kern0pt}\isaconst{cmp}\ \isaconst{num{\isacharunderscore}{\kern0pt}codeIdx}\ \isactrlbold {\isasymcdot}\ \isactrlbold {\isadigit{0}}\ \isactrlbold {\isasymcdot}\ \isavar{x}{\isacharparenright}{\kern0pt}{\isacharparenright}{\kern0pt}\ \isactrlbold {\isasymcdot}\ \isaconst{uinCert}\ \isavar{Q}\ \isavar{G}\ \isavar{x}\ \isactrlbold {\isasymnoteq}\ \isactrlbold {\isadigit{0}}}%
}%
\SNIP{ConRGA-Paper-Statements:thm:uinM-eq-0}{%
\isakeywordONE{thm}\isamarkupfalse%
\ uinM{\isacharunderscore}{\kern0pt}eq%
\isanewline%
}%
\SNIP{ConRGA-Paper-Statements:thm:uinM-eq-1}{%
\isa{\isavar{{\isasymGamma}}\ {\isasymturnstile}\ \isavar{x}\ \isactrlbold N\ {\isasymLongrightarrow}\ \isavar{{\isasymGamma}}\ {\isasymturnstile}\ \isaconst{uinM}\ \isavar{Q}\ \isavar{G}\ \isactrlbold {\isasymcdot}\ \isavar{x}\ \isactrlbold {\isasymcdot}\ \isavar{x}\ \isactrlbold {\isacharequal}{\kern0pt}\ \isaconst{uinCert}\ \isavar{Q}\ \isavar{G}\ \isavar{x}}%
}%
\SNIP{ConRGA-Paper-Statements:thm:univ-inst-nec-0}{%
\isakeywordONE{thm}\isamarkupfalse%
\ univ{\isacharunderscore}{\kern0pt}inst{\isacharunderscore}{\kern0pt}nec%
\isanewline%
}%
\SNIP{ConRGA-Paper-Statements:thm:univ-inst-nec-1}{%
\isa{{\isasymlbrakk}\isaconst{wft}\ \isavar{G}{\isacharsemicolon}{\kern0pt}\ \isaconst{Syntax{\isacharunderscore}{\kern0pt}Term{\isachardot}{\kern0pt}closed}\ \isavar{G}{\isacharsemicolon}{\kern0pt}\ {\isasymemptyset}\ {\isasymturnstile}\ \isactrlbold A\ \isactrlbold {\isasymcdot}\ {\isacharparenleft}{\kern0pt}\isactrlbold {\isasymlambda}\isactrlbold A\ \isactrlbold {\isasymcdot}\ {\isacharparenleft}{\kern0pt}\isavar{G}\ \isactrlbold {\isasymcdot}\ \isactrlbold v\isactrlsub {\isadigit{0}}{\isacharparenright}{\kern0pt}{\isacharparenright}{\kern0pt}{\isasymrbrakk}\ {\isasymLongrightarrow}\ {\isasymemptyset}\ {\isasymturnstile}\ \isactrlbold A\ \isactrlbold {\isasymcdot}\ {\isacharparenleft}{\kern0pt}\isactrlbold {\isasymlambda}\isaconst{iWtrue}\ {\isacharparenleft}{\kern0pt}\isaconst{mkAppS}\ \isaconst{mkAS}\ {\isacharparenleft}{\kern0pt}\isaconst{mkAppS}\ {\isacharparenleft}{\kern0pt}\isactrlbold n{\isasymlbrace}\isavar{G}{\isasymrbrace}{\isacharparenright}{\kern0pt}\ {\isacharparenleft}{\kern0pt}\isaconst{cmp}\ \isaconst{num{\isacharunderscore}{\kern0pt}codeIdx}\ \isactrlbold {\isasymcdot}\ \isactrlbold {\isadigit{0}}\ \isactrlbold {\isasymcdot}\ \isactrlbold v\isactrlsub {\isadigit{0}}{\isacharparenright}{\kern0pt}{\isacharparenright}{\kern0pt}{\isacharparenright}{\kern0pt}{\isacharparenright}{\kern0pt}}%
}%
\SNIP{ConRGA-Paper-Statements:thm:univ-inst-nec-1-0}{%
\isakeywordONE{thm}\isamarkupfalse%
\ univ{\isacharunderscore}{\kern0pt}inst{\isacharunderscore}{\kern0pt}nec{\isacharunderscore}{\kern0pt}{\isadigit{1}}%
\isanewline%
}%
\SNIP{ConRGA-Paper-Statements:thm:univ-inst-nec-1-1}{%
\isa{{\isasymemptyset}\ {\isasymturnstile}\ \isactrlbold A\ \isactrlbold {\isasymcdot}\ {\isacharparenleft}{\kern0pt}\isactrlbold {\isasymlambda}\isaconst{iWtrue}\ {\isacharparenleft}{\kern0pt}\isaconst{mkAppS}\ \isaconst{mkAS}\ {\isacharparenleft}{\kern0pt}\isaconst{mkAppS}\ {\isacharparenleft}{\kern0pt}\isactrlbold n{\isasymlbrace}\isaconst{uinG{\isadigit{1}}}{\isasymrbrace}{\isacharparenright}{\kern0pt}\ {\isacharparenleft}{\kern0pt}\isaconst{cmp}\ \isaconst{num{\isacharunderscore}{\kern0pt}codeIdx}\ \isactrlbold {\isasymcdot}\ \isactrlbold {\isadigit{0}}\ \isactrlbold {\isasymcdot}\ \isactrlbold v\isactrlsub {\isadigit{0}}{\isacharparenright}{\kern0pt}{\isacharparenright}{\kern0pt}{\isacharparenright}{\kern0pt}{\isacharparenright}{\kern0pt}}%
}%
\SNIP{ConRGA-Paper-Statements:thm:mpCert-chk-0}{%
\isakeywordONE{thm}\isamarkupfalse%
\ mpCert{\isacharunderscore}{\kern0pt}chk%
\isanewline%
}%
\SNIP{ConRGA-Paper-Statements:thm:mpCert-chk-1}{%
\isa{{\isasymlbrakk}\isavar{{\isasymGamma}}\ {\isasymturnstile}\ \isavar{cp}\ \isactrlbold N{\isacharsemicolon}{\kern0pt}\ \isavar{{\isasymGamma}}\ {\isasymturnstile}\ \isavar{cq}\ \isactrlbold N{\isacharsemicolon}{\kern0pt}\ \isavar{{\isasymGamma}}\ {\isasymturnstile}\ \isaconst{iIsCode}\ \isavar{cp}\ \isactrlbold {\isasymnoteq}\ \isactrlbold {\isadigit{0}}{\isacharsemicolon}{\kern0pt}\ \isavar{{\isasymGamma}}\ {\isasymturnstile}\ \isaconst{iIsCode}\ \isavar{cq}\ \isactrlbold {\isasymnoteq}\ \isactrlbold {\isadigit{0}}{\isacharsemicolon}{\kern0pt}\ \isavar{{\isasymGamma}}\ {\isasymturnstile}\ \isaconst{iWft}\ \isavar{cp}\ \isactrlbold {\isasymnoteq}\ \isactrlbold {\isadigit{0}}{\isacharsemicolon}{\kern0pt}\ \isavar{{\isasymGamma}}\ {\isasymturnstile}\ \isaconst{iWft}\ \isavar{cq}\ \isactrlbold {\isasymnoteq}\ \isactrlbold {\isadigit{0}}{\isacharsemicolon}{\kern0pt}\ \isavar{{\isasymGamma}}\ {\isasymturnstile}\ \isavar{A}\ \isactrlbold N{\isacharsemicolon}{\kern0pt}\ \isavar{{\isasymGamma}}\ {\isasymturnstile}\ \isavar{B}\ \isactrlbold N{\isacharsemicolon}{\kern0pt}\ \isavar{{\isasymGamma}}\ {\isasymturnstile}\ \isaconst{iChk}\ \isavar{A}\ {\isacharparenleft}{\kern0pt}\isaconst{mkOrS}\ {\isacharparenleft}{\kern0pt}\isaconst{mkNegS}\ \isavar{cp}{\isacharparenright}{\kern0pt}\ \isavar{cq}{\isacharparenright}{\kern0pt}{\isacharsemicolon}{\kern0pt}\ \isavar{{\isasymGamma}}\ {\isasymturnstile}\ \isaconst{iChk}\ \isavar{B}\ \isavar{cp}{\isasymrbrakk}\ {\isasymLongrightarrow}\ \isavar{{\isasymGamma}}\ {\isasymturnstile}\ \isaconst{iChk}\ {\isacharparenleft}{\kern0pt}\isaconst{mpCert}\ \isavar{cp}\ \isavar{cq}\ \isavar{A}\ \isavar{B}{\isacharparenright}{\kern0pt}\ \isavar{cq}}%
}%
\SNIP{ConRGA-Paper-Statements:thm:univ2-nec-0}{%
\isakeywordONE{thm}\isamarkupfalse%
\ univ{\isadigit{2}}{\isacharunderscore}{\kern0pt}nec%
\isanewline%
}%
\SNIP{ConRGA-Paper-Statements:thm:univ2-nec-1}{%
\isa{{\isasymlbrakk}\isaconst{wft}\ \isavar{G}{\isacharsemicolon}{\kern0pt}\ \isaconst{Syntax{\isacharunderscore}{\kern0pt}Term{\isachardot}{\kern0pt}closed}\ \isavar{G}{\isacharsemicolon}{\kern0pt}\ {\isasymemptyset}\ {\isasymturnstile}\ \isactrlbold A\ \isactrlbold {\isasymcdot}\ {\isacharparenleft}{\kern0pt}\isactrlbold {\isasymlambda}\isactrlbold A\ \isactrlbold {\isasymcdot}\ {\isacharparenleft}{\kern0pt}\isactrlbold {\isasymlambda}\isactrlbold A\ \isactrlbold {\isasymcdot}\ {\isacharparenleft}{\kern0pt}\isavar{G}\ \isactrlbold {\isasymcdot}\ \isaconst{trmVar}\ {\isadigit{1}}\ \isactrlbold {\isasymcdot}\ \isactrlbold v\isactrlsub {\isadigit{0}}{\isacharparenright}{\kern0pt}{\isacharparenright}{\kern0pt}{\isacharparenright}{\kern0pt}{\isasymrbrakk}\ {\isasymLongrightarrow}\ {\isasymemptyset}\ {\isasymturnstile}\ \isactrlbold A\ \isactrlbold {\isasymcdot}\ {\isacharparenleft}{\kern0pt}\isactrlbold {\isasymlambda}\isactrlbold A\ \isactrlbold {\isasymcdot}\ {\isacharparenleft}{\kern0pt}\isactrlbold {\isasymlambda}\isaconst{iWtrue}\ {\isacharparenleft}{\kern0pt}\isaconst{mkAppS}\ \isaconst{mkAS}\ {\isacharparenleft}{\kern0pt}\isaconst{mkAppS}\ {\isacharparenleft}{\kern0pt}\isaconst{mkAppS}\ {\isacharparenleft}{\kern0pt}\isactrlbold n{\isasymlbrace}\isavar{G}{\isasymrbrace}{\isacharparenright}{\kern0pt}\ {\isacharparenleft}{\kern0pt}\isaconst{cmp}\ \isaconst{num{\isacharunderscore}{\kern0pt}codeIdx}\ \isactrlbold {\isasymcdot}\ \isactrlbold {\isadigit{0}}\ \isactrlbold {\isasymcdot}\ \isaconst{trmVar}\ {\isadigit{1}}{\isacharparenright}{\kern0pt}{\isacharparenright}{\kern0pt}\ {\isacharparenleft}{\kern0pt}\isaconst{cmp}\ \isaconst{num{\isacharunderscore}{\kern0pt}codeIdx}\ \isactrlbold {\isasymcdot}\ \isactrlbold {\isadigit{0}}\ \isactrlbold {\isasymcdot}\ \isactrlbold v\isactrlsub {\isadigit{0}}{\isacharparenright}{\kern0pt}{\isacharparenright}{\kern0pt}{\isacharparenright}{\kern0pt}{\isacharparenright}{\kern0pt}{\isacharparenright}{\kern0pt}}%
}%
\SNIP{ConRGA-Paper-Statements:thm:univ2-disch-0}{%
\isakeywordONE{thm}\isamarkupfalse%
\ univ{\isadigit{2}}{\isacharunderscore}{\kern0pt}disch%
\isanewline%
}%
\SNIP{ConRGA-Paper-Statements:thm:univ2-disch-1}{%
\isa{{\isasymlbrakk}\isaconst{wft}\ \isavar{G}{\isacharsemicolon}{\kern0pt}\ \isaconst{Syntax{\isacharunderscore}{\kern0pt}Term{\isachardot}{\kern0pt}closed}\ \isavar{G}{\isacharsemicolon}{\kern0pt}\ {\isasymemptyset}\ {\isasymturnstile}\ \isactrlbold A\ \isactrlbold {\isasymcdot}\ {\isacharparenleft}{\kern0pt}\isactrlbold {\isasymlambda}\isactrlbold A\ \isactrlbold {\isasymcdot}\ {\isacharparenleft}{\kern0pt}\isactrlbold {\isasymlambda}{\isacharparenleft}{\kern0pt}\isaconst{cmp}\ \isavar{f}\ \isactrlbold {\isasymcdot}\ \isaconst{trmVar}\ {\isadigit{1}}\ \isactrlbold {\isasymcdot}\ \isactrlbold v\isactrlsub {\isadigit{0}}\ \isactrlbold {\isacharequal}{\kern0pt}\ \isactrlbold {\isadigit{0}}\ \isactrlbold {\isasymlongrightarrow}\ \isactrlbold A\ \isactrlbold {\isasymcdot}\ {\isacharparenleft}{\kern0pt}\isavar{G}\ \isactrlbold {\isasymcdot}\ \isaconst{trmVar}\ {\isadigit{1}}\ \isactrlbold {\isasymcdot}\ \isactrlbold v\isactrlsub {\isadigit{0}}{\isacharparenright}{\kern0pt}{\isacharparenright}{\kern0pt}{\isacharparenright}{\kern0pt}{\isacharparenright}{\kern0pt}{\isasymrbrakk}\ {\isasymLongrightarrow}\ {\isasymemptyset}\ {\isasymturnstile}\ \isactrlbold A\ \isactrlbold {\isasymcdot}\ {\isacharparenleft}{\kern0pt}\isactrlbold {\isasymlambda}\isactrlbold A\ \isactrlbold {\isasymcdot}\ {\isacharparenleft}{\kern0pt}\isactrlbold {\isasymlambda}{\isacharparenleft}{\kern0pt}\isaconst{cmp}\ \isavar{f}\ \isactrlbold {\isasymcdot}\ \isaconst{trmVar}\ {\isadigit{1}}\ \isactrlbold {\isasymcdot}\ \isactrlbold v\isactrlsub {\isadigit{0}}\ \isactrlbold {\isacharequal}{\kern0pt}\ \isactrlbold {\isadigit{0}}\ \isactrlbold {\isasymlongrightarrow}\ \isaconst{iWtrue}\ {\isacharparenleft}{\kern0pt}\isaconst{mkAppS}\ \isaconst{mkAS}\ {\isacharparenleft}{\kern0pt}\isaconst{mkAppS}\ {\isacharparenleft}{\kern0pt}\isaconst{mkAppS}\ {\isacharparenleft}{\kern0pt}\isactrlbold n{\isasymlbrace}\isavar{G}{\isasymrbrace}{\isacharparenright}{\kern0pt}\ {\isacharparenleft}{\kern0pt}\isaconst{cmp}\ \isaconst{num{\isacharunderscore}{\kern0pt}codeIdx}\ \isactrlbold {\isasymcdot}\ \isactrlbold {\isadigit{0}}\ \isactrlbold {\isasymcdot}\ \isaconst{trmVar}\ {\isadigit{1}}{\isacharparenright}{\kern0pt}{\isacharparenright}{\kern0pt}\ {\isacharparenleft}{\kern0pt}\isaconst{cmp}\ \isaconst{num{\isacharunderscore}{\kern0pt}codeIdx}\ \isactrlbold {\isasymcdot}\ \isactrlbold {\isadigit{0}}\ \isactrlbold {\isasymcdot}\ \isactrlbold v\isactrlsub {\isadigit{0}}{\isacharparenright}{\kern0pt}{\isacharparenright}{\kern0pt}{\isacharparenright}{\kern0pt}{\isacharparenright}{\kern0pt}{\isacharparenright}{\kern0pt}{\isacharparenright}{\kern0pt}}%
}%
\SNIP{ConRGA-Paper-Statements:thm:poly-all-sentence-0}{%
\isakeywordONE{thm}\isamarkupfalse%
\ poly{\isacharunderscore}{\kern0pt}all{\isacharunderscore}{\kern0pt}sentence%
\isanewline%
}%
\SNIP{ConRGA-Paper-Statements:thm:poly-all-sentence-1}{%
\isa{{\isasymemptyset}\ {\isasymturnstile}\ \isactrlbold A\ \isactrlbold {\isasymcdot}\ {\isacharparenleft}{\kern0pt}\isactrlbold {\isasymlambda}\isactrlbold A\ \isactrlbold {\isasymcdot}\ {\isacharparenleft}{\kern0pt}\isactrlbold {\isasymlambda}\isaconst{impCh}\ {\isacharbrackleft}{\kern0pt}\isaconst{iWft}\ {\isacharparenleft}{\kern0pt}\isaconst{trmVar}\ {\isadigit{1}}{\isacharparenright}{\kern0pt}\ \isactrlbold {\isasymnoteq}\ \isactrlbold {\isadigit{0}}{\isacharcomma}{\kern0pt}\ \isaconst{iIsCode}\ {\isacharparenleft}{\kern0pt}\isaconst{trmVar}\ {\isadigit{1}}{\isacharparenright}{\kern0pt}\ \isactrlbold {\isasymnoteq}\ \isactrlbold {\isadigit{0}}{\isacharcomma}{\kern0pt}\ \isaconst{iQf}\ {\isacharparenleft}{\kern0pt}\isaconst{trmVar}\ {\isadigit{1}}{\isacharparenright}{\kern0pt}\ \isactrlbold {\isasymnoteq}\ \isactrlbold {\isadigit{0}}{\isacharcomma}{\kern0pt}\ \isaconst{iWft}\ \isactrlbold v\isactrlsub {\isadigit{0}}\ \isactrlbold {\isasymnoteq}\ \isactrlbold {\isadigit{0}}{\isacharcomma}{\kern0pt}\ \isaconst{iIsCode}\ \isactrlbold v\isactrlsub {\isadigit{0}}\ \isactrlbold {\isasymnoteq}\ \isactrlbold {\isadigit{0}}{\isacharcomma}{\kern0pt}\ \isaconst{iQf}\ \isactrlbold v\isactrlsub {\isadigit{0}}\ \isactrlbold {\isasymnoteq}\ \isactrlbold {\isadigit{0}}{\isacharcomma}{\kern0pt}\ \isaconst{cmp}\ \isaconst{polyEqIdx}\ \isactrlbold {\isasymcdot}\ \isaconst{trmVar}\ {\isadigit{1}}\ \isactrlbold {\isasymcdot}\ \isactrlbold v\isactrlsub {\isadigit{0}}\ \isactrlbold {\isasymnoteq}\ \isactrlbold {\isadigit{0}}{\isacharbrackright}{\kern0pt}\ {\isacharparenleft}{\kern0pt}\isaconst{iWtrue}\ {\isacharparenleft}{\kern0pt}\isaconst{mkAppS}\ \isaconst{mkAS}\ {\isacharparenleft}{\kern0pt}\isaconst{mkAppS}\ {\isacharparenleft}{\kern0pt}\isaconst{mkAppS}\ {\isacharparenleft}{\kern0pt}\isactrlbold n{\isasymlbrace}\isaconst{polyQ}{\isasymrbrace}{\isacharparenright}{\kern0pt}\ {\isacharparenleft}{\kern0pt}\isaconst{cmp}\ \isaconst{num{\isacharunderscore}{\kern0pt}codeIdx}\ \isactrlbold {\isasymcdot}\ \isactrlbold {\isadigit{0}}\ \isactrlbold {\isasymcdot}\ \isaconst{trmVar}\ {\isadigit{1}}{\isacharparenright}{\kern0pt}{\isacharparenright}{\kern0pt}\ {\isacharparenleft}{\kern0pt}\isaconst{cmp}\ \isaconst{num{\isacharunderscore}{\kern0pt}codeIdx}\ \isactrlbold {\isasymcdot}\ \isactrlbold {\isadigit{0}}\ \isactrlbold {\isasymcdot}\ \isactrlbold v\isactrlsub {\isadigit{0}}{\isacharparenright}{\kern0pt}{\isacharparenright}{\kern0pt}{\isacharparenright}{\kern0pt}{\isacharparenright}{\kern0pt}{\isacharparenright}{\kern0pt}{\isacharparenright}{\kern0pt}}%
}%
\SNIP{ConRGA-Paper-Statements:thm:tru-univ-by-reflection-0}{%
\isakeywordONE{thm}\isamarkupfalse%
\ tru{\isacharunderscore}{\kern0pt}univ{\isacharunderscore}{\kern0pt}by{\isacharunderscore}{\kern0pt}reflection%
\isanewline%
}%
\SNIP{ConRGA-Paper-Statements:thm:tru-univ-by-reflection-1}{%
\isa{{\isasymlbrakk}\isaconst{ptmA}\ \isavar{a}{\isacharsemicolon}{\kern0pt}\ \isaconst{ptmA}\ \isavar{b}{\isacharsemicolon}{\kern0pt}\ \isaconst{trimH}\ {\isacharparenleft}{\kern0pt}\isaconst{normH}\ \isavar{a}{\isacharparenright}{\kern0pt}\ {\isacharequal}{\kern0pt}\ \isaconst{trimH}\ {\isacharparenleft}{\kern0pt}\isaconst{normH}\ \isavar{b}{\isacharparenright}{\kern0pt}{\isasymrbrakk}\ {\isasymLongrightarrow}\ {\isasymemptyset}\ {\isasymturnstile}\ \isaconst{iWtrue}\ {\isacharparenleft}{\kern0pt}\isactrlbold n{\isasymlbrace}\isactrlbold A\ \isactrlbold {\isasymcdot}\ {\isacharparenleft}{\kern0pt}\isaconst{polyQ}\ \isactrlbold {\isasymcdot}\ {\isacharparenleft}{\kern0pt}\isactrlbold n{\isasymlbrace}\isavar{a}{\isasymrbrace}{\isacharparenright}{\kern0pt}\ \isactrlbold {\isasymcdot}\ {\isacharparenleft}{\kern0pt}\isactrlbold n{\isasymlbrace}\isavar{b}{\isasymrbrace}{\isacharparenright}{\kern0pt}{\isacharparenright}{\kern0pt}{\isasymrbrace}{\isacharparenright}{\kern0pt}}%
}%
}{}

\newcommand{\GetSnip}[1]{%
    \ifcsname snippet--#1\endcsname%
        \csname snippet--#1\endcsname%
    \else%
        \PackageWarning{snips}{Snippet ``#1'' is undefined.}%
        \emph{Warning: Snippet ``#1'' is undefined.}%
    \fi%
}

\newcommand{\Snippet}[1]{{%
  \newcount\i
  \i=0
  \loop
    \GetSnip{#1-\the\i}%
    \advance \i 1
  \ifcsname snippet--#1-\the\i\endcsname
  \repeat
}}

\newcommand{\SnippetPart}[3]{{%
  \newcount\i
  \i=#1
  \loop
    \ifnum \i=#2
      \renewcommand{\isanewline}{}%
    \fi
    \GetSnip{#3-\the\i}%
    \advance \i 1
    \ifnum \i>#2 {}
    \else \repeat
}}

\hypersetup{breaklinks=true}

\newcommand{\com}[1]{}

\long\def\note#1{}	%

\newcommand{\pga}{PGA\xspace}		%
\newcommand{\rga}{RGA\xspace}		%

\newcommand{\gdl}{G\"odel\xspace}

\newcommand{\limp}{\rightarrow}		%
\newcommand{\liff}{\leftrightarrow}	%
\newcommand{\ldef}{\equiv}		%

\newcommand{\typestyle}[1]{\textsf{#1}\xspace}
\com{	%
\newcommand{\tbool}{\typestyle{bool}}
\newcommand{\ttrue}{\typestyle{true}}
\newcommand{\tfalse}{\typestyle{false}}

\newcommand{\tnat}{\typestyle{nat}}

}%
\newcommand{\tbool}{\typestyle{B}}	%
\newcommand{\ttrue}{\typestyle{T}}
\newcommand{\tfalse}{\typestyle{F}}

\newcommand{\tnat}{\typestyle{N}}

\newcommand{\ctrue}{\ttrue}
\newcommand{\cfalse}{\tfalse}

\newcommand{\judgment}[1]{\textsf{ #1}}

\newcommand{\jbool}{\judgment{\tbool}}

\newcommand{\jnat}{\judgment{\tnat}}

\makeatletter
\newlength{\doublefracgap}
\DeclareRobustCommand{\doublefrac}[2]{%
  \mathinner{\mathpalette\doublefrac@{{#1}{#2}}}%
}
\newcommand{\doublefrac@}[2]{\doublefrac@@#1#2}
\newcommand{\doublefrac@@}[3]{%
  \ooalign{%
    \raisebox{\doublefracgap}{$\m@th#1\frac{#2}{\phantom{#3}}$}\cr
    \raisebox{-\doublefracgap}{$\m@th#1\frac{\phantom{#2}}{#3}$}\cr
  }%
}
\newcommand{\ddoublefrac}[2]{{\displaystyle\doublefrac{#1}{#2}}}

\makeatother

\newcommand{\irl}[1]{\ensuremath{\mathit{#1}}}		%

\newcommand{\infrule}[3][]{\cfrac{#2}{#3}\IfStrEq{#1}{}{}{\ \irl{#1}}}
\newcommand{\infeqv}[3][]{\ddoublefrac{#2}{#3}\IfStrEq{#1}{}{}{\ \irl{#1}}}
\newcommand{\infceqv}[4][]{\cfrac{#2\qquad}{}\ddoublefrac{#3}{#4}
				\IfStrEq{#1}{}{}{\ \irl{#1}}}

\newcommand{\z}{\mathbf 0}

\newcommand{\suc}{\mathbf S}

\newcommand{\cR}{\mathcal R}

\newcommand{\quo}[1]{\ulcorner{#1}\urcorner}

\newcommand{\tofn}[1]{\underline{#1}}

\SetSymbolFont{stmry}{bold}{U}{stmry}{m}{n}

\theoremstyle{definition}
\newtheorem{defi}[thm]{Definition}
\newtheorem{rem}[thm]{Remark}	%
\theoremstyle{plain}
\newcommand{\cmp}{\mathit{cmp}}		%
\newcommand{\den}[1]{\llbracket #1 \rrbracket}	%
\newcommand{\cA}{\mathcal A}		%
\newcommand{\lft}[1]{{\uparrow}#1}	%
\newcommand{\vz}{v_0}			%
\newcommand{\wbev}[1]{\Downarrow^{\mathrm{w}}_{#1}}	%
\newcommand{\whev}[1]{\Downarrow^{\mathrm{w,hd}}_{#1}}	%
\newcommand{\wf}{\mathrm{wf}}		%
\newcommand{\proves}{\mathrel{:}}	%
\newcommand{\cpair}[2]{\langle #1, #2 \rangle}	%
\newcommand{\Some}[1]{\mathsf{Some}\,#1}
\newcommand{\None}{\mathsf{None}}
\newcommand{\iRec}{\mathsf{Rec}}	%
\newcommand{\iSub}{\mathsf{Sub}}	%
\newcommand{\sof}[2]{\hat{#1}({#2})}	%
\newcommand{\ssub}[3]{[{#1}{\mapsto}{#2}]{#3}}	%
\newcommand{\bquo}[1]{\underline{\vrule width0pt depth0.08em\relax\smash[b]{\quo{#1}}}}

\newcommand{\decV}{\mathsf{dec}}		%
\newcommand{\decP}{\mathsf{dec}^{\#}}	%
\newcommand{\dtab}{\mathsf{tab}}		%
\newcommand{\taut}{\tau_{\mathsf{tab}}}	%
\newcommand{\Dec}{\mathit{Dec}}		%
\newcommand{\Aff}{\mathit{Aff}}		%
\newcommand{\Tru}{\mathit{True}}	%
\newcommand{\Prv}{\mathit{Prv}}		%
\newcommand{\Chk}{\mathit{Chk}}		%
\newcommand{\ChkE}{\mathit{Chk}^{\varepsilon}}	%
\newcommand{\Cert}{\mathit{Cert}}	%
\newcommand{\Subst}{\mathit{Subst}}	%
\newcommand{\Num}{\mathit{Num}}		%
\newcommand{\Wfc}{\mathit{IsWf}}		%
\newcommand{\IsCode}{\mathit{IsCode}}	%
\newcommand{\Diag}{\mathit{Diag}}	%
\newcommand{\ConR}{\mathit{Con}}	%
\newcommand{\ChkQf}{\Chk_{\mathrm{qf}}}	%
\newcommand{\QfScan}{\mathit{QfScan}}	%
\newcommand{\qfs}{\mathit{qf}}		%
\newcommand{\Rcd}{\mathit{Val}}		%
\newcommand{\Frag}{\mathit{Frag}}	%
\newcommand{\LvC}{\mathit{Lv}}		%
\newcommand{\LF}{\mathit{LF}}		%
\newcommand{\Dep}{\mathit{dp}}		%
\newcommand{\Ok}{\mathit{Ok}}		%
\newcommand{\Res}{\mathit{Res}}		%
\newcommand{\DynFam}{\mathit{DynFam}}	%
\newcommand{\Stat}{\mathit{Stat}}	%
\newcommand{\vdashO}{\vdash_{\mathrm O}}	%
\newcommand{\oga}{OGA\xspace}		%
\newcommand{\eqv}{\mathrel{\dashv\vdash}}	%
\newcommand{\PrC}{\mathit{Pr}}		%

\newcommand{\artifacttag}{\texttt{rga-}\allowbreak\texttt{2026-11-11}\xspace}

\begin{document}

\title{Internal Reasoning and Consistency \\ in Reflective Grounded Arithmetic}
\author{Bryan Ford}
\date{}
\maketitle

\begin{abstract}
Reflective Grounded Arithmetic (RGA) is a paracomplete arithmetic
allowing direct expression of any general-recursive computation.
Every RGA term is a program and truth is grounded evaluation, with a
proof certificate at each universal quantifier.
This paper
shows, in a machine-checked Isabelle/HOL development, that such a system
is powerful enough to carry out its own metamathematics.
A fuel-indexed decider of RGA's grounded semantics
and checkers for RGA's proofs
are compiled primitive-recursive RGA terms.
RGA proves internally that its decider never affirms $0 = 1$,
and that every well-formed term it affirms
has an RGA proof accepted by a transparent test.
RGA proves that, for each valid RGA proof, its conclusion is true.
RGA proves that each grounded-true well-formed term is provable.
RGA proves that no proof in either of two decided classes
of RGA proofs derives $0 = 1$.

RGA's truth predicate also serves as a medium for proof by reflection,
of the kind that adds nothing to the system.
RGA derives disquotation for quantifier-free formulas
as rules in every context, at a free variable,
and proves, as one sentence, that whenever its normalizer for polynomials
in one variable gives two terms the same normal form,
the equation between them is true at every value of the variable.
Each such polynomial identity then follows from that sentence
by closed evaluation and disquotation, with no reflection rule.
A companion paper proves \gdl's incompleteness theorems
for RGA and its extensions.
\end{abstract}

\section{Introduction}
\label{sec:intro}

G\"odel's second incompleteness theorem is usually read as the point
where a formal system's knowledge of itself runs out: a consistent
system of sufficient strength cannot prove its own consistency.  The theorem's
textbook proof is valid even in intuitionistic arithmetic, but the proof
relies on two key elements.
The proof asserts a self-referential sentence to be
equivalent to that sentence's own unprovability, and the proof derives implications by
discharging hypotheses, provability and consistency statements, that
the system can neither prove nor refute.
This paper asks how much a formal system can prove about itself when
neither element is available, and when everything the system can say
is a computation.

Reflective Grounded Arithmetic, or \rga~\cite{rga-paper}, is such a
system.  \rga's terms are programs: numerals, $\lambda$-terms, a
recursor, and a quantifier combinator $\cA$.  A term is true when its
evaluation yields $1$, with a proof certificate at each universal
quantifier, false when it yields $0$, and \emph{ungrounded} when it
yields neither value.  So \rga is paracomplete: excluded middle fails at
ungrounded terms.  \rga takes its logic, strong Kleene, and the idea of
grounding from Kripke's theory of truth~\cite{kripke75outline}, but
Kripke's sentences are classically grounded, in a bivalent base
arithmetic, while \rga's are computationally grounded, arithmetic
included (\S\ref{sec:bg-rga}).  \rga discharges a hypothesis only when \rga can prove that hypothesis
decided, whether the discharge is by implication introduction or by
case analysis.  The companion paper~\cite{rga-paper} proves in
Isabelle/HOL that \rga is sound and consistent, and that \rga is
complete for grounded truth: provability coincides with grounded truth
for well-formed terms.  This completeness is not \emph{negation
completeness} (also called syntactic completeness: every sentence is
provable or refutable), nor is it the completeness of G\"odel's completeness
theorem, which is relative to every model of a theory.
\rga's completeness is completeness for one intended semantics.
Classically, completeness for one intended semantics is impossible for arithmetic,
since true arithmetic is not recursively enumerable.
Grounded truth is recursively enumerable, however,
which is what makes possible RGA's completeness for grounded truth.

\gdl's incompleteness theorems hold for \rga,
as a separate paper, the G\"odel paper~\cite{goedel-paper}, proves.
The sentences that classical reasoning calls true but
unprovable are ungrounded in \rga, because \rga's universal quantifier is settled
only by a proof certificate, not instance by instance.
Read instance by instance,
in the $\omega$-semantics (\S\ref{sec:inside-omega}),
these sentences are true.
These ungrounded sentences do not make \rga weak.
For example, in the formalization that the two papers share,
\rga proves, of each closed ground-type term of \gdl's System~T,
that this term denotes a natural number~\cite[Thm.~5.7]{rga-paper}.
The proof gives each System~T term its own \rga derivation,
which a HOL proof of normalization for System~T constructs.

Everything \rga can express is computation.
Syntactically, the code of an \rga term is a natural number,
and a natural number is a term.  \rga's semantics coincides
with a fuel-indexed computation (Theorem~\ref{thm:in-dec}).  \rga is
therefore a closed universe of computational reasoning.
This paper's main exhibit is that
this computational universe can reason about itself.
Coding, compiling primitive-recursive functions, deciding
\rga's own semantics, and checking \rga's own proofs all happen inside
\rga, as \rga terms that \rga itself reasons about.
This paper's main insight is thus, in essence,
that pure computation can reason about itself.

Paracomplete logics have a reputation for being too weak for this task.
Developing a closely related paracomplete theory of truth, Maudlin
concluded that ``banning the use of T-Inferences in \emph{reductio}
arguments eliminates all of the profound results of
metalogic''~\cite[p.~134]{maudlin06truth}, and earlier work on
grounded deduction left open whether grounded reasoning loses G\"odel's
diagonal arguments or gains from them~\cite{ford24reasoning}.
\rga's restriction is of the same kind: \rga discharges a hypothesis
only when the hypothesis is provably decided.  For \rga at least, the
profound results survive, in a precisely stated form.
Tarski's theorem takes the form that
the material T-biconditional at a Liar is explosive:
at the Liar of \rga's truth predicate, and at the diagonal
Liar of every closed predicate term.
So no consistent extension of
\rga by hypotheses adds such a biconditional as an axiom
(\S\ref{sec:pillars-tarski}).
G\"odel's and L\"ob's theorems survive as well, in the forms that the
G\"odel paper states (``G\"odel's theorems'' below).

\paragraph{Takeaways.}
Three results organize the paper.  All three results
are mechanized in Isabelle/HOL.  A term is \emph{decided} when \rga
proves the term true or false, $\vdash t \lor \neg t$.

\begin{enumerate}
\item[(a)] \emph{\rga is a working system for mathematics, including
  mathematics about itself} (Sections~\ref{sec:inside},
  \ref{sec:pillars} and~\ref{sec:con}).  A compiled decider for \rga's
  own fuel-indexed semantics, a truth predicate $\Tru$ built from the decider,
  and a faithful proof checker $\Chk$ are \rga terms, and \rga proves
  substantial facts about these machines internally, as single sentences: that
  the decider never affirms $\z = \tofn 1$ at any fuel
  (Theorem~\ref{thm:in-never}); that whatever well-formed term the
  decider affirms has a certificate, a number that a transparent proof
  test $\pi$ accepts as a proof of that term (Theorem~\ref{thm:pl-int}); and the
  consistency results of~(b) below.  Around these machines, \rga proves, for each
  \rga proof, that the proof's conclusion is true, and for each
  grounded-true well-formed term, that the term is provable. On the
  classical reading these two schemata are the analogues of
  the Hilbert--Bernays derivability condition D1 for $\Tru$ and for provability
  (Section~\ref{sec:pillars}).
\item[(b)] \emph{Unconditional internal consistency for two decided
  classes of proofs} (Section~\ref{sec:con}).  \rga proves, in the
  empty context and with no hypothesis, that no accepted proof without
  universal quantifiers concludes $\z = \tofn 1$, and the same for a
  larger decided class of proofs that admits quantifiers away from a
  few kinds of positions; \rga proves that the first class is contained
  in the second.  A third
  decided class admits quantified hypotheses whose quantifiers nest to
  a certified depth.  \rga derives the consistency sentence of this third
  class under one closed hypothesis that \rga neither proves nor refutes, and
  \oga~\cite{oga-paper}, \rga with one more rule, proves that
  hypothesis.  These results are \rga's counterparts of the classical
  fact that Peano arithmetic proves the consistency of each of its
  fragments $\mathrm{I}\Sigma_n$~\cite{hajek93metamathematics}.  These
  results do not escape the second theorem, and the proof-theoretic
  strength of the three classes is open.
\item[(c)] \emph{The truth predicate is a medium for proof by reflection}
  (\S\ref{sec:pillars-reflect}).  In proof by reflection, a decision
  procedure proved correct inside a system turns each goal of a class
  into a computation.  Of the three kinds of reflection that
  \S\ref{sec:pillars-reflect} distinguishes, logical reflection
  principles add strength, and computational reflection through
  provability yields no new theorems when the system is 1-consistent,
  and by G\"odel's second theorem the system does not prove its own
  1-consistency.  Computational
  reflection through a denotation adds nothing to the system.  \rga's
  construction is of this third kind, with the decider's values and
  truth as the denotation.  \rga derives T-In and T-Out for
  quantifier-free formulas as rules in every context, at the images of
  quoted templates and at an arbitrary term, a free variable included
  (Theorem~\ref{thm:pl-rf-tin}).  \rga also contains a normalizer for
  polynomials in one variable, built from explicit primitive-recursive
  machines and proved correct by one internal induction
  (Theorem~\ref{thm:pl-rf-norm}).  From no hypotheses, \rga proves one
  sentence of reflection: for all codes $x$ and $y$ that pass the code
  guards and the normalizer's equality test, and for every $v$, the
  image at $v$ of the equation between $x$ and $y$ is true
  (Theorem~\ref{thm:pl-rf-sentence}).  Each polynomial identity in one
  variable then follows at a free variable from that sentence, with
  the test's value and the guards by closed evaluation and the step
  from truth to the identity by open T-Out, and with no reflection rule
  (Theorem~\ref{thm:pl-rf-reflect}).  The construction adds no
  theorems, since the construction adds no rule.  Nor does the
  construction give small certificates.  \rga's conversion rules do not
  evaluate, so the closed evaluations reach \rga through open
  completeness.  T-In at a universal is also a proof builder in \rga:
  an explicit primitive-recursive builder, whose
  acceptance \rga verifies, yields a sentence one quantifier level above
  the sentence of reflection, with the guards outside the inner
  universal over~$v$ (Theorem~\ref{thm:pl-tin-disch}).
\end{enumerate}

\paragraph{G\"odel's theorems.}
The G\"odel paper~\cite{goedel-paper} proves G\"odel's theorems for
\rga with the machines of Section~\ref{sec:inside}, and
\S\ref{sec:inside-codes} lists the five results of the G\"odel paper
that this paper uses.  In \rga, G\"odel's proof survives as
entailments.  By the disjunction property every implication that the
textbook proof uses fails, while every entailment and every rule that
the proof uses holds [\citen{goedel-paper}, Section ``G\"odel's proof in four forms''].  \rga does not derive
its consistency sentence $\ConR$ [\citen{goedel-paper}, theorem~\textsf{rga\_goedel2}].  \rga proves
its own second theorem as the entailment $\ConR \vdash
\neg\PrC(\bquo{\ConR})$, and in \emph{decided-antecedent form}, a
sentence $\forall s.\,(D(s) \limp B)$ with $D$ decided that stands for
the conditional from $\exists s.\,D(s)$ to $B$
(\S\ref{sec:inside-codes}).  This sentence is $\forall q.\,(\Chk(q,
\bquo{\ConR}) \limp \neg\ConR)$, with an explicit transformer from
accepted proofs of $\ConR$ to accepted proofs of $\z = \tofn 1$
([\citen{goedel-paper}, theorem~\textsf{rga\_goedel2\_internal\_contra}], [\citen{goedel-paper}, theorem~\textsf{rga\_goedel2\_sentence}]).  The
G\"odel paper also proves the second theorem once for every extension
of \rga whose checker meets the G\"odel conditions
{}[\citen{goedel-paper}, Section ``The second theorem once: the G\"odel conditions''], with instances at \rga and at the
inductive grounded arithmetic IGA~\cite{iga-paper}.  Adoption, a
reflection rule keyed to a system's own checker, gives consistent
extensions of \rga and of IGA that derive $\ConR$
{}[\citen{goedel-paper}, Section ``Adoption: the reflection rule, and the consistency family''].  Whether IGA alone derives $\ConR$ is open
{}[\citen{goedel-paper}, Section ``Does IGA alone prove \rga's consistency?''].

\paragraph{Classical versus grounded readings.}
This paper's results admit two alternative readings.
On the classical reading, \rga's
grounded semantics is $\Sigma_1$-shaped: a sentence is grounded-true
when some finite evaluation with finite certificates establishes it,
and for well-formed sentences grounded truth coincides with
provability.  $\Tru$ is then a provability predicate, and the paper's
results are theorems about a system that reasons correctly about its
own derivations.  On \rga's own reading, truth is \emph{grounded
truth}, truth as establishment: a sentence is true when a finite
evaluation establishes it, with a proof certificate at each universal
quantifier.  This is a \emph{constructive} notion of truth in a narrow
sense, truth only through a construction, in the tradition that
identifies truth with proof (Brouwer and
Heyting~\cite{brouwer75philosophy,heyting71intuitionism}; Dummett,
Prawitz, Martin-L\"of~\cite{dummett78truth,prawitz77meaning,martin-lof80intuitionistic}),
which philosophers often call anti-realist or verificationist.  Grounded
truth does not bring intuitionistic logic with it: \rga's logic is strong Kleene, with
$\neg\neg A \eqv A$.  On this reading, the fact that $\Tru$ coincides with
provability is not a disguise but the thesis, here a theorem.
Sentences that no evaluation settles are ungrounded: the Liar for truth,
but also the G\"odel sentence and the consistency sentence $\ConR$,
which do not refer to their own truth at all.  These ungrounded
sentences are not paradoxical, only unestablished.  \rga even proves
every numeral instance of the universal quantifier inside the G\"odel
sentence and inside $\ConR$.
Section~\ref{sec:pillars} develops both readings, including the sense
in which grounded truth is establishment-transparent, and the classical
sense in which grounded truth is not transparent (\S\ref{sec:pillars-tarski}).  A companion
essay~\cite{ford-diag} discusses the grounded reading at length; this
paper does not depend on that essay.

\paragraph{What is new.}
The arithmetization is classical in kind.  Coded syntax, a verified
compiler, an evaluator and a proof checker are the standard machinery
behind any treatment of G\"odel's theorems (Section~\ref{sec:inside}).
What is new is where this machinery lives and what the system proves
about the machinery.  In \rga the machines are explicit
primitive-recursive trees inside a concrete arithmetic whose
arithmetic sentences themselves can lack a truth value, and a check
that runs with every build confirms that the machines contain no
opaque index (\S\ref{sec:inside-audit}).  \rga proves the machines'
laws at abstract arguments.  The substantial pieces are the decider's
internal soundness (Theorem~\ref{thm:in-never}), internal completeness
as one sentence (Theorem~\ref{thm:pl-int}), the consistency theorems
for decided classes with the location of the obstruction to full
internal consistency (Section~\ref{sec:con}), and the sentence of
reflection through the truth predicate
(Theorem~\ref{thm:pl-rf-sentence}).  The reflection construction
follows Boutin's pattern of reflection through a
denotation~\cite{boutin97using}.  What differs is the medium.  The
denotation is \rga's own truth predicate, so the construction needs no
separate syntax type, no separately defined denotation and no separate
lemma that relates the two (\S\ref{sec:pillars-reflect}).  The
metatheory is classical HOL; what is new lies in the objects inside
\rga and in what \rga derives about those objects.

\paragraph{Kinds of claim.}
Every theorem in this paper is a theorem of HOL about \rga, and the
paper says of each theorem which kind of claim the theorem makes: a schema of \rga derivations, one
internal sentence, an entailment, or a HOL fact about \rga.  For every
conditional claim the paper also says which of the four forms of a
conditional is meant.  \S\ref{sec:inside-codes} sets out both
distinctions and maps the kinds of claim onto the forms of a conditional
(Table~\ref{tab:in-forms}).  The soundness theorem of
takeaway~(a), for example, is a schema, and
Theorem~\ref{thm:pl-no-uniform} shows that the single-sentence form of
that schema is not an \rga theorem.  Every theorem carries its Isabelle
name in its header, and Appendix~\ref{app:endpoints} maps each theorem
to its theory file
in the development.

\paragraph{Outline.}
Section~\ref{sec:bg} recalls \rga.  Section~\ref{sec:inside} builds the coding machinery,
the decider, the truth predicate and the checker inside \rga, sets
out the four forms of a conditional, and lists the results that this
paper takes from the G\"odel paper.  Section~\ref{sec:pillars}
proves the soundness and completeness pillars, studies the truth
predicate and its Liars, sets the truth predicate against the
axioms of truth, and uses the truth predicate as a medium for proof by
reflection (\S\ref{sec:pillars-reflect}).
Section~\ref{sec:con} presents what \rga does prove about its own
consistency: unconditional internal consistency for two decided
classes of proofs, a third class under one closed hypothesis, and
where the obstruction to full consistency lies.  Section~\ref{sec:mech}
describes the mechanization, its trusted base and its build-time
checks, Section~\ref{sec:related} discusses related work, and
Section~\ref{sec:concl} concludes.  The appendices give the next level
of detail: a glossary of the paper's notation
(Appendix~\ref{app:notation}); the coding and the machines with their
build-time checks (Appendix~\ref{app:machines}); the proofs of the
pillars, of the truth predicate against the axioms of truth and of the
reflection results (Appendix~\ref{app:pillars}); the consistency folds
(Appendix~\ref{app:fold}); the level theorems and \oga
(Appendix~\ref{app:oga}); a catalog of the development's traps,
the tempting steps that fail (Appendix~\ref{app:traps}); and a table
mapping every theorem to its Isabelle name
(Appendix~\ref{app:endpoints}).

\paragraph{AI disclosure.}
The machine-checked proofs and the first drafts of this paper's text
were produced with substantial assistance from Claude models from
Anthropic, working under the direction and review of the author
(Section~\ref{sec:mech}).  The definitions that a reader must trust
are listed in Section~\ref{sec:mech}, and the endpoint theorems with
their Isabelle names in Appendix~\ref{app:endpoints}; the author
reviewed these definitions and theorems.  Appendix~\ref{app:endpoints-statements}
prints the headline statements and the core definitions as Isabelle
generates them from the build.

\section{Background}
\label{sec:bg}

This section recalls \rga from the companion paper~\cite{rga-paper}.
G\"odel's theorems for \rga, with the classical proof of the second
theorem as background, are the subject of a separate paper, ``G\"odel's
Incompleteness Theorems in Grounded Arithmetic''~\cite{goedel-paper}, the G\"odel
paper for short.  The section assumes
familiarity with operational semantics and the untyped
$\lambda$-calculus, but no particular background in proof theory.
Appendix~\ref{app:notation} collects the paper's notation in one
table, for reference while reading.

\subsection{Reflective Grounded Arithmetic}
\label{sec:bg-rga}

\paragraph{Grounded truth.}
\rga belongs to the grounded-deduction family of
logics~\cite{ford24reasoning}.  A sentence is a computation, and it is
true when its evaluation settles on true, false when its evaluation
settles on false, and \emph{ungrounded} when its evaluation never
settles.  Excluded middle plays the role of the assumption that every
program halts, and \rga declines it for the same reason a programmer
would.  Ungroundedness is not a third truth value, just as divergence
is not a third return value; it is the absence of a verdict.  Logics
with such gaps are \emph{paracomplete}.  \rga takes its logic, strong
Kleene, and the idea of grounding from Kripke's theory of
truth~\cite{kripke75outline}, but it grounds differently.  Kripke's
sentences are \emph{classically grounded}: his base arithmetic is
bivalent, so every arithmetic sentence, the consistency sentence
included, has a classical value, and gaps arise only among sentences
that contain the truth predicate, such as the Liar and the
truth-teller.  \rga's sentences are
\emph{computationally grounded}, all the way down: arithmetic is
computation too, and a universal quantifier is settled only by a proof
certificate.  So \rga's gaps sit at the arithmetic level.  This
grounded truth is truth as establishment, in the sense of
Section~\ref{sec:intro}.  A \emph{sentence} is a closed term.  A
sentence is ungrounded (or gapped) when no evaluation settles it,
whether or not it refers to itself.  On well-formed sentences whose values are truth
values (the negations, disjunctions, equations and universals,
\S\ref{sec:inside-dec}), the gaps are exactly the sentences that \rga
neither proves nor refutes (Theorem~\ref{thm:pl-prov-iff}, applied to
a sentence and to its negation); terms with a value that is not a
truth value, such as $\suc\suc\z$, are neither proved nor refuted
either.  The gaps include the Liar, but also the consistency sentence,
which is not self-referential.  HOL proves \rga consistent, so the
consistency sentence is classically true, and true in the classical
$\omega$-semantics of \S\ref{sec:inside-omega}; in \rga's own semantics
the consistency sentence is ungrounded, never grounded-true, because no certificate establishes
it.

\paragraph{The language.}
\rga's terms form an untyped $\lambda$-calculus with ten operators:
the divergent constant $\bot$, a term with no value, zero $\z$, successor $\suc$, equality of naturals
$=$, negation $\neg$, disjunction $\lor$, $\lambda$-abstraction,
application $\cdot$, a recursor $\mathcal R$, and the universal
quantifier $\cA$.  There is no distinction between terms and formulas:
a formula is a term, evaluated like any other.  Variables are de
Bruijn indices, and $\lambda$ is the only binder, so the quantifier is
a combinator: $\cA\cdot u$ says that $u$ holds of every natural
number, and $\forall x.\,b$ abbreviates $\cA\cdot(\lambda b)$.  The
paper writes binders with names, $\forall x$, $\lambda x$, for
readability; formally the bound variable is a de Bruijn index.
Application is call by name.  Conjunction, implication and the
biconditional are the classical abbreviations from $\neg$ and $\lor$;
in particular $a \limp b$ is $\neg a \lor b$.  Self-reference needs no
special machinery.  The Liar is a term $L$ whose evaluation unfolds to
$\neg L$, and \rga's semantics never grounds $L$.

\paragraph{The semantics.}
The semantics is a step-indexed evaluation relation: at fuel $s$ a
term $t$ evaluates to a value $r$, a natural number.  (The fuel is
what the companion paper calls the step index.)  The truth
values are the numerals $\ctrue = 1$ and $\cfalse = 0$.  Nine operators
have ordinary structural clauses.  The universal quantifier cannot be
evaluated by visiting its infinitely many instances, and \rga grounds
it instead in a finite object: a proof.  In the companion paper's
semantics, $\cA\cdot u$ is true at fuel $s+1$ when a proof in \rga's
own proof system, with code at most $s$, derives the open instance of
$u$, and false when a proof, with code at most $s$, refutes some
numeral instance of $u$.  A proof certificate lives inside the truth
clause, and that placement of the certificate gives \rga its name.  The companion paper also
defines a \emph{witnessed} evaluation relation $t \wbev{s} r$, whose
universal clauses add the semantic facts that the certificates
promise: every numeral instance of $u$ is true, an infinite family of
premises, and the refuted instance is false.  The companion paper
proves that on closed terms the two relations coincide.  The proof of
that coincidence takes \rga's soundness: the
certificates deliver what the added premises assert.  This paper works
with $\wbev{s}$ throughout.  A term is \emph{grounded-true} when $t
\wbev{s} \ctrue$ for some $s$, and grounded-false when $t \wbev{s} \cfalse$
for some $s$.  There is no circularity
here: the proof system is defined first, syntactically, and the
semantics only checks certificates against it.  The certificate also
lets a fuel-indexed decider approximate the semantics: the decider
never evaluates the infinite family, it is sound at each fuel,
complete up to fuel and monotone in the fuel, and so grounded truth is
a $\Sigma_1$ condition, once the certificates are counted as part of
the evaluation (\S\ref{sec:inside-dec}, ``The certificate grain'').  This paper's
internalization rests on that decider.  Grounded truth is recursively
enumerable, like the set of \rga's theorems.

\paragraph{The proof system.}
The judgment $\Gamma \vdash c$ says that $c$ is derivable from the
finite set of hypotheses $\Gamma$.  The rules are recognizable, but
three differences from classical logic shape every internal argument
in this paper.
\begin{itemize}
\item \emph{Case splits and discharge need decidedness.}  There is no
  rule concluding $a \lor \neg a$ for arbitrary $a$.  \rga reasons by
  cases on a decided sentence, and discharges a hypothesis $a$ into an
  implication $a \limp b$ only when $a$ is provably decided, that is,
  when \rga derives $a \jbool$, which abbreviates $a \lor \neg a$.  \rga introduces $a \limp b$ freely from
  $\neg a$ or from $b$, and it eliminates any derived disjunction by
  cases.  Hypothesizing is free; discharging is
  not.  So the deduction theorem, which turns a derivation of $b$ under
  the hypothesis $a$ into the implication $a \limp b$, holds exactly at
  decided hypotheses (\S\ref{sec:inside-codes}), and fails at
  ungrounded ones: the Liar entails itself, $L \vdash L$, but
  $\emptyset \nvdash L \limp L$.
  Double-negation elimination holds, so \rga is not an intuitionistic
  system either.
\item \emph{Being a number is earned.}  The judgment $t \jnat$ (``$t$
  is a natural number'') abbreviates $t = t$ and holds when $t$
  evaluates to a numeral; it plays the role of the definedness
  predicate of the logic of partial terms~\cite{beeson85foundations}.
  Universal elimination and induction carry $\jnat$ premises where
  classical rules assume definedness silently.
\item \emph{Witnesses are real.}  A derivable disjunction has a
  derivable disjunct, and a derivable negated universal has a derivable
  counterexample at some numeral (Theorem~\ref{thm:br-mech}).  The
  G\"odel paper shows which forms of G\"odel's proof survive these two
  witness properties ([\citen{goedel-paper}, Section ``G\"odel's proof in four forms'']).
\end{itemize}
Figure~\ref{fig:bg-rules} shows the main rules that later proofs use,
each with its Isabelle name.  The box in \S\ref{sec:inside-codes} lists the
companion paper's facts about the proof system as a whole, such as
soundness and the disjunction property.  The conclusion of a
derivable judgment is always well-formed, so in the rules only side
premises carry a well-formedness condition.

\begin{figure}[tp]
\centering
\small
\newcommand{\fr}[3]{\ensuremath{\dfrac{#2}{#3}}\;{\scriptsize\textsf{#1}}}
\newcommand{\fd}[3]{\ensuremath{\dfrac{#2}{#3}}\;{\scriptsize\textsf{#1}$^\dagger$}}
\newcommand{\fhd}[1]{\multicolumn{2}{@{}l@{}}{\textit{#1}}}
\renewcommand{\arraystretch}{1.0}
\begin{tabular}{@{}>{\centering\arraybackslash}p{0.48\linewidth}>{\centering\arraybackslash}p{0.48\linewidth}@{}}
\hline
\multicolumn{2}{@{}c@{}}{$a\jnat \ldef a = a \qquad
  a\jbool \ldef a \lor \neg a \qquad
  a \limp b \ldef \neg a \lor b \qquad
  \forall.\,b \ldef \cA\cdot(\lambda b)$} \\
\hline
\fhd{Structural} \\[0.5ex]
\fr{pr\_hyp}{\wf\,c}{c, \Gamma \vdash c} &
\fr{pr\_cut}{\Gamma \vdash a \quad a, \Gamma \vdash c}{\Gamma \vdash c} \\[2.6ex]
\multicolumn{2}{@{}c@{}}{\fr{pr\_sub1}{\Gamma \vdash c \quad
  \wf\,\Gamma,\ c,\ s}{\ssub{v_i}{s}{\Gamma} \vdash \ssub{v_i}{s}{c}}} \\[2.6ex]
\fhd{Negation and disjunction} \\[0.5ex]
\fr{pr\_negE}{\Gamma \vdash \neg a \quad \Gamma \vdash a \quad \wf\,c}{\Gamma \vdash c} &
\fr{pr\_dnegE}{\Gamma \vdash \neg\neg a}{\Gamma \vdash a} \\[2.6ex]
\fr{pr\_disjI1}{\Gamma \vdash a \quad \wf\,b}{\Gamma \vdash a \lor b} &
\fr{pr\_disjI2}{\Gamma \vdash b \quad \wf\,a}{\Gamma \vdash a \lor b} \\[2.6ex]
\fr{pr\_disjE1}{\Gamma \vdash a \lor b \quad a, \Gamma \vdash c \quad
  b, \Gamma \vdash c}{\Gamma \vdash c} &
\fr{pr\_disjE1 at $b = \neg a$}{\Gamma \vdash a\jbool \quad a, \Gamma \vdash c \quad
  \neg a, \Gamma \vdash c}{\Gamma \vdash c} \\[2.6ex]
\fhd{Implication (all derived)} \\[0.5ex]
\fd{ir\_impI\_bool}{\Gamma \vdash a\jbool \quad a, \Gamma \vdash c}{\Gamma \vdash a \limp c} &
\fd{ir\_mp}{\Gamma \vdash a \limp c \quad \Gamma \vdash a \quad \wf\,c}{\Gamma \vdash c} \\[2.6ex]
\fd{ir\_impI1}{\Gamma \vdash \neg a \quad \wf\,c}{\Gamma \vdash a \limp c} &
\fd{ir\_impI2}{\Gamma \vdash c \quad \wf\,a}{\Gamma \vdash a \limp c} \\[2.6ex]
\fhd{Numbers, equality and conversion} \\[0.5ex]
\fr{pr\_0}{}{\Gamma \vdash \z\jnat} \qquad
\fr{pr\_suc\_TI}{\Gamma \vdash a\jnat}{\Gamma \vdash \suc a\jnat} &
\fr{pr\_eqTI}{\Gamma \vdash a\jnat \quad \Gamma \vdash b\jnat}{\Gamma \vdash (a = b)\jbool} \\[2.6ex]
\fr{pr\_conv}{\Gamma \vdash c' \quad c \leadsto c'}{\Gamma \vdash c} &
\fr{pr\_convE}{\Gamma \vdash c \quad c \leadsto c'}{\Gamma \vdash c'} \\[2.6ex]
\multicolumn{2}{@{}c@{}}{\fr{pr\_eq\_subs}{\Gamma \vdash a = b \quad
  \Gamma \vdash \ssub{v_i}{a}{p} \quad \wf\,p}{\Gamma \vdash \ssub{v_i}{b}{p}}} \\[2.6ex]
\multicolumn{2}{@{}c@{}}{\fr{pr\_R\_suc}{\Gamma \vdash m\jnat \quad
  \Gamma \vdash (h\cdot m\cdot(\cR\cdot g\cdot h\cdot m))\jnat \quad \wf\,g}
  {\Gamma \vdash \cR\cdot g\cdot h\cdot \suc m = h\cdot m\cdot(\cR\cdot g\cdot h\cdot m)}} \\[2.6ex]
\fhd{Quantifiers} \\[0.5ex]
\fr{pr\_allI}{\vz\jnat, \lft\Gamma \vdash (\lft u)\cdot\vz}{\Gamma \vdash \cA\cdot u} &
\fd{ir\_allI\_open}{\vz\jnat, \lft\Gamma \vdash b}{\Gamma \vdash \forall.\,b} \\[2.6ex]
\fr{pr\_allE}{\Gamma \vdash \cA\cdot u \quad \Gamma \vdash a\jnat}{\Gamma \vdash u\cdot a} &
\fd{gd2\_allE\_beta}{\Gamma \vdash \forall.\,b \quad \Gamma \vdash a\jnat}{\Gamma \vdash b[a]} \\[2.6ex]
\fr{pr\_not\_allI}{\Gamma \vdash \neg(u\cdot a) \quad \Gamma \vdash a\jnat}{\Gamma \vdash \neg(\cA\cdot u)} &
\fr{pr\_not\_allE}{\Gamma \vdash \neg(\cA\cdot u) \quad
  \vz\jnat, \neg((\lft u)\cdot\vz), \lft\Gamma \vdash \lft c}{\Gamma \vdash c} \\[2.6ex]
\fhd{Induction} \\[0.5ex]
\multicolumn{2}{@{}c@{}}{\fr{pr\_nat\_ind}{\Gamma \vdash \ssub{v_i}{\z}{c} \quad
  v_i\jnat, c, \Gamma \vdash \ssub{v_i}{\suc v_i}{c} \quad
  \Gamma \vdash a\jnat \quad \wf\,c \quad v_i \text{ not free in } \Gamma}
  {\Gamma \vdash \ssub{v_i}{a}{c}}} \\[1ex]
\hline
\end{tabular}
\caption{The main rules of \rga used in later proofs, with their
  Isabelle names.  \rga has 42 primitive rules in all, two of which,
  the directions of $=$refl, are identities, since $a\jnat$ abbreviates
  $a = a$; the few other rules
  cited (weakening \textsf{pr\_wk1}, the truth-value rules
  \textsf{pr\_true\_eqI}, \textsf{pr\_true\_eqE} and
  \textsf{pr\_false\_eqE}, and \textsf{pr\_conv\_negE}) are named where
  used.  Rules marked $\dagger$ are derived; the others are primitive
  (\texttt{RGA/RGA\_Proof.thy}, where the derived rules are proved, with
  \textsf{ir\_impI\_bool} and \textsf{ir\_mp} in
  \texttt{RGA/RGA\_Sound.thy}).  A primitive rule \textsf{pr\_}$x$ is
  stated for proof lists, and its form for derivability is
  \textsf{ir\_}$x$ (\textsf{ir\_0N} for \textsf{pr\_0}).  $\Gamma$
  is a finite set of hypotheses, and $a, \Gamma$ adds $a$ to $\Gamma$.
  $\wf\,c$ says that $c$ is well-formed.  $\vz$ is
  the de Bruijn variable~0, and ${\lft t}$ lifts the free variables of
  $t$ by one, so $\vz$ is fresh in $\lft\Gamma$, $\lft u$ and $\lft c$:
  the lift does the work of the classical eigenvariable condition.
  $\ssub{v_i}{a}{c}$ substitutes $a$ for $v_i$ in $c$, and $b[a]$ is
  the $\beta$-instance.  $c \leadsto c'$ says that $c$ head-steps to
  $c'$ (\textsf{hstep}): $(\lambda b)\cdot t \leadsto b[t]$,
  $\cR\cdot g\cdot h\cdot\z \leadsto g$ and $\cR\cdot g\cdot h\cdot\suc\tofn
  n \leadsto h\cdot\tofn n\cdot(\cR\cdot g\cdot h\cdot\tofn n)$ at a
  numeral, each step when $t$, $h$ and $g$ respectively are well-formed; and
  $t\cdot u \leadsto t'\cdot u$ when $t \leadsto t'$.  The rule \textsf{pr\_disjE1} at $b = \neg a$ is the
  case split on a decided sentence.}
\label{fig:bg-rules}
\end{figure}

\paragraph{The companion paper's results.}
The companion paper~\cite{rga-paper} proves in Isabelle/HOL that every
derivable sentence is grounded-true (soundness), so that \rga does not
derive $\z = \tofn 1$ (consistency); that for well-formed terms
derivability from no hypotheses coincides with grounded truth (open
completeness, where an open term counts as grounded-true when its
evaluation, with its free variables stuck, yields true); that every
derivable totality claim is backed by a value; and that the sets \rga
represents (those $p$ with $p(n)$ iff $\vdash t\cdot\tofn n$) are
exactly the recursively enumerable sets.  Open completeness is not
negation completeness, and it is weak: it fails under hypotheses.  The
false equation $\z = \tofn 1$ is vacuously a grounded consequence of
the consistency sentence, which is never grounded-true, yet \rga does
not derive $\z = \tofn 1$ from the consistency sentence.  For a closed hypothesis $a$, every
closed well-formed grounded consequence of $a$ is derivable from $a$
exactly when $a$ is grounded-true or $a \vdash \z = \tofn 1$
(an argued remark, not a mechanized theorem~\cite[Sect.~6]{rga-paper}).  \rga is thus consistent and complete
for grounded truth.  That completeness does not escape G\"odel's theorems,
which hold for \rga ([\citen{goedel-paper}, Section ``G\"odel's theorems for \rga'']).  The sentences that classical reasoning
calls true but unprovable, such as the G\"odel sentence and the
consistency sentence, are not grounded-true, so completeness for
grounded truth does not reach them.  The classical argument concludes
that the G\"odel sentence is true because each of its numeral
instances is true.  Grounded truth makes a universal true only when a
proof certifies the universal's schematic instance, so the G\"odel and consistency
sentences, each of whose numeral instances \rga proves, are ungrounded
rather than true; read instance by instance, in an $\omega$-semantics
for which \rga is also sound (\S\ref{sec:inside-omega}), these two
sentences are true and unprovable, as classical reasoning says.  In the formalization
that the two papers share, \rga proves, of each closed ground-type
term of G\"odel's System~T, that this term denotes a natural number, one
derivation per term~\cite[Thm.~5.7]{rga-paper}, by normalization of
System~T in HOL; that result measures the reach of \rga's derivations,
not \rga's deductive strength.  \rga is not weak in
the sense of the self-verifying tradition~\cite{willard01self}: it
proves multiplication, each compiled primitive-recursive function and
Ackermann's function total by quantified
theorems~\cite[Thms.~5.4--5.6]{rga-paper}.  What \rga changes is the
reading of the universal quantifier, which a proof certificate
settles, and the conditional, which discharges only decided
hypotheses.  \rga does not
contain Heyting arithmetic, though: an implication between sentences
that \rga neither proves nor refutes, even $A \limp A$, is not a
theorem (Theorem~\ref{thm:br-mech}(3)).  In this paper, ``in HOL'' refers to
statements about \rga proved in Isabelle/HOL, and ``in \rga'' to
derivations in \rga's own proof system.

\section{RGA inside RGA}
\label{sec:inside}

Section~\ref{sec:intro} described \rga as a closed universe of
computational reasoning, in which syntax and semantics are both
computation.  This section and the next exhibit what that description means: the
arithmetization behind G\"odel's theorems (coding, compiling, deciding
and checking) happens inside \rga, as \rga terms that \rga reasons
about.

This section builds, inside \rga, the objects that the rest of the
paper reasons about.  The objects are codes of terms and proofs
(\S\ref{sec:inside-codes}), primitive-recursive functions compiled
into \rga terms (\S\ref{sec:inside-cmp}), a decider for \rga's
step-indexed semantics (\S\ref{sec:inside-dec}), a truth predicate
built around the compiled decider (\S\ref{sec:inside-truth}), and a
proof checker (\S\ref{sec:inside-chk}).  For each object the section
states what \rga proves about it, and at which arguments.  Two short
subsections then state what the development proves about the machines
at closed numerals (\S\ref{sec:inside-d3}), and describe a mechanical
check that no machine contains an opaque index
(\S\ref{sec:inside-audit}).  Appendix~\ref{app:machines} gives the
details.

The arithmetization itself is ordinary.  Coded syntax, a verified
compiler, an evaluator and a proof checker are the standard machinery
behind any treatment of G\"odel's theorems, and Isabelle/HOL checks
the construction carried out here.  Two features are specific to \rga.
The first feature is the proof certificate in the semantics' clause for the
universal quantifier, which lets a fuel-indexed decider approximate
truth one fuel bound at a time (\S\ref{sec:inside-dec}).  The second feature is paracompleteness: \rga
discharges a hypothesis (implication introduction from a derivation
under it, or case analysis on it) only when the hypothesis is provably
decided,
and every internal argument below is shaped by that discipline.

Theorem headers carry Isabelle names, and Appendix~\ref{app:endpoints}
maps every result of the paper to its name and theory file.

\subsection{Codes, numerals, and four kinds of statement}
\label{sec:inside-codes}

\paragraph{Codes.}
Every \rga term $t$ has a \emph{code} $\quo{t}$, a natural number
built with the Cantor pairing function
$\cpair{x}{y} = (x+y)(x+y+1)/2 + y$, with tuples nested to the right.
A variable $v_i$ has code $\cpair{0}{i}$.  An operator node has code
$\cpair{1}{\cpair{o}{\cpair{b}{\ell}}}$, where $o$ numbers the
operator (the ten operators $\bot, \z, \suc, {=}, \neg, \lor,
\lambda, \cdot, \cR, \cA$ in that order, from 0 to 9), $b$ counts
the variables the operator binds, and $\ell$ codes the list of the
node's arguments.  A list has code $0$ when empty and
$\cpair{x}{l} + 1$ for an element with code $x$ in front of a list
with code $l$.  Proofs receive codes the same way.  A proof is a
list of judgments, and the code of a judgment $H \vdash c$ is the pair
of the codes of the hypothesis set $H$ (as a sorted list) and of $c$.
The coding is injective but not onto: for example, no term has a
code whose first Cantor component is $2$.

\paragraph{Numerals and code numerals.}
The numeral $\tofn{n}$ is the term $\suc^n\z$.  The \emph{code
numeral} of a term $t$ is the numeral of its code,
$\bquo{t} = \tofn{\quo{t}}$.  The paper keeps codes and code numerals apart
typographically.  Plain corners $\quo{t}$ denote a natural number,
an object of the metalanguage HOL.  Underlined corners $\bquo{t}$, underlined like
any numeral, denote a closed \rga term, which \rga's own derivations
can mention.  Every argument of an internal predicate below is a term,
so codes appear there underlined.  A HOL-level function such as the packed
verdict of \S\ref{sec:inside-dec} takes a plain code.

\paragraph{Example.}
The term $\z$ has code
$\quo{\z} = \cpair{1}{\cpair{1}{\cpair{0}{0}}} = 4$, and
$\quo{\suc\z} = \cpair{1}{\cpair{2}{\cpair{0}{\quo{[4]}}}} =
5\,247\,178$, where $[4]$ is the one-element list holding the code of
$\z$.  The code of the false equation $\z = \tofn 1$ has 206
decimal digits.  The equation's code numeral $\bquo{\z = \tofn 1}$ is a term
with about $5\cdot 10^{205}$ successors (the code itself), so the
paper never writes a code out.

\paragraph{Four kinds of statement.}
Every theorem in this paper is a theorem of HOL about \rga.  What
varies is where the quantifier over codes sits.
\begin{itemize}
\item A \emph{schematic} theorem quantifies in HOL.  For every HOL
  object $t$ (a term, a proof, a number) that meets the theorem's HOL
  hypotheses, \rga derives one judgment about the code numeral
  $\bquo{t}$.  The theorem supplies one derivation per instance, not a
  single sentence that covers the instances.  Derived rules, which say
  that \rga derives a conclusion whenever it derives the premises
  (Theorems~\ref{thm:in-cmp-rec}, \ref{thm:in-wtrue-intro}
  and~\ref{thm:in-chk-dec}), also count as schematic.
  \S\ref{sec:pillars-axioms} separates two kinds of derived rule: a
  \emph{rule in every context} derives its conclusion in every context
  $\Gamma$ in which \rga derives its premises, and a \emph{rule on
  theorems} applies only to closed theorems $\emptyset \vdash c$.
\item An \emph{internal} theorem states that \rga derives one
  particular sentence, in which the quantifier over codes is \rga's
  own $\cA$.  A theorem that derives one fixed sentence in every
  context $\Gamma$ (Theorem~\ref{thm:in-never}) is internal.
\item A \emph{derivation under hypotheses} states that \rga derives
  one judgment $\Gamma \vdash c$ in a particular non-empty context
  $\Gamma$, such as $\vz\jnat,\ H(\vz) \vdash \neg\Chk(\vz,
  \bquo{\z = \tofn 1})$.  Such a derivation is not the
  implication $\emptyset \vdash H(\vz) \limp \neg\Chk(\vz, \bquo{\z =
  \tofn 1})$, and becomes one only when \rga proves the hypothesis
  decided (below).  Section~\ref{sec:con} on the restricted consistency
  theorems uses this kind of derivation twice: for the consistency of the
  level class under the closed hypothesis $\LF$ (\S\ref{sec:con-lv}), and
  for the consistency of every proof under the residual hypothesis on the
  proof code (\S\ref{sec:con-res}).
\item A \emph{meta} theorem concludes a HOL proposition about \rga
  other than a derivation: $\mathit{False}$, the existence of an
  evaluation, or that \rga does not derive some sentence.
\end{itemize}
Compare the schema ``for each $n$, $\emptyset \vdash \tofn{n}\jnat$''
with the single sentence $\forall s.\, s \jnat$: the schema is a family
of derivations, the sentence one derivation of one sentence.
The distinction matters for every claim in Sections~\ref{sec:pillars}
to~\ref{sec:con}.  ``For each proof of $c$, \rga derives that $c$
is true'' is schematic, and \S\ref{sec:pillars-sound} proves this schema,
proof by proof.
``\rga proves that everything it proves is true'' would be internal.
\S\ref{sec:pillars-con} shows that \rga does not prove even the
instance of that sentence at the false equation $\z = \tofn 1$.

\paragraph{Implication, entailment and rule.}
A conditional claim, ``$A$ gives $B$'', has four forms: implication,
decided-antecedent form, entailment and rule.  The decided-antecedent
form applies only
to a hypothesis of a common shape and is treated in the next
paragraph.  The paper always says which form it means.
\begin{itemize}
\item An \emph{implication} $\vdash A \limp B$ is one sentence.  Read
  computationally, for closed $A$ and $B$, it is one closed computation
  that terminates with the value true, and as a value it can be
  negated, disjoined and, once decided, split on inside \rga.  \rga
  discharges a hypothesis into an implication only when it can prove
  the hypothesis decided; \rga introduces $a \limp b$ freely from $\neg
  a$ or from $b$ (the box below).  Evaluating $A \limp B$ decides
  $\neg A \lor B$; the implication carries no transformation of
  evidence, which lives only in entailments and decided-antecedent
  forms.
\item An \emph{entailment} $A \vdash B$ is one derivation under the
  hypothesis $A$, the third kind of statement above; its schematic
  version is what the literature calls a derivable rule.  An entailment has two
  readings.  As a
  statement about derivations it always has content: cut turns any
  derivation of $A$, in any context $\Gamma$, into a derivation of $B$
  in $\Gamma$, and for hypotheses about proof codes this transformation
  also exists inside \rga, as an explicit primitive-recursive map of
  accepted proof codes ([\citen{goedel-paper}, theorem~\textsf{cut1M\_chk}]).  As a statement about
  truth it is soundness in context: when $A$ is grounded-true, so is
  $B$.  The development proves this preservation in HOL, by the
  soundness theorem; it extracts from the derivation no
  primitive-recursive bound on $B$'s fuel in terms of $A$'s.  At a
  hypothesis
  that is never grounded-true, such as $\ConR$ or $G$
  (\S\ref{sec:pillars-con}), the truth reading holds vacuously.  The
  entailment's content there is the derivational one, a theorem of the
  theory \mbox{\rga{}${}+A$}, and truth preservation in the $\omega$-semantics
  (\S\ref{sec:inside-omega}), in which $\ConR$ and $G$ are true and \rga
  is sound in every context.  Two-way entailment, $A \vdash B$ and $B
  \vdash A$, is written $A \eqv B$.
\item A \emph{rule}, ``if $\vdash A$ then $\vdash B$'', is a closure
  property of the theorems, a fact of HOL; its schematic version is
  what the literature calls an admissible rule.
\end{itemize}
Each form gives the next: an implication gives the entailment by
modus ponens, and an entailment gives the rule by cut.  The step back
from entailment to implication is the deduction theorem, and
Theorem~\ref{thm:in-dt} gives its exact boundary.

\begin{thm}[\textsf{deduction\_theorem\_exact}, \textsf{sigma\_deduction\_theorem}]
\label{thm:in-dt}
Let $\Gamma$ be a context.
\begin{enumerate}
\item For every well-formed $A$: $A, \Gamma \vdash B$ gives $\Gamma
  \vdash A \limp B$ for every $B$ if and only if $\Gamma \vdash A \lor
  \neg A$.  For the ordinary Liar $L$, $L \vdash L$ while $\emptyset
  \nvdash L \limp L$.
\item If $\vz\jnat, \lft\Gamma \vdash D(\vz)\jbool$, then for every
  $C$: $\exists s.\,D(s), \Gamma \vdash C$ if and only if $\Gamma
  \vdash \forall s.\,(D(s) \limp C)$.  In particular, at a closed code
  term $c$ with $\Gamma \vdash c\jnat$ and a closed $C$,
  $\PrC(c), \Gamma \vdash C$ if and only if $\Gamma \vdash \forall
  q.\,(\Chk(q, c) \limp C)$.
\end{enumerate}
\end{thm}

\noindent
In part~1 the ``only if'' half is the case $B \ldef A$, since $A
\vdash A$ always holds.  Both parts are HOL equivalences between
derivability facts.

\begin{cor}[\textsf{dec\_contrapos}, \textsf{sigma\_contrapos}]
\label{cor:in-contrapos}
Contraposition, from $A, \Gamma \vdash B$ to $\neg B, \Gamma \vdash
\neg A$, holds when $\Gamma \vdash A \lor \neg A$ and when $A$ is
$\exists s.\,D(s)$ with $D$ decided as in Theorem~\ref{thm:in-dt}(2);
in particular a hypothesis of that shape that entails $\z = \tofn 1$
is refuted.  Contraposition fails in general: $L_{\Tru} \wedge \neg L_{\Tru}
\vdash \z = \tofn 1$, while $\neg(\z = \tofn 1) \nvdash \neg(L_{\Tru}
\wedge \neg L_{\Tru})$.
\end{cor}

\noindent
Here $L_{\Tru}$ is the Liar for the truth predicate
(\S\ref{sec:pillars-tarski}).  The conclusion $\z = \tofn 1$ of the
failing instance is decided, so what matters is the hypothesis.

\paragraph{The decided-antecedent form, a fourth form.}
For a hypothesis $\exists s.\, D(s)$ with $D$ decided, an entailment
$\exists s.\, D(s) \vdash B$ is equivalent to the single sentence
$\forall s.\, (D(s) \limp B)$ (Theorem~\ref{thm:in-dt}(2)), the
\emph{decided-antecedent form} of the conditional.  Provability, truth
and certification are hypotheses of this shape: $\PrC(c) = \exists
q.\,\Chk(q, c)$ has the decided checker verdict as its matrix.  When
the entailment holds, the sentence in decided-antecedent form is an
\rga theorem, hence grounded-true in its own right, whatever the
status of $\exists s.\,D(s)$: a fourth form of the conditional,
between the implication and the entailment.  The implication gives
the decided-antecedent sentence, by cases on the disjunction
$\neg\exists s.\,D(s) \lor B$, and that sentence gives the entailment.  What fails at an undecided $\exists s.\,D(s)$
is only the step from the decided-antecedent form to the implication
$(\exists s.\,D(s)) \limp B$, the extraction of the quantifier.  At
$\ConR$ the step fails although the sentence holds:
$\vdash \forall q.\,(\Chk(q, \bquo{\ConR}) \limp \neg\ConR)$, while
$\nvdash (\exists q.\,\Chk(q, \bquo{\ConR})) \limp \neg\ConR$
([\citen{goedel-paper}, theorem~\textsf{rga\_goedel2\_sentence}]).  The decided forms of
{}[\citen{goedel-paper}, Section ``G\"odel's theorems for \rga''] are sentences of this shape, strengthened by
naming the consequent's witness with an explicit primitive-recursive
transformer of proofs.

At a hypothesis that is neither decided nor of the shape $\exists
s.\,D(s)$, such as $\ConR$ or $G$, entailment
and implication come apart: $\ConR \vdash G$ holds, while $\ConR \limp
G$ is not a theorem ([\citen{goedel-paper}, Section ``G\"odel's proof in four forms'']).  The G\"odel
entailments from $\ConR$ or $G$ still follow from sentences in
decided-antecedent form, by contraposition.  For example, $\neg G$ is
literally $\PrC(S_G)$, the sentence $\forall q.\,(\Chk(q, S_G) \limp
\neg\ConR)$ holds ([\citen{goedel-paper}, Appendix ``G\"odel II in full'']; at the code numeral
$\bquo{G}$ it is [\citen{goedel-paper}, theorem~\textsf{rga\_goedel2\_sentence}](3)), and $\Sigma$-contraposition (Corollary~\ref{cor:in-contrapos})
and double negation turn that sentence into $\ConR \vdash G$
([\citen{goedel-paper}, Appendix ``G\"odel II in full'']).  The converse direction, from such an
entailment back to a sentence, is not generic.

\paragraph{Why discharge needs decidedness.}
\rga's terms include fixed points, so with an unrestricted deduction
theorem Curry's sentence $C$, with $C$ equal to $C \limp \z = \tofn 1$,
would derive $\z = \tofn 1$ and with it every
sentence~\cite{curry42inconsistency}.  Since $\limp$ is material, $C
\limp \z = \tofn 1$ is $\neg C \lor \z = \tofn 1$, and $C$ is
Liar-like: the restriction to decided hypotheses is what blocks both
paradoxes.  Kleene's strong three-valued conditional is material in the
same way, so that $A \limp A$ is not valid in his
logic~\cite{kleene52introduction}; Field adds a non-material
conditional to recover laws such as $A \limp A$~\cite{field08saving}, and \rga keeps the
material one.  In Curry--Howard terms, an entailment $A \vdash B$ is a
hypothetical judgment, an open term $x{:}A \vdash M{:}B$.  \rga's
implication is not a function type but the boolean implication $\neg A
\lor B$, which \rga introduces by discharge exactly when $A$ is
decided, as a proof assistant reflects a decidable proposition into a
boolean.  The decided-antecedent form $\forall s.\,(D(s) \limp B)$ is
what \rga has of the BHK reading of $(\exists s.\,D(s)) \limp B$, a
function from witnesses: one sentence that discharges the decided
$D(s)$ at every candidate witness $s$.

\paragraph{Why implications between undecided sentences fail.}
Two witness properties of \rga decide which implications \rga proves.
Both properties follow from the companion paper's soundness theorem,
and the disjunction property uses open completeness as well.
\begin{thm}[\textsf{ir\_disj\_prop}, \textsf{rga\_not\_all\_witness}]
\label{thm:br-mech}
\leavevmode
\begin{enumerate}
\item \emph{The disjunction property.}  If $\emptyset \vdash a \lor b$,
  then $\emptyset \vdash a$ or $\emptyset \vdash b$.
\item \emph{The witness property.}  If $\emptyset \vdash
  \neg\forall x.\,\varphi$, then $\emptyset \vdash \neg\varphi(\tofn n)$
  for some numeral $\tofn n$.
\item If $\emptyset \nvdash \delta$ and $\emptyset \nvdash \neg\delta$,
  then $\emptyset \nvdash \delta \limp \delta$ and $\emptyset \nvdash
  \delta \liff \delta$.
\end{enumerate}
\end{thm}

\noindent
Since $a \limp b$ is $\neg a \lor b$, part~1 makes an implication $a
\limp b$ an \rga theorem only when \rga refutes $a$ or proves $b$.  So
no implication between two sentences that \rga neither proves nor
refutes is a theorem, not even $\delta \limp \delta$ (part~3), while
the entailment $\delta \vdash \delta$ holds at every well-formed
$\delta$ by the hypothesis rule.  The G\"odel paper follows this
mechanism through each step of G\"odel's proof
([\citen{goedel-paper}, \S\,``One mechanism'']).

\paragraph{A caution.}
Interderivability is weaker than intersubstitutivity, which would let
$A$ and $B$ replace each other inside any sentence; for the truth
predicate, intersubstitutivity fails under negation
(\S\ref{sec:pillars-axioms}).

\paragraph{One map of the distinctions.}
Two questions classify every conditional claim in this paper.  The
first question is the claim's form: implication, decided-antecedent form, entailment,
or rule, where a rule either holds in every context (a derivable
rule) or applies only to theorems (an admissible rule).  The second question
is where the quantifier over codes sits: in HOL, as a schema with one
instance per HOL object (the decided D3 at each $\varphi$,
{}[\citen{goedel-paper}, theorem~\textsf{D3\_decided\_any}]), or inside \rga, as one sentence over
codes (the decided D3 over the code, [\citen{goedel-paper}, theorem~\textsf{D3\_uniform}]).
Meta facts of HOL, such as $\nvdash \ConR$, are not a form of the
conditional; Table~\ref{tab:in-forms} lists such meta facts below a rule, beside
two conditionals followed through the four forms.

\begin{table}[tbp]
\centering
\small
\setlength{\tabcolsep}{3pt}
\begin{tabular}{@{}>{\raggedright\arraybackslash}p{0.15\linewidth}>{\raggedright\arraybackslash}p{0.17\linewidth}>{\raggedright\arraybackslash}p{0.31\linewidth}>{\raggedright\arraybackslash}p{0.31\linewidth}@{}}
\hline
form & written & from $\ConR$ to $G$ & from $\PrC(\bquo{\ConR})$ to $\neg\ConR$ \\
\hline
implication & $\vdash A \limp B$ & fails ([\citen{goedel-paper}, theorem~\textsf{no\_formalized\_goedel1}]) & fails ([\citen{goedel-paper}, theorem~\textsf{rga\_goedel2\_sentence}](5)) \\
decided-antecedent form & $\vdash \forall s.\,(D(s) \limp B)$ & not of this shape ($\ConR$ is not $\exists s.\,D(s)$); the conditional's contrapositive is: $\forall q.\,(\Chk(q, \bquo{G}) \limp \neg\ConR)$ holds ([\citen{goedel-paper}, theorem~\textsf{rga\_goedel2\_sentence}](3)) & holds: $\forall q.\,(\Chk(q, \bquo{\ConR}) \limp \neg\ConR)$ ([\citen{goedel-paper}, theorem~\textsf{rga\_goedel2\_sentence}](1)) \\
entailment & $A \vdash B$ & holds ([\citen{goedel-paper}, theorem~\textsf{goedelG\_under\_ConRGA}]) & holds ([\citen{goedel-paper}, theorem~\textsf{rga\_goedel2\_internal\_contra}]) \\
rule in every context & $\Gamma \vdash A \Rightarrow \Gamma \vdash B$ & holds, as the entailment does & holds, as the entailment does \\
rule on theorems & $\vdash A \Rightarrow {}\vdash B$ & holds, vacuously: $\nvdash \ConR$ & holds, vacuously: $\nvdash \PrC(\bquo{\ConR})$ \\
\hline
meta facts & HOL facts & $\nvdash \ConR$, $\nvdash G$ ([\citen{goedel-paper}, Section ``G\"odel's theorems for \rga'']) & $\nvdash \ConR$ \\
\hline
\end{tabular}
\caption{The four forms of a conditional (\S\ref{sec:inside-codes}),
  with two worked conditionals; a rule holds either in every context
  (a derivable rule) or on theorems only (an admissible rule).  The
  decided-antecedent form applies to a hypothesis $\exists s.\,D(s)$
  with $D$ decided, such as $\PrC(c)$; $\neg G$ is literally
  $\PrC(S_G)$, and $\vdash S_G = \bquo{G}$.  $\nvdash \PrC(\bquo{\ConR})$ follows from $\nvdash
  \ConR$ by the reflection rule ([\citen{goedel-paper}, theorem~\textsf{D1\_rule\_pr}]).  Meta
  facts are not a form of the conditional.  Each form can be a HOL
  schema or one internal sentence over codes.}
\label{tab:in-forms}
\end{table}

\paragraph{Notation.}
The paper uses the notation of the companion
paper~\cite{rga-paper}: the dot $\cdot$ is \rga application, $\vz$ is the de
Bruijn variable $0$, and ${\lft t}$ lifts the free variables of $t$
by one.  The quantifier shorthands are $\forall.\,b \ldef \cA \cdot
(\lambda b)$ and $\exists.\,b \ldef \neg(\cA \cdot (\lambda \neg
b))$; with a named variable, $\forall s.\,b$ means the same with $s$
standing for the bound index.  The term $a \jnat$ (``$a$ is a
natural number'') abbreviates $a = a$, and $a \jbool$ (``$a$ is
decided'') abbreviates $a \lor \neg a$.  An implication $a \limp b$
abbreviates $\neg a \lor b$.  \rga can prove such an implication by discharging the
hypothesis $a$ (the deduction step) only when $a$ is decided; the box
below lists this restriction and the other companion-paper facts that later proofs
use.

\begin{center}
\fbox{\begin{minipage}{0.95\linewidth}
\small
\textbf{Facts imported from the companion paper~\cite{rga-paper}}
(Isabelle names and theory files in Appendix~\ref{app:machines}).
\begin{itemize}
\item \emph{Values and equations.}  Values are natural numbers
  (Table~\ref{tab:in-sem}).  An equation $a = b$ is grounded-true when
  its two sides evaluate to the same number, and grounded-false when
  they evaluate to different numbers.  So $a \jnat$, which is $a = a$,
  is grounded-true exactly when $a$ has a value, and $a \jbool$ exactly
  when $a$ evaluates to $\ctrue$ or $\cfalse$.
\item \emph{Implication.}  \rga introduces $a \limp b$ from
  $\Gamma \vdash \neg a$ or from $\Gamma \vdash b$, with no condition
  beyond the well-formedness of the other side.  Discharging a
  hypothesis, from $a, \Gamma \vdash b$ to $\Gamma \vdash a \limp b$,
  needs $\Gamma \vdash a \jbool$.
\item \emph{Soundness}: if $\emptyset \vdash c$, then $c$ is
  grounded-true.  \emph{Consistency}: \rga never derives both $c$ and
  $\neg c$.
\item \emph{Open completeness}: every well-formed term, open or
  closed, that evaluates to $\ctrue$ is derivable from $\emptyset$ (an
  open term evaluates with its free variables stuck).  The
  value form of open completeness: for a well-formed $t$, $\emptyset \vdash t = \tofn{r}$
  exactly when $t$ evaluates to $r$.  So for well-formed terms,
  derivability from $\emptyset$ and grounded truth coincide.
\item \emph{Unique values}: a term has at most one value, whatever
  the fuels.
\item \emph{Weakening}: $\emptyset \vdash c$ implies $\Gamma \vdash
  c$.
\item \emph{Disjunction property}: if
  $\emptyset \vdash a \lor b$, then $\emptyset \vdash a$ or $\emptyset
  \vdash b$.
\end{itemize}
\end{minipage}}
\end{center}

\noindent
One more fact, proved in this paper rather than imported, is used
throughout.  \emph{Transport} for the truth predicate of
\S\ref{sec:inside-truth} (Appendix~\ref{app:machines-truth}): if $\Gamma$
derives $a \jnat$, $b \jnat$, $a = b$ and $\Tru(a)$, then $\Gamma
\vdash \Tru(b)$.

\begin{center}
\fbox{\begin{minipage}{0.95\linewidth}
\small
\textbf{Results taken from the G\"odel paper~\cite{goedel-paper}.}
The G\"odel paper proves G\"odel's theorems for \rga with the checker
$\Chk$ of \S\ref{sec:inside-chk}, the consistency sentence $\ConR =
\forall p.\ \neg\Chk(p, \bquo{\z = \tofn 1})$ and the provability
formula $\PrC(c) = \exists q.\ \Chk(q, c)$ of
\S\ref{sec:pillars-complete}.  This paper uses five of the G\"odel
paper's results.
\begin{itemize}
\item \emph{G\"odel's second theorem for \rga}
  (\textsf{rga\_goedel2}, [\citen{goedel-paper}, theorem~\textsf{rga\_goedel2}]):
  $\emptyset \nvdash \ConR$.
\item \emph{The consistency sentence is unrefutable}
  (\textsf{ConRGA\_unrefutable}, [\citen{goedel-paper}, theorem~\textsf{ConRGA\_unrefutable}]):
  $\emptyset \nvdash \neg\ConR$.
\item \emph{The reflection rule} (\textsf{rga\_pr\_reflect},
  {}[\citen{goedel-paper}, theorem~\textsf{D1\_rule\_pr}](4)): if $\emptyset \vdash
  \PrC(\bquo{\psi})$, then $\emptyset \vdash \psi$.
\item \emph{The formalized second theorem, as an entailment}
  (\textsf{rga\_goedel2\_internal}, [\citen{goedel-paper}, theorem~\textsf{rga\_goedel2\_internal\_contra}](2)):
  $\ConR \vdash \neg\PrC(\bquo{\ConR})$.
\item \emph{The G\"odel sentence}: $G \ldef \forall q.\ \neg\Chk(q,
  S_G)$ for a closed code term $S_G$ with $\emptyset \vdash S_G =
  \bquo{G}$ (\textsf{goedelS\_eq}, [\citen{goedel-paper}, theorem~\textsf{goedel\_diag}]), and
  $\emptyset \nvdash G$ (\textsf{rga\_goedel1},
  {}[\citen{goedel-paper}, theorem~\textsf{rga\_goedel1}]).
\end{itemize}
Two of the three unprovability claims of Theorem~\ref{thm:br-omega}
on $\omega$-incompleteness, at $G$ and at $\ConR$, are the last and the
first of these results.  The G\"odel paper's proof of the second
theorem uses neither the decider, nor the truth predicate, nor the
internal completeness of \S\ref{sec:pillars-int}, so the proof of
Theorem~\ref{thm:pl-no-uniform} in \S\ref{sec:pillars-con} can use the
second theorem without circularity.
\end{minipage}}
\end{center}

\subsection{Compiled primitive-recursive machines}
\label{sec:inside-cmp}

\paragraph{Index trees and the compiler.}
The development describes binary primitive-recursive (PR) functions
by \emph{index trees}.  The leaves are four initial functions (zero,
successor, and the left and right projections).  The node
$\iSub\,g\,h_1\,h_2$ composes, computing $g(h_1(x,y), h_2(x,y))$.
The node $\iRec\,g\,h\,h'$ recurses on the second argument:
$f(x,0) = g(x,0)$ and $f(x,k{+}1) = h(h'(x,k), f(x,k))$.  An index
$f$ denotes the binary function $\den{f}$.  The companion paper's
compiler $\cmp$ maps each index tree $f$ to a closed \rga term
$\cmp(f)$, one $\lambda$-term per node, with \rga's recursor $\cR$ at
each $\iRec$ node~\cite{rga-paper}.  Evaluating $\cmp(f)\cdot a\cdot
b$ in the semantics computes $\den{f}$ (Appendix~\ref{app:machines}).
\rga can then prove two things about every compiled function.

\begin{thm}[\textsf{cmp\_eq\_num}, \textsf{cmp\_total}]
\label{thm:in-cmp-num}
\label{thm:in-cmp-total}
For every index $f$:
\begin{enumerate}
\item for all naturals $x$ and $y$, $\emptyset \vdash \cmp(f) \cdot
  \tofn{x} \cdot \tofn{y} = \tofn{\den{f}(x,y)}$;
\item $\emptyset \vdash \forall.\,\forall.\,(\cmp(f)\cdot v_1\cdot
  v_0) \jnat$.
\end{enumerate}
\end{thm}

\noindent
Both parts are schemata over $f$: for each PR index $f$, \rga proves
$\cmp(f)$ total.

The proof of Theorem~\ref{thm:in-cmp-num} builds no internal
computation trace.  The proof evaluates the closed term in the HOL semantics
and turns the evaluation into a derivation with the companion paper's
open completeness theorem.  For this reason the theorem holds for
every index, however the development obtained the index.

\paragraph{Transparent and opaque indices.}
An index is \emph{explicit} when the development writes it out as a
tree of initial functions, $\iSub$ and $\iRec$ nodes;
\S\ref{sec:inside-audit} defines the checked property that the paper
calls transparent.  An index is \emph{opaque} when Hilbert's choice operator
picks it: for a function $F$ known to be PR, ``some index $i$ with
$\den{i} = F$''.  The only fact about such an index is that it
denotes $F$.  Both kinds compile, and at numerals both obey
Theorem~\ref{thm:in-cmp-num}.  At abstract arguments they differ.
Reasoning about $\cmp(f)\cdot m\cdot k$ for terms $m$ and $k$ that
are not numerals unfolds the compiled term one node of the tree at a
time, so the node must be visible.  ``Explicit'' means choice-free and
definitionally transparent, not of feasible size.  A code enters a tree
as a unary numeral, built by one successor node per unit, and each
node of that construction duplicates its argument, so the tree for
the code $k$ has about $2^k$ nodes (a doubling construction would make
the size logarithmic in $k$; no statement depends on the size).  The
transformers of
{}[\citen{goedel-paper}, Section ``G\"odel's theorems for \rga''] embed codes with hundreds of decimal digits,
so their trees are astronomically large.  Like the table machine of
\S\ref{sec:inside-dec}, these trees serve definability and proof, not
execution.  The unfolding laws read:

\begin{thm}[\textsf{cmp\_prRec\_zero}, \textsf{cmp\_prRec\_suc\_any}]
\label{thm:in-cmp-rec}
Let $g$, $h$, $h'$ be indices and $f = \iRec\,g\,h\,h'$.  If
$\Gamma \vdash m \jnat$, then
$\Gamma \vdash \cmp(f)\cdot m\cdot \z = \cmp(g)\cdot m\cdot \z$.
If moreover $\Gamma \vdash k \jnat$, then
\[
  \Gamma \vdash \cmp(f)\cdot m\cdot \suc k
    = \cmp(h)\cdot(\cmp(h')\cdot m\cdot k)\cdot(\cmp(f)\cdot m\cdot k) .
\]
\end{thm}

\noindent
An opaque index is still a definite tree, and $\cmp$ of it is a
definite closed term.  What is hidden is the tree's identity from the
development's HOL proofs, which cannot name the tree.  The development
knows only that $\cmp(f)$ is total and what it computes at each pair
of numerals, so a proof at abstract arguments has no node to unfold.
A programmer will recognize a sealed module: the interface (totality,
numeral values) is available, the implementation is not, and any
argument that depends on the implementation cannot be written.  No
theorem of the development says that the unfolding laws fail for an
opaque index; they are unavailable, not refuted.

The distinction decides which internal arguments are available.
Theorem~\ref{thm:in-never} below proves a statement about the decider
at every fuel, by an induction whose step unfolds the decider one level
at an abstract fuel.  For a decider compiled from an opaque index,
totality and numeral values give the induction step nothing to work
with, although that decider and the explicit-tree decider agree at every numeral fuel and every
term-code numeral (Appendix~\ref{app:machines}).  The machines that
later arguments reason about at abstract arguments are therefore
explicit trees; \S\ref{sec:inside-audit} describes the check that
confirms this property.

Logicians know this phenomenon as the intensionality of
arithmetization.  Feferman showed that a theory's consistency
statement, and with it G\"odel's second theorem, depends on how the
theory is presented inside arithmetic, not only on which theorems it
has~\cite{feferman60arithmetization}.  In this paper explicit trees
fix the presentation, and \S\ref{sec:inside-audit} checks mechanically
that $\ConR$, the truth predicate and the other transparent machines
are built from them.  The check is necessary but not sufficient for a
natural consistency statement.  A checker that rejects the code
$\bquo{\z = \tofn 1}$ outright would be explicit, and faithful and
genuine in the sense of Theorem~\ref{thm:in-chk}, yet \rga would prove
the consistency statement built from that checker trivially.  What makes $\ConR$ the
consistency statement is that \rga proves the checker's structural
laws at abstract arguments: the laws by which an accepted table is
extended, cited and combined (Appendix~\ref{app:machines} and
{}[\citen{goedel-paper}, Appendix ``G\"odel II in full'']), the analogue of the classical requirement that
the derivability conditions be provable for the provability
predicate.  \rga proves those conditions for $\Chk$, and the
development collects them in one theorem ([\citen{goedel-paper}, \S\,``The derivability conditions \rga has'']).
The G\"odel paper's section on G\"odel's theorems for \rga
{}[\citen{goedel-paper}, Section ``G\"odel's theorems for \rga''] uses exactly
these laws.

\subsection{A fuel-indexed decider of the semantics}
\label{sec:inside-dec}

\paragraph{The semantics.}
The companion paper defines truth by a step-indexed evaluation
relation whose universal clauses carry only the proof certificate,
and proves that relation equal on closed terms to a \emph{witnessed} relation
$t \wbev{s} r$, which adds the instance family (requirement~(iii)
below)~\cite{rga-paper}.  In the development, $t \wbev{s} r$ is the
instance at \rga's own proofs of an evaluation parametric in the system
whose proofs serve as certificates.  The facts this paper imports
are stated for $t \wbev{s} r$ itself, under the names given here, and
a build-time check compares the statement of every named fact of
\rga's soundness, completeness and Ackermann theories with a recorded
table.
The same construction at the proofs of IGA, the extension of \rga by
ground induction, makes IGA sound and complete for a semantics of the
same kind~\cite{iga-paper}.  This paper works with $t \wbev{s} r$: within
fuel $s$, the term $t$ evaluates to the value $r$, a natural number.
A second relation $t \whev{s} h$ reports the \emph{head form} $h$ of $t$: an
abstraction $\lambda b$, the recursor $\cR$ applied to at most two
arguments, or $\cA$.  Table~\ref{tab:in-sem} lists the clauses.
Evaluation is call by name: a $\beta$-step substitutes the argument
unevaluated, and only the recursor's counter argument is evaluated to
a numeral.  A variable has neither a head form nor a value, so a term
whose evaluation needs one of its free variables has no value.

Truth values are numerals, $\ctrue = 1$ and $\cfalse = 0$.  A term is
\emph{grounded} when $t \wbev{s} r$ for some $s$ and $r$, and
\emph{ungrounded} when it has no value at all.  A term is
\emph{grounded-true} when its value is $\ctrue$ and
\emph{grounded-false} when its value is $\cfalse$.  A negation, a
disjunction, an equation or a universal can only take the value $\ctrue$
or $\cfalse$, so such a term is ungrounded exactly when it is neither
grounded-true nor grounded-false.

\begin{table}[tbp]
\centering
\small
\begin{tabular}{@{}>{\raggedright\arraybackslash}p{0.27\linewidth}>{\raggedright\arraybackslash}p{0.25\linewidth}>{\raggedright\arraybackslash}p{0.40\linewidth}@{}}
\hline
term & head form, $\whev{s+1}$ & value, $\wbev{s+1}$ \\
\hline
$v_i$, $\bot$ & none & none \\
$\z$ & none & $0$ \\
$\suc a$ & none & $r+1$ if $a \wbev{s} r$ \\
$a = b$ & none & $\ctrue$ if $a$ and $b$ have equal values at fuel $s$,
  $\cfalse$ if their values differ \\
$\neg a$ & none & $\ctrue$ if $a \wbev{s} \cfalse$; $\cfalse$ if
  $a \wbev{s} \ctrue$ \\
$a \lor b$ & none & $\ctrue$ if $a \wbev{s} \ctrue$ or $b \wbev{s} \ctrue$;
  $\cfalse$ if both give $\cfalse$ \\
$\lambda b$, $\cR$, $\cA$ & the term itself & none \\
$t\cdot u$, head of $t$ is $\lambda b$ & the head form of $b[u]$ at
  fuel $s$ & the value of $b[u]$ at fuel $s$ \\
$t\cdot u$, head of $t$ is $\cR$ or $\cR\cdot g$ & $\cR\cdot u$,
  resp.\ $\cR\cdot g\cdot u$ & none \\
$t\cdot m$, head of $t$ is $\cR\cdot g\cdot h$ & \multicolumn{2}{>{\raggedright\arraybackslash}p{0.66\linewidth}@{}}{as
  $g$ if $m \wbev{s} 0$; as $h\cdot\tofn{k}\cdot(\cR\cdot g\cdot
  h\cdot\tofn{k})$ if $m \wbev{s} k{+}1$ (head form and value alike)} \\
$t\cdot u$, head of $t$ is $\cA$ & none & the two certificate clauses
  displayed in the text \\
\hline
\end{tabular}
\caption{The companion paper's witnessed semantics, clause by clause.
  Every clause concludes at fuel $s+1$
  from premises at fuel $s$; ``head of $t$'' means $t \whev{s} h$.
  $b[u]$ is the $\beta$-instance, $u$ substituted for the variable that
  $\lambda$ binds.}
\label{tab:in-sem}
\end{table}

The disjunction is parallel, in the strong Kleene manner, like
Plotkin's parallel-or~\cite{plotkin77lcf}: $a \lor b$ is true as soon
as either disjunct is, whichever the other is.  The
decider of this section realizes the parallel disjunction without any
interleaving, because at fuel $s+1$ it reads both disjuncts' verdicts
at fuel $s$ and affirms the disjunction when either is $\ctrue$.  That a
term has at most one value concerns values, not which disjunct
supplies a true one.

\paragraph{The $\omega$-semantics.}
\label{sec:inside-omega}
Beside grounded evaluation, the paper uses a classical
$\omega$-semantics, $\mathit{omv}$.  The $\omega$-semantics assigns values by the clauses
of Table~\ref{tab:in-sem}, except at a universal: $\cA\cdot u$ has the
value true when every numeral instance $u\cdot\tofn n$ is true, and
false when some instance is false, with no certificate.  This semantics is a
partial classical valuation, classically grounded in the sense of
\S\ref{sec:bg-rga}: a term whose evaluation passes through a universal
and returns to itself, such as $W\cdot W$ with $W = \lambda w.\
\forall x.\ w\,w$, has no $\omega$-value.  A term has at most one
$\omega$-value, and every grounded value is an
$\omega$-value~\cite[Sect.~6]{rga-paper}.  The
converse fails at true universals without a certificate: $\ConR$, $G$
and $L_{\Tru}$ are $\omega$-true and have no grounded value.  \rga is
sound for the $\omega$-semantics in every
context~\cite[Sect.~6]{rga-paper}: if $\Gamma \vdash c$
with $\Gamma$ well-formed and every member of $\Gamma$ $\omega$-true,
then $c$ is $\omega$-true.  So an entailment $A \vdash B$ fails
whenever $A$ is $\omega$-true and $B$ is not, the argument behind most
non-derivability results of Section~\ref{sec:con} and of the
G\"odel paper ([\citen{goedel-paper}, Section ``G\"odel's proof in four forms'']).  Appendix~\ref{app:machines} gives the Isabelle
names.

\paragraph{The quantifier clauses.}
The universal quantifier has two clauses, each carrying a proof
certificate.  Writing $P \proves H \vdash c$ for ``$P$ is a valid
proof whose conclusion is the judgment $H \vdash c$'', the clauses read
\[
\infrule{t \whev{s} \cA \quad
         \quo{P} \le s \quad
         P \proves \{\vz \jnat\} \vdash \lft u \cdot \vz \quad
         \forall n.\ \exists s_n.\ u\cdot\tofn{n} \wbev{s_n} \ctrue}
        {t \cdot u \wbev{s+1} \ctrue}
\]
\[
\infrule{t \whev{s} \cA \quad
         \quo{P} \le s \quad
         P \proves \emptyset \vdash \neg(u\cdot\tofn{i}) \quad
         \exists s_i.\ u\cdot\tofn{i} \wbev{s_i} \cfalse}
        {t \cdot u \wbev{s+1} \cfalse}
\]
A universal is true when (i) its head evaluates to $\cA$, (ii) a proof
$P$ with code at most $s$ proves the open instance ``$u$ holds of an
arbitrary natural $\vz$'', and (iii) every numeral instance is true.
The lift in $\lft u \cdot \vz$ is the usual shift under a binder: it
keeps $u$'s free variables apart from the new variable $\vz$.  A
universal is false when a refuted instance comes with a certificate of
its refutation.

\paragraph{The certificate grain.}
Requirement (iii) is an infinite family, and no procedure can check
it.  Requirements (i) and (ii) are bounded: running the head for $s$
steps is finite work, and so is checking each of the finitely many
numbers $p \le s$ with a PR proof checker.  The decider checks (i) and
(ii) and never evaluates an instance.  Requirement (iii) returns in
the decider's correctness proof, which applies \rga's soundness
theorem to the certificate: whatever $P$ proves about an arbitrary
natural holds of every numeral.  The machine checks certificates, and
the theorem relating machine and semantics spends them.  This paper
calls the pattern the \emph{certificate grain} of the semantics.  The
companion paper's own truth clauses omit requirement~(iii), and the
same soundness theorem is what makes that semantics coincide with the
witnessed one.

\paragraph{The decider.}
Concretely, the development defines two mutually recursive,
fuel-indexed functions in HOL: a head decider, which computes the head
form of Table~\ref{tab:in-sem}, and a value decider $\decV$.  The
verdict $\decV_s(t)$ is either $\None$ or $\Some{r}$.  At fuel $0$ the
verdict is $\None$; every recursive call at fuel $s+1$ uses fuel $s$.
At an application $t\cdot u$ whose head decider reports $\cA$ at fuel
$s$, the verdict at fuel $s+1$ is
\begin{itemize}
\item $\Some{\ctrue}$ if some $p \le s$ codes a proof $P$ with $P
  \proves \{\vz\jnat\} \vdash \lft u \cdot \vz$;
\item otherwise $\Some{\cfalse}$ if some $p \le s$ codes a proof of
  $\emptyset \vdash \neg(u\cdot\tofn{i})$ for some $i \le s$;
\item otherwise $\None$.
\end{itemize}
Appendix~\ref{app:machines-dec} gives the remaining clauses and the
range check that rejects the numbers coding no proof.

\paragraph{Fuel.}
The fuel $s$ bounds three things.  In every clause it bounds the depth
of the evaluation, since each clause spends one unit.  In the two
universal clauses it also bounds the code of the certificate,
$\quo{P} \le s$.  In the decider it bounds, in addition, the refuted
instance $i \le s$, which the semantic clause leaves unbounded.  The
family premises in (iii) carry their own fuels $s_n$, with no bound.
The bounds on depth and on certificates make each fuel slice a finite
search, and the certificate bound is why a true universal can need a
fuel far larger than any evaluation it involves (the example below).
More fuel never changes a verdict (Theorem~\ref{thm:in-dec}, part~3).
The decider is sound at its exact fuel, because each of its clauses
reads the matching semantic clause at the same fuel.  Completeness
holds only up to fuel: the semantics does not bound the refuted
instance $i$ by the fuel, and the proof aligns the fuels of
sub-verdicts by monotonicity.  So the decider does not decide the
exact fuel slice $t \wbev{s} r$; it is a decidable, fuel-monotone
approximation of the semantics, which agrees with the semantics on
whether a term
is grounded-true at some fuel.

\begin{thm}[\textsf{prDecV\_sound}, \textsf{prDecV\_complete}]
\label{thm:in-dec}
For every term $t$, value $r$ and fuels $s \le s'$:
\begin{enumerate}
\item if $\decV_s(t) = \Some{r}$, then $t \wbev{s} r$
  (sound at its exact fuel);
\item if $t \wbev{s} r$, then $\decV_{s''}(t) = \Some{r}$ for some
  fuel $s''$ (complete up to fuel);
\item if $\decV_s(t) = \Some{r}$, then $\decV_{s'}(t) = \Some{r}$
  (monotone).
\end{enumerate}
\end{thm}

\noindent
\emph{Proof idea.}  Soundness goes by induction on the fuel.  At a
quantifier node the found certificate supplies requirement (iii)
through the companion paper's soundness theorem, applied in HOL.
Completeness goes by induction on the evaluation, and monotonicity by
induction on the fuel.  The three statements are the lemmas that relate
a clocked, fuel-indexed definitional interpreter to a relational
semantics in the step-indexed treatment of functional
languages~\cite{owens16functional,amin17type}.

\emph{Classical versus \rga.}  Classically, the truth of a universal
sentence is an infinite conjunction over its instances, and arithmetic
truth is not even arithmetical, let alone decidable in bounded
slices.  \rga's
semantics demands a certificate at each true universal, and
Theorem~\ref{thm:in-dec} turns that demand into a decidable,
fuel-monotone approximation.  So grounded truth, $\exists s.\ t \wbev{s} \ctrue$,
coincides with $\exists s.\ \decV_s(t) = \Some{\ctrue}$, a $\Sigma_1$
condition.  Every later internalization spends this fact.

\paragraph{A PR table machine.}
The mutual recursion above is not PR by its shape: a $\beta$-step
substitutes into a body, so arguments grow.  The development
therefore defines a table machine $\dtab$ by primitive recursion on
the fuel.  Level $s$ of the table holds, for every code up to a bound
$B(s,x)$, the pair of the two deciders' verdicts at fuel $s$, packed
into numbers.  The packed value verdict is
$\decP_s(x) = 0$ when $\decV_s(\mathrm{decode}\,x) = \None$ and
$\decP_s(x) = r+1$ when $\decV_s(\mathrm{decode}\,x) = \Some{r}$.  So
$2$ means $\Some{\ctrue}$ and $1$ means $\Some{\cfalse}$.  The machine
$\dtab(s,x)$ reads the cell of $x$ at level $s$.

\begin{thm}[\textsf{n2pr\_prDecTab}, \textsf{prDecTab\_correct}]
\label{thm:in-tab}
The function $\dtab$ is primitive recursive, with no hypothesis.
For every term $t$ and fuel $s$, $\dtab(s, \quo{t})$ is the Cantor
pair of the packed head verdict and the packed value verdict
$\decP_s(\quo{t})$.
\end{thm}

\noindent
The correctness equation holds at term codes, not at every number.
The reasoning about how far the deciders reach lives only in the
correctness proof, so the machine itself is PR outright.  The bound
$B$ grows like a tower of exponentials in $s$.  The table is a device
for PR definability, not an algorithm anyone would run.

\paragraph{Example.}
The false equation $\z = \tofn 1$ needs fuel 1 for $\z$, fuel 2 for
$\tofn 1 = \suc\z$, and fuel 3 for the equation.  Accordingly
$\decV_2(\z = \tofn 1) = \None$ and $\decV_3(\z = \tofn 1) =
\Some{\cfalse}$, so the packed verdict is 0 at fuel 2 and 1 at fuel 3,
and by monotonicity 1 at every larger fuel.  The packed verdict is never 2.
A universal shows the grain.  For $u = \lambda(\vz \jnat)$, the
decider affirms $\cA\cdot u$, ``every natural is a natural'', at every
fuel above the code of a two-line proof of $\{\vz\jnat\} \vdash \lft
u\cdot\vz$, a code of about 3,400 decimal digits, and never evaluates
an instance $u\cdot\tofn{n}$ (Appendix~\ref{app:machines-dec}).

\subsection{The compiled decider and the truth predicate}
\label{sec:inside-truth}

\paragraph{The compiled decider.}
The development writes the table machine out as one explicit index
tree $\taut$, with $\den{\taut}(x, \quo{t}) = \dtab(x, \quo{t})$ at
every term code (Appendix~\ref{app:machines}).  The \emph{compiled
decider} is the \rga term
\[
  \Dec(s, c) \;\ldef\; \mathit{snd}(\cmp(\taut)\cdot s\cdot c) ,
\]
where $\mathit{snd}$ is the compiled second Cantor projection, applied
to the term $\cmp(\taut)\cdot s\cdot c$.  Table~\ref{tab:in-layers}
shows the layers from the semantics down to $\Dec$ and the predicates
built on it, with the theorem that links each layer to the one above;
the checkers of \S\ref{sec:inside-chk} have the same structure.

\begin{table}[t]
\centering
\small
\begin{tabular}{@{}>{\raggedright\arraybackslash}p{0.15\linewidth}>{\raggedright\arraybackslash}p{0.22\linewidth}>{\raggedright\arraybackslash}p{0.17\linewidth}>{\raggedright\arraybackslash}p{0.36\linewidth}@{}}
\hline
layer & object & lives in & link to the layer above \\
\hline
semantics & $t \wbev{s} r$ & HOL relation & (companion paper) \\
decider & $\decV_s(t)$ & HOL function & Thm.~\ref{thm:in-dec}: sound at its
  fuel, complete up to fuel, monotone \\
table machine & $\dtab(s,x)$ & HOL, PR function & Thm.~\ref{thm:in-tab}:
  PR outright; its cell at $\quo{t}$ packs $\decP_s(\quo{t})$ \\
index tree & $\taut$ & explicit PR index & $\den{\taut}(x,\quo{t}) =
  \dtab(x,\quo{t})$ at term codes (App.~\ref{app:machines}) \\
compiled decider & $\Dec(s,c)$ & closed \rga term & Thm.~\ref{thm:in-dec-closed}:
  $\emptyset\vdash\Dec(\tofn{s},\bquo{t}) = \tofn{\decP_s(\quo{t})}$,
  by Thm.~\ref{thm:in-cmp-num} \\
predicates & $\Aff(s,c)$, $\Tru(c)$ & \rga terms & Def.~\ref{def:in-truth};
  Thm.~\ref{thm:in-wtrue-intro}; adequacy, Thm.~\ref{thm:pl-adeq} \\
\hline
proofs & $P \proves H\vdash c$ & HOL relation & (companion paper) \\
detectors & rule detectors (\S\ref{sec:inside-chk}), validity test, $\pi$ & explicit PR indices &
  correct at proof codes, genuine at every number (App.~\ref{app:machines};
  Thm.~\ref{thm:pl-genuine} for $\pi$) \\
tree checker & $\Chk(p,c)$ & \rga term & Thm.~\ref{thm:in-chk}: faithful at
  proof codes, genuine at every numeral \\
predicates & $\ConR$, $\Cert(h)$ & \rga terms & \S\ref{sec:inside-chk},
  Def.~\ref{def:pl-cert} \\
\hline
\end{tabular}
\caption{The layers of the decider (top) and of the checker (bottom),
  each with the result that links it to the layer above.  ``Genuine''
  means that every accepted number codes a real proof.}
\label{tab:in-layers}
\end{table}

\begin{thm}[\textsf{iPrDecV\_closed}]
\label{thm:in-dec-closed}
For every fuel $s$ and term $t$,
$\emptyset \vdash \Dec(\tofn{s}, \bquo{t}) = \tofn{\decP_s(\quo{t})}$.
\end{thm}

\noindent
The theorem is a schema over $s$ and $t$, with no side condition.  Its
proof is Theorem~\ref{thm:in-cmp-num} at $\taut$, followed by
Theorem~\ref{thm:in-tab} and the projection.  For each $s$ and $t$,
\rga proves the value of the compiled decider at $(\tofn s, \bquo{t})$; the
right-hand side is a HOL function, so inside \rga the decider is its
own specification.

\paragraph{The truth predicate.}
Two definitions turn the compiled decider into a truth predicate.

\begin{defi}[\textsf{iPrDecTrueAt}, \textsf{iWtrue}]
\label{def:in-truth}
The \emph{affirmation test} $\Aff(s, c)$ is the 0/1-valued equality
test of $\Dec(s, c)$ against $\tofn 2$: its value is $1$ (true) when
the compiled decider affirms $c$ at fuel $s$, and $0$ otherwise.  The
\emph{truth predicate} is
\[
  \Tru(c) \;\ldef\; \neg(\cA\cdot(\lambda\neg\Aff(\vz, \lft c))) ,
  \qquad\text{that is,}\qquad
  \Tru(c) \;\ldef\; \exists s.\ \Aff(s, c) .
\]
\end{defi}

\noindent
A 0/1-valued test $t$ such as $\Aff(s, c)$, written where a formula is
expected, is true when its value is $\ctrue = 1$.  When \rga proves such
a test to be a natural number, $t$, $t = \tofn 1$ and $t \neq \z$ are
interderivable, and the paper writes whichever form the development
uses.

$\Tru(c)$ says that for some fuel the compiled decider affirms $c$.
\rga's existential is the paracomplete $\neg\forall\neg$.  This existential is true
when a witness comes with a certificate (a proof refuting the
witness's negated instance), false when the universal $\forall\neg$ is
certified (a proof of the open negated instance, and every instance
fails), and ungrounded otherwise.  Introduction takes an affirming
equation at a fuel term, which may be abstract:

\begin{thm}[\textsf{iWtrueI}]
\label{thm:in-wtrue-intro}
If $\Gamma \vdash c \jnat$, $\Gamma \vdash a \jnat$ and
$\Gamma \vdash \Dec(a, c) = \tofn 2$, then $\Gamma \vdash \Tru(c)$.
\end{thm}

\noindent
Elimination opens a fresh fuel variable with the affirmation as a
hypothesis, and a transport rule moves $\Tru$ along a provable equation of codes
(Appendix~\ref{app:machines}).  The adequacy of $\Tru$ at closed
codes, in both directions, is the first theorem of
Section~\ref{sec:pillars} on the reflection pillars.

\paragraph{An internal theorem about one code.}
The decider never affirms the false equation, and \rga proves this fact
as a single sentence.

\begin{thm}[\textsf{iPrDecTrueAt\_cEq01\_never}]
\label{thm:in-never}
For every context $\Gamma$,
$\Gamma \vdash \forall s.\ \neg\Aff(s, \bquo{\z = \tofn 1})$.
\end{thm}

\noindent
\emph{Proof idea.}  The obvious motive is not inductive; the proof
instead establishes two facts about every level by a case split, and
the motive follows without the induction hypothesis.  At an abstract
fuel the proof splits on whether the fuel is zero or a successor, a
split that \rga derives by induction, and settles each case directly; formally the proof is an induction
inside \rga on the fuel, with the decided statement $\Dec(v,
\bquo{\z=\tofn 1}) \ne \tofn 2$, whose step does not use the induction
hypothesis.
The base case is a closed numeral, settled by
Theorem~\ref{thm:in-dec-closed} and the example of
\S\ref{sec:inside-dec}.  The step unfolds $\Dec(\suc v, c)$, at an
abstract fuel $v$, into one application of the table's step function
to the level-$v$ table.  This unfolding exists because $\taut$ is an
explicit tree (Theorem~\ref{thm:in-cmp-rec} is the pattern).  The new
cell of $\quo{\z = \tofn 1}$ depends only on the level-$v$ cells of
$\quo{\z}$ and $\quo{\suc\z}$, the two sides of the equation.  The
obvious motive on tables, ``this table does not affirm the equation'',
is not preserved by a step.  A table whose value cell holds $\Some 0$
at every code up to that of the equation does not affirm the equation, but one step from that table does, since both
sides then read $0$ (Appendix~\ref{app:machines}).  So the proof does
not strengthen the motive.  Instead, the proof establishes two facts about the table at
every level, by a case split on the level (zero or a successor).  The
step function always writes $\Some 0$ into the cell of $\quo{\z}$, so
at any level that cell holds $\None$ or $\Some 0$, and therefore the
cell of $\quo{\suc\z}$ holds $\None$ or $\Some 1$.  With these facts
the new cell of $\quo{\z = \tofn 1}$ never holds $\Some{\ctrue}$, the
packed value 2.

Theorem~\ref{thm:in-never} is the paper's example of transparency at
work.  The opaque twin of $\Dec$ agrees with $\Dec$ at every numeral
fuel and every term-code numeral, but the development can unfold only
the explicit tree.  Section~\ref{sec:pillars} on the reflection
pillars uses Theorem~\ref{thm:in-never} as the internal companion of the
consistency route of Theorem~\ref{thm:pl-con}.

\subsection{The proof checker}
\label{sec:inside-chk}

\paragraph{Proofs and rules.}
A proof is a finite list of judgments.  Its first entry is its
conclusion, and each entry follows from later entries by one of
\rga's 42 rules.  $P \proves H \vdash c$ says that every entry of
$P$ is justified and that the conclusion is $H \vdash c$.

\paragraph{The tree checker.}
The checker consists of 41 \emph{rule detectors}, explicit index trees
that each recognize one rule's instances at codes.  One detector
serves the two directions of the reflexivity rule: since $a\jnat$
abbreviates $a = a$, both directions conclude $a = a$ from $a = a$.
A compiled validity test $\mathit{TreeValid}(p)$ requires, at every
position of the list coded by $p$, that the entry is a genuine code of
a judgment and that some detector fires on it against the entries
after it.  The tree checker adds requirements on the conclusion.

\begin{defi}[\textsf{iChk}]
\label{def:in-chk}
$\Chk(p, c)$ is the conjunction of four decided tests:
$\mathit{TreeValid}(p)$; $p \ne 0$; the conclusion of the list coded
by $p$ has no hypotheses; and that conclusion's code equals $c$.
\end{defi}

\begin{thm}[\textsf{iChk\_faithful}, \textsf{iChk\_genuine}]
\label{thm:in-chk}
\leavevmode
\begin{enumerate}
\item For every proof $P$ and term $c$:
  $\emptyset \vdash \Chk(\bquo{P}, \bquo{c})$ if and only if
  $P \proves \emptyset \vdash c$.
\item For every natural $p$ and term $c$: if
  $\emptyset \vdash \Chk(\tofn{p}, \bquo{c})$, then
  $p = \quo{P}$ for some proof $P$ with $P \proves \emptyset
  \vdash c$.
\end{enumerate}
\end{thm}

\noindent
Part~1 (\emph{faithfulness}) is an equivalence at proof codes, for
every rule, quantifier rules included.  Part~2 (\emph{genuineness})
strengthens one direction to every natural $p$.  Part~1 is
schematic; part~2 is meta, since it concludes a HOL fact about proofs
from a derivation.  An internal sentence
that quantifies over all numerals needs part~2 to be a statement about
proofs.  The consistency sentence $\ConR \ldef \forall p.\
\neg\Chk(p, \bquo{\z = \tofn 1})$ of [\citen{goedel-paper}, Section ``G\"odel's theorems for \rga''] is
such a sentence.

\emph{Proof idea.}  Both parts rest on one fact from the companion
paper: for a well-formed term $t$ and a natural $r$, \rga
derives $t = \tofn{r}$ exactly when $t$ evaluates to $r$
(the value form of open completeness, \S\ref{sec:inside-codes}).  For part~1, the compiled validity test at a
valid proof's code evaluates to 1, so \rga derives that the test equals 1;
conversely a derived acceptance forces the value 1, and the value 1
forces a valid proof.  Part~2 holds because the validity test checks
at every position that the entry is a genuine judgment code, so any
accepted number decodes to a proof list.  Without that check the test
accepts a number that codes no proof (Appendix~\ref{app:machines}).

\emph{Classical versus \rga.}  Nothing in the result differs from a
classical arithmetization.  The \rga-specific point is the route:
both directions come from \rga's HOL soundness and completeness at
closed codes, not from an internal analysis of the checker.

The checker's verdict is decided at abstract arguments, which
{}[\citen{goedel-paper}, Section ``G\"odel's theorems for \rga''] on G\"odel's theorems for \rga uses for a case split under
$\forall$-introduction:

\begin{thm}[\textsf{hq\_iChk}]
\label{thm:in-chk-dec}
If $\Gamma \vdash p \jnat$ and $\Gamma \vdash c \jnat$, then
$\Gamma \vdash \Chk(p, c) \jbool$.
\end{thm}

\paragraph{Example.}
The one-entry list $P = [(\emptyset, \z\jnat)]$ is a valid proof of
$\z \jnat$, by the rule that zero is a natural.  Its conclusion
$\z\jnat$ is the term $\z = \z$, with code
$789\,401\,477\,180\,779$, and $\quo{P}$ has 59 decimal digits.
Theorem~\ref{thm:in-chk} gives $\emptyset \vdash \Chk(\bquo{P},
\bquo{\z\jnat})$.  Conversely, any numeral that \rga proves accepted
for $\z \jnat$ is the code of a proof like this one.

\paragraph{Three checkers.}
Two more checkers appear later, and each theorem names the checker it
uses.  Table~\ref{tab:in-checkers} lists all three checkers.  The certificate
predicate of \S\ref{sec:pillars-int} wraps the transparent test $\pi$,
``$p$ codes a valid proof with conclusion judgment $h$''.  The
quantifier-free checker of Section~\ref{sec:con} is $\Chk$ together
with a decided scan for the absence of $\cA$.  Beside these checkers the
development has an opaque twin $\ChkE$ of $\Chk$, the companion
paper's PR proof-certificate test compiled through an index chosen by
Hilbert choice; $\ChkE$ serves as a control of the check of
\S\ref{sec:inside-audit}, and Appendix~\ref{app:machines-opaque}
describes it with the other opaque twins.

\begin{table}[t]
\centering
\small
\begin{tabular}{@{}>{\raggedright\arraybackslash}p{0.09\linewidth}>{\raggedright\arraybackslash}p{0.19\linewidth}>{\raggedright\arraybackslash}p{0.11\linewidth}>{\raggedright\arraybackslash}p{0.16\linewidth}>{\raggedright\arraybackslash}p{0.14\linewidth}>{\raggedright\arraybackslash}p{0.16\linewidth}@{}}
\hline
checker & built from & transparent & faithful at proof codes & genuine at every number & used in \\
\hline
$\Chk$ & the rule detectors and the validity test & yes, audited &
  yes, Thm.~\ref{thm:in-chk}(1) & yes, Thm.~\ref{thm:in-chk}(2) &
  $\ConR$, $\PrC$, Thms.~\ref{thm:pl-no-uniform}, \ref{thm:pl-complete}, Cors.~\ref{thm:pl-reflect}, \ref{thm:pl-bridge}, Sect.~\ref{sec:con}, [\citen{goedel-paper}, Section ``G\"odel's theorems for \rga''] \\
$\Chk_{\mathrm{qf}}$ & $\Chk$ and a decided scan
  for $\cA$ & yes, audited & yes, for $\forall$-free proofs
  (Thm.~\ref{thm:con-qf-faithful}) & yes (Thm.~\ref{thm:c-chkqf-genuine}) &
  Sect.~\ref{sec:con} \\
$\pi$ & detectors, validity test, head check &
  yes, audited inside $\Cert$ & yes, in HOL (Thm.~\ref{thm:a-chk-correct}) & yes,
  Thm.~\ref{thm:pl-genuine} & $\Cert$, Thm.~\ref{thm:pl-int} \\
\hline
\end{tabular}
\caption{The three proof checkers.  ``Transparent'' refers to the
  checks of \S\ref{sec:inside-audit}; ``genuine'' means that every
  accepted number codes a real proof.  The opaque twin $\ChkE$ of
  $\Chk$ is in Appendix~\ref{app:machines-opaque}.}
\label{tab:in-checkers}
\end{table}

\subsection{Closed-numeral representability}
\label{sec:inside-d3}

Classically, the third Hilbert--Bernays--L\"ob derivability condition,
an instance of formalized $\Sigma_1$-completeness, reads
\[
  \vdash \PrC(\quo{\varphi}) \limp \PrC(\quo{\PrC(\quo{\varphi})}) .
\]
The development proves a weaker condition of a different kind, which
this paper calls \emph{closed-numeral representability}: each machine,
run at closed numerals, is provably equal to the value of its HOL
counterpart.  (The Isabelle names keep the development's label D3.  In
this paper D3 is the classical condition; [\citen{goedel-paper}, Section ``G\"odel's theorems for \rga''] on
G\"odel's theorems for \rga proves a decided form of it and a family form, D3G.)  The condition
uses the numeral-code encoder $\Num \ldef
\cmp(\tau_{\mathrm{num}})\cdot\z$, compiled from an explicit tree
$\tau_{\mathrm{num}}$, which maps $\tofn{n}$ to $\bquo{\tofn{n}}$:
$\emptyset \vdash \Num\cdot\tofn{n} = \bquo{\tofn{n}}$ for every $n$,
an instance of Theorem~\ref{thm:in-cmp-num}(1).  Closed-numeral
representability collects this law with three schemata stated above:
the compiler's closed values (Theorem~\ref{thm:in-cmp-num}(1)), the
compiled decider's closed values (Theorem~\ref{thm:in-dec-closed}), and
the checker's acceptance of the code of every proof
(Theorem~\ref{thm:in-chk}(1)); and the decider is monotone in the fuel
(Theorem~\ref{thm:in-dec}(3)), so more fuel never undoes a verdict.
Unlike the classical condition, closed-numeral representability says
nothing about the provability
of statements about provability.

Closed-numeral laws hold for any compiled index
(Theorem~\ref{thm:in-cmp-num}), so the same laws hold for the opaque
twins of these machines.  The twins of the encoder and the decider
agree with the transparent encoder and decider at closed arguments (term codes, for the decider); the
checker $\Chk$ and its twin $\ChkE$ agree at proof codes, since both are faithful there
(Appendix~\ref{app:machines}).  Closed-numeral representability is
stated for the transparent machines because the internal arguments that need more than closed
values run on those machines: Theorem~\ref{thm:in-never}, the internal
completeness of \S\ref{sec:pillars-int}, and D3G.  D3G is a sentence
in decided-antecedent form over the checker's verdict, proved in
{}[\citen{goedel-paper}, Section ``G\"odel's theorems for \rga''] at the tree checker alone; its proof uses
neither the decider nor the truth predicate.  The development also
states forms of D1 and D2, which [\citen{goedel-paper}, Section ``G\"odel's proof in four forms''] on the four
forms of G\"odel's proof compares with
the classical conditions.

\subsection{Machines without opaque indices}
\label{sec:inside-audit}

The check described here protects the meaning of the paper's
sentences, not the soundness of its proofs: Isabelle checks every
proof, whatever the machines contain.  By the intensionality noted in
\S\ref{sec:inside-cmp}, what a sentence such as $\ConR$ says depends
on the machine terms it contains, and the paper claims that those
terms are the explicit trees.

G\"odel's second theorem for \rga ([\citen{goedel-paper}, Section ``G\"odel's theorems for \rga'']) needs a
primitive-recursive witness function built as an explicit tree, so
that its compiled form unfolds at abstract arguments
(\S\ref{sec:inside-cmp}).  No theorem about the value of an index can
state this requirement.  The statement ``$\den{i}$
is PR'' holds of every index $i$, chosen by Hilbert choice or not, so
the theorem of [\citen{goedel-paper}, Section ``G\"odel's theorems for \rga''] that asserts that the witness
is PR does not show that the witness contains no opaque index
(Appendix~\ref{app:machines-audit}).

The development makes the requirement mechanical.  A check that runs
with every build walks the definitional closure of every machine the
paper calls transparent, and of every explicit witness and transformer
of [\citen{goedel-paper}, Section ``G\"odel's theorems for \rga''], through every specification of every
constant it meets.  The check fails the build if a closure reaches an index
chosen by Hilbert choice, any other constant defined by choice outside
a fixed whitelist of HOL library constants (such as the conditional and
well-founded recursion), a whitelisted choice constant applied at the
index type, or a constant with no defining specification at the index
type.  Positive controls, among them the opaque twins of
\S\ref{sec:pillars-complete}, must be flagged, and a size floor rules
out a vacuous pass.  Appendix~\ref{app:machines-audit} gives the
roots, the whitelist, the controls and the closure sizes.  A second
check asks Isabelle's code generator to translate every index-valued
function in these closures into Standard ML from its defining
equations.  All these functions translate with no call that aborts (the
check skips five internal constants of Isabelle's packages, which the
generated code never calls), while a control index chosen by a case split on $\forall n.\, n =
0$, which the first check passes, fails.  The second check relies on
the correctness of the code generator, and nothing is executed.  In
this paper a machine is \emph{transparent} when these checks pass on it.  For what
a sentence such as $\ConR$ says, transparency is necessary but not
sufficient (\S\ref{sec:inside-cmp}).

\paragraph{Mechanization.}
Every theorem in this paper is a theorem of Isabelle/HOL, and the
validation of the development reports no use of oracles.
Section~\ref{sec:mech} gives the Isabelle version, the size of the
development, its trusted base, and the plan for the published
artifact.

\section{The reflection pillars}
\label{sec:pillars}

This section states what \rga proves about its own truth and its own
provability; Appendix~\ref{app:pillars} gives the proofs.  The first
pillar is soundness, proof by proof: for each proof of $c$, \rga
derives that $c$ is true by its own truth predicate
(\S\ref{sec:pillars-sound}).  The second pillar is internalized
completeness: for each well-formed grounded-true term, \rga derives
that the term is provable (\S\ref{sec:pillars-complete}), and \rga
also proves one internal sentence of completeness and the entailment
form of that sentence, from truth to a certificate (\S\ref{sec:pillars-int}).  The soundness pillar rests on the adequacy
of the truth predicate at closed codes, which comes first
(\S\ref{sec:pillars-adeq}); the completeness pillar runs on the
checker, and adequacy links the two pillars (Corollary~\ref{thm:pl-bridge}).
\S\ref{sec:pillars-con} says what the soundness pillar gives for
consistency, and proves that \rga does not prove the pillar's
single-sentence form at the false equation $\z = \tofn 1$.
\S\ref{sec:pillars-tarski} explains why a truth predicate that \rga
defines for itself does not run into Tarski's theorem, in what sense
grounded truth is transparent, and what \rga establishes about the
untruth of the Liars and of refuted sentences, and how, under the
hypothesis $\ConR$, $\Tru$ tracks \rga's own verdicts;
\S\ref{sec:breaks-liars} says where the T-biconditional holds and
compares the truth-Liar with the G\"odel sentence and the consistency
sentence; \S\ref{sec:pillars-axioms} sets the predicate against the axioms of
truth, with the full comparison in Appendix~\ref{app:pillars-axioms};
and \S\ref{sec:pillars-reflect} uses $\Tru$ as the medium of proof by
reflection, with a normalizer for polynomial identities inside \rga.

The pillars internalize \rga's own semantics.  The predicate $\Tru$
runs a compiled decider of the step-indexed evaluation that defines
grounded truth, so the pillars state, inside \rga, that \rga is sound
and complete with respect to that evaluation.  Two readings of this
result are possible, and the paper keeps both in view.

On the classical reading, grounded truth is a $\Sigma_1$ notion, and
for well-formed terms it coincides with provability
(Theorem~\ref{thm:pl-prov-iff}).  $\Tru$ is then a provability-like
predicate, and both pillars are forms of necessitation, the first
derivability condition: ``if $\vdash \varphi$ then $\vdash
\Tru(\bquo{\varphi})$'' and ``if $\varphi$ is true, hence provable,
then $\vdash \PrC(\bquo{\varphi})$'', with $\PrC$ the provability
formula of \S\ref{sec:pillars-complete}.  The converses of these two rules are
reflection rules: for $\Tru$, the rule CONEC of Friedman and
Sheard~\cite{friedman87axiomatic}, from $\vdash \Tru(\bquo{\varphi})$
to $\vdash \varphi$ (Theorem~\ref{thm:pl-nec}); for $\PrC$, the
reflection rule of [\citen{goedel-paper}, theorem~\textsf{D1\_rule\_pr}](4).  A $\Sigma_1$-sound theory has the
reflection rule for its provability predicate.  On this
reading the internal content of the section lies elsewhere: in the
internal completeness sentence (Theorem~\ref{thm:pl-int}), in
Theorem~\ref{thm:in-never}, and in the fact that a compiled decider of
\rga's own semantics exists inside \rga at all.

On \rga's own reading, grounded truth is not an approximation of a
completed totality of truths.  A sentence is true when computation
establishes it, false when computation refutes it, and otherwise
neither.  $\Tru$ expresses this notion for \rga's whole language,
quantified sentences included; \S\ref{sec:pillars-axioms} says which
axioms of truth $\Tru$ satisfies, and in which form.  Sentences that
classical reasoning counts as true but no evaluation settles, the
G\"odel sentence and the consistency sentence
(\S\ref{sec:pillars-con}), are ungrounded, and so is the Liar for truth
(\S\ref{sec:pillars-tarski}), for the same reason: none of these sentences is ever
established.  The consistency sentence does not refer to itself at
all.  A companion essay develops this grounded view of truth and
diagonalization~\cite{ford-diag}; nothing in this paper depends on that essay.
Both readings agree on the fact that matters for the paper's theme:
\rga's semantics is itself a computation, so \rga can carry its own
semantics inside itself.

The distinction of \S\ref{sec:inside-codes} between schematic,
internal and meta theorems governs every statement here;
Table~\ref{tab:pillars} in Appendix~\ref{app:pillars} lists each
result with its kind.  The pillars use the truth predicate $\Tru$ of
\S\ref{sec:inside-truth} and the compiled decider $\Dec$ beneath it,
and checker-based predicates defined below.

\subsection{The truth predicate is adequate at closed codes}
\label{sec:pillars-adeq}

\begin{thm}[\textsf{iWtrue\_of\_wosr}, \textsf{wosr\_of\_iWtrue}]
\label{thm:pl-adeq}
For every term $t$:
\begin{enumerate}
\item if $t \wbev{s} \ctrue$ for some $s$, then
  $\emptyset \vdash \Tru(\bquo{t})$;
\item if $\emptyset \vdash \Tru(\bquo{t})$, then $t \wbev{s} \ctrue$ for
  some $s$.
\end{enumerate}
\end{thm}

\noindent
Part~1 is schematic and part~2 is meta, since part~2 concludes the
existence of an evaluation.  Both parts quantify over $t$, and neither
assumes that $t$ is well-formed.  Together the two parts say that \rga derives
$\Tru(\bquo{t})$ exactly when $t$ is grounded-true.  The pair is a
statement about which closed judgments \rga derives, not an
internal biconditional $\Tru(\bquo{t}) \liff t$.  \rga proves that
biconditional at every well-formed $t$ that is grounded-true, or
grounded-false and $\forall$-free, and at no ungrounded $t$
(\S\ref{sec:breaks-liars}).

\emph{Proof idea.}  Part 1: completeness up to fuel
(Theorem~\ref{thm:in-dec}) gives a fuel $s'$ with $\decV_{s'}(t) =
\Some{\ctrue}$, so the packed verdict is 2.  Theorem~\ref{thm:in-dec-closed}
turns this verdict into $\emptyset \vdash \Dec(\tofn{s'}, \bquo{t}) = \tofn 2$,
and the introduction rule (Theorem~\ref{thm:in-wtrue-intro}) gives
$\Tru(\bquo{t})$.  Part 2: by the companion paper's soundness theorem,
the derived sentence $\Tru(\bquo{t}) \ldef \neg\forall s.\
\neg\Aff(s, \bquo{t})$ is grounded-true.  So the universal inside $\Tru(\bquo{t})$
is grounded-false, and its falsity clause supplies a numeral $i$ at
which $\Aff(\tofn{i}, \bquo{t})$ is grounded-true.  Open completeness
makes that affirmation derivable, and Theorem~\ref{thm:in-dec-closed}
then forces the packed verdict at $(i, \quo{t})$ to be 2.  Decider soundness then gives $t \wbev{i}
\ctrue$.  Neither part needs well-formedness of $t$, because the
decider's laws at closed codes have no side condition.  The absence
of a side condition is a design constraint, not luck.  A table whose
$\beta$-step substitutes through coders that read only two argument
slots affirms the ill-formed redex $(\lambda\vz)\cdot(\z = \z = \z)$,
whose equation node has three arguments.  Appendix~\ref{app:machines}
proves in HOL that such a decider cannot have both an introduction
rule like Theorem~\ref{thm:in-wtrue-intro} and a side-condition-free
part~2.

\emph{Classical versus \rga.}  In a classical recursively enumerable
theory of arithmetic, no formula $T$ satisfies ``$\vdash
T(\quo{\varphi})$ if and only if $\varphi$ is true'' for all sentences
$\varphi$, because arithmetic truth is not recursively enumerable.  The
classical analogue of Theorem~\ref{thm:pl-adeq} is the provability
predicate: a $\Sigma_1$-sound theory proves $\PrC(\quo{\varphi})$
exactly when it proves $\varphi$.  \rga's grounded truth is
recursively enumerable, because the semantics requires a certificate
at every true universal (\S\ref{sec:inside-dec}), and for well-formed
terms grounded truth coincides with provability (Theorem~\ref{thm:pl-prov-iff}).
So Theorem~\ref{thm:pl-adeq} is that analogue, stated about \rga's
semantics rather than about its proofs.  Classically, $\PrC(\quo{\varphi})$
and the truth of $\varphi$ come apart.  In \rga provability and truth do not, because
truth is defined by the same certificate-carrying computation that
the proof checker verifies.

\subsection{Soundness, proof by proof}
\label{sec:pillars-sound}

The soundness pillar covers every instance of a proved conclusion.
Instances come from a substitution: a finite HOL list $\zeta$ of pairs
of a variable index and a term.  $\sof{\zeta}{c}$ is the
simultaneous substitution of $\zeta$ into $c$.  The compiled
substitution $\Subst$ computes this substitution on codes: for a well-formed term $c$
and a list $\zeta$ of well-formed terms, $\emptyset \vdash
\Subst(\bquo{\zeta}, \bquo{c}) = \bquo{\sof{\zeta}{c}}$
(Theorem~\ref{thm:b-subst}).

\begin{thm}[\textsf{rga\_sound\_sch}]
\label{thm:pl-sound}
Let $P$ be a proof with $P \proves \emptyset \vdash c$, and let
$\zeta$ be a list of pairs whose terms are all well-formed.  Then for
every context $\Delta$,
\[
  \Delta \vdash \Tru(\Subst(\bquo{\zeta}, \bquo{c})) .
\]
\end{thm}

\noindent
The theorem is a schema over $P$, $c$, $\zeta$ and $\Delta$.  The theorem has
no well-formedness hypothesis on $c$; well-formedness follows from the
proof.  With $\zeta$ empty, the closed equation for $\Subst$ and
transport reduce the conclusion to $\emptyset \vdash \Tru(\bquo{c})$
(Theorem~\ref{thm:b-of-proof}).  Quantifiers
in $c$ need no special treatment.  $\Tru$ runs the decider, and the
decider handles each universal through its certificate.

\emph{Proof idea: pointwise collapse.}  The proof does not redo the
soundness induction over proofs inside \rga.  That induction runs
once, in HOL, as the companion paper's soundness theorem: for each entry of a valid proof
and each well-formed substitution that makes the entry's hypotheses
true, the instantiated conclusion is grounded-true.  Internalization
is then one composition per conclusion
(Theorem~\ref{thm:b-collapse}).
Theorem~\ref{thm:pl-adeq}(1) turns the grounded truth of
$\sof{\zeta}{c}$ into $\emptyset \vdash \Tru(\bquo{\sof{\zeta}{c}})$.
The closed equation for $\Subst$ and a transport rule, which moves
$\Tru$ along a provable equation of codes, rewrite the argument to
$\Subst(\bquo{\zeta}, \bquo{c})$.  Weakening moves the judgment to
$\Delta$.

\emph{Classical versus \rga.}  Classically the schema ``if $\vdash
\varphi$ then $\vdash T(\quo{\varphi})$'' alone is cheap: the
provability predicate $\PrC$ satisfies it, by the derivability condition
D1.  The content of Theorem~\ref{thm:pl-sound} is that \rga's predicate
satisfies the schema and is also adequate (Theorem~\ref{thm:pl-adeq}).
For a classical recursively enumerable theory no predicate is adequate
to truth across the whole language, since arithmetic truth is not
recursively enumerable~\cite{tarski83concept}.  What \rga does not
have is the single-sentence form of soundness, and
\S\ref{sec:pillars-con} proves that \rga does not derive even that form's
instance at $\z = \tofn 1$.

\paragraph{Example.}
Take the one-entry proof $P = [(\emptyset, \z\jnat)]$ of
\S\ref{sec:inside-chk}, and $\zeta$ empty.
Theorem~\ref{thm:pl-sound} with $\zeta$ empty, collapsed by
Theorem~\ref{thm:b-of-proof}, gives $\emptyset \vdash
\Tru(\bquo{\z\jnat})$.  Traced: soundness makes $\z = \z$
grounded-true; the decider has no verdict at fuel 1 and affirms the
equation at fuel 2, $\decV_2(\z = \z) = \Some{\ctrue}$; \rga proves
$\Dec(\tofn 2, \bquo{\z\jnat}) = \tofn 2$
(Theorem~\ref{thm:in-dec-closed}); and the introduction rule applies,
for instance at the fuel $\tofn 2$.  The mechanized proof takes
whatever fuel decider completeness supplies; any fuel from 2 on
works, by monotonicity.

\subsection{Consistency through the truth predicate, and its limit}
\label{sec:pillars-con}

\begin{thm}[\textsf{rga\_consistent\_via\_wtrue}]
\label{thm:pl-con}
No proof $P$ satisfies $P \proves \emptyset \vdash \z = \tofn 1$.
\end{thm}

\noindent
Theorem~\ref{thm:pl-con} is a meta theorem.  Its statement is \rga's consistency as a
theorem of HOL, which the companion paper's soundness theorem already
gives.  What is new is the route through the truth predicate, the HOL
shadow of the internal route through Theorem~\ref{thm:in-never}.  A
proof of $\z = \tofn 1$ would give $\emptyset \vdash \Tru(\bquo{\z =
\tofn 1})$ by Theorem~\ref{thm:pl-sound}.  A second HOL argument
refutes that judgment, again through soundness
(Appendix~\ref{app:pillars-con}).  The judgment would make the universal inside $\Tru$
grounded-false, so its falsity clause would supply a numeral fuel at
which the decider affirms $\z = \tofn 1$.  The decider's closed value
there is never 2 (the example of \S\ref{sec:inside-dec}).  So \rga
would derive a false equation between numerals, which soundness
excludes.

Theorem~\ref{thm:in-never} refutes the same judgment inside \rga, in
every context: $\Gamma \vdash \forall s.\ \neg\Aff(s, \bquo{\z = \tofn
1})$, which is exactly the universal that
$\Tru(\bquo{\z = \tofn 1})$ negates.  Theorem~\ref{thm:pl-con} does
not use this internal refutation.  The two facts would combine into an
internal proof of consistency if \rga proved soundness as one sentence.
Classically the same gap separates a schema from a sentence.  For each
specific number $p$, Peano arithmetic proves that $p$ is not the code
of a proof of $\z = \tofn 1$, since that claim is a finite check.  Peano
arithmetic does not prove that no $p$ is such a code, since that
sentence is Peano arithmetic's own consistency statement.  \rga does not prove the single sentence either:

\begin{thm}[\textsf{rga\_no\_uniform\_soundness}]
\label{thm:pl-no-uniform}
Let $U$ be the sentence $\forall p.\ \Chk(p, \bquo{\z = \tofn 1}) \limp
\Tru(\bquo{\z = \tofn 1})$.  Then $\emptyset \nvdash U$.
\end{thm}

\noindent
\emph{Proof idea.}  Suppose \rga proved $U$.  Open a fresh
variable $p$ with $p \jnat$.  The verdict $\Chk(p, \bquo{\z = \tofn
1})$ is decided (Theorem~\ref{thm:in-chk-dec}), so \rga may split on
this verdict.  If the checker accepts, $U$ gives $\Tru(\bquo{\z =
\tofn 1})$, which contradicts Theorem~\ref{thm:in-never}, so this case
yields anything, in particular $\neg\Chk(p, \bquo{\z = \tofn 1})$.  If
the checker rejects, $\neg\Chk(p, \bquo{\z = \tofn 1})$ is the case's
hypothesis.  Universal introduction then gives $\forall p.\
\neg\Chk(p, \bquo{\z = \tofn 1})$, the consistency sentence $\ConR$.
\rga does not prove $\ConR$, by G\"odel's second theorem for \rga
([\citen{goedel-paper}, theorem~\textsf{rga\_goedel2}]).  The proof uses the
assumed implication but introduces none.  The proof has the shape of the
internal consistency proof for the quantifier-free fragment
(Section~\ref{sec:con}), with $\Chk$ in place of the quantifier-free
checker.  No circularity arises: the proof of [\citen{goedel-paper}, theorem~\textsf{rga\_goedel2}]
uses neither the decider, nor the truth predicate, nor the internal
completeness of \S\ref{sec:pillars-int}.

Theorems~\ref{thm:pl-sound} and~\ref{thm:pl-no-uniform} together mark
the boundary.  \rga proves its soundness for each proof.  \rga does not
prove the single-sentence form of soundness, which fails already at
the conclusion $\z = \tofn 1$.  The entailment form of soundness fails there too:
$\PrC(\bquo{\z = \tofn 1}) \nvdash \Tru(\bquo{\z = \tofn 1})$
([\citen{goedel-paper}, \S\,``L\"ob's theorem and reflection'']), and so
does the entailment at an open code (\S\ref{sec:pillars-int}).  (The Isabelle name says ``uniform'' in
the sense of ``over all proof codes''; the reflection literature uses
``uniform'' for a different principle.)  The obstruction is G\"odel's
second theorem, as in the classical case.  There the local reflection
instance $\PrC(\quo{\z = \tofn 1}) \limp \z = \tofn 1$ is equivalent to
the consistency statement, which G\"odel's second theorem rules out,
and L\"ob's theorem~\cite{lob55solution} extends the obstruction to
every sentence the theory does not prove.  The \rga case and the classical case are not literally the same:
Theorem~\ref{thm:pl-no-uniform}'s sentence has $\Tru(\bquo{\z = \tofn
1})$ in its consequent, not $\z = \tofn 1$, and
Theorem~\ref{thm:in-never} is what links the two consequents.

\paragraph{Consistency is ungrounded.}
In \rga's own semantics, $\ConR$ and the sentence $U$ of
Theorem~\ref{thm:pl-no-uniform}, single-sentence soundness at $\z =
\tofn 1$, have the same status: \rga proves every numeral instance,
and the universal has no value.  For $\ConR$, \rga proves the instance
$\neg\Chk(\tofn{p}, \bquo{\z = \tofn 1})$ at every numeral $\tofn{p}$,
and $\ConR$ has no value at any fuel, because a well-formed universal
that \rga neither proves nor refutes has no value, and $\ConR$ is
unprovable ([\citen{goedel-paper}, theorem~\textsf{rga\_goedel2}]) and unrefutable
([\citen{goedel-paper}, theorem~\textsf{ConRGA\_unrefutable}]).  \rga does not refute $U$ either, $U$ has
no value, and each numeral instance of $U$ follows from the instance
of $\ConR$ at the same numeral, as the first disjunct of the
implication.  Like the instances of $\ConR$, these proofs use \rga's
soundness (Appendix~\ref{app:pillars-con}).  The $\omega$-semantics of
\S\ref{sec:inside-omega} makes $\ConR$ true, while \rga does not
derive $\Tru(\bquo{\ConR})$: a universal all of whose instances are
grounded-true is not thereby grounded-true, because its truth needs a
certificate.  So \rga is $\omega$-incomplete at $\ConR$
(Theorem~\ref{thm:br-omega}), and \S\ref{sec:breaks-liars}, on the
T-schema and the Liars, states that theorem together with the same
$\omega$-incompleteness at the Liar of \S\ref{sec:pillars-tarski}.

\subsection{Internalized completeness}
\label{sec:pillars-complete}

The completeness pillar runs on the checker, not on the decider.  The pillar's
natural provability predicate is the existential over the transparent
tree checker of \S\ref{sec:inside-chk}:
\[
  \PrC(c) \;\ldef\; \neg\forall q.\ \neg\Chk(q, c) .
\]

\begin{thm}[\textsf{rga\_complete\_sch\_chk}]
\label{thm:pl-complete}
For every well-formed term $t$: if $t \wbev{s} \ctrue$, then
$\emptyset \vdash \PrC(\bquo{t})$.
\end{thm}

\noindent
The proof is three steps: open completeness gives a
proof $P$ of $t$, faithfulness (Theorem~\ref{thm:in-chk}) makes
$\Chk(\bquo{P}, \bquo{t})$ derivable, and the rule that introduces a
negated universal from one refuted instance concludes $\PrC(\bquo{t})$
at the numeral $\bquo{P}$.

\begin{cor}[\textsf{rga\_truth\_reflect\_sch\_chk}]
\label{thm:pl-reflect}
For every well-formed term $t$: if $t \wbev{s} \ctrue$, then
$\emptyset \vdash \Tru(\bquo{t})$ and $\emptyset \vdash
\PrC(\bquo{t})$.
\end{cor}

\begin{cor}[\textsf{iPrC\_of\_iWtrue}]
\label{thm:pl-bridge}
For every well-formed term $t$: if $\emptyset \vdash \Tru(\bquo{t})$,
then $\emptyset \vdash \PrC(\bquo{t})$.
\end{cor}

\noindent
All three results are schemata with the hypothesis that $t$ is well-formed,
which open completeness needs.  Corollary~\ref{thm:pl-bridge} is a
rule on theorems, not an internal implication.  The corollary's proof composes
Theorem~\ref{thm:pl-adeq}(2) with Theorem~\ref{thm:pl-complete}, and
Corollary~\ref{thm:pl-reflect} pairs Theorem~\ref{thm:pl-complete}
with Theorem~\ref{thm:pl-adeq}(1).

\emph{Proof idea.}  Open completeness turns the grounded truth of a
well-formed $t$ into a derivation, and so into a HOL proof object $P$
of $t$.  The checker accepts $P$'s code, and negated universal
introduction at that numeral gives $\PrC(\bquo{t})$.  The witness is
whatever proof open completeness supplies, not one the paper chooses.
The decider plays no part.
Soundness enters only through the companion paper's
soundness--completeness pair, of which the proof uses the completeness
direction.  So the two judgments of Corollary~\ref{thm:pl-reflect} come
from one semantic witness through two machines: the decider for
$\Tru$, the checker for $\PrC$.  Neither derivation uses the other
machine.  The same results hold for the opaque twin $\Prv$ of $\PrC$,
the existential over the opaque checker $\ChkE$ of
Appendix~\ref{app:machines-opaque}; $\Prv$ serves as a control of the
check of \S\ref{sec:inside-audit} (Definition~\ref{def:pl-prv},
Appendix~\ref{app:pillars-complete}).

\emph{Classical versus \rga.}  For a $\Sigma_1$-sound classical
theory, ``every true sentence is provably provable'' would make every
true sentence provable, contradicting G\"odel's first theorem.  In
\rga, ``true'' means grounded-true, and grounded truth coincides with
provability:

\begin{thm}[\textsf{rga\_prov\_iff\_true}]
\label{thm:pl-prov-iff}
For every well-formed term $t$: some proof $P$ satisfies $P \proves
\emptyset \vdash t$ if and only if $t \wbev{s} \ctrue$ for some $s$.
\end{thm}

\noindent
Theorem~\ref{thm:pl-prov-iff} is a meta theorem: both sides are HOL
propositions.  Provability in \rga is certified groundedness.  The Liar
for the truth predicate (\S\ref{sec:pillars-tarski}) and \rga's
G\"odel sentence ([\citen{goedel-paper}, Section ``G\"odel's theorems for \rga'']) are well-formed, unprovable
and unrefutable.  So by Theorem~\ref{thm:pl-prov-iff}, applied to each
sentence and to its negation, neither sentence is grounded-true or
grounded-false.  The Liar is a negation and the G\"odel sentence a
universal, so neither can take another value: both are ungrounded
(Theorem~\ref{thm:pl-no-value} and [\citen{goedel-paper}, \S\,``G\"odel's first theorem'']).
The completeness pillar says nothing about these two sentences.  The pillar is
completeness for grounded truth, not negation completeness.

The pillar has a converse: if $\emptyset \vdash \PrC(\bquo{t})$, then
$\emptyset \vdash t$ ([\citen{goedel-paper}, theorem~\textsf{D1\_rule\_pr}](4)).  The witness
property gives a numeral that $\Chk$ accepts for $\bquo{t}$, and
genuineness (Theorem~\ref{thm:in-chk}(2)) turns that numeral into a proof.
The G\"odel paper compares \rga's derivability conditions with the
classical ones ([\citen{goedel-paper}, \S\,``The derivability conditions \rga has'']).

\paragraph{What the truth predicate satisfies.}
$\Tru$ satisfies the two rules that Friedman and Sheard call NEC and
CONEC~\cite{friedman87axiomatic}.

\begin{thm}[\textsf{iWtrue\_nec}, \textsf{iWtrue\_conec}]
\label{thm:pl-nec}
NEC: if $\emptyset \vdash t$, then $\emptyset \vdash
\Tru(\bquo{t})$, for every term $t$.  CONEC: if $\emptyset \vdash
\Tru(\bquo{t})$ and $t$ is well-formed, then $\emptyset \vdash t$.  So
for well-formed $t$, $\emptyset \vdash t$ holds exactly when
$\emptyset \vdash \Tru(\bquo{t})$.
\end{thm}

\noindent
NEC is Theorem~\ref{thm:pl-sound} with $\zeta$ empty, and CONEC is
Theorem~\ref{thm:pl-adeq}(2) followed by
Theorem~\ref{thm:pl-prov-iff}.  CONEC needs its well-formedness
hypothesis: \rga derives $\Tru(\bquo{t})$ for a grounded-true term $t$
that is not well-formed, and hence not derivable
(Theorem~\ref{thm:pl-not-sentence}).  \S\ref{sec:pillars-axioms} sets
$\Tru$ against the other axioms and rules of axiomatic truth theories,
and \S\ref{sec:breaks-liars} states where the internal T-biconditional
holds and where it fails.

\subsection{Internal completeness as one sentence}
\label{sec:pillars-int}

The schematic pillar of \S\ref{sec:pillars-complete} gives one
derivation per grounded-true term.  \rga also proves a single internal
sentence of completeness.  The sentence's consequent uses the transparent test
$\pi$, through a certificate predicate.

\begin{defi}[\textsf{iHasCert}]
\label{def:pl-cert}
Let $\pi$ be the explicit index of the test ``$p$ codes a valid
proof whose conclusion judgment has code $h$''.  The \emph{certificate
predicate} is $\Cert(h) \ldef \exists p.\ \cmp(\pi)\cdot h\cdot p \ne
\z$.  For a code term $c$, write $\Cert_\emptyset(c)$ for $\Cert$
applied to the compiled pair of $\bquo{\emptyset}$ and $c$, the code
of the judgment $\emptyset \vdash c$.
\end{defi}

\noindent
Write $\Wfc(c)$ and $\IsCode(c)$ for the compiled, 0/1-valued tests
``$c$ codes a well-formed term'' and ``$c$ is the code of a term''.

\begin{thm}[\textsf{rga\_complete\_int2}]
\label{thm:pl-int}
\[
  \emptyset \vdash \forall s.\ \forall c.\
    \bigl(\Wfc(c) \wedge \IsCode(c) \wedge \Aff(s, c)\bigr)
    \limp \Cert_\emptyset(c) .
\]
\end{thm}

\noindent
The theorem is internal: one sentence, quantified by \rga's own $\cA$.
The conjunction
in the antecedent is the 0/1-valued conjunction of three decided
tests, so the antecedent is decided, and \rga can introduce the
implication by discharging the antecedent.  The sentence says: whenever the decider affirms a
well-formed term code at some fuel, the test $\pi$ accepts some number
as a proof of that term.  The guards $\Wfc$ and $\IsCode$ sit in the
antecedent because $c$ ranges over all numerals, including numbers
that code nothing.

\emph{Proof idea of Theorem~\ref{thm:pl-int}.}  One induction inside \rga, on the decider's level,
proves a motive with a value half and a head half, mirroring the
decider's two parts.  The value half says: if $\Dec$ at the level
gives the code $c$ of a well-formed term $t$ the verdict $v$, then
$\Cert_\emptyset$ holds of the code of the equation $t = \tofn{v}$,
built on codes with $\Num$.  The head
half carries the analogous statement for head verdicts, with its own
certificate witness.  Each arm of the induction step mirrors one
clause of the decider, and builds the matching proof step with the
checker's laws at abstract arguments.  At a universal, the decider's
affirmation is itself a found certificate of the open instance.  The
checker accepts that certificate, and one step of $\forall$-introduction
on top of it gives a certificate for the universal.  A final step
turns a certificate of $t = \tofn 1$ into a certificate of $t$.  Every
antecedent in the motive is decided, so the proof discharges every
hypothesis legitimately.  Appendix~\ref{app:pillars} follows the
induction arm by arm.

``Everything true is provable'' has the forms of
\S\ref{sec:inside-codes}, and Theorem~\ref{thm:pl-int} is one of
those forms.  The natural grounded reading is the entailment: from the
hypothesis that a numeral $c$ coding a well-formed term is true, \rga
derives that $c$ has a certificate.

\begin{thm}[\textsf{complete\_\allowbreak ent}, \textsf{complete\_\allowbreak imp\_\allowbreak fails}]
\label{thm:pl-int-ent}
\leavevmode
\begin{enumerate}
\item $\Tru(\vz),\ \Wfc(\vz),\ \IsCode(\vz),\ \vz\jnat \vdash
  \Cert_\emptyset(\vz)$.  As a rule in every context: if $\Gamma$
  derives $c\jnat$, $\Wfc(c)$, $\IsCode(c)$ and $\Tru(c)$, then
  $\Gamma \vdash \Cert_\emptyset(c)$.
\item $\emptyset \nvdash \forall c.\ \Tru(c) \limp
  \Cert_\emptyset(c)$, and $\emptyset \nvdash \forall c.\
  (\Wfc(c) \wedge \IsCode(c)) \limp (\Tru(c) \limp
  \Cert_\emptyset(c))$.
\item For every well-formed $t$ and every context $\Gamma$,
  $\Tru(\bquo{t}), \Gamma \vdash \PrC(\bquo{t})$.
\end{enumerate}
\end{thm}

\noindent
The development states the guards of parts~1 and~2 as $\Wfc(c) \ne
\z$ and $\IsCode(c) \ne \z$, interderivable with $\Wfc(c)$ and
$\IsCode(c)$ by the convention of \S\ref{sec:inside-truth}.

\noindent
\emph{Proof idea.}  (1) $\Tru(c)$ is $\exists s.\ \Aff(s, c)$, and
$\Aff$ is decided.  Eliminate the existential at a fresh fuel $s$,
with $\Aff(s, c)$ as a hypothesis, and apply the rule form of
Theorem~\ref{thm:pl-int} at $s$ and $c$
(Theorem~\ref{thm:b-cert-rule}).  (2) Instantiate at the code of
the truth-predicate Liar $L_{\Tru}$ of \S\ref{sec:pillars-tarski}.
\rga proves both guards there, and by the disjunction property the
instance needs \rga to refute $\Tru(\bquo{L_{\Tru}})$ or to prove
$\Cert_\emptyset(\bquo{L_{\Tru}})$; \rga does neither.  (3) At the code
of a term, a certificate yields an accepted proof in every context;
compose with part~1.  \qed

\medskip
\noindent
The forms of ``everything true is provable'' then line up.  The
implication, over all codes, fails at the Liar's code, where $\Tru$
is undecided (part~2).  The entailment holds (part~1), with an
existential consequent and no explicit transformer.
Theorem~\ref{thm:pl-int} is that entailment in its decided-antecedent
form: by Theorem~\ref{thm:in-dt}(2), an entailment from $\Tru(c) = \exists
s.\ \Aff(s, c)$ is equivalent to the sentence $\forall s.\ \Aff(s,
c) \limp \Cert_\emptyset(c)$, and Theorem~\ref{thm:pl-int} states
that sentence for all codes at once, with the decided guards moved
into the antecedent.  So the decided affirmation $\Aff(s, c)$ can be
an antecedent, and the undecided $\Tru(c)$ only a hypothesis.  The
converse, soundness, fails even as an entailment: $\PrC(\vz),\
\vz\jnat \nvdash \Tru(\vz)$, since at $\vz = \bquo{\z = \tofn 1}$, by
Theorem~\ref{thm:in-dt}(2), that entailment would give the sentence of
Theorem~\ref{thm:pl-no-uniform} (Appendix~\ref{app:pillars-int}).

The consequent of parts~1 and~2 is the certificate test $\pi$, not the
checker $\Chk$; part~3 is the completeness entailment with $\PrC$, a
schema over $t$.  Whether \rga proves the passage from a certificate
to an accepted proof as one sentence over all codes, in the
decided-antecedent form of \S\ref{sec:inside-codes}, is open.

In HOL, the transparent test is genuine at every natural:

\begin{thm}[\textsf{provesChk\_genuine}]
\label{thm:pl-genuine}
For every natural $p$, hypothesis set $H$ and term $c$: if
$\den{\pi}(\cpair{\quo{H}}{\quo{c}}, p) \ne 0$, then $p = \quo{P}$ for
some proof $P$ with $P \proves H \vdash c$.
\end{thm}

\noindent
So every number that the certificate predicate can be talking about
is the code of a real proof.

\emph{Classical versus \rga.}  Theorem~\ref{thm:pl-int} plays the role
of formalized $\Sigma_1$-completeness, ``every true $\Sigma_1$
sentence is provably provable''.  Theorem~\ref{thm:pl-int} differs in two ways.  The
antecedent is the decider's verdict, which covers quantified sentences
because their certificates are part of the verdict.  And the
implication is \rga's implication obtained by discharging a decided
antecedent, the only antecedent \rga can discharge.

G\"odel's second theorem for \rga does not need
Theorem~\ref{thm:pl-int}.  The proof of that theorem ([\citen{goedel-paper}, theorem~\textsf{rga\_goedel2}]) and the proof
of the internal second theorem ([\citen{goedel-paper}, theorem~\textsf{rga\_goedel2\_internal\_contra}]) use
neither the decider nor the truth predicate nor internal completeness.
Both proofs build the family form D3G of the third derivability condition at
the tree checker.  A second proof of the same statement
([\citen{goedel-paper}, Appendix ``The search route'']) builds its D3G witness by a search,
and uses there the rule form of Theorem~\ref{thm:pl-int}
(Theorem~\ref{thm:b-cert-rule}), not the sentence itself.

\paragraph{Example.}
Instantiate Theorem~\ref{thm:pl-int} at the fuel $\tofn 2$ and the
code $\bquo{\z\jnat}$.  The antecedent is a closed, decided term.
$\Wfc$ and $\IsCode$ hold of $\bquo{\z\jnat}$, and
$\Aff(\tofn 2, \bquo{\z\jnat})$ holds by the example of
\S\ref{sec:pillars-sound}.  So \rga derives
$\Cert_\emptyset(\bquo{\z\jnat})$: not every number fails to be an
accepted proof of $\z\jnat$.  The development proves this instance
directly from the rule form (Theorem~\ref{thm:b-cert-rule}).  In HOL,
Theorem~\ref{thm:pl-genuine} says that every accepted number is the
code of a proof of $\z \jnat$, such as the one-entry proof $P$.

\subsection{The truth predicate and the Liar}
\label{sec:pillars-tarski}

This subsection and the next use results on G\"odel's theorems for
\rga from the G\"odel paper~\cite{goedel-paper}; the box at the end of
\S\ref{sec:inside-codes} states the ones this paper relies on, and the
text cites the others where it uses them.
Tarski's theorem for a theory needs the diagonal biconditional $L
\liff \neg T(\quo{L})$ and the T-biconditionals $T(\quo{\varphi}) \liff
\varphi$ as provable sentences.  From these biconditionals a contradiction follows by
reasoning that is valid even intuitionistically.  \rga has the
diagonal only at the level of codes.  The compiled diagonal function
$\Diag$ satisfies $\emptyset \vdash \Diag(\bquo{t}, \tofn{k}) =
\bquo{\ssub{0}{\tofn{k}}{t}}$ for every well-formed $t$
(Appendix~\ref{app:pillars-liar}).  Define
\[
  D_{\Tru} \ldef \neg\Tru(\Diag(\vz, \vz)) , \qquad
  S_{\Tru} \ldef \Diag(\bquo{D_{\Tru}}, \bquo{D_{\Tru}}) , \qquad
  L_{\Tru} \ldef \neg\Tru(S_{\Tru}) .
\]

\begin{thm}[\textsf{goedelW\_diag}, \textsf{goedelWS\_eq}]
\label{thm:pl-diag}
$\ssub{0}{\bquo{D_{\Tru}}}{D_{\Tru}} = L_{\Tru}$ (an identity of terms), and
$\emptyset \vdash S_{\Tru} = \bquo{L_{\Tru}}$.
\end{thm}

\begin{thm}[\textsf{goedelW\_no\_value}]
\label{thm:pl-no-value}
For all $s$ and $r$, it is not the case that $L_{\Tru} \wbev{s} r$.
\end{thm}

\begin{thm}[\textsf{goedelW\_unprovable}, \textsf{goedelW\_unrefutable}]
\label{thm:pl-liar-indep}
$\emptyset \nvdash L_{\Tru}$, $\;\emptyset \nvdash \neg L_{\Tru}$, and
$\emptyset \nvdash L_{\Tru} \jbool$.
\end{thm}

\noindent
\emph{Proof of Theorem~\ref{thm:pl-no-value}.}  Suppose $L_{\Tru} =
\neg\Tru(S_{\Tru})$ had a value.  Then $\Tru(S_{\Tru})$ has a value, true or false.
If $\Tru(S_{\Tru})$ is grounded-true, open completeness makes $\Tru(S_{\Tru})$
derivable, and transport along $S_{\Tru} = \bquo{L_{\Tru}}$ makes
$\Tru(\bquo{L_{\Tru}})$ derivable.  By Theorem~\ref{thm:pl-adeq}(2),
$L_{\Tru}$ is then grounded-true, while $L_{\Tru} = \neg\Tru(S_{\Tru})$ is
grounded-false; values are unique, a contradiction.  If $\Tru(S_{\Tru})$ is
grounded-false, then $L_{\Tru}$ is grounded-true and so derivable.
Theorem~\ref{thm:pl-adeq}(1) makes $\Tru(\bquo{L_{\Tru}})$ derivable,
and transport makes $\Tru(S_{\Tru})$ derivable.  \rga would then derive both
$\Tru(S_{\Tru})$ and its negation $L_{\Tru}$, against consistency.
Theorem~\ref{thm:pl-liar-indep} follows: a well-formed term with no
value is neither provable nor refutable, by soundness, and so by
\rga's disjunction property such a term is not provably decided either.

Tarski's semantic argument needs, in addition, bivalence at
$L_{\Tru}$: the argument's two cases are exhaustive classically, and the argument
ends in a contradiction.  In \rga's grounded semantics the two cases are not
exhaustive, and Theorem~\ref{thm:pl-no-value} shows that $L_{\Tru}$
falls in the gap.  The argument ends in a diagnosis: $L_{\Tru}$ is
ungrounded.  On \rga's grounded reading this diagnosis is the intended outcome,
not an escape hatch.  A sentence that only speaks of its own truth
gives evaluation nothing to compute, so it receives no value, much as
a program that only calls itself returns none.  This view is close in
spirit to Kripke's fixed-point construction~\cite{kripke75outline},
with one difference: \rga's grounded truth is defined by a computation
that \rga carries inside itself, not by a fixed point over a classical
model.  The syntactic version of Tarski's argument needs the internal
T-biconditional at the Liar.  \rga proves that biconditional at every
well-formed sentence that is grounded-true, or grounded-false and
$\forall$-free, but at no ungrounded sentence, $L_{\Tru}$ included
(\S\ref{sec:breaks-liars}).  Taken as a hypothesis, the material
T-biconditional at the Liar is explosive, which is Tarski's syntactic
argument run in \rga's logic.

\begin{thm}[\textsf{tbicond\_liar\_explosive}, \textsf{tbicond\_liar\_unprovable}]
\label{thm:pl-tb-explosive}
For every context $\Gamma$, $\bigl(\Tru(\bquo{L_{\Tru}}) \liff
L_{\Tru}\bigr), \Gamma \vdash \z = \tofn 1$; so the hypothesis entails
both $L_{\Tru}$ and $\neg L_{\Tru}$, and
$\emptyset \nvdash \Tru(\bquo{L_{\Tru}}) \liff L_{\Tru}$.
\end{thm}

\noindent
\emph{Proof idea.}  The biconditional is the material conjunction of
$\neg\Tru(\bquo{L_{\Tru}}) \lor L_{\Tru}$ and $\neg L_{\Tru} \lor
\Tru(\bquo{L_{\Tru}})$.  Eliminate the first disjunction: in its left
case, transport along $S_{\Tru} = \bquo{L_{\Tru}}$ turns
$\neg\Tru(\bquo{L_{\Tru}})$ into $\neg\Tru(S_{\Tru})$, which is
$L_{\Tru}$; in its right case $L_{\Tru}$ is the case hypothesis.  So
the hypothesis entails $L_{\Tru}$, and with it
$\neg\Tru(\bquo{L_{\Tru}})$.  In both cases of the second disjunction
negation elimination then gives $\z = \tofn 1$.  The proof uses only
the elimination of a derived disjunction, which needs no decidedness,
conjunction elimination and transport.  The non-derivability follows
by a cut, since \rga does not derive $\z = \tofn 1$.  \qed

\medskip
\noindent
The argument uses nothing of the truth predicate: only the elimination
of disjunctions, equality substitution and the code equation
$S_{\Tru} = \bquo{L_{\Tru}}$ of the predicate's own diagonal Liar.  So
the argument is generic in the predicate.

\begin{thm}[\textsf{tbicond\_liar\_explosive\_gen}, \textsf{tbicond\_liar\_explosive\_diag}]
\label{thm:pl-tb-generic}
\leavevmode
\begin{enumerate}
\item Let $\varphi$ be a well-formed term, read as the predicate
  $T(u) \ldef \ssub{v}{u}{\varphi}$ of its variable $v$, and let $S$
  be a term and $\Gamma$ a context with $\Gamma \vdash S = \bquo{L}$,
  where $L \ldef \neg T(S)$.  Then $\bigl(T(\bquo{L}) \liff L\bigr),
  \Gamma \vdash \z = \tofn 1$.
\item Every closed well-formed term $T$, read as the predicate
  $u \mapsto T\cdot u$, has a diagonal Liar: for $D \ldef \neg(T\cdot
  \Diag(\vz, \vz))$ and $S \ldef \Diag(\bquo{D}, \bquo{D})$,
  $\emptyset \vdash S = \bquo{L}$ with $L \ldef \neg(T\cdot S)$.  So
  $\bigl(T\cdot\bquo{L} \liff L\bigr), \Gamma \vdash \z = \tofn 1$
  for every context $\Gamma$, and $\emptyset \nvdash T\cdot\bquo{L}
  \liff L$.
\end{enumerate}
\end{thm}

\noindent
Theorem~\ref{thm:pl-tb-explosive} is part~1 at $\varphi =
\Tru(\vz)$ and $S = S_{\Tru}$, and the development proves that theorem that way.
Part~2 builds the Liar of $T$ as $L_{\Tru}$ is built.  An extension of
\rga by hypotheses is a context $\Gamma$ of axioms.  By part~1, an
extension that derives the code equation of a predicate's Liar and
contains the material T-biconditional at that Liar derives $\z = \tofn
1$; by part~2, the code equation is a theorem of \rga for every closed
predicate term.  So no consistent extension of \rga by hypotheses adds
the material T-biconditional at the Liar of the truth predicate $\Tru$, or of
any other predicate whose Liar's code equation the extension derives, as an axiom.
What remains open
is narrower (Section~\ref{sec:concl}): whether a non-material
conditional in Field's style~\cite{field08saving} could carry the
biconditional consistently, and whether a primitive truth predicate
with PKF's rules T-Intro and T-Elim~\cite{halbach06axiomatizing},
interpreted for instance by a Kripke construction over the
$\omega$-semantics, could sit beside $\Tru$.  The predicate $\Tru$
exists for a reason that the classical setting also shows.  Classical
arithmetic truth is not recursively enumerable, so no classical
recursively enumerable theory can derive $T(\quo{\varphi})$ for
exactly the true $\varphi$.  Grounded truth is recursively enumerable,
because of the certificate requirement, and $\Tru$ is adequate to grounded truth.

\paragraph{Two kinds of transparency.}
In Field's theory and in PKF~\cite{field08saving,halbach06axiomatizing},
truth is transparent in a classical sense: $T\quo{\varphi}$ and
$\varphi$ are intersubstitutable in every context, under hypotheses
and under negation, so $\varphi \eqv T\quo{\varphi}$ and $\neg\varphi
\eqv \neg T\quo{\varphi}$ for every sentence, and strong Kleene logic
keeps this intersubstitutivity consistent at the Liar.  ``$\varphi$ is true'' then says no
more than $\varphi$.  Grounded truth is transparent in a different
sense, the one the grounded reading calls for: establishing $\varphi$
and establishing that $\varphi$ is true are the same achievement.  For
well-formed $\varphi$, $\vdash \varphi$ holds exactly when $\vdash
\Tru(\bquo{\varphi})$ does; this \emph{establishment transparency} is
the pair of rules NEC and CONEC (Theorem~\ref{thm:pl-nec}).  The pair
is not specific to grounded truth.  The Friedman--Sheard theory FS has
both rules~\cite{friedman87axiomatic,halbach11axiomatic}, and so does
the provability predicate of any $\Sigma_1$-sound theory, by D1 and
$\Sigma_1$-soundness.  What is specific to $\Tru$ is what it lets \rga
say.  $\Tru$ is adequate to grounded truth at closed codes
(Theorem~\ref{thm:pl-adeq}), and $\neg\Tru$ says ``not established'',
which differs from ``refuted''.  So \rga proves that the looping term
$\Omega \ldef (\lambda(\vz\cdot\vz))\cdot(\lambda(\vz\cdot\vz))$, which
head-steps to itself, is not established, and that the ordinary Liar
is neither true nor false (Theorems~\ref{thm:pl-omega}
and~\ref{thm:pl-liar-neither}), statements that a classically
transparent predicate cannot make without a separate determinacy
operator.

What establishment transparency does not give, and should not, are
the two classical intersubstitutions.  First, $\varphi \vdash
\Tru(\bquo{\varphi})$ would make a hypothesis count as an
establishment, and assuming is not establishing: the entailment fails
at the truth-Liar $L_{\Tru}$, at $G$ and at $\ConR$, and
$\Tru(\bquo{L_{\Tru}}) \nvdash L_{\Tru}$ (Theorem~\ref{thm:br-tin-ent}).
Second, $\neg\varphi \eqv \neg\Tru(\bquo{\varphi})$ would identify
``refuted'' with ``not established'', which \rga keeps apart: \rga proves
$\neg\Tru(\bquo{\Omega})$ while $\Omega$ itself stays ungrounded.  In a
classically transparent theory ``the Liar is not true'' just is the
Liar, which is why Field adds a separate determinacy operator, with
its own hierarchy and revenge problems~\cite{field08saving}.

Classical transparency is not available to $\Tru$ even if one wanted
it.  $\Tru$ is not a primitive predicate but a term defined from
\rga's decider, so $\Tru$ has a classical interpretation in the
$\omega$-semantics (\S\ref{sec:inside-omega}), where its extension is
the set of codes of grounded-true sentences and \rga's derivations,
under any $\omega$-true hypotheses, are sound.  There $L_{\Tru}$ is
true and $\Tru(\bquo{L_{\Tru}})$ is false, so no entailment can carry
one to the other.  This failure is Tarski's theorem in the $\omega$-semantics:
no defined predicate whose entailments are sound for a classical model
is classically transparent.  Nor can transparency be added by fiat, by
Theorem~\ref{thm:pl-tb-explosive}.

\begin{table}[tbp]
\centering
\small
\begin{tabular}{@{}>{\raggedright\arraybackslash}p{0.22\linewidth}>{\raggedright\arraybackslash}p{0.27\linewidth}>{\raggedright\arraybackslash}p{0.45\linewidth}@{}}
\hline
& classical transparency (Field, PKF) & establishment transparency (grounded truth) \\
\hline
idea & ``$\varphi$ is true'' says no more than $\varphi$, in every context & establishing $\varphi$ is establishing ``$\varphi$ is true'' \\
$\vdash \varphi$ iff $\vdash \Tru(\bquo{\varphi})$ & holds & holds (NEC and CONEC, Thm.~\ref{thm:pl-nec}), as for FS and for the provability predicate of a $\Sigma_1$-sound theory \\
$\varphi \eqv \Tru(\bquo{\varphi})$ under hypotheses, at well-formed $\varphi$ & holds & not expected (assuming is not establishing); fails at $L_{\Tru}$, $G$ and $\ConR$.  Under $\ConR$, $\Tru(\bquo{\varphi}) \liff \varphi$ is derivable at every well-formed $\varphi$ that \rga proves or refutes, and refuted at $G$ \\
refuted vs unestablished: $\neg\varphi \eqv \neg\Tru(\bquo{\varphi})$ & holds & not expected (refuted is not unestablished).  A refuted $\forall$-free sentence is provably unestablished; under the hypothesis $\ConR$ every refuted sentence is; at a refuted universal, $\neg\Tru$ is interderivable with $\ConR$.  \rga proves $\neg\Tru(\bquo{\Omega})$ but not $\neg\Omega$ \\
can state ``$\varphi$ is not true'' of an ungrounded sentence & no, at the Liar (needs a determinacy operator) & yes: proved at $\Omega$ and at the ordinary Liar $L$, and of $\neg L$; under $\ConR$, at $G$ and at $\ConR$; not at $L_{\Tru}$, whose untruth is $L_{\Tru}$ itself \\
\hline
\end{tabular}
\caption{Two kinds of transparency (\S\ref{sec:pillars-tarski}).  The
  second row does not discriminate; the last two rows do.  The fourth
  row is Theorem~\ref{thm:pl-refuted}; the last row is
  Theorems~\ref{thm:pl-liar-neither} and~\ref{thm:pl-omega}; the
  entries under $\ConR$ are Theorem~\ref{thm:pl-tracking}.}
\label{tab:pl-transparency}
\end{table}

Table~\ref{tab:pl-transparency} sums up the comparison.  Two theorems
illustrate the table's last two rows: one about the ordinary Liar and the
truth-Liar, one about refuted sentences.

\begin{thm}[\textsf{ordinary\_liar\_neither}, \textsf{truth\_liar\_not\_untrue\_unprovable}]
\label{thm:pl-liar-neither}
\leavevmode
\begin{enumerate}
\item For the ordinary Liar $L$ of \S\ref{sec:bg-rga},
  $\emptyset \vdash \neg\Tru(\bquo{L})$ and $\emptyset \vdash
  \neg\Tru(\bquo{\neg L})$, while $\emptyset \nvdash L$ and $\emptyset
  \nvdash \neg L$.
\item $\emptyset \nvdash \neg\Tru(\bquo{L_{\Tru}})$.
\end{enumerate}
\end{thm}

\noindent
Part~1 says that \rga establishes, internally and with no determinacy
operator, that the ordinary Liar is neither true nor false.  The
internal halves are proved by one induction inside \rga on the
decider's level (Appendix~\ref{app:pillars-liar}): a decided value of $L$ at one
level becomes, by two $\beta$-laws and the negation law, the opposite
decided value three levels down, so a window of consecutive levels at
which $L$ has no decided value carries the induction.  These internal
halves use no soundness theorem; the two non-derivability facts do.
Part~2 is the contrast.
The untruth of the truth-Liar, $\neg\Tru(\bquo{L_{\Tru}})$, is
$L_{\Tru}$ itself up to the code equation $S_{\Tru} =
\bquo{L_{\Tru}}$, and this untruth behaves like a G\"odel sentence: true in the
$\omega$-semantics, and not established.

\begin{thm}[\textsf{refuted\_unestablished\_trichotomy}, \textsf{refuted\_universal\_unestablished\_eqv\_con}]
\label{thm:pl-refuted}
Let $\Gamma$ be any context.
\begin{enumerate}
\item If $t$ is $\forall$-free and $\emptyset \vdash \neg t$, then
  $\emptyset \vdash \neg\Tru(\bquo{t})$.
\item If $\emptyset \vdash \neg t$, then $\ConR, \Gamma \vdash
  \neg\Tru(\bquo{t})$ and $\ConR, \Gamma \vdash \neg\PrC(\bquo{t})$.
\item For every well-formed $u$, $\neg\Tru(\bquo{\cA\cdot u}), \Gamma
  \vdash \ConR$; so $\emptyset \nvdash \neg\Tru(\bquo{\cA \cdot u})$.
\item If $\emptyset \vdash \neg(\cA\cdot u)$, then $\ConR \eqv
  \neg\Tru(\bquo{\cA\cdot u})$, in every context.  Also $\ConR \eqv
  \neg\Tru(\bquo{\ConR})$.  In particular \rga refutes $\forall x.\, x
  = \tofn 1$ but does not derive $\neg\Tru(\bquo{\forall x.\, x =
  \tofn 1})$.
\end{enumerate}
\end{thm}

\noindent
At a refuted universal, unestablishment is interderivable with
consistency.  At a universal $\cA\cdot u$, the decider affirms exactly
when it finds a certificate, so $\neg\Tru(\bquo{\cA\cdot u})$ says that
no proof certificate exists at any fuel, an unprovability statement.
Part~3 is the argument of [\citen{goedel-paper}, theorem~\textsf{no\_code\_unprovable}](1) with
$\Tru(\bquo{\cA\cdot u})$ in place of $\PrC(\bquo{\varphi})$: under
$\neg\ConR$, ex falso turns a proof of $\z = \tofn 1$ into a
certificate, and \rga proves internally that the decider then affirms
the universal.  The hypothesis of part~3 is used only through weakening, so the
entailment holds at every well-formed $u$, refuted or not, and the
meta form $\emptyset \nvdash \neg\Tru(\bquo{\cA\cdot u})$ follows by
G\"odel's second theorem.  Part~2, a rule whose conclusion is an
entailment, runs the other way.  Inside \rga,
establishment yields an accepted proof
(Theorem~\ref{thm:pl-int-ent}(3)); with the refutation, D1 and EXPL
give $\neg\ConR$; and since $\Tru(\bquo{t})$ is a hypothesis of the
$\Sigma$ shape, contraposition (Corollary~\ref{cor:in-contrapos})
gives $\ConR \vdash \neg\Tru(\bquo{t})$.  The decided-antecedent forms
hold too: $\forall q.\,(\Chk(q, \bquo{t}) \limp \neg\ConR)$ and
$\forall s.\,(\Aff(s, \bquo{t}) \limp \neg\ConR)$ for every refuted $t$
(Appendix~\ref{app:pillars-liar}).  Part~4 combines parts~2 and~3;
$\ConR$ is itself a universal $\cA\cdot u$, so part~3 and the tracking
result of Theorem~\ref{thm:pl-tracking}(3) give $\ConR \eqv
\neg\Tru(\bquo{\ConR})$.  Parts~1 to~3 are the three cases of the
fourth row of Table~\ref{tab:pl-transparency}: a refuted
$\forall$-free sentence is unestablished outright, every refuted
sentence is unestablished under the hypothesis $\ConR$, and the
unestablishment of a universal entails $\ConR$.  \rga can establish
that $\Omega$ and $L$ are never established, because their failure is
visible in their evaluation, but not that a refuted universal is never
established, because that claim would need \rga's own consistency.  Part~1,
the non-derivability in part~3 and the last sentence of part~4 rest on
soundness or on G\"odel's second theorem; the entailments of parts~2
to~4 do not.

\paragraph{Tracking under $\ConR$.}
With $\ConR$ as a hypothesis, the three cases close on the sentences
that \rga proves or refutes: $\Tru$ and $\PrC$ then agree with \rga's
own verdicts exactly, and on every sentence what \mbox{\rga{}${}+\ConR$}
derives about grounded truth is right.  These facts are grounded counterparts
of standard classical facts.  A theory $T$ as in [\citen{goedel-paper}, \S\,``The classical proof of G\"odel's second theorem'']
proves $\ConR_T \limp \neg\PrC(\quo{\varphi})$ for every refutable
$\varphi$, and if it proves $\neg\PrC(\quo{\varphi})$ for any
$\varphi$, it proves $\ConR_T$~\cite{boolos93logic,smorynski77incompleteness}; and the
Kripke--Feferman theory KF proves the T-schema on its bivalent base
language~\cite{halbach11axiomatic}.  Three things are specific to \rga.
The tracking covers grounded truth for \rga's whole language, its own
quantified sentences included, and the bivalent part of that language is fixed by what
\rga proves or refutes, not by syntax.  The dividing line falls at
universals: that a sentence is unestablished is established outright
at refuted $\forall$-free sentences and at visibly divergent ones such
as $\Omega$ and the ordinary Liar, but at a refuted universal that claim is
equivalent to consistency (Theorem~\ref{thm:pl-refuted}(4)), because a
universal is true only with a proof.  And the classical arguments use
implications over undecided antecedents, while here the same arguments go through as
entailments, by contraposition and reductio at hypotheses of the
$\Sigma$ shape.

\begin{thm}[\textsf{con\_tschema\_decided}, \textsf{con\_tracking\_sound}]
\label{thm:pl-tracking}
Let $\Gamma$ be any context.
\begin{enumerate}
\item If $t$ is well-formed and $\emptyset \vdash t$ or $\emptyset
  \vdash \neg t$, then $\ConR, \Gamma \vdash \Tru(\bquo{t}) \liff t$
  and $\ConR, \Gamma \vdash \PrC(\bquo{t}) \liff t$, and $\ConR, \Gamma$
  derives $\Tru(\bquo{t}) \jbool$ and $\PrC(\bquo{t}) \jbool$.
\item For every $t$, if $\ConR \vdash \Tru(\bquo{t})$ then $t$ is
  grounded-true, and if $\ConR \vdash \neg\Tru(\bquo{t})$ then $t$ is
  not grounded-true.
\item $\ConR, \Gamma \vdash G \wedge \neg\Tru(\bquo{G})$, so $\ConR,
  \Gamma \vdash \neg(\Tru(\bquo{G}) \liff G)$; and $\ConR, \Gamma
  \vdash \neg\Tru(\bquo{\ConR})$.
\end{enumerate}
\end{thm}

\noindent
Part~1 holds case by case: at a proved $t$ by D1 and NEC, at a refuted
$t$ by Theorem~\ref{thm:pl-refuted}(2).  Part~1's $\PrC$ forms and its
$\Tru$ biconditional at refuted sentences use no soundness theorem;
part~1's $\Tru$ forms at proved sentences rest on soundness, through NEC.
Part~2 is the $\omega$-soundness of \mbox{\rga{}${}+\ConR$}
([\citen{goedel-paper}, \S\,``The formalized theorems, as entailments'']) applied to $\Tru$.  Part~3 shows that
beyond the sentences \rga proves or refutes the agreement fails, and
\mbox{\rga{}${}+\ConR$} proves that the agreement fails.  Both of part~3's refutations
of truth come from one rule that holds in every context: if $\Gamma
\vdash \neg\PrC(\bquo{t})$, then $\Gamma \vdash \neg\Tru(\bquo{t})$,
because establishment yields an accepted proof.  Under $\ConR$ the
premise holds at $G$ ([\citen{goedel-paper}, theorem~\textsf{goedelG\_under\_ConRGA}]) and at $\ConR$
(formalized G\"odel~II, [\citen{goedel-paper}, theorem~\textsf{rga\_goedel2\_internal\_contra}](2)).  So under
$\ConR$, $G$ is derivable but unestablished, the T-biconditional at
$G$ is refuted, and $\ConR$ itself is unestablished: a hypothesis is
not an establishment.  Appendix~\ref{app:pillars-liar} gives the
proofs.

\begin{rem}
\label{rem:pl-gapped}
\mbox{\rga{}${}+\ConR$} cannot track the ungrounded sentences exactly.
The theorems of \mbox{\rga{}${}+\ConR$} are computably enumerable, and
so are the sentences that \mbox{\rga{}${}+\ConR$} would mark as ungrounded, by refuting the truth
both of the sentence and of its negation.  The ungrounded sentences are
not: the sentences with a value are computably enumerable, and whether
a term has a value is undecidable.  This argument is not mechanized.
\end{rem}

\paragraph{Asserting what is not established.}
Theorem~\ref{thm:pl-tracking}(3) makes \mbox{\rga{}${}+\ConR$} derive a
sentence together with its untruth: $\ConR \vdash G \wedge
\neg\Tru(\bquo{G})$.  Classically this pattern is the feature of the
Kripke--Feferman theory that is often counted as its central cost:
KF proves $L \wedge \neg T\quo{L}$ for its Liar $L$, a sentence
together with the claim that this sentence is not
true~\cite{reinhardt86remarks,field08saving,halbach11axiomatic}, and
PKF and Field's theory are designed to avoid this pattern.  Here the
situation differs in two ways.  First, \rga itself derives no such
pair (Theorem~\ref{thm:b-no-moorean}): if $\emptyset \vdash \varphi$,
then $\emptyset \vdash \Tru(\bquo{\varphi})$ by NEC
(Theorem~\ref{thm:pl-nec}), so by consistency $\emptyset \nvdash
\neg\Tru(\bquo{\varphi})$.  The
conjunction is derivable only under the hypothesis $\ConR$, which is
never grounded-true; on \rga's own reading the derivation establishes
nothing, since assuming is not establishing and $\Tru$ says
``established''.  Second, on the classical reading the conjunction is
true: \mbox{\rga{}${}+\ConR$} is sound for the $\omega$-semantics, in
which $G$ is true and $\Tru(\bquo{G})$ is false, because $G$ is true
but has no certificate.  KF's $L$ and $\neg T\quo{L}$ concern one
notion of truth; $G$ and $\neg\Tru(\bquo{G})$ concern classical truth
and establishment, which this paper keeps apart.

\paragraph{Revenge.}
From the classical side, $L_{\Tru}$ is true: the $\omega$-semantics
says so, and $\Tru$ does not.  Kripke met the same situation in his
own theory, where the metalanguage says that the Liar is not true,
which is what the Liar says~\cite{kripke75outline}.  On the grounded
reading this situation is a disagreement between two notions of truth, not a
defect in one.  That $L_{\Tru}$ is never affirmed is established in
the classical metatheory, in HOL, not in \rga, just as \rga's
consistency is.  \rga neither proves $L_{\Tru}$ nor has a predicate
that calls $L_{\Tru}$ true.  (Contrast $\Omega$ and the ordinary Liar
$L$, whose non-establishment \rga does establish.)

\S\ref{sec:breaks-liars}, next, returns
to $L_{\Tru}$: \rga proves
$\neg\Aff(\tofn{n}, S_{\Tru})$ for every $n$ but not the universal
$\forall s.\ \neg\Aff(s, S_{\Tru})$, and that subsection compares $L_{\Tru}$
with the Liar at the level of evaluation records.

\subsection{Truth, the T-schema, and the Liars}
\label{sec:breaks-liars}

\S\ref{sec:pillars-tarski} showed that Tarski's argument does not apply
to \rga's truth predicate, and that the material T-biconditional at
the Liar is explosive as a hypothesis.  This subsection asks where the
T-biconditional and its halves T-In and T-Out hold as theorems, in the
forms of a conditional of \S\ref{sec:inside-codes}, and then compares
the truth-Liar $L_{\Tru}$ with the G\"odel sentence and the consistency
sentence.  As an implication, the T-biconditional fails wherever
grounded truth gives out.

\begin{thm}[\textsf{tschema\_unprovable\_ungrounded}, \textsf{tin\_iff\_decided}]
\label{thm:br-tschema}
\leavevmode
\begin{enumerate}
\item If $t$ has no value at any fuel, then $\emptyset \nvdash t \limp
  \Tru(\bquo{t})$ and $\emptyset \nvdash \Tru(\bquo{t}) \liff t$.
\end{enumerate}
For well-formed $t$:
\begin{enumerate}
\setcounter{enumi}{1}
\item (T-In) $\emptyset \vdash t \limp \Tru(\bquo{t})$ if and only if
  \rga proves or refutes $t$, that is, if and only if $t$ is
  grounded-true or grounded-false;
\item (T-Out) $\emptyset \vdash \Tru(\bquo{t}) \limp t$ if and only if
  \rga proves $t$ or refutes $\Tru(\bquo{t})$; at a $\forall$-free $t$
  that is grounded-true or grounded-false, $\Gamma \vdash \Tru(\bquo{t})
  \limp t$ in every context $\Gamma$, among them at $\z = \tofn 1$;
\item (TB) $\emptyset \vdash \Tru(\bquo{t}) \liff t$ at every
  grounded-true $t$ and at every grounded-false $\forall$-free $t$.
\end{enumerate}
\end{thm}

\noindent
Part~1 is the disjunction property at work
(Theorem~\ref{thm:br-mech}), as in [\citen{goedel-paper}, theorem~\textsf{no\_goedel\_biconditional}]: by
part~3 of Theorem~\ref{thm:br-mech}, $t \liff t$ fails at such a $t$ too.  The
rule forms are NEC and CONEC (\S\ref{sec:pillars-complete}).  The
$\forall$-free cases rest on the internal refutation of $\Tru$ at every
grounded-false $\forall$-free term (\S\ref{sec:pillars-axioms}).  At
grounded-false sentences with quantifiers \rga does not refute
$\Tru(\bquo{t})$ in general: \rga does not refute that truth at the refuted universal
$\forall x.\, x = \tofn 1$, and refuting the truth of a universal at
all would give $\vdash \ConR$; under the hypothesis $\ConR$ \rga refutes
the truth of every refuted sentence (Theorem~\ref{thm:pl-refuted}).
So the T-biconditional holds at every grounded-true sentence and every
grounded-false $\forall$-free one, fails at every ungrounded sentence,
and fails at that refuted universal (by part~3, and
Table~\ref{tab:truth-axioms}).  T-Out alone has no uniform answer at
ungrounded sentences: \rga proves T-Out in every context at the
self-reducing term $\Omega$ (Theorem~\ref{thm:pl-omega}), but not at
$L_{\Tru}$ (Table~\ref{tab:truth-axioms}).

Where an implication holds, so does the entailment, by modus ponens.
At the ungrounded sentences the entailments are a separate question,
and at the Liar they fail too.

\begin{thm}[\textsf{tin\_\allowbreak ent\_\allowbreak fails\_\allowbreak liar}, \textsf{tout\_\allowbreak ent\_\allowbreak fails\_\allowbreak liar}]
\label{thm:br-tin-ent}
\leavevmode
\begin{enumerate}
\item T-In fails as an entailment at $L_{\Tru}$, $G$ and $\ConR$:
  $L_{\Tru} \nvdash \Tru(\bquo{L_{\Tru}})$, $G \nvdash \Tru(\bquo{G})$
  and $\ConR \nvdash \Tru(\bquo{\ConR})$.
\item T-Out fails as an entailment at $L_{\Tru}$:
  $\Tru(\bquo{L_{\Tru}}) \nvdash L_{\Tru}$.
\item $\neg L_{\Tru} \nvdash \neg\Tru(\bquo{L_{\Tru}})$.
\end{enumerate}
\end{thm}

\noindent
\emph{Proof idea.}  (1) Each hypothesis is true in the
$\omega$-semantics.  The conclusion is not: \rga refutes every instance
$\Aff(\tofn n, c)$ at the code $c$ of a sentence that is never
grounded-true, so $\Tru(c)$ is false there.  (2) Transport along $S_{\Tru} =
\bquo{L_{\Tru}}$ turns the entailment into $\Tru(S_{\Tru}) \vdash
\neg\Tru(S_{\Tru})$, and reductio at this hypothesis of the $\Sigma$
shape would give $\vdash \neg\Tru(S_{\Tru})$, which is $L_{\Tru}$,
against Theorem~\ref{thm:pl-liar-indep}.  (3) runs the same way,
through double negation.  \qed

\medskip
\noindent
The truth rules of PKF~\cite{halbach06axiomatizing}, T-Intro and
T-Elim, are NEC and CONEC with side formulas in the context, and in
\rga the two rules are these entailments; so the two rules fail at the Liar.  This failure is
the expected outcome.  $\Tru$ is grounded truth, and the
$\omega$-semantics, which \rga respects in every context, makes
$L_{\Tru}$ true and $\Tru(\bquo{L_{\Tru}})$ false.
Table~\ref{tab:truth-axioms} gives the entailment form of every rule
of truth.

The two sentences built by diagonalizing over \rga's own predicates
behave alike.  The Liar for the truth predicate, $L_{\Tru}$ of
\S\ref{sec:pillars-tarski}, and the G\"odel sentence $G$ of
{}[\citen{goedel-paper}, Section ``G\"odel's theorems for \rga''] are both $\omega$-incomplete.  \rga proves
every numeral instance $\neg\Aff(\tofn n, S_{\Tru})$ and
$\neg\Chk(\tofn n, S_G)$ of the two universals, yet proves neither
universal, and the $\omega$-semantics of \S\ref{sec:inside-omega},
which applies an $\omega$-rule at each universal, makes both of them
true (Theorem~\ref{thm:br-omega}).  The consistency sentence behaves
the same way.

\begin{thm}[\textsf{goedelW\_omega\_incomplete},
  \textsf{goedelG\_omega\_incomplete}]
\label{thm:br-omega}
Let $S_{\Tru}$ be the code term of $L_{\Tru}$
(\S\ref{sec:pillars-tarski}).  For every $n$,
\[
  \emptyset \vdash \neg\Aff(\tofn n, S_{\Tru}) , \quad
  \emptyset \vdash \neg\Chk(\tofn n, S_G) \quad\text{and}\quad
  \emptyset \vdash \neg\Chk(\tofn n, \bquo{\z = \tofn 1}) ,
\]
but
\[
  \emptyset \nvdash \forall s.\ \neg\Aff(s, S_{\Tru}) , \quad
  \emptyset \nvdash \forall q.\ \neg\Chk(q, S_G) \quad\text{and}\quad
  \emptyset \nvdash \ConR .
\]
The $\omega$-semantics makes $L_{\Tru}$, $G$ and $\ConR$ true.
\end{thm}

\noindent
The companion paper proves that \rga is $\omega$-incomplete~\cite{rga-paper}; the
three sentences of Theorem~\ref{thm:br-omega} are explicit instances.
The ordinary Liar $L$, a term-level fixed point that head-steps to
$\neg L$, with no codes involved, behaves differently.  $L$ has no
grounded value and is neither provable nor refutable
(Appendix~\ref{app:traps-logic}).  Yet \rga proves, by one internal
induction on the fuel, that $L$'s evaluation never yields a value: its
value record $\Rcd_v(\bquo{L})$, the decided fuel-bounded record of a
term's evaluation that Section~\ref{sec:con} defines
(Definition~\ref{def:con-rcd}), is refuted at every value $v$, so
\rga decides that record (Appendix~\ref{app:fold}).  The sentence is ungrounded, and
\rga refutes, at every value, the statement that the sentence has that value.
With the truth predicate the same holds: \rga proves that neither $L$
nor $\neg L$ is true (Theorem~\ref{thm:pl-liar-neither}).  A
Liar built over the value record itself, $L_{\Rcd}$, ``my own record at
the value true fails'', behaves differently again.  \rga decides that
sentence's record exactly when it decides the sentence, and whether
\rga does so is open.
The record evaluates $\forall$-free codes within its fuel bound, but
its quantifier rows read the numeral families outside that bound.  For
that reason the development regards $L_{\Rcd}$ as not comparable with
the ordinary Liar, and does not claim that $L_{\Rcd}$ is a strengthened Liar.

For grounded truth itself, the revenge situation is $L_{\Tru}$
(\S\ref{sec:pillars-tarski}).  This revenge situation has the familiar shape of revenge in
paracomplete theories of truth~\cite{field08saving,beall08revenge}, one
level down: the gap is in \rga's own arithmetic, and the standpoint
that classifies it is HOL.  Which value record, if any, would yield a
strengthened Liar inside \rga is open.

\subsection{The truth predicate against the axioms of truth}
\label{sec:pillars-axioms}

Table~\ref{tab:truth-axioms} sets $\Tru$ against the axioms and rules
of the standard axiomatic theories of truth~\cite{halbach11axiomatic}:
the compositional clauses of the Kripke--Feferman theory
KF~\cite{feferman91reflecting} and of the Friedman--Sheard theory
FS~\cite{friedman87axiomatic}, the iteration axioms T-Rep and T-Del,
the disquotation schemata T-In and T-Out and the T-biconditionals
(TB), the rules NEC and CONEC and their negative forms, the axioms
CONS and COMP, and the intersubstitutivity of PKF~\cite{halbach06axiomatizing}
and of Field's theory~\cite{field08saving}.  An entry is
\emph{internal} when \rga derives one sentence, quantified over codes
by $\cA$ (for T-In, T-Out and TB, one sentence at each $t$).  A
\emph{rule in every context} derives its conclusion $\Gamma \vdash
\dots$ from its premises $\Gamma \vdash \dots$ in every context
$\Gamma$; this is the form that PKF's rules take.  A \emph{rule on
theorems} takes closed theorems $\emptyset \vdash \dots$ as its
premises, as NEC and CONEC do.  The column ``entailment'' gives the
form $A \vdash B$ of \S\ref{sec:inside-codes}, with the variables over
codes free and typed as numbers in the context.  A rule in every
context gives the entailment, by instantiating its context with the
premises, and the entailment gives the rule on theorems, by
substitution and cut; so where a rule on theorems fails, the other two
forms fail as well.
With hypotheses, NEC and CONEC are PKF's rules T-Intro and T-Elim with
side formulas, and their entailment forms are T-In and T-Out.  A form
\emph{fails} when a theorem of HOL says that \rga does not derive the form or
exhibits a sentence at which the form fails, and a form is \emph{open} when the
development has no theorem either way.  In the table, $\dagger$ means
``at codes of $\forall$-free terms only'', $\Omega$ is the looping
term of \S\ref{sec:pillars-tarski}, and $U$ is the refuted universal
$\forall x.\, x = \tofn 1$ of Theorem~\ref{thm:pl-refuted}.
Appendix~\ref{app:pillars-axioms} gives every cell with its proof and
its Isabelle names, the fuel-indexed forms of the clauses for $\Aff$
that hold internally, and the structural laws of $\Tru$.

\begin{table}[!htb]
\centering
\footnotesize
\setlength{\tabcolsep}{2.5pt}
\begin{tabular}{@{}>{\raggedright\arraybackslash}p{0.12\linewidth}>{\raggedright\arraybackslash}p{0.17\linewidth}>{\raggedright\arraybackslash}p{0.18\linewidth}>{\raggedright\arraybackslash}p{0.22\linewidth}>{\raggedright\arraybackslash}p{0.22\linewidth}@{}}
\hline
principle & internal & rule in every context & rule on theorems & entailment \\
\hline
KF: $=$, $\neq$ & holds & holds & holds & holds \\
KF: $\neg\neg$, $\lor$, $\neg\lor$ & fails & holds ($\neg\lor$, $\leftarrow$: $\dagger$; open beyond) & holds & holds ($\neg\lor$, $\leftarrow$: $\dagger$; open beyond) \\
FS: $\neg$ & fails & $\rightarrow$: holds${}^\dagger$, fails at $U$; $\leftarrow$: fails at $\Omega$ & $\rightarrow$: holds${}^\dagger$, fails at $U$; $\leftarrow$: fails at $\Omega$ & $\rightarrow$: holds${}^\dagger$, fails at $U$; $\leftarrow$: fails at $\Omega$ \\
KF, FS: $\forall$; KF: $\neg\forall$ & $\forall$: fails; certificate form of $\leftarrow$: holds (Thm.~\ref{thm:pl-tin-cert}); $\neg\forall$: open & $\forall$, $\rightarrow$: fails at $U$; $\leftarrow$: fails at $L_{\Tru}$; $\neg\forall$: open & $\forall$, $\rightarrow$: holds; $\leftarrow$: only from a proof of the open instance, no $\omega$-rule; $\neg\forall$: holds & $\forall$, $\rightarrow$: fails at $U$; $\leftarrow$: fails at $L_{\Tru}$; $\neg\forall$: open \\
T-Rep, T-Del; KF: $\neg T$ & T-Rep: fails; otherwise open & T-Rep, T-Del: open; KF $\neg T$: fails, as on theorems & hold, except KF $\neg T$: $\rightarrow$ fails at $\Omega$, $\leftarrow$ at $U$ & T-Rep, T-Del: open; KF $\neg T$: $\rightarrow$ fails at $\Omega$, $\leftarrow$ at $U$ \\
head step & open & holds${}^\dagger$ & holds & holds${}^\dagger$ \\
T-In & holds exactly at $t$ with a value & --- & holds (NEC) & fails at $L_{\Tru}$, $G$, $\ConR$ \\
T-Out & holds exactly where \rga proves $t$ or refutes $\Tru(\bquo{t})$; fails at $L_{\Tru}$ and $U$ & --- & holds at well-formed $t$ (CONEC) & fails at $L_{\Tru}$ and $U$; open at $G$, $\ConR$, $L$ \\
TB & holds exactly at grounded-true and $\forall$-free grounded-false $t$; fails at $U$ & --- & --- & fails at $L_{\Tru}$ and $U$ \\
NEC, CONEC & NEC over all proof codes: fails at $\z = \tofn 1$ & NEC: holds for substitution instances; NEC, CONEC as T-Intro, T-Elim with side formulas: hold at images of quoted quantifier-free templates (Thm.~\ref{thm:pl-rf-tin}) & hold (CONEC: at well-formed $t$) & as T-In and T-Out \\
$\neg$-NEC, $\neg$-CONEC & --- & fail where they fail on theorems & $\neg$-NEC: holds${}^\dagger$, fails at $U$, holds under $\ConR$; $\neg$-CONEC: holds at $t$ with a value, fails at $\Omega$ & $\neg$-NEC: fails at $L_{\Tru}$; $\neg$-CONEC: fails at $\Omega$ \\
CONS & fails & as a refutation: holds${}^\dagger$, fails at $U$ & holds (never both theorems) & $\Tru(\mathit{neg}(x)) \vdash \neg\Tru(x)$: holds${}^\dagger$, fails at $U$; the converse refutation: holds${}^\dagger$, open beyond \\
COMP & fails & fails at $L_{\Tru}$ & fails at $L_{\Tru}$ & fails at $L_{\Tru}$ \\
inter\-sub\-sti\-tu\-tivity & --- & fails under $\neg$, at $\Omega$ & holds under $\lor$, $\wedge$, $\neg\neg$; fails under $\neg$, at $\Omega$ & fails at $L_{\Tru}$, and under $\neg$ at $\Omega$ \\
truths are sentences & fails & fails & fails & fails \\
\hline
\end{tabular}
\caption{$\Tru$ against the axioms and rules of truth
  (\S\ref{sec:pillars-axioms}); ``---'': not applicable.
  $\dagger$: at codes of $\forall$-free terms only; $U$: the refuted
  universal $\forall x.\, x = \tofn 1$.  Every cell with its theorem
  is in Table~\ref{tab:c-truth-axioms} of
  Appendix~\ref{app:pillars-axioms}.}
\label{tab:truth-axioms}
\end{table}

\emph{Why the internal forms fail.}  The internal forms, implications
over codes, fail by the disjunction property (Theorem~\ref{thm:br-mech}).  Each failing sentence is, or
implies, a universal over codes of an implication.  The sentence's instance at a code built from the Liar
$L_{\Tru}$ of \S\ref{sec:pillars-tarski}, usually the code term
$S_{\Tru}$ itself, is an implication $a \limp b$ for which \rga
neither refutes $a$ nor proves $b$, and by the disjunction property
such an implication is not a theorem.  The failures persist when the
quantifier is restricted to codes of well-formed terms, by the guard
$\Wfc(x) \neq \z \wedge \IsCode(x) \neq \z$, since \rga proves both
guards at $S_{\Tru}$.  ``Truths are sentences'' fails for a different
reason, at a term that is true but not well-formed.

\begin{thm}[\textsf{true\_\allowbreak not\_\allowbreak sentence}, \textsf{true\_\allowbreak not\_\allowbreak sentence\_\allowbreak internal}]
\label{thm:pl-not-sentence}
Let $j$ be the term whose root is the symbol $\z$ with one argument
$\z$ ($\z$ takes no arguments), and let $t \ldef (\z = \z) \lor j$.
Then $t$ is not well-formed, $\emptyset \nvdash t$, and $\emptyset
\vdash \Tru(\bquo{t})$.  \rga does not prove $\forall x.\ \Tru(x)
\limp \Wfc(x) \neq \z$, nor the same sentence with the consequent
$\Wfc(x) \neq \z \wedge \IsCode(x) \neq \z$.
\end{thm}

\noindent
\emph{Proof idea.}  The semantics evaluates a disjunction in parallel,
in the strong Kleene manner: $a \lor b$ is true as soon as $a$ is
(Table~\ref{tab:in-sem}), whatever $b$ is.  So $t$ is grounded-true, and Theorem~\ref{thm:pl-adeq}(1),
which has no well-formedness hypothesis, gives $\emptyset \vdash
\Tru(\bquo{t})$.  A derivable term is well-formed, so $\emptyset
\nvdash t$.  For the sentence, \rga proves $\Wfc(\bquo{t}) = \z$, so
the instance at $\bquo{t}$ is an implication whose antecedent \rga
proves and whose consequent it refutes.  \qed

\medskip
The rules that fail on theorems fail at the refuted universal $U$, at
two sentences with no value, or at the universal inside one of those two sentences.
At $U$, each failing form would let \rga refute the truth of a
universal, which entails $\ConR$ (Theorem~\ref{thm:pl-refuted}(3)).
At the Liar, \rga proves neither $\Tru(\bquo{L_{\Tru}})$ nor
$\Tru(\bquo{\neg L_{\Tru}})$, so COMP fails as a rule, and \rga does not
prove T-Out at $L_{\Tru}$.  The other sentence is $\Omega$, whose lack
of a value \rga can see, unlike that of $L_{\Tru}$:

\begin{thm}[\textsf{omega\_\allowbreak never\_\allowbreak affirmed}, \textsf{neg\_\allowbreak conec\_\allowbreak fails}]
\label{thm:pl-omega}
The term $\Omega$ satisfies the following.
\begin{enumerate}
\item For every context $\Gamma$, $\Gamma \vdash \forall s.\
  \neg\Aff(s, \bquo{\Omega})$ and $\Gamma \vdash \Tru(\bquo{\Omega})
  \limp \Omega$.
\item $\emptyset \vdash \neg\Tru(\bquo{\Omega})$ and $\emptyset \vdash
  \Tru(\bquo{\neg\Tru(\bquo{\Omega})})$.
\item $\emptyset \nvdash \neg\Omega$ and $\emptyset \nvdash
  \Tru(\bquo{\neg\Omega})$.
\end{enumerate}
\end{thm}

\noindent
\emph{Proof idea.}  Part~1 is an induction inside \rga on the fuel,
like that of Theorem~\ref{thm:in-never}.  The decider reads the cell
of the redex $\Omega$ two levels up as the cell of its reduct, which
is $\Omega$ again, one level up.  So an affirmation of $\Omega$ at
level $v+2$ gives one at level $v+1$, and the motive ``no affirmation
at levels $v$ and $v+1$'' carries the induction; the base cases are
closed computations.  The implication $\Tru(\bquo{\Omega}) \limp
\Omega$ is $\neg\Tru(\bquo{\Omega}) \lor \Omega$, and its left disjunct
is the universal just proved, up to double negation.  Part~2 follows
by double negation and NEC.  Part~3 holds because $\Omega$ has no
value, by soundness and Theorem~\ref{thm:pl-adeq}(2).  The refutation
$\neg\Tru(\bquo{\Omega})$ of part~2 also has a proof that does not use
\rga's soundness theorem: that proof's base case uses only the decider's
correctness on $\forall$-free terms (Appendix~\ref{app:pillars-axioms}).  \qed

\medskip
\noindent
Parts~2 and~3 refute four rules on theorems at $\Omega$: $\neg$-CONEC
($\neg\Tru(\bquo{\Omega})$ without $\neg\Omega$), intersubstitutivity
under $\neg$ (the same pair), the $\leftarrow$ direction of FS's clause
for $\neg$ ($\neg\Tru(\bquo{\Omega})$ without $\Tru(\bquo{\neg\Omega})$),
and KF's clause for $\neg T$ ($\Tru(\bquo{\neg\Tru(\bquo{\Omega})})$
without $\Tru(\bquo{\neg\Omega})$).  Part~1 is T-Out at a sentence with
no value, where T-In fails (Theorem~\ref{thm:br-tschema}).  The
contrast with $L_{\Tru}$ is in part~1: \rga proves each instance
$\neg\Aff(\tofn n, S_{\Tru})$ but not their universal
(Theorem~\ref{thm:br-omega}), while at $\Omega$ \rga proves the
universal.

The table repeats the pattern of [\citen{goedel-paper}, Section ``G\"odel's proof in four forms'']: internal
forms fail where rule forms hold.  The exceptions run both ways.  The
clauses for equations between numerals hold internally, and apart from
``truths are sentences'' and the failures at $U$, the rules that fail
on theorems fail at $L_{\Tru}$ (or the universal inside it) or at
$\Omega$, two sentences with no value.  The entailment column sits
between the internal forms and the rule forms.  KF's clauses for $\neg\neg$ and $\lor$ hold as
entailments in both directions, although their internal forms fail;
there the hypothesis $\Tru(x)$ is the obstacle to the implication, not
to the entailment.  The entailments that fail fail at $L_{\Tru}$, at
$G$ or $\ConR$, at $\Omega$, at $U$, or at the true term that is not a
sentence, and the $\leftarrow$ half of the clause for $\forall$, a
formalized $\omega$-rule, fails even as an entailment, and so in every
context (\S\ref{sec:pillars-reflect}, ``T-In at a universal'').  So PKF's truth
rules with side formulas fail at the Liar, as \S\ref{sec:breaks-liars}
explains.  A failure at $U$ is a witness, not a claim about every
grounded-false sentence with quantifiers; Section~\ref{sec:concl}
lists the cells that remain open.

\subsection{The truth predicate as a reflection medium}
\label{sec:pillars-reflect}

In proof by reflection, a decision procedure written and proved
correct inside a proof assistant turns each goal of a class into a
computation~\cite{boutin97using}.  Coq's tactic for ring identities is
the standard example~\cite{gregoire05proving}.  This subsection shows
that \rga supports the same pattern, with its own truth predicate as
the medium.  \rga derives disquotation at a free variable for a class
of quantifier-free formulas: from a derivation of an instance to the
truth of the template's image, the code that \rga computes for the
instance, and back.  \rga also contains a
normalizer for polynomials in one variable, written as explicit
primitive-recursive machines and proved correct by one internal
induction.  The two parts combine into one \rga sentence of
reflection, and into a derived rule that proves polynomial identities
at a free variable by computation.  At a universal, \rga proves that
truth is the existence of an accepted proof of the open instance, and
an explicit machine builds such proofs for the instances of a closed
universal theorem.  Appendix~\ref{app:pillars-reflect} gives the
proofs.

\paragraph{Three kinds of reflection.}
The word ``reflection'' names three different things, with different
costs.
\begin{enumerate}
\item \emph{Logical reflection principles}, such as local reflection
  $\PrC(\quo{\varphi}) \limp \varphi$ for every sentence $\varphi$, add
  strength.  A theory that proves the reflection principle for another
  theory is properly stronger than that
  theory~\cite{kreisel68reflection,beklemishev05reflection}\cite[p.~19]{harrison95metatheory}.
  The G\"odel paper studies the corresponding step for \rga, the
  adoption of \rga's own checker~[\citen{goedel-paper}, Section ``Adoption: the reflection rule, and the consistency family''].
\item \emph{Computational reflection through provability}, in
  Harrison's formulation~\cite[pp.~19--20]{harrison95metatheory}, proves inside the
  system that a decision procedure $\mathit{DP}$ is correct for
  provability, $\vdash \mathit{DP}(\quo{\varphi}) \limp
  \PrC(\quo{\varphi})$, and then adds the rule from $\vdash
  \mathit{DP}(\quo{\varphi})$ to $\vdash \varphi$.  The added rule
  yields no new theorems when the system is 1-consistent, that is,
  when every $\Sigma_1$ sentence it proves is true; by G\"odel's second
  theorem the system does not prove its own 1-consistency.  Davis and Schwartz's
  metamathematical extensibility accepts a new rule once the system
  proves the rule's justifying sentence, which says that the rule takes
  theorems to theorems, and the system then applies the rule directly:
  the medium is provability, no reflection principle joins the system,
  and no proof is replayed~\cite{davis79metamathematical}.  Davis and
  Schwartz's proof that the extended system has no new theorems assumes
  weak $\omega$-consistency, a form of 1-consistency, so that, as they
  note, the proof cannot be formalized in the system, and their compiled
  evaluators rest instead on correctness for a truth predicate of a
  decidable fragment~\cite[\S6]{davis79metamathematical}.  Harrison
  counts Davis and Schwartz among the earliest advocates of
  computational reflection~\cite[p.~20]{harrison95metatheory}.
\item \emph{Computational reflection through a denotation}, as in
  Boutin's work, whose own case is the normalization of commutative-ring
  expressions in Coq, proves inside the system that the procedure is
  correct for a denotation function, and lets the system's own
  computation rules run the procedure.  Nothing is added to the
  system, so this kind of reflection is conservative outright.
\end{enumerate}
The three-way split is this paper's.  Harrison defines computational
reflection by the provability form and files the denotation form under
the same name in his account of the Boyer--Moore prover's
metafunctions, and he writes that the idea of computational reflection
``is not to make a formal system stronger, but rather to make its
deductive process more efficient''~\cite[pp.~20 and~27]{harrison95metatheory}.
\rga's construction is of the third kind, with the decider's values
and truth as the denotation.  The medium is the truth predicate $\Tru$, a search
over the compiled decider, $\exists s.\ \Aff(s, c)$
(Definition~\ref{def:in-truth}), not itself a compiled program.

\paragraph{No reify/denote pair.}
In Coq a reflective tactic needs, besides the procedure, an inductive
type of syntax, a denotation function from that syntax to the terms of
the goal, and a lemma that relates the two, together with a
metalevel step, reification, that maps each goal to its syntax.  \rga
needs no separate syntax type, no separately defined denotation and
no separate lemma that relates the two.  The syntax is \rga's own term
codes.  The denotation of a term code is the value that the compiled
decider reads at the code's image under an environment, and the
denotation of a formula code is $\Tru$ at the image.  The step from
the denotation back to the goal is open disquotation
(Theorems~\ref{thm:pl-rf-tin} and~\ref{thm:c-rf-eq}), theorems about
$\Tru$ that hold for every template, not a lemma written for the
procedure.  The metalevel
still quotes: HOL supplies the code numerals of the two sides of an
identity and checks that the two sides belong to the procedure's
class.

\paragraph{Templates and images.}
The classes are HOL classes of terms with no free variable other
than $\vz$.

\begin{defi}[\textsf{ptm}, \textsf{qfm}]
\label{def:pl-rf-tpl}
The \emph{term templates} are built from the variable $\vz$, the
numerals $\tofn k$, successor $\suc t$, and the compiled applications
$\cmp(f)\cdot a\cdot b$ of term templates $a$ and $b$, for every index
$f$.  The \emph{formula templates} are the equations $a = b$ between
term templates, closed under $\neg$ and $\lor$.
\end{defi}

\noindent
Since $f$ ranges over all indices, the term templates contain every
binary primitive-recursive function, not only addition and
multiplication; a template may contain the compiled decider at a fixed
fuel term.  \rga's conjunction and implication are defined from $\neg$ and
$\lor$, so the formula templates contain them as well.  Templates are
quantifier-free in the sense that they contain no $\cA$.  The compiled
constants $\cmp(f)$ inside them are closed $\lambda$-terms, with bound
variables of their own.  A template never contains $\Tru$, whose
definition contains $\cA$.

For a term $v$, the environment code $\mathit{env}(v)$ maps the
variable $\vz$ to the code of the numeral of $v$.  \rga computes
$\mathit{env}(v)$ from the term $v$, which may be a free variable.  The
\emph{image} of a template $t$ at $v$ is $\Subst(\mathit{env}(v),
\bquo{t})$, and the \emph{instance} $\ssub{\vz}{v}{t}$ is HOL's
substitution of $v$ for $\vz$ in $t$.  The difference from the
soundness pillar lies in the environment.  Theorem~\ref{thm:pl-sound}
concerns proofs from no hypotheses, at a quoted environment
$\bquo{\zeta}$, the code numeral of a HOL list.  Here the environment
is an \rga term that is open whenever $v$ is.  The value that the decider
reads at a code is a $\Sigma$ statement of the same shape as $\Tru$:
\[
  \mathit{DecVal}(c, x) \;\ldef\; \exists l.\ \Dec(l, c) = \suc x .
\]
The packed verdict~2 is an affirmation, so $\mathit{DecVal}(c, \tofn
1)$ gives $\Tru(c)$ (Appendix~\ref{app:pillars-reflect}).

\paragraph{Open disquotation.}
The decider reads at the image of every template the value of its
instance, with $v$ free.

\begin{thm}[\textsf{ptm\_val}, \textsf{qfm\_val}]
\label{thm:pl-rf-val}
Let $\Gamma \vdash v\jnat$.
\begin{enumerate}
\item For every term template $t$, $\Gamma \vdash
  \mathit{DecVal}(\Subst(\mathit{env}(v), \bquo{t}),
  \ssub{\vz}{v}{t})$.
\item For every formula template $p$, \rga decides the instance:
  $\Gamma \vdash \ssub{\vz}{v}{p} \lor \neg\ssub{\vz}{v}{p}$.  If $\Gamma \vdash
  \ssub{\vz}{v}{p}$, then $\Gamma \vdash
  \mathit{DecVal}(\Subst(\mathit{env}(v), \bquo{p}), \tofn 1)$, and if
  $\Gamma \vdash \neg\ssub{\vz}{v}{p}$, then $\Gamma \vdash
  \mathit{DecVal}(\Subst(\mathit{env}(v), \bquo{p}), \z)$.
\end{enumerate}
\end{thm}

\noindent
Both parts are schemata, proved by HOL induction over the template.
The application case is a schema in the index $f$ as well, because the
value law of a compiled application (Theorem~\ref{thm:c-rf-eval}) is
proved by HOL induction over the index tree.  So the theorem is a
family of \rga derivations, one per template, and not one sentence
over codes.  In each derivation $v$ stays free.  In part~2, the
value $\tofn 1$ gives T-In and the decidedness gives T-Out.

\begin{thm}[\textsf{qfm\_tin}, \textsf{qfm\_tout}]
\label{thm:pl-rf-tin}
For every formula template $p$, every context $\Gamma$ and every term
$v$ with $\Gamma \vdash v\jnat$:
\begin{enumerate}
\item (T-In) if $\Gamma \vdash \ssub{\vz}{v}{p}$, then $\Gamma \vdash
  \Tru(\Subst(\mathit{env}(v), \bquo{p}))$;
\item (T-Out) if $\Gamma \vdash \Tru(\Subst(\mathit{env}(v),
  \bquo{p}))$, then $\Gamma \vdash \ssub{\vz}{v}{p}$.
\end{enumerate}
\end{thm}

\noindent
On the classical reading, both parts are HOL statements about
derivability, rules in every context.  On \rga's own reading, each
instance is an \rga derivation of the disquotation step at a free
variable, since $v$ may be one.  \emph{Proof idea.}  T-In is the value
$\tofn 1$ of Theorem~\ref{thm:pl-rf-val}(2), turned into truth.  T-Out
splits on the decided instance.  Under the hypothesis $\neg\ssub{\vz}{v}{p}$,
the image has the value $\z$, so the code of the image's negation is true.  The
image is the code of a $\forall$-free term, and the decider never
affirms both such a code and its negation, so the two truths give any
well-formed conclusion by ex falso, here $\ssub{\vz}{v}{p}$
(Appendix~\ref{app:pillars-reflect}).

T-Out holds at the images of quoted templates, at an arbitrary term
$v$, and not at an arbitrary open code.  The internal sentence below
quantifies over codes, but the sentence's consequent stays $\Tru$.
HOL takes the step to the identity afterwards, per template.  The
development does not state the implication form $\Gamma \vdash
\ssub{\vz}{v}{p} \limp \Tru(\Subst(\mathit{env}(v), \bquo{p}))$
either, although every instance is decided.
Table~\ref{tab:truth-axioms} lists T-In and T-Out as one sentence at
each $t$, with no entry in the column ``rule in every context''.
Theorem~\ref{thm:pl-rf-tin} bears on the row NEC, CONEC instead.  With
side formulas, as in PKF, both rules hold in every context at the
images of quoted quantifier-free templates, at an arbitrary term $v$,
and the table's cell ``rule in every context'' for that row records
this class.  Theorem~\ref{thm:pl-rf-tin} changes no other cell.

The normalizer below reasons about arbitrary codes, not only about
images of templates.  Its application step needs the value law at
argument codes that \rga knows only through guards.  The guard
$\mathit{IsQf}(c)$ is the compiled test that $c$ codes a
$\forall$-free term, nonzero exactly at such codes, beside the tests
$\Wfc$ and $\IsCode$ of \S\ref{sec:pillars-int}.  The code constructor
$\mathit{ap}(a, b)$ computes the code of an application from the codes
of its two parts, and $\mathit{eq}(a, b)$ that of an equation.

\begin{thm}[\textsf{iPrDecVal\_\allowbreak app\_\allowbreak open},
  \textsf{iPrDecVal\_\allowbreak app\_\allowbreak img}]
\label{thm:pl-rf-app}
Let $f$ be an index.  Let $\Gamma$ derive $c_i\jnat$, $\Wfc(c_i) \neq
\z$, $\IsCode(c_i) \neq \z$ and $\mathit{IsQf}(c_i) \neq \z$ for $i
= 1, 2$, and $x\jnat$ and $y\jnat$.  Let $A \ldef
\mathit{ap}(\mathit{ap}(\bquo{\cmp(f)}, c_1), c_2)$.
\begin{enumerate}
\item If $\Gamma \vdash \mathit{DecVal}(c_1, x)$ and $\Gamma \vdash
  \mathit{DecVal}(c_2, y)$, then $\Gamma \vdash \mathit{DecVal}(A,
  \cmp(f)\cdot x\cdot y)$.
\item If $\Gamma \vdash v\jnat$, $\Gamma \vdash
  \mathit{DecVal}(\Subst(\mathit{env}(v), c_1), x)$ and $\Gamma \vdash
  \mathit{DecVal}(\Subst(\mathit{env}(v), c_2), y)$, then $\Gamma
  \vdash \mathit{DecVal}(\Subst(\mathit{env}(v), A), \cmp(f)\cdot
  x\cdot y)$.
\end{enumerate}
\end{thm}

\noindent
The codes $c_1$, $c_2$ and the values $x$, $y$ are arbitrary terms,
open ones included.  The theorem is still a schema in $f$, not one
sentence over indices, and the theorem needs the $\forall$-free guard.  The one
new ingredient beyond the value law of Theorem~\ref{thm:c-rf-eval} is
an internal law: lifting by zero is the identity on codes that pass
the tests $\Wfc$ and $\IsCode$ (Appendix~\ref{app:pillars-reflect}).

\paragraph{A polynomial normalizer inside \rga.}
A polynomial in one variable is its list of coefficients, lowest
degree first, in the list coding of \S\ref{sec:inside-codes}.  Six
machines act on these lists: Horner evaluation $\mathit{ev}(v, l)$,
scaling by a constant, addition, multiplication by shift and add,
trimming of trailing zeros, and equality after trimming.  Each machine
is an explicit primitive-recursive tree built on one skeleton, a
course-of-values recursion, whose step may read any earlier value.
The skeleton recurs on the second argument with two lookups into the
history, and the equality machine applies an equality test to two
trimmings.  Each machine has a HOL list function as its semantics.
\rga derives the skeleton's unfolding law once, at open codes, and
from it the machines' laws.  \rga also derives, at open values, the
ring laws of the compiled multiplication and of an addition defined
by recursion, which \rga proves equal to the compiled addition at
numbers.  Among these laws are commutativity and associativity of
multiplication and both distributive laws.  The evaluation homomorphisms come next: at open
codes and values, in every context, Horner evaluation commutes with
scaling, addition, multiplication and trimming, each by one internal
course-of-values induction.  The ring laws enter there, once, and not
per identity.

The normalizer $\mathit{Norm}(x)$ is an explicit tree of the same
skeleton over a term code.  Its value is a pair $\langle
\mathit{tag}, \mathit{poly} \rangle$, where tag $0$ means ``a
polynomial term''.  At the code $0$, the code of $\vz$, the value is
$\langle 0, [0, 1] \rangle$.  At a nonzero code the normalizer runs
four tests in order and otherwise rejects: the code of $\z$ gives $\langle 0, [\,]
\rangle$; a successor adds $[1]$ to its argument's polynomial; an
application of the compiled addition or multiplication adds or
multiplies the polynomials of its two arguments, with the sum of
their tags as its tag; any other code is rejected with the tag~3.  Each
test is a round trip: the code equals the constructor applied to the
code's own decoded parts.  The tests compare the head of an application with the fixed code
numerals of the compiled addition and multiplication, and the
normalizer never decodes those heads.  On the HOL
side the \emph{polynomial templates} are the term templates built from
$\vz$, $\z$, successor and the compiled addition and multiplication
alone, a narrower class than the term templates, with numerals only as
towers of successors.  On every polynomial template $t$, HOL computes
that $\mathit{Norm}$ yields tag $0$ and the HOL normal form of $t$
(Theorem~\ref{thm:c-rf-hol}).

\begin{thm}[\textsf{norm\_correct}, \textsf{nmCor\_allI}]
\label{thm:pl-rf-norm}
Let $\Gamma$ derive $v\jnat$, $x\jnat$, $\Wfc(x) \neq \z$,
$\IsCode(x) \neq \z$, $\mathit{IsQf}(x) \neq \z$ and
$\mathit{fst}(\mathit{Norm}(x)) = \z$.  Then
\[
  \Gamma \vdash \mathit{DecVal}(\Subst(\mathit{env}(v), x),\,
  \mathit{ev}(v, \mathit{snd}(\mathit{Norm}(x)))) .
\]
The judgment follows from one sentence $\forall u.\ C(u)$ that \rga
proves from no hypotheses.  $C(u)$ states the conclusion at the two
components $v = \mathit{fst}(u)$ and $x = \mathit{snd}(u)$, under the
three code guards and the tag guard as decided antecedents of
implications.  The judgment is the sentence's instance at the pair of
$v$ and $x$, with the four guards discharged.
\end{thm}

\noindent
Both $v$ and $x$ are arbitrary terms, open ones included, and the theorem
has no HOL premise on $x$.  The class enters only through the tag: a
code that passes no test of the normalizer gets the tag 3, and a
compound code adds the tags of its parts.  Correctness is stated
relative to the decider's value at the image, not relative to an
equality of terms.  The normalizer recognizes only the compiled
addition and multiplication, and only one variable.  \emph{Proof
idea.}  One internal course-of-values induction over the pair
$\langle v, x \rangle$ proves the sentence.  Its step dispatches the
code $0$, the four tests and the default.  The application case uses
Theorem~\ref{thm:pl-rf-app}(2) and the homomorphisms, and the
rejecting case is ex falso, since the tag would be 3 while the guard
says 0.

\paragraph{Proof by reflection.}
The test $\mathit{polyEq}(x, y)$ accepts two codes when both
normal forms have the tag $0$ and their coefficient lists agree after
trimming.

\begin{thm}[\textsf{poly\_\allowbreak reflect\_\allowbreak sentence}]
\label{thm:pl-rf-sentence}
\rga proves, from no hypotheses, the sentence
\[
  \begin{array}{@{}l@{}}
  \forall x.\,\forall y.\,\forall v.\ \Wfc(x) \neq \z \limp \IsCode(x)
  \neq \z \limp \mathit{IsQf}(x) \neq \z \limp \Wfc(y) \neq \z \limp
  \IsCode(y) \neq \z \limp {} \\
  \qquad \mathit{IsQf}(y) \neq \z \limp \mathit{polyEq}(x, y) \neq \z
  \limp \Tru(\Subst(\mathit{env}(v), \mathit{eq}(x, y))) ,
  \end{array}
\]
with the implications nested to the right.
\end{thm}

\noindent
Every antecedent is a decided test, a compiled term compared with
$\z$, so \rga's restriction of implication to decided antecedents
(Theorem~\ref{thm:in-dt}) is met.  \emph{Proof idea.}  After the
three universal introductions, the seven implications are introduced
one by one, each by the deduction theorem at a decided hypothesis.  The body
is the correctness of the test: from the guards and the test's
acceptance, Theorem~\ref{thm:pl-rf-norm} at $x$ and at $y$ and the
homomorphism for trimming give the decider's value $\tofn 1$ at the
image of $\mathit{eq}(x, y)$, and that value gives truth.

\begin{thm}[\textsf{poly\_reflect}]
\label{thm:pl-rf-reflect}
Let $a$ and $b$ be polynomial templates whose HOL normal forms agree
after trimming.  For every context $\Gamma$ and every $v$ with $\Gamma
\vdash v\jnat$, $\Gamma \vdash \ssub{\vz}{v}{a} = \ssub{\vz}{v}{b}$.
\end{thm}

\noindent
The premise on the normal forms is an equation between two HOL lists
of numbers, which HOL checks by evaluation.  \emph{Proof idea.}  \rga's
derivation of an identity is assembled anew for each identity.  The test's value at
the code numerals $\bquo{a}$ and $\bquo{b}$, $\emptyset \vdash
\mathit{polyEq}(\bquo{a}, \bquo{b}) = \tofn 1$, comes from a closed evaluation through open completeness
(Theorem~\ref{thm:in-cmp-num}), and the code guards at these numerals
are closed computations as well.  The proof weakens the
sentence of Theorem~\ref{thm:pl-rf-sentence} to $\Gamma$, instantiates
the sentence at $\bquo{a}$, $\bquo{b}$ and $v$, and discharges its
antecedents.  Transport along $\mathit{eq}(\bquo{a}, \bquo{b}) =
\bquo{a = b}$ gives $\Tru(\Subst(\mathit{env}(v), \bquo{a = b}))$, and
the step back from truth to the identity is open T-Out for one
equation (Theorem~\ref{thm:c-rf-eq}), a HOL induction over the
template.  Only the sentence in the middle is one fixed internal
theorem, and the proof applies the sentence at the unary numerals of the
codes.

\begin{cor}[\textsf{rfl\_binom\_by\_reflection}]
\label{cor:pl-rf-binom}
For every context $\Gamma$ and every $v$ with $\Gamma \vdash v\jnat$,
\[
  \Gamma \vdash \suc v \times \suc v = (v \times v + \tofn 2 \times v) +
  \tofn 1 ,
\]
where $\times$ and $+$ are the compiled multiplication and addition.
\end{cor}

\noindent
The proof checks membership in the class of polynomial templates by
the class's introduction rules, checks the agreement of the two normal
forms by one evaluation in HOL, and applies
Theorem~\ref{thm:pl-rf-reflect}.  No algebraic derivation of the
identity is built.  The identities $(v+1)^3 = v^3 + 3v^2 + 3v + 1$ and
$(v+2)(v+3) = v^2 + 5v + 6$, where $v+k$ abbreviates $\suc^k v$, $v^2$
abbreviates $v \times v$ and $v^3$ abbreviates $v \times (v \times
v)$, follow the same way (Theorem~\ref{thm:c-rf-more}).  Open T-Out also derives Corollary~\ref{cor:pl-rf-binom}'s
identity by hand, but only from the identity's truth taken as a
hypothesis (Theorem~\ref{thm:c-rf-eq}(2)).  A dependency check of the
build confirms that the proof of Corollary~\ref{cor:pl-rf-binom} does
not use that derivation.

\paragraph{What the construction buys, and what it does not.}
On the classical reading, the reflection construction is free in a
precise sense.  It costs no proof effort per identity beyond
computation, and no new principle.  The construction uses no
reflection rule, no proof checker and no adoption of a checker, and
the sentence of Theorem~\ref{thm:pl-rf-sentence} is a theorem of \rga
itself, whose medium is truth, not provability.  On \rga's own
reading, the work per identity is computation, not an algebraic
derivation.  HOL compares the two normal forms by evaluation, and \rga
evaluates the compiled test.
The ring laws enter once, in the homomorphisms behind
Theorem~\ref{thm:pl-rf-norm}.  What reflection through $\Tru$ buys is
one internal uniform sentence and a decision procedure inside \rga.

This construction does not buy new theorems.  The identities are \rga
theorems already, and the reflection construction adds no rule to
\rga, so it adds no theorems.  Nor does reflection through $\Tru$ buy
small certificates, the feature that makes proof by reflection efficient in
type theory.  \rga's conversion rules move only along single head
steps and do not evaluate, so the closed evaluations of the test and
of the guards reach \rga through open completeness.  Read informally,
the derivation that open completeness builds follows the evaluation
step by step; the development states no bound on that derivation's
size.

\begin{rem}[not mechanized]
Theorem~\ref{thm:pl-rf-reflect}, at the context $\vz\jnat$ and $v =
\vz$, followed by universal introduction, gives $\emptyset \vdash
\forall \vz.\,a = b$.  NEC (Theorem~\ref{thm:pl-nec}) then gives $\Tru$
of that sentence's code.  This step is NEC at a universal.  The step
adds no theorem, and the reflected class stays quantifier-free, so the
construction does not cross a quantifier level.
\end{rem}

\emph{Classical versus \rga.}  Classically, truth for quantifier-free
formulas over every primitive-recursive function symbol is not itself
primitive recursive, since no primitive-recursive function is
universal for the primitive-recursive functions.  That truth is
$\Sigma_1$: a formula is true when a terminating computation says so.
For each formula the theory $\mathrm{I}\Sigma_1$ proves the T-schema of
this $\Sigma_1$ truth definition, as part of the classical theory of
partial truth definitions~\cite{hajek93metamathematics}.  $\Tru$ has
exactly this shape at the codes the construction uses, a fuel and an
affirming decider run, and the value law of
Theorem~\ref{thm:c-rf-eval} is a schema in the index, matching the
per-formula nature of the classical T-schema.  Read classically, the
reflection construction is this partial truth definition at work, and
nothing in it needs grounded truth.  The guard $\mathit{IsQf}$ confines
$\Tru$ to codes of $\forall$-free terms.

On \rga's own reading, the grounded notion would show only across a
quantifier level: one type-free truth predicate serving a reflected
class with quantifiers, for instance a decision procedure for a class
of universal sentences.  Such a procedure needs T-In at a universal,
the $\leftarrow$ direction of the clause for $\forall$, in a form that
holds in every context.

\paragraph{T-In at a universal.}
For a code $c$, let $J(c)$ be the code, computed inside \rga from $c$,
of the open instance of the universal whose body has the code $c$: the
judgment with the hypothesis $\vz\jnat$ and the conclusion
$\mathit{ap}(\mathit{lift}(c), \bquo{\vz})$, where $\mathit{lift}$ is
the compiled de Bruijn lift by one.  As in
Table~\ref{tab:c-truth-axioms}, $\mathit{ap}(\bquo{\cA}, c)$ is the
code of that universal, and $\pi$ is the proof test of
Definition~\ref{def:pl-cert}.

\begin{thm}[\textsf{tru\_\allowbreak univ\_\allowbreak cert\_\allowbreak int}, \textsf{tru\_\allowbreak univ\_\allowbreak cert\_\allowbreak conv\_\allowbreak int}]
\label{thm:pl-tin-cert}
\rga proves, from no hypotheses, the two sentences
\[
  \forall c.\,\forall p.\ \cmp(\pi)\cdot J(c)\cdot p \neq \z \limp
  \Tru(\mathit{ap}(\bquo{\cA}, c))
  \qquad\text{and}\qquad
  \forall c.\,\forall s.\ \Aff(s, \mathit{ap}(\bquo{\cA}, c)) \neq \z
  \limp \Cert(J(c)) .
\]
\end{thm}

\noindent
Both antecedents are decided tests, and the body code is free.  On
\rga's own reading, inside \rga the truth of a universal is the
existence of an accepted proof of the universal's open instance.  The
$\omega$-form of this direction, with the truths of all numeral
instances in place of the proof, fails internally and as an entailment
(Table~\ref{tab:truth-axioms}).  On the classical reading, the
$\omega$-form fails as a rule in every context too, at the matrix $u_L
\ldef \lambda\neg\Aff(\vz, S_{\Tru})$ of $L_{\Tru}$
(Theorem~\ref{thm:c-tin-omega}), which settles the table's cell ``rule
in every context'' for the direction $\leftarrow$.  On \rga's own
reading, local reflection into $\Tru$ holds at every closed universal:
$\Chk(\vz, \bquo{\cA\cdot u}),\ \vz\jnat \vdash \Tru(\bquo{\cA\cdot
u})$ for every well-formed closed $u$.  On the classical reading, the
same entailment at $\z = \tofn 1$ is not derivable, since it would
give the sentence that Theorem~\ref{thm:pl-no-uniform} excludes
(Theorem~\ref{thm:c-tin-reflect}).  From a closed theorem
about all instances, \rga builds the proofs that
Theorem~\ref{thm:pl-tin-cert} asks for:

\begin{thm}[\textsf{univ\_\allowbreak inst\_\allowbreak nec}]
\label{thm:pl-tin-nec}
Let $G$ be a well-formed closed term with $\emptyset \vdash \forall
x.\,\forall v.\ G\cdot x\cdot v$.  Then
\[
  \emptyset \vdash \forall x.\ \Tru(\mathit{ap}(\bquo{\cA},
  \mathit{ap}(\bquo{G}, \Num\cdot x))) .
\]
\end{thm}

\noindent
At each numeral $\tofn k$, \rga proves the code
$\mathit{ap}(\bquo{\cA}, \mathit{ap}(\bquo{G}, \Num\cdot\tofn k))$
equal to $\bquo{\cA\cdot(G\cdot\tofn k)}$, so the theorem is NEC at the
universal instances, as one sentence; NEC itself
(Theorem~\ref{thm:pl-nec}) gives one theorem per numeral.  On \rga's
own reading, the proof is a construction: one explicit
primitive-recursive machine, with no index chosen by Hilbert choice,
computes from $x$ a proof of the open instance, and \rga verifies, at
an arbitrary $x$, that the checker accepts that proof
(Theorem~\ref{thm:c-tin-builder}).  The only new internal law is that
the compiled lift fixes numeral codes.  On the classical reading, the
machine makes the passage from the closed theorem to the truth of its
universal instances an explicitly verifiable rule in Artemov's
sense~\cite[Def.~13]{artemov99explicit}: a rule from premises $\Gamma$
to a conclusion $\varphi$ together with a total computable function
$f$ on proof codes such that the system proves that $f$ maps every
proof of $\Gamma$ to a proof of $\varphi$.  Such a rule need only be
admissible, not derivable, and Artemov shows that adding one gives a
system whose proofs the original system can translate into its
own~\cite[Thms.~9 and~11]{artemov99explicit}.  The
built proof contains a certificate that the numeral of $x$ is a
number, so the proof grows with $x$, and the construction is not
reflection with small certificates.  The sentence of
Theorem~\ref{thm:pl-tin-nec} keeps any guards of the matrix inside the
universal.  A decision procedure for a class of universal sentences
needs the decided guards outside the inner universal, discharged inside
the built proof.

\begin{thm}[\textsf{poly\_\allowbreak all\_\allowbreak sentence}, \textsf{tru\_\allowbreak univ\_\allowbreak by\_\allowbreak reflection}]
\label{thm:pl-tin-disch}
Let $Q \ldef \lambda\lambda\lambda\,\Tru(\Subst(\mathit{env}(\vz),
\mathit{eq}(v_2, v_1)))$, so that $Q\cdot x\cdot y\cdot v$ states the
truth of the image of the equation $\mathit{eq}(x, y)$ at $v$.
\begin{enumerate}
\item \rga proves, from no hypotheses, the sentence
\[
  \begin{array}{@{}l@{}}
  \forall x.\,\forall y.\ \Wfc(x) \neq \z \limp \IsCode(x)
  \neq \z \limp \mathit{IsQf}(x) \neq \z \limp \Wfc(y) \neq \z \limp
  \IsCode(y) \neq \z \limp {} \\
  \qquad \mathit{IsQf}(y) \neq \z \limp \mathit{polyEq}(x, y) \neq \z
  \limp \Tru(\mathit{ap}(\bquo{\cA}, \mathit{ap}(\mathit{ap}(\bquo{Q},
  \Num\cdot x), \Num\cdot y))) ,
  \end{array}
\]
with the implications nested to the right.
\item For polynomial templates $a$ and $b$ whose HOL normal forms agree
  after trimming, \rga proves, from no hypotheses,
  $\Tru(\bquo{\cA\cdot(Q\cdot\bquo{a}\cdot\bquo{b})})$.
\end{enumerate}
\end{thm}

\noindent
The sentence of part~1 stands one quantifier level above
Theorem~\ref{thm:pl-rf-sentence}, with the same seven guards outside
the inner universal.  On the classical reading, the $\omega$-form of
T-In cannot give this sentence, since that form fails as an entailment
(Table~\ref{tab:truth-axioms}).  On \rga's own reading, the proof is a
construction at two parameters, with generic forms in
Theorems~\ref{thm:c-tin-univ2} and~\ref{thm:c-tin-disch}.  One explicit
primitive-recursive builder computes from $x$ and $y$ a proof of the
open instance, and \rga verifies, under the internal hypothesis that
the test is zero, that the checker accepts that proof.  The builder
chains the eliminations of the universal at the two numerals, the
value certificate of the guard test, moved to zero by that hypothesis
as an equality of codes, modus ponens at abstract formula codes with
one-element contexts only (Theorem~\ref{thm:c-tin-mp}), and the open
top.  Part~2 is the guarded form at $\bquo{a}$ and $\bquo{b}$,
discharged by one closed evaluation of the guards and the test, and
holds at the three identities of Corollary~\ref{cor:pl-rf-binom} and
Theorem~\ref{thm:c-rf-more}.  The built proof grows with $x$ and $y$
and contains the value certificate of the test.  Part~2 is not the
object identity inside $\Tru$, $\emptyset \vdash
\Tru(\bquo{\cA\cdot\lambda(a = b)})$, which needs a builder that also
emits the steps of open T-Out; the development does not have one
(Section~\ref{sec:concl}).

\section{Consistency, from unconditional to conditional}
\label{sec:con}

\rga does not prove its own consistency sentence
$\ConR \ldef \forall p.\ \neg\Chk(p, \bquo{\z = \tofn 1})$
([\citen{goedel-paper}, theorem~\textsf{rga\_goedel2}]).  This section states what \rga does prove
about the consistency of its proofs: a family of restricted
consistency theorems (Table~\ref{tab:con-tower}).
Appendix~\ref{app:fold} gives the folds and the fragment and residual
results in full, Appendix~\ref{app:oga} the level theorems and their
discharge in \oga, and Appendix~\ref{app:traps} the tempting designs
that fail.

Two of the theorems are internal.  \rga proves, in the empty context
and with no hypothesis, one consistency sentence for its
$\forall$-free proofs (\S\ref{sec:con-qf}) and one for a second
decided class of proofs, the \emph{fragment} $\Frag$, which admits
quantifiers away from a few kinds of positions (\S\ref{sec:con-frag});
the $\forall$-free theorem follows from the fragment theorem inside
\rga.  The \emph{level theorem} covers a third decided class, $\LvC$, which
admits quantified case hypotheses and critical formulas
(\S\ref{sec:con-frag}) whose quantifiers nest, through head steps, to
a certified depth.  \rga derives that class's consistency sentence
under one closed hypothesis, $\LF$, that \rga neither proves nor
refutes, and the stronger system \oga proves $\LF$ and with $\LF$ the
sentence (\S\ref{sec:con-lv}).  A last result covers every proof under
one hypothesis about the proof code, the \emph{residual}; that result locates
the obstruction to $\ConR$ (\S\ref{sec:con-res}), and no system sound
for the $\omega$-semantics can run the fold to remove that obstruction
(\S\ref{sec:con-ext}).

The classes do not form a chain.  Among accepted codes, \rga proves
these containments, and no others are claimed:
\begin{itemize}
\item every code that the $\forall$-free checker accepts lies in
  $\Frag$, in $\LvC$ and in each per-level class $\LvC_N$
  (Theorems~\ref{thm:con-qf-in-frag} and~\ref{thm:con-lv-contain});
  the containment in $\Frag$ is strict (Theorem~\ref{thm:con-witness});
\item an accepted code in $\Frag$ lies in $\LvC$ when it passes one
  more decided scan, the head-step scan, and every code in $\LvC$
  passes that scan (Theorem~\ref{thm:con-lv-contain}); otherwise
  $\Frag$ is incomparable with $\LvC$ and with each $\LvC_N$
  (\S\ref{sec:con-break});
\item at closed codes every $\LvC_N$ lies in $\LvC$
  (Theorem~\ref{thm:con-lv-closed}).
\end{itemize}

Every theorem of the family is $\ConR$ restricted in one of two ways.
Either the universal ranges only over a decided class of proof codes,
or the sentence is the conclusion of an entailment, derived under
hypotheses rather than in the empty context.

One method proves the theorems of the family: an induction inside \rga
over the positions of a coded proof, which this paper calls the
\emph{fold}.  The fold shows, position by position, that every accepted
proof has a true conclusion, so that no accepted proof concludes $\z =
\tofn 1$.  The $\forall$-free theorem runs the fold with the truth
predicate $\Tru$ of \S\ref{sec:inside-truth}.  The other theorems run
one shared fold with a certificate-free notion of truth, and differ
only in what that fold is given at two \emph{interfaces}: deciding the
case hypotheses on a path, and moving a value between two terms with
the same value (\S\ref{sec:con-frag}).  \S\ref{sec:con-break} explains
why only these interfaces are at stake.

The conditional results are entailments, derivations under hypotheses.
An entailment becomes an implication exactly when its hypothesis is
provably decided, or has the form $\exists s.\,D(s)$ with $D$ decided
(Theorem~\ref{thm:in-dt}).  The hypotheses here are of neither kind,
so these results are stated as entailments; whether they hold as
implications is open, and for the level theorem the implication holds
exactly when the conclusion holds outright
(Theorem~\ref{thm:con-lf-indep}).  As entailments, the conditional results turn any
derivation of the hypothesis into a derivation of the consistency
statement.

\subsection{The shape of every consistency sentence}
\label{sec:con-shape}

Every theorem of this section contains the same short argument, the
\emph{decided-matrix split}.  Fix a fresh variable $p$ with $p \jnat$.
The checker's verdict $\Chk(p, \bquo{\z = \tofn 1})$ is decided at
the abstract $p$ (Theorem~\ref{thm:in-chk-dec}), and so is the verdict
of the $\forall$-free checker $\ChkQf$ (\S\ref{sec:con-qf}), which
the $\forall$-free theorem uses in place of $\Chk$
(Theorem~\ref{thm:c-chkqf-dec}).  So \rga may split on the verdict.
\begin{itemize}
\item If the checker rejects, the goal $\neg\Chk(p, \bquo{\z = \tofn
  1})$ is the case's hypothesis.
\item If the checker accepts, the fold runs from the root of $p$ and
  concludes that $\z = \tofn 1$ is true.  \rga refutes that truth
  internally, so this case derives anything, in particular the goal.
\end{itemize}
Disjunction elimination closes the split.  For the internal theorems,
and for those under a closed hypothesis ($\LF$, $\LF_{\le N}$ of
\S\ref{sec:con-lv}), $\forall$-introduction at $p$ then gives the
universal.  The theorems whose hypotheses mention the proof code
(the reach-indexed family formula $\DynFam$ of \S\ref{sec:con-lv} and
the residual $\Res$ of \S\ref{sec:con-res}) stop at the free
variable.  Where the sentence
restricts $p$ to a decided class, such as $\Frag(p)$, the implication
$\Frag(p) \limp \cdots$ is introduced by discharging the decided class
hypothesis.  No proof in this section discharges a hypothesis that is
not provably decided.  The hypotheses of the conditional theorems sit
in the derivation's context and are only eliminated.  The hypotheses
about the proof code ($\DynFam$, $\Res$) are not shown decided.  The
closed ones, $\LF$, $\LF_{\forall}$ and $\LF_{\le N}$, are independent
of \rga: \rga neither proves nor refutes them, so \rga does not prove
these hypotheses decided (Theorems~\ref{thm:con-lf-indep} and~\ref{thm:con-lvN}(3)).

\begin{table}[!hb]
\centering
\small
\begin{tabular}{@{}>{\raggedright\arraybackslash}p{0.17\linewidth}>{\raggedright\arraybackslash}p{0.07\linewidth}>{\raggedright\arraybackslash}p{0.19\linewidth}>{\raggedright\arraybackslash}p{0.13\linewidth}>{\raggedright\arraybackslash}p{0.31\linewidth}@{}}
\hline
theorem & system & hypotheses (context) & class of $p$ & status \\
\hline
$\forall$-free, Thm.~\ref{thm:con-qf} & \rga & none & ($\ChkQf$ includes the scan) &
  unconditional, internal \\
fragment, Thm.~\ref{thm:con-frag} & \rga & none & $\Frag(p)$ &
  unconditional, internal \\
level, Thm.~\ref{thm:con-lv} & \rga & $\LF$ &
  $\LvC(p)$ & $\LF$ neither proved nor refuted
  (Thm.~\ref{thm:con-lf-indep}); one sentence for all levels \\
level in \oga, Thm.~\ref{thm:con-oga} & \oga & none & $\LvC(p)$ &
  unconditional in \oga, which proves $\LF$ (Thm.~\ref{thm:con-olf}) \\
per level, Thm.~\ref{thm:con-lvN} & \rga; \oga & $\LF_{\le N}$; none &
  $\LvC_N(p)$ & as the two rows above, one result for each $N$ \\
reach-guarded, Thm.~\ref{thm:con-e1} & \rga & $\vz\jnat$, $\DynFam(\vz)
  \wedge \LF$ (per level: $\wedge \LF_{\le N}$) & $\Stat(\vz)$
  ($\Stat_N(\vz)$) & $\DynFam$ speaks about the proof code; contains
  the level and per-level theorems \\
residual, Thm.~\ref{thm:con-res} & \rga & $\vz\jnat$, $\Res(\vz)$ &
  every code & locates the obstruction: $\Res$ is not derivable
  (Cor.~\ref{thm:con-res-underivable}), nor from the classes
  (Thms.~\ref{thm:con-res-not-frag}, \ref{thm:con-res-lf}), but holds
  at swap-free codes in $\Frag$ (Thm.~\ref{thm:con-res-frag}) \\
\hline
\end{tabular}
\caption{Restricted consistency theorems for \rga's proofs.  Each
  theorem concludes that $p$ is not accepted as a proof of $\z = \tofn 1$, by
  $\ChkQf$ in the first row and by the full checker $\Chk$ of $\ConR$
  elsewhere.  Columns: the hypotheses sit in the context and are only
  eliminated; the class of $p$ is a decided antecedent, introduced as an
  implication; rows without $\vz\jnat$ quantify over $p$ with $\cA$,
  the other rows are about the free variable $\vz$.  $\LvC$, $\LvC_N$,
  $\LF$, $\Stat$ and $\DynFam$ are defined in \S\ref{sec:con-lv}, and
  $\Res$ in \S\ref{sec:con-res}.  The containments between the classes
  are listed at the start of this section.}
\label{tab:con-tower}
\end{table}

\S\ref{sec:pillars-con} used this shape to show that \rga does not
prove the single-sentence form of soundness
(Theorem~\ref{thm:pl-no-uniform}), and [\citen{goedel-paper}, Section ``G\"odel's theorems for \rga''] used
the same shape in the central step of G\"odel's second theorem
([\citen{goedel-paper}, theorem~\textsf{goedelG\_of\_expl\_d3g\_ctx}]).  Here the accepting
case has real content: the fold.

\emph{Classical versus \rga.}  Classically one proves ``if $p$ proves
$\z = \tofn 1$ then false'' by assuming the antecedent and deriving a
contradiction, whatever kind of formula the antecedent is.  \rga can
introduce such an implication only because the checker's verdict is a
decided 0/1-valued term.  The same restriction shapes the fold's
induction motive, below.

\subsection{The \texorpdfstring{$\forall$}{A}-free fragment}
\label{sec:con-qf}

A term is \emph{$\forall$-free} when the quantifier combinator $\cA$
occurs nowhere in it; this is \rga's form of ``quantifier-free''.  A
proof is $\forall$-free when every hypothesis and every conclusion of
every judgment in it is $\forall$-free.  A compiled test
$\qfs(c)$ checks a code for $\forall$-freeness, with $0$ meaning
``pass''.  The compiled whole-proof scan $\QfScan(p)$ applies the test
to every judgment of the list coded by $p$.

\begin{defi}[\textsf{iChkQf}]
\label{def:con-chkqf}
The \emph{$\forall$-free checker} $\ChkQf(p, c)$ is the conjunction
of the four decided tests of Definition~\ref{def:in-chk} and the scan
$\QfScan(p)$.  In every context, \rga derives $\ChkQf(p, c)$ exactly
when it derives $\Chk(p, c) \wedge \QfScan(p)$.
\end{defi}

\begin{thm}[\textsf{iChkQf\_faithful}]
\label{thm:con-qf-faithful}
For every proof $P$ and term $c$: $\emptyset \vdash
\ChkQf(\bquo{P}, \bquo{c})$ if and only if $P \proves \emptyset
\vdash c$ and $P$ is $\forall$-free.
\end{thm}

\noindent
The theorem is schematic, like Theorem~\ref{thm:in-chk}(1), and the
checker is genuine at every number (Theorem~\ref{thm:c-chkqf-genuine}).

\begin{thm}[\textsf{rga\_con\_qf\_int}, \textsf{rga\_sound\_qf\_int2}]
\label{thm:con-qf}
\label{thm:con-qf-sound}
\leavevmode
\begin{enumerate}
\item $\emptyset \vdash \forall p.\ \neg\ChkQf(p, \bquo{\z = \tofn 1})$.
\item $\emptyset \vdash \forall p.\ \forall c.\ \ChkQf(p, c) \limp
  \Tru(\Subst(\z, c))$.
\end{enumerate}
\end{thm}

\noindent
Both parts are internal: one sentence each, derived in the empty
context with no hypothesis.  Part~2 says that every accepted
$\forall$-free proof has a true conclusion, with truth in the sense of
the decider-based $\Tru$ of \S\ref{sec:inside-truth}.  $\Subst(\z,
c)$ is the image of $c$ under the empty substitution, whose code is
$\z$.  The implication of part~2 is introduced over the decided verdict
$\ChkQf(p, c)$.

\emph{Proof idea: the decider fold.}  The induction runs over the
entries of the proof list coded by $p$, as a course-of-values
induction inside \rga.  The motive at an entry $x$ says: if the entry
is accepted and passes the scan, then for every well-formed
$\forall$-free substitution $z$ and every fuel $s$, if the decider
affirms at fuel $s$ every hypothesis of the entry's judgment under
$z$, then the entry's conclusion under $z$ is $\Tru$.  This motive is the
companion paper's soundness invariant for a proof entry (every
substitution that makes the hypotheses true makes the conclusion
true), with grounded truth replaced by the decider's verdict.  The
fuel quantifier sits in front of the hypothesis test, and that
placement is essential.  ``The hypotheses are $\Tru$'' is not
decided in general, since $\Tru$ is an existential over fuels: \rga
does not prove $\Tru(S_{\Tru}) \jbool$ for the code $S_{\Tru}$ of the
Liar $L_{\Tru} = \neg\Tru(S_{\Tru})$ (Theorem~\ref{thm:pl-liar-indep}).  But
``the decider affirms the hypotheses at fuel $s$'' is a decided 0/1
test, so every implication in the motive is introduced legitimately.

The induction step dispatches on the checker's node table, one arm
per rule, 41 arms in all (the two directions of reflexivity share a
row).  Most arms need one law of the compiled decider at the node's
codes.  Eight arms, among them cut and one of the disjunction
eliminations, also need monotonicity of the decider in the fuel at
$\forall$-free codes (Table~\ref{tab:c-qf-arms} lists them).  The
arms for conversion, for the recursor's successor step, for equality
substitution, for induction and for substitution need further internal
laws about the decider at $\forall$-free codes: invariance under head
steps, a value-congruence (``swap'') engine, laws for substitution and
for the freeness test, and an induction on a stage numeral inside
\rga.  The four quantifier rules cannot occur,
since their premises or conclusions contain $\cA$ and the scan rejects
them.  At the root the substitution is empty and the conclusion has
no hypotheses, so the motive yields $\Tru(\Subst(\z, \bquo{\z = \tofn
1}))$.  The empty substitution leaves the code unchanged, and
Theorem~\ref{thm:in-never} refutes $\Tru(\bquo{\z = \tofn 1})$
inside \rga, which closes the accepting case of the decided-matrix
split.
Appendix~\ref{app:fold} groups the 41 arms and exhibits two of them.

\emph{Why $\forall$-free makes the fold go through.}  At a
$\forall$-free code the decider never enters its quantifier cell, the
bounded certificate search of \S\ref{sec:inside-dec}.  Everywhere else
the decider is exact, and its laws (monotonicity, head-step
invariance) are provable inside \rga.  At a quantifier node, ``the
decider affirms the universal'' means ``a certificate was found''.
Turning that finding into the truth of the universal's instances is the
soundness of the certificate, which is the companion paper's soundness
theorem itself, run in HOL.  No internal sentence states that step;
the step's whole-proof form at $\z = \tofn 1$ is not derivable
(Theorem~\ref{thm:pl-no-uniform}).

\paragraph{What the fragment covers.}
Every theorem of \rga that is itself $\forall$-free has a
$\forall$-free proof, which $\ChkQf$ accepts:

\begin{thm}[\textsf{rga\_afree\_normal\_chk}]
\label{thm:con-normal}
If $P \proves \emptyset \vdash c$ and $c$ is $\forall$-free, then
some $\forall$-free proof $Q$ satisfies $Q \proves \emptyset \vdash
c$ and $\emptyset \vdash \ChkQf(\bquo{Q}, \bquo{c})$.
\end{thm}

\noindent
Theorem~\ref{thm:con-normal} is a meta theorem, and its proof is
semantic, not a syntactic normalization: soundness gives the
conclusion's evaluation, and the completeness of the $\forall$-free
fragment rebuilds a $\forall$-free derivation from that evaluation.  Combined with
Theorem~\ref{thm:con-qf}, Theorem~\ref{thm:con-normal} shows in HOL that no proof of $\z =
\tofn 1$ exists at all, which soundness already gives.  \rga cannot
make this combination internally.  \rga does not prove the internal
form of Theorem~\ref{thm:con-normal} at $\z = \tofn 1$, $\forall p.\
\Chk(p, \bquo{\z = \tofn 1}) \limp \neg\forall q.\ \neg\ChkQf(q,
\bquo{\z = \tofn 1})$ ([\citen{goedel-paper}, Appendix ``G\"odel II in full'']): with
Theorem~\ref{thm:con-qf} that internal form would give $\ConR$, which G\"odel's second
theorem excludes.  Classically the same obstruction appears in Gentzen's
analysis of Peano arithmetic~\cite{gentzen36widerspruchsfreiheit}:
Peano arithmetic does not prove that Gentzen's reduction of its proofs
terminates, which needs induction up to $\varepsilon_0$, since Peano
arithmetic would then prove its own consistency.

\subsection{The fragment theorem}
\label{sec:con-frag}

The second internal theorem concerns a class of proofs that admits
quantifiers except at four kinds of positions
(Definition~\ref{def:con-frag}).  The fold for this class needs a
different notion of truth.

\paragraph{The value record.}
The decider's $\Tru$ handles a universal through its certificate,
which is exactly what the fold cannot use (\S\ref{sec:con-qf}).  The
\emph{value record} drops certificates.  The value record is an internal formula
$\Rcd_v(c)$, ``the term with code $c$ has the value $v$ at some
fuel'', defined inside \rga by self-reference over a record body with
one row per clause of the semantics (Table~\ref{tab:in-sem}).  The
rows for $\z$, $\suc$, ${=}$, $\neg$, $\lor$, head steps and the
recursor mirror the semantics.  The two quantifier rows read the
numeral instances directly: $\cA\cdot u$ has value $\ctrue$ when every
numeral instance $u\cdot\tofn n$ has value $\ctrue$, and value $\cfalse$
when some instance has value $\cfalse$.  These two \emph{numeral families}
of $u$ are the formulas $\forall n.\ \Rcd_{1}(u\cdot\bquo{\tofn n})$
and $\exists n.\ \Rcd_{0}(u\cdot\bquo{\tofn n})$, where the instance
code is built inside \rga.  The fuel-bounded core $\Rcd^s_v(c)$
records the value $v$ at fuel $s$; the families inside the core are not
bounded by $s$.

\begin{defi}[\textsf{iRec}, \textsf{iRecAt}]
\label{def:con-rcd}
$\Rcd_v(c)$ is the value record at some fuel, and $\Rcd^s_v(c)$
its core at fuel $s$, with the rows described above.
\end{defi}

\noindent
The value record is an $\omega$-style truth.  Its HOL counterpart
contains grounded truth at well-formed terms: every grounded value of a well-formed
term is a recorded value.  At a universal the record asks for every instance,
not for a certificate, so the record can hold values that grounded truth
does not have.  At closed $\forall$-free codes of well-formed terms the
record and grounded truth agree (Theorem~\ref{thm:c-orec}).  At
quantified codes the two notions come apart: the record's HOL
counterpart gives $\ConR$ the value $\ctrue$, while $\ConR$ has no
grounded value and \rga does not derive $\Tru(\bquo{\ConR})$
(\S\ref{sec:pillars-con}, [\citen{goedel-paper}, Appendix ``G\"odel II in full'']).  The
consistency theorems below do not depend on the choice of truth
notion, since both notions refute the false root; the soundness
theorems do, and each soundness theorem names its truth notion.

\paragraph{The fragment.}
Two kinds of proof steps matter for the fold.  Three rules are
\emph{case-dependent}: a premise carries an added hypothesis that no
other premise establishes, a case assumption.  The three case-dependent rules
are the elimination of a negated universal (the hypothesis is the
witness instance), induction (the induction hypothesis at the stage),
and disjunction elimination (the branch).  Five rules are \emph{swap}
rules: a premise and the conclusion are instances of one formula at
two terms that have the same value.  The five swap rules are the recursor's successor
step, universal elimination, introduction of a negated universal,
induction, and equality substitution.  The formula whose instances a
swap rule compares is its \emph{critical formula}.

In one sentence, the fragment admits quantifiers everywhere except at
the positions where the fold would have to decide a family of
infinitely many numeral instances.  Precisely:

\begin{defi}[\textsf{iFrag}]
\label{def:con-frag}
$\Frag(p)$ is a decided scan of the proof code $p$.  The scan requires that
at every substitution node the substituted term that the checker reads
(the result of the checker's witness search) is $\forall$-free;
that every hypothesis added by a case-dependent rule is
$\forall$-free; that at every swap node the critical formula is
$\forall$-free, well-formed and the code of a term; and that at three
of the swap rules the swapped term is $\forall$-free.  Nothing else is
constrained.
\end{defi}

\noindent
Conclusions, cut formulas and every other judgment of a proof in
$\Frag$ may contain quantifiers.

\begin{thm}[\textsf{rga\_con\_frag\_int}, \textsf{rga\_sound\_frag\_int}]
\label{thm:con-frag}
\label{thm:con-frag-sound}
\leavevmode
\begin{enumerate}
\item $\emptyset \vdash \forall p.\ \Frag(p) \limp \neg\Chk(p, \bquo{\z =
  \tofn 1})$.
\item $\emptyset \vdash \forall p.\ \Frag(p) \limp \forall c.\ \Chk(p, c)
  \limp \Rcd_1(\Subst(\z, c))$.
\end{enumerate}
\end{thm}

\noindent
Both parts are internal, in the empty context, with no hypothesis.
Both parts use the full tree checker $\Chk$, the checker of $\ConR$.  The
only implications introduced are over decided antecedents: the scan
$\Frag(p)$ and the verdict $\Chk(p, c)$.  Part~2 covers every
conclusion $c$, quantified conclusions included, with truth in the
sense of the value record.

\paragraph{The decided-record lemma.}
The fragment is unconditional because of one internal lemma, the
decided-record lemma (Lemma~Q in the development).

\begin{thm}[the decided-record lemma, Lemma Q; \textsf{lemQ\_int}]
\label{thm:con-lemq}
Let $\Gamma$ be a well-formed context.  If $\Gamma$ derives $s\jnat$,
$c\jnat$, $k\jnat$ and $\qfs(c) = \z$, then $\Gamma \vdash
\Rcd^s_k(c) \jbool$.
\end{thm}

\noindent
At every $\forall$-free code, even an abstract one, \rga proves that
the fuel-bounded record is decided.  The proof is one internal
induction on the fuel.  A $\forall$-free code never reaches a
quantifier row, so the core consults no numeral family.

\paragraph{Example: $\forall$-free is not the same as decided.}
The Liar $L$ of the companion paper, the fixed point of
$\lambda\neg\vz$, unfolds to its own negation.  $L$ is $\forall$-free
and has no value.  Lemma~Q decides the core of $L$'s record at every fuel.  A separate
internal induction on the fuel refutes the core at every fuel, and
from that refutation \rga refutes $\Rcd_v(\bquo{L})$ for every value $v$
(Appendix~\ref{app:fold-frag}).  So a proof whose case hypothesis is the
Liar still has a decided guard, in the sense below, and can belong to
$\Frag$.  What Lemma~Q decides is the bounded record, not the term's
value.

\emph{Proof idea of Theorem~\ref{thm:con-frag}: the record fold.}  A
proof list cites earlier entries, and one entry can be cited along
several paths from the root.  The fold walks the tree that this
citation structure unfolds into, by a course-of-values induction
inside \rga over the nodes' post-order ranks.  The walk is per path
because the eigenvariables of $\forall$-introduction are instantiated
along a path, so one list entry carries a different substitution, the
\emph{path environment}, on each path.  The motive at a rank $r$
says: for all instantiations of the path's eigenvariables and all fuel
bounds,
\[
  \neg\,\mathit{Guard}_r \;\lor\;
  \Rcd_1(\text{the conclusion at $r$ under its path environment}) .
\]
The \emph{guard} $\mathit{Guard}_r$ is the conjunction, over the
case-dependent hypotheses added along the path to $r$, of their
fuel-bounded value-$\ctrue$ records.  The disjunction is how a
hypothesis enters the motive without an implication over a possibly
undecided antecedent.

The step has one arm per row of the node table.  Most arms hold for
every class of proofs at once (Appendix~\ref{app:fold}).  Seven arms
need help at two \emph{interfaces}.  The hypothesis-leaf arm and the
$\forall$-introduction arm need the \emph{guard interface}: the guard
is decided at every point of the path, so that the arm can split on
the guard.  The five swap arms need the \emph{swap interface}: a fuel-bounded
record of the critical formula's instance at one term gives the record
of the formula's instance at any other term with the same value.  At the root the path
environment is empty and the guard is vacuous, so the root's
conclusion has value $\ctrue$; at $\bquo{\z = \tofn 1}$ no record has
that value (Theorem~\ref{thm:c-false-root}).  For $\Frag$, Lemma~Q
decides every guard item at its $\forall$-free instance, and a swap
engine proved inside \rga (Theorem~\ref{thm:c-swap}) supplies the swap
at a $\forall$-free critical formula.  No numeral family is consulted
anywhere, so the fold needs nothing beyond $p\jnat$ and $\Frag(p)$.
(A third, minor interface asks that the swapped term at three of the
swap rules be $\forall$-free under its path environment.  The class
scans supply this interface, and the guarded swap of \S\ref{sec:con-res}
does not need it.)

\emph{Classical versus \rga.}  Classically a theory proves the
consistency of a fragment by defining a partial truth predicate for
the fragment's formulas and showing by induction over proofs that
every line is true.  Peano arithmetic, for example, proves the
consistency of each fragment $\mathrm{I}\Sigma_n$ by a partial truth
(satisfaction) predicate of bounded
complexity~\cite{hajek93metamathematics}.  \rga follows the same plan,
with two differences.  First, the fold can use its truth notion inside
\rga only where the relevant tests are decided: the decider's verdict
at $\forall$-free codes ($\Tru$), and the fuel-bounded record where
the numeral families are decided ($\Rcd_1$).  Both notions are
adequate in HOL more widely: $\Tru$ at every closed code
(Theorem~\ref{thm:pl-adeq}), the record at closed $\forall$-free
codes of well-formed terms.  Second, the motive's implications must be
introduced by discharging decided hypotheses, so ``the hypotheses are true'' must
become a decided test, the fuel-prefixed hypothesis test of the
$\forall$-free fold and the guarded disjunction of the record fold.
What breaks beyond the fragments is not the truth predicate but the
decidedness of the guard and of the swap (\S\ref{sec:con-break}).
The fragments restrict where quantifiers may occur in a proof, not
which axioms the proof may use.  The comparison with
$\mathrm{I}\Sigma_n$ concerns the method only, not the strength of the
classes (below).

\paragraph{The fragment and the $\forall$-free class.}
The two internal theorems use different checkers, $\ChkQf$ and $\Chk$,
and different scans.  Every code that the $\forall$-free checker
accepts lies in the fragment nevertheless, and \rga proves this containment as one
sentence.

\begin{thm}[\textsf{rga\_qf\_in\_frag\_int}]
\label{thm:con-qf-in-frag}
$\emptyset \vdash \forall p.\ \forall c.\ \ChkQf(p, c) \limp \Frag(p)$.
\end{thm}

\noindent
The theorem is internal, in the empty context with no hypothesis, and
its implication is introduced over the decided verdict $\ChkQf(p, c)$.
The proof reads nothing back from a particular code and is not an
induction over the proof.  $\Frag(p)$ is a conjunction of two scans
over the nodes of the code $p$, and \rga introduces each scan from a
fact about the node at an arbitrary position.  An accepted
$\forall$-free node meets every clause of $\Frag$, because each term a
clause reads comes out of the node's own $\forall$-free judgments; two
internal facts about the substitution coder, each a course-of-values
induction inside \rga, handle the substitution rules
(Appendix~\ref{app:fold-contain}).

With Theorem~\ref{thm:con-frag}, the containment gives the
$\forall$-free consistency theorem a second derivation, inside \rga.

\begin{thm}[\textsf{rga\_con\_qf\_via\_frag\_int}]
\label{thm:con-qf-via-frag}
$\emptyset \vdash \forall p.\ \neg\ChkQf(p, \bquo{\z = \tofn 1})$,
derived inside \rga from the sentences of
Theorems~\ref{thm:con-qf-in-frag} and~\ref{thm:con-frag}.
\end{thm}

\noindent
The statement is that of Theorem~\ref{thm:con-qf}.  Under the decided
hypothesis $\ChkQf(p, \bquo{\z = \tofn 1})$,
Theorem~\ref{thm:con-qf-in-frag} instantiated at $p$ and $\bquo{\z =
\tofn 1}$ gives $\Frag(p)$, Theorem~\ref{thm:con-frag} at $p$ gives
$\neg\Chk(p, \bquo{\z = \tofn 1})$, and the full checker accepts every
code that the $\forall$-free checker accepts.  The derivation uses only first-order
rules and the two sentences, no fact about codes proved in HOL.  So
Theorem~\ref{thm:con-qf} is a consequence of
Theorem~\ref{thm:con-frag} inside \rga.

The containment is strict.  Let $P_{\mathrm w}$ be the following
three-line proof of $\forall x.\ x = x$.  One line is the hypothesis
$\vz\jnat$, that is $\vz = \vz$, in the context $\vz\jnat$.  One
conversion step gives $(\lambda x.\ x = x)\cdot\vz$ in the same
context.  $\forall$-introduction discharges $\vz\jnat$ and concludes
$\cA\cdot(\lambda x.\ x = x)$, which is $\forall x.\ x = x$, in the
empty context.

\begin{thm}[\textsf{fwP\_in\_frag}, \textsf{fwP\_not\_qf}]
\label{thm:con-witness}
\leavevmode
\begin{enumerate}
\item $\emptyset \vdash \Chk(\bquo{P_{\mathrm w}}, \bquo{\forall x.\ x =
  x})$.
\item $\emptyset \vdash \Frag(\bquo{P_{\mathrm w}})$.
\item $P_{\mathrm w}$ is not $\forall$-free, $\emptyset \vdash
  \neg\QfScan(\bquo{P_{\mathrm w}})$, and $\emptyset \nvdash
  \ChkQf(\bquo{P_{\mathrm w}}, c)$ for every term $c$.
\end{enumerate}
\end{thm}

\noindent
$P_{\mathrm w}$ also lies in $\LvC$ and in every per-level class
(Theorem~\ref{thm:c-witness}).  So, with
Theorem~\ref{thm:con-qf-in-frag}, the class of
Theorem~\ref{thm:con-frag} strictly contains the class of
Theorem~\ref{thm:con-qf}.  $P_{\mathrm w}$ belongs to these classes
because none of the classes' clauses applies at the three nodes of
$P_{\mathrm w}$: the quantifier of $P_{\mathrm w}$ sits in judgments that the classes leave unconstrained.  A
second witness, $P_3$, exercises the swap clause.  With $g = \lambda x.\
(\lambda y.\ y = y)\cdot x$, $P_3$ derives $\forall x.\ g\cdot x$ as
$P_{\mathrm w}$ derives $\forall x.\ x = x$ (with one more conversion
step), and concludes $g\cdot\z$ by universal elimination at $\z$.  The
checker reads the root of $P_3$ as a universal elimination, and $P_3$ lies in
$\Frag$, in $\LvC$ and in every $\LvC_N$
(Appendix~\ref{app:fold-contain}).

\paragraph{Strength, and the relation to \pga.}
The proof-theoretic strength of the two internal classes is open.
These results are \rga's counterparts of the classical fact that Peano
arithmetic proves the consistency of each of its fragments
$\mathrm{I}\Sigma_n$~\cite{hajek93metamathematics}; these results do not escape
G\"odel's second theorem.  Both classes admit induction only on
$\forall$-free formulas: the $\forall$-free class trivially, and
$\Frag$ because the induction rule's critical formula and its case
hypothesis must be $\forall$-free.  A $\forall$-free formula is not
quantifier-free in the classical sense: through untyped $\lambda$ such a
formula can run an unbounded search, so its truth can be a $\Sigma_1$
condition, and it can lack a value.  Candidate comparisons therefore
range from quantifier-free or $\Sigma_1$ induction (primitive recursive
arithmetic, and $\mathrm{I}\Sigma_1$, which is $\Pi_2$-conservative
over primitive recursive arithmetic) to quantifier-free formulations of System~T, which the classes
resemble because they contain $\lambda$ and System~T's recursor, and
which is far stronger.  Full $\mathrm{I}\Sigma_n$ for $n \ge 2$ is not
the right comparison, since both classes restrict induction to
$\forall$-free formulas.  Both classes are in the family of
\pga~\cite{ford25have}, propositional grounded arithmetic, the
quantifier-free member of the grounded-deduction family: they share
\pga's propositional, equational and inductive core, with $\lambda$
and System~T's recursor in place of \pga's coded application and
combinators.  A mechanized proof that \rga proves \pga itself
consistent, which nothing in G\"odel's theorem forbids, is future
work.

\paragraph{The fragment as a search that never halts.}
\rga also proves the content of Theorem~\ref{thm:con-frag} as a
divergence.  Let $\mathit{fragSearch}$ be the term that searches for a
code $x$ in $\Frag$ that $\Chk$ accepts as a proof of $\z = \tofn 1$.
\rga proves $\neg\Rcd_v(\bquo{\mathit{fragSearch}})$ for every $v$: the
search has no recorded value.  For the unrestricted search over all
codes, \rga does not refute the value-$\ctrue$ record, since a
refutation would yield $\emptyset \vdash \ConR$
(Theorem~\ref{thm:c-search}).

\subsection{Where the fold stops}
\label{sec:con-break}

For full \rga the guard interface of the record fold fails.  Deciding
a guard item means deciding a fuel-bounded record.  At a quantified
case hypothesis, the record's quantifier rows consult the numeral
families of the hypothesis's quantified subterms, and deciding a
family $\forall n.\ \Rcd_1(u\cdot\bquo{\tofn n})$ is an $\omega$-type
question.  If \rga decided such families uniformly, the complement of
the halting set would be recursively enumerable
(Theorem~\ref{thm:con-not-lf} below).  For the guard this is a
theorem: at proofs whose case hypothesis says that a program does not
halt, \rga does not derive that the guard is decided
(\S\ref{sec:con-res}).  Everything else in the fold holds for every
proof.

\begin{rem}
\label{rem:con-swap}
The swap at a quantified critical formula would need the same
families, so the swap interface plausibly fails there as the guard
interface does.  This reasoning is an argument, not a theorem; whether the swap
interface fails there is an open question.  The swap-half failure that is proved
(\S\ref{sec:con-res}) sits at a $\forall$-free critical formula and has
a different cause.
\end{rem}

The remaining theorems therefore supply the two interfaces in
different ways:
\begin{itemize}
\item the \emph{level theorem} assumes the families at certified
  positions, as one closed formula $\LF$, which \oga proves
  (\S\ref{sec:con-lv});
\item the \emph{reach-guarded theorem} assumes the families for the
  guard only where the evaluation of a hypothesis reaches them, and
  takes the swap from $\LF$ (\S\ref{sec:con-lv});
\item the \emph{residual theorem} assumes the two interfaces
  themselves, as one formula about the proof code
  (\S\ref{sec:con-res}).
\end{itemize}
The family is one fold at several assumptions, not a chain of
theorems in which each implies the next; the containments among the
classes are those listed at the start of this section.  $\Frag$ is
incomparable with $\LvC$ and with the per-level classes: \rga refutes
membership in $\LvC$ (or in every $\LvC_N$) at an accepted code in
$\Frag$, and membership in $\Frag$ at an accepted code in $\LvC$ (or in
every $\LvC_N$) (Theorem~\ref{thm:e-frag-lv}).

On accepted codes in $\Frag$, $\LvC$ adds one decided clause.  Write
$\mathit{HS}(p)$ for the \emph{head-step scan}: at every node of $p$
that the checker reads as a universal elimination or a
negated-universal introduction, the compiled head-step reader finds a
step of the node's critical formula.  Like the scans of $\Frag$,
$\mathit{HS}$ is a bounded scan over the nodes of $p$, and it is
decided.  Let $P_{\mathrm R}$ be the accepted proof in $\Frag$ whose
root eliminates $\forall x.\ \cR\cdot g\cdot h\cdot x$, a universal
over a recursor matrix (Theorem~\ref{thm:e-frag-lv}).

\begin{thm}[\textsf{rga\_qf\_in\_lvi\_int}, \textsf{rga\_frag\_in\_lvi\_int}]
\label{thm:con-lv-contain}
\leavevmode
\begin{enumerate}
\item $\emptyset \vdash \forall p.\ \forall c.\ \ChkQf(p, c) \limp
  \LvC(p)$.
\item $\emptyset \vdash \forall p.\ \forall c.\ \Chk(p, c) \wedge
  \Frag(p) \wedge \mathit{HS}(p) \limp \LvC(p)$.
\item If $\Gamma$ is well-formed, $\Gamma \vdash p\jnat$ and $\Gamma
  \vdash \LvC(p)$, then $\Gamma \vdash \mathit{HS}(p)$.
\item $\emptyset \vdash \neg\mathit{HS}(\bquo{P_{\mathrm R}})$.
\item $\emptyset \vdash \neg\LvC(\bquo{P_{\mathrm R}})$, and
  $\emptyset \nvdash \LvC(\bquo{P_{\mathrm R}})$.
\end{enumerate}
Parts~1 to~3 hold for every per-level class $\LvC_N$ in place of
$\LvC$.
\end{thm}

\noindent
Parts~1 and~2 are internal, in the empty context with no hypothesis,
and their implications are introduced over decided antecedents.  The
proofs of parts~1 and~2 work node by node at the abstract code $p$, as for
Theorem~\ref{thm:con-qf-in-frag}, and read nothing back.  A
$\forall$-free code is certified, at depth~1.  At a universal
elimination or a negated-universal introduction, $\LvC$ asks that the
head reduct of the critical formula be static, certified, and a
well-formed term code; the proof of part~2 rests on one internal law,
that the head reduct of a well-formed term code is again a well-formed
term code.  By parts~2 and~3, on accepted codes in $\Frag$ the
head-step scan is exactly what $\LvC$ adds.  The scan is needed even
at the root: $P_{\mathrm R}$ fails it, because $\cR\cdot g\cdot
h\cdot\vz$ has no head step, so $P_{\mathrm R}$ lies outside $\LvC$
and every $\LvC_N$.  Part~5 is proved at the closed code, where the
disjunction property turns the non-derivability into a refutation,
since $\LvC$ is decided at numerals.  Appendix~\ref{app:oga-frag} gives
the proofs.

\subsection{Family hypotheses, and the discharge in OGA}
\label{sec:con-lv}

\paragraph{Certified codes.}
One compiled course-of-values machine computes the \emph{certified
depth} $\Dep(c)$ of a code $c$.  The depth is $0$ if $c$ is not
certified; otherwise the depth is the length of the longest chain of
quantifiers in $c$ nested through head steps, and $1$ if $c$ has none.
A $\forall$-free code has depth~$1$.  A spine node, built by $\suc$,
$\neg$, $\lor$ or ${=}$, has the largest depth of its arguments if none
of them has depth~$0$.  A quantifier application $\cA\cdot m$ has
depth~$1$ if its matrix $m$ is $\forall$-free, and depth $1 + \Dep(R)$
if $m$ is a $\lambda$ whose instance at $\vz$ head-steps to a static,
well-formed term code $R$ of nonzero depth.  Every other code has
depth~$0$.  A code $c$ is \emph{certified} when $\Dep(c) \ne \z$ and
$c$ is a well-formed term code, and \emph{certified at level $n$}, for
a term $n$, when moreover $\Dep(c) \le n + 1$.  Both tests are decided.
The machine has no level parameter; Appendix~\ref{app:oga-level} gives
the machine's definition and laws.

\begin{defi}[\textsf{iLvClassI}]
\label{def:con-lv}
The class $\LvC(p)$ is the decided class of proof codes whose
case-dependent hypotheses are certified, whose substituted terms are
$\forall$-free, and whose swap nodes satisfy three conditions: the swap
target is static and certified; at universal elimination and at
introduction of a negated universal, the critical formula has a head
step; and at three of the swap rules the swapped term is
$\forall$-free.
\end{defi}

\noindent
The swap target is the critical formula, or its head reduct at the two
rules that need a head step; Appendix~\ref{app:oga} gives the scans in
full.  No level appears in the class: each code is certified at its
own depth.  The witness $P_{\mathrm w}$ of \S\ref{sec:con-frag} lies in
the class, since no clause of the class applies at the nodes of $P_{\mathrm w}$.

\begin{defi}[\textsf{iLFI}, \textsf{iLFA}]
\label{def:con-lf}
$\LF$ is the closed sentence: at every certified quantifier
application code, under every well-formed $\forall$-free substitution,
both numeral families of the instance's matrix are decided.  For a
term $n$, $\LF(n)$ is the same sentence at the codes certified at
level $n$, and $\LF_{\forall} \ldef \forall n.\ \LF(n)$.
\end{defi}

\noindent
$\LF$ is closed: it mentions no proof code and no level.  \rga derives
each of $\LF$ and $\LF_{\forall}$ from the other, in every
well-formed context, so $\LF \eqv \LF_{\forall}$, since a certified
code is certified at the level of its own depth.

\begin{thm}[\textsf{rga\_con\_lvi\_under\_lf\_all}, \textsf{rga\_con\_lvi\_under\_lf}]
\label{thm:con-lv}
\[
  \LF \vdash \forall p.\ \LvC(p) \limp \neg\Chk(p,
  \bquo{\z = \tofn 1}) ,
\]
and the same with $\LF_{\forall}$ in place of $\LF$.
\end{thm}

\noindent
The theorem is an entailment: one sentence, derived under one closed
hypothesis, which the derivation only eliminates.  No meta-level number
appears.  The implication is introduced over the decided class
$\LvC(p)$.  The fold is the same as for $\Frag$; the families that the
guard and the swap consult at certified positions come from $\LF$.  The
fold gives the theorem at a free variable, $\vz\jnat,\ \LF \vdash
\LvC(\vz) \limp \neg\Chk(\vz, \bquo{\z = \tofn 1})$, and since $\LF$ is
closed, $\forall$-introduction at $\vz$ gives the stated form.  The
swap at a nested quantifier needs the swap at its reduct, one level
down.  \rga supplies the swap at every level by one internal induction on the level
(Appendix~\ref{app:oga-lv}), whose motive quantifies the code, the
environment and the fuel bound and has only decided antecedents.

\begin{thm}[\textsf{rga\_not\_iLFI}]
\label{thm:con-not-lf}
$\emptyset \nvdash \LF$.
\end{thm}

\noindent
The proof is a computability argument, not G\"odel's second theorem.
The quantifier application of the halting family is certified, at
depth~1.  If \rga proved $\LF$, the numeral families that $\LF$ decides there
would make the complement of the halting set recursively enumerable
(Theorem~\ref{thm:e-halting}).

\paragraph{The discharge in OGA.}
\oga, Omega Grounded Arithmetic~\cite{oga-paper}, extends \rga's 42
rules by one rule, ATI, which introduces the decidedness of a
universal: from $\forall x.\ (u\cdot x)\jbool$ conclude $(\cA\cdot
u)\jbool$.  A universal all of whose instances are decided is
decided.  ATI is sound for an $\omega$-semantics, and every \rga
derivation is an \oga derivation.  Write $\vdashO$ for derivability in
\oga.  This paper uses only these facts about \oga: its one extra rule
and its $\omega$-soundness (Appendix~\ref{app:oga-oga}).

\begin{thm}[\textsf{o\_iLFI}, \textsf{rga\_oga\_separation\_iLFI}]
\label{thm:con-olf}
$\emptyset \vdashO \LF$ and $\emptyset \vdashO \LF_{\forall}$, while
$\emptyset \nvdash \LF$.
\end{thm}

\noindent
\oga proves $\LF_{\forall}$ by one internal induction on the level.  At
level~$0$ no nested quantifier is certified.  The step from $\LF(x)$, a
hypothesis of the context, to $\LF(x+1)$ decides the families at a
nested quantifier from those of its reduct, which is certified at level
$x$.  \rga's derivation of $\LF$ from $\LF_{\forall}$ then gives $\LF$.
What \oga proves, \rga does not refute.  So $\LF$ is independent of
\rga, and \rga cannot discharge $\LF$ in either direction.
IGA, the extension of \rga by ground induction~\cite{iga-paper}, does
not derive $\LF$ either, by the same computability argument run on the
counterpart of \rga's grounded semantics (the witnessed evaluation)
whose certificates are IGA's proofs.

\begin{thm}[\textsf{iLFI\_independent}, \textsf{iLFI\_imp\_iff}]
\label{thm:con-lf-indep}
If $\emptyset \vdashO c$, then $\emptyset \nvdash \neg c$.  Further,
$\emptyset \nvdash \LF$ and $\emptyset \nvdash \neg\LF$, hence
$\emptyset \nvdash (\LF)\jbool$.  For every well-formed $S$:
$\emptyset \vdash \LF \limp S$ if and only if $\emptyset \vdash S$.
\end{thm}

\noindent
\emph{Proof idea.}  The theorems of \oga, and those of \rga, are true
in the $\omega$-semantics.  A term has at most one $\omega$-value, so
$c$ and $\neg c$ cannot both be true there.  Apply this fact to
Theorem~\ref{thm:con-olf}, and combine with
Theorem~\ref{thm:con-not-lf}.  A proof of $(\LF)\jbool$ would, by the
disjunction property, give a proof of $\LF$ or of its negation.  The
implication $\LF \limp S$ abbreviates the disjunction $\neg\LF \lor S$,
so by the disjunction property a proof of that implication gives a proof of
$\neg\LF$, which is excluded, or of $S$; the converse is
$\lor$-introduction.  The same holds of $\LF_{\forall}$
(Appendix~\ref{app:oga-notlf}).  \qed

\medskip
So the entailment is the only informative conditional form of
Theorem~\ref{thm:con-lv}.  \rga derives the implication from $\LF$ to
the consistency sentence over $\LvC$ exactly when it derives that
sentence outright, and the rule ``if $\vdash \LF$ then $\vdash
\forall p.\ \LvC(p) \limp \neg\Chk(p, \bquo{\z = \tofn 1})$'' is
vacuous, since $\nvdash \LF$.  Whether \rga derives that sentence
outright is an open question.

\begin{thm}[\textsf{rga\_con\_lvi\_oga}]
\label{thm:con-oga}
$\emptyset \vdashO \forall p.\ \LvC(p) \limp \neg\Chk(p, \bquo{\z
= \tofn 1})$.
\end{thm}

\noindent
The proof is one cut of the free-variable form of
Theorem~\ref{thm:con-lv}, embedded in \oga, against $\LF$
(Theorem~\ref{thm:con-olf}), then $\forall$-introduction.  ATI enters
only through $\LF$.  Theorem~\ref{thm:con-oga} is \oga proving
something about \rga's checker: that theorem is neither \oga's
consistency nor \rga's own.  The theorem's soundness companion has the same form, with
$\forall c.\ \Chk(p, c) \limp \Rcd_1(\Subst(\z, c))$ as conclusion
(Theorem~\ref{thm:e-sound-oga}).  Appendix~\ref{app:oga} gives the
proof of $\LF$ in \oga.  The $\forall$-free base uses Lemma~Q and two
ATI steps, one over the fuel and one over the numeral.  The nesting
step reuses the same two ATI steps at a nested matrix.  Whether \oga
proves $\ConR$ itself, the consistency sentence for all of \rga's
proofs, is an open question.

\paragraph{Per-level and reach-guarded forms.}
The level theorem also has per-level forms.  For each meta-level
number $N$, a compiled scan $\Ok_N$ certifies quantifier nesting up to
depth $N+1$; with $\Ok_N$ in place of the certificate, the per-level class
$\LvC_N$ and the closed hypothesis $\LF_{\le N}$ satisfy the analogues
of Theorems~\ref{thm:con-lv} to~\ref{thm:con-oga}: $\LF_{\le N}$
entails the consistency sentence for $\LvC_N$, \oga proves that
sentence, and \rga neither proves nor refutes $\LF_{\le N}$
(Theorem~\ref{thm:con-lvN}).  At closed codes every $\LvC_N$ lies in
$\LvC$ (Theorem~\ref{thm:con-lv-closed}).  A reach-guarded form asks
for the families only where the evaluation of a case hypothesis's
instance actually reaches a quantifier application, through a
hypothesis $\DynFam(p)$ about the proof code; with the static half
$\Stat(p)$ of $\LvC$ (its swap-node and substitution clauses) as the
class, the level and per-level theorems each follow from the reach-guarded
form by one cut (Theorem~\ref{thm:con-e1}).  Appendix~\ref{app:oga}
gives all three forms.

\subsection{Full RGA under the residual: locating the obstruction}
\label{sec:con-res}

The last result covers every proof code.  This result is a remark
rather than a consistency theorem in the sense of the previous
subsections: the result locates the obstruction to $\ConR$.  The
result assumes the fold's two
interfaces in their general form, stated as one formula about the
given code.

\begin{defi}[\textsf{iRes}]
\label{def:con-res}
The \emph{residual} $\Res(p)$ is the conjunction of two formulas
about the proof code $p$.  In both formulas, ``$p$ is accepted'' means that the
checker's validity test $\mathit{TreeValid}(p)$ holds.  The guard half
says: if $p$ is accepted, the guard is decided at every occupied path
point, for every instantiation of the eigenvariables and every fuel
bound.  The swap half says: if $p$ is accepted, at every occupied path
point of the five swap rules, the \emph{guarded swap body} holds at
the rule's critical formula $C$ and the arm's environment $z$.  The
body says: if $X$ and $Y$ are well-formed codes of terms (four decided
tests), then if every entry of $z$ is a well-formed term code, and $X$
and $Y$ have records of the same value, the fuel-bounded record of
$C$'s instance at $X$ gives the record of $C$'s instance at $Y$.
\end{defi}

\noindent
The four guards on $X$ and $Y$, and the test on $z$, are decided.  The
implications behind them have records as antecedents, which are not
shown decided.  Such implications can be stated, since an implication
abbreviates a disjunction, but \rga's discharge rule does not apply to
them.  These implications are usable here because the residual is a context
hypothesis, only ever eliminated.

\begin{thm}[\textsf{rga\_con\_under\_res}]
\label{thm:con-res}
$\vz\jnat,\ \Res(\vz) \vdash \neg\Chk(\vz, \bquo{\z = \tofn 1})$.
\end{thm}

\noindent
Theorem~\ref{thm:con-res} is an entailment from the residual, not the
implication $\emptyset \vdash \Res(\vz) \limp \neg\Chk(\vz, \bquo{\z =
\tofn 1})$.  $\Res$ is neither shown decided nor of the form $\exists
s.\, D(s)$ with $D$ decided, so neither deduction theorem applies
(Theorem~\ref{thm:in-dt}).  Whether the implication holds is open.
The theorem's soundness companion has the same shape, with $\forall c.\
\Chk(\lft\vz, c) \limp \Rcd_1(\Subst(\z, c))$ as conclusion
(Theorem~\ref{thm:c-res-sound}).  Where \rga refutes acceptance, it
derives the residual: both halves are implications from
$\mathit{TreeValid}(p)$, so one $\lor$-introduction gives each half, and at
the numeral of every rejected code closed evaluation supplies the
refutation (Appendix~\ref{app:fold-res}).

\rga does not derive the residual at an arbitrary code:

\begin{cor}[\textsf{iRes\_underivable}]
\label{thm:con-res-underivable}
$\vz\jnat \;\nvdash\; \Res(\vz)$.
\end{cor}

\noindent
\emph{Proof.}  If \rga derived $\Res(\vz)$ from $\vz\jnat$ alone, a
cut with Theorem~\ref{thm:con-res} would give $\vz\jnat \vdash
\neg\Chk(\vz, \bquo{\z = \tofn 1})$, and $\forall$-introduction would
give $\emptyset \vdash \ConR$, which G\"odel's second theorem
([\citen{goedel-paper}, theorem~\textsf{rga\_goedel2}]) excludes.  \qed

\medskip
With Theorem~\ref{thm:con-res}, the corollary shows where the
obstruction to $\ConR$ sits.  Given the two interfaces, the fold goes
through for every proof, and what \rga cannot derive for an arbitrary
code is the residual itself: the guard's decidedness and the guarded
swap, at every point of the proof.  The next results locate the
obstruction inside the fragment.

\paragraph{The residual at the codes of the fragment.}
Write $\mathit{ResGuard}(p)$ and $\mathit{ResSwap}(p)$ for the guard
half and the swap half of $\Res(p)$, and $\mathit{NoSwap}(p)$ for the
decided scan ``no node of $p$ is read as one of the five swap rules''.

\begin{thm}[\textsf{iFrag\_iResGuard}, \textsf{iFrag\_noswap\_iRes}]
\label{thm:con-res-frag}
\leavevmode
\begin{enumerate}
\item $\emptyset \vdash \forall p.\ \Frag(p) \limp
  \mathit{ResGuard}(p)$.
\item $\emptyset \vdash \forall p.\ \Frag(p) \wedge \mathit{NoSwap}(p)
  \limp \Res(p)$.
\item $\emptyset \vdash \Res(\bquo{P_{\mathrm w}})$.
\end{enumerate}
\end{thm}

\noindent
The guard half is the fragment's guard interface stated as a formula;
each of its implications has a decided antecedent, so \rga can
introduce the guard half.  Where no node is read as a swap rule, the swap half
holds because the swap half's antecedent, that the node has a swap tag, is refuted
at every node.  Part~3 is part~2 at the code of $P_{\mathrm w}$, which
passes the scan.  So \rga derives the residual at the code of a proof
with a quantifier, and Theorem~\ref{thm:con-res} is not vacuous.

\rga does not derive the residual from membership in the fragment,
however.  The residual's swap half fails at the code of an accepted three-node
proof $P_3$, which lies in $\Frag$, in $\LvC$ and in every $\LvC_N$
and whose root is a swap node, by the disjunction property and
G\"odel's second theorem at an unbounded proof search
(Theorem~\ref{thm:con-res-not-frag}).  So the residual follows from
none of these classes, as an implication or as an entailment, and no
set of hypotheses true in the $\omega$-semantics, $\LF$ and
$\LF_{\le N}$ included, changes this (Theorem~\ref{thm:con-res-lf}).
Beyond $\Frag$ the guard half fails too, at proofs with a quantified
case hypothesis, by Turing's argument
(Theorem~\ref{thm:con-guard-not}).  Within the fragment, then, the swap
positions are responsible: the guard half holds on all of $\Frag$, the
whole residual at the codes of $\Frag$ without swap nodes, and the swap half
fails at $P_3$.  Appendix~\ref{app:fold-res} gives these theorems and
their proofs, and Appendix~\ref{app:traps-t01} shows why the swap body
needs its guards.

\subsection{The residual in stronger systems}
\label{sec:con-ext}

A stronger system might supply what \rga cannot: an extension of \rga
whose further principles derive the residual would run the same fold
and prove $\ConR$.  The natural candidates are sound for the
$\omega$-semantics: \rga under hypotheses true there, such as $\ConR$;
\oga; and \rga extended by induction over strictly positive fixed
points.  The development proves the $\omega$-soundness of all three
candidates, and none of these candidates can run the fold.  The witness is an accepted proof $P_{\mathrm h}$ of nine
judgments with an empty root context.  $P_{\mathrm h}$ concludes $Z \lor Y$, where
$Z \ldef (\lambda.\,\z\jnat)\cdot\z$ and $Y \ldef W\cdot W$ is the
ungrounded term of Theorem~\ref{thm:con-res-lf}, by eliminating the
disjunction $Y \lor Z$, so one of the proof's case hypotheses is $Y$ and one of
its nodes is the hypothesis judgment $Y \vdash Y$.

\begin{thm}[\textsf{xder\_not\_fold\_residual}, \textsf{phP\_fMot\_not\_omega}]
\label{thm:con-res-ext}
$\emptyset \vdash \Chk(\bquo{P_{\mathrm h}}, \bquo{Z \lor Y})$, and:
\begin{enumerate}
\item the guard half $\mathit{ResGuard}(\bquo{P_{\mathrm h}})$ is not
  true in the $\omega$-semantics;
\item the fold's motive at the hypothesis node of $P_{\mathrm h}$ is
  not true in the $\omega$-semantics, although, for every residual $R$
  with $\mathit{arms}(R)$, the fold derives the motive in every context that
  derives $R(\bquo{P_{\mathrm h}})$ and in which $v_1$ is not free;
\item let $\vdash_X$ be any entailment over \rga's terms that is sound
  for the $\omega$-semantics: whenever $\Gamma \vdash_X c$ and $\Gamma$
  is a well-formed context of hypotheses true in the
  $\omega$-semantics, $c$ is true there.  Then no such $\Gamma$ has
  $\Gamma \vdash_X \mathit{ResGuard}(\bquo{P_{\mathrm h}})$, and no such
  $\Gamma$ has $\Gamma \vdash_X R(\bquo{P_{\mathrm h}})$ for any residual $R$
  with $\mathit{arms}(R)$ for which $R(\bquo{P_{\mathrm h}})$ is
  well-formed and does not contain $v_1$ free.
\end{enumerate}
\end{thm}

\noindent
Part~3 is stated for an abstract entailment, a locale whose one
assumption is $\omega$-soundness.  \rga's entailment is an instance,
and the $\omega$-soundness theorems of \oga and of the extension by
positive induction have exactly the assumed shape.  The
obstruction is the fold's motive, not the choice of residual.  At the
node $Y \vdash Y$ the motive, $\neg\mathit{Guard} \lor
\Rcd_1(\text{conclusion})$, yields inside \rga the disjunction
$\neg\Rcd^2_1(\bquo{Y}) \lor \Rcd_1(\bquo{Y})$: the hypothetical
judgment ``assuming $Y$, $Y$'' internalized as the implication ``if $Y$
is true, $Y$ is true''.  Neither record of $Y$ has a value in the
$\omega$-semantics, so the implication is not true there, and no
$\omega$-sound system derives it.  That internalization is the step that grounded
reasoning refuses throughout: an entailment from a hypothesis becomes
an implication only at a decided hypothesis
(Theorem~\ref{thm:in-dt}), assuming is not establishing
(\S\ref{sec:pillars-tarski}), and conditional proof at arbitrary
antecedents makes every extension with modus ponens inconsistent, by
Curry's paradox~\cite{rga-paper}.  Part~1 strengthens
Theorem~\ref{thm:con-guard-not}, which shows by Turing's argument and
the disjunction property that \rga does not derive the guard half at
every code; here the guard half fails at one code, semantically, in
every $\omega$-sound system.  Appendix~\ref{app:fold-res} gives the
proof $P_{\mathrm h}$ and the proofs.

\begin{rem}[not mechanized]
\label{rem:con-any-motive}
Part~3 covers every residual that the paper's fold accepts.  That no
internal-soundness fold whatsoever proves $\ConR$ in an
$\omega$-sound extension is not a theorem: stating it needs a general
notion of such a fold, of what its motive must give at the root and
how the steps for the rules consume the motive.  The argument for that
claim is the one above.  At a hypothesis node $h \vdash h$, any motive that yields
the truth of the conclusion from the truth of the hypotheses renders
``if $h$ is true, $h$ is true'', which at an ungrounded $h$ such as $Y$
is excluded middle for a sentence without a value, and so is not true
in the $\omega$-semantics.
\end{rem}

\subsection{What the theorems show}
\label{sec:con-summary}

Each theorem of the family is $\ConR$ with one restriction made
explicit.  The two internal theorems restrict the class of proofs, by
decided scans, and are unconditional; the $\forall$-free theorem follows from
the fragment theorem inside \rga, and their proof-theoretic strength is open.  The
level theorem uses a third decided class, $\LvC$, whose quantified case
hypotheses and critical formulas nest to any certified depth; the level theorem is an
entailment from one closed hypothesis, $\LF$, which \rga neither proves
nor refutes and \oga proves.  The residual result drops the class
restriction and keeps a hypothesis about the proof code, which \rga
cannot derive in general, nor from membership in any of the classes,
even under the level formulas; within $\Frag$ the residual's swap half is what
fails, and beyond $\Frag$ its guard half fails too.  The result locates
the obstruction rather than covering more proofs, and a stronger
system does not remove the obstruction: in every $\omega$-sound extension the fold's
motive fails at a hypothesis node whose hypothesis is ungrounded
(Theorem~\ref{thm:con-res-ext}).  G\"odel's second
theorem ([\citen{goedel-paper}, theorem~\textsf{rga\_goedel2}]) shows that the restrictions cannot all
be removed, and its central step ([\citen{goedel-paper}, theorem~\textsf{goedelG\_of\_expl\_d3g\_ctx}]) is the
decided-matrix split of \S\ref{sec:con-shape}.

\paragraph{No Gentzen-style proof of this shape.}
A Gentzen-style consistency proof needs a reduction on proof codes that
\rga proves correct, and a proof in \rga that the reduction's descent
stops.  Here ``correct'' means that every accepted proof of a
$\forall$-free conclusion either lies in the fragment of
\S\ref{sec:con-frag} or reduces to an accepted proof of the same
conclusion with a smaller notation.  Such a pair would give $\emptyset
\vdash \ConR$, through Theorem~\ref{thm:con-frag}.  So, by
{}[\citen{goedel-paper}, theorem~\textsf{rga\_goedel2}], for no reduction that \rga proves correct does
\rga also prove that the reduction's descent stops
([\citen{goedel-paper}, Appendix ``Consequences'']).  The descent statement here is a
single no-descent sentence about one sequence, not a
transfinite-induction schema.  Whether that sentence can be derived from such a
schema in \rga, whose implications there would have undecided
antecedents, is open; the development does not formalize such a schema.

\paragraph{The proof search for $\z = \tofn 1$.}
Let $\mathit{Stops}$ be the value record
$\Rcd_1(\bquo{\mathit{conSearch}})$ (Definition~\ref{def:con-rcd}),
which says that the search through all proof codes for an accepted
proof of $\z = \tofn 1$ returns.  The search never returns (a HOL
fact), but \rga neither proves nor refutes $\mathit{Stops}$, and does
not prove it decided.  As entailments, a proof of $\z = \tofn 1$ makes
the search stop, $\PrC(\bquo{\z = \tofn 1}) \vdash \mathit{Stops}$, and
by contraposition $\neg\mathit{Stops} \vdash \ConR$.  With
the theorems \textsf{goedelG\_under\_ConRGA} and
\textsf{rga\_goedel2\_internal\_contra} of [\citen{goedel-paper}], $\neg\mathit{Stops}$
also entails $G$, $\neg\PrC(\bquo{G})$ and $\neg\PrC(\bquo{\ConR})$.  In
the other direction, $\mathit{Stops} \nvdash \z = \tofn 1$; whether
$\mathit{Stops} \vdash \PrC(\bquo{\z = \tofn 1})$ and whether $\ConR
\vdash \neg\mathit{Stops}$ are open ([\citen{goedel-paper}, Appendix ``The entailments and the cut transformers'']).
So ``the search never stops'' is true but never established in \rga.
\rga can use that statement only as a hypothesis, and from that
hypothesis \rga derives its own consistency, its G\"odel sentence, and the unprovability
of $\ConR$ and of $G$.

\paragraph{What the internal derivations use.}
A natural worry is circularity: the internal theorems might borrow
\rga's soundness from HOL.  An \rga derivation cannot appeal to a HOL
theorem, so the question is what the HOL constructions of these
derivations use.  A check that runs with every build reads the
theorem-dependency graph that Isabelle's kernel records.  The check confirms
that the constructions of this section's internal derivations (the
consistency and soundness sentences of the two decided classes, the
containments, Lemma~Q, the level theorems with their \oga forms, and
the residual theorem) use neither \rga's soundness induction nor an
$\omega$-soundness induction, identified by name and by the shape of
their statements, nor the HOL shadow models.  These constructions may use completeness, which turns a HOL evaluation of a
closed term into a derivation and adds no reasoning of its own.  Facts
about particular closed codes, every non-derivability claim and the
faithfulness of the checker do rest on soundness; for non-derivability
this dependence is unavoidable, since a claim that \rga does not derive some
sentence implies \rga's consistency.  Each internal theorem asserts
that a derivation exists; extracting one concrete derivation is left
to future work (Section~\ref{sec:concl}).
Appendix~\ref{app:machines-audit} lists the roots and the controls and
says what the check cannot see.

\section{Mechanization}
\label{sec:mech}

Every theorem in this paper is proved in Isabelle/HOL (Isabelle2025-2),
and Appendix~\ref{app:endpoints} maps each theorem to its Isabelle name and
theory file.  The development builds on the companion paper's
formalization of \rga~\cite{rga-paper}, the \emph{companion
formalization}, which the development extends with ten
sessions: the HOL models of the design, the coding kit, the decider,
the checker, the $\forall$-free consistency proof, the value record,
G\"odel's theorems, the fragment and residual theorems, the Liar and
divergence results, and the level theorems and \oga results.  These
sessions hold about 267{,}000 lines of theory text (\texttt{wc -l}) and
about 13{,}600 named items in 97 theory files, measured, like every
size and timing figure in this paper, at the point named under
``Artifact'' below.

\paragraph{Trusted base.}
Besides Isabelle/HOL's kernel (the validation pins its Pure and HOL
heaps by hash), a reader must trust the definitions in which the
theorems are stated, in two tiers.  The headline results are stated
in \rga's derivability relation, with its 42 rules, and the proof and
entailment relations over it; the coding of terms, $\quo{\cdot}$; the compilation
$\cmp$; and this paper's $\ConR$, $\Chk$ (with its validity test),
$\Dec$, $\Aff$, $\Tru$ and the code $\bquo{\z = \tofn 1}$, the
transparent machines of Section~\ref{sec:inside}.  The other results
are also stated in the companion paper's grounded semantics and
$\omega$-semantics, the same paper's well-formedness of terms and contexts,
$\forall$-freeness, and the derivability relation of \oga.  No
headline result mentions an opaque twin.  The locale \textsf{rga} in
which all these definitions are stated has no parameters and no assumptions.
$\Chk$ need not be read line by line.  Faithfulness and genuineness
(Theorem~\ref{thm:in-chk}) pin $\Chk$ down extensionally, at proof codes
and at every numeral, and depend only on the derivability relation and
the coding; that $\ConR$ is the natural consistency sentence and not
merely one with the right extension is the intensional criterion of
\S\ref{sec:inside-cmp}, which the derivability conditions of
{}[\citen{goedel-paper}, \S\,``The derivability conditions \rga has''] meet.  Validation fingerprints each trusted
definition of both tiers, and every constant that
Appendix~\ref{app:endpoints} names for a numbered result, over its
whole definitional closure, every project definition it depends on
transitively, so a change to any definition in such a closure fails
validation.  Appendix~\ref{app:endpoints-statements} prints the statements of
the headline results and the definitions of $\ConR$, $\Chk$, $\Aff$,
$\Tru$ and $\bquo{\z = \tofn 1}$, generated by Isabelle from the
build rather than transcribed by hand.

\paragraph{Assurance.}
The development contains no \textsf{sorry}, and an oracle audit over
every fact of the final session's heap finds no oracle.  An axiom audit
confirms that every axiom of the final session beyond those of
Isabelle/HOL's \texttt{Main} is a definition, a type definition or the
arity axiom that declaring a type constructor adds (5{,}858 such axioms;
no axiom of any other kind).
Three further kinds of check run with every build and, like the axiom
audit, fail the build on a violation.
The choice audit and the code-generation check of
\S\ref{sec:inside-audit} make transparency a checked property of the
machines, witnesses and transformers that the paper names.  Dependency
checks read Isabelle's theorem-dependency graph.  One dependency check confirms that the
HOL constructions of the internal derivations, 163 theorems in all, do
not use \rga's soundness theorem or a HOL shadow model
(\S\ref{sec:con-summary}).  This check forbids the soundness inductions by
name.  In the dependency cones of these theorems it finds no theorem
that cites an induction rule of the derivation relations, and it reads
the statement of every named theorem: none has a derivability premise
and a semantic conclusion, apart from the eight introduction rules of the
semantics at the universal quantifier, which define a universal's value
by an \rga proof.  Another dependency check confirms that G\"odel's second
theorem, its internal form, the theorems derived from them and
explicit reflection at each numeral (31 roots) use neither the decider
nor the truth predicate: no theorem
in the dependency cones of these roots has a statement that mentions the truth
predicate, the decider's verdicts or its compiled table, and neither
the search witness nor the decider's internal completeness is in a
cone.  A third dependency check confirms that on the path of
{}[\citen{goedel-paper}, theorem~\textsf{rga\_goedel2}] the only implication introductions in the
G\"odel theories are those of EXPL and D3G, over decided antecedents;
\rga's rules guarantee that every derivation discharges only decided
hypotheses, and the check records the shape of this particular proof.
A validation script compares statement and definitional fingerprints
with a recorded baseline, so that a change to a fingerprinted
statement, or to any definition beneath it, fails validation.  The
statement fingerprints cover the development's endpoints and every
theorem that a numbered result of this paper cites; the definitional
fingerprints cover the closures of the endpoints and of the trusted
definitions.  The script also
confirms that every statement and definition of the companion paper's
formalization is unchanged.  Appendix~\ref{app:machines-audit} gives
the roots, the controls, the closure and cone sizes, and the limits of
each check.  Statements marked as remarks are not mechanized.  The
whole chain, from the grounded-deduction base through the last
session, builds in about eight minutes of wall-clock time on one
workstation (an Apple M4 Max with 36\,GB of memory, ten threads);
Appendix~\ref{app:machines-audit} on the build-time check gives the
command and the conditions.

\paragraph{Size of the G\"odel proof.}
Measured on Isabelle's theorem-dependency graph
(Appendix~\ref{app:machines-audit} gives the method), the dependency
cone of [\citen{goedel-paper}, theorem~\textsf{rga\_goedel2}] has about 62{,}000 lines of theory
text, most of it machinery the paper needs anyway: the coding kit, the
checker and the companion formalization.  The part in the session
\texttt{RGA\_Goedel} is about 3{,}100 lines.  The cone of the search
route of [\citen{goedel-paper}, \S\,``G\"odel's second theorem''] is larger, about 99{,}000 lines, and
that of the internal second theorem smaller, about 60{,}000 lines.  The
figures depend on proof style and are quoted only as relative sizes.

\paragraph{How it was built.}
The development was produced over about two months by the author
directing Claude models from Anthropic, followed by a cleanup that
removed what the results do not need and renamed the library for
readability.  No endpoint statement was weakened, and no definition
beneath the endpoints changed its meaning
(Appendix~\ref{app:machines-audit}, ``Fingerprints'').  The workflow,
its error profile and its costs are reported separately.

\paragraph{Artifact.}
\label{sec:mech-artifact}
The development is recorded at the tag \artifacttag of the project
repository, whose README gives the build command and the sessions.
The reflection results of \S\ref{sec:pillars-reflect} are in three
further sessions, \texttt{RGA\_Reflect}, \texttt{RGA\_Poly} and
\texttt{RGA\_Tin}, which are not among the ten sessions above.  The
first two build on the session of G\"odel's theorems, and
\texttt{RGA\_Tin} builds on the session of the Liar and divergence
results and uses the generic checker of the session
\texttt{RGA\_Checker\_Gen}.  The command
\begin{center}
\small\texttt{isabelle build -b -v -o threads=10 -d hol ConOGA RGA\_Poly RGA\_Tin}
\end{center}
builds the ten sessions, the three further sessions and
\texttt{RGA\_Checker\_Gen}.  The
theorem names and theory files that this paper gives are those of the
development at that tag.  The size, timing and audit figures (this
section and Appendix~\ref{app:machines-audit}) were measured at the tag
\texttt{conrga-}\allowbreak\texttt{paper-v3}, the measured endpoint of
this paper's proof campaign, including its cleanup.  Later commits
moved some of the theorems the paper cites into the companion
formalization and added other theorems; the figures were not measured again after
those commits.

\section{Related work}
\label{sec:related}

The design genealogy of \rga itself (Fitch's basic logic, illative
combinatory logic, the recursive $\omega$-rule, Kripke's fixed points)
is discussed in the companion paper~\cite{rga-paper}, as are \rga's
completeness for grounded truth, the recursive enumerability of
grounded truth and \rga's $\omega$-incompleteness.  This section
places the present paper's results.

\paragraph{Mechanized incompleteness.}
Paulson gave the first mechanized proof of G\"odel's second theorem,
in Isabelle/HOL for the theory of hereditarily finite
sets~\cite{paulson14machine,paulson15mechanised}, following
\'Swierczkowski~\cite{swierczkowski03finite}.  Popescu and Traytel
distilled abstract versions of both theorems, with the derivability
conditions as assumptions, and instantiated them using Paulson's
proofs~\cite{popescu19formally,popescu21distilling}.  Popescu and Traytel's framework
assumes only intuitionistic reasoning, but its instances are
classical, and it assumes implication introduction, which \rga lacks.  Saitou and Noguchi proved the second
theorem in Lean~4 for every $\Delta_1$ theory extending
$\mathrm{I}\Sigma_1$, establishing D3 by a model-theoretic argument
through the completeness theorem~\cite{saitou26mechanizing}.  Coquand
formalized in Agda the second theorem for Church's basic recursive
arithmetic, following Guard's unpublished lecture notes of 1963, which
follow Church's Princeton lectures~\cite{coquand26auto}; the report's
title has ``binary'' where its text has ``basic''.
All of these mechanizations concern classical, or free-variable and hence everywhere decided,
systems.  Carnielli and Fuenmayor analyzed the incompleteness theorems
under weak logics of formal inconsistency in a shallow embedding in
Isabelle/HOL, with the fixed point taken as a premise and the
derivability conditions cited rather than proved, and found that the
second theorem then needs extra premises asserting the consistency of
particular formulas~\cite{carnielli20goedel}; Carnielli and Fuenmayor's setting is
abstract, with no arithmetic.  Gross, Gallagher and Fallenstein
formalize L\"ob's theorem in Agda for a postulated, syntactic
provability predicate, and read the theorem, through Curry--Howard, as a limit
on self-interpreters for total languages~\cite{gross16lob}.  The mechanized proofs of the second
theorem known to the author are for classical
theories~\cite{paulson15mechanised,saitou26mechanizing,hydras-goedel2}
or for quantifier-free classical arithmetic~\cite{coquand26auto}; the
non-classical ones are abstract, with the fixed point or the
derivability conditions assumed~\cite{popescu21distilling,carnielli20goedel}.
No mechanized proof of the second theorem is known for Heyting arithmetic, and none for a concrete
arithmetic with quantifiers whose object logic is not
classical.\footnote{The search covered the Lean~4 Foundation library,
the Rocq developments hydra-battles and the Coq library for first-order
logic, Agda, the Isabelle Archive of Formal Proofs and HOL Light, and,
by search, Metamath, Mizar, HOL4, Nuprl and Minlog, as of October
2026.}
Several developments mechanize the first theorem only: Shankar's in
the Boyer--Moore prover~\cite{shankar94metamathematics}, Harrison's in
HOL Light~\cite{harrison-hol-light-godel}, and the synthetic
treatments of undecidability and essential incompleteness by Kirst and
coauthors, also in Coq~\cite{kirst23synthetic,kirst23godel}.
O'Connor's development in Coq~\cite{oconnor05essential} proved the
first theorem; its maintained version also derives the second for
extensions of Peano arithmetic, from D2 and D3 taken as
hypotheses~\cite{hydras-goedel2}.
The present development also differs in architecture: Isabelle/HOL
verifies statements about \rga, and \rga in turn reasons about its own
codes, through a checker and a decider written as \rga terms.

\paragraph{Self-verification and reflection in proof assistants.}
Several proof assistants have been verified in themselves or in a
close relative.  Harrison verified HOL Light's logic and a model of its
kernel inside HOL Light, proving the consistency of HOL without the
axiom of infinity, and of full HOL under an additional
axiom~\cite{harrison06towards}.  Kumar, Arthan, Myreen and Owens
formalized higher-order logic and a verified implementation in
HOL4~\cite{kumar16self}, and Candle, a clone of HOL Light, is verified
down to its machine code~\cite{abrahamsson22candle}.  Fallenstein and
Kumar obtained a reflection principle for HOL in HOL4 by assuming a
large cardinal~\cite{fallenstein15proof}.  Milawa extends itself by
reflection, and its kernel is verified in HOL4 down to the machine
code that runs it~\cite{davis15reflective}.  MetaCoq verifies a type
checker for Coq's kernel in Coq, with strong normalization as an
axiom~\cite{sozeau20coq}, and Lean4Lean verifies parts of a type
checker for Lean in Lean~\cite{carneiro26lean4lean}.  Where these
developments prove a system consistent, they work in a stronger theory
or under an added assumption, as the second theorem requires.  This
paper is no exception: HOL proves \rga consistent, and \rga does not
prove its own consistency ([\citen{goedel-paper}, theorem~\textsf{rga\_goedel2}]).  What the paper
adds is the object theory's own reasoning: \rga proves theorems about
its own checker and decider, the decided steps of G\"odel's proof
are derivations inside \rga, and \rga derives its own second
incompleteness theorem as an entailment, $\ConR \vdash
\neg\PrC(\bquo{\ConR})$, and in decided-antecedent form, $\forall
q.\,(\Chk(q, \bquo{\ConR}) \limp \neg\ConR)$
([\citen{goedel-paper}, Section ``G\"odel's theorems for \rga'']).

Brown and Palsberg built a self-interpreter for System~F$_\omega$, a
strongly normalizing language, breaking through what they call the
normalization barrier: the conventional wisdom that a strongly
normalizing language has no self-interpreter~\cite{brown16breaking}.
\rga's terms are untyped, so its situation is that of the untyped
$\lambda$-calculus, where self-interpreters are
classical~\cite{mogensen92efficient}; its decider is moreover a
primitive-recursive function of an explicit fuel bound
(Theorem~\ref{thm:in-tab}), so that no barrier of this kind arises.
The content here is what \rga proves about its interpreter, not that
the interpreter exists.  In proof by reflection, a decision procedure
written and proved correct inside a proof assistant turns a closed
computation into a proof, through the conversion
rule~\cite{boutin97using}.  \rga's conversion rules move only along
single head steps and do not evaluate.  Theorem~\ref{thm:in-cmp-num}
takes the conversion rule's place: HOL evaluates the closed term, and
open completeness turns the evaluation into a derivation.  Read
informally, that derivation follows the evaluation step by step, so
\rga lacks reflection's small certificate.
Section~\ref{sec:pillars-reflect} on the truth predicate as a
reflection medium builds the rest of the pattern inside \rga.  It
gives a polynomial normalizer proved correct by one internal sentence
over arbitrary codes, with $\Tru$ as the denotation, so that
polynomial identities at a free variable follow from one closed
evaluation of the test (Theorem~\ref{thm:pl-rf-reflect}).  The other
internal theorems at abstract arguments are derivations by internal
induction, not computations turned into proofs: the checker's
structural laws, Theorem~\ref{thm:in-never}, and the certificate
builder of D3G ([\citen{goedel-paper}, \S\,``D3G with an explicit primitive-recursive witness'']).
Step-indexed models of recursive types~\cite{appel01indexed} and
Nakano's modality for recursion~\cite{nakano00modality} underlie
logics such as LSLR~\cite{dreyer11logical} and Iris~\cite{jung18iris}.
These logics have a ``later'' modality $\triangleright$, which LSLR
traces to G\"odel--L\"ob logic, and a L\"ob rule: from $\triangleright
P \limp P$ infer $P$, valid by induction on the step index.  \rga has no later
modality, and its fuel indexes evaluations, not the truth of
propositions.  \rga's L\"ob rule ([\citen{goedel-paper}, theorem~\textsf{no\_code\_unprovable}]) is the rule
of provability logic, and it follows from G\"odel's second theorem and
the disjunction property, not from induction on fuel.  Fuel does enter
\rga's internal theorems in another way: the affirmation test
$\Aff(s, c)$ is a step-indexed judgment written as an \rga term, and
Theorems~\ref{thm:in-never} and~\ref{thm:pl-int}, which quantify over
the test's fuel inside \rga, are internalized step-indexed reasoning.

\paragraph{Assurance in proof assistants.}
The build-time checks of Section~\ref{sec:mech} have precedents.  Coq's
\texttt{Print Assumptions} and Lean's \texttt{\#print axioms} list the
axioms a theorem depends on~\cite{huet16coq,moura21lean}, as the
oracle and axiom audits do here; Coq's \texttt{Qed} and \texttt{Defined}, and
Lean's \texttt{opaque} and \texttt{irreducible}, control which
definitions unfold, the property that the choice audit makes checkable
for the machines of this paper; Isabelle records theorem dependencies
in its proof bodies, which the dependency checks read.  Pollack
separates two questions about a machine-checked proof, whether the proof is a
derivation and whether its statement means what is
claimed~\cite{pollack98believe}, and Wiedijk asks that a prover read
back what it prints~\cite{wiedijk12pollack}; the fingerprints and the
transparency check of this paper address the second question for the paper's
trusted definitions.

\paragraph{Hilbert--Bernays and free-variable arithmetic.}
In the first detailed proof of the second
theorem~\cite{hilbert39grundlagen}, the derivability conditions are
stated on the proof predicate, with the third restricted to
primitive-recursive equations as antecedents~\cite{kurahashi20note};
\rga's decided D3, which \rga proves for every $\varphi$ by a
construction whose instance at the code term $S_G$ is D3G
([\citen{goedel-paper}, theorem~\textsf{D3\_decided\_any}]), and also as one sentence over the code
([\citen{goedel-paper}, theorem~\textsf{D3\_uniform}]), has a premise certificate of the shape
of that third condition.  The first two conditions
and Hilbert and Bernays' consistency statement are implications over
provability formulas, and Hilbert and Bernays' proof reasons with such
implications~\cite{kurahashi20note}.  The decided form of the route
used in [\citen{goedel-paper}, Section ``G\"odel's theorems for \rga''] is that of free-variable arithmetic.
In Coquand's formalization, after Guard and Church, self-reference is a
term equation, D3 and ex falso are primitive-recursive proof
transformers with decided antecedents, the G\"odel formula is derived
from the consistency formula by propositional reasoning on decidable
formulas followed by generalization, and the theorem is a
transformation of derivations~\cite{coquand26auto}; the decided case
split of [\citen{goedel-paper}, theorem~\textsf{goedelG\_of\_expl\_d3g\_ctx}] is that propositional step.  In those settings every formula
used is decided, so the restriction to decided antecedents costs
nothing.  In \rga the restriction is forced, because every implication of the
textbook proof fails while the corresponding entailments hold
([\citen{goedel-paper}, Section ``G\"odel's proof in four forms'']), and the decided forms are those
entailments written as sentences.  These decided forms must be built inside a logic
with quantifiers and undecided sentences,
where each term used as a number must be shown to be one.  Jeroslow
showed that D2 can be dispensed with for Hilbert and Bernays'
consistency statement once D3 holds for every primitive-recursive
term~\cite{jeroslow73redundancies,kurahashi20note},
Kurahashi compared the sufficient sets of
conditions~\cite{kurahashi20note}, and Beklemishev and Shamkanov proved
abstract versions of the second theorem over non-classical consequence
relations (below).
Smory\'nski's handbook chapter is the standard account of the
derivability conditions~\cite{smorynski77incompleteness}, and Visser
discusses how to state the second theorem in a general, coordinate-free
form~\cite{visser16second}.

Proofs of the second theorem that avoid the textbook route are not
new.  Free-variable arithmetic proves it with decided proof
transformers.  Jeroslow's proof avoids D2 through a fixed point that
asserts its own refutability~\cite{jeroslow73redundancies,kurahashi20note}.
Kritchman and Raz's proof argues from Chaitin's incompleteness theorem
and uses no self-referential sentence~\cite{kritchman10surprise};
Coquand has mechanized Kritchman and Raz's proof for Basic Recursive
Arithmetic~\cite{coquand26chaitin}.  That every implication of the
textbook proof fails in \rga while every entailment holds
([\citen{goedel-paper}, Section ``G\"odel's proof in four forms'']) is the material conditional and the
disjunction property at work, not a new route.  What this paper adds
is the decided route carried out inside a gappy arithmetic: the D3G
certificate builder, the cut transformers, and formalized G\"odel~II in
decided-antecedent form with its explicit witness.  Artemov's schematic
consistency, ``$n$ is not a proof of $\z = \tofn 1$'' for each
$n$~\cite{artemov19provability}, is the classical analogue of the
instance-by-instance consistency that \rga derives on its decided
matrix ([\citen{goedel-paper}, \S\,``Formalized G\"odel~I and the consistency sentence'']).
Heyting arithmetic does not prove its consistency, by the textbook
proof.  Heyting arithmetic does not refute its own G\"odel sentence either, most simply
because it is contained in Peano arithmetic, which is
$\Sigma_1$-sound; inside intuitionistic arithmetic, Heyting arithmetic's existence
property and closure under Markov's rule together play the role that
the witness property plays
in \rga~\cite{troelstra73metamathematical,friedman78classically}.  The
textbook proof goes through in Heyting arithmetic, which derives
$\mathit{Con} \limp G$ for its own consistency sentence and G\"odel
sentence by implication introduction over an antecedent it does not
decide.  Implication introduction over such an antecedent is the step \rga lacks: \rga derives the entailment
$\ConR \vdash G$ ([\citen{goedel-paper}, theorem~\textsf{goedelG\_under\_ConRGA}]), and formalized G\"odel~II
as the entailment $\ConR \vdash \neg\PrC(\bquo{\ConR})$
([\citen{goedel-paper}, theorem~\textsf{rga\_goedel2\_internal\_contra}]), but not the implications.  The two
systems also obtain
their closure properties differently.  Heyting arithmetic's closure properties come from
realizability and from Friedman's
translation~\cite{troelstra73metamathematical,friedman78classically};
\rga's witness property comes from the soundness of its
certificate-carrying semantics.

\paragraph{Self-reference by equation.}
The equation form of self-reference used in [\citen{goedel-paper}, \S\,``The sentences'']
is standard: the strong diagonal lemma provides, for each formula
$A(x)$, a closed term $t$ with $\vdash t = \quo{A(t)}$.  The history of this equation form
runs through the Kreisel--Henkin exchange and is surveyed by Halbach
and Visser~\cite{halbach14self} and by Grabmayr and
Visser~\cite{grabmayr20self}; Kripke obtains the same effect with
constants~\cite{kripke23goedel}.  In computability terms the equation
form is Kleene's second recursion theorem~\cite{kleene38notation}, and
$S_G = \Diag(\bquo{D}, \bquo{D})$ is the usual construction of a
quine.  The device is not new.  What is
specific to \rga is that the equation form is forced:
{}[\citen{goedel-paper}, theorem~\textsf{no\_goedel\_biconditional}] shows that the biconditional form is
unprovable for every sentence.

\paragraph{Forms of consequence.}
That a principle can hold as a rule or an entailment while failing as
an implication is familiar.  The admissible rules of Heyting arithmetic
outrun its derivable ones, and its admissible propositional rules are
exactly those of intuitionistic propositional
logic~\cite{visser99rules,visser02substitutions,iemhoff01admissible,rybakov97admissibility}.
In modal logic, necessitation holds for global but not for local
consequence, which is why the deduction theorem needs a restricted
necessitation rule~\cite{blackburn01modal,hakli12deduction}.  \rga's
necessitation for entailments ([\citen{goedel-paper}, theorem~\textsf{cut1M\_chk}]) is of the
global kind.  Kleene's strong three-valued conditional is material,
and $A \limp A$ is not valid in his logic~\cite{kleene52introduction}.
Field adds a non-material conditional to recover $A \limp A$, and does
not declare all of his rules truth-preserving~\cite{field08saving}.
\rga keeps the material conditional, so it proves no implication
between independent sentences, and it recovers the textbook steps as
entailments and, for provability hypotheses, in decided-antecedent
form.
Beall and Murzi's validity Curry~\cite{beall13two} turns on a validity
predicate obeying ``if $A \vdash B$ then $\vdash \mathrm{Val}(A, B)$''
and ``$A, \mathrm{Val}(A, B) \vdash B$'', and they blame structural
contraction.  Field replies that for validity read as provability the
second principle is local reflection, which the second theorem already
excludes~\cite{field17disarming}.  \rga sides with Field: the first
principle holds for $\PrC$ as [\citen{goedel-paper}, theorem~\textsf{cut1M\_chk}], the second
fails ([\citen{goedel-paper}, theorem~\textsf{formalized\_loeb\_ent}]), and contraction is kept.

Beklemishev and Shamkanov~\cite{beklemishev16abstract} prove the
second theorem, and its formalized version, for abstract provability
structures: a consequence relation with provability and refutability
operators $\Box$ and $\boxtimes$ and no connectives, under conditions
C1--C4, with a Jeroslow fixed point $p \eqv \boxtimes p$.  Condition
C1 asks that $\Box$ be monotone and $\boxtimes$ anti-monotone in the
consequence order.  Beklemishev and Shamkanov's $\top$ is read here as $\z = \z$ and their
$\bot$ as $\neg(\z = \z)$, a classical falsum, not \rga's divergent
$\bot$.  Beklemishev and Shamkanov's formalized theorem is itself an entailment,
$\boxtimes\boxtimes\top \leqslant \boxtimes\top$, the shape of
{}[\citen{goedel-paper}, theorem~\textsf{rga\_goedel2\_internal\_contra}](1).  In \rga, with $x \leqslant y$ read as
$x \vdash y$ and $\boxtimes x$ as $\PrC(\bquo{\neg x})$, conditions C2 to
C4, Beklemishev and Shamkanov's C5 and the monotone half of C1 hold at every sentence.  The
anti-monotone half of C1, contraposition under provability, fails in
general, at the Liar, but holds at decided sentences and at $\Sigma$
sentences with a decided matrix.  That class is closed under $\Box$ and
$\boxtimes$ and contains the Jeroslow fixed point $\neg G$, so on that class
\rga is an abstract provability structure, and Beklemishev and Shamkanov's Theorem~1 holds
for \rga by their proof, mechanized ([\citen{goedel-paper}, Appendix ``The entailments and the cut transformers'']
gives the statements).  So the contrapositive form of \rga's
formalized second theorem, $\PrC(\bquo{\ConR}) \vdash \neg\ConR$, is an
instance of Beklemishev and Shamkanov's Theorem~1 on that fragment: the $\boxtimes\top$ of that theorem is
\rga's $\neg\ConR$ up to interderivability, and necessitation of
entailments ([\citen{goedel-paper}, theorem~\textsf{cut1M\_chk}]) with $\Sigma$-contraposition
moves the instance under $\PrC$.  With an implication obeying the deduction
theorem, Beklemishev and Shamkanov also show that contraction restricted to provability
formulas is what G\"odel's argument needs, and they give a
contraction-free modal system in which formalized G\"odel~II fails.
\rga lacks that deduction theorem and keeps contraction.  Beklemishev and Shamkanov's
framework does not provide a concrete arithmetic with gaps that meets
its conditions; \rga is one, built and checked inside the logic, and
the internal decided route, the cut transformers and the
decided-antecedent form with its explicit witness go beyond the
abstract theorem.  The closest existing logic of this entailment shape
is positive provability logic: sequents $A \vdash B$ between strictly
positive modal formulas, as in Dashkov's positive fragment of
GLP~\cite{dashkov12positive} and Beklemishev's reflection
calculus~\cite{beklemishev14positive}.  The judgmental reading of
modal logic gives the same split between hypothetical and categorical
judgments: necessitation is a rule on categorical judgments, $\vdash A$
gives $\vdash \Box A$, not an implication~\cite{pfenning01judgmental},
and \rga's necessitation for entailments ([\citen{goedel-paper}, theorem~\textsf{cut1M\_chk}])
lifts that rule to hypothetical judgments.

Artemov's Logic of Proofs~\cite{artemov01explicit} makes the proof
operations of provability explicit.  LP's application operation corresponds to
\rga's decided D2 transformer $m_{\varphi\psi}$
([\citen{goedel-paper}, theorem~\textsf{D2\_decided}]), and LP's proof checker ${!}$, with $t{:}A
\limp {!t}{:}(t{:}A)$, to the certificate builder of D3G
([\citen{goedel-paper}, \S\,``D3G with an explicit primitive-recursive witness'']); the witness $\sigma$ is that builder
followed by $\exists$-introduction into $\PrC$.  LP's explicit
reflection $t{:}A \limp A$ is a schema over proof terms, and its \rga
counterpart, $\vdash \Chk(\tofn n, \bquo{A}) \limp A$ for every
numeral $n$ and well-formed $A$, holds, because the checker's verdict
at a numeral is decided and genuine ([\citen{goedel-paper}, theorem~\textsf{chk\_reflect\_num}]).
Like LP's explicit reflection, the \rga counterpart is a schema: one theorem for each numeral and term,
proved in HOL by cases on the decided verdict.  The single sentence over all proofs,
$\forall p.\,(\Chk(p, \bquo{\z = \tofn 1}) \limp \z = \tofn 1)$, fails
([\citen{goedel-paper}, the table of the ingredients of G\"odel's proof in the four forms]), as implicit reflection does.
Artemov's selector proofs~\cite{artemov19provability,artemov25serial} let
Peano arithmetic prove each instance of its consistency scheme
uniformly without proving $\mathrm{Con}(\mathrm{PA})$.

\paragraph{Non-classical logics and the second theorem.}
Every consistent system known to the author that proves its own consistency does
so in one of two ways.  Either it is arithmetically weak, as Willard's
self-verifying theories, which do not prove multiplication
total~\cite{willard01self},
and Pakhomov's set theory, which does not
interpret Robinson arithmetic~\cite{pakhomov19weak}; or it uses an
intensional, nonstandard provability predicate, as in Feferman's
numerations~\cite{feferman60arithmetization}.  Meyer's relevant
arithmetic $R^\#$ proves that $0 = 1$ is unprovable, but only because
$R^\#$ lacks explosion, so that the unprovability of $0 = 1$ no longer expresses
negation consistency; Meyer notes that the classical barriers remain
for negation consistency~\cite{meyer21consistency}.  Friedman and Meyer showed
that $R^\#$ is properly weaker than Peano arithmetic: a strictly
positive theorem of Peano arithmetic is unprovable in
$R^\#$~\cite{friedman92whither}.  Paraconsistent arithmetics
avoid the theorem by being inconsistent~\cite{priest06in}.  Under weak
logics of formal inconsistency, the second theorem needs extra
premises asserting the consistency of particular
formulas~\cite{carnielli20goedel}, the paraconsistent counterpart of
\rga's restriction of discharge to decided hypotheses.

\rga takes
none of these routes.  It proves multiplication, each compiled
primitive-recursive function and Ackermann's function total by
quantified theorems, its
consistency statement uses a checker whose structural laws it proves,
it has explosion, and it is consistent.  Its non-classical logic
changes the form in which the steps of the second theorem's proof
hold, entailments and decided-antecedent forms rather than
implications, but not the theorem.  Nor does \rga's completeness for
grounded truth~\cite{rga-paper} conflict with the theorem: the
G\"odel sentence and $\ConR$ are not grounded-true
(\S\ref{sec:bg-rga}).  Shapiro argued that a dialetheic
theory inherits incompleteness rather than escaping
it~\cite{shapiro02incompleteness}.  Arithmetics with a many-valued truth
predicate keep their arithmetic atoms crisp: Restall studies arithmetic
and truth in \L ukasiewicz's infinitely valued
logic~\cite{restall93arithmetic}, and H\'ajek, Paris and Shepherdson add
a many-valued truth predicate to crisp Peano
arithmetic~\cite{hajek00liar}.  Substructural approaches keep a
transparent truth predicate by restricting a structural rule:
contraction~\cite{zardini11truth} or the transitivity of
consequence~\cite{cobreros13reaching}.  \rga keeps both structural rules, and its gaps
reach the arithmetic sentences themselves; no earlier proof of the
second theorem is known to the author for an arithmetic of that
kind.

\paragraph{Paracomplete theories of truth.}
Kripke's theory of truth~\cite{kripke75outline} and its
axiomatizations, KF~\cite{feferman91reflecting} and
PKF~\cite{halbach06axiomatizing}, and Field's
theory~\cite{field08saving}, are paracomplete about truth but keep
classical arithmetic underneath; see Halbach~\cite{halbach11axiomatic}
for a survey, and Table~\ref{tab:truth-axioms} for the status of their
axioms and rules in \rga.  Undecided arithmetic sentences, including consistency
statements, still obey excluded middle in those theories, and the classical proof
of the second theorem applies.  This classical behavior holds even for PKF, whose logic
comes closest to \rga's.  PKF's partial, strong-Kleene logic governs
its whole language~\cite{halbach06axiomatizing}, but PKF extends Peano
arithmetic, and classical logic still applies to its arithmetic
sentences~\cite{halbach18costs}.  In \rga the arithmetic itself is
paracomplete, and that paracompleteness is why the diagonal biconditional and D2 and
D3 fail as implications (the theorems \textsf{no\_goedel\_biconditional} and
\textsf{no\_D2\_schema\_goedelG} of [\citen{goedel-paper}]), while their entailment and rule forms hold.
PKF's truth rules with side formulas, T-Intro and T-Elim, are
entailments.  For \rga's defined predicate $\Tru$ they fail at the
truth-predicate Liar (Theorem~\ref{thm:br-tin-ent}): the
$\omega$-semantics, for which \rga is sound in every context, makes
that Liar true and its image under $\Tru$ false.  This failure is a fact about
the defined $\Tru$, and the material T-biconditional at its Liar is
explosive (Theorem~\ref{thm:pl-tb-explosive}).  Whether \rga tolerates
a primitive truth predicate with PKF's rules beside $\Tru$ is open
(Section~\ref{sec:concl}).  \rga is also weaker than KF with the axiom
CONS here: that theory proves T-Out at every
sentence~\cite{feferman91reflecting,halbach11axiomatic}, while \rga
fails T-Out as an entailment at its Liar.

PKF also shows what a
partial logic can cost in strength: PKF has the proof-theoretic
strength of ramified analysis up to $\omega^\omega$, against ramified
analysis up to $\varepsilon_0$ for
KF~\cite{halbach06axiomatizing,feferman91reflecting}.  \rga's
proof-theoretic strength is open; \rga proves Ackermann's function
total by a quantified theorem, and, of each closed ground-type
System~T term, that the term denotes a natural, one derivation per
term~\cite[Thm.~5.7]{rga-paper}.

Three further contrasts
are worth noting.  First, Field's theory respects the second theorem
because it does not declare all of its logical rules
truth-preserving~\cite{field08saving}; \rga proves, for each of its
proofs, that the proof's conclusion is true, and Theorem~\ref{thm:pl-no-uniform} shows that the
single-sentence form of that claim fails, as does the claim's entailment form
(\S\ref{sec:pillars-int}).  Classical reflection follows the same
pattern.  Kreisel and L\'evy's local reflection principle for a theory
consists of the sentences $\PrC(\quo{\varphi}) \limp \varphi$, one for
each sentence~$\varphi$~\cite{kreisel68reflection}.  The principle's instance at
$\z = \tofn 1$ yields the theory's consistency, so by the second
theorem a consistent theory that meets the derivability conditions
does not prove that instance.  In \rga even the entailment form of
that instance, $\PrC(\bquo{\z = \tofn 1}) \vdash \z = \tofn 1$, is
equivalent to $\vdash \ConR$ and fails
([\citen{goedel-paper}, theorem~\textsf{formalized\_loeb\_ent}]).  Beklemishev surveys these reflection
principles~\cite{beklemishev05reflection}.  Second, the grounded truths of Kripke's
least fixed point are $\Pi^1_1$-complete~\cite{burgess86truth},
whereas \rga's grounded truth is recursively enumerable by design,
which is what lets a compiled decider define it.  Reinhardt's proposal
to use KF instrumentally and assert only what it proves
true~\cite{reinhardt86remarks,castaldo23reinhardt} yields a recursively
enumerable, gappy set of significant sentences; \rga's grounded truth
is a set of that kind, given directly.  Third, McGee's
theorem~\cite{mcgee85truthlike} and the Friedman--Sheard
analysis~\cite{friedman87axiomatic} show how easily truth principles
yield $\omega$-inconsistency or outright inconsistency.  \rga's
witness property and its consistency make it $\omega$-consistent
([\citen{goedel-paper}, \S\,``G\"odel's first theorem'']); its $\omega$-incompleteness
(Theorem~\ref{thm:br-omega}) is the expected counterpart.

\paragraph{Truth over applicative theories.}
Untyped $\lambda$-calculus and combinatory logic with a partial truth
predicate are the closest relatives of \rga's slogan that terms are
programs and truth is grounded evaluation.  Aczel's Frege structures
add predicates for proposition and truth to a model of the untyped
$\lambda$-calculus~\cite{aczel80frege}.  Feferman's type-free theories
of partial predicates and partial operations~\cite{feferman84toward},
Cantini's axiomatic theories of truth and
abstraction~\cite{cantini96logical}, and the truth theories over
applicative theories surveyed by Kahle~\cite{kahle01truth} develop the
idea axiomatically.  In these theories truth is a predicate added to
an applicative base and governed by axioms.  In \rga truth is not
added: a term is true when its evaluation yields $1$, and the predicate
$\Tru$ is defined inside \rga, from a compiled decider of that
evaluation (Section~\ref{sec:inside}).

\paragraph{Consistency of fragments.}
Classically, Peano arithmetic proves the consistency of each fragment
$\mathrm{I}\Sigma_n$~\cite{hajek93metamathematics}, and Gentzen proved
the consistency of Peano arithmetic by transfinite induction up to
$\varepsilon_0$~\cite{gentzen36widerspruchsfreiheit}.  The two
internal results of Section~\ref{sec:con} resemble the first of these classical results: a
theory proves, outright, the consistency of a restricted class of its
own proofs, here decided classes of \rga proofs whose strength is
open.  Where the analogy stops shows where \rga is weaker.  The
level theorems hold only under the closed hypothesis $\LF$, which \rga
neither proves nor refutes, while Peano arithmetic needs no hypothesis
for $\mathrm{I}\Sigma_n$; and \rga has no excluded middle at the
universal families that the level hypothesis is about.
\S\ref{sec:con-summary} shows that \rga proves no Gentzen-style
consistency proof of the shape described there: for no reduction that
\rga proves correct does \rga also prove that the reduction's descent stops.

\section{Conclusion}
\label{sec:concl}

\rga set out to be a logic in which everything is computation.  This
paper shows that such a logic can do the metamathematics of itself.
A fuel-indexed decider of \rga's semantics and checkers of \rga's
proofs are compiled primitive-recursive machines and \rga terms, and
\rga proves internally that the decider never affirms $\z = \tofn 1$,
that whatever well-formed term the decider affirms has a proof
certificate, and that two decided classes of \rga's own proofs are
consistent, the second class
containing the first.  For a third decided class, under the
certified-depth condition on the class's quantified case hypotheses
and critical formulas, \rga derives the same consistency sentence
under one closed hypothesis that \rga neither proves nor refutes, and
\oga, \rga with one more rule, proves that hypothesis and with that
hypothesis the sentence.  Around the machines, \rga proves, for each of its proofs,
that the proof's conclusion is true, and for each grounded-true
well-formed term, that the term is provable.

The truth predicate is also a medium for proof by reflection, of the
kind that adds nothing to the system.  \rga derives T-In and T-Out as
rules in every context at the images of quoted quantifier-free
templates, and proves one sentence of reflection for a normalizer of
polynomials in one variable, from which each such polynomial identity
follows at a free variable by computation, with no reflection rule
(\S\ref{sec:pillars-reflect}).

Tarski's theorem takes the form that the material T-biconditional at
\rga's Liar is explosive (\S\ref{sec:pillars-tarski}).  G\"odel's
theorems for \rga and its extensions, and what they say about the
consistency of \rga, are the subject of the G\"odel
paper~\cite{goedel-paper}.

\paragraph{Future work.}
Several questions remain open.  These questions fall into five groups.

\medskip\noindent\emph{Consistency and strength.}\par\nopagebreak
\begin{itemize}
\item \emph{Consistency of \pga inside \rga.}  The decided classes of
  Section~\ref{sec:con} share the propositional, equational and
  inductive core of \pga~\cite{ford25have}.  An internal proof that
  \rga proves \pga consistent is not ruled out by G\"odel's theorem,
  and is a substantial next step.
\item \emph{Strength.}  The proof-theoretic strength of the decided
  classes of Section~\ref{sec:con} (the $\forall$-free class, $\Frag$
  and $\LvC$) is open; the induction of the first two classes is restricted to
  $\forall$-free formulas.  Whether the $\forall$-free class interprets
  the quantifier-free equations of System~T, and whether \rga proves
  the consistency of primitive recursive arithmetic, are open.  \rga
  proves Ackermann's function total by a quantified theorem, and no reduction that
  \rga proves correct has a descent that \rga proves to stop
  (\S\ref{sec:con-summary}); locating \rga's proof-theoretic ordinal
  would make both facts precise.  Whether \rga derives that no-descent
  sentence from a transfinite-induction schema, whose implications
  would have undecided antecedents, is open.
\item \emph{The residual hypothesis.}  The residual result of
  \S\ref{sec:con-res} derives $\neg\Chk(p, \bquo{\z = \tofn 1})$ at
  every proof code under a residual hypothesis.  \rga derives the
  residual at every code in the fragment without swap nodes, but not
  at an arbitrary code, and not at the code of one genuine proof in
  the fragment, where the residual's swap half fails
  (Corollary~\ref{thm:con-res-underivable} and
  Theorem~\ref{thm:con-res-not-frag}).  The residual's guard half holds
  on the fragment, but \rga does not derive that guard half at every
  code
  (Theorem~\ref{thm:con-guard-not}).  Neither the level formulas nor
  any set of hypotheses true in the $\omega$-semantics makes the
  residual derivable on the level classes, and at the code of that
  genuine proof the residual is not even true
  (Theorem~\ref{thm:con-res-lf}).  No extension of \rga sound for the
  $\omega$-semantics derives any residual that the fold accepts at the
  code of an accepted proof with the ungrounded case hypothesis $W\cdot
  W$, since the fold's motive is not true there
  (Theorem~\ref{thm:con-res-ext}); whether some other internal-soundness
  argument avoids that motive is open (Remark~\ref{rem:con-any-motive}).
  Which genuine proofs satisfy the
  residual semantically is not characterized, and whether the fold's
  swap interface fails at a quantified critical formula, as its guard
  interface does at a quantified case hypothesis, is open
  (\S\ref{sec:con-break}).
\item \emph{The level classes.}  Whether \rga derives the consistency
  sentence over $\LvC$ outright, without the hypothesis $\LF$, is
  open (\S\ref{sec:con-lv}); by Theorem~\ref{thm:con-lf-indep}, \rga
  does so exactly when \rga derives the implication from $\LF$ to that sentence.
  The same holds per level, already at the lowest level, $N = 1$.  At
  closed codes every code in a per-level class $\LvC_N$ lies in $\LvC$
  (Theorem~\ref{thm:con-lv-closed}); whether \rga derives that
  containment as one sentence, $\forall p.\ \LvC_N(p) \limp \LvC(p)$,
  is open.
\item \emph{\oga and $\ConR$.}  Whether \oga proves $\ConR$, the
  consistency sentence for all of \rga's proofs, is open; such a proof
  would need the families at every position at once, not only at certified
  ones.  \oga does decide $\ConR$: it proves $\ConR \lor \neg\ConR$,
  since its proof checker's verdict at every number is decided and
  \oga's universal rule gathers those verdicts (\textsf{oga\_ConRGA\_decided},
  session \textsf{OGA\_ConDec}; a classical theorem), and \oga does not
  refute $\ConR$, by \oga's $\omega$-soundness
  (\textsf{oga\_not\_refutes\_ConRGA}).  Systems with a reflection rule do derive $\ConR$.
  Adoption, a rule keyed to a system's own checker, gives consistent
  extensions of \rga and of IGA, the inductive grounded
  arithmetic~\cite{iga-paper}, that derive $\ConR$
  {}[\citen{goedel-paper}, Section ``Adoption: the reflection rule, and the consistency family''].  Whether IGA alone derives $\ConR$ is
  open, and the G\"odel paper names what each route to an answer would
  have to supply [\citen{goedel-paper}, Section ``Does IGA alone prove \rga's consistency?''].
\end{itemize}

\medskip\noindent\emph{Provability.}\par\nopagebreak
\begin{itemize}
\item \emph{G\"odel's theorems.}  The G\"odel paper~\cite{goedel-paper}
  takes up the questions about G\"odel's theorems for \rga and its
  extensions: the second theorem for extensions of \rga, uniform proof
  transformers, and the provability logic of $\PrC$.
\item \emph{Open entailments about the search.}  Whether the search for
  a proof of $\z = \tofn 1$, if the search stops, yields such a proof inside
  \rga, and whether $\ConR$ entails that the search does not stop, are
  open (\S\ref{sec:con-summary}).
\item \emph{Internal completeness with a transformer.}  The internal
  completeness sentence (Theorem~\ref{thm:pl-int}) has an existential
  consequent; whether that sentence has a form with an explicit
  primitive-recursive proof builder from the decider's table, like the
  decided D3, is open.
\end{itemize}

\medskip\noindent\emph{Truth.}\par\nopagebreak
\begin{itemize}
\item \emph{Rules of truth.}  Whether \rga proves $\neg\Tru(\bquo{\neg
  L_{\Tru}})$, that the negated truth-Liar is not established, is open.
  Most of the open cells below are expected to fail by G\"odel's second
  theorem.
\item \emph{Open cells of the truth-axiom table} (Table~\ref{tab:truth-axioms},
  Appendix~\ref{app:pillars-axioms}).  The clause for
  $\neg\forall$ in every context and as an entailment; the rules marked
  $\dagger$ beyond $\forall$-free codes where no failure is known, with
  fuel monotonicity of $\Aff$ at codes of universals; T-Rep and T-Del in
  every context and as entailments; the internal forms of T-Del and of
  KF's clauses for $\neg\forall$ and $\neg T$; T-Out, TB and KF's clause
  for $\neg T$ at grounded-false sentences with quantifiers other than
  refuted universals, such as grounded-false existentials; T-Out at
  ungrounded sentences other than $L_{\Tru}$ and $\Omega$, such as
  $\ConR$, as an implication and as an entailment; the internal
  head-step clause; and the uniform T-biconditionals beyond the
  equations between numerals.  \rga refutes the truth of no universal
  (Theorem~\ref{thm:pl-refuted}(3)), so the refutations of $\Tru$ that
  fail at grounded-false universals fail for this reason; under the
  hypothesis $\ConR$ every refuted sentence is provably unestablished.
  Two cells have entries from \S\ref{sec:pillars-reflect}.  With side
  formulas, as in PKF, T-In and T-Out hold as rules in every context at
  the images of quoted quantifier-free templates, at an arbitrary term
  (Theorem~\ref{thm:pl-rf-tin}), in the row NEC, CONEC.  The
  $\leftarrow$ direction of the clause for $\forall$ fails as a rule in
  every context, at the matrix of $L_{\Tru}$
  (Theorem~\ref{thm:c-tin-omega}), and its certificate form holds
  internally (Theorem~\ref{thm:pl-tin-cert}).  Every cell listed above
  stays open.
\item \emph{Reasoning under ungrounded hypotheses.}
  Theorem~\ref{thm:pl-tracking} assumes that \rga proves or refutes
  $t$; whether $\ConR, t \jbool \vdash \Tru(\bquo{t}) \jbool$ holds
  with $t$ only assumed decided is open, and both semantics are blind
  there: grounded soundness says nothing under the hypothesis $\ConR$,
  and in the $\omega$-semantics the conclusion is true.  The natural
  candidate counterexample is the vacuous universal $V =
  \forall x.\,\neg\ConR$, ungrounded in \rga and refuted under
  $\ConR$, at which the question reduces to G\"odel's second theorem
  for \mbox{\rga{}${}+\ConR$}.
\item \emph{Classically transparent truth beside $\Tru$.}  The material
  T-biconditional is explosive at the Liar of $\Tru$ and at the diagonal
  Liar of every closed predicate term (Theorems~\ref{thm:pl-tb-explosive}
  and~\ref{thm:pl-tb-generic}).  Whether a non-material
  conditional in Field's style could carry the T-biconditionals, and
  whether a primitive truth predicate with PKF's rules T-Intro and
  T-Elim, interpreted for instance by a Kripke construction over the
  $\omega$-semantics, could sit consistently beside $\Tru$, are open.
\item \emph{Reflection across a quantifier level.}  The reflected
  class of \S\ref{sec:pillars-reflect} is quantifier-free, and read
  classically the reflection construction is the partial truth
  definition for that class at work.  At a universal, \rga proves that
  truth is the existence of an accepted proof of the open instance
  (Theorem~\ref{thm:pl-tin-cert}).  An explicit builder discharges
  decided guards inside the proof it builds, so \rga proves a sentence
  one quantifier level above the reflection sentence, with the guards
  outside the inner universal (Theorem~\ref{thm:pl-tin-disch}).  Three steps remain open.  The
  first is the identity itself inside $\Tru$: a sentence whose
  conclusion is the truth of the universal of the reflected identity,
  not of the universal of the identity's truth at an image.  Such a
  sentence needs a builder that also emits the steps of open T-Out.
  The second step is modus ponens at abstract codes in larger contexts.
  The development's modus ponens uses one-element contexts only
  (Theorem~\ref{thm:c-tin-mp}), and a context of two or more symbolic
  codes needs an internal law that the context normalizer's sort pass
  yields a sorted list.  Finally,
  the certificate side is generic in the checker table
  (Appendix~\ref{app:pillars-reflect}, ``T-In at a universal''), but
  $\Tru$ searches only \rga's certificates, so the development derives
  truth only from certificates at \rga's own table.  A truth predicate
  indexed by the checker table is not defined in the development.
\item \emph{The record-level Liar.}  Whether \rga decides the record
  of the Liar $L_{\Rcd}$, built over that sentence's own value record, is open
  (\S\ref{sec:breaks-liars}), as is which record, if any, yields a
  strengthened Liar inside \rga.  For grounded truth itself the
  revenge situation is $L_{\Tru}$, which the metatheory and the
  $\omega$-semantics classify and \rga does not
  (\S\ref{sec:pillars-tarski}).
\end{itemize}

\medskip\noindent\emph{Entailment and classical reasoning.}\par\nopagebreak
\begin{itemize}
\item \emph{Content of entailments.}  Whether a derivation of $A
  \vdash B$ yields, uniformly, a primitive-recursive bound on the fuel
  of $B$'s grounded truth in terms of the fuel of $A$'s, or one that \rga
  proves correct, is open: a realizability reading of entailment, and
  a cut transformer at the level of fuel and certificates.  (A
  computable bound alone is immediate: by soundness, a search over
  fuel finds $B$'s.)  Whether a discharge rule that avoids contraction
  could admit implications over some undecided hypotheses is open too.
\item \emph{Grounded reasoning as a refinement of classical
  reasoning.}  The negative translation makes Heyting arithmetic prove
  the translation of every theorem of Peano arithmetic, while Heyting arithmetic stays
  more discriminating outside the negative fragment.  This paper
  conjectures an analogue for \rga, a decidedness relativization: if
  Peano arithmetic proves $\varphi$, then \rga derives $\varphi$ under
  hypotheses that the relevant subformulas of $\varphi$ are decided,
  with separation theorems that say which decidedness hypotheses \rga
  can discharge.  The entailments of G\"odel's proof in the G\"odel
  paper~\cite{goedel-paper}, and the level theorem's hypothesis $\LF$, are
  instances of the pattern.
\end{itemize}

\medskip\noindent\emph{Proof objects.}\par\nopagebreak
\begin{itemize}
\item \emph{Certificates and accepted proofs.}  Whether \rga proves, as
  one sentence over all codes, that a certificate of the test $\pi$
  yields an accepted proof for $\Chk$ is open (\S\ref{sec:pillars-int}).
  Every internal theorem of this paper asserts that a derivation
  exists; extracting a concrete \rga proof object for one internal
  theorem, checked independently by closed evaluation of $\Chk$, is
  future work.
\item \emph{Wider normalizers.}  The normalizer of
  \S\ref{sec:pillars-reflect} recognizes only the compiled addition and
  multiplication, and only one variable.  A normalizer for polynomials
  in several variables is future work, as is a variant that emits the
  derivation of each identity, which would follow the LCF pattern
  rather than reflection.
\end{itemize}

\appendix

\section{Notation}
\label{app:notation}

Table~\ref{tab:notation} lists the named objects that the paper uses,
grouped by role.  For each object the table gives the symbol, a one-line meaning,
where the paper defines the object, the Isabelle constant, and how the object
stands with respect to the build-time check of
\S\ref{sec:inside-audit} (column ``audit'').  In that column, \emph{audited} means
that the object is a root of that check or lies in the definitional
closure of a root (the root is named); \emph{opaque} means that the
object is compiled from an index chosen by Hilbert choice, and the
check requires that the object be found; \emph{HOL} means that the object is a function, relation or
notation of the metalanguage, for which the question does not arise;
``over $\Chk$'' marks a formula with no constant of its own, built from
the audited $\Chk$; and --- marks an \rga notation outside the check's
subject (numerals, code numerals, typing abbreviations, the Liar
and $\Omega$); \emph{not checked} marks an object introduced in
\S\ref{sec:pillars-reflect}, whose results are not roots of that
check.
Isabelle names are those of the theory text; the theory file of each constant is
in the appendix that the ``defined'' column points to.

{\footnotesize
\setlength{\tabcolsep}{3pt}
\renewcommand{\arraystretch}{1.12}
\begin{longtable}{@{}>{\raggedright\arraybackslash}p{0.15\linewidth}>{\raggedright\arraybackslash}p{0.35\linewidth}>{\raggedright\arraybackslash}p{0.12\linewidth}>{\raggedright\arraybackslash}p{0.20\linewidth}>{\raggedright\arraybackslash}p{0.12\linewidth}@{}}
\caption{The paper's named machines, predicates and classes.}
\label{tab:notation}\\
\hline
symbol & meaning & defined & Isabelle & audit \\
\hline
\endfirsthead
\hline
symbol & meaning & defined & Isabelle & audit \\
\hline
\endhead
\hline
\endfoot
\multicolumn{5}{@{}l}{\textit{Coding and compiling}} \\
$\quo{t}$ & the code of $t$, a natural number (plain corners) &
  \S\ref{sec:inside-codes} & \textsf{encode}, written $\{t\}$ (\textsf{trm\_encode} on terms) & HOL \\
$\bquo{t}$ & the code numeral of $t$, the closed \rga term
  $\tofn{\quo{t}}$ (underlined corners, like any numeral) & \S\ref{sec:inside-codes} &
  \textsf{iCeil} & --- \\
$\tofn n$ & the numeral $\suc^n\z$ & \S\ref{sec:inside-codes} &
  \textsf{trm\_num} & --- \\
$\cpair{x}{y}$ & Cantor pair of two naturals & \S\ref{sec:inside-codes}
  & \textsf{cpair} & HOL \\
$\iSub$, $\iRec$ & composition and recursion nodes of a PR index tree &
  \S\ref{sec:inside-cmp} & \textsf{prSub}, \textsf{prRec} & HOL \\
$\den{f}$ & the binary function that the index $f$ denotes &
  \S\ref{sec:inside-cmp} & \textsf{pr2fi} & HOL \\
$\cmp(f)$ & the compiler: the closed \rga term for the index $f$ &
  \S\ref{sec:inside-cmp} & \textsf{cmp} & HOL \\
$\Num$ & numeral-code encoder, $\Num\cdot\tofn n = \bquo{\tofn n}$,
  compiled from the tree $\tau_{\mathrm{num}}$ & \S\ref{sec:inside-d3} &
  \textsf{iNumT}; \textsf{num\_codeIdx} & audited \\
$\Subst(z, c)$ & compiled simultaneous substitution on codes; $\zeta$
  is a HOL substitution list and $\sof{\zeta}{c}$ its application &
  \S\ref{sec:pillars-sound} & \textsf{iSubM} & audited \\
$\Diag(x, k)$ & compiled diagonal substitution: the code of $t$ with
  $\tofn k$ for $\vz$, at $x = \bquo{t}$ & \S\ref{sec:pillars-tarski},
  {}[\citen{goedel-paper}, \S\,``The sentences''] & \textsf{iDiagSub};
  \textsf{diagSubIdx} & audited (root; also in $\sigma$'s tree) \\
$\mathit{neg}$, $\mathit{or}$, $\mathit{eq}$, $\mathit{ap}$,
  $\mathit{succ}$ & compiled code constructors: the code of a negation,
  disjunction, equation, application or successor from the codes of its
  parts & {}[\citen{goedel-paper}, Appendix ``G\"odel II in full''] & \textsf{mkNegS}, \textsf{mkOrS},
  \textsf{mkEqS}, \textsf{mkAppS}, \textsf{mkSucS} & audited \\
$\mathit{tr}(x)$ & $\Diag(\bquo{\Tru(\vz)}, x)$, the code of $\Tru$
  applied to the numeral $x$ & \S\ref{sec:pillars-axioms} & no
  constant (\textsf{iDiagSub}) & audited (\textsf{iDiagSub}) \\
$\mathit{fst}$, $\mathit{snd}$ & compiled Cantor projections &
  \S\ref{sec:inside-truth}, \S\ref{sec:pillars-int} & \textsf{cpFst},
  \textsf{cpSnd} & $\mathit{snd}$: audited (in $\Dec$) \\
$\Wfc(c)$ & compiled 0/1 test: $c$ codes a well-formed term &
  \S\ref{sec:pillars-int} & \textsf{iWft} & audited (in $\Res$) \\
$\IsCode(c)$ & compiled 0/1 test: $c$ is the code of a term &
  \S\ref{sec:pillars-int} & \textsf{iIsCode} & audited (in $\Res$) \\[0.6ex]
\multicolumn{5}{@{}l}{\textit{Semantics, deciders and truth}} \\
$t \wbev{s} r$ & witnessed evaluation: within fuel $s$ (the companion
  paper's step index), $t$ has the
  value $r$; grounded-true when $r = \ctrue$ for some $s$ &
  \S\ref{sec:bg-rga}, \S\ref{sec:inside-dec} & \textsf{wosr} & HOL \\
$t \whev{s} h$ & head evaluation: within fuel $s$, $t$ has the head
  form $h$ & \S\ref{sec:inside-dec} & \textsf{whev} & HOL \\
$\ctrue$, $\cfalse$ & the truth values, the numerals 1 and 0 &
  \S\ref{sec:inside-dec} & \textsf{encode True}, \textsf{encode False} &
  HOL \\
$c \leadsto c'$ & $c$ head-steps to $c'$ & Fig.~\ref{fig:bg-rules} &
  \textsf{hstep} & HOL \\
$\decV_s(t)$ & the value decider: $\None$ or $\Some r$ at fuel $s$ (with
  a head decider beside it) & \S\ref{sec:inside-dec} & \textsf{prDecV},
  \textsf{prDecH} & HOL \\
$\decP_s(x)$ & the packed verdict: 0 for $\None$, $r+1$ for $\Some r$ &
  \S\ref{sec:inside-dec} & \textsf{prDecVEnc} & HOL \\
$\dtab(s, x)$ & the PR table machine; $B(s, x)$ is its code bound &
  \S\ref{sec:inside-dec} & \textsf{prDecTab}, \textsf{prBound} & HOL \\
$\taut$ & the explicit index tree of $\dtab$ & \S\ref{sec:inside-truth}
  & \textsf{prDecTabIdx} & audited (in $\Dec$) \\
$\Dec(s, c)$ & the compiled decider, $\mathit{snd}(\cmp(\taut)\cdot
  s\cdot c)$ & \S\ref{sec:inside-truth} & \textsf{iPrDecV} & audited \\
$\Aff(s, c)$ & affirmation test: 1 when $\Dec(s, c) = \tofn 2$, else 0
  & Def.~\ref{def:in-truth} & \textsf{iPrDecTrueAt} (arguments $c, s$) &
  audited (in $\Tru$) \\
$\Tru(c)$ & the truth predicate, $\exists s.\ \Aff(s, c)$ &
  Def.~\ref{def:in-truth} & \textsf{iWtrue} & audited \\
(no symbol) & the opaque twins of $\Dec$ and $\Num$ &
  \S\ref{sec:inside-cmp}, \S\ref{sec:inside-d3} &
  \textsf{iPrDecOpqT}, \textsf{iNum} & opaque \\
$\Rcd_v(c)$, $\Rcd^s_v(c)$ & the value record, certificate-free: $c$ has
  the value $v$ at some fuel, resp.\ at fuel $s$ &
  Def.~\ref{def:con-rcd} & \textsf{iRec}, \textsf{iRecAt} & audited (roots \textsf{RF}, \textsf{recSelf}) \\
$\mathit{Guard}_r$ & the guard of the record fold at rank $r$ &
  \S\ref{sec:con-frag} & \textsf{iGuardAt} & audited (in $\Res$) \\
$\omega$-semantics & the classical semantics: a universal is true
  when every instance is, with no certificate & \S\ref{sec:inside-omega} &
  \textsf{omv} & HOL \\
$\omega$-record & the certificate-free value record of the fold, for
  $\forall$-free terms & Def.~\ref{def:c-orec} & \textsf{orec} &
  HOL \\[0.6ex]
\multicolumn{5}{@{}l}{\textit{Proofs and checkers}} \\
$\Gamma \vdash c$ & \rga derives $c$ from the hypotheses $\Gamma$ &
  \S\ref{sec:bg-rga} & \textsf{entails} & HOL \\
$a, \Gamma$ & hypothesis insertion: the finite set $\Gamma$ with the
  hypothesis $a$ added; Isabelle, and the statements of
  \S\ref{app:endpoints-statements}, write $a \triangleright \Gamma$ &
  Fig.~\ref{fig:bg-rules} & \textsf{finsert}, written $\triangleright$
  & HOL \\
$\Gamma \vdashO c$ & the same in \oga & \S\ref{sec:con-lv} &
  \textsf{oent} & HOL \\
$A \vdash B$, $A \eqv B$ & entailment: \rga derives $B$ under the
  hypothesis $A$; interderivability, $A \vdash B$ and $B \vdash A$ (not
  intersubstitutivity) & \S\ref{sec:inside-codes} & \textsf{entails}
  ($A \triangleright \emptyset \vdash B$) & HOL \\
$\exists s.\, D(s)$ & $\neg\forall s.\, \neg D(s)$; with $D$ decided, a
  hypothesis of the $\Sigma$ shape (\S\ref{sec:inside-codes}) &
  \S\ref{sec:inside-codes} & $\exists$ of \texttt{Cov\_Induct};
  \textsf{iPrc} for $\PrC$ & --- \\
$\forall s.\,(D(s) \limp B)$ & the decided-antecedent form of a
  conditional from the hypothesis $\exists s.\,D(s)$, $D$ decided: one
  sentence, equivalent to the entailment $\exists s.\,D(s) \vdash B$
  (Theorem~\ref{thm:in-dt}) & \S\ref{sec:inside-codes} & (no constant;
  the Isabelle names say ``prenex'') & --- \\
$\bot$ & the divergent constant, a term with no value
  (Table~\ref{tab:in-sem}); the false sentence of the paper is $\z =
  \tofn 1$ & \S\ref{sec:bg-rga} & \textsf{roBot} & --- \\
$P \proves H \vdash c$ & $P$ is a valid proof whose conclusion is
  $H \vdash c$ (read ``$P$ proves $H \vdash c$'') & \S\ref{sec:inside-dec} &
  \textsf{proves} & HOL \\
$a\jnat$, $a\jbool$ & $a = a$ ($a$ is a natural); $a \lor \neg a$ ($a$
  is decided) & \S\ref{sec:inside-codes} & \textsf{trmNat},
  \textsf{trmTypeBool} & --- \\
$\mathit{TreeValid}(p)$ & compiled validity test: every entry of the
  list $p$ is a judgment code that a rule detector accepts &
  \S\ref{sec:inside-chk} & \textsf{iTreeValid} & audited (in $\Chk$) \\
$\Chk(p, c)$ & the tree checker: $p$ codes a valid proof of $\emptyset
  \vdash c$ & Def.~\ref{def:in-chk} & \textsf{iChk} & audited \\
$\ChkQf(p, c)$ & the $\forall$-free checker, $\Chk$ with the scan
  $\QfScan$ & Def.~\ref{def:con-chkqf} & \textsf{iChkQf} & audited \\
$\QfScan(p)$ & every judgment of the list $p$ is $\forall$-free &
  \S\ref{sec:con-qf} & \textsf{qfScanAcc} & audited \\
$\qfs(c)$ & $\forall$-freeness test of one code, 0 meaning ``pass'' &
  \S\ref{sec:con-qf} & \textsf{iQfScan} & audited \\
$\ChkE(p, c)$ & the companion paper's PR proof checker, compiled
  through a Hilbert-choice index & \S\ref{sec:pillars-complete} &
  \textsf{iChkOpq}; \textsf{chkIdx} & opaque \\
$\pi$ & the explicit index of ``$p$ codes a valid proof whose
  conclusion judgment has code $h$'' & Def.~\ref{def:pl-cert} &
  \textsf{provesChkIdx} & audited (in $\Cert$) \\[0.6ex]
\multicolumn{5}{@{}l}{\textit{Provability}} \\
$\PrC(x)$ & \rga's provability formula, $\neg\forall q.\ \neg\Chk(q,
  x)$ & \S\ref{sec:pillars-complete} & no constant & over $\Chk$ \\
$\Prv(c)$ & the provability predicate over $\ChkE$, $\exists p.\ \ChkE(p,
  c)$ & Def.~\ref{def:pl-prv} & \textsf{iPrv} & opaque \\
$\Cert(h)$, $\Cert_\emptyset(c)$ & certificate predicate, $\exists p.\
  \cmp(\pi)\cdot h\cdot p \ne \z$; at the code of $\emptyset \vdash c$ &
  Def.~\ref{def:pl-cert} & \textsf{iHasCert}; \textsf{emptyJdg} &
  audited \\[0.6ex]
\multicolumn{5}{@{}l}{\textit{Consistency and the classes of Section~\ref{sec:con}}} \\
$\ConR$ & the consistency sentence, $\forall p.\ \neg\Chk(p, \bquo{\z =
  \tofn 1})$ & \S\ref{sec:inside-chk} & \textsf{ConRGA} & audited \\
$\Frag(p)$ & the fragment: a decided scan of the proof code &
  Def.~\ref{def:con-frag} & \textsf{iFrag} & audited \\
$\Dep(c)$ & the certified depth: $0$ if $c$ is not certified;
  otherwise the length of the longest chain of quantifiers nested
  through head steps, and $1$ if $c$ has none &
  \S\ref{sec:con-lv} & \textsf{iDep} & audited \\
$\Ok(c)$, $\Ok(n, c)$ & $c$ is certified ($\Dep(c) \ne \z$, a
  well-formed term code); certified at level $n$ ($\Dep(c) \le n+1$) &
  \S\ref{sec:con-lv} & \textsf{OkI}, \textsf{OkU} & audited \\
$\LvC(p)$ & the class of proof codes whose case hypotheses and swap
  targets are certified & Def.~\ref{def:con-lv} & \textsf{iLvClassI} &
  audited \\
$\LF$ & closed sentence: the numeral families are decided at every
  certified quantifier application & Def.~\ref{def:con-lf} &
  \textsf{iLFI} & audited \\
$\LF(n)$, $\LF_{\forall}$ & the same at the codes certified at level
  $n$; $\LF_{\forall} = \forall n.\ \LF(n)$ & Def.~\ref{def:con-lf} &
  \textsf{iLFU}, \textsf{iLFA} & audited \\
$\LvC(n, p)$ & the class at a level $n$ given as an \rga term:
  $\LvC$ with ``certified at level $n$'' in place of ``certified'' &
  \S\ref{sec:con-lv} & \textsf{iLvClassU} & audited \\
$\Stat(p)$ & the static half of $\LvC$: its swap-node and
  substitution clauses & \S\ref{sec:con-lv} & \textsf{iStatLvI},
  \textsf{iSub1Qf} & audited (in $\LvC$) \\
$\LvC_N(p)$ & the per-level class: $\LvC$ with the level-$N$ scan as
  certificate & \S\ref{sec:con-lv} & \textsf{iLvClass} & audited \\
$\LF_N$, $\LF_{\le N}$ & the per-level sentences: $\LF$ with the
  level-$N$ scan as certificate; $\LF_{\le N} = \LF_N \wedge \cdots
  \wedge \LF_0$ & \S\ref{sec:con-lv} & \textsf{iLF}, \textsf{iLFs} &
  audited \\
$\Stat_N(p)$ & the static half of the per-level class &
  \S\ref{sec:con-lv} & \textsf{iStatLv}, \textsf{iSub1Qf} & audited \\
$\DynFam(p)$ & the reach-indexed family hypothesis about the proof code
  $p$ & \S\ref{sec:con-lv} & \textsf{iDynWFD} & audited \\
$\Res(p)$ & the guarded residual: the fold's two interfaces at $p$ &
  Def.~\ref{def:con-res} & \textsf{iRes} & audited \\
$\mathit{ResGuard}(p)$, $\mathit{ResSwap}(p)$ & the guard half and the
  swap half of $\Res(p)$ & \S\ref{sec:con-res} & \textsf{iResGuard},
  \textsf{iResSwap} & audited (in $\Res$) \\
$\mathit{NoSwap}(p)$ & no node of $p$ is read as one of the five swap
  rules & \S\ref{sec:con-res} & \textsf{iNoSwap} & audited \\
$\mathit{HS}(p)$ & the head-step scan: at every universal elimination
  and negated-universal introduction of $p$, the critical formula has a
  head step & \S\ref{sec:con-break} & \textsf{iHdStep} & audited \\
$P_x$ & a valid proof whose case hypothesis $a_x = \forall n.\
  \mathit{ofam}(x)\cdot n$ says that program $x$ never halts &
  \S\ref{sec:con-res} & \textsf{gxP} & HOL \\
$\mathit{fragSearch}$ & search for a code in $\Frag$ that $\Chk$ accepts
  for $\z = \tofn 1$ & \S\ref{sec:con-frag} & \textsf{fragSearch} &
  audited \\
$\mathit{conSearch}$ & search through all proof codes for one that
  $\Chk$ accepts for $\z = \tofn 1$ & \S\ref{sec:con-res} &
  \textsf{conSearch} & audited \\
$\mathit{Stops}$ & $\Rcd_1(\bquo{\mathit{conSearch}})$: the search
  returns & {}[\citen{goedel-paper}, \S\,``Consequences''] & \textsf{conStops} (abbreviation)
  & audited (in $\Rcd$) \\[0.6ex]
\multicolumn{5}{@{}l}{\textit{Self-reference and G\"odel's theorems}} \\
$L$ & the ordinary Liar: a term that head-steps in two steps to $\neg
  L$, with no codes & \S\ref{sec:bg-rga}, \S\ref{sec:breaks-liars} &
  \textsf{trmLiar} (\texttt{RGA/RGA\_Paradox}) & --- \\
$\Omega$ & the self-reducing term $(\lambda(\vz\cdot\vz))\cdot
  (\lambda(\vz\cdot\vz))$, which has no value &
  \S\ref{sec:pillars-axioms} & \textsf{trmOmega}
  (\texttt{RGA/RGA\_Paradox}) & --- \\
$D_{\Tru}$ & the Liar template,
  $\neg\Tru(\Diag(\vz, \vz))$ & \S\ref{sec:pillars-tarski} &
  \textsf{lwD} & audited \\
$S_{\Tru}$ & $\Diag(\bquo{D_{\Tru}}, \bquo{D_{\Tru}})$, the code term of
  $L_{\Tru}$ & \S\ref{sec:pillars-tarski} &
  \textsf{goedelWS} & audited \\
$L_{\Tru}$ & the Liar for the truth predicate, $\neg\Tru(S)$ &
  \S\ref{sec:pillars-tarski} & \textsf{goedelW} & audited \\
$L_{\Rcd}$ & the Liar for the value record, ``my own record at the
  value true fails'' & \S\ref{sec:breaks-liars} & \textsf{liarRec} &
  audited \\
$U$ & single-sentence soundness at $\z = \tofn 1$ &
  Thm.~\ref{thm:pl-no-uniform} & (no constant) & over $\Chk$ and $\Tru$ \\
\multicolumn{5}{@{}l}{\textit{Reflection through the truth predicate}} \\
$\mathit{DecVal}(c, x)$ & the value that the decider reads at $c$ is
  $x$: $\exists l.\ \Dec(l, c) = \suc x$ & \S\ref{sec:pillars-reflect}
  & \textsf{iPrDecVal} & not checked \\
$\mathit{env}(v)$ & environment code that maps $\vz$ to the code of
  the numeral of $v$ & \S\ref{sec:pillars-reflect} & \textsf{envVar0}
  (abbreviation over \textsf{updSubM}, \textsf{numC}) & not checked \\
$\mathit{IsQf}(c)$ & compiled test: $c$ codes a $\forall$-free term,
  nonzero meaning ``pass'' (unlike $\qfs$, where 0 means ``pass'') &
  \S\ref{sec:pillars-reflect} & \textsf{iQf} & not checked \\
$\mathit{Norm}(x)$ & the polynomial normalizer: a pair of a tag (0: a
  polynomial term) and a coefficient list & \S\ref{sec:pillars-reflect}
  & \textsf{plNormT} (tree \textsf{plNormIdx}) & not checked \\
$\mathit{polyEq}(x, y)$ & the reflection test: both tags are 0 and the
  trimmed coefficient lists agree & \S\ref{sec:pillars-reflect} &
  \textsf{polyEqT} (tree \textsf{polyEqIdx}) & not checked \\
$\mathit{ev}(v, l)$, $\mathit{sc}$, $\oplus$, $\otimes$,
  $\mathit{trim}$ & Horner evaluation, and the scaling, addition,
  multiplication and trimming machines on coefficient lists &
  \S\ref{sec:pillars-reflect}, \S\ref{app:pillars-reflect} &
  \textsf{evT}; \textsf{plScIdx}, \textsf{plAddIdx}, \textsf{plMulIdx},
  \textsf{plTrimIdx} & not checked \\
$\mathit{Rec}(m, n)$ & compiled value of the course-of-values skeleton
  & \S\ref{app:pillars-reflect} & \textsf{recV} (tree \textsf{recIdx})
  & not checked \\
$J(c)$ & code of the open instance $\{\vz\jnat\} \vdash
  \mathit{ap}(\mathit{lift}(c), \bquo{\vz})$ of the universal with body
  code $c$ & \S\ref{sec:pillars-reflect} & \textsf{openJ} & not checked \\
$\mathit{lift}(c)$ & compiled de Bruijn lift by one at cutoff 0 &
  \S\ref{sec:pillars-reflect} & \textsf{iLiftS} & not checked \\
$u_L$ & the matrix $\lambda\neg\Aff(\vz, S_{\Tru})$ of $L_{\Tru}$ &
  \S\ref{sec:pillars-reflect} & (no constant) & --- \\
\end{longtable}
}

\paragraph{Symbols that name two things.}
A few letters carry two meanings, which the context separates.
\begin{itemize}
\item The templates $D_{\Tru}$ and $D_G$ are not related to the
  conditions D1--D3 or D3G.
\item The italic $S_{\Tru}$ and $S_G$ are code terms; the bold $\suc$
  is the successor.
\item $L$ is the ordinary Liar; $L_{\Tru}$ and $L_{\Rcd}$ are the Liars
  built from codes over $\Tru$ and over the value record.
\item $\PrC$ is the classical provability predicate in
  {}[\citen{goedel-paper}, \S\,``The classical proof of G\"odel's second theorem''] and \rga's provability formula over $\Chk$ from
  \S\ref{sec:pillars-complete} on.
\item $T$ is a theory in [\citen{goedel-paper}, \S\,``The classical proof of G\"odel's second theorem''] and a classical truth
  predicate in Section~\ref{sec:pillars}.
\item $B(s, x)$ is the table's code bound; the sans-serif $\jbool$ in
  $a\jbool$ is decidedness.  Likewise the italic level $N$ of $\LvC_N$
  and $\LF_N$ is a number, and the sans-serif $\jnat$ in $a\jnat$ is
  the natural-number judgment.
\item A subscript level, as in $\LvC_N$, $\LF_N$ and $\Ok_N$, is a
  meta-level number; an argument level, as in $\LF(n)$ and $\Ok(n,
  c)$, is a term of \rga.  $\LF_{\forall}$ quantifies the argument
  level inside \rga.
\item The certified depth is written $\Dep(c)$; the letter $D$ is
  kept for the decided matrix $D(s)$ of a hypothesis $\exists s.\,D(s)$,
  for the conditions D1--D3 and D3G, and for the templates $D_{\Tru}$
  and $D_G$, which are not related to the conditions.
\item $\Sigma$ appears only in $\Sigma_1$, the arithmetical class, and
  in ``the $\Sigma$ deduction theorem'' and ``a hypothesis of the
  $\Sigma$ shape'', which name hypotheses $\exists s.\, D(s)$ with $D$
  decided.  No machine is written $\Sigma$ or $\Pi$.
\item $\bot$ is \rga's divergent constant, never ``false'';
  Beklemishev and Shamkanov's $\bot$ (\S\ref{sec:related}) is the
  classical falsum $\neg(\z = \z)$.
\item $\vdash$ with a hypothesis on its left, $A \vdash B$, is an
  entailment; $\eqv$ is interderivability, not the biconditional
  $\liff$, which is a single sentence.
\item $L_\psi$ is the L\"ob sentence for $\psi$, not a Liar.
\end{itemize}

\section{Coding and machines}
\label{app:machines}
\begingroup\sloppy

This appendix supplies the details that Section~\ref{sec:inside},
on the objects built inside \rga, defers: the theory file of each name, the coding, the compiler's
clauses and side conditions, the decider and its table machine, the
introduction, elimination and transport rules of the truth predicate,
the opaque twins, the checker's construction, the negative results
that explain three design choices, and the build-time check of
\S\ref{sec:inside-audit}.  Statements carry every hypothesis of their
Isabelle form.  Theory files are given relative to \texttt{hol/}.
Names resolve in the development
(Section~\ref{sec:mech}, ``Artifact'').

\subsection{Theory files}
\label{app:machines-files}

Table~\ref{tab:a-files} gives the theory file of every name that
Section~\ref{sec:inside} cites, including the companion-paper facts of
the box in \S\ref{sec:inside-codes}, the semantics of
Table~\ref{tab:in-sem} and the checkers of
Table~\ref{tab:in-checkers}.  The endpoint table of
Appendix~\ref{app:endpoints} lists the same names by theorem.

\begin{table}[ht]
\centering
\footnotesize
\begin{tabular}{@{}>{\raggedright\arraybackslash}p{0.30\linewidth}>{\raggedright\arraybackslash}p{0.66\linewidth}@{}}
\hline
theory file & names \\
\hline
\texttt{RGA/RGA\_Compile} & \textsf{cmp}, \textsf{cmp\_correct} \\
\texttt{RGA/RGA\_SystemT\_PR} & \textsf{cmp\_eq\_num}, \textsf{cmp\_prSub\_unfold} \\
\texttt{RGA/RGA\_Total} & \textsf{cmp\_total}, \textsf{cmp\_total\_app} \\
\texttt{RGA/RGA\_Cantor\_Internal} & \textsf{cmp\_prRec\_zero}, \textsf{cmp\_prRec\_suc},
  \textsf{cmp\_prRec\_suc\_any} \\
\texttt{RGA/RGA\_Proof} & \textsf{derives} (the 42 rules), \textsf{ir\_impI1},
  \textsf{ir\_impI2} \\
\texttt{RGA/RGA\_GCalc} & \textsf{ir\_weaken0}, \textsf{wsubf} \\
\texttt{RGA/RGA\_Sound} & \textsf{wosr}, \textsf{whev} (Table~\ref{tab:in-sem}),
  \textsf{ir\_impI\_bool}, \textsf{rga\_provable\_wosr}, \textsf{rga\_consistent},
  \textsf{rga\_open\_iff}, \textsf{rga\_val\_iff}, \textsf{wosr\_unique},
  \textsf{ir\_disj\_prop}, \textsf{prDecV\_allT\_family},
  \textsf{prDecV\_allF\_family} \\
\texttt{RGA/RGA\_Complete} & \textsf{rga\_complete\_open} \\
\texttt{RGA/RGA\_Paradox} & \textsf{rga\_deduction\_theorem\_fails} \\
\texttt{RGA/RGA\_Decide} & \textsf{prDecH}, \textsf{prDecV}, \textsf{prDecV\_sound},
  \textsf{prDecV\_complete}, \textsf{prDecH\_complete}, \textsf{prDecV\_mono},
  \textsf{prDecVEnc} \\
\texttt{RGA\_Decider/RGA\_Decide\_PR} & \textsf{prDecTab}, \textsf{n2pr\_prDecTab},
  \textsf{prDecTab\_correct}, \textsf{prDecIdx}, \textsf{iPrDecVOpq\_closed},
  \textsf{allT\_example\_proof}, \textsf{prDecV\_allT\_example},
  \textsf{prDecV\_allT\_fuel\_lower} \\
\texttt{RGA\_Kit/Code\_Index} & \textsf{prDecTabIdx} ($\taut$), \textsf{iPrDecV},
  \textsf{iPrDecTrueAt}, \textsf{iWtrue}, \textsf{hq\_iPrDecV},
  \textsf{iPrDecV\_suc\_unfold}, \textsf{iPrDecV\_suc\_unfold\_any}, \textsf{drvChkList}, \textsf{validChkIdx},
  \textsf{provesChkIdx} ($\pi$) \\
\texttt{RGA\_Decider/Decider\_Truth} & \textsf{pr2fi\_prDecTabIdx},
  \textsf{iPrDecV\_closed}, \textsf{iPrDecV\_agree}, \textsf{iWtrueI},
  \textsf{iWtrueI\_abs}, \textsf{rga\_wtrueE}, \textsf{iWtrue\_transport},
  \textsf{iPrDecTrueAt\_cEq01\_never}, \textsf{iNumT}, \textsf{iNumT\_closed},
  \textsf{iNum}, \textsf{iNum\_agree} \\
\texttt{RGA\_Decider/Chk\_Index} & \textsf{drvChk\_correct}, \textsf{validChk\_correct},
  \textsf{provesChk\_correct} \\
\texttt{RGA\_Checker/Chk\_Base} & \textsf{validChk\_genuine} \\
\texttt{RGA\_Checker/Chk\_Formers} & \textsf{iTreeValid} \\
\texttt{RGA\_Checker/Chk\_Faithful} & \textsf{iChk}, \textsf{iChk\_faithful},
  \textsf{iChk\_genuine}, \textsf{iChk\_genuine\_proves}, \textsf{hq\_iChk},
  \textsf{iTreeValid\_faithful}, \textsf{iTreeValid\_genuine},
  \textsf{iChkOpq\_iChk\_iff}, \textsf{iChkOpq\_iChk\_agree} \\
\texttt{RGA\_Checker/Chk\_QF\_Faithful} & \textsf{iChkQf}, \textsf{iChkQf\_faithful},
  \textsf{iChkQf\_genuine\_proves} \\
\texttt{RGA\_Kit/Term\_Coders} & \textsf{cEq01}, \textsf{chkIdx}, \textsf{iChkOpq},
  \textsf{iChkOpq\_true}, \textsf{iChkOpq\_iff\_proves}, \textsf{iChkOpq\_decodes},
  \textsf{iChkOpq\_genuine\_iff}, \textsf{prf\_proof\_noncode} \\
\texttt{RGA\_Checker/RGA\_Deriv} & \textsf{D3\_closed}, \textsf{D3\_monotone} \\
\texttt{RGA\_Goedel/Goedel\_II} & \textsf{Choice\_Audit}, the walk of the
  build-time check of \S\ref{app:machines-audit}; \textsf{Dep\_Audit}, the
  dependency checks \\
\texttt{RGA\_Goedel/Goedel\_II\_PR} & \textsf{goedelSigma\_is\_pr}; the first run of
  the build-time check; the G\"odel-path dependency check; the index-instance
  audit \textsf{Choice\_Audit\_WL} and its control \textsf{wlCtlIdx}, and the
  unspecified-constant audit \textsf{Choice\_Audit\_UD} and its control
  \textsf{udCtlIdx} \\
\texttt{RGA\_Goedel/Goedel\_D3} & the fourth run of the build-time check \\
\texttt{RGA\_Goedel/Goedel\_D3\_Uniform} & the fifth run of the
  build-time check \\
\texttt{RGA\_Goedel/Goedel\_D2} & \textsf{d2\_valid\_accepted},
  \textsf{D2\_decided}; the sixth run of the build-time check \\
\texttt{RGA\_Goedel/Goedel\_Cut} & the seventh and eighth runs of the
  build-time check \\
\texttt{ConOGA/ConOGA\_Traps} & the second run of the build-time check \\
\texttt{ConOGA/ConOGA\_Audit} & the third run of the build-time check;
  the whitelist and unspecified-constant audits of the later runs; the
  soundness dependency check and the semantic-fold check \\
\texttt{RGA\_Decider/RGA\_Decider\_Traps}, \texttt{RGA\_Checker/RGA\_Checker\_Traps}
  & the negative results of \S\ref{app:machines-traps} \\
\hline
\end{tabular}
\caption{Theory files of the names cited in Section~\ref{sec:inside}.}
\label{tab:a-files}
\end{table}

\subsection{Coding}
\label{app:machines-coding}

The term coding is \textsf{trm\_encode} (\texttt{Syntax\_Term}): a
variable $v_i$ has code $\cpair{0}{i}$, and an operator node with
operator number $o$, binder count $b$ and argument list $ts$ has code
$\cpair{1}{\cpair{o}{\cpair{b}{\quo{ts}}}}$, where $\quo{ts}$ codes the
list of the arguments' codes.  The pair is the Cantor pair
(\texttt{Cantor\_Pair}).  Lists are coded by \textsf{encode\_list}
(\texttt{Coding}): the empty list has code $0$, and $x$ in front of
$l$ has code $\cpair{\quo{x}}{\quo{l}} + 1$.  A finite set of terms,
such as the hypothesis set of a judgment, is coded as the sorted list
of its elements' codes (\textsf{encode\_fset\_def}).  The coding is
injective (\textsf{trm\_encode\_inj}) and not onto.  A number whose
first Cantor component is $2$ codes no term
(\texttt{RGA\_Kit/Code\_Total}), and likewise not every number codes a
proof, which is why the decider and the checker test that a number is
in range, $\quo{\mathrm{decode}\,p} = p$, where $\mathrm{decode}$
inverts the coding on its range.

The underlined corners are the abbreviation \textsf{iCeil}
(\texttt{RGA\_Kit/Cov\_Induct}): $\bquo{t}$ is the numeral of
$\quo{t}$.  The code numeral of the false equation, $\bquo{\z = \tofn
1}$, is the constant \textsf{cEq01}, and the equation's plain code $\quo{\z =
\tofn 1}$ is \textsf{sEq01} (both in \texttt{RGA\_Kit/Term\_Coders}).

\subsection{The compiler}
\label{app:machines-cmp}

\paragraph{Index trees.}
Binary PR functions are described by index trees of type
\textsf{pridx} (\texttt{Indexed\_Recursion}).  The denotation
$\den{f}$ is \textsf{pr2fi}~$f$, and \textsf{n2pr}~$F$ says that a
binary function $F$ is primitive recursive.  The operator
\textsf{priof}~$F$ picks, by Hilbert's choice operator, some index
whose denotation is $F$.  The only fact the development has about
such an index is
\begin{thm}[\textsf{pr2fi\_priof}]
\label{thm:a-priof}
If $\mathit{n2pr}\,F$, then $\den{\mathit{priof}\,F} = F$.
\end{thm}
\noindent
The statement ``$\den{i}$ is PR'' holds of every index $i$
(\textsf{pr2fi\_n2pr}), so a theorem of that form cannot tell an
explicit tree from one chosen by \textsf{priof}.

\paragraph{The compiled terms.}
The compiler \textsf{cmp} (\texttt{RGA/RGA\_Compile}) maps the four
initial functions to $\lambda\lambda.\,\z$, $\lambda\lambda.\,\suc
v_1$, $\lambda\lambda.\,v_1$ and $\lambda\lambda.\,v_0$ (zero,
successor of the first argument, the two projections), and the two
node forms to
\begin{align*}
  \cmp(\iSub\,g\,h_1\,h_2) &= \lambda\lambda.\ \cmp(g)\cdot(\cmp(h_1)\cdot v_1\cdot v_0)\cdot(\cmp(h_2)\cdot v_1\cdot v_0) , \\
  \cmp(\iRec\,g\,h\,h') &= \lambda\lambda.\ \cR\cdot(\cmp(g)\cdot v_1\cdot\z)\cdot\bigl(\lambda\lambda.\ \cmp(h)\cdot(\cmp(h')\cdot v_3\cdot v_1)\cdot v_0\bigr)\cdot v_0 .
\end{align*}
In the step abstraction of the recursion, $v_1$ is the counter, $v_0$
the previous value, and $v_3$ the outer first argument.  Correctness
is stated against the companion paper's evaluation relation $A
\Downarrow_s a$ (Isabelle \textsf{osr}), the semantics whose quantifier
clauses keep the certificate but drop the evaluations of instances
(the family of requirement~(iii) of \S\ref{sec:inside-dec}, and the
refuted instance's evaluation in the falsity clause), and in which the
companion paper's completeness theorem is proved.

\begin{thm}[\textsf{cmp\_correct}]
\label{thm:a-cmp-correct}
Let $A$ and $B$ be closed terms with $A \Downarrow_{s_a} a$ and $B
\Downarrow_{s_b} b$.  Then for every index $f$ there is an $s$ with
$\cmp(f)\cdot A\cdot B \Downarrow_s \den{f}(a, b)$.
\end{thm}

\noindent
\emph{Proof sketch.}  Induction on $f$.  The composition case runs the
two inner calls and then $g$; the recursion case evaluates the counter
to a numeral and unrolls $\cR$ by an inner induction on the counter's value.
\qed

\medskip
Theorem~\ref{thm:in-cmp-num} is this theorem at numerals
(\textsf{cmp\_correct\_num}), followed by the companion paper's value
form of open completeness, which turns an evaluation of a well-formed
term into a derivable equation (\textsf{ir\_osr\_eq}, behind it
\textsf{rga\_complete\_open}).  Totality has a rule form at abstract
arguments, of which Theorem~\ref{thm:in-cmp-num}(2) is the closed
universal instance:

\begin{thm}[\textsf{cmp\_total\_app}]
\label{thm:a-cmp-total-app}
For every index $f$: if $\Gamma \vdash a \jnat$ and $\Gamma \vdash b
\jnat$, then $\Gamma \vdash (\cmp(f)\cdot a\cdot b) \jnat$.
\end{thm}

\paragraph{Unfolding at abstract arguments.}
Theorem~\ref{thm:in-cmp-rec} has the following exact hypotheses.  The
zero law, \textsf{cmp\_prRec\_zero}, needs only $\Gamma \vdash m
\jnat$.  The successor law, \textsf{cmp\_prRec\_suc}, needs $\Gamma
\vdash m \jnat$, $\Gamma \vdash k \jnat$, and a variable index $i$
that does not occur in $m$ and is not free in $\Gamma$; the variable
does not occur in the conclusion.  Such an $i$ always exists, because
$m$ and $\Gamma$ are finite, so the side condition can always be met.
The variant \textsf{cmp\_prRec\_suc\_any} states the law without the
side condition, by choosing $i$ inside the variant's proof.  The composition node unfolds as
well.

\begin{thm}[\textsf{cmp\_prSub\_unfold}]
\label{thm:a-cmp-sub}
Let $m$ and $n$ be well-formed, and suppose $\Gamma \vdash
(\cmp(g)\cdot(\cmp(h_1)\cdot m\cdot n)\cdot(\cmp(h_2)\cdot m\cdot n))
\jnat$.  Then
\[
  \Gamma \vdash \cmp(\iSub\,g\,h_1\,h_2)\cdot m\cdot n
    = \cmp(g)\cdot(\cmp(h_1)\cdot m\cdot n)\cdot(\cmp(h_2)\cdot m\cdot n) .
\]
\end{thm}

\noindent
The totality premise is the one Theorem~\ref{thm:a-cmp-total-app}
supplies whenever $m$ and $n$ are naturals.  All three laws read the
node at the root of the tree, so they apply only to an index whose
root the development can name.

\subsection{The decider and its table machine}
\label{app:machines-dec}

\paragraph{The clauses.}
The head decider \textsf{prDecH} and the value decider
\textsf{prDecV} (\texttt{RGA/RGA\_Decide}) are defined by one mutual
recursion on the fuel.  At fuel $0$ and at a variable both deciders return
$\None$.  At fuel $s+1$ every recursive call uses fuel $s$, and the
clauses follow Table~\ref{tab:in-sem}.  The head decider reports
$\lambda b$, $\cR$ and $\cA$ for those terms, $\cR\cdot u$ and
$\cR\cdot g\cdot u$ for an application whose head is $\cR$ or
$\cR\cdot g$, the head of $b[u]$ for an application whose head is
$\lambda b$, and, for an application $t\cdot m$ whose head is
$\cR\cdot g\cdot h$, the head of $g$ when $m$ has value $0$ and of
$h\cdot\tofn k\cdot(\cR\cdot g\cdot h\cdot\tofn k)$ when $m$ has value
$k+1$.  The value decider computes $0$, successor, equality,
negation and disjunction as in the table (disjunction is parallel: one
true side suffices), follows $\beta$ and the recursor as the head
decider does, and handles an application whose head is $\cA$ by the
bounded certificate search of \S\ref{sec:inside-dec}.  No other
clause applies, so every other term gets $\None$.

\paragraph{The range check.}
At a quantifier node the certificate search tests, for each candidate
$p \le s$, that $\quo{\mathrm{decode}\,p} = p$ before it reads $p$ as
a proof: the true verdict needs some such $p \le s$ that codes a proof
$P$ with $P \proves \{\vz\jnat\} \vdash \lft u \cdot \vz$, and the
false verdict some such $p \le s$ that codes a proof of $\emptyset
\vdash \neg(u\cdot\tofn{i})$ for some $i \le s$.  The test is a range
check.  The proof coding is not onto, so a search over numbers must
reject the numbers that code no proof.

\paragraph{A universal at the decider.}
For $u = \lambda(\vz \jnat)$, the sentence $\cA\cdot u$ says that every
natural is a natural.  The decider affirms $\cA\cdot u$ at the first
fuel that evaluates the head and exceeds the code of some proof of
$\{\vz\jnat\} \vdash \lft u\cdot\vz$.  A two-line proof suffices: the
hypothesis $\vz \jnat$, then the instance by a conversion step
(\textsf{allT\_example\_proof}), and the decider affirms $\cA\cdot u$
at every fuel above that proof's code
(\textsf{prDecV\_allT\_example}).  The decider never evaluates
$u\cdot\tofn{n}$ for any $n$.  Every affirmation of a universal
$\cA\cdot u$ needs a fuel above the code of the judgment
$\{\vz\jnat\} \vdash \lft u\cdot\vz$
(\textsf{prDecV\_allT\_fuel\_lower}).  Here the two-line proof's code
has about 3,400 decimal digits, and the judgment's code about 1,700.

\paragraph{The certificate grain in HOL.}
At a quantifier node, decider soundness needs requirement~(iii) of
the semantic clause, the infinite family.  The certificate supplies the family
through the companion paper's soundness invariant \textsf{prf\_owsjt}.

\begin{thm}[\textsf{prDecV\_allT\_family}, \textsf{prDecV\_allF\_family}]
\label{thm:a-family}
Let $p$ be a number and $P = \mathrm{decode}\,p$.
\begin{enumerate}
\item If $P \proves \{\vz\jnat\} \vdash \lft u\cdot\vz$, then for
  every $n$ there is an $s_n$ with $u\cdot\tofn n \wbev{s_n} \ctrue$.
\item If $P \proves \emptyset \vdash \neg(u\cdot\tofn i)$, then
  $u\cdot\tofn i \wbev{s_i} \cfalse$ for some $s_i$.
\end{enumerate}
\end{thm}

\noindent
\emph{Proof.}  For part~1, apply \textsf{prf\_owsjt} to the valid
proof $P$, which contains the judgment $\{\vz\jnat\} \vdash \lft
u\cdot\vz$.  Instantiate the resulting soundness statement at the
substitution that sends $\vz$ to $\tofn n$ and each $v_{k+1}$ to
$v_k$.  The substitution satisfies the hypothesis $\vz\jnat$, since
$\tofn n$ has a value.  The substitution's shift undoes the lift, so the
substitution maps $\lft u\cdot\vz$ to $u\cdot\tofn n$.  So $u\cdot\tofn n$ is grounded-true.  Part~2 instantiates the
same invariant at the identity substitution: the negation is
grounded-true, so its argument is grounded-false.  \qed

\medskip
Completeness has a head part, \textsf{prDecH\_complete}: if $t
\whev{s} h$, then $\mathsf{prDecH}_{s'}(t) = \Some{h}$ for some $s'$.
Head completeness is proved together with Theorem~\ref{thm:in-dec}(2), by one
induction over the two semantic relations.  Monotonicity holds for
both deciders.

\paragraph{The packed verdict.}
\textsf{prDecVEnc} packs the value verdict at a code:
$\decP_s(x) = 0$ when $\decV_s(\mathrm{decode}\,x) = \None$, and
$\decP_s(x) = r+1$ when the verdict is $\Some r$.  The head verdict
is packed the same way (\textsf{prDecHEnc}).

\paragraph{The table machine.}
The table machine (\texttt{RGA\_Decider/RGA\_Decide\_PR}) is
\[
  \dtab(s, x) \;=\; \mathit{nth}\bigl(\mathit{levels}(B(s,x), s),\ x\bigr) ,
\]
where $\mathit{levels}(B, s)$ is the level-$s$ table over the codes
$0, \ldots, B$, built by primitive recursion on $s$ from a base table
and a step function that computes each cell of level $s+1$ from the
cells of level $s$ (\textsf{prTabLevels}).  The bound is $B(s,x) =
W^{c(s)}(x)$ (\textsf{prBound}), where $W$ is a growth step that
dominates every code a decider step can read from a code (the
$\beta$-instance, the recursor's reduct, and the subterms;
\textsf{Wstep}), and $c(s) = 2^s - 1$ counts the steps
(\textsf{cnt}: $c(0) = 0$, $c(s+1) = 2c(s) + 1$).  Iterating $W$
a number of times exponential in $s$ gives the tower-like growth noted
in \S\ref{sec:inside-dec}.

\emph{Proof sketch of Theorem~\ref{thm:in-tab}.}  Primitive
recursiveness (\textsf{n2pr\_prDecTab}) is closure of PR functions
under composition and the list operations, with no hypothesis.
Correctness reads one cell.  An invariant on levels
(\textsf{table\_inv}) states that, for every term $t'$ whose
\emph{nested majorant} $\mathit{Mnest}(lv, \quo{t'})$ lies below the
bound, the level-$lv$ cells of $\quo{t'}$ hold the two packed verdicts
at fuel $lv$.  The majorant bounds every code that the deciders read
from $t'$ within $lv$ steps, re-entering at the $\beta$-instance and
at the recursor's reduct, which can be larger than the input.  The
invariant goes by induction on the level, with step fidelity for the
head and value step functions under bounded readback.  The theorem
then follows from $\mathit{Mnest}(s, x) \le B(s, x)$
(\textsf{Mnest\_reenter}).  \qed

\subsection{The compiled decider and the truth predicate}
\label{app:machines-truth}

\paragraph{The explicit tree.}
The tree $\taut$ (\textsf{prDecTabIdx}, \texttt{RGA\_Kit/Code\_Index})
writes the table machine out from the initial functions, $\iSub$ and
$\iRec$: the level recursion is one $\iRec$ node over the step tree,
and the bound is an explicit tree for $W^{c(s)}$.

\begin{thm}[\textsf{pr2fi\_prDecTabIdx}]
\label{thm:a-taut}
For every natural $x$ and term $t$, $\den{\taut}(x, \quo{t}) =
\dtab(x, \quo{t})$.
\end{thm}

\noindent
\emph{Proof.}  The denotations of the subtrees are the HOL functions
that the subtrees transcribe (the bound tree computes \textsf{prBound}; the level
tree's table agrees with \textsf{prTabLevels} on the cells the
invariant covers, \textsf{table\_invT}).  The two cell readings at
$\quo{t}$ are the packed verdicts, by the same majorant argument as
Theorem~\ref{thm:in-tab}, and \textsf{prDecTab\_correct} closes the
equation.  \qed

\medskip
Theorem~\ref{thm:in-dec-closed} composes Theorem~\ref{thm:in-cmp-num}
at $\taut$, Theorem~\ref{thm:a-taut}, Theorem~\ref{thm:in-tab} and the
closed law of the compiled second projection (\textsf{cpv\_snd}).  At
abstract arguments, the compiled decider is a natural number:

\begin{thm}[\textsf{hq\_iPrDecV}]
\label{thm:a-dec-nat}
If $\Gamma \vdash s \jnat$ and $\Gamma \vdash c \jnat$, then $\Gamma
\vdash \Dec(s, c) \jnat$.
\end{thm}

\noindent
So $\Aff(s, c)$, an equality test of $\Dec(s, c)$ against $\tofn 2$,
is decided whenever its arguments are naturals.  The proof is
Theorem~\ref{thm:a-cmp-total-app} at $\taut$ and the totality of the
compiled projection.

\paragraph{One fuel level at an abstract fuel.}
Write $\tau_{\mathrm{lev}}$, $\tau_{\mathrm{step}}$ and $\tau_B$ for
the explicit trees of the level table, of its step function and of the
bound (\textsf{prTabLevelsIdx}, \textsf{prStepTabIdx},
\textsf{prBoundIdx}).  Write $\mathit{Table}(b, v) =
\cmp(\tau_{\mathrm{lev}})\cdot b\cdot v$ for the compiled level-$v$
table at bound $b$, $\mathit{nth}(T, c)$ for the compiled list lookup,
and $\mathit{pair}$ for the compiled Cantor pair.

\begin{thm}[\textsf{iPrDecV\_suc\_unfold\_any}]
\label{thm:a-suc-unfold}
Suppose $\Gamma \vdash v \jnat$ and $\Gamma \vdash c \jnat$, and write
$b = \cmp(\tau_B)\cdot\suc v\cdot c$.  Then
\[
  \Gamma \vdash \Dec(\suc v, c) =
    \mathit{snd}\bigl(\mathit{nth}\bigl(\cmp(\tau_{\mathrm{step}})\cdot
      \mathit{pair}(b, v)\cdot \mathit{Table}(b, v),\ c\bigr)\bigr) .
\]
\end{thm}

\noindent
The right side is the value half of the cell of $c$ after one
application of the step function to the level-$v$ table.  The bound
$b$ stays a compiled application, because the development has only
the denotation law of $\tau_B$, not an unfolding of it at abstract
arguments.  \emph{Proof.}  The successor law of
Theorem~\ref{thm:in-cmp-rec} unfolds the level recursion at $\suc v$,
the read-out at the top of $\taut$ unfolds by the composition law, and
congruence under the compiled projection finishes.  The successor law
needs a variable that does not occur in $b$ and is not free in
$\Gamma$; the proof chooses one.  The underlying lemma
\textsf{iPrDecV\_suc\_unfold} takes any $b$ with $\Gamma \vdash
\cmp(\tau_B)\cdot\suc v\cdot c = b$ and such a variable as
hypotheses.  \qed

\paragraph{Introduction, elimination and transport.}
Theorem~\ref{thm:in-wtrue-intro} (\textsf{iWtrueI}) is the
equation form of the following rule.

\begin{thm}[\textsf{iWtrueI\_abs}]
\label{thm:a-wtrue-intro}
If $\Gamma \vdash c \jnat$, $\Gamma \vdash a \jnat$ and $\Gamma \vdash
\Aff(a, c) \ne \z$, then $\Gamma \vdash \Tru(c)$.
\end{thm}

\noindent
\emph{Proof.}  The test is then derivable, so its negation equals
$\z$.  One $\beta$-step at the fuel $a$ shows that the body
$\lambda\neg\Aff(\vz, \lft c)$ applied to $a$ is that negation, so the
instance of the universal at $a$ is refuted.  \rga's rule for refuting
a universal from a refuted instance at a natural gives
$\neg(\cA\cdot(\lambda\neg\Aff(\vz,\lft c)))$, which is $\Tru(c)$.
\qed

\begin{thm}[\textsf{rga\_wtrueE}]
\label{thm:a-wtrue-elim}
Suppose $\Gamma \vdash \Tru(c)$, and
\[
  \vz\jnat,\ \Aff(\vz, \lft c),\ \lft\Gamma \;\vdash\; \lft G ,
\]
where $\lft\Gamma$ lifts every hypothesis of $\Gamma$.  Then $\Gamma
\vdash G$.
\end{thm}

\noindent
The elimination opens the fuel as the fresh variable $\vz$ and hands
the consumer the affirmation at that fuel.  The lifts keep the
context's and the goal's free variables apart from $\vz$, so the rule
needs no freshness side condition.  \emph{Proof.}  $\Tru(c)$ is a
negated universal, and \rga's elimination rule for a negated universal
opens a fresh natural $\vz$ with the negated instance
$\neg\neg\Aff(\vz, \lft c)$ as hypothesis.  Double-negation
elimination at that hypothesis and a cut give the premise's context.
\qed

\begin{thm}[\textsf{iWtrue\_transport}]
\label{thm:a-transport}
If $\Gamma \vdash a \jnat$, $\Gamma \vdash b \jnat$, $\Gamma \vdash a
= b$ and $\Gamma \vdash \Tru(a)$, then $\Gamma \vdash \Tru(b)$.
\end{thm}

\noindent
\emph{Proof.}  Eliminate $\Tru(a)$ by Theorem~\ref{thm:a-wtrue-elim}.
In the opened context, the lifted equation $\lft a = \lft b$ rewrites
the code argument inside $\Aff(\vz, \lft a)$ by congruence of compiled
application, of the projection and of the equality test
(\textsf{cmp\_argI}, \textsf{cpSnd\_argI}, \textsf{eqT\_cong}), which
gives $\Aff(\vz, \lft b) \ne \z$.  Theorem~\ref{thm:a-wtrue-intro} at
the fuel $\vz$ gives $\Tru(\lft b)$, the lift of the goal.  \qed

\paragraph{The numeral-code encoder.}
$\Num = \cmp(\tau_{\mathrm{num}})\cdot\z$ (\textsf{iNumT}), where
$\tau_{\mathrm{num}}$ (\textsf{num\_codeIdx}) is one $\iRec$ node that
starts from the code of $\z$ and applies the explicit tree of the
successor-code constructor at each step.  The numeral-code law of
\S\ref{sec:inside-d3} is \textsf{iNumT\_closed}, proved by
Theorem~\ref{thm:in-cmp-num}, and $\Gamma \vdash a\jnat$ implies
$\Gamma \vdash (\Num\cdot a)\jnat$ (\textsf{iNumT\_N}).

\paragraph{The proof of Theorem~\ref{thm:in-never}.}
The proof runs in four steps, all in
\texttt{RGA\_Decider/Decider\_Truth}.  Recall that a packed cell
holds $0$ for $\None$ and $r+1$ for $\Some r$, so ``affirmed'' means
the packed value $2$ in the cell of an equation.
\begin{enumerate}
\item \emph{Slot facts.}  At every level $l$ and bound $B$ with
  $\Gamma \vdash B \jnat$ and $\Gamma \vdash l \jnat$, \rga proves
  that the value cell of $\quo{\z}$ holds the packed value $0$ or $1$
  ($\None$ or $\Some 0$; \textsf{ir\_tabL\_slot0}), and that the value
  cell of $\quo{\suc\z}$ holds $0$ or $2$ ($\None$ or $\Some 1$;
  \textsf{ir\_tabL\_slot1}).  Both slot facts go by a case split on the level
  (\textsf{ir\_tabL\_slot\_skel}).  At level zero the base table
  holds $\None$ everywhere.  A successor level is one application of
  the step function, which writes $\Some 0$ into the cell of $\z$, and
  into the cell of $\suc\z$ either $\None$ or the successor of the
  value in the previous cell of $\z$.
\item \emph{The new cell.}  The step function computes the cell of
  $\quo{\z = \tofn 1}$ from the cells of the equation's two sides: $\None$ if
  either side is $\None$, otherwise the truth value of the equality of
  the two values.  With the slot facts the two values are $0$ and $1$,
  so the cell is $\None$ or $\Some{\cfalse}$, never the packed value $2$
  (\textsf{ir\_cellV\_cEq01\_ne2}, \textsf{stepTab\_cEq01\_cell\_ne2}).
  The HOL mirror of this step is \textsf{stepV\_sEq01\_ne2\_any}.
\item \emph{One fuel level.}  Theorem~\ref{thm:a-suc-unfold} at
  $c = \bquo{\z = \tofn 1}$ turns step~2 into $\Gamma \vdash
  \Dec(\suc v, \bquo{\z = \tofn 1}) \ne \tofn 2$ for every $v$ with
  $\Gamma \vdash v\jnat$ (\textsf{iPrDecV\_cEq01\_step\_reduced}, with
  a variable that is fresh for $v$ and not free in $\Gamma$, which
  always exists).
  This step uses neither an induction hypothesis nor any fact about $v$.
\item \emph{Assembly.}  The base case at the numeral $\z$, and indeed
  at every numeral, is Theorem~\ref{thm:in-dec-closed} with the HOL
  fact $\decP_s(\quo{\z = \tofn 1}) \ne 2$
  (\textsf{prDecVEnc\_sEq01\_ne2}, \textsf{iPrDecV\_cEq01\_num\_ne2}).
  \rga's induction rule carries the decided statement $\Dec(v,
  \bquo{\z=\tofn 1}) \ne \tofn 2$ to an abstract natural
  (\textsf{iPrDecV\_cEq01\_ne2\_abs}), and to the eigenvariable $\vz$
  (\textsf{iPrDecV\_cEq01\_ne2\_eigen}).  $\forall$-introduction then
  gives the universal of Theorem~\ref{thm:in-never}, in every context.
\end{enumerate}
Steps~1 and~3 together cover every successor fuel without using the
induction hypothesis, and fuel $\z$ is the base case; the induction
serves only to reach an abstract fuel.  (The
zero-or-successor case split of step~1 on an abstract natural,
\textsf{ir\_nat\_zero\_or\_suc}, is itself derived in \rga by
induction.)
\S\ref{app:machines-traps} on negative results shows why the proof does not use the
statement itself as the induction motive.

\subsection{The opaque twins}
\label{app:machines-opaque}

In the ten sessions that this development adds to the companion
formalization (Section~\ref{sec:mech}), exactly four definitions pick an
index with \textsf{priof}:
\begin{itemize}
\item \textsf{prDecIdx} $=$ \textsf{priof}~\textsf{prDecTab}
  (\texttt{RGA\_Decider/RGA\_Decide\_PR}), compiled to the opaque
  decider \textsf{iPrDecOpq} and its value projection
  \textsf{iPrDecVOpq};
\item \textsf{numCodeIdx} (\texttt{RGA\_Decider/Decider\_Truth}),
  compiled to the opaque numeral-code encoder \textsf{iNum};
\item \textsf{chkIdx} $=$ \textsf{priof}~\textsf{chk2}
  (\texttt{RGA\_Kit/Term\_Coders}), where \textsf{chk2} is the
  companion paper's PR proof-certificate test at the empty context,
  compiled to $\ChkE(p, c) \ldef (\cmp(\mathit{chkIdx})\cdot p\cdot c
  = \tofn 1)$ (\textsf{iChkOpq});
\item \textsf{qfStepChkIdx} (\texttt{RGA\_Checker/Chk\_Base}), the
  per-node test of an opaque quantifier-free checker that no result of
  the paper uses.  The definition of the $\forall$-free checker
  $\ChkQf$ of Section~\ref{sec:con} does not mention it.
\end{itemize}
By Theorem~\ref{thm:a-priof} and Theorem~\ref{thm:in-cmp-num}, each
twin has the closed laws of its transparent counterpart:

\begin{thm}[\textsf{iPrDecVOpq\_closed}, \textsf{iPrDecV\_agree}, \textsf{iNum\_closed}, \textsf{iNum\_agree}]
\label{thm:a-agree}
For every fuel $s$, term $t$ and natural $n$,
\begin{gather*}
  \emptyset \vdash \Dec^{\mathrm{opq}}(\tofn s, \bquo{t}) = \tofn{\decP_s(\quo{t})} , \qquad
  \emptyset \vdash \Dec(\tofn s, \bquo{t}) = \Dec^{\mathrm{opq}}(\tofn s, \bquo{t}) , \\
  \emptyset \vdash \mathit{iNum}\cdot\tofn n = \bquo{\tofn n} , \qquad
  \emptyset \vdash \Num\cdot\tofn n = \mathit{iNum}\cdot\tofn n ,
\end{gather*}
where $\Dec^{\mathrm{opq}}$ is \textsf{iPrDecVOpq}.
\end{thm}

\noindent
The development proves the agreement equations at closed numeral
arguments, and for the decider at term-code numerals; at abstract
arguments the agreement equations would need an unfolding of the
opaque index.  The
opaque checker is faithful at proof codes:

\begin{thm}[\textsf{iChkOpq\_true}, \textsf{iChkOpq\_iff\_proves}]
\label{thm:a-chkopq}
For every proof $P$ and term $C$: if $P \proves \emptyset \vdash C$,
then $\emptyset \vdash \ChkE(\bquo{P}, \bquo{C})$; and $\emptyset
\vdash \ChkE(\bquo{P}, \bquo{C})$ holds if and only if $P \proves
\emptyset \vdash C$.
\end{thm}

\noindent
At an arbitrary natural $p$ the opaque checker accepts exactly what
the HOL decoding of $p$ proves.

\begin{thm}[\textsf{iChkOpq\_decodes}, \textsf{iChkOpq\_genuine\_iff}]
\label{thm:a-chkopq-decode}
For every natural $p$ and term $C$, $\emptyset \vdash \ChkE(\tofn p,
\bquo{C})$ holds if and only if the decoding of $p$ proves $\emptyset
\vdash C$.  And $\ChkE$ is genuine at every number (every $p$ that it
accepts for $\bquo{C}$ is $\quo{P}$ for a proof $P$ of $\emptyset
\vdash C$) if and only if the decoder's default value proves no
judgment $\emptyset \vdash C$.
\end{thm}

\noindent
Decoding is the inverse of the coding of proofs, and at a natural
outside the coding's range it returns a default value, a HOL constant
whose value HOL leaves unspecified.  The range is not all of the
naturals: some natural codes no proof (\textsf{prf\_proof\_noncode}).
So genuineness of $\ChkE$ at every number is equivalent to a
property of an unspecified constant, which HOL leaves open; the
statement that holds is Theorem~\ref{thm:a-chkopq-decode}.  The
transparent checker $\Chk$ is genuine at every number
(Theorem~\ref{thm:in-chk}(2)).  Since both checkers are faithful at
proof codes, the two checkers agree there as a HOL fact: \rga derives
$\ChkE(\bquo{P}, \bquo{C})$ exactly when it derives
$\Chk(\bquo{P}, \bquo{C})$ (\textsf{iChkOpq\_iChk\_iff}, a rule in both
directions).  The two checkers also agree inside \rga at these codes: both checkers are
decided at numerals, so for every proof $P$ and term $C$, \rga proves
the two implications $\ChkE(\bquo{P}, \bquo{C}) \limp \Chk(\bquo{P},
\bquo{C})$ and $\Chk(\bquo{P}, \bquo{C}) \limp \ChkE(\bquo{P},
\bquo{C})$ as one internal sentence (\textsf{iChkOpq\_iChk\_agree}, a
schema over the closed codes).  At abstract arguments no such sentence
is available, since the opaque checker's unfolding laws are not
derivable.

\subsection{The checker}
\label{app:machines-chk}

\paragraph{Detectors and the validity test.}
\rga's proof system \textsf{derives} (\texttt{RGA/RGA\_Proof}) has 42
rules.  The list \textsf{drvChkList} (\texttt{RGA\_Kit/Code\_Index})
holds 41 detector trees, one per rule, except that the two directions
of reflexivity share \textsf{chkEqReflIdx}: since $a\jnat$ abbreviates
$a = a$ (\textsf{trmNat}, \texttt{Arithmetic}), both directions conclude $a = a$
from $a = a$.  The disjunction of the 41 detectors is the rule test
\textsf{drvChkIdx}.  The validity tree \textsf{validChkIdx} checks, at
every position $z$ below the length of the coded list, the body
\textsf{validBody}: the entry at $z$ is a genuine judgment code (its
hypothesis slot codes a finite set of terms and its conclusion slot
codes a term, \textsf{nodeCodeIdx}), and the rule test fires on that
entry against the entries after $z$.  The compiled test is
$\mathit{TreeValid}(p) \ldef (\cmp(\mathit{validChkIdx})\cdot p\cdot p
= \tofn 1)$ (\textsf{iTreeValid}, \texttt{RGA\_Checker/Chk\_Formers}).
The tree $\pi$ (\textsf{provesChkIdx}) conjoins $p \ne 0$, ``the head
judgment's code is $h$'', and validity.

\begin{thm}[\textsf{drvChk\_correct}, \textsf{validChk\_correct}, \textsf{provesChk\_correct}]
\label{thm:a-chk-correct}
At codes, the rule test computes the companion paper's
primitive-recursive one-step rule check; for every proof $P$ and
every $w$, $\den{\mathit{validChkIdx}}(\quo{P}, w)$ is $1$ if $P$ is
valid and $0$ otherwise; and for every $H$, $c$ and $P$,
$\den{\pi}(\cpair{\quo{H}}{\quo{c}}, \quo{P})$ is $1$ if $P \proves H
\vdash c$ and $0$ otherwise.
\end{thm}

\begin{thm}[\textsf{validChk\_genuine}]
\label{thm:a-valid-genuine}
For all naturals $p$ and $w$: if $\den{\mathit{validChkIdx}}(p, w) \ne
0$, then $p = \quo{P}$ for some proof $P$.
\end{thm}

\noindent
\emph{Proof.}  Decode $p$ as a list of numbers.  At every position the
body has fired, so the node-code test has accepted the entry, which is
therefore the code of a judgment.  A list of judgment codes is the
code of the list of those judgments.  \qed

\medskip
The node-code test is what makes Theorem~\ref{thm:a-valid-genuine}
true; \S\ref{app:machines-traps} on negative results shows a validity
test without the node-code test accepting a code of no proof.

\paragraph{From HOL values to derivations.}
\textsf{iTreeValid\_faithful} states that $\emptyset \vdash
\mathit{TreeValid}(\bquo{P})$ holds exactly when $P$ is valid, and
\textsf{iTreeValid\_genuine} that $\emptyset \vdash
\mathit{TreeValid}(\tofn p)$ implies $p = \quo{P}$ for some valid $P$.
Both rest on \textsf{rga\_val\_iff}: for the well-formed closed term
$\cmp(\mathit{validChkIdx})\cdot\tofn p\cdot\tofn p$, \rga derives
that this term equals $\tofn 1$ exactly when the term evaluates to $1$,
which by
Theorem~\ref{thm:a-cmp-correct} is the HOL value of the tree
(\textsf{iTreeValid\_iff\_value}).  Theorem~\ref{thm:in-chk} follows.
Part~2 (\textsf{iChk\_genuine\_proves}, packaged with validity as
\textsf{iChk\_genuine}) reads the four conjuncts of $\Chk$ off the
derivation, applies \textsf{iTreeValid\_genuine}, and reads the head
judgment through the compiled accessors (\textsf{iConcOf\_encode},
\textsf{iMcLen\_iHypsOfT\_encode}).  Part~1 needs in addition the
converse, \textsf{iChk\_accept\_of\_parts}: for a valid proof
$(\emptyset, c) :: Q$, \rga derives $\Chk(\bquo{(\emptyset, c) :: Q},
\bquo{c})$.
Since \textsf{rga\_val\_iff} rests
on \rga's soundness (\textsf{rga\_val\_iff} is a control of the dependency check of
\S\ref{app:machines-audit}), both directions of
\textsf{iTreeValid\_faithful} do too.  The direction that a valid
proof's code passes the test needs no soundness: \rga derives
$\mathit{TreeValid}(\bquo{P})$ for every valid $P$ by closed
evaluation of the compiled test alone (\textsf{d2\_valid\_accepted},
\texttt{RGA\_Goedel/Goedel\_D2}, a schema over $P$).  The decided D2
(\textsf{D2\_decided}), which the dependency check lists as
soundness-free, uses that form.

\paragraph{A meta corollary.}
Genuineness and \rga's consistency give \textsf{rga\_con\_meta\_num}:
if $\emptyset \vdash \Chk(\tofn p, \bquo{\z = \tofn 1})$, then
$\mathit{False}$.  The proof extracts a proof of $\z = \tofn 1$ by
genuineness and refutes it by Theorem~\ref{thm:pl-con}.  The corollary
is a HOL theorem about which judgments \rga derives at numerals.  It
is not internal consistency, which is $\emptyset \vdash \ConR$ and
which \rga does not prove ([\citen{goedel-paper}, Section ``G\"odel's theorems for \rga'']).

\subsection{Negative results behind three design choices}
\label{app:machines-traps}

\paragraph{Coders exact at every arity.}
The table's $\beta$-step substitutes into a body through coders of
lifting and substitution.  An earlier family of coders dispatches on
at most two argument slots of an operator node.  That covers every
well-formed term, since no \rga operator takes more than two
arguments.  The trap is ill-formed input.  Let $E_3$ be the equation
node with three arguments $\z, \z, \z$, written $\z = \z = \z$, and
$R_3 = (\lambda\vz)\cdot E_3$.  Neither term is well-formed.  The two-slot
lift coder truncates the code of $E_3$ to the code of $\z = \z$, and a
table whose $\beta$-arm substitutes through that lift coder (the local twin
\textsf{prDecTabTr}) therefore affirms $R_3$.  The facts are, in
\texttt{RGA\_Decider/RGA\_Decider\_Traps}: the truncating table has
cell $2$ at $(3, \quo{R_3})$ (\textsf{prDecTabTr\_redexAr3\_affirms});
the real table and the compiled decider have cell $0$ there
(\textsf{prDecTab\_redexAr3\_refuses},
\textsf{iPrDecV\_redexAr3\_refuses}); and $R_3$ has no value at any
fuel (\textsf{redexAr3\_no\_wosr}).

\begin{thm}[\textsf{truncating\_decider\_contradicts\_converse}]
\label{thm:a-trunc}
Let $D$ be any function from two terms to a term, and $W$ any function
from terms to terms.  The following three assumptions are jointly
contradictory.
\begin{enumerate}
\item \emph{$D$ computes the truncating table:} for every $s$ and
  $t$, $\emptyset \vdash D(\tofn s, \bquo{t}) = \tofn{d}$, where $d$
  is the value half of the truncating table's cell at $(s, \quo{t})$.
\item \emph{$W$ has an introduction rule like
  Theorem~\ref{thm:in-wtrue-intro}:} for all terms $c$ and $a$, if
  $\emptyset \vdash c\jnat$, $\emptyset \vdash a\jnat$ and $\emptyset
  \vdash D(a, c) = \tofn 2$, then $\emptyset \vdash W(c)$.
\item \emph{$W$ has a converse without side conditions, like
  Theorem~\ref{thm:pl-adeq}(2):} for every term $t$, if $\emptyset
  \vdash W(\bquo{t})$, then $t \wbev{s} \ctrue$ for some $s$.
\end{enumerate}
\end{thm}

\noindent
\emph{Proof.}  Assumption~1 at $(3, R_3)$ gives $\emptyset \vdash
D(\tofn 3, \bquo{R_3}) = \tofn 2$.  Assumption~2 at numerals gives
$\emptyset \vdash W(\bquo{R_3})$, assumption~3 a value of $R_3$, and
$R_3$ has none.  \qed

\medskip
The theorem is a HOL theorem about three hypothetical rules; no
derivation of $\mathit{False}$ in \rga is involved.  The development
needs the converse of Theorem~\ref{thm:pl-adeq} without a
well-formedness hypothesis, so it builds the table on coders that are
exact at every arity (\textsf{lift0\_code}, \textsf{bsub\_code}).  The
corollary \textsf{iPrDecV\_not\_truncating} records that the compiled
decider $\Dec$ does not compute the truncating table: instantiating
the theorem at $D = \Dec$ and $W = \Tru$, assumptions~2 and~3 are
Theorem~\ref{thm:in-wtrue-intro} and Theorem~\ref{thm:pl-adeq}(2).

\paragraph{A motive that is not inductive.}
The obvious motive for Theorem~\ref{thm:in-never} is ``the table does
not affirm $\z = \tofn 1$''.  Let $T_1$ be the table whose value cell
holds the packed value $1$ ($\Some 0$) at every code up to $\quo{\z =
\tofn 1}$.  $T_1$ satisfies the motive.  One step of the value step
function from $T_1$ reads the value $0$ on both sides of the equation
and writes $\Some{\ctrue}$, the packed value $2$
(\textsf{sEq01\_motive\_not\_inductive}).  So the motive is not
preserved by a step, and the proof of
\S\ref{app:machines-truth} proves the two slot facts at every level
instead.

\paragraph{The node-code test.}
A validity test without the node-code test of
\S\ref{app:machines-chk} checks only that some detector fires at each
position.  The detector for the axiom $\z\jnat$ reads only the
conclusion slot of its entry.  So the two-entry code $p_0$
whose entries both conclude $\z\jnat$, the second with a hypothesis
slot that codes no set of terms, passes that test.
\textsf{p0\_frame\_refuted} (\texttt{RGA\_Checker/RGA\_Checker\_Traps})
states that $p_0 > 0$, that the head entry of $p_0$ is the code of
$\emptyset \vdash \z\jnat$, that the weaker test accepts $p_0$, and
that the in-range test returns $0$ at $p_0$: no proof has code $p_0$.
With the node-code test, Theorem~\ref{thm:a-valid-genuine} excludes
such codes.

\subsection{The build-time check}
\label{app:machines-audit}

The walk of the check of \S\ref{sec:inside-audit} is the Isabelle/ML
structure \textsf{Choice\_Audit} of
\texttt{RGA\_Goedel/Goedel\_II} (subsection ``The choice audit, as an
ML structure''), a theory that every audit site imports ($\ConR$ is
defined there too).  A block near the end of
\texttt{RGA\_Goedel/Goedel\_II\_PR} runs the walk on the roots of the first
run, a block at the end of \texttt{ConOGA/ConOGA\_Traps} on the
machines of Section~\ref{sec:con}, a block in
\texttt{ConOGA/ConOGA\_Audit} on the remaining objects of the notation
table, a block at the end of \texttt{RGA\_Goedel/Goedel\_D3} on
the transformers of the decided D3, and blocks at the ends of
\texttt{RGA\_Goedel/Goedel\_D3\_Uniform} and
\texttt{RGA\_Goedel/Goedel\_D2} on the transformers of the decided D3
as one sentence over the code and of the decided D2, and two blocks in
\texttt{RGA\_Goedel/Goedel\_Cut}, on the cut transformers and on the
explicit witnesses of the G\"odel sentences in decided-antecedent form.  The walk reads the theory's record
of definitional dependencies (\textsf{Defs.specifications\_of}): from
a constant it visits every specification of that constant, and every
constant and type on each specification's right-hand side, without
unfolding or reducing anything, so an overloaded constant contributes
all its instances.  The closure of a set of roots is the set of items
visited.

A closure is \emph{flagged} if it contains the constant
\textsf{Indexed\_Recursion.priof}, or a constant whose own definition
applies HOL's choice operator \textsf{SOME} or definite-description
operator \textsf{THE} directly and that is not on a fixed whitelist.
The whitelist has 28 entries: 27 constants of the HOL library defined
that way (among them \textsf{If}, \textsf{wfrec}, \textsf{Least},
\textsf{THE\_default}, \textsf{the\_inv\_into}, \textsf{Finite\_Set.fold},
the case combinators of the library's datatypes, their folds for lists
and binary numerals, and the cardinal constants of the BNF package),
and the fold \textsf{ctor\_fold\_pridx} that the datatype package
generates for the index type \textsf{pridx}.  Every audit checks that
each whitelisted constant declared in the audit's theory is defined by a
direct \textsf{SOME} or \textsf{THE} and is not index-valued, and that
every transparent closure contains each of the 27 library constants
that is declared in the audit's theory.  So the library entries are not
exemptions chosen per root: when the whitelist was fixed, the
direct-choice constants common to the 40 closures then audited by the
first three runs were exactly the 27, and the only other one met was
the \textsf{pridx} fold, in the witness's closure.  The eight roots
added to the third run later, and the roots of the fourth to eighth
runs, pass the same checks: the closures of these roots contain the
library constants declared in the audited theory, and no other
choice-defined constant outside the whitelist.  This flag does not
fire at a constant that reaches choice only through a whitelisted
library constant, for instance one defined by \textsf{Least}.  A
second audit closes that gap for indices (described after the list).
An audit raises an error, which fails the build of the theory, if any
of the following holds:
\begin{itemize}
\item a positive control of the choice flags fails: the index
  \textsf{prConstFunSuc} of the grounded-deduction base
  (\texttt{Indexed\_Recursion}), defined by a direct \textsf{SOME} outside
  \textsf{priof}, must be flagged with no \textsf{priof} in its
  closure; the closure of \textsf{inv\_into}, defined by a direct
  \textsf{SOME} and not index-valued, must be flagged; and the
  direct-use test must fire at \textsf{inv\_into} (\textsf{SOME}) and
  at \textsf{the\_inv\_into} (\textsf{THE});
\item the joint closure of the witness constants
  \textsf{goedelSigmaOf}, \textsf{goedelSigmaIdx} and
  \textsf{premCertIdx} is flagged, or has fewer than 100 constants;
\item the closure of one of the thirteen transparent roots (the seven
  machines $\Tru$, $\Dec$, $\Chk$, $\ConR$, $\Num$, $\Cert$ and
  $\Subst$, the EXPL and ex-falso trees
  \textsf{goedelEIdx} and \textsf{exfIdx}, the proof object
  \textsf{neq01Prf} of $\neg(\z = \tofn 1)$, and $\ChkQf$ with its two
  scans \textsf{qfScanAcc} and \textsf{iQfScan}) is flagged, or has
  fewer than 100 constants;
\item the closure of \textsf{prDecIdx}, or of the opaque compiled
  decider \textsf{iPrDecOpqT}, which reaches \textsf{priof} only
  through \textsf{prDecIdx}, does not contain \textsf{priof};
\item the closure of one of the three opaque twins, or of the index
  \textsf{chkIdx} of $\ChkE$, does not contain \textsf{priof};
\item some whitelisted constant is not defined by a direct choice, or is
  index-valued, or some transparent closure lacks one of the 27 library
  constants declared in the audit's theory.
\end{itemize}
Otherwise the block prints its report.  Since the twins must reach
\textsf{priof}, listing one of them among the transparent roots fails
the build.

The second audit, \textsf{Choice\_Audit\_WL} (near the end of
\texttt{Goedel\_II\_PR}), flags every constant in a closure, of any
type, whose own definition applies one of the whitelisted constants
that use choice in an essential way (\textsf{the\_inv\_into},
\textsf{THE\_default}, \textsf{Least}, \textsf{Finite\_Set.fold}, the
cardinal constants of the BNF package and the quotient representation
\textsf{quot\_type.rep}) at a type instance that involves the index
type \textsf{pridx}.  The case combinators, \textsf{If}, the recursors
and \textsf{wfrec} are not in this list: they use \textsf{THE} only
to name a determined value.  For the same reason one definition shape
is exempt: the function package defines \textsf{f\_sumC} as
\textsf{THE\_default} applied to the graph \textsf{f\_graph} of $f$.
The audit fails the build unless no transparent closure has a flagged
constant and every control has one.  The control of this audit,
\textsf{wlCtlIdx}, is
an index defined by \textsf{the\_inv\_into} at the index type; the
build also checks that the first audit does not flag \textsf{wlCtlIdx}, so the
control lies in exactly the gap this audit closes.  The audit runs in
\texttt{Goedel\_II\_PR} on the roots of the first run, and in
\texttt{ConOGA\_Audit} on the roots of all the other runs.  All the
whitelist audit's transparent closures are clean, and its control is
flagged.  Choice
reached through a whitelisted constant at a type that does not
involve indices is still flagged by neither audit; the code-generation
check below meets such choice wherever the code of an index-valued
constant calls the whitelisted constant.

The third audit, \textsf{Choice\_Audit\_UD} (at the end of
\texttt{Goedel\_II\_PR}), flags, in the closure of a root, every
occurrence inside the specification of a closure constant of a
constant that has no defining specification, at a type involving the
index type: HOL's unspecified value \textsf{undefined}, for instance,
or an axiomatized constant.  The audit also flags every closure constant of
such a type that has no defining specification itself.  The
primitives of the logic, \textsf{The} and \textsf{Eps} (which the
first audit handles), the morphisms of type definitions (the index
type is built through them), and \textsf{undefined} as the default
value of a terminating function's definition (the exemption of the
whitelist audit) are exempt; each run checks that the exempt names
are declared and have no defining specification.  The build fails
unless no transparent closure has a flag, the control
\textsf{udCtlIdx}, an index defined as \textsf{undefined}, is flagged,
and the first audit leaves that control clean.  The audit runs at the
end of \texttt{Goedel\_II\_PR} on its own roots and in
\texttt{ConOGA\_Audit} on the roots of the other runs, the same list
as the whitelist audit's.  No root is flagged.  The audit is not
vacuous: each closure meets between 114 and 129 constants without a
definition, and in the union of the closures exactly 15 occurrences at
a type involving the index type, all of them exempt.  Like the first
audit, the unspecified-constant audit walks constants by name, so a constant without a definition
that is used polymorphically and instantiated at the index type only
downstream is seen only through the types of the definitions it
walks.

The other runs have the same pass conditions, except the positive
controls of the choice flags, which run only in the first run.  The second
run, at the end
of \texttt{ConOGA\_Traps}, the only theory that sees all of the
machines of Section~\ref{sec:con} that this run audits, has as roots the
value record (\textsf{RF}, \textsf{recSelf}), the fragment and level
classes (\textsf{iFrag}, \textsf{iStat}, \textsf{iFragV},
\textsf{iStatLv}, \textsf{iLvClass}, \textsf{iSub1Qf}),
\textsf{iDynWFD}, \textsf{iLF}, \textsf{iLFs} and \textsf{iRes}; the
second run's control is the opaque twin \textsf{iPrv}.  The third run, in
\texttt{ConOGA\_Audit}, a last theory that imports every theory of the
development defining an object of the notation table, has as roots
the code constructors \textsf{mkNegS}, \textsf{mkOrS},
\textsf{mkEqS}, \textsf{mkAppS} and \textsf{mkSucS}, the diagonal
substitution \textsf{iDiagSub}, the no-swap scan \textsf{iNoSwap}, the
searches \textsf{fragSearch}, \textsf{conSearch} and
\textsf{goedelSigmaSearch}, and the Liar objects \textsf{lwD},
\textsf{goedelWS}, \textsf{goedelW} and \textsf{liarRec}, the
head-step scan \textsf{iHdStep}, and the objects of the level theorem
of \S\ref{sec:con-lv}: the depth machine \textsf{iDep}, the
certificates \textsf{OkI} and \textsf{OkU}, the level sentences
\textsf{iLFA} and \textsf{iLFI}, the class \textsf{iLvClassI}, and
\textsf{iLvClassU}, the class at a level given as an object term; the
third run's controls are \textsf{iPrv} and \textsf{chkIdx}.  The fourth run, at the
end of \texttt{Goedel\_D3}, has as roots the transformer of the
decided D3 and its builders (\textsf{d3SigmaOf}, \textsf{d3SigmaAt},
\textsf{d3SigmaIdx}, \textsf{d3CertIdx}).  The fifth run, at the end of
\texttt{Goedel\_D3\_Uniform}, has as roots the transformer of the
decided D3 as one sentence over the code, its certificate builder and
the code of its consequent (\textsf{uSigIdx}, \textsf{uCertIdx},
\textsf{prCU}, \textsf{prCUIdx}).  The sixth run, at the end of
\texttt{Goedel\_D2}, has as roots the transformer of the decided D2
and its tree (\textsf{d2M}, \textsf{d2Idx}).  The seventh run, at the
end of \texttt{Goedel\_Cut}, has as roots the cut transformers of
{}[\citen{goedel-paper}, theorem~\textsf{cut1M\_chk}] and their trees (\textsf{cut1M},
\textsf{cut1Idx}, \textsf{cut2M}, \textsf{cut2Idx}).  The eighth run, also
at the end of \texttt{Goedel\_Cut}, has as roots the explicit
witnesses of [\citen{goedel-paper}, theorem~\textsf{rga\_goedel2\_sentence}] and their trees
(\textsf{goedel1T}, \textsf{goedel1Idx}, \textsf{goedel2T},
\textsf{goedel2Idx}).  The whitelist audit and the
unspecified-constant audit in \texttt{ConOGA\_Audit} have the roots of
the last two runs among theirs too.  The control of the last five runs
is \textsf{prDecIdx}.  All transparent closures are clean, and all
controls are flagged.  Table~\ref{tab:a-audit} gives the closure sizes
of all eight runs.  Batch builds do not
keep the printed report, so these figures come from scratch theories that
import exactly the audit site and print each row.  Sizes depend on the
theory in which the closure is taken, since an overloaded constant
gains instances: \textsf{iPrv} has 928 constants at
\texttt{Goedel\_II\_PR} and 983 at \texttt{ConOGA\_Traps}, and
\textsf{chkIdx} has 919 at \texttt{Goedel\_II\_PR} and 974 at
\texttt{ConOGA\_Audit}.

\begingroup
\small
\begin{longtable}{@{}>{\raggedright\arraybackslash}p{0.13\linewidth}>{\raggedright\arraybackslash}p{0.37\linewidth}>{\raggedright\arraybackslash}p{0.21\linewidth}>{\raggedleft\arraybackslash}p{0.19\linewidth}@{}}
\caption{Closure sizes of the build-time check, per audit site,
  measured by scratch theories that import exactly
  the audit site.  No transparent root and not the witness is flagged;
  every control is.  Sizes depend on the
  theory in which the closure is taken, so later imports can change
  them.}
\label{tab:a-audit}\\
\hline
role & root & paper & constants \\
\hline
\endfirsthead
\hline
role & root & paper & constants \\
\hline
\endhead
\hline
\endfoot
\multicolumn{4}{@{}l}{\emph{first run, end of \texttt{Goedel\_II\_PR}}} \\
witness & \textsf{goedelSigmaOf}, \textsf{goedelSigmaIdx}, \textsf{premCertIdx}
  & $\sigma$ ([\citen{goedel-paper}, Section ``G\"odel's theorems for \rga'']) & 1{,}182 \\
transparent & \textsf{iWtrue} & $\Tru$ & 1{,}189 \\
transparent & \textsf{iPrDecV} & $\Dec$ & 1{,}126 \\
transparent & \textsf{iChk} & $\Chk$ & 1{,}085 \\
transparent & \textsf{ConRGA} & $\ConR$ & 1{,}087 \\
transparent & \textsf{iNumT} & $\Num$ & 808 \\
transparent & \textsf{iHasCert} & $\Cert$ & 1{,}138 \\
transparent & \textsf{iSubM} & $\Subst$ & 846 \\
transparent & \textsf{goedelEIdx} & tree of $e$ (EXPL) & 823 \\
transparent & \textsf{exfIdx} & ex-falso tree & 824 \\
transparent & \textsf{neq01Prf} & proof of $\neg(\z=\tofn 1)$ & 791 \\
transparent & \textsf{iChkQf} & $\ChkQf$ & 1{,}095 \\
transparent & \textsf{qfScanAcc} & $\QfScan$ & 861 \\
transparent & \textsf{iQfScan} & $\qfs$ & 840 \\
control & \textsf{prDecIdx} & & 962 \\
control & \textsf{iPrDecOpqT} & & 969 \\
control & \textsf{chkIdx} & index of $\ChkE$ & 919 \\
opaque twin & \textsf{iPrv} & $\Prv$ & 928 \\
opaque twin & \textsf{iChkOpq} & $\ChkE$ & 927 \\
opaque twin & \textsf{iNum} & & 812 \\
control (second flag) & \textsf{prConstFunSuc} & & 796 \\
control (third flag) & \textsf{inv\_into} & & 8 \\
\hline
\multicolumn{4}{@{}l}{\emph{second run, end of \texttt{ConOGA\_Traps}}} \\
transparent & \textsf{RF} / \textsf{recSelf} & $\Rcd$ & 915 / 1{,}059 \\
transparent & \textsf{iFrag} / \textsf{iStat} / \textsf{iFragV} & $\Frag$ & 1{,}266 / 1{,}234 / 1{,}245 \\
transparent & \textsf{iStatLv} / \textsf{iLvClass} / \textsf{iSub1Qf} & $\LvC_N$ & 1{,}265 / 1{,}287 / 1{,}219 \\
transparent & \textsf{iDynWFD} & $\DynFam$ & 1{,}403 \\
transparent & \textsf{iLF} / \textsf{iLFs} & $\LF_N$, $\LF_{\le N}$ & 1{,}134 / 1{,}135 \\
transparent & \textsf{iRes} & $\Res$ & 1{,}422 \\
control & \textsf{iPrv} & $\Prv$ & 983 \\
\hline
\multicolumn{4}{@{}l}{\emph{third run, \texttt{ConOGA\_Audit}}} \\
transparent & \textsf{mkNegS}, \textsf{mkOrS}, \textsf{mkEqS}, \textsf{mkAppS}, \textsf{mkSucS} & code constructors & 826 each \\
transparent & \textsf{iDiagSub} & $\Diag$ & 872 \\
transparent & \textsf{iNoSwap} & $\mathit{NoSwap}$ & 1{,}215 \\
transparent & \textsf{fragSearch} / \textsf{conSearch} & searches (\S\ref{sec:con-frag}) & 1{,}285 / 1{,}114 \\
transparent & \textsf{goedelSigmaSearch} & search witness ([\citen{goedel-paper}, Appendix ``G\"odel II in full'']) & 1{,}232 \\
transparent & \textsf{lwD} / \textsf{goedelWS} / \textsf{goedelW} & $D_{\Tru}$, $S_{\Tru}$, $L_{\Tru}$ & 1{,}261 / 1{,}262 / 1{,}263 \\
transparent & \textsf{liarRec} & $L_{\Rcd}$ & 1{,}080 \\
transparent & \textsf{iHdStep} & head-step scan $\mathit{HS}$ & 1{,}227 \\
transparent & \textsf{iDep} / \textsf{OkI} / \textsf{OkU} & $\Dep$, certificates & 937 / 946 / 950 \\
transparent & \textsf{iLFA} / \textsf{iLFI} & level sentences & 1{,}130 / 1{,}128 \\
transparent & \textsf{iLvClassI} / \textsf{iLvClassU} & $\LvC$, $\LvC$ at a term level & 1{,}283 / 1{,}281 \\
control & \textsf{iPrv} / \textsf{chkIdx} & $\Prv$, index of $\ChkE$ & 983 / 974 \\
\hline
\multicolumn{4}{@{}l}{\emph{fourth run, end of \texttt{Goedel\_D3}}} \\
transparent & \textsf{d3SigmaOf} / \textsf{d3SigmaAt} / \textsf{d3SigmaIdx} / \textsf{d3CertIdx} & decided D3 & 1{,}175 / 1{,}174 / 1{,}173 / 1{,}172 \\
control & \textsf{prDecIdx} & & 962 \\
\hline
\multicolumn{4}{@{}l}{\emph{fifth run, end of \texttt{Goedel\_D3\_Uniform}}} \\
transparent & \textsf{uSigIdx} / \textsf{uCertIdx} / \textsf{prCU} / \textsf{prCUIdx} & D3 over the code & 1{,}174 / 1{,}171 / 1{,}110 / 1{,}105 \\
control & \textsf{prDecIdx} & & 962 \\
\hline
\multicolumn{4}{@{}l}{\emph{sixth run, end of \texttt{Goedel\_D2}}} \\
transparent & \textsf{d2Idx} / \textsf{d2M} & decided D2 & 826 / 838 \\
control & \textsf{prDecIdx} & & 962 \\
\hline
\multicolumn{4}{@{}l}{\emph{seventh run, end of \texttt{Goedel\_Cut}}} \\
transparent & \textsf{cut1Idx} / \textsf{cut1M} / \textsf{cut2Idx} / \textsf{cut2M} & cut transformers & 820 / 832 / 821 / 833 \\
control & \textsf{prDecIdx} & & 962 \\
\hline
\multicolumn{4}{@{}l}{\emph{eighth run, end of \texttt{Goedel\_Cut}}} \\
transparent & \textsf{goedel1Idx} / \textsf{goedel1T} / \textsf{goedel2Idx} / \textsf{goedel2T} & witnesses of [\citen{goedel-paper}, theorem~\textsf{rga\_\allowbreak goedel2\_\allowbreak sentence}] & 1{,}182 / 1{,}183 / 1{,}185 / 1{,}186 \\
control & \textsf{prDecIdx} & & 962 \\
\end{longtable}
\endgroup

The tree $\pi$ is not a separate root, but it lies inside the closure
of $\Cert$.

\paragraph{The code-generation check.}
The three audits above read definitions.  None of these audits sees an index
that a definition chooses between explicit trees by a case split on a
condition that HOL cannot evaluate: the control \textsf{ifCtlIdx},
defined as $\mathsf{if}\ (\forall n{::}\mathsf{nat}.\ n = 0)\
\mathsf{then}\ \mathsf{prSucL}\ \mathsf{else}\ \mathsf{prSucR}$ in
\texttt{Goedel\_II\_PR}, passes all three audits, and a block after the
control's definition checks that the control passes them.  A fourth
check, the ML structure
\textsf{Code\_Audit} at the end of \texttt{ConOGA/ConOGA\_Audit},
closes this gap.  The code-generation check asks Isabelle's code
generator to translate the
index-valued constants of the transparent closures into Standard ML
from their defining equations (\textsf{Code\_Target.produce\_code},
which translates and serializes; nothing is executed).  The check's roots are
those of all eight runs of the choice audit, 68 root names with 86
root constants, and their closures contain 625 constants of the index
type \textsf{pridx}.  Among them are the fifteen index-valued roots:
the trees of the witness $\sigma$ and of its premise certificate, of
EXPL and ex falso, of the decided D3 (two), of the decided D3 over the
code (three), of the decided D2, of the two cut transformers and of
the witnesses of [\citen{goedel-paper}, theorem~\textsf{rga\_goedel2\_sentence}], and the index of the
whole-proof $\forall$-free scan.

The constants of the development are defined inside the locale
\textsf{rga}, whose definitions Isabelle does not register as code
equations, and a datatype declared inside a locale gets no
code-datatype setup.  These two gaps are packaging artifacts, not
non-executability.  So the check works on a copy of the theory that is
never installed.  On the copy the check registers the constructors of each
free datatype of the development that has no code-datatype setup (twelve locale
datatypes) as code constructors, declares the recursion equations, or
else the definition, of each constant of the development in the
closures as code equations (908 of them), and declares further
definitions on demand when the translation reports a constant without
code equations.  The 625 constants then translate in one translation,
except five internal constants of the definition packages that the
generated code does not call and that the check skips: the function
package's \textsf{ifNestIdx\_sumC}, \textsf{orNestIdx\_sumC} and
\textsf{bldInit\_sumC}, and the morphisms \textsf{Abs\_pridx} and
\textsf{ctor\_pridx} of the index type.  The generated code contains
no aborting call (\texttt{raise Fail}).  The fifteen roots are among
the 620 translated constants.  The check fails the build unless the
translation succeeds with no aborting call, at least 100 constants are
translated, and each control fails.  The controls fail as they must:
\textsf{ifCtlIdx} (the universal quantifier over the naturals has no
code), \textsf{wlCtlIdx} (\textsf{the\_inv\_into}, defined by
\textsf{THE}), \textsf{udCtlIdx} (\textsf{undefined}: it translates,
with one aborting call, which counts as failure) and
\textsf{prDecIdx} (\textsf{priof}, defined by \textsf{SOME}).  The
check takes about a second and a half.

The check's notion of explicitness is that the code generator can
translate the function from its defining equations into ML with no
aborting call.  The check relies on the correctness of Isabelle's code
generator, that the generated code computes the equations it was
given; the equations themselves are theorems of the development.  The
check sees a conditional on an undecidable condition, \textsf{SOME} and
\textsf{THE}, \textsf{undefined}, an unbounded \textsf{LEAST} over the
naturals and quantifiers over infinite types.  A function that is
computable but slow, such as a bounded search, passes, as it should.

\paragraph{The axiom audit.}
A block of \texttt{ConOGA/ConOGA\_Audit} takes every axiom of the
theory's ancestry (\textsf{Theory.all\_axioms\_of}) that is not an
axiom of Isabelle/HOL's \texttt{Main}.  Each must be a definition
(the axiom of a specification that the theory records in
\textsf{Theory.defs\_of}), the axiom of a type definition
(\textsf{type\_definition}), or the arity axiom that declaring a type
constructor adds for the class \textsf{type}.  The build fails
otherwise, and fails unless there are at least 1{,}000 such axioms and
\texttt{Main}'s own non-definitional axioms, among them
\textsf{HOL.ext} and \textsf{Hilbert\_Choice.someI}, are found.  Of the
9{,}309 axioms, 3{,}451 are \texttt{Main}'s, 20 of them
non-definitional (the axioms of HOL and Pure, \textsf{someI}, the
axioms of the natural numbers' representation and of set
comprehension).  The other 5{,}858 axioms, from 807 theory and locale name
spaces that include the parts of HOL-Library beyond \texttt{Main} and
the whole development, are 5{,}724 definitions, 67 type definitions
and 67 type-constructor arities; none is of another kind.

\paragraph{The dependency checks.}
A second ML structure at the end of \texttt{Goedel\_II},
\textsf{Dep\_Audit}, reads the theorem-dependency graph that the
kernel records in the proof bodies.  The \emph{cone} of a theorem is
every theorem node reachable from its proof body; a theorem
\emph{cites} another theorem directly when that other theorem is a named node reachable
through unnamed nodes only.  Three checks use this graph, and each check fails the
build on a violation.
\begin{itemize}
\item In \texttt{ConOGA/ConOGA\_Audit}, no soundness induction lies
  in the cone of an internal theorem of Section~\ref{sec:con}
  (\S\ref{sec:con-summary}).  The roots are
  \textsf{rga\_con\_qf\_int}, \textsf{rga\_sound\_qf\_int2},
  \textsf{rga\_con\_frag\_int}, \textsf{rga\_sound\_frag\_int},
  \textsf{rga\_qf\_in\_frag\_int},
  \textsf{rga\_con\_qf\_via\_frag\_int},
  \textsf{rga\_con\_lv\_under\_lf},
  \textsf{rga\_con\_lv\_under\_lf\_all}, \textsf{rga\_con\_lv\_oga},
  \textsf{rga\_con\_under\_res}, \textsf{rga\_sound\_}\allowbreak\textsf{under\_res},
  \textsf{rga\_sound\_}\allowbreak\textsf{under\_res\_v0}, \textsf{iFrag\_}\allowbreak\textsf{iResGuard}
  and \textsf{iFrag\_}\allowbreak\textsf{noswap\_iRes}, with cones of 62{,}532 to
  105{,}271 theorem nodes.  Two further roots are the per-level forms of the internal
  containments of Theorem~\ref{thm:con-lv-contain},
  \textsf{rga\_qf\_in\_lv\_int} and \textsf{rga\_frag\_in\_lv\_int}.
  Nine more roots are the level-free results of \S\ref{sec:con-lv}: the
  \oga derivations of the level sentences, \textsf{o\_iLFA} and
  \textsf{o\_iLFI}; the level theorem in its three variants,
  \textsf{rga\_con\_lvi\_}\allowbreak\textsf{under\_lf},
  \textsf{rga\_con\_lvi\_}\allowbreak\textsf{under\_lf\_all}
  and \textsf{rga\_con\_lvi\_}\allowbreak\textsf{under\_lfa\_all}; its \oga form and soundness
  companion, \textsf{rga\_con\_lvi\_oga} and \textsf{rga\_sound\_}\allowbreak\textsf{lvi\_oga};
  and the containments \textsf{rga\_qf\_}\allowbreak\textsf{in\_lvi\_int} and
  \textsf{rga\_frag\_}\allowbreak\textsf{in\_lvi\_int}.  These eleven roots have cones
  of 59{,}547 to 97{,}249 nodes.  Twenty-one more roots were added after these eleven:
  the reach-guarded theorem of Theorem~\ref{thm:con-e1} and its
  siblings, per level (\textsf{rga\_con\_}\allowbreak\textsf{stat\_under\_dwfd\_lf},
  \textsf{rga\_con\_}\allowbreak\textsf{stat\_oga\_under\_dwfd},
  \textsf{rga\_con\_}\allowbreak\textsf{lv\_under\_lf\_via\_stat},
  \textsf{rga\_con\_}\allowbreak\textsf{lv\_oga\_via\_stat}) and over $\LvC$
  (\textsf{rga\_con\_}\allowbreak\textsf{statI\_under\_dwfd\_lf},
  \textsf{rga\_con\_}\allowbreak\textsf{statI\_oga\_under\_dwfd},
  \textsf{iHypLvI\_}\allowbreak\textsf{iDynWFD},
  \textsf{rga\_con\_}\allowbreak\textsf{lvi\_under\_lf\_via\_stat},
  \textsf{rga\_con\_}\allowbreak\textsf{lvi\_oga\_via\_stat}); the internal
  Lemma~Q (\textsf{lemQ\_int}, Theorem~\ref{thm:con-lemq}) and
  \textsf{iHdStep\_of\_lvi}; the class at a term level and its
  consistency sentences (\textsf{lvU\_in\_lvI},
  \textsf{rga\_con\_}\allowbreak\textsf{lvu\_under\_lf\_all},
  \textsf{rga\_con\_}\allowbreak\textsf{lvu\_oga}); the tree form of the
  $\forall$-free containment (\textsf{qf\_accepted\_}\allowbreak\textsf{in\_frag\_tree});
  the residual at rejected codes (\textsf{iRes\_of\_}\allowbreak\textsf{rejected\_num});
  and, on the G\"odel side, \textsf{negInstC\_num}, the decided D3 as
  one sentence over the code (\textsf{D3\_uniform},
  \textsf{D3\_uniform\_num}) and the decided D2 (\textsf{D2\_decided},
  \textsf{d2M\_chk}).  These 46 roots have cones of 24{,}216 to
  105{,}271 nodes.  Forty-one more roots were added with the entailment
  results of the sections ``G\"odel's proof in four forms'' and
  ``G\"odel's theorems for \rga'' of [\citen{goedel-paper}]: the
  entailment kit (\textsf{sigma\_}\allowbreak\textsf{deduction\_theorem},
  \textsf{sigma\_reductio}, \textsf{sigma\_contrapos},
  \textsf{sigma\_}\allowbreak\textsf{refute\_bot},
  \textsf{deduction\_}\allowbreak\textsf{theorem\_exact},
  \textsf{dec\_contrapos}, \textsf{ent\_of\_ctx\_rule}); the G\"odel ingredients as entailments
  (\textsf{iPrc\_transfer}, \textsf{iPrc\_transfer2},
  \textsf{goedel\_diag\_ent}, \textsf{goedel\_}\allowbreak\textsf{diag\_ent\_all},
  \textsf{D3\_ent}, \textsf{D3\_at\_ent}, \textsf{D3G\_ent},
  \textsf{expl\_ent}, \textsf{exf\_ent}, \textsf{exfIdx\_chk\_sf},
  \textsf{con\_notprc\_}\allowbreak\textsf{goedelG\_ent},
  \textsf{liar\_diag\_ent}); internal completeness and the truth
  clauses as entailments (\textsf{iHasCert\_}\allowbreak\textsf{of\_iWtrue},
  \textsf{complete\_ent}, \textsf{kf\_dneg\_ent}, \textsf{kf\_or\_ent});
  the search for a proof of $\z = \tofn 1$
  (\textsf{prc\_stops\_ent}, \textsf{not\_stops\_con\_ent},
  \textsf{not\_stops\_}\allowbreak\textsf{goedelG\_ent},
  \textsf{not\_stops\_}\allowbreak\textsf{notprc\_goedelG\_ent},
  \textsf{not\_stops\_}\allowbreak\textsf{notprc\_con\_ent}); and the cut
  transformers and what they give (\textsf{cut1M\_chk},
  \textsf{cut2M\_chk}, \textsf{ent\_nec}, \textsf{ent\_nec\_contra},
  \textsf{ent\_nec2}, \textsf{D2\_ent},
  \textsf{rga\_goedel2\_}\allowbreak\textsf{internal\_contra} (a root
  under both of its names; Appendix~\ref{app:endpoints}),
  \textsf{rga\_goedel2\_}\allowbreak\textsf{internal},
  \textsf{goedel\_diag\_}\allowbreak\textsf{neg\_ent},
  \textsf{not\_goedelG\_}\allowbreak\textsf{not\_con\_ent},
  \textsf{formalized\_}\allowbreak\textsf{loeb\_ent},
  \textsf{formalized\_}\allowbreak\textsf{loeb\_0eq1\_ent}).  Their cones
  have 12{,}932 to 86{,}520 nodes, in 53 to 146 theories; the smallest
  is that of \textsf{ent\_of\_ctx\_rule}.  Sixty-three more roots were added
  with the results of \S\ref{sec:pillars-tarski},
  {}[\citen{goedel-paper}, \S\,``The formalized theorems, as entailments''] and the comparison with Beklemishev and
  Shamkanov (\S\ref{sec:related}): the decided-antecedent forms and G\"odel's
  proof read with entailments (\textsf{prc\_sentence\_iff},
  \textsf{rga\_goedel2\_}\allowbreak\textsf{sentence},
  \textsf{goedel1\_sentence}, \textsf{goedel1\_}\allowbreak\textsf{sentence\_diag},
  \textsf{formalized\_}\allowbreak\textsf{loeb\_sentence},
  \textsf{goedelG\_under\_}\allowbreak\textsf{ConRGA\_textbook}); the
  explicit witnesses (\textsf{cmp\_diag\_beta},
  \textsf{goedel1Idx\_chk}, \textsf{goedel2Idx\_chk},
  \textsf{chk\_imp\_}\allowbreak\textsf{sentence},
  \textsf{chk\_imp\_}\allowbreak\textsf{matrix\_B},
  \textsf{goedel1\_}\allowbreak\textsf{sentence\_witness},
  \textsf{goedel2T\_}\allowbreak\textsf{sentence},
  \textsf{rga\_goedel2\_}\allowbreak\textsf{sentence\_witness}); the
  explosive T-biconditional (\textsf{tbicond\_}\allowbreak\textsf{liar\_explosive},
  \textsf{tbicond\_}\allowbreak\textsf{liar\_both}); D1 without soundness
  and the interderivabilities (\textsf{prc\_of\_}\allowbreak\textsf{derivable},
  \textsf{notcon\_}\allowbreak\textsf{prc\_con\_ent}, \textsf{goedelG\_con\_ent},
  \textsf{goedelG\_con\_}\allowbreak\textsf{interderivable},
  \textsf{notcon\_prc\_con\_}\allowbreak\textsf{interderivable}); the
  derivability conditions (\textsf{exf\_int\_sf},
  \textsf{chk\_derivability\_}\allowbreak\textsf{conditions}); the
  Beklemishev--Shamkanov conditions and theorem
  (\textsf{bs\_structure\_}\allowbreak\textsf{conditions},
  \textsf{bs\_C1\_anti\_sigma}, \textsf{bs\_C1\_anti\_dec},
  \textsf{bsSD\_prc}, \textsf{bs\_aps\_sigma},
  \textsf{bs\_jeroslow\_fp}, \textsf{bs\_theorem1},
  \textsf{bs\_theorem1\_i\_rule}, \textsf{neg\_goedelG\_dec},
  \textsf{bs\_goedel2}, \textsf{bs\_incon\_con}); the limits of
  establishment and the ordinary Liar
  (\textsf{ir\_A\_cert\_affirm}, \textsf{ir\_inRangePf\_}\allowbreak\textsf{of\_valid},
  \textsf{ir\_certTrue\_}\allowbreak\textsf{intro},
  \textsf{iWtrue\_univ\_}\allowbreak\textsf{of\_chk01},
  \textsf{notcon\_tru\_}\allowbreak\textsf{univ\_ent},
  \textsf{universal\_}\allowbreak\textsf{unestablished\_needs\_con},
  \textsf{refuted\_universal\_}\allowbreak\textsf{unestablished\_needs\_con},
  \textsf{liar\_never\_}\allowbreak\textsf{decided},
  \textsf{liar\_never\_}\allowbreak\textsf{affirmed},
  \textsf{neg\_liar\_never\_}\allowbreak\textsf{affirmed}); refuted
  sentences under consistency (\textsf{iPrc\_of\_}\allowbreak\textsf{iHasCert\_code},
  \textsf{iPrc\_of\_iWtrue}, \textsf{tru\_prc\_ent},
  \textsf{refuted\_prc\_}\allowbreak\textsf{notcon\_ent},
  \textsf{refuted\_tru\_}\allowbreak\textsf{notcon\_ent},
  \textsf{con\_refuted\_}\allowbreak\textsf{unprovable},
  \textsf{con\_refuted\_}\allowbreak\textsf{unestablished},
  \textsf{refuted\_chk\_}\allowbreak\textsf{inconsistent},
  \textsf{refuted\_affirmed\_}\allowbreak\textsf{inconsistent},
  \textsf{refuted\_universal\_}\allowbreak\textsf{unestablished\_iff\_con}); and
  $\Omega$ without soundness (\textsf{omega\_never\_}\allowbreak\textsf{affirmed\_sf},
  \textsf{omega\_}\allowbreak\textsf{unestablished\_sf}); and tracking
  under $\ConR$ (\textsf{unestablished\_}\allowbreak\textsf{of\_unprovable},
  \textsf{con\_tschema\_}\allowbreak\textsf{refuted},
  \textsf{con\_prc\_}\allowbreak\textsf{schema\_decided},
  \textsf{con\_prc\_}\allowbreak\textsf{bool\_decided},
  \textsf{con\_goedelG\_}\allowbreak\textsf{derivable\_unestablished},
  \textsf{con\_tschema\_}\allowbreak\textsf{goedelG\_refuted},
  \textsf{con\_con\_}\allowbreak\textsf{unestablished}).  Their cones
  have 22{,}542 to 94{,}039 nodes, in 73 to 143 theories.  Nine more roots were added with the entailment forms
  of \S\ref{sec:pillars-tarski} and [\citen{goedel-paper}, \S\,``The formalized theorems, as entailments'']: local
  reflection at $\z = \tofn 1$ in every context
  (\textsf{reflect\_ent\_}\allowbreak\textsf{iff\_con\_ctx},
  \textsf{con\_reflect\_}\allowbreak\textsf{0eq1\_ent}); the G\"odel
  sentence from the consistency hypothesis
  (\textsf{goedelG\_under\_}\allowbreak\textsf{ConRGA}); verdict tracking
  as one conjunction (\textsf{con\_goedelG\_}\allowbreak\textsf{and\_unestablished});
  and refuted universals as entailments
  (\textsf{universal\_unestablished\_}\allowbreak\textsf{needs\_con\_ctx},
  \textsf{universal\_unestablished\_}\allowbreak\textsf{con\_ent},
  \textsf{refuted\_universal\_}\allowbreak\textsf{unestablished\_eqv\_con},
  \textsf{refuted\_universal\_}\allowbreak\textsf{unestablished\_iff\_con\_ctx},
  \textsf{con\_eqv\_}\allowbreak\textsf{not\_true\_con}), with cones of
  25{,}913 to 94{,}260 nodes.  Four more roots are
  the generic T-biconditional of Theorem~\ref{thm:pl-tb-generic}
  (\textsf{tbicond\_liar\_}\allowbreak\textsf{explosive\_gen},
  \textsf{tbicond\_liar\_}\allowbreak\textsf{explosive\_pred},
  \textsf{pred\_liar\_eq},
  \textsf{tbicond\_liar\_}\allowbreak\textsf{explosive\_diag}), with
  cones of 18{,}881 to 30{,}198 nodes.  So the check has 163 roots in
  all, with cones of 12{,}932 to 105{,}271 nodes in 53 to 159
  theories.  \textsf{exfIdx\_chk\_sf} re-proves the
  ex-falso law of [\citen{goedel-paper}, theorem~\textsf{exfIdx\_chk}] without soundness: the
  acceptance of the fixed proof of $\neg(\z = \tofn 1)$ inside
  \textsf{exfIdx\_chk} reaches \textsf{prf\_owsjt} through the
  checker's faithfulness, and the new proof evaluates that acceptance
  as a closed computation instead (\textsf{iChk\_accept\_sf},
  \textsf{iChk\_proves\_sf}).  The numeral instances of the uniform-soundness
  sentence (\textsf{uniform\_soundness\_}\allowbreak\textsf{instances}) are not
  a root: they follow from those of $\ConR$, which rest on soundness.
  The forbidden theorems are \rga's
  soundness induction \textsf{prf\_owsjt}, \rga's $\omega$-soundness
  induction \textsf{prf\_oomjt} and \oga's
  \textsf{oga\_prf\_oomjt}.  The controls must reach the forbidden theorems:
  \textsf{rga\_provable\_wosr}, \textsf{rga\_val\_iff},
  \textsf{rga\_open\_iff}, \textsf{rga\_consistent},
  \textsf{rga\_goedel2}, \textsf{rga\_afree\_normal\_chk} and
  \textsf{qf\_accepted\_in\_frag} reach \textsf{prf\_owsjt};
  \textsf{rga\_omega\_sound} and
  \textsf{rga\_omega\_true\_consistent} reach \textsf{prf\_oomjt};
  \oga's \textsf{omega\_sound} reaches \textsf{oga\_prf\_oomjt}.
  The check also fails if a root's cone has fewer than 1{,}000 nodes,
  or if some forbidden theorem is reached by no control.  A second
  run over the same roots forbids the HOL shadow models of
  Appendix~\ref{app:traps}: the four base folds \textsf{path\_fold},
  \textsf{path\_fold\_tot}, \textsf{qf\_fold} and
  \textsf{shadow\_fold}, through which every shadow soundness or
  consistency theorem passes, with the controls
  \textsf{path\_consistent}, \textsf{shadow\_consistent\_}\allowbreak\textsf{of\_tot},
  \textsf{shadow\_sound\_qf} and \textsf{shadow\_sound}, one for each
  fold.  Some cones contain lemmas of the shadow theories that are not
  model theorems (shape lemmas), and the build prints, for each root,
  how many such nodes its cone has.  The folds
  \textsf{fold\_consistent\_ctx} and \textsf{fold\_sound\_ctx} of
  Appendix~\ref{app:fold} lie in the cones of the consistency theorems
  proved through the internal fold, 20 and 4 roots respectively.  The
  conclusions of these two folds are \rga derivations, built by the fold
  inside \rga, so the two folds are the internal argument itself and
  are not forbidden.  The
  build runs at proofs level~0, at which the kernel still records every
  theorem dependency (only proof terms depend on the level), and the
  walk reads each proof body through \textsf{Thm.proof\_body\_of},
  which joins the theorem's pending proofs, so the graph it reads is
  complete.  All the node counts of this appendix come from one build.
  These counts include unnamed nodes, whose number varies
  slightly between builds of the same theory text (by at most
  0.15\,\% between the builds measured, by up to 154 nodes); the
  counts of named nodes and of theories do not vary.

  Two further passes over the same roots read statements rather than
  names (\textsf{Dep\_Audit.shape} and
  \textsf{Dep\_Audit.induction\_report}).  The \emph{shape check}
  splits the statement of a node into premise--conclusion pairs:
  the meta-level premises and conclusion, then through HOL's $\forall$,
  $\longrightarrow$, $\wedge$ and both directions of $\longleftrightarrow$.
  A \emph{derivability constant} is one of the eight constants named
  \textsf{derives} (of \rga, \oga, PGA, BGA and the classical system),
  \textsf{proves}, \textsf{entails} or \textsf{valid}, mentioned
  directly in a premise.  A \emph{semantic constant} is one of
  \textsf{wosr}, \textsf{whev}, their generic counterparts
  \textsf{gwosr} and \textsf{gwhev}, \textsf{osr}, \textsf{hev},
  \textsf{omv}, \textsf{omh} and the shadow models' \textsf{qfsnd},
  \textsf{famTrue} and \textsf{siM}, or any constant whose definitional
  closure contains one, through constants with a single specification;
  so wrappers such as \textsf{owsjt} count without being listed.  The
  derivability side is read by direct mention, not by closure, because
  \rga's semantics is defined from its proof system: a universal's value
  needs a certificate, an \rga proof of its open instance, so every
  semantic constant has \textsf{derives} in its closure, and a closure
  test would flag every semantic lemma.  A node is
  \emph{soundness-shaped} when a premise mentions a derivability
  constant and the node's conclusion a semantic constant.  The check fails the
  build if a root or a named node of the union of the root cones is
  soundness-shaped, except the thirteen certificate clauses of the
  semantics, the semantics' introduction rules at the universal
  quantifier, or if a control
  has no soundness-shaped node outside that list.  The controls are
  \textsf{rga\_provable\_wosr}, \textsf{rga\_consistent},
  \textsf{rga\_omega\_sound} and \oga's \textsf{omega\_sound}; each is
  flagged, through itself or nodes such as \textsf{prf\_owsjt},
  \textsf{contains\_owsjt} and \textsf{prf\_oomjt}.  Eight of the
  thirteen clauses belong to \rga's two evaluations
  (\textsf{osr\_allT}, \textsf{osr\_allF}, \textsf{hev\_osr.intros},
  \textsf{wosr\_allT}, \textsf{wosr\_allF}, \textsf{whev\_wosr.intros},
  \textsf{wosr\_allT\_i} and \textsf{wosr\_allF\_i}).  The other five are
  their generic counterparts in the evaluation parametric in the
  system whose proofs serve as certificates (\textsf{gwosr\_allT},
  \textsf{gwosr\_allF}, \textsf{gwhev\_gwosr.intros},
  \textsf{gwosr\_allT\_i} and \textsf{gwosr\_allF\_i},
  \texttt{RGA/RGA\_GSem}).  The facts about \textsf{wosr} proved as
  instances of that evaluation reach them through the equation identifying
  \textsf{wosr} with that evaluation's instance at \rga's proofs
  (\textsf{gR\_rw\_whev\_wosr}).  The union of the
  163 cones has 159{,}317 nodes.
  Its soundness-shaped named nodes are eleven of the thirteen clauses,
  all but \textsf{wosr\_allT\_i} and \textsf{wosr\_allF\_i}.  Eighteen
  unnamed nodes, ten in \texttt{RGA/RGA\_GSem} and eight in
  \texttt{RGA/RGA\_Sound}, are soundness-shaped too: internal forms of
  the same clauses that the inductive package generates.  The check
  covers named nodes and reports the unnamed count.

  The \emph{induction report} lists the theorems of the union cone that
  cite, directly or through unnamed nodes, an induction rule of
  \textsf{derives}, \textsf{proves}, \textsf{entails} or
  \textsf{valid}; the build requires that there be none, and there are
  none.  \rga's soundness theorem \textsf{prf\_owsjt} is not proved by
  such a rule but by induction on the proof list, so the report also
  lists the theorems that cite list induction and whose statement
  mentions a derivability constant: 96 nodes, 24 of them named.  All 24
  are proof-list bookkeeping and syntactic transformations
  (\textsf{contains\_proves}, \textsf{valid\_appendI},
  \textsf{prf\_strip} of \rga and \oga, \textsf{valid\_embed}, the
  or-nest and if-nest lemmas, the spine lemmas of the emission), and the
  unnamed ones are their internal steps; none of the 96 nodes has a semantic conclusion,
  which the shape check would flag.

  Together the name-based check and the shape check catch a re-proof of
  soundness as a named lemma, under any name and with any wrapper over
  the semantics.  The shape check does not count the contrapositive
  form, a negated semantic premise with a negated derivability
  conclusion, as soundness-shaped, and a re-proof inside one proof, with
  no named intermediate statement, is visible only in the list of
  list-induction citers.

  What the check certifies is that \rga does the reasoning: no HOL
  proof that builds one of the root derivations passes through ``true,
  hence derivable'' by way of \rga's soundness.  The check allows completeness
  (\textsf{rga\_complete\_val}, \textsf{rga\_complete\_open}), which
  turns a HOL evaluation of a closed term into a derivation and adds no
  reasoning of its own, as in Theorem~\ref{thm:in-cmp-num}; where a
  construction needs the decider's verdict at $\z = \tofn 1$ or at
  $\z\jnat$, the construction computes that verdict.  The name-based part forbids the three
  soundness inductions and the four shadow folds; the shape check and
  the induction report above cover a re-proof under another name, within
  the limits stated there.  Some results
  of Section~\ref{sec:con} are deliberately not roots.
  Theorem~\ref{thm:con-normal} and the meta corollaries of the
  $\forall$-free theorem (Appendix~\ref{app:fold}) use soundness, and
  the check confirms the dependence for Theorem~\ref{thm:con-normal}
  and \textsf{qf\_accepted\_in\_frag}, two of the check's controls.  Most facts
  at closed codes (Theorem~\ref{thm:con-witness}(2),
  \textsf{fwP\_in\_lvi}, the memberships of $P_3$ and $P_{\mathrm R}$,
  Theorem~\ref{thm:con-lv-contain}(5),
  Theorem~\ref{thm:con-lv-closed}(1) and
  Theorem~\ref{thm:con-res-frag}(3)) are HOL statements about closed
  codes, and these facts use soundness too: to close a scan at a closed
  code their proofs apply the disjunction property and \rga's consistency, two HOL
  consequences of soundness (\textsf{nodeScan\_closedI}).  The numeral
  instances of Theorem~\ref{thm:con-qf-in-frag} have two proofs.  The
  first proof (\textsf{qf\_accepted\_in\_frag}) closes the scans in this
  way; the second proof instantiates the internal tree theorem of
  Appendix~\ref{app:fold} in the empty context, is a root, and uses no
  soundness (\textsf{qf\_accepted\_}\allowbreak\textsf{in\_frag\_tree}).
  The internal containments (Theorems~\ref{thm:con-qf-in-frag}
  and~\ref{thm:con-lv-contain}) and the second derivation of the
  $\forall$-free theorem (Theorem~\ref{thm:con-qf-via-frag}) close their
  scans by introduction rules inside \rga instead.  Faithfulness of the
  checker (Theorem~\ref{thm:in-chk}) rests on soundness, as does every
  non-derivability claim of the paper.  That dependence concerns the
  adequacy of the statements, a meta-level matter, and is unavoidable:
  a claim that \rga does not derive some sentence implies \rga's
  consistency.  Every internal theorem is a HOL statement that a
  derivation exists; the check does not extract one concrete derivation
  (Section~\ref{sec:concl}).
\item In \texttt{ConOGA/ConOGA\_Audit}, the G\"odel results built at the
  checker use neither the decider nor the truth predicate.  The roots
  are 31 theorems: \textsf{rga\_goedel2}, \textsf{goedelG\_of\_con},
  \textsf{D3G\_goedelSigma}, \textsf{goedelG\_under\_}\allowbreak\textsf{ConRGA}
  and the internal second theorem with the results derived from it
  (\textsf{con\_notprc\_}\allowbreak\textsf{goedelG\_ent},
  \textsf{prc\_con\_ent\_bot}, \textsf{rga\_goedel2\_internal},
  \textsf{rga\_goedel2\_}\allowbreak\textsf{internal\_contra},
  \textsf{rga\_goedel2\_}\allowbreak\textsf{sentence},
  \textsf{rga\_goedel2\_}\allowbreak\textsf{sentence\_witness},
  \textsf{goedel2T\_}\allowbreak\textsf{sentence},
  \textsf{goedelG\_con\_}\allowbreak\textsf{interderivable},
  \textsf{notcon\_prc\_con\_}\allowbreak\textsf{interderivable},
  \textsf{not\_stops\_}\allowbreak\textsf{notprc\_con\_ent}); the
  textbook route (\textsf{goedelG\_under\_}\allowbreak\textsf{ConRGA\_textbook},
  \textsf{rga\_goedel2\_}\allowbreak\textsf{textbook}); local reflection
  at $\z = \tofn 1$ (\textsf{reflect\_ent\_}\allowbreak\textsf{iff\_con\_ctx},
  \textsf{con\_reflect\_}\allowbreak\textsf{0eq1\_ent}); the ingredients
  (\textsf{D3\_uniform}, \textsf{D3G\_ent}, \textsf{D2\_ent},
  \textsf{expl\_ent}, \textsf{ent\_nec}, \textsf{ent\_nec2},
  \textsf{chk\_derivability\_}\allowbreak\textsf{conditions});
  \textsf{goedel1\_}\allowbreak\textsf{sentence\_witness},
  \textsf{formalized\_}\allowbreak\textsf{loeb\_ent} and
  \textsf{bs\_goedel2}; the numeral lemma \textsf{isNum\_subM},
  that substitution maps numeral codes to numeral codes; and explicit
  reflection at each numeral (\textsf{chk\_reflect\_num},
  \textsf{chk\_reflect\_num\_ent}, [\citen{goedel-paper}, theorem~\textsf{chk\_reflect\_num}]),
  which uses soundness but not the decider.  Two checks
  run on each root.  The dependency check forbids the search witness's
  D3G (\textsf{D3G\_}\allowbreak\textsf{goedelSigmaSearch}), the
  decided-search theorem behind it (\textsf{D3Gdec\_holds}) and the
  decider's internal completeness
  (\textsf{iHasCert\_of\_}\allowbreak\textsf{iPrDecTrueAt}), through
  which every use of $\Tru \vdash \PrC$ passes; the dependency check's controls
  \textsf{goedelG\_of\_con\_}\allowbreak\textsf{via\_search},
  \textsf{rga\_goedel2\_}\allowbreak\textsf{via\_search} and
  \textsf{iHasCert\_of\_iWtrue} reach these forbidden theorems.  A statement check reads
  the statement of every node of each cone, named or unnamed: no node
  may mention the truth predicate, the affirmation test, the decided
  value, the search, the HOL-level decider or the constants of the
  compiled decider's table.  The constants are resolved by base name,
  and the check fails if one of these constants is missing, so a renaming breaks
  the check loudly.  The statement check's controls, flagged, are the
  search route, the decider's
  internal completeness at $\Tru$ and \textsf{ncG\_zero}, whose
  statement reads the table's value cell at the code of zero.  The
  statement check excludes every lemma about the verdicts, whatever its
  name.  What the statement check cannot see is a decider constant that occurs only as a
  term inside one proof step with no lemma about it.  A third run, with
  23 roots, forbids only the search: the five theorems that track \rga's
  verdicts under $\ConR$ state $\Tru$ and use the decider through $\Tru$,
  the two theorems behind Theorem~\ref{thm:b-no-moorean} use NEC and consistency,
  and the sixteen meta theorems that cite G\"odel's second theorem (among them
  \textsf{ConRGA\_undecided}, \textsf{ConRGA\_omega\_}\allowbreak\textsf{incomplete},
  \textsf{no\_formalized\_}\allowbreak\textsf{goedel2\_nonprenex} and
  \textsf{prenex\_extraction\_}\allowbreak\textsf{fails\_con}) rest on
  soundness.  Of the 17{,}548 theorems of the locale \textsf{rga} in the
  final session's heap, exactly five have the search in their cone, all
  of them search routes by design.
\item At the end of \texttt{Goedel\_II\_PR}, in the cones of
  \textsf{rga\_goedel2} and \textsf{rga\_goedel2\_via\_search}, the
  theorems of \texttt{Goedel\_II} and \texttt{Goedel\_II\_PR} that
  cite one of the development's ten implication-introduction rules or
  four case-split rules are exactly the expected ones, and EXPL and
  D3G cite only \textsf{ir\_impI\_bool}
  ([\citen{goedel-paper}, \S\,``The assembly'']).  Over every theory these two cones
  reach (127 and 139 theories), no theorem cites
  \textsf{ir\_impI\_provable}, whose antecedent need only be provable,
  or \textsf{ir\_impI\_Bw}.
\end{itemize}
Walking a cone of 100{,}000 nodes takes well under a second: the 173
cones of the first check (163 roots and 10 controls) take 7.0 seconds
together.  The checks
were also run, in scratch theories, on inputs that must fail (a root
that uses soundness, a control that does not, a wrong set of citing
theorems), and failed.

The second conjunct of \textsf{goedelSigma\_is\_pr} states
that the witness's denotation is PR.  That conjunct holds of every
index (\textsf{pr2fi\_n2pr}), so it does not show that the witness
contains no index chosen by \textsf{priof}.  The build-time check does.

\paragraph{Cone sizes.}
The cone of G\"odel's second theorem ([\citen{goedel-paper}, theorem~\textsf{rga\_goedel2}]) has
85{,}966 nodes (11{,}929 named) in 130 theories, and the cone of
its search route (\textsf{rga\_goedel2\_via\_search},
{}[\citen{goedel-paper}, Appendix ``The search route'']) 102{,}222 nodes in 142
theories; both cones contain the whole trunk.  The cone of the internal second
theorem ([\citen{goedel-paper}, theorem~\textsf{rga\_goedel2\_internal\_contra}]) has 81{,}534 nodes
(11{,}652 named) in 134 theories.  This cone is the smallest of the three cones,
since it contains no soundness theorem, and spans four more theories
than the headline cone.  The internal theorem's cone contains the
one-premise cut transformer,
necessitation, the decided D3, D3G and EXPL, all at the checker.
Measured in theory lines, the lines of the named items in the cone
without the grounded-deduction base, the headline cone has 62{,}399
lines, the search route's 98{,}946 and the internal theorem's 60{,}181.
The part of the headline cone in the session \texttt{RGA\_Goedel} has
3{,}144 lines, so about 95\,\% of the cone lies outside that
session.

\paragraph{Fingerprints.}
Two further checks run in the validation script that accepts each
batch of commits, not in the build.  The statement fingerprints hash a
canonical rendering of a theorem's statement (de Bruijn bound
variables, numbered free variables, names without theory qualifiers,
mapped back through the cleanup's renamings); a batch passes only if
every hash is unchanged.  The statement fingerprints cover the 131 endpoint theorems of
Appendix~\ref{app:endpoints}, fingerprinted when the proof campaign
and its cleanup closed, sixteen endpoints added after that (among them
the internal second theorem, the decided-antecedent forms and the
derivability conditions), and the 401 theorem names that
Appendix~\ref{app:endpoints} lists for the paper's numbered results,
read from that appendix's table.
The definitional fingerprints hash, for each endpoint, its
definitional closure: from the constants and types of its statement,
the walk of \textsf{Choice\_Audit} through every specification, kept
to the items defined in the development's own sessions.  Each kept item
contributes the canonical text of its defining axioms, and the
fingerprint is a hash of the sorted per-item hashes, so a difference
names the changed item.  The HOL library is not hashed but pinned: the
Isabelle version and the hashes of the Pure and HOL heaps form one
row, so a library change shows as one visible difference.  The
definitional fingerprints cover the 131 endpoints, the sixteen later
endpoints, and 99 trusted definitions, each in its own row, in the two
tiers of Section~\ref{sec:mech} and a third.  The first tier has 13
rows: \rga's derivability relation and the proof, entailment and
validity relations over it, the coding of terms, $\cmp$, $\ConR$, the
code of $\z = \tofn 1$, $\Chk$, the existential and $\PrC$, $\Tru$ and
$\Dec$.  The second tier has 7 rows: the grounded semantics (\textsf{wosr} and
\textsf{whev}), the $\omega$-semantics (\textsf{omv}), the
well-formedness of terms and contexts, $\forall$-freeness and \oga's
derivability relation.  The third tier has the 75 further constants that
Appendix~\ref{app:endpoints} lists for numbered results.  Four more
rows hold the derivability relations of the other systems of the
development (\pga, BGA, primitive recursive arithmetic and classical
logic with equality), which no statement of this paper uses.  The
statements of the cited theorems also mention auxiliary constants,
witness terms and helpers that Appendix~\ref{app:endpoints} does not
list; those statements are fingerprinted, and the definitions are
given in the cited theories but not separately fingerprinted.  A batch passes only if every
definitional hash equals the recorded baseline, or the batch lists the
change as an approved redefinition.  The walk is by name, locale
parameters are not constants and are not hashed, and the trusted rows
are matched by name.  On an older tree, the comparison first resolves
each trusted name through the cleanup's renamings, so that a name which
a transparent machine took over from its opaque twin
(Section~\ref{sec:mech}) denotes the same machine on both sides.
Compared that way with the freeze of the campaign's proofs, the trusted
definitions that existed there (all but $\PrC$) differ only by the
naming of one subterm of the primitive-recursive index that appends two
coded lists, and the build proves the earlier defining equation from
the current definitions (\textsf{appendEncIdx\_unfold}).

\paragraph{Build and measurement.}
Every figure in this paper was measured at the point named in
Section~\ref{sec:mech} (``Artifact), on one workstation: an Apple M4 Max (10 performance and 4 efficiency
cores, 36\,GB of memory) under macOS, with Isabelle2025-2 (Poly/ML
5.9.2 in its default 32-bit mode, Java heap 12\,GB) and ten threads,
with no other build running and no swapping.  The command
\begin{center}
\small\texttt{isabelle build -b -v -o threads=10 -d hol ConOGA}
\end{center}
builds the whole chain of 13 sessions, from the grounded-deduction
base through the companion formalization and \oga to the ten sessions
of this paper; with the Pure and HOL images prebuilt and every other
session rebuilt, the build took 8 minutes 11 seconds
of wall-clock time (7 minutes 34 seconds summed over the sessions).
The build-time checks of this appendix run inside this build.
  The validation
script, which adds the fingerprints and the oracle audit, passed twice
there.  The development's README describes both commands and the
scripts that measure closures and cones.

\endgroup

\section{The pillars' proofs}
\label{app:pillars}
\begingroup\sloppy

This appendix proves the results of Section~\ref{sec:pillars}, except
those of \S\ref{sec:pillars-axioms}, which are proved in place or in
\S\ref{app:pillars-axioms}, along
the dependency spines of the results' Isabelle endpoints.  Each step is stated at
the granularity of the mechanized lemma that carries it, with every
hypothesis.  The internal completeness theorem, whose proof is the longest
induction inside \rga in Section~\ref{sec:pillars}
(Theorem~\ref{thm:pl-omega}(1) is a short fuel induction like
Theorem~\ref{thm:in-never}), receives the most space (\S\ref{app:pillars-int}): its motive is given in
full, one value arm and the arm for the universal are exhibited, and the
remaining arms are tabulated.  Facts below the spines are cited by theory
file.  Theory files are given relative to \texttt{hol/}.  Every name
resolves in the development (Section~\ref{sec:mech}, ``Artifact''; for
\S\ref{app:pillars-reflect}, see that subsection's first paragraph).

Table~\ref{tab:pillars} lists the results of Section~\ref{sec:pillars}
by kind, by machine and by well-formedness hypothesis.

\begin{table}[htbp]
\centering
\footnotesize
\begin{tabular}{@{}>{\raggedright\arraybackslash}p{0.38\linewidth}>{\raggedright\arraybackslash}p{0.20\linewidth}>{\raggedright\arraybackslash}p{0.16\linewidth}>{\raggedright\arraybackslash}p{0.18\linewidth}@{}}
\hline
result & kind & machine & well-formedness hypothesis \\
\hline
Thm.~\ref{thm:pl-adeq} (adequacy) & schematic (1); meta (2) & decider & none \\
Thm.~\ref{thm:pl-sound} (soundness) & schematic & decider & terms of $\zeta$ only \\
Thm.~\ref{thm:pl-con} (consistency) & meta & decider & none \\
Thm.~\ref{thm:in-never} (never affirms $\z=\tofn 1$) & internal & decider & none \\
Thm.~\ref{thm:pl-no-uniform} (no single-sentence soundness) & meta (non-derivability) & tree checker, decider & none \\
Thm.~\ref{thm:pl-complete} (completeness) & schematic & tree checker & $t$ \\
Cor.~\ref{thm:pl-reflect} (both judgments) & schematic & decider; tree checker & $t$ \\
Cor.~\ref{thm:pl-bridge} (truth to provability) & rule on theorems & decider; tree checker & $t$ \\
Thm.~\ref{thm:pl-prov-iff} (provable iff grounded-true) & meta & none & $t$ \\
Thm.~\ref{thm:pl-nec} (NEC and CONEC) & rules on theorems & decider & $t$ (CONEC) \\
Thm.~\ref{thm:pl-int} (internal completeness) & internal & decider; $\pi$ & guards $\Wfc$, $\IsCode$ in the antecedent \\
Thm.~\ref{thm:pl-int-ent} (completeness as an entailment) & entailment, rule in every context (1, 3); meta (non-derivability) (2) & decider; $\pi$ & guards $\Wfc$, $\IsCode$ as hypotheses (1) \\
Thm.~\ref{thm:pl-genuine} ($\pi$ is genuine) & meta & $\pi$ & none \\
Thm.~\ref{thm:pl-diag} (diagonal) & meta (1); internal (2) & $\Diag$, decider & none \\
Thm.~\ref{thm:pl-no-value} ($L_{\Tru}$ has no value) & meta & decider & none \\
Thm.~\ref{thm:pl-liar-indep} ($L_{\Tru}$ independent) & meta (non-derivability) & decider & none \\
Thm.~\ref{thm:pl-tb-explosive} (the T-biconditional at $L_{\Tru}$ is explosive) & entailment; meta (non-derivability) & decider & none \\
Thm.~\ref{thm:pl-tb-generic} (the same for every predicate) & entailment (1, 2); internal (the code equation, 2); meta (non-derivability) (2) & $\Diag$ & $\varphi$ (1); $T$, closed (2) \\
Thm.~\ref{thm:pl-liar-neither} (the ordinary Liar is neither true nor false) & internal (the two refutations); meta (non-derivability) & decider & none \\
Thm.~\ref{thm:pl-refuted} (refuted and unestablished) & rule on theorems (1); rule with an entailment as conclusion (2); entailments (3, 4); meta (non-derivability) (3, 4) & decider; tree checker & none \\
Thm.~\ref{thm:pl-tracking} (tracking under $\ConR$) & entailments (1, 3); meta (2) & decider; tree checker & $t$ (1) \\
Thm.~\ref{thm:pl-not-sentence} (a truth that is no sentence) & internal ($\emptyset \vdash \Tru(\bquo{t})$); meta (the rest) & decider & none \\
Thm.~\ref{thm:pl-omega} ($\Omega$ is never affirmed) & internal (1, 2); meta (non-derivability) (3) & decider & none \\
Table~\ref{tab:truth-axioms} (axioms of truth) & by cell & decider & by cell \\
Thm.~\ref{thm:pl-rf-val} (values at open images) & schematic (rules in every context) & decider & templates (HOL class) \\
Thm.~\ref{thm:pl-rf-tin} (open T-In and T-Out) & schematic (rules in every context) & decider & templates (HOL class) \\
Thm.~\ref{thm:pl-rf-app} (value law at open codes) & schematic (rule in every context) & decider & guards $\Wfc$, $\IsCode$, $\mathit{IsQf}$ as hypotheses \\
Thm.~\ref{thm:pl-rf-norm} (the normalizer is correct) & internal (the sentence); rule in every context (the instance) & decider; normalizer & guards as decided antecedents (the sentence); as hypotheses (the instance) \\
Thm.~\ref{thm:pl-rf-sentence} (the reflection sentence) & internal & decider; normalizer & guards as decided antecedents \\
Thm.~\ref{thm:pl-rf-reflect} (reflection) & schematic (rule in every context) & decider; normalizer & polynomial templates (HOL class) \\
Cor.~\ref{cor:pl-rf-binom} (a binomial identity) & schematic (rule in every context) & decider; normalizer & none \\
Thm.~\ref{thm:pl-tin-cert} (T-In at a universal, certificate form) & internal & decider; $\pi$ & none \\
Thm.~\ref{thm:pl-tin-nec} (NEC at universal instances) & rule on theorems, with an internal sentence as conclusion & decider; $\pi$; the builder & $G$, closed \\
\hline
\end{tabular}
\caption{The results of Section~\ref{sec:pillars} and
  Theorem~\ref{thm:in-never}, by kind (\S\ref{sec:inside-codes}), by
  the machine each result runs on, and by the well-formedness hypothesis each
  result takes.}
\label{tab:pillars}
\end{table}

\subsection{Theory files}
\label{app:pillars-files}

Table~\ref{tab:b-files} lists the theory file of each name behind the
results of Section~\ref{sec:pillars}.

\begin{table}[htbp]
\centering
\footnotesize
\begin{tabular}{@{}>{\raggedright\arraybackslash}p{0.27\linewidth}>{\raggedright\arraybackslash}p{0.70\linewidth}@{}}
\hline
theory file & names \\
\hline
\texttt{RGA\_Kit/\allowbreak Term\_Coders} & \textsf{iSubM\_closed},
  \textsf{iChkOpq\_decodes} \\
\texttt{RGA\_Decider/\allowbreak Decider\_Truth} &
  \textsf{iWtrue\_of\_wosr}, \textsf{wosr\_of\_iWtrue},
  \textsf{pointwise\_collapse}, \textsf{rga\_sound\_sch},
  \textsf{iWtrue\_of\_proof}, \textsf{rga\_consistent\_via\_wtrue},
  \textsf{iWtrue\_false\_root\_of\_proof}, \textsf{iWtrue\_no\_false\_root},
  \textsf{iPrv}, \textsf{rga\_complete\_sch}, \textsf{rga\_truth\_reflect\_sch},
  \textsf{iPrv\_of\_iWtrue}, \textsf{iWtrue\_nec},
  \textsf{tschema\_\allowbreak unprovable\_\allowbreak ungrounded},
  \textsf{iWtrue\_conec}, \textsf{iWtrue\_nec\_conec}, \textsf{trep\_rule},
  \textsf{tdel\_rule}, \textsf{true\_cons\_rule},
  \textsf{iWtrue\_neg\_transport}, \textsf{iPrv\_reflect},
  \textsf{tin\_iff\_decided}, \textsf{tschema\_comp\_iff\_grounded},
  \textsf{tout\_iff}, \textsf{tout\_false\_root},
  \textsf{tschema\_grounded\_true},
  \textsf{true\_\allowbreak allE\_\allowbreak rule}, \textsf{true\_\allowbreak notall\_\allowbreak rule}, \textsf{true\_\allowbreak hstep\_\allowbreak rule}, \textsf{isub\_\allowbreak or}, \textsf{isub\_\allowbreak and}, \textsf{isub\_\allowbreak dneg}, \textsf{aff\_\allowbreak true\_\allowbreak int}, \textsf{true\_\allowbreak not\_\allowbreak sentence}, \textsf{neg\_\allowbreak conec\_\allowbreak grounded}, \textsf{tb\_grounded}, \textsf{tb\_\allowbreak fails\_\allowbreak ungrounded}, \textsf{iWtrue\_\allowbreak transport}, \textsf{rga\_\allowbreak wtrueE}, \textsf{iWtrue\_\allowbreak allI\_\allowbreak of\_\allowbreak proof} \\
\texttt{RGA\_Checker/\allowbreak RGA\_Deriv} & \textsf{rga\_prov\_iff\_true}, \textsf{Dctx\_proj}, \textsf{Dctx\_mono} \\
\texttt{RGA\_Checker/\allowbreak Chk\_HasCert} & \textsf{iHasCert} \\
\texttt{RGA\_Checker/\allowbreak Chk\_Completeness} &
  \textsf{rga\_complete\_int}, \textsf{iHasCert\_of\_iPrDecTrueAt},
  \textsf{rga\_complete\_int2}, \textsf{rga\_complete\_int\_zeroN} \\
\texttt{RGA\_Checker/\allowbreak Chk\_Base} & \textsf{provesChk\_genuine} \\
\texttt{RGA\_Checker/\allowbreak Chk\_Faithful} & \textsf{rga\_complete\_sch\_chk}, \textsf{rga\_truth\_reflect\_sch\_chk},
  \textsf{iPrC\_of\_iWtrue}, \textsf{rga\_prov\_reflect\_sch\_chk} \\
\texttt{RGA\_Checker/\allowbreak Chk\_QF\_Node\_Laws} &
  \textsf{iWtrue\_\allowbreak dneg\_\allowbreak intro}, \textsf{iWtrue\_\allowbreak dneg\_\allowbreak elim}, \textsf{iWtrue\_\allowbreak or\_\allowbreak intro1}, \textsf{iWtrue\_\allowbreak or\_\allowbreak intro2}, \textsf{iWtrue\_\allowbreak disjE2\_\allowbreak elim}, \textsf{iWtrue\_\allowbreak disjE3\_\allowbreak elim}, \textsf{iWtrue\_\allowbreak sucEqE\_\allowbreak elim}, \textsf{iWtrue\_\allowbreak sucne0\_\allowbreak intro}, \textsf{iWtrue\_\allowbreak sucneI\_\allowbreak intro}, \textsf{iWtrue\_\allowbreak sucNeE\_\allowbreak elim}, \textsf{iWtrue\_\allowbreak eqsym\_\allowbreak intro}, \textsf{iWtrue\_\allowbreak neSym\_\allowbreak intro},
  \textsf{iWtrue\_or\_elim} \\
\texttt{RGA\_Checker/\allowbreak Chk\_Emission} &
  \textsf{true\_\allowbreak not\_\allowbreak sentence\_\allowbreak internal} \\
\texttt{RGA\_QF/\allowbreak QF\_Arms} & \textsf{iWtrue\_\allowbreak or\_\allowbreak intro3\_\allowbreak qf}, \textsf{iWtrue\_refute\_afree}, \textsf{negnec\_qf}, \textsf{tout\_qf},
  \textsf{tschema\_afree\_grounded}, \textsf{iWtrue\_\allowbreak refute\_\allowbreak of\_\allowbreak neg\_\allowbreak qf}, \textsf{iWtrue\_\allowbreak neg\_\allowbreak refute\_\allowbreak qf}, \textsf{aff\_\allowbreak cons\_\allowbreak qf\_\allowbreak int}, \textsf{tneg\_\allowbreak rule\_\allowbreak grounded} \\
\texttt{RGA\_QF/\allowbreak QF\_Master\_Key} & \textsf{hstep\_\allowbreak iWtrue\_\allowbreak fwd}, \textsf{hstep\_\allowbreak iWtrue\_\allowbreak rev}, \textsf{omega\_\allowbreak never\_\allowbreak affirmed}, \textsf{tout\_\allowbreak omega}, \textsf{neg\_\allowbreak conec\_\allowbreak fails}, \textsf{fs\_\allowbreak neg\_\allowbreak rule\_\allowbreak fails}, \textsf{tneg\_\allowbreak rule\_\allowbreak fails}, \textsf{aff\_\allowbreak or\_\allowbreak iff\_\allowbreak int} \\
\texttt{RGA\_QF/\allowbreak QF\_Mono} & \textsf{imono\_\allowbreak qf\_\allowbreak le} \\
\texttt{RGA/\allowbreak RGA\_Paradox} & \textsf{trmOmega}, \textsf{omega\_\allowbreak no\_\allowbreak value}, \textsf{rga\_liar\_\allowbreak unprovable}, \textsf{rga\_liar\_\allowbreak irrefutable},
  \textsf{indep\_allapp\_no\_value} \\
\texttt{RGA\_Goedel/\allowbreak Goedel\_II} & \textsf{iDiagSub\_closed},
  \textsf{rga\_goedel2\_via\_search}, \textsf{rga\_pr\_reflect}, \textsf{ConRGA\_instances\_refuted},
  \textsf{goedelG\_no\_value}, \textsf{kf\_\allowbreak eq\_\allowbreak int}, \textsf{kf\_\allowbreak ne\_\allowbreak int}, \textsf{aff\_\allowbreak dneg\_\allowbreak iff\_\allowbreak int} \\
\texttt{RGA\_Goedel/\allowbreak Goedel\_II\_PR} & \textsf{rga\_goedel2}, \textsf{ConRGA\_omega\_incomplete}, \textsf{ConRGA\_no\_value},
  \textsf{ConRGA\_not\_true}, \textsf{goedelG\_under\_ConRGA},
  \textsf{rga\_no\_uniform\_}\allowbreak\textsf{soundness},
  \textsf{uniform\_soundness\_}\allowbreak\textsf{unrefutable},
  \textsf{uniform\_soundness\_}\allowbreak\textsf{no\_value},
  \textsf{uniform\_soundness\_}\allowbreak\textsf{instances} \\
\texttt{RGA\_Goedel/\allowbreak Entailment\_Goedel} & \textsf{complete\_ent}, \textsf{iHasCert\_\allowbreak of\_iWtrue},
  \textsf{prc\_true\_\allowbreak ent\_fails}, \textsf{omv\_\allowbreak ConRGA\_true},
  \textsf{omv\_\allowbreak goedelG\_true}, \textsf{con\_ctx\_\allowbreak omega\_sound} \\
\texttt{RGA\_Liar\_Div/\allowbreak Entailment\_Liar} & \textsf{complete\_\allowbreak imp\_fails}, \textsf{liar\_code\_\allowbreak not\_refuted\_true},
  \textsf{liar\_code\_\allowbreak no\_cert}, \textsf{tbicond\_\allowbreak liar\_\allowbreak explosive\_gen},
  \textsf{tbicond\_\allowbreak liar\_\allowbreak explosive\_pred}, \textsf{pred\_liar\_eq},
  \textsf{tbicond\_\allowbreak liar\_\allowbreak explosive\_diag},
  \textsf{tbicond\_\allowbreak liar\_\allowbreak unprovable\_diag} \\
\texttt{RGA\_Liar\_Div/\allowbreak Entailment\_Truth} & \textsf{ent\_of\_\allowbreak ctx\_rule}, \textsf{kf\_dneg\_ent}, \textsf{kf\_or\_ent},
  \textsf{fs\_neg\_\allowbreak ent\_fails}, \textsf{all\_omega\_\allowbreak ent\_fails},
  \textsf{tneg\_ent\_\allowbreak fails\_omega}, \textsf{negconec\_\allowbreak ent\_fails\_omega},
  \textsf{comp\_ent\_\allowbreak fails}, \textsf{true\_sent\_\allowbreak ent\_fails},
  \textsf{omv\_iWtrue\_\allowbreak not\_true}, \textsf{tin\_ent\_\allowbreak fails\_liar},
  \textsf{tin\_ent\_\allowbreak fails\_goedelG}, \textsf{tin\_ent\_\allowbreak fails\_con},
  \textsf{tout\_ent\_\allowbreak fails\_liar}, \textsf{negnec\_ent\_\allowbreak fails\_liar},
  \textsf{omv\_\allowbreak goedelW\_true}, \textsf{tbicond\_\allowbreak liar\_explosive},
  \textsf{tbicond\_\allowbreak liar\_both}, \textsf{tbicond\_\allowbreak liar\_unprovable},
  \textsf{establishment\_\allowbreak transparency},
  \textsf{refuted\_\allowbreak unestablished\_qf},
  \textsf{omega\_\allowbreak unestablished} \\
\texttt{RGA\_Liar\_Div/\allowbreak Truth\_Limits} & \textsf{ordinary\_\allowbreak liar\_neither}, \textsf{truth\_liar\_\allowbreak not\_untrue\_\allowbreak unprovable},
  \textsf{liar\_never\_\allowbreak affirmed}, \textsf{neg\_liar\_\allowbreak never\_affirmed},
  \textsf{refuted\_\allowbreak unestablished\_\allowbreak trichotomy},
  \textsf{universal\_\allowbreak unestablished\_\allowbreak needs\_con},
  \textsf{universal\_\allowbreak unestablished\_\allowbreak unprovable},
  \textsf{con\_refuted\_\allowbreak unestablished}, \textsf{con\_refuted\_\allowbreak unprovable},
  \textsf{refuted\_\allowbreak universal\_\allowbreak unestablished\_\allowbreak iff\_con},
  \textsf{refuted\_\allowbreak unestablished\_\allowbreak fails},
  \textsf{refuted\_chk\_\allowbreak inconsistent}, \textsf{refuted\_\allowbreak affirmed\_\allowbreak inconsistent},
  \textsf{ir\_A\_cert\_\allowbreak affirm}, \textsf{notcon\_tru\_\allowbreak univ\_ent},
  \textsf{iPrc\_of\_\allowbreak iHasCert\_code}, \textsf{tru\_prc\_ent},
  \textsf{omega\_\allowbreak unestablished\_sf},
  \textsf{con\_tschema\_\allowbreak decided}, \textsf{con\_tracking\_\allowbreak sound},
  \textsf{con\_prc\_\allowbreak schema\_decided}, \textsf{con\_prc\_\allowbreak bool\_decided},
  \textsf{con\_true\_\allowbreak bool\_decided}, \textsf{con\_tschema\_\allowbreak refuted},
  \textsf{con\_goedelG\_\allowbreak derivable\_\allowbreak unestablished},
  \textsf{con\_tschema\_\allowbreak goedelG\_refuted},
  \textsf{con\_con\_\allowbreak unestablished},
  \textsf{unestablished\_\allowbreak of\_unprovable},
  \textsf{con\_goedelG\_\allowbreak and\_\allowbreak unestablished},
  \textsf{universal\_\allowbreak unestablished\_\allowbreak needs\_con\_ctx},
  \textsf{universal\_\allowbreak unestablished\_\allowbreak con\_ent},
  \textsf{refuted\_\allowbreak universal\_\allowbreak unestablished\_\allowbreak eqv\_con},
  \textsf{refuted\_\allowbreak universal\_\allowbreak unestablished\_\allowbreak iff\_con\_ctx},
  \textsf{con\_eqv\_\allowbreak not\_true\_con},
  \textsf{no\_moorean\_\allowbreak theorem}, \textsf{no\_moorean\_\allowbreak pair},
  and the eleven theorems of the cells at $U$ in
  Table~\ref{tab:c-truth-axioms} \\
\texttt{RGA\_Liar\_Div/\allowbreak Liar\_Record} & \textsf{goedelW\_diag},
  \textsf{goedelWS\_eq}, \textsf{goedelW\_no\_value},
  \textsf{goedelW\_unprovable}, \textsf{goedelW\_unrefutable},
  \textsf{goedelW\_undecided}, \textsf{goedelW\_\allowbreak instances\_\allowbreak refuted},
  \textsf{goedelW\_\allowbreak omega\_\allowbreak incomplete}, \textsf{comp\_rule\_fails},
  \textsf{iWtrue\_\allowbreak undecided}, \textsf{guards\_\allowbreak at\_\allowbreak liar}, \textsf{kf\_\allowbreak dneg\_\allowbreak not\_\allowbreak internal}, \textsf{kf\_\allowbreak dneg\_\allowbreak iff\_\allowbreak not\_\allowbreak internal}, \textsf{kf\_\allowbreak or\_\allowbreak not\_\allowbreak internal}, \textsf{kf\_\allowbreak or\_\allowbreak iff\_\allowbreak not\_\allowbreak internal}, \textsf{kf\_\allowbreak negor\_\allowbreak not\_\allowbreak internal}, \textsf{kf\_\allowbreak negor\_\allowbreak iff\_\allowbreak not\_\allowbreak internal}, \textsf{fs\_\allowbreak neg\_\allowbreak not\_\allowbreak internal}, \textsf{fs\_\allowbreak neg\_\allowbreak iff\_\allowbreak not\_\allowbreak internal}, \textsf{all\_\allowbreak omega\_\allowbreak rule\_\allowbreak fails}, \textsf{kf\_\allowbreak all\_\allowbreak not\_\allowbreak internal}, \textsf{kf\_\allowbreak all\_\allowbreak iff\_\allowbreak not\_\allowbreak internal}, \textsf{trep\_\allowbreak not\_\allowbreak internal}, \textsf{comp\_\allowbreak not\_\allowbreak internal}, \textsf{tout\_\allowbreak fails\_\allowbreak liar} \\
\hline
\end{tabular}
\caption{The theory files of the names cited in Section~\ref{sec:pillars}
  and of the truth-predicate results of \S\ref{sec:breaks-liars}; the
  names of \S\ref{sec:pillars-reflect} are in Table~\ref{tab:c-rf-files}.}
\label{tab:b-files}
\end{table}

\subsection{The companion paper's soundness invariant}
\label{app:pillars-sem}

The soundness pillar reuses, rather than repeats, the induction over
proofs in the companion paper's soundness theorem~\cite{rga-paper}.  That
theorem is stated with three semantic notions, all in
\texttt{RGA/\allowbreak RGA\_Sound.thy}.  A \emph{substitution} $\sigma$ maps every
variable to a term, and it is \emph{well-formed} when every image is
well-formed (\textsf{wsubf}, in \texttt{RGA/\allowbreak RGA\_GCalc.thy}).  A substitution may leave variables free.

\begin{defi}[\textsf{wstrue}, \textsf{wshsat}, \textsf{owsjt}]
\label{def:b-owsjt}
A term $t$ is \emph{true under} $\sigma$ when the instance $\sigma(t)$
is grounded-true, that is, $\sigma(t) \wbev{s} \ctrue$ for some $s$.  The
substitution $\sigma$ \emph{satisfies} a hypothesis set $H$ when every
member of $H$ is true under $\sigma$.  The \emph{open semantic judgment}
$H \models c$ holds when $c$ is true under every well-formed
substitution that satisfies $H$, closing or not.
\end{defi}

\noindent
A proof $P$ \emph{contains} the judgment $H \vdash c$ when the pair
$(H, c)$ is one of the entries of $P$ (\textsf{prf\_contains}).  This relation is
weaker than $P \proves H \vdash c$, which also requires $P$ to be valid
and $H \vdash c$ to be the conclusion of $P$ (the first entry of $P$).
The soundness proof also uses two facts about the satisfaction of a
context, projection and monotonicity (\textsf{wshsatD},
\textsf{wshsat\_mono}, restated for this paper as \textsf{Dctx\_proj} and
\textsf{Dctx\_mono}).  Since \rga's deduction theorem holds only at
decided hypotheses, an entailment cannot in general be traded for an
implication, and a context is part of every derivability judgment.

\begin{thm}[\textsf{prf\_owsjt}]
\label{thm:b-prf-owsjt}
If $P$ is valid, $P$ contains $H \vdash c$, and every member of $H$ is
well-formed, then $H \models c$.
\end{thm}

\noindent
The theorem is the companion paper's, proved by induction on $P$.  In
the development it is the instance, at \rga's proofs, of the same
theorem for an evaluation parametric in the system whose proofs serve
as certificates (\textsf{gsem\_sound.prf\_gowsjt},
\texttt{RGA/\allowbreak RGA\_GSem.thy}).  The
companion paper derives from this theorem soundness at the empty context
(\textsf{rga\_provable\_wosr}) and consistency
(\textsf{rga\_consistent}).  Open completeness
(\textsf{rga\_complete\_open}, \texttt{RGA/\allowbreak RGA\_Complete.thy}) takes a
well-formed $c$ and an evaluation of $c$ to $\ctrue$ in the companion
paper's evaluation relation $\Downarrow$ (Isabelle \textsf{osr},
\S\ref{app:machines-cmp}); the lemma \textsf{wosr\_osr} converts a witnessed evaluation
$c \wbev{s} \ctrue$ into an evaluation of $c$ to $\ctrue$ in $\Downarrow$.

\subsection{The soundness pillar}
\label{app:pillars-sound}

The instances of a conclusion come from a finite list $\zeta$ of pairs
of a variable and a term.  The substitution $\hat\zeta$ of the list sends each
variable listed in $\zeta$ to the term of the variable's first entry in
$\zeta$, and every other variable to itself (\textsf{subm\_to\_subf}).  The code $\bquo{\zeta}$ is the code
numeral of the list.

\begin{thm}[\textsf{iSubM\_closed}]
\label{thm:b-subst}
For every well-formed term $t$ and every list $\zeta$ whose terms are all
well-formed,
$\emptyset \vdash \Subst(\bquo{\zeta}, \bquo{t}) = \bquo{\sof{\zeta}{t}}$.
\end{thm}

\noindent
\emph{Proof sketch.}  $\Subst$ is the compiled form of an explicit index
tree whose denotation computes simultaneous substitution on codes.  The
coders inside that tree read at most two argument slots of an operator
node, so the denotation is exact at every code whose operator nodes have
at most two arguments.  The lemma \textsf{iSubM\_closed\_ar} states the
equation under that weaker hypothesis on the terms of $\zeta$, and every
well-formed term satisfies that hypothesis (\textsf{wft\_ar2}).  The evaluation of the
closed term turns into a derivation by Theorem~\ref{thm:in-cmp-num}.
\qed

\begin{thm}[\textsf{pointwise\_collapse}]
\label{thm:b-collapse}
Let $P$ be a valid proof that contains $H \vdash c$, where every member
of $H$ is well-formed and $c$ is well-formed.  Let $\zeta$ be a list whose
terms are all well-formed, such that $\hat\zeta$ satisfies $H$.  Then for
every context $\Delta$,
$\Delta \vdash \Tru(\Subst(\bquo{\zeta}, \bquo{c}))$.
\end{thm}

\noindent
The theorem asks only that $P$ contain the judgment, not that $P$ prove
the judgment, and the theorem asks for the well-formedness of $c$.
Theorem~\ref{thm:pl-sound} needs neither hypothesis in that form,
because a proof of $c$ contains its own conclusion and makes $c$
well-formed.

\emph{Proof.}  The proof is one composition, with no induction of its
own.
\begin{enumerate}
\item The substitution $\hat\zeta$ is well-formed, because the terms of
  $\zeta$ are.  Theorem~\ref{thm:b-prf-owsjt} and the hypothesis that
  $\hat\zeta$ satisfies $H$ make $c$ true under $\hat\zeta$
  (\textsf{owsjtD}): $\sof{\zeta}{c} \wbev{s} \ctrue$ for some $s$.
\item Theorem~\ref{thm:pl-adeq}(1) turns this evaluation into
  $\emptyset \vdash \Tru(\bquo{\sof{\zeta}{c}})$
  (\textsf{iWtrue\_of\_wosr}).
\item Theorem~\ref{thm:b-subst} gives
  $\emptyset \vdash \Subst(\bquo{\zeta}, \bquo{c}) = \bquo{\sof{\zeta}{c}}$.
  The transport rule (Appendix~\ref{app:machines},
  \textsf{iWtrue\_transport}) moves $\Tru$ along this equation, and
  weakening (\textsf{ir\_weaken0}) moves the judgment to $\Delta$.
  \qed
\end{enumerate}

\noindent
\emph{Proof of Theorem~\ref{thm:pl-sound}} (\textsf{rga\_sound\_sch}).
A proof $P$ with $P \proves \emptyset \vdash c$ is valid
(\textsf{proves\_valid}), contains $\emptyset \vdash c$
(\textsf{proves\_contains}), and has a well-formed conclusion
(\textsf{proves\_wff}).  The empty hypothesis set is well-formed and every
substitution satisfies it, so Theorem~\ref{thm:b-collapse} applies.
\qed

\medskip
The same composition gives the pillar for proofs from hypotheses, under
the truth of the hypotheses.

\begin{thm}[\textsf{rga\_sound\_sch\_ctx}]
\label{thm:b-sound-ctx}
Let $P \proves H \vdash c$, where every member of $H$ is well-formed,
and let $\zeta$ be a list whose terms are all well-formed, such that
$\hat\zeta$ satisfies $H$.  Then for every context $\Delta$,
$\Delta \vdash \Tru(\Subst(\bquo{\zeta}, \bquo{c}))$.
\end{thm}

\begin{thm}[\textsf{iWtrue\_of\_proof}]
\label{thm:b-of-proof}
If $P \proves \emptyset \vdash c$, then
$\emptyset \vdash \Tru(\bquo{c})$.
\end{thm}

\noindent
\emph{Proof.}  Apply Theorem~\ref{thm:pl-sound} with $\zeta$ empty.
The term $c$ is well-formed because it is proved (\textsf{proves\_wff}),
so Theorem~\ref{thm:b-subst} at the empty list gives
$\emptyset \vdash \Subst(\bquo{[\,]}, \bquo{c}) = \bquo{c}$
(\textsf{subm\_trm\_empty}), and transport rewrites the argument.
\qed

\medskip
The truth predicate also admits $\forall$-introduction at the level of
codes.  The code of $\cA\cdot u$ is computed from $\bquo{u}$ by the code
constructor $\mathit{ap}$ of [\citen{goedel-paper}, Appendix ``G\"odel II in full''] (Isabelle
\textsf{mkAppS}), applied to the compiled-pair spelling of the code of
$\cA$ (Isabelle \textsf{mkAS}).  In this appendix $\bquo{\cA}$ in
$\mathit{ap}(\bquo{\cA}, \cdot)$ denotes that spelling, which \rga
proves equal to the numeral code of $\cA$ (\textsf{mkAS\_val}).

\begin{thm}[\textsf{iWtrue\_allI\_of\_proof}]
\label{thm:b-allI}
If the number $p$ decodes to a proof of
$\{\vz\jnat\} \vdash \lft u \cdot \vz$, then for every context $\Delta$,
$\Delta \vdash \Tru(\mathit{ap}(\bquo{\cA}, \bquo{u}))$.
\end{thm}

\noindent
\emph{Proof.}  The proof stays in HOL until its last step.  The decoded
proof is a certificate for the universal, and the companion paper's
soundness theorem rebuilds from the certificate the family of true numeral instances
(\textsf{prDecV\_allT\_family}).  The truth clause of the semantics then
evaluates $\cA\cdot u$ to $\ctrue$ (\textsf{wosr\_allT\_i}), at a fuel
above the certificate's code.  Theorem~\ref{thm:pl-adeq}(1) crosses to
$\emptyset \vdash \Tru(\bquo{\cA\cdot u})$, and transport along
$\emptyset \vdash \mathit{ap}(\bquo{\cA}, \bquo{u}) = \bquo{\cA\cdot u}$
(\textsf{mkAS\_app\_val}) and weakening finish.
\qed

\subsection{The consistency route}
\label{app:pillars-con}

Theorem~\ref{thm:pl-con} (\textsf{rga\_consistent\_via\_wtrue}) composes
the following two facts, both in
\texttt{RGA\_Decider/\allowbreak Decider\_Truth.thy}.

\begin{thm}[\textsf{iWtrue\_false\_root\_of\_proof}, \textsf{iWtrue\_no\_false\_root}]
\label{thm:b-false-root}
\leavevmode
\begin{enumerate}
\item If $P \proves \emptyset \vdash \z = \tofn 1$, then
  $\emptyset \vdash \Tru(\bquo{\z = \tofn 1})$.
\item If $\emptyset \vdash \Tru(\bquo{\z = \tofn 1})$, then
  $\mathit{False}$.
\end{enumerate}
\end{thm}

\noindent
\emph{Proof.}  Part~1 is Theorem~\ref{thm:b-of-proof} at $c$ the false
equation.  For part~2, write $u$ for the body
$\lambda\neg\Aff(\vz, \lft\bquo{\z = \tofn 1})$, so that
$\Tru(\bquo{\z = \tofn 1})$ is $\neg(\cA\cdot u)$.
\begin{enumerate}
\item By the companion paper's soundness theorem
  (\textsf{rga\_provable\_wosr}), $\neg(\cA\cdot u)$ is grounded-true.  So
  $\cA\cdot u$ evaluates to $\cfalse$ (\textsf{wosr\_not\_inv}).
\item The falsity clause of the universal (\textsf{wosr\_Aapp\_inv}; the
  truth clause would give $\ctrue$) supplies a numeral $\tofn i$ and a proof
  of $\emptyset \vdash \neg(u\cdot\tofn i)$, so
  $\emptyset \vdash \neg(u\cdot\tofn i)$.
\item The $\beta$-step from $u\cdot\tofn i$ to
  $\neg\Aff(\tofn i, \bquo{\z = \tofn 1})$ (\textsf{hs\_beta},
  \textsf{bsub\_iPrDecTrueAt\_fuel}) transports the refutation along the
  head step (\textsf{ir\_conv\_negE}), and double-negation elimination
  gives $\emptyset \vdash \Aff(\tofn i, \bquo{\z = \tofn 1})$.
\item At the numeral fuel $\tofn i$, Theorem~\ref{thm:in-dec-closed} gives
  $\emptyset \vdash \Dec(\tofn i, \bquo{\z = \tofn 1}) =
  \tofn{\decP_i(\quo{\z = \tofn 1})}$, and the packed verdict is never 2
  (\textsf{prDecVEnc\_sEq01\_ne2}).  So \rga derives
  $\Aff(\tofn i, \bquo{\z = \tofn 1}) = \z$.
\item The affirmation is then both nonzero and zero, so \rga derives
  $\z = \tofn 1$.  By the companion paper's value characterization
  (\textsf{rga\_val\_iff}), the term $\z$ would evaluate to 1, which
  \textsf{osr\_num\_inv} excludes.  \qed
\end{enumerate}

\noindent
The packed verdict of step~4 changes with the fuel ($0$ at fuel 2 and
$1$ from fuel 3 on), so the refutation does not go through a fixed point
of the decider.  The refutation goes through HOL soundness and back.  The internal
companion of part~2 is Theorem~\ref{thm:in-never}, which derives the
universal $\forall s.\ \neg\Aff(s, \bquo{\z = \tofn 1})$ that
$\Tru(\bquo{\z = \tofn 1})$ negates.

\paragraph{Proof of Theorem~\ref{thm:pl-no-uniform}.}
The theorem \textsf{rga\_no\_uniform\_}\allowbreak\textsf{soundness}
(\texttt{RGA\_Goedel/\allowbreak Goedel\_II\_PR.thy}) states
$\emptyset \nvdash \forall p.\ \Chk(p, \bquo{\z = \tofn 1}) \limp
\Tru(\bquo{\z = \tofn 1})$.  Suppose \rga derived that sentence, call it
$U$, and write $c$ for $\Chk(\vz, \bquo{\z = \tofn 1})$ and
$\Gamma_0$ for the context $\{\vz\jnat\}$.  The derivation runs as
follows.
\begin{enumerate}
\item $\Gamma_0 \vdash \vz\jnat$ by the hypothesis rule.
\item $\Gamma_0 \vdash c \lor \neg c$, since $c$ is decided at the
  abstract $\vz$ (Theorem~\ref{thm:in-chk-dec}, \textsf{hq\_iChk}, read
  through \textsf{ir\_boolE}).
\item In the accepting branch, with context $c, \Gamma_0$: weakening and
  universal elimination of $U$ at $\vz$, followed by the $\beta$-step,
  give $c \limp \Tru(\bquo{\z = \tofn 1})$
  (\textsf{gd2\_allE\_beta}).
\item Modus ponens with the branch hypothesis $c$ gives
  $\Tru(\bquo{\z = \tofn 1})$ (\textsf{ir\_mp}).
\item Unfolding $\Tru$, this is $\neg(\cA\cdot(\lambda\neg\Aff(\vz,
  \lft\bquo{\z=\tofn 1})))$.  Theorem~\ref{thm:in-never} derives the
  universal itself, so negation elimination (\textsf{ir\_negE}) gives any
  term, in particular $\neg c$.
\item In the rejecting branch, with context $\neg c, \Gamma_0$, the goal
  $\neg c$ is the branch hypothesis.
\item Disjunction elimination (\textsf{ir\_disjE1}) on step~2 gives
  $\Gamma_0 \vdash \neg c$.
\item Universal introduction at the eigenvariable $\vz$
  (\textsf{ir\_allI\_open}) gives $\emptyset \vdash \ConR$, which
  {}[\citen{goedel-paper}, theorem~\textsf{rga\_goedel2}] (\textsf{rga\_goedel2}) excludes.  \qed
\end{enumerate}

\noindent
The derivation eliminates the assumed implication in step~4 and
introduces none.  The derivation has the shape of the decided-matrix split of
\S\ref{sec:con-shape}.

Each numeral instance of $U$'s matrix is a theorem: for every $n$,
\rga derives $(\lambda(\Chk(\vz, \bquo{\z = \tofn 1}) \limp
\Tru(\bquo{\z = \tofn 1})))\cdot\tofn n$
(\textsf{uniform\_soundness\_}\allowbreak\textsf{instances}, a schema
over $n$).  At each $n$ \rga refutes $\Chk(\tofn n, \bquo{\z = \tofn
1})$ (\textsf{ConRGA\_instances\_refuted}), and the implication, which
unfolds to a disjunction, follows by introducing its first disjunct.
The schema rests on soundness, since \textsf{ConRGA\_instances\_refuted}
does.

\subsection{The completeness pillar}
\label{app:pillars-complete}

\emph{Proof of Theorem~\ref{thm:pl-complete}} (\textsf{rga\_complete\_sch\_chk},
\texttt{RGA\_Checker/\allowbreak Chk\_Faithful.thy}).
Let $t$ be well-formed with $t \wbev{s} \ctrue$.
\begin{enumerate}
\item The evaluation converts to the companion paper's evaluation
  relation $\Downarrow$ (\textsf{wosr\_osr}), and open completeness gives
  $\emptyset \vdash t$ (\textsf{rga\_complete\_open}).
\item A derivation has a proof object: some $P$ satisfies
  $P \proves \emptyset \vdash t$ (\textsf{entails\_provesE}).
\item The checker accepts $P$'s code:
  $\emptyset \vdash \Chk(\bquo{P}, \bquo{t})$ (\textsf{iChk\_faithful}).
\item A $\beta$ step and negated universal introduction
  (\textsf{ir\_not\_allI}) at the numeral $\bquo{P}$ conclude
  $\emptyset \vdash \PrC(\bquo{t})$.  \qed
\end{enumerate}

\noindent
Corollary~\ref{thm:pl-reflect} (\textsf{rga\_truth\_reflect\_sch\_chk})
pairs Theorem~\ref{thm:pl-complete} with Theorem~\ref{thm:pl-adeq}(1).  Corollary~\ref{thm:pl-bridge}
(\textsf{iPrC\_of\_iWtrue}) composes Theorem~\ref{thm:pl-adeq}(2)
(\textsf{wosr\_of\_iWtrue}) with Theorem~\ref{thm:pl-complete}.  The
proof of Theorem~\ref{thm:pl-adeq}(2) is the one sketched in
\S\ref{sec:pillars-adeq}.  In Isabelle that proof runs through
\textsf{rga\_provable\_wosr}, \textsf{wosr\_not\_inv},
\textsf{wosr\_Aapp\_inv} (the falsity clause supplies a numeral fuel
$\tofn i$), \textsf{hs\_beta} with \textsf{whstep\_rev\_wosr} (the
affirmation at $\tofn i$ is grounded-true), \textsf{rga\_open\_iff},
\textsf{iPrDecTrueAt\_true\_elim} with \textsf{iPrDecV\_closed} (the
packed verdict at fuel $i$ is 2), and \textsf{prDecVEnc\_sound}
(decider soundness).  None of these steps needs $t$ to be well-formed.

Faithfulness, used in step~3, itself rests on \rga's soundness and
open completeness in HOL (Theorem~\ref{thm:in-chk}).  That dependence
concerns the adequacy of the checker, a meta-level fact; where an
internal derivation needs the acceptance of a fixed proof, the
development obtains that acceptance instead by closed evaluation of the checker,
without soundness (\textsf{iChk\_proves\_sf}, used by
\textsf{prc\_of\_derivable}).

Theorem~\ref{thm:pl-prov-iff} (\textsf{rga\_prov\_iff\_true},
\texttt{RGA\_Checker/\allowbreak RGA\_Deriv.thy}) has two directions.  The forward
direction is soundness at the empty context; the backward direction is
steps~1 and~2 above.  Composed with Corollary~\ref{thm:pl-reflect},
Theorem~\ref{thm:pl-prov-iff} gives a third form of the pillar, stated from a proof rather than from
an evaluation.

\begin{thm}[\textsf{rga\_prov\_reflect\_sch\_chk}]
\label{thm:b-prov-reflect}
For every term $t$: if some proof $P$ satisfies
$P \proves \emptyset \vdash t$, then $\emptyset \vdash \Tru(\bquo{t})$
and $\emptyset \vdash \PrC(\bquo{t})$.
\end{thm}

\noindent
The theorem needs no well-formedness hypothesis, since a proof of $t$
makes $t$ well-formed (\textsf{entails\_wff}).

\paragraph{The opaque predicate.}
The development also has a second provability predicate, the opaque
twin of $\PrC$.  The opaque predicate wraps the opaque checker $\ChkE$, compiled from an
index chosen by Hilbert choice for the companion paper's
proof-certificate test (\texttt{RGA\_Kit/\allowbreak Term\_Coders.thy}:
\textsf{chk2}, \textsf{chkIdx}, \textsf{iChkOpq};
Appendix~\ref{app:machines-opaque}).  \rga derives $\ChkE(\tofn p,
\bquo{c})$ exactly when the HOL decoding of $p$ is a proof of
$\emptyset \vdash c$ (\textsf{iChkOpq\_decodes}).

\begin{defi}[\textsf{iPrv}]
\label{def:pl-prv}
$\Prv(c) \ldef \neg(\cA\cdot(\lambda\neg\ChkE(\vz, \lft c)))$, that
is, $\Prv(c) \ldef \exists p.\ \ChkE(p, c)$.
\end{defi}

\noindent
$\Prv$ serves as the opaque control of the build-time check of
\S\ref{sec:inside-audit}.  The introduction rule of $\Prv$ mirrors that
of $\Tru$.

\begin{thm}[\textsf{iPrvI}]
\label{thm:b-prvI}
If $\Gamma \vdash c\jnat$, $\Gamma \vdash a\jnat$ and
$\Gamma \vdash \ChkE(a, c)$, then $\Gamma \vdash \Prv(c)$.
\end{thm}

\noindent
The forms of the pillar over $\Prv$ (\textsf{rga\_complete\_sch},
\textsf{rga\_truth\_reflect\_sch} and \textsf{iPrv\_of\_iWtrue} in
\texttt{RGA\_Decider/\allowbreak Decider\_Truth.thy}, and
\textsf{rga\_prov\_reflect\_sch} in \texttt{RGA\_Checker/\allowbreak
RGA\_Deriv.thy}, which carries a well-formedness hypothesis) have the
same proofs, with two changes.  Step~3 becomes $\emptyset \vdash
\ChkE(\bquo{P}, \bquo{t})$ (\textsf{iChkOpq\_true}).  The proof
evaluates the closed checker application in HOL
(\textsf{cmp\_correct\_num}).  The value of that application is the value
of the denotation of the chosen index, which is the proof-certificate test itself
(\textsf{pr2fi\_priof} with \textsf{n2pr\_chk2}), and that test accepts
the code of a proof (\textsf{prfCertCode\_encode}).  Open completeness
(\textsf{rga\_open\_iff}) turns the evaluation into a derivation.
Step~4 becomes Theorem~\ref{thm:b-prvI} at the numerals $\bquo{t}$ and
$\bquo{P}$.

\paragraph{Converses.}
For the transparent $\PrC$ the converse is a theorem: if $\emptyset
\vdash \PrC(\bquo{\psi})$, then $\emptyset \vdash \psi$
(\textsf{rga\_pr\_reflect}, [\citen{goedel-paper}, theorem~\textsf{D1\_rule\_pr}](4),
\texttt{RGA\_Goedel/\allowbreak Goedel\_II.thy}).  The witness property
(\textsf{rga\_not\_all\_witness}) gives a numeral $\tofn n$ with
$\emptyset \vdash \neg((\lambda\neg\Chk(\vz, \bquo{\psi}))\cdot\tofn n)$.
A $\beta$ step (\textsf{ir\_conv\_negE}) and double-negation elimination
(\textsf{ir\_dnegE}) give $\emptyset \vdash \Chk(\tofn n, \bquo{\psi})$,
and genuineness (\textsf{iChk\_genuine\_proves},
Theorem~\ref{thm:in-chk}(2)) turns that numeral into a proof of $\psi$.

$\Prv$ has the same converse: if $\emptyset \vdash \Prv(\bquo{\psi})$,
then $\emptyset \vdash \psi$ (\textsf{iPrv\_reflect},
\texttt{RGA\_Decider/\allowbreak Decider\_Truth.thy}).  The same steps,
with \textsf{bsub\_iChkOpq\_proof} for the $\beta$ step, give
$\emptyset \vdash \ChkE(\tofn n, \bquo{\psi})$.  In place of
genuineness, Theorem~\ref{thm:a-chkopq-decode}
(\textsf{iChkOpq\_decodes}) says that the decoding of $n$ is a proof of
$\psi$, and \textsf{entails\_provesI} turns that proof into a
derivation.

\subsection{Internal completeness}
\label{app:pillars-int}

This subsection proves Theorem~\ref{thm:pl-int}
(\textsf{rga\_complete\_int2}, through its packed-pair form
\textsf{rga\_complete\_int}).  Unless noted, names are in
\texttt{RGA\_Checker/\allowbreak Chk\_Completeness.thy}; the value half of the
motive and the value arms are in \texttt{RGA\_Checker/\allowbreak Chk\_Transfers.thy}, and
the certificate steps in \texttt{RGA\_Checker/\allowbreak Chk\_HasCert.thy}.
The read-down law \textsf{readdownV\_oprNot} is in
\texttt{Chk\_Readdown}, \textsf{iWft\_neg1\_kid} in
\texttt{Chk\_Wf\_Invariant}, \textsf{iOprHd1\_neg\_eta} in
\texttt{Chk\_QF\_Faces}, the step laws
\textsf{ir\_stepV\_\allowbreak app\_cert\_hitT},
\textsf{ir\_stepV\_\allowbreak app\_cert\_hitF} and
\textsf{ir\_stepV\_\allowbreak app\_cert\_miss} in \texttt{Chk\_Coder\_Bounds}, and
\textsf{ir\_certTrue\_elim} in \texttt{Chk\_Formers} (all in
\texttt{RGA\_Checker}).

\paragraph{Notation.}
The code constructors are those of [\citen{goedel-paper}, Appendix ``G\"odel II in full''] on
G\"odel's second theorem in full:
$\mathit{eq}(a,b)$, $\mathit{neg}(a)$ and $\mathit{ap}(a,b)$ build the
codes of an equation, a negation and an application from the codes of
the parts, and $\mathit{jdg}(c)$ is the code of the judgment
$\emptyset \vdash t$ for the term $t$ with code $c$ (\textsf{emptyJdg}).
So $\Cert_\emptyset(c)$ is $\Cert(\mathit{jdg}(c))$.  The value code of a
verdict $v$ is $\Num\cdot v$, the code of the numeral $\tofn v$.  The
compiled decider $\Dec$ packs a value verdict $\Some{r}$ as $r+1$, so
$\Dec(\ell, c) = \suc v$ says that the value decider at level $\ell$ gives
the code $c$ the value $v$.  The compiled head decider
$\Dec^{\mathrm{hd}}(\ell, c)$ (\textsf{iPrDecH}) reads the head half of the
same table cell.

\paragraph{The certificate layer.}
$\Cert(h)$ is an existential over tree codes: some $p$ with
$\cmp(\pi)\cdot h\cdot p \ne \z$.  Every step of the proof below has the
same shape (\textsf{iHasCert\_step1}, \textsf{iHasCert\_step2}).  Each step
eliminates the existential of each premise, which opens a fresh
variable for the premise's tree.  The step appends the new proof nodes to the
coded tree, and proves with the checker's laws at abstract arguments that
$\pi$ accepts the extended tree.  The step then introduces the existential
again.  The checker laws that append one node for one rule are the
\emph{emission cells} of \texttt{RGA\_Checker/\allowbreak Chk\_Emission.thy} and
\texttt{Chk\_Emission\_Cells.thy}.

\begin{defi}[\textsf{decValBody}, \textsf{decHdCore}, \textsf{decHdBody}, \textsf{decMot}]
\label{def:b-motive}
The motive at level $\ell$ is the internal universal
$\mathit{M}(\ell) \ldef \forall q.\ V(\ell, \mathit{fst}(q)) \wedge
\mathit{Hd}(\ell, \mathit{snd}(q))$ of two clauses.
\begin{itemize}
\item The \emph{value clause} $V(\ell, x)$, for $x$ the pair of a code
  $c$ and a value $v$, is
  \[
    \Dec(\ell, c) = \suc v \;\limp\; \Wfc(c) \ne \z \;\limp\;
    \IsCode(c) \ne \z \;\limp\;
    \Cert_\emptyset(\mathit{eq}(c, \Num\cdot v)) .
  \]
\item The \emph{head clause} $\mathit{Hd}(\ell, y)$ reads six components
  from $y$: a code $t$, a head verdict $h$, a spine list $L$ with a
  length $n$, a value code $w$, and a certificate witness $P$.  Write
  $\mathit{sp}(a)$ for the code of $a$ applied to the arguments of the
  spine $L, n$ (\textsf{iAppSpineCI}), and $\mathit{canon}(h)$ for the
  code of the canonical term with head verdict $h$
  (\textsf{iCanonOf}).  The clause first assumes four guards: the spine's
  codes are well-formed and genuine, and $\Wfc(w) \ne \z$ and
  $\IsCode(w) \ne \z$.  It then says
  \[
  \begin{aligned}
    &\Dec^{\mathrm{hd}}(\ell, t) = \suc h \limp \Wfc(t) \ne \z \limp
      \IsCode(t) \ne \z \\
    &\quad\limp \cmp(\pi)\cdot\mathit{jdg}(\mathit{eq}(\mathit{sp}(\mathit{canon}(h)), w))
      \cdot P \ne \z
    \limp \Cert_\emptyset(\mathit{eq}(\mathit{sp}(t), w)) .
  \end{aligned}
  \]
\end{itemize}
\end{defi}

\noindent
Every antecedent is a decided test: an equation between numerals, a
disequation with zero, or an application of a compiled checker.  So each
implication is introduced by \textsf{ir\_impI\_bool}.  The guards
$\Wfc$ and $\IsCode$ are there because the decider dispatches on shape
tags alone, so a verdict carries neither fact.  The head clause names
its certificate witness $P$ and hypothesises the compiled test, rather
than hypothesising $\Cert$ of the canonical equation: \rga discharges
a hypothesis into an implication only when the hypothesis is decided
(\textsf{ir\_impI\_bool}), and $\Cert$ is an existential that \rga
does not in general prove decided.  At every well-formed $\varphi$
that \rga does not prove, \rga neither proves nor refutes
$\Cert_\emptyset(\bquo{\varphi})$, so $\emptyset \nvdash
\Cert_\emptyset(\bquo{\varphi})\jbool$ (\textsf{iHasCert\_undecided},
a metatheorem; the refutation half goes through
{}[\citen{goedel-paper}, theorem~\textsf{rga\_goedel2}]).  In particular
$\Cert_\emptyset(\bquo{\ConR})$ is not provably decided
(\textsf{iHasCert\_}\allowbreak\textsf{ConRGA\_undecided}).

\begin{thm}[\textsf{decMot\_ladder}]
\label{thm:b-ladder}
If $\Gamma \vdash \ell\jnat$, then $\Gamma \vdash \mathit{M}(\ell)$.
\end{thm}

\noindent
\emph{Proof.}  By \rga's induction rule on the level
(\textsf{decMot\_ladderI}, through \textsf{decMot\_natind}).  At level
$\z$ the decider gives no verdict, so every value and head antecedent is
refuted (\textsf{decMot\_base}).  The step \textsf{decMot\_stepI} derives
$\mathit{M}(\suc v_i)$ from the induction hypothesis $\mathit{M}(v_i)$ at
a fresh variable.  The step unfolds the decider one level at the abstract
$v_i$ and dispatches on the node's shape tag.  The value side goes
through the router \textsf{valRouter} (six shapes), the head side
through \textsf{hdRouter}.  Five arms are premises of the step, named as
slots, because these arms need the full motive at the inner level: the
numeral-disequality slot, the application slot of the value router, and
the three contraction slots ($\beta$, recursor at zero, recursor at a
successor).  The ladder fills these slots with \textsf{iHasCert\_num\_neq},
\textsf{valArm\_oprApp} (whose certificate branch is
Theorem~\ref{thm:b-arm-A}), \textsf{decSlot\_beta}, \textsf{decSlot\_recg}
and \textsf{decSlot\_recr}.  \qed

\medskip
The two exhibited arms follow.  In both arms, the premise ``the operator test
fires'' is the decided test that the node at $c$ has the operator in
question (\textsf{iOprHd}).  The negation arm takes $\Gamma \vdash
V^{\forall}(\ell)$, the value clause at level $\ell$ as an internal
universal (\textsf{decValMot}), which the full motive provides
(\textsf{decValMot\_of\_decMot}).  The arm for the universal takes the
full motive $\mathit{M}(\ell)$, because that arm's push spends the head
clause.

\begin{thm}[\textsf{valArm\_oprNot}]
\label{thm:b-arm-not}
Suppose $\Gamma$ derives $\ell\jnat$, $c\jnat$ and $v\jnat$; that the
negation test fires at $c$; that $\Gamma \vdash \Dec(\suc\ell, c) =
\suc v$; that $\Gamma \vdash V^{\forall}(\ell)$; and that $\Gamma$
derives $\Wfc(c) \ne \z$ and $\IsCode(c) \ne \z$.  Then
$\Gamma \vdash \Cert_\emptyset(\mathit{eq}(c, \Num\cdot v))$.
\end{thm}

\noindent
\emph{Proof.}  Write $a$ for the code of the negated subterm
(\textsf{iKid1}).  The negation clause of the decider reads down
internally (\textsf{readdownV\_oprNot}): either the decider gives $a$ the
value $\cfalse$ at level $\ell$ and $v = \ctrue$, or it gives $a$ the value
$\ctrue$ and $v = \cfalse$.  Each case extends $\Gamma$ by the case's facts.
Take the first case.
\begin{enumerate}
\item The guards pass from the node to its child: $\Wfc(a) \ne \z$
  (\textsf{iWft\_neg1\_kid}) and $\IsCode(a) \ne \z$
  (\textsf{iIsCode\_gK\_one}).
\item The value clause at level $\ell$, instantiated at $a$ and $\cfalse$,
  discharges its three decided antecedents and gives
  $\Cert_\emptyset(\mathit{eq}(a, \Num\cdot\tofn 0))$
  (\textsf{decValMot\_inst}), which transports to
  $\Cert_\emptyset(\mathit{eq}(a, \bquo{\tofn 0}))$
  (\textsf{iHasCert\_verd\_F\_in}).
\item Two proof nodes are appended (\textsf{iHasCert\_notT}, through the
  emission cell \textsf{emit\_cell\_notT}): the rule from $a = \tofn 0$ to
  $\neg a$ (\textsf{pr\_false\_eqE}), then the rule from $\neg a$ to
  $\neg a = \tofn 1$ (\textsf{pr\_true\_eqI}).  This gives
  $\Cert_\emptyset(\mathit{eq}(\mathit{neg}(a), \bquo{\tofn 1}))$.
\item \rga derives $c = \mathit{neg}(a)$ from the firing test
  (\textsf{iOprHd1\_neg\_eta}), and the certificate moves along this
  equation (\textsf{iHasCert\_head\_swap}) and along $v = \ctrue$
  (\textsf{iHasCert\_verd\_T\_out}).
\end{enumerate}
The second case is symmetric, with three appended nodes
(\textsf{iHasCert\_notF}: from $a = \tofn 1$ to $a$, then $\neg\neg a$,
then $\neg a = \tofn 0$).  \qed

\medskip
The arm for the universal is where the decider's certificate search meets
the checker.  At an application node whose head reads as $\cA$, the
decider at level $\suc\ell$ runs two bounded searches on the argument
code $u$ (\textsf{iKid2}): the compiled truth search
$\cmp(\mathit{certTrueIdx})\cdot\ell\cdot u$, for a number up to $\ell$
that codes a proof of $\{\vz\jnat\} \vdash \lft u'\cdot\vz$ ($u'$ the term
with code $u$), and the falsity search
$\cmp(\mathit{certFalseIdx})\cdot\ell\cdot u$, for a proof of a refuted
instance with index up to $\ell$ (\texttt{RGA\_Kit/\allowbreak Code\_Index.thy}).

\begin{thm}[\textsf{aNode\_cert}]
\label{thm:b-arm-A}
Suppose $\Gamma$ derives $\ell\jnat$, $c\jnat$, $v\jnat$ and $h\jnat$;
that the application test fires at $c$; that $\Gamma \vdash
\Dec(\suc\ell, c) = \suc v$ and $\Gamma \vdash \mathit{M}(\ell)$; that
$\Gamma$ derives $\Wfc(c) \ne \z$ and $\IsCode(c) \ne \z$; that the head
cell of $c$'s operator child at level $\ell$ holds $\suc h$
(\textsf{tabCellH}); and that the shape tag of $h$ is that of $\cA$ (tag
$1$, \textsf{hdTag}).  Then
$\Gamma \vdash \Cert_\emptyset(\mathit{eq}(c, \Num\cdot v))$.
\end{thm}

\noindent
\emph{Proof.}  Both search results are numerals, so \rga splits on
whether each search result is zero (\textsf{ir\_zero\_dich}, twice).
\begin{enumerate}
\item \emph{The truth search hits.}  The step function then writes
  $\Some{\ctrue}$ (\textsf{ir\_stepV\_app\_cert\_hitT}), so $v = \ctrue$.  The
  hit unpacks to a tree $p$ below $\ell$ that passes the range check and
  the validity test and concludes the open instance
  (\textsf{ir\_certTrue\_elim}).  These are exactly the data of
  \textsf{iHasCertI}, so \rga derives $\Cert$ of the judgment
  $\{\vz\jnat\} \vdash \lft u'\cdot\vz$ (\textsf{certTrue\_hasCert}).  One
  appended $\forall$-introduction node (\textsf{pr\_allI}) gives
  $\Cert_\emptyset(\mathit{ap}(\bquo{\cA}, u))$
  (\textsf{iHasCert\_allT}), and a true-equation node gives
  $\Cert_\emptyset(\mathit{eq}(\mathit{ap}(\bquo{\cA}, u), \bquo{\tofn 1}))$
  (\textsf{iHasCert\_trueEqI}; together \mbox{\textsf{aNode\_cert\_true}}).
\item \emph{The truth search misses and the falsity search hits.}  Then
  $v = \cfalse$ (\textsf{ir\_stepV\_app\_cert\_hitF}).  The hit unpacks to a
  tree certifying $\emptyset \vdash \neg(u'\cdot\tofn i)$ for some
  $i \le \ell$.  \rga also certifies $\tofn i\jnat$ at the numeral code
  (\textsf{iHasCert\_numRefl\_code}, from an internal induction,
  \textsf{iHasCert\_numRefl}).  A negated-universal introduction node
  (\textsf{pr\_not\_allI}) and a false-equation node give
  $\Cert_\emptyset(\mathit{eq}(\mathit{ap}(\bquo{\cA}, u),
  \bquo{\tofn 0}))$ (\textsf{iHasCert\_allF}, \mbox{\textsf{aNode\_cert\_false}}).
\item \emph{Both searches miss.}  The step function then writes $0$
  (\textsf{ir\_stepV\_app\_cert\_miss}), which contradicts
  $\Dec(\suc\ell, c) = \suc v$.  The case closes by \textsf{ir\_neE}.
\end{enumerate}
In the first two cases the certificate is about the canonical node
$\mathit{ap}(\bquo{\cA}, u)$, while the hypothesis is about $c$, whose
operator child only reads as $\cA$ at its head.  The push
\textsf{aNode\_push} closes this gap.  The push spends the head clause of
$\mathit{M}(\ell)$ at the operator child and transports the certificate
along the resulting equation of codes.  \qed

\medskip
The decider's verdict at a universal is thus the certificate the decider found,
and the checker accepts that certificate under one more rule.  No
instance of the universal is evaluated, inside or outside \rga.
Table~\ref{tab:b-arms} lists the remaining arms.

\begin{table}[ht]
\centering
\small
\begin{tabular}{@{}>{\raggedright\arraybackslash}p{0.27\linewidth}>{\raggedright\arraybackslash}p{0.30\linewidth}>{\raggedright\arraybackslash}p{0.37\linewidth}@{}}
\hline
arm & decider clause mirrored & certificate step \\
\hline
\textsf{valArm\_opr0} & value of $\z$ is $0$ & the rule $\z\jnat$, i.e.\ $\z = \z$ (\textsf{iHasCert\_zero\_nat}) \\
\textsf{valArm\_oprSuc} & value of $\suc a$ is $r+1$ & successor step (\textsf{iHasCert\_sucI}) \\
\textsf{valArm\_oprEq} & equation, equal or different values & \textsf{iHasCert\_eqI\_true}, \textsf{iHasCert\_eqI\_false}; numeral disequality from the slot filled by \textsf{iHasCert\_num\_neq} (internal double induction) \\
\textsf{valArm\_oprOr} & disjunction, strong Kleene & \textsf{iHasCert\_orT1}, \textsf{iHasCert\_orT2}, \textsf{iHasCert\_orF} \\
\textsf{valArm\_oprNot} & negation & Theorem~\ref{thm:b-arm-not} \\
\textsf{valArm\_oprApp}, $\beta$ branch & head $\lambda b$: value of $b[u]$ & slot \textsf{decSlot\_beta}; head-reduction node (\textsf{iHasCert\_hev\_appLam}) \\
\textsf{valArm\_oprApp}, $\cR$ at zero & head $\cR\cdot g\cdot h$, counter $0$ & slot \textsf{decSlot\_recg} (\textsf{iHasCert\_hev\_appR0}) \\
\textsf{valArm\_oprApp}, $\cR$ at successor & head $\cR\cdot g\cdot h$, counter $k+1$ & slot \textsf{decSlot\_recr} (\textsf{iHasCert\_hev\_appRS}) \\
\textsf{valArm\_oprApp}, $\cA$ branch & the two certificate clauses & Theorem~\ref{thm:b-arm-A} \\
\textsf{hdArm\_base\_lam}, \textsf{\_R}, \textsf{\_A} & $\lambda b$, $\cR$, $\cA$ are their own head forms & transport along the canonical code (\textsf{hdArm\_canon\_close}) \\
\textsf{hdArm\_oprApp} & application: head of the operator, spine grown & spine growth (\textsf{hdArm\_grow}); contractions via the three slots \\
\hline
\end{tabular}
\caption{The arms of the internal completeness induction.  The value arms
  are in \texttt{RGA\_Checker/\allowbreak Chk\_Transfers.thy} except the application
  arm; the application arm, the slots and the head arms are in
  \texttt{RGA\_Checker/\allowbreak Chk\_Completeness.thy}.}
\label{tab:b-arms}
\end{table}

\paragraph{From the ladder to the sentence.}
The value projection \textsf{iHasCert\_of\_iPrDecV} instantiates the
ladder: if $\Gamma$ derives $\ell\jnat$, $c\jnat$, $v\jnat$,
$\Dec(\ell, c) = \suc v$, $\Wfc(c) \ne \z$ and $\IsCode(c) \ne \z$, then
$\Gamma \vdash \Cert_\emptyset(\mathit{eq}(c, \Num\cdot v))$.  At
$v = \ctrue$ one more appended node, from $t = \tofn 1$ to $t$
(\textsf{pr\_true\_eqE}), gives $\Cert_\emptyset(c)$
(\textsf{iHasCert\_reflect\_true}, \textsf{iHasCert\_of\_iPrDecV\_true}).

\begin{thm}[\textsf{iHasCert\_of\_iPrDecTrueAt}]
\label{thm:b-cert-rule}
If $\Gamma$ derives $L\jnat$, $c\jnat$, $\Aff(L, c) \ne \z$,
$\Wfc(c) \ne \z$ and $\IsCode(c) \ne \z$, then
$\Gamma \vdash \Cert_\emptyset(c)$.
\end{thm}

\noindent
\emph{Proof.}  The affirmation test reads back as
$\Dec(L, c) = \tofn 2 = \suc\ctrue$ (\textsf{iPrDecTrueAt\_true\_elim}), and
\textsf{iHasCert\_of\_iPrDecV\_true} applies.  \qed

\medskip
\noindent
\emph{Proof of Theorem~\ref{thm:pl-int}} (\textsf{rga\_complete\_int2}).
Open the eigenvariables $s$ and $c$ with $s\jnat$ and $c\jnat$.  The
antecedent is the 0/1-valued conjunction of $\Wfc(c)$, $\IsCode(c)$
and $\Aff(s, c)$, and that conjunction is decided (\textsf{andT\_B}).
The decided antecedent is discharged (\textsf{ir\_impI\_bool}); as a
hypothesis, the antecedent gives the three nonzero tests (\textsf{andT2\_ne0D}), and
Theorem~\ref{thm:b-cert-rule} gives $\Cert_\emptyset(c)$.  Two
universal introductions (\textsf{ir\_allI\_open}) close the sentence.
\qed

\noindent
The development also states the sentence over one packed pair $k$ of
a fuel and a code (\textsf{rga\_complete\_int}), with the same proof at
one eigenvariable and the packed test \textsf{iPrDecTruePair}.  The
instance at the code of $\z\jnat$ (\textsf{rga\_complete\_int\_zeroN})
is proved directly from Theorem~\ref{thm:b-cert-rule}, at a numeral
fuel at which the decider affirms $\z\jnat$ (\textsf{prDecV\_complete},
\textsf{iPrDecV\_closed}).

\medskip
\noindent
\emph{Proof of Theorem~\ref{thm:pl-int-ent}} (\textsf{iHasCert\_of\_iWtrue},
\textsf{complete\_ent}, \textsf{complete\_imp\_fails}).  For the rule
in every context, eliminate the hypothesis $\Tru(c) = \exists s.\
\Aff(s, c)$ at a fresh fuel (\textsf{rga\_wtrueE}), which adds
$\vz\jnat$ and the affirmation at $\vz$ to the lifted context; turn
the affirmation into the nonzero test, lift the guards
(\textsf{ir\_lift\_sloppy}), and apply Theorem~\ref{thm:b-cert-rule}
at $L = \vz$ and the lifted code; the conclusion is lifted back
(\textsf{lift\_iHasCert}).  The entailment is the rule at the context
of its hypotheses.  The implication fails at the code of $L_{\Tru}$:
\rga proves both guards there, does not refute
$\Tru(\bquo{L_{\Tru}})$ (\textsf{liar\_code\_not\_refuted\_true},
from the transport of $\neg\Tru$ and the independence of the Liar),
and does not prove $\Cert_\emptyset(\bquo{L_{\Tru}})$
(\textsf{liar\_code\_no\_cert}, since a certificate gives a proof),
so the disjunction property refutes the instance, and a universal with
an underivable instance is underivable (\textsf{not\_all\_of\_inst},
and \textsf{not\_all\_guarded\_of\_inst} for the guarded form).
\textsf{iHasCert\_of\_iWtrue} and \textsf{complete\_ent} are in
\texttt{RGA\_Goedel/}\allowbreak\texttt{Entailment\_Goedel};
\textsf{complete\_imp\_fails}, \textsf{liar\_code\_not\_refuted\_true} and
\textsf{liar\_code\_no\_cert} are in
\texttt{RGA\_Liar\_Div/}\allowbreak\texttt{Entailment\_Liar}.  \qed

\medskip
Two further facts relate $\Cert$ to the other machines.  In HOL, every
witness that $\pi$ accepts is the code of a proof
(Theorem~\ref{thm:pl-genuine}).  Inside \rga, an accepted tree of the tree
checker is a certificate: if $\Gamma \vdash p\jnat$ and $\Gamma \vdash
\Chk(p, c)$, then $\Gamma \vdash \Cert_\emptyset(c)$
(\textsf{iHasCert\_of\_iChk}, \texttt{RGA\_Goedel/\allowbreak Goedel\_II.thy}).
Theorem~\ref{thm:b-cert-rule} is used by the search route to G\"odel's
second theorem (\textsf{rga\_goedel2\_via\_search}, through
\textsf{iHasCert\_negInstC}, [\citen{goedel-paper}, Appendix ``The search route'']).  The
headline proof \textsf{rga\_goedel2} uses neither
Theorem~\ref{thm:b-cert-rule} nor Theorem~\ref{thm:pl-int}, nor the
decider or the truth predicate.

\subsection{The Liar, unestablishment and tracking}
\label{app:pillars-liar}

The definitions of $D_{\Tru}$, $S_{\Tru}$ and $L_{\Tru}$ in
\S\ref{sec:pillars-tarski} are \textsf{lwD}, \textsf{goedelWS} and
\textsf{goedelW} in \texttt{RGA\_Liar\_Div/\allowbreak Liar\_Record.thy}.  The
compiled diagonal function is \textsf{iDiagSub}
(\texttt{RGA\_Goedel/\allowbreak Goedel\_II.thy}), whose closed law
\textsf{iDiagSub\_closed} assumes that $t$ is well-formed.  The identity
$\ssub{0}{\bquo{D_{\Tru}}}{D_{\Tru}} = L_{\Tru}$ (\textsf{goedelW\_diag}) holds because
substitution passes through $\Tru$ and the compiled diagonal unchanged,
and the equation $\emptyset \vdash S_{\Tru} = \bquo{L_{\Tru}}$
(\textsf{goedelWS\_eq}) is \textsf{iDiagSub\_closed} at $t = D_{\Tru}$ and
$k = \quo{D_{\Tru}}$.  Transport along that equation, in both directions, gives
the two lemmas \textsf{goedelW\_to\_code} and \textsf{goedelW\_from\_code}:
$\emptyset \vdash \Tru(S_{\Tru})$ exactly when $\emptyset \vdash
\Tru(\bquo{L_{\Tru}})$.

The proof of Theorem~\ref{thm:pl-no-value} (\textsf{goedelW\_no\_value})
in \S\ref{sec:pillars-tarski} follows the Isabelle proof step by step.
If $L_{\Tru} = \neg\Tru(S_{\Tru})$ had a value, then $\Tru(S_{\Tru})$ would have a
value $b$ (\textsf{wosr\_not\_inv}).  For $b = \ctrue$, open completeness
(\textsf{rga\_open\_iff}, \textsf{wosr\_osr}) derives $\Tru(S_{\Tru})$,
\textsf{goedelW\_to\_code} and Theorem~\ref{thm:pl-adeq}(2) make
$L_{\Tru}$ grounded-true, and \textsf{wosr\_unique} contradicts the value
$\cfalse$.  For $b = \cfalse$, $L_{\Tru}$ is grounded-true and derivable,
Theorem~\ref{thm:pl-adeq}(1) and \textsf{goedelW\_from\_code} derive
$\Tru(S_{\Tru})$, and \textsf{rga\_consistent} refutes the pair.  The
independence results of Theorem~\ref{thm:pl-liar-indep} follow from
\textsf{no\_value\_independent} (\texttt{RGA/\allowbreak RGA\_Paradox.thy}): a
well-formed term with no value is neither provable nor refutable.  The
undecidedness \textsf{goedelW\_undecided} adds the disjunction property
(\textsf{ir\_disj\_prop}).  [\citen{goedel-paper}, Section ``G\"odel's proof in four forms''], on G\"odel's proof in four forms, states the
$\omega$-incompleteness of $L_{\Tru}$ (Theorem~\ref{thm:br-omega}).

\paragraph{The explosive T-biconditional.}
Theorem~\ref{thm:pl-tb-generic}(1) is \textsf{tbicond\_liar\_explosive\_gen}
in \texttt{RGA\_Liar\_Div/\allowbreak Entailment\_Liar.thy}.  The
biconditional is \rga's material one (\textsf{trmIff\_def}): the
conjunction of $T(\bquo{L}) \limp L$ and $L \limp T(\bquo{L})$, each a
disjunction.  Conjunction elimination gives both disjunctions under
the hypothesis.  In the first disjunction, both cases give $L$: the left case by
transport of $\neg T$ along $S = \bquo{L}$, since $L$ is $\neg T(S)$;
the right case is the case hypothesis.  Transport is full equality
substitution (\textsf{ir\_eq\_subs}) at the term $\neg\varphi$, whose
instances at $\bquo{L}$ and at $S$ are $\neg T(\bquo{L})$ and $L$.  In
the second disjunction, both cases contradict $L$, the right one after
transport in the other direction, and negation elimination concludes
$\z = \tofn 1$.  No decidedness and no property of $\varphi$ beyond
well-formedness is used.  The form for a predicate term applied to its
argument, \textsf{tbicond\_liar\_explosive\_pred}, is the instance at
$\varphi = T\cdot v$ with $v$ fresh for $T$.
Theorem~\ref{thm:pl-tb-explosive} is \textsf{tbicond\_liar\_explosive}
in \texttt{RGA\_Liar\_Div/\allowbreak Entailment\_Truth.thy}, proved
as the instance of \textsf{tbicond\_liar\_explosive\_gen} at $\varphi =
\Tru(\vz)$ and $S = S_{\Tru}$, with the code equation
\textsf{goedelWS\_eq} (substitution passes through $\Tru$,
\textsf{sub1\_iWtrue}).  \textsf{tbicond\_liar\_both} adds the two
consequences $L_{\Tru}$ and $\neg L_{\Tru}$, and
\textsf{tbicond\_liar\_unprovable} cuts with
\textsf{rga\_no\_proof\_0eq1}.  For part~2, \textsf{pred\_liar\_eq}
builds the diagonal Liar of a closed well-formed $T$ as $L_{\Tru}$ is
built, by \textsf{iDiagSub\_closed} at $t = D$; closedness of $T$ makes
the substitution of $\bquo{D}$ for $\vz$ in $D$ pass through $T$.
\textsf{tbicond\_liar\_explosive\_diag} weakens that equation into
$\Gamma$ and applies the predicate form, and
\textsf{tbicond\_liar\_unprovable\_diag} cuts with
\textsf{rga\_no\_proof\_0eq1}, so that theorem uses soundness, while the
entailments of Theorems~\ref{thm:pl-tb-explosive}
and~\ref{thm:pl-tb-generic} and the code equation do not (those facts are
roots of the dependency check of Appendix~\ref{app:machines-audit}).

\paragraph{The ordinary Liar.}
Theorem~\ref{thm:pl-liar-neither}(1) is \textsf{ordinary\_liar\_neither}
in \texttt{RGA\_Liar\_Div/\allowbreak Truth\_Limits.thy}, for the
ordinary Liar $L = D\cdot D$, with $D = \lambda x.\,F\cdot(x\cdot x)$
and $F = \lambda y.\,\neg y$ (\textsf{trmLiar}).  The term $L$ head-steps
through $F\cdot(D\cdot D)$ to $\neg L$.  The
two refutations come from \textsf{liar\_never\_affirmed} and
\textsf{neg\_liar\_never\_affirmed}, both proved by one internal
induction on the fuel.  A decided value of $L$ at fuel $w+4$ becomes,
by two $\beta$-laws of the decider's table read downwards
(\textsf{beta\_val\_down}, \textsf{liar\_val\_down}), a decided value
of $\neg L$ at fuel $w+2$, and by the negation law
(\textsf{neg\_val\_down\_T}, \textsf{neg\_val\_down\_F}) the opposite
decided value of $L$ at fuel $w+1$.  So the motive is a window of four
consecutive levels at which $L$ is neither 1 nor 2
(\textsf{liar\_never\_decided}), and the base levels 0 to 3 are closed
computations.  The base uses the decider's correctness on
$\forall$-free terms only (\textsf{liar\_slice\_None}), so the internal
halves use no soundness theorem.  An affirmation of $\neg L$ at some
fuel decides $L$ false one level down, which the same induction
excludes.  The two non-derivability facts use soundness.  Part~2 is
\textsf{truth\_liar\_not\_untrue\_unprovable}, an alias of
\textsf{liar\_code\_not\_refuted\_true}.

\paragraph{Refuted and unestablished.}
Theorem~\ref{thm:pl-refuted} is \textsf{refuted\_unestablished\_trichotomy}
with its parts, in \texttt{RGA\_Liar\_Div/\allowbreak Truth\_Limits.thy}.
Part~1 is \textsf{refuted\_unestablished\_qf}, the $\forall$-free
decider collision.  Part~3 rests on the internal law that the decider
affirms a universal at a fuel above a certificate of its open instance
(\textsf{ir\_A\_cert\_affirm}, with the hit introduction
\textsf{ir\_certTrue\_intro}).  From an accepted proof $a$ of $\z =
\tofn 1$, ex falso and one weakening step build that certificate, so
$\Gamma \vdash \Tru(\bquo{\cA\cdot u})$ at an abstract $a$
(\textsf{iWtrue\_univ\_of\_chk01}), hence $\neg\ConR \vdash
\Tru(\bquo{\cA\cdot u})$ (\textsf{notcon\_tru\_univ\_ent}).  Then the
argument of \textsf{ConRGA\_of\_code\_unprovable}, a decided split on
$\Chk(v, \bquo{\z = \tofn 1})$ at an eigenvariable, turns $\vdash
\neg\Tru(\bquo{\cA\cdot u})$ into $\vdash \ConR$
(\textsf{universal\_unestablished\_needs\_con}).  The premise is used
only through weakening, so the argument runs in every context
(\textsf{universal\_unestablished\_needs\_con\_ctx}): lifted under the
eigenvariable's binder, $\neg\Tru(\bquo{\cA\cdot u})$ is unchanged,
since the argument of $\Tru$ is a closed code numeral.  Instantiating the context
with the premise gives the entailment of part~3,
$\neg\Tru(\bquo{\cA\cdot u}), \Gamma \vdash \ConR$
(\textsf{universal\_unestablished\_con\_ent}), at every well-formed $u$,
refuted or not.  Part~2 uses internal
completeness: a certificate of $t$ is an accepted proof
(\textsf{iPrc\_of\_iHasCert\_code}), so $\Tru(\bquo{t}) \vdash
\PrC(\bquo{t})$ (\textsf{tru\_prc\_ent}); with the refutation, D1
without soundness and EXPL give $\Tru(\bquo{t}) \vdash \neg\ConR$
(\textsf{refuted\_tru\_notcon\_ent}); $\Sigma$-contraposition, at the
decided matrix of $\Tru$, and double negation give $\ConR \vdash
\neg\Tru(\bquo{t})$ (\textsf{con\_refuted\_unestablished}).  Part~4
combines parts~2 and~3 into the two entailments
(\textsf{refuted\_universal\_unestablished\_eqv\_con}); in derivability
form, in every context $\Gamma$, $\Gamma \vdash \neg\Tru(\bquo{\cA\cdot
u})$ holds exactly when $\Gamma \vdash \ConR$ does
(\textsf{refuted\_universal\_unestablished\_iff\_con\_ctx}).  The
instance of this derivability form at the empty context
(\textsf{refuted\_universal\_unestablished\_iff\_con}) relates two
derivability facts that both fail, by G\"odel's second theorem, so the
paper states the entailments; at $\Gamma = \{\ConR\}$ both sides hold.
$\ConR \eqv \neg\Tru(\bquo{\ConR})$ (\textsf{con\_eqv\_not\_true\_con})
has \textsf{con\_con\_unestablished} as its first half, and
\textsf{universal\_unestablished\_con\_ent} at $u = \lambda\neg\Chk(\vz,
\bquo{\z = \tofn 1})$ as its second, since $\ConR$ is literally
$\cA\cdot u$ for that $u$.  \textsf{refuted\_unestablished\_fails}
instantiates part~4 at $\forall x.\, x = \tofn 1$ with G\"odel's
second theorem.

\paragraph{Tracking under $\ConR$.}
Theorem~\ref{thm:pl-tracking} collects results of the same file.  In
part~1 the biconditionals are built from two derivations by the
introduction of a material conjunction of implications: at a proved
$t$, both sides are derived and each implication is introduced over its
derived consequent; at a refuted $t$, both sides are refuted and each
implication is introduced over its refuted antecedent.  Neither
introduction discharges a hypothesis, so no decidedness is needed.  At
a proved $t$, $\vdash \PrC(\bquo{t})$ is D1 without soundness
(\textsf{prc\_of\_derivable}) and $\vdash \Tru(\bquo{t})$ is NEC
(\textsf{iWtrue\_nec}), which goes through soundness and the decider's
completeness.  At a refuted $t$, $\ConR \vdash \neg\PrC(\bquo{t})$ and
$\ConR \vdash \neg\Tru(\bquo{t})$ are Theorem~\ref{thm:pl-refuted}(2)
(\textsf{con\_refuted\_unprovable}, \textsf{con\_refuted\_unestablished}).
These derivations give \textsf{con\_prc\_schema\_decided} and
\textsf{con\_tschema\_decided}, with \textsf{con\_tschema\_refuted} the
refuted half of \textsf{con\_tschema\_decided}; the Boolean forms
\textsf{con\_prc\_bool\_decided} and \textsf{con\_true\_bool\_decided}
take the same case split, with $\lor$-introduction in place of the
biconditional.  Part~2 is \textsf{con\_tracking\_sound}.  If $\ConR
\vdash \Tru(\bquo{t})$, the $\omega$-soundness of \mbox{\rga{}${}+\ConR$}
(\textsf{con\_ctx\_omega\_sound}) makes $\Tru(\bquo{t})$ true in the
$\omega$-semantics, which excludes a $t$ that is never grounded-true
(\textsf{omv\_iWtrue\_not\_true}).  If $\ConR \vdash \neg\Tru(\bquo{t})$
and $t$ were grounded-true, adequacy would give $\emptyset \vdash
\Tru(\bquo{t})$ (\textsf{iWtrue\_of\_wosr}), and with that derivation $\ConR \vdash
\z = \tofn 1$, which \textsf{con\_ent\_0eq1\_fails} excludes.  For
part~3, \textsf{unestablished\_of\_unprovable} takes $\Gamma \vdash
\neg\PrC(\bquo{t})$; internal completeness gives $\Tru(\bquo{t}),
\Gamma \vdash \PrC(\bquo{t})$ (\textsf{tru\_prc\_ent}), so
$\Tru(\bquo{t}), \Gamma \vdash \neg\Tru(\bquo{t})$, and reductio at a
hypothesis of the $\Sigma$ shape (\textsf{sigma\_reductio}, at the
decided matrix of $\Tru$) gives $\Gamma \vdash \neg\Tru(\bquo{t})$; no
$\ConR$ is used in this step.  At $G$, the premise is
\textsf{con\_notprc\_goedelG\_ent} and the derivation of $G$ is
\textsf{goedelG\_under\_ConRGA}
(\textsf{con\_goedelG\_derivable\_unestablished}), and conjunction
introduction makes the two facts one \rga sentence, $\ConR, \Gamma \vdash G
\wedge \neg\Tru(\bquo{G})$ (\textsf{con\_goedelG\_and\_unestablished}).
That conjunction is true in the $\omega$-semantics: \mbox{\rga{}${}+\ConR$}
is sound for that conjunction (\textsf{con\_ctx\_omega\_sound}), $G$ is $\omega$-true
and $\Tru(\bquo{G})$ is $\omega$-false, since $G$ is never grounded-true
(\S\ref{sec:pillars-tarski}, ``Asserting what is not established'').  The refuted
biconditional (\textsf{con\_tschema\_goedelG\_refuted}) refutes its
second implication, $G \limp \Tru(\bquo{G})$, from the two conjuncts.  At $\ConR$,
the premise is \textsf{rga\_goedel2\_internal}
(\textsf{con\_con\_unestablished}).  These unestablishment results are
consistent with part~2:
$G$ has no value (\textsf{goedelG\_no\_value}), and \mbox{\rga{}${}+\ConR$}
refutes the truth of $G$.  Part~2 and the $\Tru$ forms of part~1 at proved
sentences rest on $\omega$-soundness or soundness.  The $\PrC$ forms of
part~1, \textsf{con\_tschema\_refuted},
\textsf{unestablished\_of\_unprovable} and the three results of part~3
are roots of the dependency check of
Appendix~\ref{app:machines-audit}.

\paragraph{No sentence together with its untruth.}
\rga itself derives no conjunction of the kind that part~3 of
Theorem~\ref{thm:pl-tracking} derives under $\ConR$.

\begin{thm}[\textsf{no\_moorean\_theorem}, \textsf{no\_moorean\_pair}]
\label{thm:b-no-moorean}
For every term $\varphi$: if $\emptyset \vdash \varphi$, then
$\emptyset \nvdash \neg\Tru(\bquo{\varphi})$; and $\emptyset \nvdash
\varphi \wedge \neg\Tru(\bquo{\varphi})$.
\end{thm}

\noindent
\emph{Proof.}  If $\emptyset \vdash \varphi$, then $\emptyset \vdash
\Tru(\bquo{\varphi})$ by NEC (Theorem~\ref{thm:pl-nec},
\textsf{iWtrue\_nec}), and a derivation of $\neg\Tru(\bquo{\varphi})$
would contradict \rga's consistency (\textsf{rga\_consistent}).  The
second statement follows by conjunction elimination.  Both theorems are in
\texttt{RGA\_Liar\_Div/\allowbreak Truth\_Limits.thy}; they are meta
theorems and use soundness, through consistency, and the decider's
completeness, through NEC.  \qed

\subsection{$\Tru$ against the axioms of truth}
\label{app:pillars-axioms}

This subsection gives the comparison of \S\ref{sec:pillars-axioms}
cell by cell (Table~\ref{tab:c-truth-axioms}), with the theorems that
settle each cell.  The column meanings are those of
\S\ref{sec:pillars-axioms}.  A rule in every context is an
entailment, by instantiating its context with the premises
(\textsf{ent\_of\_ctx\_rule}), so those cells hold as entailments too.

In the table the variables $x$, $y$, $f$, $a$ and $b$ range over
numerals, and their universal quantifiers are left implicit.
$\mathit{neg}(c)$, $\mathit{or}(a, b)$, $\mathit{eq}(a, b)$ and
$\mathit{ap}(a, b)$ are the compiled code constructors of
{}[\citen{goedel-paper}, Appendix ``G\"odel II in full''], $\bquo{\cA}$ is the compiled spelling of the
code of $\cA$, which \rga proves equal to the code numeral of $\cA$
(\textsf{mkAS\_val}), and
$\Num\cdot n$ is the code of the numeral of $n$.  The code
$\mathit{tr}(x) \ldef \Diag(\bquo{\Tru(\vz)}, x)$ is the code of $\Tru$
applied to the numeral $x$, so that $\emptyset \vdash
\mathit{tr}(\bquo{t}) = \bquo{\Tru(\bquo{t})}$.  \rga's conjunction is
defined from $\neg$ and $\lor$, so KF's clause for $\wedge$ appears as
the clause for $\neg\lor$.  For a clause $\varphi \liff \psi$, the
direction $\rightarrow$ goes from $\varphi$ to $\psi$, and
$\leftarrow$ back.  The mark $\dagger$ means ``at codes of
$\forall$-free terms only''.  Finally, $\Omega \ldef
(\lambda(\vz\cdot\vz))\cdot(\lambda(\vz\cdot\vz))$ is a well-formed,
$\forall$-free term that head-steps to itself and has no value
(\textsf{omega\_\allowbreak no\_\allowbreak value}).  The refuted
universal $U \ldef \cA\cdot(\lambda(\vz = \tofn 1))$, that is,
$\forall x.\, x = \tofn 1$, is the sentence of
Theorem~\ref{thm:pl-refuted}(4).

{\footnotesize
\setlength{\tabcolsep}{2.5pt}
\renewcommand{\arraystretch}{1.0}
\begin{longtable}{@{}>{\raggedright\arraybackslash}p{0.085\linewidth}>{\raggedright\arraybackslash}p{0.21\linewidth}>{\raggedright\arraybackslash}p{0.235\linewidth}>{\raggedright\arraybackslash}p{0.225\linewidth}>{\raggedright\arraybackslash}p{0.19\linewidth}@{}}
\caption{$\Tru$ against the axioms and rules of truth, in full
  (\S\ref{sec:pillars-axioms}; theory files in
  Table~\ref{tab:b-files}).  The column ``entailment'' gives the form
  $A \vdash B$, with the code variables typed in the context;
  $\dagger$: at codes of $\forall$-free terms only; $U$ is the
  refuted universal $\forall x.\, x = \tofn 1$.}
\label{tab:c-truth-axioms}\\
\hline
principle & \rga form & status for $\Tru$ & entailment & Isabelle \\
\hline
\endfirsthead
\hline
principle & \rga form & status for $\Tru$ & entailment & Isabelle \\
\hline
\endhead
\hline
\endfoot
\multicolumn{5}{@{}l}{\textit{Compositional clauses}} \\*
KF: $=$, $\neq$ &
  $\Tru(\mathit{eq}(\Num\cdot a, \Num\cdot b)) \liff{}$\linebreak
  \mbox{$a = b$}; with
  $\mathit{neg}(\mathit{eq}(\dots))$, $a \neq b$ &
  internal &
  holds, from the internal form &
  \textsf{kf\_\allowbreak eq\_\allowbreak int}, \textsf{kf\_\allowbreak ne\_\allowbreak int} \\
KF: $\neg\neg$ &
  $\Tru(\mathit{neg}(\mathit{neg}(x))) \liff \Tru(x)$ &
  rule in every context ($\rightarrow$, $\leftarrow$); internal: fails, as does its
  $\leftarrow$ half &
  holds ($\rightarrow$, $\leftarrow$) (\textsf{kf\_\allowbreak dneg\_\allowbreak ent}) &
  \textsf{iWtrue\_\allowbreak dneg\_\allowbreak elim}, \textsf{iWtrue\_\allowbreak dneg\_\allowbreak intro};
  \textsf{kf\_\allowbreak dneg\_\allowbreak iff\_\allowbreak not\_\allowbreak internal}, \textsf{kf\_\allowbreak dneg\_\allowbreak not\_\allowbreak internal} \\
KF: $\lor$ &
  $\Tru(\mathit{or}(x, y)) \liff \Tru(x) \lor \Tru(y)$ &
  rule in every context ($\rightarrow$, $\leftarrow$); internal: fails, as does its
  $\leftarrow$ half &
  holds ($\rightarrow$, $\leftarrow$) (\textsf{kf\_\allowbreak or\_\allowbreak ent}) &
  \textsf{iWtrue\_\allowbreak or\_\allowbreak elim}, \textsf{iWtrue\_\allowbreak or\_\allowbreak intro1}, \textsf{iWtrue\_\allowbreak or\_\allowbreak intro2};
  \textsf{kf\_\allowbreak or\_\allowbreak iff\_\allowbreak not\_\allowbreak internal}, \textsf{kf\_\allowbreak or\_\allowbreak not\_\allowbreak internal} \\
KF: $\neg\lor$ &
  $\Tru(\mathit{neg}(\mathit{or}(x, y))) \liff \Tru(\mathit{neg}(x))
  \wedge \Tru(\mathit{neg}(y))$ &
  $\rightarrow$: rule in every context; $\leftarrow$: rule in every
  context${}^\dagger$, open beyond $\dagger$; internal: fails, as does
  its $\leftarrow$ half &
  $\rightarrow$: holds; $\leftarrow$: holds${}^\dagger$, open beyond
  $\dagger$ &
  \textsf{iWtrue\_\allowbreak disjE2\_\allowbreak elim}, \textsf{iWtrue\_\allowbreak disjE3\_\allowbreak elim};
  \textsf{iWtrue\_\allowbreak or\_\allowbreak intro3\_\allowbreak qf}; \textsf{kf\_\allowbreak negor\_\allowbreak iff\_\allowbreak not\_\allowbreak internal},
  \textsf{kf\_\allowbreak negor\_\allowbreak not\_\allowbreak internal} \\
FS: $\neg$ &
  $\Tru(\mathit{neg}(x)) \liff \neg\Tru(x)$ &
  $\rightarrow$: rule in every context${}^\dagger$; beyond $\dagger$
  fails as a rule on theorems (at $U$), hence in every context;
  $\leftarrow$: fails as a rule on theorems (at $\Omega$);
  internal: fails, as does its
  $\leftarrow$ half &
  $\rightarrow$: holds${}^\dagger$; fails at $U$, and with the code
  variable free (\textsf{fs\_\allowbreak neg\_\allowbreak fwd\_\allowbreak ent\_\allowbreak fails\_\allowbreak univ},
  \textsf{fs\_\allowbreak neg\_\allowbreak fwd\_\allowbreak ent\_\allowbreak fails});
  $\leftarrow$: fails (at $\Omega$) (\textsf{fs\_\allowbreak neg\_\allowbreak ent\_\allowbreak fails}) &
  \textsf{iWtrue\_\allowbreak refute\_\allowbreak of\_\allowbreak neg\_\allowbreak qf}; \textsf{fs\_\allowbreak neg\_\allowbreak fwd\_\allowbreak rule\_\allowbreak fails}; \textsf{fs\_\allowbreak neg\_\allowbreak rule\_\allowbreak fails};
  \textsf{fs\_\allowbreak neg\_\allowbreak iff\_\allowbreak not\_\allowbreak internal}, \textsf{fs\_\allowbreak neg\_\allowbreak not\_\allowbreak internal} \\
KF, FS: $\forall$ &
  $\Tru(\mathit{ap}(\bquo{\cA}, f)) \liff \forall n.\
  \Tru(\mathit{ap}(f, \Num\cdot n))$ &
  $\rightarrow$: rule on theorems, instance by instance; $\leftarrow$:
  from a proof of the open instance (Thm.~\ref{thm:b-allI}), but no
  $\omega$-rule; in every context: $\rightarrow$ fails, since that rule would
  give the entailment; $\leftarrow$ fails, at the matrix of $L_{\Tru}$
  (Thm.~\ref{thm:c-tin-omega}); internal: fails, as does its
  $\leftarrow$ half; the certificate form of $\leftarrow$, with the
  premise $\cmp(\pi)\cdot J(f)\cdot p \neq \z$, holds internally, and
  so does its converse (Thm.~\ref{thm:pl-tin-cert}) &
  $\rightarrow$: fails at $U$, refuted at $n = 0$
  (\textsf{kf\_\allowbreak all\_\allowbreak fwd\_\allowbreak ent\_\allowbreak fails});
  $\leftarrow$, a formalized $\omega$-rule from
  one universal hypothesis: fails, at the matrix of $L_{\Tru}$
  (\textsf{all\_\allowbreak omega\_\allowbreak ent\_\allowbreak fails}) &
  \textsf{true\_\allowbreak allE\_\allowbreak rule}; \textsf{iWtrue\_\allowbreak allI\_\allowbreak of\_\allowbreak proof};
  \textsf{all\_\allowbreak omega\_\allowbreak rule\_\allowbreak fails}; \textsf{all\_\allowbreak omega\_\allowbreak ctx\_\allowbreak rule\_\allowbreak fails}; \textsf{kf\_\allowbreak all\_\allowbreak not\_\allowbreak internal},
  \textsf{kf\_\allowbreak all\_\allowbreak iff\_\allowbreak not\_\allowbreak internal};
  \textsf{tru\_\allowbreak univ\_\allowbreak cert\_\allowbreak int}, \textsf{tru\_\allowbreak univ\_\allowbreak cert\_\allowbreak conv\_\allowbreak int} \\
KF: $\neg\forall$ &
  $\Tru(\mathit{neg}(\mathit{ap}(\bquo{\cA}, f))) \liff{}$\linebreak
  $\exists n.\ \Tru(\mathit{neg}(\mathit{ap}(f, \Num\cdot n)))$ &
  rule on theorems ($\rightarrow$, $\leftarrow$) at $f = \bquo{u}$, $u$
  well-formed; in every context and internal: open &
  open &
  \textsf{true\_\allowbreak notall\_\allowbreak rule} \\
head step &
  $\Tru(\bquo{c}) \liff \Tru(\bquo{c'})$ for $c \leadsto c'$ &
  rule on theorems ($\rightarrow$, $\leftarrow$); rule in every
  context${}^\dagger$ ($\rightarrow$, $\leftarrow$) for
  $\Tru(\Subst(z, u)) \liff \Tru(\Subst(z, v))$ at code terms $u$,
  $v$ with $\Gamma \vdash \mathit{hred}(u) = \suc v$ (the compiled head
  step), under an environment $z$ of well-formed $\forall$-free terms;
  internal: open &
  holds${}^\dagger$ ($\rightarrow$, $\leftarrow$), in the compiled
  form; open beyond $\dagger$ &
  \textsf{true\_\allowbreak hstep\_\allowbreak rule}; \textsf{hstep\_\allowbreak iWtrue\_\allowbreak fwd}, \textsf{hstep\_\allowbreak iWtrue\_\allowbreak rev} \\
\multicolumn{5}{@{}l}{\textit{Iterated truth}} \\*
T-Rep, T-Del &
  $\Tru(x) \limp \Tru(\mathit{tr}(x))$;
  $\Tru(\mathit{tr}(x)) \limp \Tru(x)$ &
  rules on theorems; in every context: open; internal T-Rep: fails;
  internal T-Del: open &
  open (T-Rep and T-Del) &
  \textsf{trep\_\allowbreak rule}, \textsf{tdel\_\allowbreak rule}; \textsf{trep\_\allowbreak not\_\allowbreak internal} \\
KF: $\neg T$ &
  $\Tru(\mathit{neg}(\mathit{tr}(x))) \liff \Tru(\mathit{neg}(x))$ &
  rule on theorems ($\rightarrow$, $\leftarrow$) at every well-formed
  $t$ that is grounded-true, or grounded-false and $\forall$-free;
  $\rightarrow$ fails as a rule on theorems (at $\Omega$);
  $\leftarrow$ fails as a rule on theorems at $U$ ($\rightarrow$ is
  vacuous there); other grounded-false $t$, and internal: open &
  $\rightarrow$: fails (at $\Omega$) (\textsf{tneg\_\allowbreak ent\_\allowbreak fails\_\allowbreak omega});
  $\leftarrow$: fails at $U$
  (\textsf{tneg\_\allowbreak rev\_\allowbreak ent\_\allowbreak fails}) &
  \textsf{tneg\_\allowbreak rule\_\allowbreak grounded}; \textsf{tneg\_\allowbreak rule\_\allowbreak fails};
  \textsf{tneg\_\allowbreak rev\_\allowbreak rule\_\allowbreak fails} \\
\multicolumn{5}{@{}l}{\textit{Disquotation, at well-formed $t$}} \\*
T-In &
  $t \limp \Tru(\bquo{t})$ &
  proved exactly at the $t$ that are grounded-true or grounded-false;
  fails at every $t$ with no value &
  holds where the implication does; fails at $L_{\Tru}$, $G$ and
  $\ConR$ (Thm.~\ref{thm:br-tin-ent}); other ungrounded $t$: open &
  \textsf{tin\_\allowbreak iff\_\allowbreak decided}, \textsf{tschema\_\allowbreak comp\_\allowbreak iff\_\allowbreak grounded};
  \textsf{tschema\_\allowbreak unprovable\_\allowbreak ungrounded} \\
T-Out &
  $\Tru(\bquo{t}) \limp t$ &
  proved exactly when \rga proves $t$ or refutes $\Tru(\bquo{t})$; in
  every context at $\forall$-free $t$ that are grounded-true or
  grounded-false, and at $\Omega$; fails at $L_{\Tru}$ and at $U$;
  other ungrounded $t$, and other grounded-false $t$ with quantifiers:
  open &
  holds where the implication does; fails at $L_{\Tru}$
  (Thm.~\ref{thm:br-tin-ent}) and at $U$
  (\textsf{tout\_\allowbreak ent\_\allowbreak fails\_\allowbreak univ});
  at $G$, $\ConR$, $L$: open &
  \textsf{tout\_\allowbreak iff}, \textsf{tout\_\allowbreak qf},
  \textsf{tout\_\allowbreak omega}; \textsf{tout\_\allowbreak fails\_\allowbreak liar},
  \textsf{tout\_\allowbreak fails\_\allowbreak univ} \\
TB &
  $\Tru(\bquo{t}) \liff t$ &
  proved only at grounded-true and grounded-false $t$: at every
  grounded-true $t$ and every $\forall$-free grounded-false $t$; fails
  at every other $t$, whether ungrounded or with a value that is not
  a truth value, such as $\suc\suc\z$; fails at $U$; other
  grounded-false $t$ with quantifiers: open &
  $t \eqv \Tru(\bquo{t})$ as T-In and T-Out: fails at $L_{\Tru}$,
  and at $U$ by its T-Out half &
  \textsf{tschema\_\allowbreak grounded\_\allowbreak true}, \textsf{tschema\_\allowbreak afree\_\allowbreak grounded};
  \textsf{tb\_\allowbreak grounded}, \textsf{tb\_\allowbreak fails\_\allowbreak ungrounded},
  \textsf{tb\_\allowbreak fails\_\allowbreak univ},
  \textsf{tschema\_\allowbreak unprovable\_\allowbreak ungrounded} \\
\multicolumn{5}{@{}l}{\textit{Rules of truth}} \\*
NEC, CONEC &
  $\vdash t \;\Rightarrow\; \vdash \Tru(\bquo{t})$, and back &
  rules on theorems, CONEC for well-formed $t$ only
  (Thm.~\ref{thm:pl-not-sentence}); NEC also for substitution
  instances, in every context (Thm.~\ref{thm:pl-sound}); NEC and
  CONEC as PKF's T-Intro and T-Elim with side formulas, in every
  context, at the images of quoted quantifier-free templates at an
  arbitrary term (Thm.~\ref{thm:pl-rf-tin}); internal NEC
  (one sentence over proof codes): fails at $\z = \tofn 1$
  (Thm.~\ref{thm:pl-no-uniform}) &
  with side formulas, as in PKF, these rules are T-In and T-Out: both fail
  at $L_{\Tru}$ &
  \textsf{iWtrue\_\allowbreak nec}, \textsf{rga\_\allowbreak sound\_\allowbreak sch}, \textsf{iWtrue\_\allowbreak conec};
  \textsf{qfm\_\allowbreak tin}, \textsf{qfm\_\allowbreak tout};
  \textsf{true\_\allowbreak not\_\allowbreak sentence} \\
$\neg$-NEC, $\neg$-CONEC &
  $\vdash \neg t \;\Rightarrow\; \vdash \neg\Tru(\bquo{t})$, and back &
  $\neg$-NEC: rule on theorems at $\forall$-free $t$, conclusion in
  every context (at $\z = \tofn 1$ a theorem outright,
  Thm.~\ref{thm:in-never}); fails with quantifiers, at $\forall x.\,
  x = \tofn 1$ (\textsf{refuted\_\allowbreak unestablished\_\allowbreak fails}); under the
  hypothesis $\ConR$, holds for every refuted $t$
  (Thm.~\ref{thm:pl-refuted}).
  $\neg$-CONEC: rule on theorems at well-formed grounded-true or
  grounded-false $t$; fails at $\Omega$ &
  $\neg$-NEC fails at $L_{\Tru}$ (Thm.~\ref{thm:br-tin-ent});
  $\neg$-CONEC fails at $\Omega$
  (\textsf{negconec\_\allowbreak ent\_\allowbreak fails\_\allowbreak omega}) &
  \textsf{negnec\_\allowbreak qf}; \textsf{refuted\_\allowbreak unestablished\_\allowbreak fails}, \textsf{con\_\allowbreak refuted\_\allowbreak unestablished}; \textsf{neg\_\allowbreak conec\_\allowbreak grounded}; \textsf{neg\_\allowbreak conec\_\allowbreak fails} \\
inter\-sub\-sti\-tu\-tivity &
  $C[t] \eqv C[\Tru(\bquo{t})]$ &
  rule on theorems at well-formed $t$ for $C$ a context $\lor$,
  $\wedge$ or $\neg\neg$; fails under $\neg$ (at $\Omega$) &
  fails: at $L_{\Tru}$ in the empty context $C$ (T-In), under $\neg$
  at $\Omega$ ($\neg$-CONEC) &
  \textsf{isub\_\allowbreak or}, \textsf{isub\_\allowbreak and}, \textsf{isub\_\allowbreak dneg}; \textsf{neg\_\allowbreak conec\_\allowbreak fails} \\
\multicolumn{5}{@{}l}{\textit{Consistency and completeness of the predicate}} \\*
CONS, COMP &
  $\neg(\Tru(x) \wedge \Tru(\mathit{neg}(x)))$;
  $\Tru(x) \lor \Tru(\mathit{neg}(x))$ &
  CONS: never both theorems; in every context${}^\dagger$, as a
  refutation; as a theorem, fails at $U$, and as a refutation the
  direction from $\Tru(\mathit{neg}(x))$ to $\neg\Tru(x)$ fails at
  $U$; internal: fails, plain and with the guards.  COMP: fails on
  theorems (at $L_{\Tru}$) and internally &
  CONS: $\Tru(\mathit{neg}(x)) \vdash \neg\Tru(x)$ holds${}^\dagger$
  and fails at $U$; $\Tru(x) \vdash \neg\Tru(\mathit{neg}(x))$
  holds${}^\dagger$, open beyond $\dagger$.  COMP, from
  $\Wfc(x) \ne \z$ and $\IsCode(x) \ne \z$: fails (at $L_{\Tru}$)
  (\textsf{comp\_\allowbreak ent\_\allowbreak fails}) &
  \textsf{true\_\allowbreak cons\_\allowbreak rule}; \textsf{iWtrue\_\allowbreak refute\_\allowbreak of\_\allowbreak neg\_\allowbreak qf},
  \textsf{iWtrue\_\allowbreak neg\_\allowbreak refute\_\allowbreak qf}; \textsf{cons\_\allowbreak fails\_\allowbreak univ},
  \textsf{cons\_\allowbreak not\_\allowbreak internal}; \textsf{comp\_\allowbreak rule\_\allowbreak fails},
  \textsf{comp\_\allowbreak not\_\allowbreak internal} \\
truths are sentences &
  $\Tru(x) \limp \Wfc(x) \neq \z$ &
  fails as a rule on theorems and internally
  (Thm.~\ref{thm:pl-not-sentence}) &
  fails (\textsf{true\_\allowbreak sent\_\allowbreak ent\_\allowbreak fails}) &
  \textsf{true\_\allowbreak not\_\allowbreak sentence}, \textsf{true\_\allowbreak not\_\allowbreak sentence\_\allowbreak internal} \\
\end{longtable}
}

\paragraph{The cells at the refuted universal.}
\rga refutes $U$ and proves $\Tru(\bquo{\neg U})$ by NEC, but does not
refute $\Tru(\bquo{U})$ (Theorem~\ref{thm:pl-refuted}(4),
\textsf{refuted\_unestablished\_fails}); the six facts about $U$ behind this
paragraph are collected in \textsf{allEq1\_facts} (all the theorems of this paragraph are in
\texttt{RGA\_Liar\_Div/\allowbreak Truth\_Limits.thy}).  Each form that fails at $U$ would
give $\emptyset \vdash \neg\Tru(\bquo{U})$: a rule on theorems
directly, and an entailment with the hypothesis $\Tru(\bquo{U})$ through
\textsf{tru\_reductio}, reductio at the $\Sigma$ hypothesis $\Tru$, or
through \textsf{ir\_nand\_right} ($\neg(a \wedge b), b \vdash \neg a$).
The forward clause for $\forall$ fails at $f = \bquo{\lambda(\vz = \tofn
1)}$ because the instance at $n = 0$ is $\forall$-free and refuted, so
\rga refutes the truth of that instance (\textsf{negnec\_qf}).  These
failure theorems are meta theorems: they cite G\"odel's second theorem, so soundness is in their cone.  A
failure at $U$ is a witness, not a failure at every grounded-false
sentence with quantifiers: a grounded-false existential, such as
$\neg\cA\cdot(\lambda(\vz = \vz))$, is untested.

\paragraph{The internal forms and $\Omega$.}
The internal forms fail at a code built from $L_{\Tru}$, usually
$S_{\Tru}$ (\S\ref{sec:pillars-axioms}); \rga proves both guards
$\Wfc$ and $\IsCode$ at $S_{\Tru}$ (\textsf{guards\_\allowbreak at\_\allowbreak liar}),
and each theorem of the table that ends in \textsf{not\_internal}
states the guarded form as well.  The refutation $\neg\Tru(\bquo{\Omega})$
of Theorem~\ref{thm:pl-omega}(2) has a second proof without \rga's
soundness theorem (\textsf{omega\_\allowbreak unestablished\_\allowbreak sf}), whose
base case uses only the decider's correctness on $\forall$-free
terms.

\paragraph{The affirmation test.}  At a fixed fuel the test $\Aff$ is
decided, and there, several clauses that fail for $\Tru$ are internal
sentences, with a fuel shift of one decider step per connective.  For
every context $\Gamma$, \rga proves
$\forall s.\,\forall c.\ \Aff(\suc\suc s, \mathit{neg}(\mathit{neg}(c)))
\liff \Aff(s, c)$ (\textsf{aff\_\allowbreak dneg\_\allowbreak iff\_\allowbreak int}) and
$\forall s.\,\forall a.\,\forall b.\ \Aff(\suc s, \mathit{or}(a, b))
\liff \Aff(s, a) \lor \Aff(s, b)$ (\textsf{aff\_\allowbreak or\_\allowbreak iff\_\allowbreak int}), CONS for any
two fuels at codes $c$ of $\forall$-free terms (\textsf{aff\_\allowbreak cons\_\allowbreak qf\_\allowbreak int}),
and the introduction half of the definition of $\Tru$, $\forall s.\,\forall
c.\ \Aff(s, c) \limp \Tru(c)$ (\textsf{aff\_\allowbreak true\_\allowbreak int}).  $\Tru(c)$ hides
the fuel under an existential, which \rga does not decide at every
code (\textsf{iWtrue\_\allowbreak undecided}), and there the clauses fail.

\paragraph{The rows marked $\dagger$, and the structural laws.}
The rows marked $\dagger$ for $\neg\lor$, FS's $\neg$ and CONS join
verdicts that the decider reaches at two different fuels.  The proofs of
these rows join the two verdicts by fuel monotonicity of $\Aff$, which \rga proves internally
only at codes of $\forall$-free terms (\textsf{imono\_\allowbreak qf\_\allowbreak le}).  The table also leaves out the predicate's
structural laws, which hold as rules in every context: introduction at
an abstract fuel (Theorem~\ref{thm:in-wtrue-intro}), elimination
(\textsf{rga\_\allowbreak wtrueE}), transport of $\Tru$ and of $\neg\Tru$ along provable
equations of codes (\textsf{iWtrue\_\allowbreak transport}, \textsf{iWtrue\_\allowbreak neg\_\allowbreak transport}),
and the clauses for successor, zero and the symmetry of equations, such
as: from $\Gamma \vdash \Tru(\mathit{eq}(\mathit{succ}(a),
\mathit{succ}(b)))$ derive $\Gamma \vdash \Tru(\mathit{eq}(a, b))$
(\textsf{iWtrue\_\allowbreak sucEqE\_\allowbreak elim} and its companions,
Table~\ref{tab:b-files}).  In every context, \rga also refutes
$\Tru(\bquo{t})$ at every $\forall$-free $t$ that is grounded-false
(\textsf{iWtrue\_\allowbreak refute\_\allowbreak afree}); this is the refutation on which the
$\forall$-free cases of T-Out, TB and $\neg$-NEC rest.

\subsection{Reflection through the truth predicate}
\label{app:pillars-reflect}

This subsection proves the results of \S\ref{sec:pillars-reflect}.
They live in three sessions: \texttt{RGA\_Reflect}, with the theory
\texttt{RGA\_Reflect/\allowbreak Reflect\_QF}; \texttt{RGA\_Poly},
with the six theories \texttt{Poly\_\allowbreak Skeleton}, \texttt{Poly\_\allowbreak Ring},
\texttt{Poly\_\allowbreak Machines}, \texttt{Poly\_\allowbreak Hom}, \texttt{Poly\_\allowbreak Norm} and
\texttt{Poly\_\allowbreak Reflect} in \texttt{RGA\_Poly/}; and \texttt{RGA\_Tin},
which builds on the session \texttt{RGA\_Liar\_Div}, uses the generic
checker of \texttt{RGA\_Checker\_Gen}, imports
\texttt{Poly\_\allowbreak Reflect} from \texttt{RGA\_Poly}, and holds
the two theories \texttt{Tin\_\allowbreak Univ} and
\texttt{Tin\_\allowbreak Univ2} in \texttt{RGA\_Tin/}.  The three sessions
postdate the release named in Section~\ref{sec:mech} (``Artifact''),
and every theory works inside the locale \textsf{rga} with no further
assumption.  \texttt{RGA\_Reflect} builds on the session
\texttt{RGA\_Goedel}, which holds the value wrapper and the value law
below in \texttt{RGA\_Goedel/\allowbreak Goedel\_II}.  In this
subsection $\mathit{ap}$ and $\mathit{eq}$ are the code constructors
of Table~\ref{tab:c-truth-axioms} (Isabelle \textsf{mkAppS},
\textsf{mkEqS}), $\mathit{IsQf}$ is \textsf{iQf}
(\texttt{RGA\_QF/\allowbreak QF\_Mono}).
Table~\ref{tab:c-rf-files} gives the theory file of each name of
\S\ref{sec:pillars-reflect}.

\begin{table}[htbp]
\centering
\footnotesize
\begin{tabular}{@{}>{\raggedright\arraybackslash}p{0.27\linewidth}>{\raggedright\arraybackslash}p{0.70\linewidth}@{}}
\hline
theory file & names \\
\hline
\texttt{RGA\_Goedel/\allowbreak Goedel\_II} & \textsf{iPrDecVal}, \textsf{pr\_\allowbreak dec\_\allowbreak eval}, \textsf{iWtrue\_\allowbreak of\_\allowbreak iPrDecVal}, \textsf{envVar0} \\
\texttt{RGA\_Reflect/\allowbreak Reflect\_QF} & \textsf{ptm}, \textsf{qfm}, \textsf{ptm\_val}, \textsf{qfm\_val},
  \textsf{qfm\_tin}, \textsf{qfm\_tout}, \textsf{ptm\_\allowbreak tout\_\allowbreak eq}, \textsf{rfl\_\allowbreak tout\_\allowbreak binom},
  \textsf{iLift\_\allowbreak zero\_\allowbreak code}, \textsf{iSubM\_\allowbreak updSubM\_\allowbreak hit\_\allowbreak code},
  \textsf{iPrDecVal\_\allowbreak app\_\allowbreak open}, \textsf{iPrDecVal\_\allowbreak app\_\allowbreak img} \\
\texttt{RGA\_Poly/\allowbreak Poly\_Skeleton}, \texttt{Poly\_Ring}, \texttt{Poly\_Machines} & \textsf{recIdx}, \textsf{recT\_unfold}; the ring laws; the six machines (Table~\ref{tab:c-rf-machines}) \\
\texttt{RGA\_Poly/\allowbreak Poly\_Hom} & \textsf{pl\_\allowbreak addT\_\allowbreak gapAdd}, \textsf{ev\_scale}, \textsf{ev\_addP}, \textsf{ev\_mulP}, \textsf{ev\_trim} \\
\texttt{RGA\_Poly/\allowbreak Poly\_Norm} & \textsf{plNormIdx}, \textsf{ptmA}, \textsf{normH}, \textsf{pr2fi\_\allowbreak normIdx\_\allowbreak ptmA},
  \textsf{nmCor\_allI}, \textsf{norm\_\allowbreak correct} \\
\texttt{RGA\_Poly/\allowbreak Poly\_Reflect} & \textsf{polyEqIdx}, \textsf{polyeq\_\allowbreak correct}, \textsf{poly\_\allowbreak reflect\_\allowbreak sentence},
  \textsf{pr2fi\_\allowbreak polyEqIdx}, \textsf{poly\_\allowbreak reflect}, \textsf{rfl\_\allowbreak binom\_\allowbreak by\_\allowbreak reflection},
  \textsf{rfl\_\allowbreak cube\_\allowbreak by\_\allowbreak reflection}, \textsf{rfl\_\allowbreak prod\_\allowbreak by\_\allowbreak reflection} \\
\texttt{RGA\_Tin/\allowbreak Tin\_Univ} & \textsf{openJ}, \textsf{tru\_\allowbreak univ\_\allowbreak of\_\allowbreak cert}, \textsf{tru\_\allowbreak univ\_\allowbreak cert\_\allowbreak ent},
  \textsf{tru\_\allowbreak univ\_\allowbreak cert\_\allowbreak int}, \textsf{tru\_\allowbreak univ\_\allowbreak cert\_\allowbreak of\_\allowbreak dec}, \textsf{tru\_\allowbreak univ\_\allowbreak cert\_\allowbreak conv},
  \textsf{tru\_\allowbreak univ\_\allowbreak cert\_\allowbreak conv\_\allowbreak ent}, \textsf{tru\_\allowbreak univ\_\allowbreak cert\_\allowbreak conv\_\allowbreak int}, \textsf{tru\_\allowbreak univ\_\allowbreak cert\_\allowbreak embed},
  \textsf{all\_\allowbreak omega\_\allowbreak ctx\_\allowbreak rule\_\allowbreak fails}, \textsf{reflCert}, \textsf{tru\_\allowbreak univ\_\allowbreak reflect\_\allowbreak at}, \textsf{tru\_\allowbreak univ\_\allowbreak reflect},
  \textsf{iLiftS\_\allowbreak numC}, \textsf{uinCert}, \textsf{uinCert\_\allowbreak proves}, \textsf{uinIdx}, \textsf{uinM}, \textsf{uinM\_\allowbreak eq},
  \textsf{univ\_\allowbreak inst\_\allowbreak nec\_\allowbreak at}, \textsf{univ\_\allowbreak inst\_\allowbreak nec}, \textsf{univ\_\allowbreak inst\_\allowbreak nec\_\allowbreak 1}, \textsf{univ\_\allowbreak inst\_\allowbreak code\_\allowbreak at} \\
\texttt{RGA\_Tin/\allowbreak Tin\_Univ2} (103 names) & \textsf{tin\_\allowbreak hins\_\allowbreak single}, \textsf{tin\_\allowbreak fsetcode\_\allowbreak single}, \textsf{tin\_\allowbreak cons\_\allowbreak valid}, \textsf{mpCert}, \textsf{mpCert\_\allowbreak chk},
  \textsf{tinEnv2}, \textsf{tinImg2}, \textsf{tinCa2}, \textsf{inst2Cert}, \textsf{inst2\_\allowbreak chk}, \textsf{openTop}, \textsf{openTop\_\allowbreak proves},
  \textsf{univ2Cert}, \textsf{univ2Cert\_\allowbreak proves}, \textsf{univ2Idx}, \textsf{univ2M}, \textsf{univ2M\_\allowbreak eq}, \textsf{univ2\_\allowbreak nec\_\allowbreak at}, \textsf{univ2\_\allowbreak nec},
  \textsf{polyG}, \textsf{polyG\_\allowbreak thm}, \textsf{poly\_\allowbreak univ2\_\allowbreak nec}, \textsf{guardCert}, \textsf{guardCert\_\allowbreak chk}, \textsf{dischCert}, \textsf{dischCert\_\allowbreak proves},
  \textsf{dischIdx}, \textsf{dischM}, \textsf{dischM\_\allowbreak eq}, \textsf{univ2\_\allowbreak disch\_\allowbreak at}, \textsf{univ2\_\allowbreak disch}, \textsf{tinPdIdx}, \textsf{tinPd\_\allowbreak D}, \textsf{tinPd\_\allowbreak I},
  \textsf{polyQ}, \textsf{polyQ\_\allowbreak conv}, \textsf{polyD\_\allowbreak thm}, \textsf{poly\_\allowbreak disch}, \textsf{poly\_\allowbreak all\_\allowbreak sentence}, \textsf{tru\_\allowbreak univ\_\allowbreak by\_\allowbreak reflection},
  \textsf{rfl\_\allowbreak binom\_\allowbreak tru\_\allowbreak univ}, \textsf{rfl\_\allowbreak cube\_\allowbreak tru\_\allowbreak univ}, \textsf{rfl\_\allowbreak prod\_\allowbreak tru\_\allowbreak univ} \\
\texttt{RGA\_Goedel/\allowbreak Entailment\_Goedel} & \textsf{true\_\allowbreak reflect\_\allowbreak chk\_\allowbreak ent\_\allowbreak fails} \\
\hline
\end{tabular}
\caption{The theory files of the names of \S\ref{sec:pillars-reflect}
  and of this subsection's headline facts.}
\label{tab:c-rf-files}
\end{table}

\paragraph{The value wrapper and the environment.}
The value wrapper is $\mathit{DecVal}(c, x) \ldef \neg\cA\cdot(\lambda\,
\neg(\Dec(\vz, \lft c) = \suc\lft x))$ (\textsf{iPrDecVal}), that is,
$\exists l.\ \Dec(l, c) = \suc x$.  Value one is affirmation: if
$\Gamma \vdash c\jnat$ and $\Gamma \vdash \mathit{DecVal}(c, \tofn 1)$,
then $\Gamma \vdash \Tru(c)$ (\textsf{iWtrue\_of\_iPrDecVal}).  The
environment $\mathit{env}(v)$ (\textsf{envVar0}) updates the empty
environment $\z$ at the variable $\vz$ with the compiled code of the
numeral of $v$ (\textsf{numC}).  If $\Gamma \vdash v\jnat$, then $\Gamma$
derives $\mathit{env}(v)\jnat$, the two environment checks
$\mathit{sigWfChk}\cdot\mathit{env}(v)$ and
$\mathit{envQf}\cdot\mathit{env}(v)$, which say that the entries are
codes of well-formed and of $\forall$-free terms
(\textsf{envVar0\_sigWfChk}, \textsf{envVar0\_envQf}), and the read
$\mathit{DecVal}(\Subst(\mathit{env}(v), \bquo{\vz}), v)$
(\textsf{envVar0\_var}).

\begin{thm}[\textsf{pr\_dec\_eval}]
\label{thm:c-rf-eval}
Let $f$ be an index and $T_1$, $T_2$ well-formed $\forall$-free terms.
Let $\Gamma$ derive $z\jnat$, $\mathit{sigWfChk}\cdot z$,
$\mathit{envQf}\cdot z$, $x\jnat$ and $y\jnat$.  If $\Gamma \vdash
\mathit{DecVal}(\Subst(z, \bquo{T_1}), x)$ and $\Gamma \vdash
\mathit{DecVal}(\Subst(z, \bquo{T_2}), y)$, then
$\Gamma \vdash \mathit{DecVal}(\Subst(z, \bquo{\cmp(f)\cdot T_1\cdot
T_2}), \cmp(f)\cdot x\cdot y)$.
\end{thm}

\noindent
\emph{Proof.}  The statement is the HOL predicate
\textsf{prDecEvalOK}~$f$, a universal over $\Gamma$, $z$, $T_1$,
$T_2$, $x$ and $y$.  The proof is a HOL induction on the index tree
$f$, with one case for the initial functions
(\textsf{prDecEvalOK\_init}), one for composition
(\textsf{prDecEvalOK\_sub}) and one for recursion
(\textsf{prDecEvalOK\_rec}), which carries one number induction, on
the numeral of the recursion argument; so the derivation for $f$
carries one such induction for each recursion node of $f$.  The
form \textsf{pr\_dec\_eval\_D} unpacks the predicate into a rule.  \qed

\paragraph{The templates.}
The classes of Definition~\ref{def:pl-rf-tpl} are the inductive
predicates \textsf{ptm} (constructors \textsf{pt\_var} for $\vz$,
\textsf{pt\_num} for numerals, \textsf{pt\_suc}, and \textsf{pt\_app}
for $\cmp(f)\cdot a\cdot b$ at any index $f$) and \textsf{qfm}
(\textsf{qf\_eq}, \textsf{qf\_neg}, \textsf{qf\_or}), with conjunction
derived (\textsf{qfm\_and}).  Every member of either class is
well-formed and $\forall$-free (\textsf{ptm\_wft\_afree},
\textsf{qfm\_wft\_afree}), and the member's only free variable is $\vz$, since
numerals and the compiled terms $\cmp(f)$ are closed.

\emph{Proof of Theorem~\ref{thm:pl-rf-val}}
(\textsf{ptm\_val}, \textsf{qfm\_val}).  Part~1 is a HOL induction over
the term template.  The variable case is the read
\textsf{envVar0\_var}; the numeral case is \textsf{iPrDecVal\_img\_num};
the successor case is \textsf{iPrDecVal\_suc} with the code identity
\textsf{iSubM\_code\_suc}; the application case is
Theorem~\ref{thm:c-rf-eval} at $z = \mathit{env}(v)$, whose
environment checks the two lemmas above supply.  Part~2 is one
simultaneous HOL induction over the formula template, of the
conjunction of the three claims (\textsf{qfm\_vals}, unfolded as
\textsf{qfm\_val}).  At an equation $a = b$, decidedness comes from
\textsf{ir\_eqTI} and \textsf{ir\_boolE}, and the two values from
part~1 at $a$ and $b$ through the decider's laws for equations
\textsf{iPrDecVal\_eq\_true} and \textsf{iPrDecVal\_eq\_false}, with the
$\forall$-free scan of the image (\textsf{iQfScan\_iSubM\_afree}).
Negation uses \textsf{iPrDecVal\_negF} and \textsf{iPrDecVal\_negT}.
Disjunction uses three value laws proved in \texttt{Reflect\_QF},
\textsf{iPrDecVal\_orT1}, \textsf{iPrDecVal\_orT2} and
\textsf{iPrDecVal\_orF}, with the code identity \textsf{iSubM\_code\_or}.
The two implications of part~2 are HOL implications between
derivability judgments.  \qed

\emph{Proof of Theorem~\ref{thm:pl-rf-tin}}
(\textsf{qfm\_tin}, \textsf{qfm\_tout}).  T-In: part~2 of
Theorem~\ref{thm:pl-rf-val} gives the value $\tofn 1$ at the image $X
\ldef \Subst(\mathit{env}(v), \bquo{p})$, and
\textsf{iWtrue\_of\_iPrDecVal} gives $\Tru(X)$.  T-Out: eliminate the
decided disjunction $\ssub{\vz}{v}{p} \lor \neg\ssub{\vz}{v}{p}$
(\textsf{ir\_disjE1}).  The first case is the goal.  In the second
case, with $\neg\ssub{\vz}{v}{p}$ as a hypothesis, the value at $X$ is
$\z$, \textsf{iPrDecVal\_negF} gives the value $\tofn 1$ at
$\mathit{neg}(X)$, and \textsf{iWtrue\_of\_iPrDecVal} gives
$\Tru(\mathit{neg}(X))$.  The scan gives $\qfs(X) = \z$.  The lemma
\textsf{iWtrue\_collision\_qf} (\texttt{RGA\_QF/\allowbreak QF\_Arms})
takes $X\jnat$, $\Tru(X)$, $\Tru(\mathit{neg}(X))$, $\qfs(X) = \z$ and a
well-formed $G$, and concludes $G$; here $G$ is $\ssub{\vz}{v}{p}$.
\qed

\begin{thm}[\textsf{ptm\_tout\_eq}, \textsf{rfl\_tout\_binom}]
\label{thm:c-rf-eq}
Let $\Gamma \vdash v\jnat$.
\begin{enumerate}
\item For term templates $a$ and $b$, if $\Gamma \vdash
  \Tru(\Subst(\mathit{env}(v), \bquo{a = b}))$, then $\Gamma \vdash
  \ssub{\vz}{v}{a} = \ssub{\vz}{v}{b}$.
\item Let $L(u) \ldef \suc u \times \suc u$ and $R(u) \ldef (u \times u
  + \tofn 2 \times u) + \tofn 1$, with the compiled operations
  (\textsf{rflL}, \textsf{rflR}).  If $\Gamma \vdash
  \Tru(\Subst(\mathit{env}(v), \bquo{L(\vz) = R(\vz)}))$, then $\Gamma
  \vdash L(v) = R(v)$.
\end{enumerate}
\end{thm}

\noindent
\emph{Proof.}  Part~1 splits on the decided equation
(\textsf{ir\_eq\_cases}).  In the false case,
\textsf{iPrDecVal\_eq\_false} and part~1 of Theorem~\ref{thm:pl-rf-val}
give the value $\z$ at the image, transported by
\textsf{iPrDecVal\_code\_eq}; then \textsf{iPrDecVal\_negF},
\textsf{iWtrue\_of\_iPrDecVal} and \textsf{iWtrue\_collision\_qf}
finish as in T-Out.  Part~2 is part~1 at the two sides.  The twin
\textsf{rfl\_tout\_binom\_generic} has the same statement and goes
through Theorem~\ref{thm:pl-rf-tin} instead.  \qed

\begin{thm}[\textsf{iLift\_zero\_code}, \textsf{iSubM\_updSubM\_hit\_code}]
\label{thm:c-rf-lift}
Let $\Gamma$ derive $c\jnat$, $\Wfc(c) \neq \z$ and $\IsCode(c) \neq
\z$.
\begin{enumerate}
\item The compiled de Bruijn lift by zero (Isabelle \textsf{iLift})
  maps $c$ to $c$: \rga derives that the lift of $c$ equals $c$.
\item If $\Gamma \vdash z\jnat$, then substituting the environment $z$
  updated at the variable $v_i$ with $c$ into the code $\bquo{v_i}$
  yields $c$.
\end{enumerate}
\end{thm}

\noindent
\emph{Proof.}  Part~1 composes two lemmas of the coding kit,
\textsf{lift\_bridge\_at\_zero} (\texttt{RGA\_Checker/\allowbreak
Chk\_Sigma\_Machinery}) and \textsf{lift0\_depth\_zero\_id}
(\texttt{RGA\_Checker/\allowbreak Chk\_Emission}), the second an
internal course-of-values induction over the codes of well-formed
terms.  Part~2 is the general update read
\textsf{iSubM\_updSubM\_hit\_any}, whose result is the lift of $c$ by
zero, followed by part~1.  At the code of a compiled constant, which
has bound variables, part~1 holds outright
(\textsf{iLift\_zero\_cmp}).  \qed

\emph{Proof of Theorem~\ref{thm:pl-rf-app}}
(\textsf{iPrDecVal\_app\_open}, \textsf{iPrDecVal\_app\_img}).
Part~1: let $z$ be the two-entry environment that updates $\z$ at
$\vz$ with $c_1$ and then at $v_1$ with $c_2$.  The guards on $c_1$
and $c_2$ give the environment checks of $z$.  By
Theorem~\ref{thm:c-rf-lift}(2) and the law that an update at a
variable that does not occur leaves a code's image unchanged
(\textsf{iSubM\_updSubM\_code\_nfree}), the images of $\bquo{\vz}$ and
$\bquo{v_1}$ under $z$ are $c_1$ and $c_2$, and the image of
$\bquo{\cmp(f)}$ is $\bquo{\cmp(f)}$ itself, by two such steps and
\textsf{iSubM\_zero\_id}.  Theorem~\ref{thm:c-rf-eval} at the
templates $\vz$ and $v_1$ gives the value at the image of
$\cmp(f)\cdot\vz\cdot v_1$, and the code identity for applications
(\textsf{iSubM\_code\_app}) rewrites that image to
$\mathit{ap}(\mathit{ap}(\bquo{\cmp(f)}, c_1), c_2)$.  Part~2: the
image under $\mathit{env}(v)$ of the application code is the
application of the images (\textsf{iSubM\_mkAppS}), and the guards
pass from $c_i$ to $\Subst(\mathit{env}(v), c_i)$
(\textsf{iSubM\_imgpres}, \textsf{iSubM\_gnpres},
\textsf{iSubM\_qfpres}); part~1 at those images finishes.  The
environment shape is fixed to explicit updates of $\z$; no law for
$\Subst(z, \bquo{\cmp(f)})$ at an arbitrary environment $z$ is used.
\qed

\paragraph{The skeleton and the ring laws.}
The skeleton \textsf{recIdx}~$B$~$S_1$~$S_2$~$H$
(\texttt{Poly\_Skeleton}) computes, in HOL, the function $F$ with
$F(m, 0) = \den{B}(m, 0)$ and, for $n \neq 0$,
$F(m, n) = \den{H}(\langle m, n\rangle, \langle F(m, s_1), F(m,
s_2)\rangle)$, where $s_i = \den{S_i}(m, n)$ when that value is below
$n$, and $n - 1$ otherwise (\textsf{recF}, \textsf{pr2fi\_recIdx}).

\begin{thm}[\textsf{recT\_unfold}]
\label{thm:c-rf-skel}
Let $\mathit{Rec}(m, n)$ be the compiled value of the skeleton (Isabelle
\textsf{recV}).  If $\Gamma$ derives $m\jnat$, $n\jnat$, $n \neq \z$,
and that $\cmp(S_1)\cdot m\cdot n$ and $\cmp(S_2)\cdot m\cdot n$ are
below $n$, then
$\Gamma \vdash \mathit{Rec}(m, n) = \cmp(H)\cdot\langle m, n\rangle
\cdot\langle \mathit{Rec}(m, \cmp(S_1)\cdot m\cdot n), \mathit{Rec}(m,
\cmp(S_2)\cdot m\cdot n)\rangle$, with compiled pairs.
\end{thm}

\noindent
The proof composes the step law without the bounds
(\textsf{recT\_step}) with the lookup law (\textsf{recT\_look}), which
rest on the coding kit's laws for course-of-values recursion
(\textsf{ir\_covVal\_covStep\_suc}, \textsf{ir\_covR\_lookup}).  The
list instance \textsf{lrecIdx} has the laws \textsf{recT\_lst\_nil}
and \textsf{recT\_lst\_cons}.  The skeleton is explicit, not of
feasible size, in the sense of \S\ref{sec:inside-cmp}.

In the ring laws and in Theorem~\ref{thm:c-rf-hom}, $+$ is the
coding kit's addition by recursion (\textsf{gapAdd},
\texttt{RGA/\allowbreak RGA\_Cantor\_Internal}), and $\times$ is the
compiled multiplication.  \rga proves the kit's addition equal to the
compiled addition at numbers (\textsf{pl\_addT\_gapAdd}).
The ring laws (\texttt{Poly\_Ring}) hold at open terms $x$, $k$
with $\Gamma \vdash x\jnat$ and $\Gamma \vdash k\jnat$, in every
context: $\z \times k = \z$ and $\suc x \times k = x \times k + k$
(\textsf{pl\_mul\_zeroL}, \textsf{pl\_mul\_sucL}), $x \times \tofn 1
= x$ and $\tofn 1 \times x = x$ (\textsf{pl\_mul\_oneR},
\textsf{pl\_mul\_oneL}), commutativity (\textsf{pl\_mul\_comm}), the two
distributive laws (\textsf{pl\_mul\_distL}, \textsf{pl\_mul\_distR}),
associativity and left commutativity of multiplication
(\textsf{pl\_mul\_assoc}, \textsf{pl\_mul\_lcomm}), and two
rearrangements of sums (\textsf{pl\_add\_lcomm},
\textsf{pl\_add\_swap4}).  Each induction is one application of
\textsf{pl\_eq\_ind}, a number induction with an equation as its
motive, over the kit's commutativity and associativity of addition.

\paragraph{The machines and the homomorphisms.}
Table~\ref{tab:c-rf-machines} lists the six machines
(\texttt{Poly\_Machines}).  A coefficient list is coded in the kit's
list coding, $[\,] = 0$ and $h :: t = \suc\langle h, t\rangle$.  Each
machine's HOL semantics is a HOL list function, and \rga derives its
unfolding laws at open codes through Theorem~\ref{thm:c-rf-skel} or
\textsf{recT\_lst\_cons}.

\begin{table}[htbp]
\centering
\footnotesize
\begin{tabular}{@{}>{\raggedright\arraybackslash}p{0.18\linewidth}>{\raggedright\arraybackslash}p{0.17\linewidth}>{\raggedright\arraybackslash}p{0.13\linewidth}>{\raggedright\arraybackslash}p{0.20\linewidth}>{\raggedright\arraybackslash}p{0.20\linewidth}@{}}
\hline
machine & index & HOL mirror & semantics & internal law \\
\hline
Horner evaluation $\mathit{ev}(v, l)$ & \textsf{evIdx} & \textsf{hornerH}
  & \textsf{pr2fi\_evIdx} & \textsf{ev\_nil}, \textsf{ev\_cons} \\
scaling & \textsf{plScIdx} & \textsf{scaleH} & \textsf{pr2fi\_plScIdx}
  & \textsf{pl\_sc\_cons} \\
addition, with padding & \textsf{plAddIdx} & \textsf{addH} &
  \textsf{pr2fi\_plAddIdx} & \textsf{pl\_add\_cons} \\
multiplication, shift and add & \textsf{plMulIdx} & \textsf{mulH} &
  \textsf{pr2fi\_plMulIdx} & \textsf{pl\_mul\_cons} \\
trimming of trailing zeros & \textsf{plTrimIdx} & \textsf{trimH} &
  \textsf{pr2fi\_plTrimIdx} & \textsf{pl\_trim\_cons} \\
equality after trimming & \textsf{plEqIdx} & --- &
  \textsf{pr2fi\_plEqIdx} & \textsf{pl\_eqP\_law} \\
\hline
\end{tabular}
\caption{The six machines of the normalizer
  (\texttt{RGA\_Poly/\allowbreak Poly\_Machines}): explicit index trees,
  their HOL mirrors and semantics lemmas, and their unfolding laws at
  open codes.  The equality machine's semantics is the test
  $\mathit{trimH}(\mathit{xs}) = \mathit{trimH}(\mathit{ys})$, and its
  law says that \rga derives, at open codes $p$ and $q$, that the
  machine's value is the compiled 0/1 equality test of the two trimmed
  lists.}
\label{tab:c-rf-machines}
\end{table}

\begin{thm}[\textsf{ev\_addP}, \textsf{ev\_mulP}]
\label{thm:c-rf-hom}
Let $\Gamma$ derive $v\jnat$, $p\jnat$ and $q\jnat$.  Then $\Gamma
\vdash \mathit{ev}(v, p \oplus q) = \mathit{ev}(v, p) +
\mathit{ev}(v, q)$ and $\Gamma \vdash \mathit{ev}(v, p \otimes q) =
\mathit{ev}(v, p) \times \mathit{ev}(v, q)$, where $\oplus$ and
$\otimes$ are the addition and multiplication machines.  Likewise, for
the scaling machine $\mathit{sc}$ and the trimming machine
$\mathit{trim}$, $\Gamma \vdash \mathit{ev}(v, \mathit{sc}(c, q)) = c
\times \mathit{ev}(v, q)$ when $\Gamma \vdash c\jnat$
(\textsf{ev\_scale}), and $\Gamma \vdash \mathit{ev}(v,
\mathit{trim}(p)) = \mathit{ev}(v, p)$ (\textsf{ev\_trim}).
\end{thm}

\noindent
\emph{Proof.}  Each law is one internal course-of-values induction by
\textsf{pl\_cov\_eq\_ind} (\texttt{Poly\_Hom}), an equation as motive
over the kit's \textsf{cov\_pair\_ind}
(\texttt{RGA\_Kit/\allowbreak Cov\_Induct}), with the pair of the
parameter and the code as carrier.  The cases come from
\textsf{pl\_lst\_cases}, the algebra from the ring laws, and the case
of trimming also uses \textsf{pl\_add\_eq0D}: a sum is $\z$ only if both
summands are.  \qed

\paragraph{The normalizer.}
The normalizer is $\mathit{Norm}(x) \ldef \cmp(\textsf{plNormIdx})\cdot
\z\cdot x$ (\texttt{Poly\_Norm}), with \textsf{plNormIdx} the head of a
skeleton instance whose step \textsf{nmHIdx} is a chain of four
conditional tests and a default.  The tests are, in order,
\textsf{nmIsZ} (the code of $\z$), \textsf{nmIsS} (a successor),
\textsf{nmIsH} at the code of the compiled addition, and \textsf{nmIsH}
at the code of the compiled multiplication; the default
\textsf{nmVRIdx} rejects with the pair $\langle 3, \z\rangle$.  Each test
compares the code with the constructor applied to the code's decoded
parts.  The first lookup reads the left argument of a binary
application, and the second lookup reads the argument of a successor or
the right argument (\textsf{nmS1Idx}, \textsf{nmS2Idx}).  The polynomial
templates are the inductive class \textsf{ptmA} (\textsf{pa\_var},
\textsf{pa\_zero}, \textsf{pa\_suc}, \textsf{pa\_add}, \textsf{pa\_mul}),
which contains the numerals as towers of successors (\textsf{ptmA\_num})
and lies inside the term templates (\textsf{ptmA\_ptm}).  The HOL normal
form \textsf{normH} sends $\vz$ to $[0, 1]$, $\z$ to $[\,]$, $\suc t$ to
the sum of $t$'s normal form with $[1]$, and the compiled sum and
product to the sum and product of the parts' normal forms.

\begin{thm}[\textsf{pr2fi\_normIdx\_ptmA}, \textsf{pr2fi\_polyEqIdx}]
\label{thm:c-rf-hol}
For polynomial templates $t$, $a$ and $b$:
\begin{enumerate}
\item $\den{\textsf{plNormIdx}}(0, \quo{t}) = \langle 0,
  \quo{\mathit{normH}(t)}\rangle$;
\item $\den{\textsf{polyEqIdx}}(\quo{a}, \quo{b})$ is $1$ if
  $\mathit{trimH}(\mathit{normH}(a)) = \mathit{trimH}(\mathit{normH}(b))$,
  and $0$ otherwise.
\end{enumerate}
\end{thm}

\noindent
Both parts are HOL computations about the index trees.  Part~2 follows
from part~1 and the semantics of the equality machine.

\emph{Proof of Theorem~\ref{thm:pl-rf-norm}}
(\textsf{norm\_correct}, \textsf{nmCor\_allI}).  The motive is
$\mathit{nmCorMot} = \lambda\,C(\vz)$ with
\[
  \begin{array}{@{}l@{}}
  C(u) \ldef \Wfc(x_u) \neq \z \limp \IsCode(x_u) \neq \z \limp
  \mathit{IsQf}(x_u) \neq \z \limp \mathit{fst}(\mathit{Norm}(x_u)) =
  \z \limp {} \\
  \qquad \mathit{DecVal}(\Subst(\mathit{env}(v_u), x_u),
  \mathit{ev}(v_u, \mathit{snd}(\mathit{Norm}(x_u))))
  \end{array}
\]
for $v_u = \mathit{fst}(u)$ and $x_u = \mathit{snd}(u)$, with \rga's
implication (\textsf{nmCorBody}).  The packed form
\textsf{nm\_cor\_packed} is one application of \textsf{cov\_pair\_ind}
with this motive, over the pair of $v$ and the code, and
\textsf{nmCor\_allI} closes the induction as $\emptyset \vdash \cA\cdot
\mathit{nmCorMot}$.  \textsf{norm\_correct} instantiates this sentence
at the pair of $v$ and $x$ (\textsf{ir\_allE}) and discharges the four
decided antecedents (\textsf{nmCor\_use}).  The induction step
\textsf{nm\_cor\_node} dispatches six cases by five decided splits: the
code $0$ of $\vz$ (\textsf{nm\_cor\_zero}), the code of $\z$
(\textsf{nm\_cor\_isZ}), a successor
(\textsf{nm\_cor\_isS}), the two applications by one generic lemma
(\textsf{nm\_cor\_app}, which uses Theorem~\ref{thm:pl-rf-app}(2) and
Theorem~\ref{thm:c-rf-hom}), and the default (\textsf{nm\_cor\_rej}),
where the tag would be $3$ while the guard says $\z$, so
\textsf{ir\_neE} gives any well-formed conclusion.  The kids' tags are
summed, and a zero sum splits by \textsf{pl\_add\_eq0D}.  The unfolding
with one bounded lookup is \textsf{recT\_unfold2}.  \qed

\emph{Proof of Theorem~\ref{thm:pl-rf-sentence}}
(\textsf{poly\_reflect\_sentence}).  The
test is $\mathit{polyEq}(x, y) \ldef \cmp(\textsf{polyEqIdx})\cdot
x\cdot y$, where \textsf{polyEqIdx} conjoins the two tag tests with the
equality machine at the two normal forms.  First,
\textsf{polyeq\_correct}: let $\Gamma$ derive $v\jnat$, $x\jnat$,
$y\jnat$, the three code guards at $x$ and at $y$, and
$\mathit{polyEq}(x, y) \neq \z$.  Then $\Gamma \vdash
\mathit{DecVal}(\Subst(\mathit{env}(v), \mathit{eq}(x, y)), \tofn
1)$.  Its proof splits the acceptance into the two zero tags and the
agreement of the trimmed lists (\textsf{polyEq\_ne0D}), applies
Theorem~\ref{thm:pl-rf-norm} at $x$ and at $y$, equates the two Horner
values by \textsf{ev\_trim}, and finishes with
\textsf{iPrDecVal\_eq\_true}.  The sentence (\textsf{polyReflS})
nests three universals over the chain of implications
\textsf{plRB}~$x$~$y$~$v$ (\textsf{impCh}, \texttt{RGA\_Kit/\allowbreak
Cov\_Induct}).  Its proof is three universal introductions
(\textsf{lam\_allI\_wk}) and \textsf{impCh\_intro}, whose
introductions go through \textsf{ir\_impI\_Bw} to \textsf{ir\_impI\_bool},
implication introduction at a hypothesis $a$ with $\Gamma \vdash
a\jbool$.  Each of the seven antecedents is decided
because the tested term is a total compiled term (\textsf{ne0\_B}), and
the body is \textsf{iWtrue\_of\_iPrDecVal} after
\textsf{polyeq\_correct}.  \qed

\emph{Proof of Theorem~\ref{thm:pl-rf-reflect}} (\textsf{poly\_reflect}).
The premises are \textsf{ptmA}~$a$, \textsf{ptmA}~$b$, the HOL equation
$\mathit{trimH}(\mathit{normH}(a)) = \mathit{trimH}(\mathit{normH}(b))$
and $\Gamma \vdash v\jnat$.  The steps are:
\begin{enumerate}
\item the closed evaluation $\emptyset \vdash \mathit{polyEq}(\bquo{a},
  \bquo{b}) = \tofn 1$, from Theorem~\ref{thm:in-cmp-num}
  (\textsf{cmp\_eq\_num}), rewritten by
  Theorem~\ref{thm:c-rf-hol}(2) and the HOL premise;
  \textsf{cmp\_eq\_num} rests on \textsf{ir\_osr\_eq} and through it on
  open completeness (\textsf{rga\_complete\_open});
\item the sentence of Theorem~\ref{thm:pl-rf-sentence}, weakened to
  $\Gamma$ (\textsf{ir\_weaken0}) and instantiated at $\bquo{a}$,
  $\bquo{b}$ and $v$ by three applications of \textsf{lam\_allE};
\item the guards at the code numerals (\textsf{pl\_code\_guards}): the
  well-formedness test by \textsf{iWft\_code\_ne0}, the code test by
  \textsf{iIsCode\_num\_ne0}, and the $\forall$-free test by a closed evaluation of the
  scan through completeness (\textsf{ir\_qfScan\_qf}), carried to the
  guard by \textsf{iQfScan\_iSubM\_afree} and \textsf{iSubM\_zero\_id};
\item the discharge of the seven antecedents (\textsf{impCh\_use}),
  which gives $\Tru(\Subst(\mathit{env}(v), \mathit{eq}(\bquo{a},
  \bquo{b})))$;
\item transport along $\mathit{eq}(\bquo{a}, \bquo{b}) = \bquo{a = b}$
  (\textsf{mkEqS\_num\_code}, \textsf{iWtrue\_transport});
\item Theorem~\ref{thm:c-rf-eq}(1), through \textsf{ptmA\_ptm}.  \qed
\end{enumerate}

\begin{thm}[\textsf{rfl\_cube\_by\_reflection}, \textsf{rfl\_prod\_by\_reflection}]
\label{thm:c-rf-more}
For every context $\Gamma$ and every $v$ with $\Gamma \vdash v\jnat$,
with the compiled operations,
\begin{enumerate}
\item $\Gamma \vdash \suc v \times (\suc v \times \suc v) = ((v
  \times (v \times v) + \tofn 3 \times (v \times v)) + \tofn 3 \times
  v) + \tofn 1$;
\item $\Gamma \vdash \suc\suc v \times \suc\suc\suc v = (v \times v +
  \tofn 5 \times v) + \tofn 6$.
\end{enumerate}
\end{thm}

\noindent
\emph{Proof.}  As for Corollary~\ref{cor:pl-rf-binom}
(\textsf{rfl\_binom\_by\_reflection}): membership in \textsf{ptmA} by
its introduction rules, the agreement of the normal forms by one
\textsf{simp} evaluation, and Theorem~\ref{thm:pl-rf-reflect}.  \qed

\paragraph{T-In at a universal.}
The results of \S\ref{sec:pillars-reflect} on T-In at a universal are
in \texttt{RGA\_Tin/\allowbreak Tin\_Univ}.  The open-instance code
is $J(c) \ldef \mathit{pair}(\bquo{\{\vz\jnat\}},
\mathit{ap}(\mathit{lift}(c), \bquo{\vz}))$ (\textsf{openJ}), where
$\mathit{pair}$ is the compiled pairing of a judgment's two parts and
$\mathit{lift}$ is the compiled de Bruijn lift by one at cutoff $0$
(Isabelle \textsf{iLiftS}); $\Gamma \vdash c\jnat$ gives $\Gamma
\vdash J(c)\jnat$ (\textsf{openJ\_N}).

\begin{thm}[\textsf{tru\_\allowbreak univ\_\allowbreak of\_\allowbreak cert}, \textsf{tru\_\allowbreak univ\_\allowbreak cert\_\allowbreak ent}]
\label{thm:c-tin-cert}
\leavevmode
\begin{enumerate}
\item If $\Gamma \vdash c\jnat$, $\Gamma \vdash p\jnat$ and $\Gamma
  \vdash \cmp(\pi)\cdot J(c)\cdot p \neq \z$, then $\Gamma \vdash
  \Tru(\mathit{ap}(\bquo{\cA}, c))$.
\item $\Cert(J(\vz)),\ \vz\jnat \vdash \Tru(\mathit{ap}(\bquo{\cA},
  \vz))$.
\end{enumerate}
\end{thm}

\noindent
\emph{Proof.}  Part~1: the head data of the verdict
(\textsf{iProvesChk\_headD}) give the certificate search's hit for
$c$ at the level $\suc p$ (\textsf{ir\_certTrue\_intro}), the hit gives
the affirmation $\Dec(\suc\suc p, \mathit{ap}(\bquo{\cA}, c)) = \tofn
2$ (\textsf{ir\_A\_cert\_affirm}), and the affirmation gives truth
(\textsf{iWtrueI}).  Part~2 eliminates the certificate predicate
(\textsf{iHasCertE}) to a witness tree $P$, rebuilds the verdict at $P$
(\textsf{iProvesChk\_headI}) and applies part~1.  The first sentence
of Theorem~\ref{thm:pl-tin-cert} (\textsf{tru\_univ\_cert\_int}) is
part~1 under two universal introductions (\textsf{ir\_allI\_open})
and the deduction theorem at the decided antecedent
(\textsf{ir\_impI\_bool}).  \qed

\begin{thm}[\textsf{tru\_\allowbreak univ\_\allowbreak cert\_\allowbreak of\_\allowbreak dec}, \textsf{tru\_\allowbreak univ\_\allowbreak cert\_\allowbreak conv}]
\label{thm:c-tin-conv}
Let $\Gamma \vdash c\jnat$.
\begin{enumerate}
\item If $\Gamma \vdash s\jnat$ and $\Gamma \vdash \Dec(s,
  \mathit{ap}(\bquo{\cA}, c)) = \tofn 2$, then $\Gamma \vdash
  \Cert(J(c))$.
\item If $\Gamma \vdash \Tru(\mathit{ap}(\bquo{\cA}, c))$, then
  $\Gamma \vdash \Cert(J(c))$; as an entailment,
  $\Tru(\mathit{ap}(\bquo{\cA}, \vz)),\ \vz\jnat \vdash \Cert(J(\vz))$
  (\textsf{tru\_univ\_cert\_conv\_ent}).
\end{enumerate}
\end{thm}

\noindent
\emph{Proof.}  Part~1 splits on the fuel.  At fuel $0$ the decider's
value is $0$ (\textsf{dec\_zero\_fuel}).  At fuel $1$ the application
arm reads the head cell of the base table, which is zero, so the value
is $0$ (\textsf{decA\_one\_fuel}, through \textsf{ir\_tabL\_zero\_cells}).
At a double successor the arm's value is the verdict of the
certificate searches (\textsf{ir\_stepV\_app\_cert\_miss},
\textsf{ir\_stepV\_app\_cert\_hitF}); only a hit of the search for a
true universal gives $\tofn 2$, and that hit is a certificate of
$J(c)$ (\textsf{certTrue\_hasCert}; together \textsf{decA\_ss\_cert}).
Part~2 eliminates $\Tru$ to a fuel at which the decider affirms
(\textsf{rga\_wtrueE}) and applies part~1.  The second sentence of
Theorem~\ref{thm:pl-tin-cert} (\textsf{tru\_univ\_cert\_conv\_int})
is part~1 under two universal introductions and the deduction theorem
at the decided antecedent, whose value $\tofn 1$ gives $\Dec(s,
\mathit{ap}(\bquo{\cA}, c)) = \tofn 2$
(\textsf{iPrDecTrueAt\_true\_elim}).  \qed

\noindent
The certificate side does not depend on \rga's checker table: for
every list $L$ of rule checkers that contains \rga's
(\textsf{drvChkList}), $\Gamma \vdash \Cert(J(c))$ gives the
certificate at the checker over $L$ (\textsf{tru\_univ\_cert\_embed},
by \textsf{chkG\_embed\_cert}).  The truth predicate is \rga's.

\begin{thm}[\textsf{all\_\allowbreak omega\_\allowbreak ctx\_\allowbreak rule\_\allowbreak fails}]
\label{thm:c-tin-omega}
Let $u_L \ldef \lambda\neg\Aff(\vz, S_{\Tru})$.  It is not the case
that for every context $\Gamma$, $\Gamma \vdash \forall x.\,
\Tru(\mathit{ap}(\bquo{u_L}, \Num\cdot x))$ implies $\Gamma \vdash
\Tru(\mathit{ap}(\bquo{\cA}, \bquo{u_L}))$.
\end{thm}

\noindent
\emph{Proof.}  Take for $\Gamma$ the context that consists of the
premise $H \ldef \forall x.\, \Tru(\mathit{ap}(\bquo{u_L}, \Num\cdot
x))$, which is well-formed.  The hypothesis rule (\textsf{ir\_hyp})
gives $\Gamma \vdash H$, and the rule would give $H \vdash
\Tru(\mathit{ap}(\bquo{\cA}, \bquo{u_L}))$, which
\textsf{all\_omega\_ent\_fails} refutes (the entailment cell of the
row for $\forall$ in Table~\ref{tab:c-truth-axioms}).  \qed

\begin{thm}[\textsf{tru\_\allowbreak univ\_\allowbreak reflect}, \textsf{true\_\allowbreak reflect\_\allowbreak chk\_\allowbreak ent\_\allowbreak fails}]
\label{thm:c-tin-reflect}
\leavevmode
\begin{enumerate}
\item For every well-formed closed term $u$,
  $\Chk(\vz, \bquo{\cA\cdot u}),\ \vz\jnat \vdash \Tru(\bquo{\cA\cdot
  u})$.
\item $\Chk(\vz, \bquo{\z = \tofn 1}),\ \vz\jnat \nvdash \Tru(\bquo{\z
  = \tofn 1})$.
\end{enumerate}
\end{thm}

\noindent
\emph{Proof.}  Part~1 is the instance at the hypothesis of the rule
\textsf{tru\_univ\_reflect\_at}: if $\Gamma \vdash a\jnat$ and $\Gamma
\vdash \Chk(a, \bquo{\cA\cdot u})$, then $\Gamma \vdash
\Tru(\bquo{\cA\cdot u})$.  The tree $\mathit{reflCert}(u, a)$
(\textsf{reflCert}) adds two nodes over a fixed three-node core
$\{\cA\cdot u, \vz\jnat\} \vdash u\cdot\vz$, made of two hypotheses
and one universal elimination, and over $a$'s closed proof; the
checker accepts the tree as a proof of $J(\bquo{u})$
(\textsf{reflCert\_proves}).  Theorem~\ref{thm:c-tin-cert}(1) gives
$\Tru(\mathit{ap}(\bquo{\cA}, \bquo{u}))$, and transport along
$\mathit{ap}(\bquo{\cA}, \bquo{u}) = \bquo{\cA\cdot u}$
(\textsf{mkAS\_app\_val}, \textsf{iWtrue\_transport}) finishes.
Part~2 (\texttt{RGA\_Goedel/\allowbreak Entailment\_Goedel}): the
hypothesis $\Chk(\vz, \bquo{\z = \tofn 1})$ is decided
(\textsf{hq\_iChk}), so the deduction theorem at that hypothesis and
universal introduction would give the sentence that
Theorem~\ref{thm:pl-no-uniform} excludes.  \qed

\begin{thm}[\textsf{iLiftS\_\allowbreak numC}]
\label{thm:c-tin-lift}
If $\Gamma \vdash k\jnat$, then $\Gamma \vdash \mathit{lift}(\Num\cdot
k) = \Num\cdot k$.
\end{thm}

\noindent
\emph{Proof.}  One internal induction on $k$ (\textsf{ir\_nat\_ind}).
The proof uses the laws that the compiled lift commutes with the
compiled successor and application formers
(\textsf{iLiftS\_mkSucS}, \textsf{iLiftS\_mkAppS}, through the
operator reads \textsf{liftVal\_opr\_read\_gen1} and
\textsf{liftVal\_opr\_read\_gen2}) and is the identity at closed
literal codes (\textsf{iLiftS\_num\_closed}).  \qed

\begin{thm}[\textsf{uinCert\_\allowbreak proves}, \textsf{uinM\_\allowbreak eq}]
\label{thm:c-tin-builder}
Let $G$ be a well-formed closed term and $Q$ a proof with $Q \proves
\emptyset \vdash \forall x.\,\forall v.\ G\cdot x\cdot v$, and let
$\Gamma \vdash x\jnat$.
\begin{enumerate}
\item $\Gamma \vdash \cmp(\pi)\cdot J(\mathit{ap}(\bquo{G}, \Num\cdot
  x))\cdot \mathit{uinCert}(Q, G, x) \neq \z$.
\item $\Gamma \vdash \mathit{uinM}_{Q,G}\cdot x\cdot x =
  \mathit{uinCert}(Q, G, x)$, where $\mathit{uinM}_{Q,G} \ldef
  \cmp(\mathit{uinIdx}(Q, G))$.
\end{enumerate}
\end{thm}

\noindent
The tree $\mathit{uinCert}(Q, G, x)$ (\textsf{uinCert}) puts four nodes
in front of the tree that applies $\bquo{Q}$ to the one-node proof
\textsf{hypQ0} of $\{\vz\jnat\} \vdash \vz\jnat$ and to the
certificate $\mathit{natCert}(x, x)$ of
{}[\citen{goedel-paper}, theorem~\textsf{natCert\_zero}], which certifies that the numeral of
$x$ is a number.  The four nodes are a universal elimination at
$\emptyset$ at that numeral, a forward conversion along the code
identity of the head reduct (from \textsf{iHeadRedC\_iSubM}), a
weakening to $\{\vz\jnat\}$, and a universal elimination at
$\{\vz\jnat\}$ over \textsf{hypQ0}.  The root's conclusion matches the
conclusion in $J(\mathit{ap}(\bquo{G}, \Num\cdot x))$ because the
compiled lift fixes $\mathit{ap}(\bquo{G}, \Num\cdot x)$
(\textsf{iLiftS\_inst}, from Theorem~\ref{thm:c-tin-lift}).  The index
$\mathit{uinIdx}(Q, G)$ (\textsf{uinIdx}) is one explicit
primitive-recursive tree, with no index chosen by Hilbert choice.
\emph{Proof.}  Part~1 verifies the tree's acceptance node by node; the
numeral certificate enters through \textsf{iChk\_natCert}.  Part~2
evaluates the index node by node (\textsf{ir\_mkAppSI},
\textsf{ir\_certConsIdxI}, \textsf{ir\_certAppIdxI}).  \qed

\begin{thm}[\textsf{univ\_\allowbreak inst\_\allowbreak nec}, \textsf{univ\_\allowbreak inst\_\allowbreak nec\_\allowbreak 1}]
\label{thm:c-tin-nec}
\leavevmode
\begin{enumerate}
\item Let $G$ be well-formed and closed with $\emptyset \vdash \forall
  x.\,\forall v.\ G\cdot x\cdot v$.  Then $\emptyset \vdash \forall
  x.\ \Tru(\mathit{ap}(\bquo{\cA}, \mathit{ap}(\bquo{G}, \Num\cdot
  x)))$.
\item At $G_1 \ldef \lambda\lambda(v_1 + \vz = v_1 + \vz)$, with the
  compiled addition, $\emptyset \vdash \forall x.\
  \Tru(\mathit{ap}(\bquo{\cA}, \mathit{ap}(\bquo{G_1}, \Num\cdot
  x)))$.
\end{enumerate}
\end{thm}

\noindent
\emph{Proof.}  Part~1 (Theorem~\ref{thm:pl-tin-nec}): HOL takes a
proof $Q$ of the premise (\textsf{entails\_provesE}).  At the context
$\{\vz\jnat\}$, transport of Theorem~\ref{thm:c-tin-builder}(1) along
part~2 of that theorem gives the verdict at
$\mathit{uinM}_{Q,G}\cdot\vz\cdot\vz$, and
Theorem~\ref{thm:c-tin-cert}(1) gives the truth
(\textsf{univ\_inst\_nec\_at}); universal introduction finishes.
Part~2 is part~1 at $G_1$, whose premise \textsf{uinG1\_thm} is
derived by two universal introductions and two $\beta$-conversions.
At each numeral $\tofn k$, $\emptyset \vdash \mathit{ap}(\bquo{\cA},
\mathit{ap}(\bquo{G}, \Num\cdot\tofn k)) = \bquo{\cA\cdot(G\cdot\tofn
k)}$ (\textsf{univ\_inst\_code\_at}).  \qed

\noindent
The results of \S\ref{sec:pillars-reflect} on the guard discharge
(Theorem~\ref{thm:pl-tin-disch}) are in
\texttt{RGA\_Tin/\allowbreak Tin\_Univ2}.  In what follows,
$\mathit{env}_2(x, y)$ (\textsf{tinEnv2}) is the environment that maps
$\vz$ to the code of the numeral of $x$ and $v_1$ to the code of the
numeral of $y$, and $\mathit{img}_2(x, y, t) \ldef
\Subst(\mathit{env}_2(x, y), \bquo{t})$ (\textsf{tinImg2}) is the image
of the template $t$ under that environment.  For a closed $G$, write
$c_G(x, y) \ldef \mathit{ap}(\mathit{ap}(\bquo{G}, \Num\cdot x),
\Num\cdot y)$ (\textsf{tinCa2}).

\begin{thm}[\textsf{mpCert\_\allowbreak chk}]
\label{thm:c-tin-mp}
Let $\Gamma$ derive $c_1\jnat$, $c_2\jnat$, $\IsCode(c_i) \neq \z$ and
$\Wfc(c_i) \neq \z$ for $i = 1, 2$, $A\jnat$, $B\jnat$,
$\Chk(A, \mathit{or}(\mathit{neg}(c_1), c_2))$ and $\Chk(B, c_1)$.
Then $\Gamma \vdash \Chk(\mathit{mpCert}(c_1, c_2, A, B), c_2)$.
\end{thm}

\noindent
The formula codes $c_1$ and $c_2$ of $\varphi$ and $\psi$ are arbitrary
terms, so the theorem is modus ponens at abstract formula codes, with
the implication $\varphi \limp \psi$ coded as
$\mathit{or}(\mathit{neg}(c_1), c_2)$.  The tree
$\mathit{mpCert}(c_1, c_2, A, B)$ (\textsf{mpCert}) puts five nodes
in front of the concatenation of $A$ and $B$: a disjunction
elimination at the root, a negation elimination, two hypothesis nodes
and one weakening.  The only contexts in these nodes are the
one-element contexts $\{\neg\varphi\}$ and $\{\psi\}$, coded as
one-element lists.  \emph{Proof.}  The acceptance of the tree is
verified node by node, through a law for adding a node whose context is
symbolic (\textsf{tin\_cons\_valid}).  The proof uses two new internal
laws about one-element contexts.  If $\Gamma \vdash c\jnat$, the
insertion of $c$ into the empty context is the one-element list
$[c]$ (\textsf{tin\_hins\_single}, by reading the normalizer's passes
that remove duplicates and sort at length one).  If moreover $\Gamma
\vdash \IsCode(c) \neq \z$, the code test for finite sets accepts $[c]$
(\textsf{tin\_fsetcode\_single}).  A context of two or more symbolic
codes would need an internal law that the normalizer's sort pass
yields a sorted list, which the development does not prove; the
one-element contexts avoid that obligation.  \qed

\begin{thm}[\textsf{univ2Cert\_\allowbreak proves}, \textsf{univ2M\_\allowbreak eq}]
\label{thm:c-tin-univ2b}
Let $G$ be a well-formed closed term and $Q$ a proof with $Q \proves
\emptyset \vdash \forall x.\,\forall y.\,\forall v.\ G\cdot x\cdot
y\cdot v$, and let $\Gamma \vdash x\jnat$ and $\Gamma \vdash y\jnat$.
\begin{enumerate}
\item $\Gamma \vdash \cmp(\pi)\cdot J(c_G(x, y))\cdot
  \mathit{univ2Cert}(Q, G, x, y) \neq \z$.
\item $\Gamma \vdash \mathit{univ2M}_{Q,G}\cdot x\cdot y =
  \mathit{univ2Cert}(Q, G, x, y)$, where $\mathit{univ2M}_{Q,G} \ldef
  \cmp(\mathit{univ2Idx}(Q, G))$.
\end{enumerate}
\end{thm}

\noindent
The tree $\mathit{univ2Cert}(Q, G, x, y)$ (\textsf{univ2Cert}) has
two parts.  The instance chain (\textsf{inst2Cert}) puts four nodes in
front of the tree that applies $\bquo{Q}$ to the numeral certificates
$\mathit{natCert}(x, x)$ and $\mathit{natCert}(y, y)$ of
{}[\citen{goedel-paper}, theorem~\textsf{natCert\_zero}]: a universal elimination at each numeral,
each followed by a forward conversion along the head step, with every
node's code the image under $\mathit{env}_2(x, y)$ of a template.  The
open top (\textsf{openTop}) adds a weakening to $\{\vz\jnat\}$ and a
universal elimination over the one-node proof \textsf{hypQ0}.  The
index $\mathit{univ2Idx}(Q, G)$ (\textsf{univ2Idx}) is one explicit
primitive-recursive tree.  \emph{Proof.}  Part~1 verifies the instance
chain (\textsf{inst2\_chk}, with \textsf{iHeadRedC\_iSubM} for both
head steps) and the open top (\textsf{openTop\_proves}, with
Theorem~\ref{thm:c-tin-lift} for the lift).  Part~2 evaluates the index
node by node.  \qed

\begin{thm}[\textsf{univ2\_\allowbreak nec}]
\label{thm:c-tin-univ2}
Let $G$ be well-formed and closed with $\emptyset \vdash \forall
x.\,\forall y.\,\forall v.\ G\cdot x\cdot y\cdot v$.  Then $\emptyset
\vdash \forall x.\,\forall y.\ \Tru(\mathit{ap}(\bquo{\cA}, c_G(x, y)))$.
\end{thm}

\noindent
\emph{Proof.}  As for Theorem~\ref{thm:c-tin-nec}(1), with
Theorem~\ref{thm:c-tin-univ2b} in place of
Theorem~\ref{thm:c-tin-builder} (\textsf{univ2\_nec\_at}).  \qed

\noindent
At $G \ldef \lambda\lambda\lambda\,R$ (\textsf{polyG}), where $R$ is
the body of the sentence of Theorem~\ref{thm:pl-rf-sentence} below its
three quantifiers, the premise follows from that sentence by three head
steps (\textsf{polyG\_thm}), and Theorem~\ref{thm:c-tin-univ2} gives
\textsf{poly\_univ2\_nec}, with the guards still inside the
universal.

\begin{thm}[\textsf{guardCert\_\allowbreak chk}]
\label{thm:c-tin-guard}
Let $f$ be an index, and let $\Gamma \vdash x\jnat$, $\Gamma \vdash
y\jnat$ and $\Gamma \vdash \cmp(f)\cdot x\cdot y = \z$.  Then
$\Gamma \vdash \Chk(\mathit{guardCert}(f, x, y),
\mathit{img}_2(x, y, \cmp(f)\cdot\vz\cdot v_1 = \z))$.
\end{thm}

\noindent
The tree $\mathit{guardCert}(f, x, y)$ (\textsf{guardCert}) is the
value certificate that the builder $B_f$ of
{}[\citen{goedel-paper}, theorem~\textsf{pr\_cert\_builder}] computes for the test at the two
numeral certificates.  \emph{Proof.}  [\citen{goedel-paper}, theorem~\textsf{pr\_cert\_builder}]
gives a tree accepted for the equation between the image of
$\cmp(f)\cdot\vz\cdot v_1$ and the code of the numeral of
$\cmp(f)\cdot x\cdot y$.  The hypothesis turns that numeral code into
the code of $\z$, an equality of codes that adds no node to the tree,
and transport of the checked conclusion finishes.  \qed

\begin{thm}[\textsf{dischCert\_\allowbreak proves}, \textsf{dischM\_\allowbreak eq}]
\label{thm:c-tin-dischb}
Let $G$ be well-formed and closed, $f$ an index and $Q$ a proof with
$Q \proves \emptyset \vdash \forall x.\,\forall y.\ \cmp(f)\cdot x\cdot
y = \z \limp \forall v.\ G\cdot x\cdot y\cdot v$, and let $\Gamma
\vdash x\jnat$ and $\Gamma \vdash y\jnat$.
\begin{enumerate}
\item If $\Gamma \vdash \cmp(f)\cdot x\cdot y = \z$, then $\Gamma
  \vdash \cmp(\pi)\cdot J(c_G(x, y))\cdot \mathit{dischCert}(Q, f, G,
  x, y) \neq \z$.
\item $\Gamma \vdash \mathit{dischM}_{Q,f,G}\cdot x\cdot y =
  \mathit{dischCert}(Q, f, G, x, y)$, where $\mathit{dischM}_{Q,f,G}
  \ldef \cmp(\mathit{dischIdx}(Q, f, G))$.
\end{enumerate}
\end{thm}

\noindent
The tree $\mathit{dischCert}(Q, f, G, x, y)$ (\textsf{dischCert}) is
the open top over modus ponens (Theorem~\ref{thm:c-tin-mp}) applied to
the instance chain, which proves the implication at the two numerals,
and to the value certificate of Theorem~\ref{thm:c-tin-guard}, which
proves the implication's antecedent.  The index $\mathit{dischIdx}(Q,
f, G)$ (\textsf{dischIdx}) is one explicit primitive-recursive tree.
\emph{Proof.}  Part~1 composes \textsf{inst2\_chk},
Theorem~\ref{thm:c-tin-guard}, Theorem~\ref{thm:c-tin-mp} and
\textsf{openTop\_proves}, with the image of the implication read as
$\mathit{or}(\mathit{neg}(\cdot), \cdot)$ of the two images.  Part~2
evaluates the index node by node.  \qed

\begin{thm}[\textsf{univ2\_\allowbreak disch}]
\label{thm:c-tin-disch}
Let $G$ be well-formed and closed, $f$ an index, and $\emptyset \vdash
\forall x.\,\forall y.\ \cmp(f)\cdot x\cdot y = \z \limp \forall v.\
G\cdot x\cdot y\cdot v$.  Then
\[
  \emptyset \vdash \forall x.\,\forall y.\ \cmp(f)\cdot x\cdot y = \z
  \limp \Tru(\mathit{ap}(\bquo{\cA}, c_G(x, y))) .
\]
\end{thm}

\noindent
\emph{Proof.}  HOL takes a proof $Q$ of the premise.  At the context
$\{\cmp(f)\cdot v_1\cdot\vz = \z, \vz\jnat, v_1\jnat\}$, transport of
Theorem~\ref{thm:c-tin-dischb}(1) along part~2 gives the verdict at
$\mathit{dischM}_{Q,f,G}\cdot v_1\cdot\vz$, and
Theorem~\ref{thm:c-tin-cert}(1) gives the truth
(\textsf{univ2\_disch\_at}).  The deduction theorem at the decided
hypothesis (\textsf{ir\_impI\_bool}) and two universal introductions
finish.  \qed

\begin{thm}[\textsf{poly\_\allowbreak all\_\allowbreak sentence}, \textsf{polyD\_\allowbreak thm}]
\label{thm:c-tin-poly}
Let $Q$ be as in Theorem~\ref{thm:pl-tin-disch}, and let $d$ be the
compiled test \mbox{\textsf{tinPdIdx}}, whose value is zero exactly when the
seven guards of Theorem~\ref{thm:pl-rf-sentence} are nonzero.
\begin{enumerate}
\item If $\Gamma$ derives $x\jnat$, $y\jnat$ and the seven guards at
  $x$ and $y$, then $\Gamma \vdash \cmp(d)\cdot x\cdot y = \z$
  (\textsf{tinPd\_I}).
\item $\emptyset \vdash \forall x.\,\forall y.\ \cmp(d)\cdot x\cdot y =
  \z \limp \forall v.\ Q\cdot x\cdot y\cdot v$.
\item The sentence of Theorem~\ref{thm:pl-tin-disch}(1).
\end{enumerate}
\end{thm}

\noindent
\emph{Proof.}  Part~1 evaluates the conjunction of the seven tests
inside \rga; the converse direction, from the test's value zero to the
seven guards, also holds (\textsf{tinPd\_D}, through
\textsf{eqT\_false\_E}).  Part~2 instantiates
Theorem~\ref{thm:pl-rf-sentence}, discharges its seven antecedents
through \textsf{tinPd\_D}, and reaches $Q\cdot x\cdot y\cdot v$ by
three conversions (\textsf{polyQ\_conv}).  Part~3: HOL takes a proof of part~2, and at the
context of the seven guards and the two typings, part~1 gives the
test's value zero and \textsf{univ2\_disch\_at} gives the truth; seven
deductions at decided hypotheses and two universal introductions
finish.  The generic form at $Q$ is \textsf{poly\_disch}
(Theorem~\ref{thm:c-tin-disch} at $Q$ and $d$).  \qed

\begin{thm}[\textsf{tru\_\allowbreak univ\_\allowbreak by\_\allowbreak reflection}, \textsf{rfl\_\allowbreak binom\_\allowbreak tru\_\allowbreak univ}]
\label{thm:c-tin-showcase}
\leavevmode
\begin{enumerate}
\item Let $a$ and $b$ be polynomial templates whose HOL normal forms
  agree after trimming.  Then \rga proves, from no hypotheses,
  $\Tru(\bquo{\cA\cdot(Q\cdot\bquo{a}\cdot\bquo{b})})$.
\item Part~1 holds at the two sides of each of the three identities of
  Corollary~\ref{cor:pl-rf-binom} and Theorem~\ref{thm:c-rf-more}
  (\textsf{rfl\_binom\_tru\_univ}, \textsf{rfl\_cube\_tru\_univ},
  \textsf{rfl\_prod\_tru\_univ}).
\end{enumerate}
\end{thm}

\noindent
\emph{Proof.}  Part~1: the test's value at the code numerals,
$\emptyset \vdash \mathit{polyEq}(\bquo{a}, \bquo{b}) = \tofn 1$,
comes from a closed evaluation (\textsf{cmp\_eq\_num},
\textsf{pr2fi\_polyEqIdx}), and the six code guards at the numerals
are closed computations as well (\textsf{pl\_code\_guards}, through
\textsf{iWft\_code\_ne0}).
Theorem~\ref{thm:c-tin-poly}(1) gives the test's value zero,
\textsf{univ2\_disch\_at} at the numerals gives
$\Tru(\mathit{ap}(\bquo{\cA}, c_Q(\bquo{a}, \bquo{b})))$, and transport
along $\mathit{ap}(\bquo{\cA}, c_Q(\bquo{a}, \bquo{b})) =
\bquo{\cA\cdot(Q\cdot\bquo{a}\cdot\bquo{b})}$ finishes.  Part~2: membership
in the class by its introduction rules and one \textsf{simp} evaluation
of the normal forms, as for Corollary~\ref{cor:pl-rf-binom}.  \qed

\paragraph{The dependency checks.}
Each of the four theories ends with a check that reads Isabelle's
theorem-dependency graph at build time and fails the build if a
condition fails.  In \texttt{Reflect\_QF}, no oracle lies under the
twelve headline facts.  \textsf{qfm\_tout} and \textsf{rfl\_tout\_binom}
reach \textsf{iWtrue\_collision\_qf} and \textsf{pr\_dec\_eval}.  As a
control that the test is not blind, \textsf{qfm\_tin} reaches
\textsf{pr\_dec\_eval} and does not reach
\textsf{iWtrue\_collision\_qf}.  In \texttt{Poly\_Reflect}, no oracle
lies under the fourteen headline facts.  \textsf{poly\_reflect} reaches
\textsf{norm\_correct}, the homomorphisms for addition, multiplication
and trimming, \textsf{ptm\_tout\_eq},
\textsf{recT\_unfold}, the sentence, \textsf{polyeq\_correct},
\textsf{cmp\_eq\_num} and \textsf{pr2fi\_polyEqIdx}.
\textsf{norm\_correct} reaches \textsf{cov\_pair\_ind},
\textsf{recT\_unfold}, \textsf{recT\_unfold2},
\textsf{iPrDecVal\_app\_img} and \textsf{pl\_subf\_iPrDecVal}, and each
of the three identities reaches \textsf{poly\_reflect} and
\textsf{pr2fi\_polyEqIdx}.  Neither \textsf{poly\_reflect} nor any of the three identities reaches
\textsf{rfl\_tout\_binom} or \textsf{rfl\_tout\_binom\_generic}, the
hand-derived proofs of the binomial identity.  As a control that this
avoidance check is not blind, the build also runs it on a lemma that
derives the binomial identity through \textsf{rfl\_tout\_binom},
\textsf{rfl\_binom\_by\_tout}, and the build fails unless the check
flags that lemma.  A second control, \textsf{rfl\_tout\_binom} reaches
\textsf{ptm\_tout\_eq}, shows that the dependency reading is not blind.
In \texttt{Tin\_Univ}, no oracle lies under the seventeen headline
facts.  \textsf{tru\_\allowbreak univ\_\allowbreak of\_\allowbreak cert} reaches
\textsf{ir\_\allowbreak A\_\allowbreak cert\_\allowbreak affirm}, \textsf{iWtrueI} and
\textsf{ir\_\allowbreak certTrue\_\allowbreak intro}; \textsf{all\_\allowbreak omega\_\allowbreak ctx\_\allowbreak rule\_\allowbreak fails}
reaches \textsf{all\_\allowbreak omega\_\allowbreak ent\_\allowbreak fails}; \textsf{tru\_\allowbreak univ\_\allowbreak cert\_\allowbreak conv}
reaches \textsf{certTrue\_\allowbreak hasCert}, the certificate-search reads
\textsf{ir\_\allowbreak stepV\_\allowbreak app\_\allowbreak cert\_\allowbreak miss}, \textsf{ir\_\allowbreak stepV\_\allowbreak app\_\allowbreak cert\_\allowbreak hitF}
and \textsf{ir\_\allowbreak stepV\_\allowbreak app\_\allowbreak headmiss}, and the fuel laws
\textsf{dec\_\allowbreak zero\_\allowbreak fuel}, \textsf{decA\_\allowbreak one\_\allowbreak fuel} and
\textsf{decA\_\allowbreak ss\_\allowbreak cert}; \textsf{tru\_\allowbreak univ\_\allowbreak reflect} reaches
\textsf{reflCert\_\allowbreak proves} and \textsf{tru\_\allowbreak univ\_\allowbreak of\_\allowbreak cert}; and
\textsf{univ\_\allowbreak inst\_\allowbreak nec\_\allowbreak 1} reaches \textsf{tru\_\allowbreak univ\_\allowbreak of\_\allowbreak cert},
\textsf{uinCert\_\allowbreak proves}, \textsf{uinM\_\allowbreak eq}, \textsf{iLiftS\_\allowbreak numC},
\textsf{iHeadRedC\_\allowbreak iSubM} and \textsf{iChk\_\allowbreak natCert}.  The five facts
\textsf{tru\_\allowbreak univ\_\allowbreak of\_\allowbreak cert}, \textsf{tru\_\allowbreak univ\_\allowbreak cert\_\allowbreak int},
\textsf{tru\_\allowbreak univ\_\allowbreak cert\_\allowbreak conv\_\allowbreak int}, \textsf{tru\_\allowbreak univ\_\allowbreak reflect} and
\textsf{univ\_\allowbreak inst\_\allowbreak nec\_\allowbreak 1} reach none of the soundness and
consistency theorems \textsf{rga\_\allowbreak omega\_\allowbreak sound}, \textsf{omega\_\allowbreak sound},
\textsf{prf\_\allowbreak owsjt}, \textsf{prf\_\allowbreak oomjt}, \textsf{rga\_\allowbreak consistent},
\textsf{rga\_\allowbreak no\_\allowbreak proof\_\allowbreak 0eq1}, \textsf{iChk\_\allowbreak genuine\_\allowbreak proves} and
\textsf{rga\_\allowbreak goedel2}.  As a control that this avoidance check is not
blind, the build also runs it on \textsf{all\_\allowbreak omega\_\allowbreak ctx\_\allowbreak rule\_\allowbreak fails},
which reaches one of these theorems, and the build fails unless the
check flags that fact.
In \texttt{Tin\_Univ2}, no oracle lies under the twenty-one headline
facts.  \textsf{poly\_\allowbreak univ2\_\allowbreak nec} reaches
\textsf{univ2M\_\allowbreak eq}, \textsf{univ2\_\allowbreak nec\_\allowbreak at},
\textsf{tru\_\allowbreak univ\_\allowbreak of\_\allowbreak cert}, \textsf{iLiftS\_\allowbreak numC},
\textsf{iProvesChk\_\allowbreak headI}, \textsf{d2\_\allowbreak cons\_\allowbreak valid},
\textsf{iHeadRedC\_\allowbreak iSubM}, \textsf{iChk\_\allowbreak natCert} and
\textsf{poly\_\allowbreak reflect\_\allowbreak sentence}.
\textsf{poly\_\allowbreak all\_\allowbreak sentence} reaches \textsf{mpCert\_\allowbreak chk},
\textsf{tin\_\allowbreak hins\_\allowbreak single}, \textsf{tin\_\allowbreak fsetcode\_\allowbreak single},
\textsf{guardCert\_\allowbreak chk}, \textsf{pr\_\allowbreak cert\_\allowbreak builder},
\textsf{tin\_\allowbreak cons\_\allowbreak valid}, \textsf{dischM\_\allowbreak eq},
\textsf{univ2\_\allowbreak disch\_\allowbreak at}, \textsf{tru\_\allowbreak univ\_\allowbreak of\_\allowbreak cert},
\textsf{iLiftS\_\allowbreak numC}, \textsf{polyD\_\allowbreak thm},
\textsf{poly\_\allowbreak reflect\_\allowbreak sentence}, \textsf{tinPd\_\allowbreak I} and
\textsf{eqT\_\allowbreak false\_\allowbreak E}.  \textsf{tru\_\allowbreak univ\_\allowbreak by\_\allowbreak reflection}
reaches \textsf{univ2\_\allowbreak disch\_\allowbreak at}, \textsf{mpCert\_\allowbreak chk},
\textsf{cmp\_\allowbreak eq\_\allowbreak num}, \textsf{pr2fi\_\allowbreak polyEqIdx},
\textsf{iWft\_\allowbreak code\_\allowbreak ne0} and \textsf{tinPd\_\allowbreak I}.  The nine facts
\textsf{mpCert\_\allowbreak chk}, \textsf{univ2\_\allowbreak nec}, \textsf{poly\_\allowbreak univ2\_\allowbreak nec},
\textsf{univ2\_\allowbreak disch}, \textsf{poly\_\allowbreak all\_\allowbreak sentence},
\textsf{tru\_\allowbreak univ\_\allowbreak by\_\allowbreak reflection} and the three instances
\textsf{rfl\_\allowbreak binom\_\allowbreak tru\_\allowbreak univ}, \textsf{rfl\_\allowbreak cube\_\allowbreak tru\_\allowbreak univ} and
\textsf{rfl\_\allowbreak prod\_\allowbreak tru\_\allowbreak univ} reach none of the eight soundness and
consistency theorems above, and the build fails unless the same check
flags the same control, \textsf{all\_\allowbreak omega\_\allowbreak ctx\_\allowbreak rule\_\allowbreak fails}.

\endgroup

\section{The consistency folds}
\label{app:fold}
\begingroup\sloppy

This appendix proves the consistency theorems of Section~\ref{sec:con}
that rest on the fold, and the containment and residual results around
those theorems, at the granularity of the mechanized lemmas that carry them.  The appendix covers the $\forall$-free checker
(\S\ref{app:fold-chkqf}), the decider fold behind
Theorem~\ref{thm:con-qf} (\S\ref{app:fold-qf}), the value record (\S\ref{app:fold-value}), the
record fold that the other consistency and soundness theorems of the
section share
(\S\ref{app:fold-record}), the fragment (\S\ref{app:fold-frag}), the
containment of the accepted $\forall$-free codes in the fragment and
the witness beyond it (\S\ref{app:fold-contain}), and the residual
(\S\ref{app:fold-res}).  Appendix~\ref{app:oga} on the level theorems
and \oga reuses the record fold for those theorems.  As in
{}[\citen{goedel-paper}, Appendix ``G\"odel II in full''] on G\"odel's second theorem, each
technique is exhibited once, repetitions are compressed into tables,
and facts below the main chain of lemmas are cited by theory file.  Theory files are
relative to \texttt{hol/}.  Every name resolves in the development
(Section~\ref{sec:mech}, ``Artifact'').

\paragraph{Notation at the level of codes.}
$\mathit{ap}(a, b)$ is the compiled constructor that builds the code of
an application from the codes of its parts (\textsf{mkAppS}), and
$\Num\cdot n$ is the code of the numeral of $n$ (\textsf{numC}~$n$,
compiled from the explicit tree of $\Num$, \S\ref{sec:inside-d3}).  An
\emph{environment} is a coded list of pairs of a variable and a term
code; the empty environment has code $\z$.  $\Subst(z, c)$ is the image
of the code $c$ under the environment $z$ (\S\ref{sec:pillars-sound}),
and $z[w{:=}X]$ updates $z$ at the variable $w$ (\textsf{updSubM}).
$\mathit{head}(p)$ is the first judgment of the proof list coded by
$p$, and $\mathit{hyps}$ and $\mathit{concl}$ read the codes of a
judgment's hypothesis list and conclusion.  For terms $a$ and $b$ that
\rga proves natural, $a < b$ abbreviates the decided test
$\suc a \mathbin{\dot-} b = \z$, with $\dot-$ the compiled truncated
subtraction.

\subsection{Names and theory files}
\label{app:fold-files}

Table~\ref{tab:c-files} gives the theory file of each name that
Section~\ref{sec:con} and this appendix cite.

\begingroup
\footnotesize

\begin{longtable}{@{}>{\raggedright\arraybackslash}p{0.33\linewidth}>{\raggedright\arraybackslash}p{0.63\linewidth}@{}}
\caption{The theory files of the main names of Section~\ref{sec:con}
  and Appendix~\ref{app:fold}; the names of the level classes are in
  Appendix~\ref{app:oga-files}.}
\label{tab:c-files}\\
\hline
theory file & names \\
\hline
\endfirsthead
\hline
theory file & names \\
\hline
\endhead
\hline
\endfoot
\texttt{RGA/RGA\_Omega\_QF} & \textsf{afree} \\
\texttt{RGA/RGA\_Normal\_QF} & \textsf{rga\_afree\_normal} \\
\texttt{RGA/RGA\_Paradox} & \textsf{rga\_deduction\_theorem\_fails} \\
\texttt{RGA\_Kit/Cov\_Induct} & \textsf{ir\_cov\_ind\_par}, \textsf{owsjt\_prenex}, \textsf{cov\_pair\_ind} \\
\texttt{Shadows/Shadow\_QF} & \textsf{afreeP}, \textsf{orec}, \textsf{wosr\_orec} \\
\texttt{RGA\_Checker/Diagonal\_WW} & \textsf{iQfScan} \\
\texttt{RGA\_Checker/Chk\_Formers} & \textsf{iwshsatDT}, \textsf{mTrueAtT} \\
\texttt{RGA\_Checker/Chk\_QF\_Faithful} & \textsf{iChkQf}, \textsf{qfScanAcc}, \textsf{qfScanP\_correct},
  \textsf{qfScanAcc\_faithful}, \textsf{iChkQf\_genuine\_proves}, \textsf{iChkQf\_complete},
  \textsf{iChkQf\_faithful}, \textsf{rga\_afree\_normal\_chk}, \textsf{iChk\_of\_iChkQf},
  \textsf{iChkQf\_of\_iChk} \\
\texttt{RGA\_Checker/Chk\_Base} & \textsf{iSub1S} \\
\texttt{RGA\_Checker/Chk\_Sigma\_Machinery} & \textsf{subm1\_opr\_read\_gen0} \\
\texttt{RGA\_Checker/Chk\_QF\_Node\_Laws} & \textsf{qf\_disjI1\_truth} \\
\texttt{RGA\_QF/QF\_Arms} & \textsf{iChkQfPair}, \textsf{hq\_iChkQf}, \textsf{mOwsjtTQ},
  \textsf{tree\_spineQ}, \textsf{tree\_step\_dispatch\_qf}, \textsf{tree\_ladder\_qf},
  \textsf{qf\_root\_truth}, \textsf{disjI1\_arm\_qf}, \textsf{allI\_node\_refuted},
  \textsf{allE\_node\_refuted}, \textsf{notAllI\_node\_refuted}, \textsf{notAllE\_node\_refuted} \\
\texttt{RGA\_QF/QF\_Last\_Arms} & \textsf{natInd\_stage\_wtrue}, \textsf{natInd\_arm\_qf},
  \textsf{rga\_sound\_qf\_int}, \textsf{rga\_con\_qf\_int}, \textsf{rga\_con\_qf\_meta},
  \textsf{qf\_con\_meta\_num}, \textsf{rga\_consistent\_qf\_via\_chk},
  \textsf{qfRevScan\_at0}, \textsf{rga\_sound\_qf\_int2}, \textsf{iChkQf\_transport} \\
\texttt{RGA\_Record/Value\_Record} & \textsf{iRec}, \textsf{iRecAt}, \textsf{famT}, \textsf{famF},
  \textsf{vfalse\_root3} \\
\texttt{RGA\_Record/Path\_Walk} & \textsf{iEnvAt}, \textsf{iGuardAt}, \textsf{iTreeValid\_node\_row},
  \textsf{nodeTbl\_tag} \\
\texttt{RGA\_Liar\_Div/Record\_Halting} & \textsf{iRec\_of\_wosr}, \textsf{iRec\_faithful},
  \textsf{orec\_iff\_wosr} \\
\texttt{RGA\_Liar\_Div/Liar\_Record} & \textsf{liar\_iRec\_none},
  \textsf{iWtrue\_undecided} \\
\texttt{RGA\_Liar\_Div/Divergence} & \textsf{frag\_search\_diverges}, \textsf{con\_search\_unrefutable},
  \textsf{con\_search\_found},
  \textsf{rga\_con\_int\_of\_con\_search} \\
\texttt{ConRGA/Full\_Fold} & \textsf{fMRb}, \textsf{armT}, \textsf{armsT}, \textsf{critAt},
  \textsf{fold\_root\_record}, \textsf{fold\_consistent\_ctx}, \textsf{fold\_sound\_ctx},
  \textsf{armsT\_of\_fixed}, \textsf{tag\_fire} \\
\texttt{ConRGA/Wfd\_Decision} & \textsf{famDec}, \textsf{lemQ\_int}, \textsf{gdIf},
  \textsf{arm\_Hyp\_path}, \textsf{arm\_AllI\_path} \\
\texttt{ConRGA/Wfd\_Static} & \textsf{swIfG}, \textsf{swIfQ}, \textsf{xqIf}, \textsf{xnIf},
  \textsf{xgIf}, \textsf{arm\_AllE\_pathG}, \textsf{fold\_consistent\_pathG},
  \textsf{fold\_sound\_pathG}, \textsf{poNode\_suffix} \\
\texttt{ConRGA/Swap\_Engine} & \textsf{rec\_swap\_qf} \\
\texttt{ConRGA/Wfd\_General} & \textsf{iFrag}, \textsf{hq\_iFrag}, \textsf{gdIf\_F}, \textsf{swIfQ\_F},
  \textsf{xqIf\_F}, \textsf{rga\_con\_frag\_int}, \textsf{rga\_sound\_frag\_int},
  \textsf{rga\_sound\_frag\_int\_prenex}, \textsf{con\_shell} \\
\texttt{ConRGA/Guarded\_Residual} & \textsf{iRes}, \textsf{swB}, \textsf{gdIf\_iRes}, \textsf{swIf\_iRes},
  \textsf{rga\_con\_under\_res}, \textsf{rga\_sound\_under\_res},
  \textsf{rga\_con\_under\_res\_unguarded}, \textsf{iRes\_unguarded\_iRes},
  \textsf{iRes\_underivable}, \textsf{rga\_sound\_under\_res\_v0} \\
\texttt{ConRGA/ConRGA\_Traps} & \textsf{swapCexProof\_proves}, \textsf{swap\_body\_unguarded\_refutable},
  \textsf{residual\_unguarded\_refutable}, \textsf{swap\_body\_guarded\_vacuous},
  \textsf{path\_env\_root} \\
\texttt{ConOGA/ConOGA\_Frag\_Witness} & \textsf{ruleTag\_code\_cases},
  \textsf{nodeScan\_closedI}, \textsf{fwP}, \textsf{fwP\_accepted}, \textsf{fwP\_tag\_off},
  \textsf{fwP\_in\_frag}, \textsf{fwP\_in\_lv}, \textsf{fwP\_not\_qf},
  \textsf{frag\_strictly\_beyond\_qf}, \textsf{pr2fi\_fhLtF\_least}, \textsf{fwP2},
  \textsf{fwP2\_root\_conv}, \textsf{pr2fi\_ifNestIdx\_cases}, \textsf{iAllLt\_closed\_ne0},
  \textsf{iStat\_closed\_off} \\
\texttt{ConOGA/ConOGA\_Frag\_QF} & \textsf{ruleTag\_afree\_not\_quant},
  \textsf{sub1Ns\_node}, \textsf{iStat\_afree}, \textsf{iFragV\_afree},
  \textsf{fqOkT\_afree\_node}, \textsf{qf\_accepted\_in\_frag}, \textsf{qf\_accepted\_in\_lv},
  \textsf{frag\_lv\_incomparable},
  \textsf{fwP3}, \textsf{fwP3\_tag0}, \textsf{fwP3\_in\_frag}, \textsf{fwP3\_in\_lv},
  \textsf{swap\_member\_beyond\_qf}, \textsf{root\_conv\_of\_reduct}, \textsf{swapCexProof\_root\_conv},
  \textsf{allE\_fwF\_root\_conv}, \textsf{frag\_not\_lv\_at\_root} \\
\texttt{ConOGA/ConOGA\_Frag\_Lv} & \textsf{frag\_accepted\_in\_lv} \\
\texttt{ConOGA/ConOGA\_Residual\_Frag} & \textsf{iNoSwap}, \textsf{hq\_iNoSwap}, \textsf{iResGuard\_of\_gdIf},
  \textsf{iFrag\_iResGuard}, \textsf{iResSwap\_of\_noswap}, \textsf{iFrag\_noswap\_iRes}, \textsf{fwP\_iNoSwap},
  \textsf{fwP\_iRes}, \textsf{fwP3\_iResGuard}, \textsf{crit25\_kid\_node}, \textsf{fwP3\_res\_disj\_gen}, \textsf{conSearch\_no\_value},
  \textsf{fwP3\_iRes\_underivable}, \textsf{fwP3\_iResSwap\_underivable}, \textsf{iRes\_not\_from\_frag},
  \textsf{iRes\_not\_from\_frag\_ctx}, \textsf{iRes\_not\_from\_lv}, \textsf{iRes\_not\_from\_lv\_ctx},
  \textsf{fwP3\_not\_noswap}, \textsf{t01P\_residual\_unguarded\_refutable} \\
\texttt{ConOGA/ConOGA\_Frag\_QF\_Int} & \textsf{rga\_qf\_in\_frag\_int},
  \textsf{qf\_tree\_in\_frag}, \textsf{rga\_con\_qf\_via\_frag\_int}, \textsf{iAllLt\_allI},
  \textsf{nodeScan\_I}, \textsf{qf\_node\_at}, \textsf{iStatR\_of\_qf}, \textsf{rsuc\_c24k\_qf},
  \textsf{natInd\_nc\_qf}, \textsf{disjE1\_slots\_qf}, \textsf{eqsubs\_c40k\_qf},
  \textsf{iSub1S\_qf\_reflect}, \textsf{iSub1S\_nfree\_indep}, \textsf{iSub1Fset\_indep},
  \textsf{sub1Wit\_fire\_at0}, \textsf{iFhLt\_hit0}, \textsf{sub1\_ns\_qf},
  \textsf{qf\_accepted\_in\_frag\_tree} \\
\texttt{ConOGA/ConOGA\_Lv\_Int} & \textsf{rga\_qf\_in\_lv\_int} \\
\texttt{ConOGA/ConOGA\_Level} & \textsf{rga\_qf\_in\_lvi\_int}, \textsf{fwP\_in\_lvi} \\
\texttt{ConOGA/ConOGA\_Level\_Ports} & \textsf{fwP3\_in\_lvi}, \textsf{iRes\_not\_from\_lvi},
  \textsf{iRes\_not\_from\_lvi\_ctx}, \textsf{iRes\_not\_from\_lvi\_omega},
  \textsf{iRes\_not\_from\_lvi\_oga}, \textsf{iRes\_not\_from\_lvi\_under\_lf},
  \textsf{iRes\_not\_from\_lvi\_under\_lf\_ctx} \\
\texttt{ConOGA/ConOGA\_Seam} & \textsf{oga\_omega\_true} \\
\texttt{ConOGA/ConOGA\_Fold} & \textsf{o\_iLFs} \\
\texttt{ConOGA/ConOGA\_Otto} & \textsf{famT\_ofam\_sound}, \textsf{renp\_famT\_ofam} \\
\texttt{RGA/RGA\_Omega\_Sound} & \textsf{rga\_omega\_sound} \\
\texttt{RGA/RGA\_Omega\_Model} & \textsf{wosr\_omv} \\
\texttt{ConOGA/ConOGA\_Guard\_NoGo} & \textsf{gxP}, \textsf{gxP\_valid}, \textsf{gxP\_tag0},
  \textsf{gxP\_root\_kid}, \textsf{gx\_rec1\_pos}, \textsf{gx\_rec1\_neg},
  \textsf{gx\_guard\_dec\_absurd}, \textsf{gdIf\_not\_of\_gx}, \textsf{iResGuard\_not\_uniform},
  \textsf{iResGuard\_not\_ctx}, \textsf{gxP\_not\_frag} \\
\texttt{ConOGA/ConOGA\_Residual\_LF} & \textsf{fwP3\_res\_disj\_ctx},
  \textsf{lf\_res\_disj}, \textsf{lf\_res\_omega}, \textsf{lf\_res\_omega\_ctx},
  \textsf{res\_fwP3\_omega\_disj}, \textsf{oga\_res\_omega}, \textsf{cycY},
  \textsf{omv\_cycY\_unvalued}, \textsf{omBlk}, \textsf{omBlk\_never},
  \textsf{fwP3\_iRes\_not\_oga}, \textsf{iRec\_cycY\_not\_true},
  \textsf{iRec1\_cycY\_unvalued}, \textsf{omv\_resDisj\_cycY},
  \textsf{fwP3\_iRes\_not\_omega}, \textsf{iRes\_not\_from\_lv\_omega},
  \textsf{iRes\_not\_from\_lv\_oga}, \textsf{iRes\_not\_from\_lv\_under\_lf},
  \textsf{iRes\_not\_from\_lv\_under\_lf\_ctx}, \textsf{fwP3\_iRes\_not\_grounded} \\
\texttt{ConOGA/ConOGA\_Motive\_NoGo} & \textsf{phP}, \textsf{phP\_valid}, \textsf{phP\_accepted},
  \textsf{phP\_tag0}, \textsf{phQ2\_tag}, \textsf{lowRows\_orapp\_off}, \textsf{phP\_walk},
  \textsf{phBK}, \textsf{phBK\_unvalued}, \textsf{ph\_guard\_dec\_disj},
  \textsf{phP\_iResGuard\_not\_omega}, \textsf{phP\_iResGuard\_underivable\_omega},
  \textsf{phP\_gdIf\_underivable\_omega}, \textsf{iResGuard\_open\_underivable\_omega},
  \textsf{ph\_motive\_disj}, \textsf{ph\_E3\_disj\_not\_omega}, \textsf{phP\_fMot\_not\_omega},
  \textsf{phMotAll\_not\_omega}, \textsf{fold\_motive\_at\_phP},
  \textsf{fold\_residual\_underivable\_omega}, \textsf{fold\_residual\_not\_omega},
  \textsf{fold\_residual\_open\_underivable}, \textsf{omega\_sound\_ext},
  \textsf{xder\_not\_iResGuard}, \textsf{xder\_not\_motive}, \textsf{xder\_not\_fold\_residual},
  \textsf{rga\_omega\_sound\_ext} \\
\end{longtable}
\endgroup

\subsection{The \texorpdfstring{$\forall$}{A}-free checker}
\label{app:fold-chkqf}

A HOL proof $P$ is \emph{$\forall$-free} (\textsf{afreeP}) when every
hypothesis and every conclusion of every judgment of $P$ is
$\forall$-free in the sense of \textsf{afree}: the constant $\cA$
occurs nowhere in that hypothesis or conclusion.  The compiled test $\qfs(c)$ of
Section~\ref{sec:con-qf} decides $\forall$-freeness for one code.  The whole-proof
scan has the shape of the checker's validity test: an explicit tree
$\tau_{\mathrm{qf}}$ (\textsf{qfScanPIdx}) requires, at every position
of the list coded by $p$, that $\qfs$ passes every hypothesis and the
conclusion of the judgment there.

\begin{defi}[\textsf{qfScanAcc}, \textsf{iChkQf}]
\label{def:c-chkqf}
$\QfScan(p) \ldef \cmp(\tau_{\mathrm{qf}})\cdot p\cdot p = \tofn 1$, and
\[
  \ChkQf(p, c) \ldef
  \bigl((\mathit{TreeValid}(p) \wedge p \ne \z) \wedge \QfScan(p)\bigr)
  \wedge \bigl(|\mathit{hyps}(\mathit{head}(p))| = \z \wedge
  \mathit{concl}(\mathit{head}(p)) = c\bigr) .
\]
\end{defi}

\noindent
Dropping the scan conjunct gives $\Chk(p, c)$, and adding a derived
$\QfScan(p)$ to a derived $\Chk(p, c)$ gives $\ChkQf(p, c)$, in every
context (\textsf{iChk\_of\_iChkQf}, \textsf{iChkQf\_of\_iChk}).

\begin{thm}[\textsf{qfScanP\_correct}, \textsf{qfScanAcc\_faithful}]
\label{thm:c-qfscan}
For every proof $P$ and natural $w$,
$\den{\tau_{\mathrm{qf}}}(\quo{P}, w)$ is $1$ if $P$ is $\forall$-free
and $0$ otherwise.  Hence $\emptyset \vdash \QfScan(\bquo{P})$ if and
only if $P$ is $\forall$-free.
\end{thm}

\noindent
\emph{Proof.}  The scan is a bounded universal over the list length,
like the validity test, and at each position the scan's body reads the
judgment there (\textsf{qfBody\_correct}).  The second statement
follows from the first, because $\emptyset \vdash \cmp(f)\cdot\tofn
x\cdot\tofn y = \tofn 1$ holds exactly when $\den{f}(x, y) = 1$
(\textsf{cmp\_num\_one\_iff}).  \qed

\begin{thm}[\textsf{iChkQf\_genuine\_proves}]
\label{thm:c-chkqf-genuine}
For every natural $p$ and term $c$: if $\emptyset \vdash
\ChkQf(\tofn{p}, \bquo{c})$, then $p = \quo{P}$ for some proof $P$ with
$P \proves \emptyset \vdash c$.
\end{thm}

\noindent
\emph{Proof.}  Drop the scan conjunct and apply the genuineness of the
tree checker, Theorem~\ref{thm:in-chk}(2).  \qed

\noindent
The conclusion does not say that $P$ is $\forall$-free.  At a proof
code, $\forall$-freeness follows from Theorem~\ref{thm:c-qfscan}, and the two theorems together
give the forward direction of Theorem~\ref{thm:con-qf-faithful}.  The
backward direction, \textsf{iChkQf\_complete}, assembles acceptance
from a valid proof list whose head has no hypotheses
(\textsf{iChkQf\_accept\_of\_parts}, from the checker's
\textsf{iChk\_accept\_of\_parts}) and from the scan.

\begin{thm}[\textsf{hq\_iChkQf}]
\label{thm:c-chkqf-dec}
If $\Gamma \vdash p \jnat$ and $\Gamma \vdash c \jnat$, then $\Gamma
\vdash \ChkQf(p, c) \jbool$.
\end{thm}

\noindent
\emph{Proof.}  Each of the five conjuncts is an equation or a
disequation between total compiled terms, hence decided, and a
conjunction of decided terms is decided.  \qed

\noindent
$\ChkQf$ and the other machines of Section~\ref{sec:con} are roots of
the audit of \S\ref{sec:inside-audit}, and are transparent.

\subsection{The decider fold}
\label{app:fold-qf}

The decider fold proves both parts of Theorem~\ref{thm:con-qf}.  The fold is a course-of-values induction inside
\rga over the positions of an accepted, scanned proof list.

\paragraph{The induction engine.}
Every internal induction over coded proofs in this appendix runs on one
engine.

\begin{thm}[\textsf{ir\_cov\_ind\_par}]
\label{thm:c-cov}
Let $m$ be a well-formed term and $\Gamma_0$ a context, with the
variable $v_1$ free in neither, and let $\Gamma_0 \vdash a \jnat$.
Suppose that for every context $\Delta$ that derives everything
$\Gamma_0$ derives, and every $x$ with $\Delta \vdash x \jnat$,
\[
  \Delta \vdash \forall y.\ (y < x) \limp m\cdot y
  \quad\text{implies}\quad
  \Delta \vdash m\cdot x .
\]
Then $\Gamma_0 \vdash m\cdot a$.
\end{thm}

\noindent
The motive $m$ need not be decided.  Only the guard $y < x$ of the
accumulator is decided, and the engine introduces the accumulator's
implication over that guard alone.  Because of the side condition on
$v_1$, the engine's sentences quantify over one packed pair: under a
nested universal the outer bound variable would be free in the root
context.  The two-variable form of Theorem~\ref{thm:con-qf}(2) is
recovered from the packed one at the end (\textsf{rga\_sound\_qf\_int2},
below).

\paragraph{The motive.}
The motive internalizes the companion paper's soundness invariant for
one proof entry, with the decider's verdict in place of grounded truth.

\begin{defi}[\textsf{iwshsatDT}, \textsf{mTrueAtT}, \textsf{mOwsjtTQ}]
\label{def:c-qfmot}
For codes $s$, $z$ and $x$, $\mathit{Sat}_s(z, x)$ says that
$\Aff(s, \Subst(z, h))$ holds for every hypothesis $h$ of
$\mathit{head}(x)$, a compiled bounded universal over the hypothesis
list.  For an environment predicate $W$,
\[
  \mathit{TrueAt}_W(x) \ldef \forall z.\ W\cdot z \limp
    \forall s.\ \mathit{Sat}_s(z, x) \limp
    \Tru(\Subst(z, \mathit{concl}(\mathit{head}(x)))) ,
\]
and the motive at $x$ is $\bigl(\mathit{TreeValid}(x) \wedge x \ne \z
\wedge \QfScan(x)\bigr) \limp \mathit{TrueAt}_W(x)$.
\end{defi}

\noindent
The fold runs at $W = W_{\mathrm{qf}}$ (\textsf{wscQfG}), which says
of $z$ that every entry is a well-formed term code and every image is
$\forall$-free.  $\mathit{Sat}_s(z, x)$ is decided at abstract
arguments, because $\Dec$ is total (\textsf{iwshsatDT\_B}), and so are
$W_{\mathrm{qf}}\cdot z$ and the guard of the motive.  Every implication
in the motive is therefore introduced over a decided antecedent.  The
fuel quantifier stands in front of the hypothesis test.  In HOL this
placement loses nothing: for every hypothesis set $H$ and term $c$, the
companion paper's invariant holds exactly when every well-formed
substitution under which all of $H$ evaluates to $\ctrue$ within one
common fuel $s$ makes $c$ grounded-true (\textsf{owsjt\_prenex}).  The
fold does not use this equivalence, but the equivalence explains why
the motive with the fuel quantifier in front is the natural internal
form.

\paragraph{The chain of lemmas.}
\begin{itemize}\raggedright
\item \textsf{rga\_con\_qf\_int}, \textsf{rga\_sound\_qf\_int} $\leftarrow$
  \textsf{rga\_con\_qf\_int\_of\_arms}, \textsf{rga\_sound\_qf\_int\_of\_arms}
  (with the eight premise arms of Table~\ref{tab:c-qf-arms} supplied);
\item \textsf{rga\_con\_qf\_int\_of\_arms} $\leftarrow$ \textsf{qf\_root\_truth},
  \textsf{iWtrue\_cEq01\_refute} (\texttt{RGA\_Checker/Chk\_Base}),
  Theorem~\ref{thm:in-never};
\item \textsf{qf\_root\_truth} $\leftarrow$ \textsf{tree\_ladder\_qf}
  $\leftarrow$ \textsf{tree\_spineQ\_root} (the engine at the motive)
  and \textsf{tree\_step\_dispatch\_qf} (the step).
\end{itemize}

\begin{thm}[\textsf{qf\_root\_truth}]
\label{thm:c-qf-root}
Assume the eight premise arms of Table~\ref{tab:c-qf-arms}.  Let
$\Gamma_0$ be a context in which $v_1$ is not free.  If $\Gamma_0
\vdash p \jnat$ and $\Gamma_0 \vdash \ChkQf(p, c)$, then $\Gamma_0
\vdash \Tru(\Subst(\z, c))$.
\end{thm}

\noindent
\emph{Proof.}  The first three conjuncts of $\ChkQf(p, c)$ are the
guard of the motive at $p$, so the ladder gives
$\mathit{TrueAt}_{W_{\mathrm{qf}}}(p)$.  Instantiate this formula at the empty
environment, which satisfies $W_{\mathrm{qf}}$, and at fuel $\z$.  The
head of $p$ has no hypotheses, so $\mathit{Sat}_{\z}(\z, p)$ holds
vacuously (\textsf{iwshsatDT\_vac}).  The instance gives $\Tru(\Subst(\z,
\mathit{concl}(\mathit{head}(p))))$, and transport along the last
conjunct of $\ChkQf(p, c)$ gives $\Tru(\Subst(\z, c))$.  \qed

\noindent
Theorem~\ref{thm:con-qf} follows by the decided-matrix split of
\S\ref{sec:con-shape} at the free variable $\vz$, with
Theorem~\ref{thm:c-chkqf-dec} for the split.  In the accepting case,
Theorem~\ref{thm:c-qf-root} gives $\Tru(\Subst(\z, \bquo{\z = \tofn
1}))$, and \textsf{iWtrue\_cEq01\_refute} combines that truth with
Theorem~\ref{thm:in-never} into the numeral equation $\tofn 2 = \tofn
1$, which \rga refutes.  The packed form \textsf{rga\_sound\_qf\_int},
$\emptyset \vdash \forall w.\ \ChkQf(\mathit{fst}(w), \mathit{snd}(w))
\limp \Tru(\Subst(\z, \mathit{snd}(w)))$, is the same argument with the
conclusion left abstract, at one eigenvariable $w$.
Theorem~\ref{thm:con-qf}(2) (\textsf{rga\_sound\_qf\_int2}) follows
at two eigenvariables $p$ and $c$: instantiate the packed sentence at
the compiled pair of $p$ and $c$, move the antecedent $\ChkQf(p, c)$ to
the projections of that pair along the projection equations
(\textsf{iChkQf\_transport}), apply modus ponens, move the consequent
back by transport, and introduce the implication over the decided
verdict and the two universals.

\paragraph{The dispatch and the arms.}
The step at a position $\suc k$ that is accepted and scanned
(\textsf{tree\_step\_dispatch\_qf}) receives the induction hypothesis
at every earlier accepted, scanned, non-zero position, in every
extension of the node's context.  The step dispatches on the checker's
41-row node table (\textsf{iTreeValid\_node\_dispatch}).
Table~\ref{tab:c-qf-arms} groups the 41 arms.  Every arm runs inside a
common frame (\textsf{qf\_arm\_frame1}, \textsf{qf\_arm\_frame2}): the arm
opens the rule's detector, reads each cited premise through the
induction hypothesis, lifts the premise's truth into the arm's
context, and closes through the node's conclusion.  An arm proves only
the step in the middle, from the truth of the cited conclusions to the
truth of the node's conclusion.

\begin{table}[ht]
\centering
\footnotesize
\begin{tabular}{@{}>{\raggedright\arraybackslash}p{0.17\linewidth}>{\raggedright\arraybackslash}p{0.42\linewidth}>{\raggedright\arraybackslash}p{0.35\linewidth}@{}}
\hline
group & rules (Isabelle arm names, without the suffix \textsf{\_arm\_qf}, except the refuted row) & what the step uses \\
\hline
environment-free (20) & \textsf{pr0}, \textsf{hyp}, \textsf{wk1}, \textsf{botTE}, \textsf{eqRefl},
  \textsf{eqSym}, \textsf{neSym}, \textsf{eqTE}, \textsf{sucNe0}, \textsf{sucTI}, \textsf{sucNeI},
  \textsf{sucNeE}, \textsf{sucTE}, \textsf{sucEqE}, \textsf{dnegI}, \textsf{dnegE}, \textsf{disjI1},
  \textsf{disjI2}, \textsf{disjE2}, \textsf{disjE3} &
  one node law of the decider at the node's codes (\texttt{RGA\_Checker/Chk\_QF\_Node\_Laws}) \\
negation elimination (1) & \textsf{negE} & two decider values of one code cannot differ
  (\textsf{iWtrue\_collision\_qf}) \\
fuel monotonicity (8) & \textsf{cut}, \textsf{trueEqI}, \textsf{trueEqE}, \textsf{falseEqI},
  \textsf{falseEqE}, \textsf{eqTI}, \textsf{disjI3}, \textsf{disjE1} &
  $\Dec$ is monotone in the fuel at $\forall$-free codes (\texttt{RGA\_QF/QF\_Mono}) \\
premise arms (8) & \textsf{conv}, \textsf{convE}, \textsf{convEq}, \textsf{convNegE} &
  $\Subst$ commutes with one head step, and $\Tru$ moves along head steps
  (\texttt{RGA\_QF/QF\_Master\_Key}) \\
 & \textsf{RSuc} & the swap engine on the decider's table (\texttt{RGA\_QF/QF\_Swap\_Engine}) \\
 & \textsf{eqSubs}, \textsf{natInd}, \textsf{sub1} & the swap engine, an internal induction on
  the stage numeral, the non-occurrence law and the reverse scan (\texttt{RGA\_QF/QF\_Last\_Arms}) \\
refuted (4) & \textsf{allI}, \textsf{notAllI}, \textsf{allE}, \textsf{notAllE}
  (\textsf{allI\_node\_refuted} and siblings) & the scan of the node finds $\cA$ \\
\hline
\end{tabular}
\caption{The 41 arms of the decider fold.  The 29 arms of the first
  three groups are proved in \texttt{RGA\_QF/QF\_Arms}.  The eight
  premise arms are premises of \textsf{rga\_con\_qf\_int\_of\_arms} and
  its two siblings, discharged in the theories named.}
\label{tab:c-qf-arms}
\end{table}

\paragraph{An environment-free arm: disjunction introduction.}
The rule concludes $\Gamma \vdash a \lor b$ from $\Gamma \vdash a$ and
the well-formedness of $b$.  The frame supplies, for the environment
$z$ and the fuel $s$ of the motive, $\Tru(\Subst(z, c))$ for the code
$c$ of the cited conclusion; the premise has the node's own context,
so the node's $\mathit{Sat}_s(z, \cdot)$ serves for both the premise
and the node.  The
detector supplies three facts about the code $y$ of the node's
conclusion: $y$ is a disjunction node, its first argument is $c$, and
its second argument is the code of a well-formed term.  The step is
\textsf{qf\_disjI1\_truth}.  \rga proves that $y$ equals the
disjunction code built from its two arguments, and that $\Subst(z,
\cdot)$ distributes over that constructor.  The decider's introduction
law for a disjunction (\textsf{iWtrue\_or\_intro1}) gives the truth of
the disjunction of the two images, and transport moves that truth back to
$\Subst(z, y)$.

\paragraph{A premise arm: induction.}
The rule \textsf{pr\_nat\_ind} concludes $\Gamma \vdash
\ssub{i}{a}{c}$ from the base $\Gamma \vdash \ssub{i}{\z}{c}$, the step
$v_i \jnat,\ c,\ \Gamma \vdash \ssub{i}{\suc v_i}{c}$ and $\Gamma
\vdash a \jnat$, where $c$ is well-formed and $v_i$ is not free in
$\Gamma$.  The arm \textsf{natInd\_arm\_qf} follows the companion
paper's HOL soundness proof of the rule, inside \rga.
\begin{enumerate}
\item The induction hypothesis at the third premise gives the truth of
  the image of $a = a$.  So the decider assigns the image of $a$ a
  value $w$ (\textsf{eq\_wtrue\_valE}).
\item The internal induction \textsf{natInd\_stage\_wtrue} is \rga's
  induction rule on a stage numeral $n$, with the motive
  $\Tru(\Subst(z[i{:=}\Num\cdot n], c))$.  The motive is not shown
  decided, and \rga's induction rule does not ask the motive to be.  The base is the
  induction hypothesis at the first premise, moved by the law that
  updating the environment at $i$ with an image equals substituting
  first (\textsf{iSubM\_updSubM\_iSub1S\_G}).
\item In the induction step the arm instantiates the step premise's
  motive at the updated environment.  Three facts make that premise's
  hypotheses hold at one fuel.  The old hypotheses stay satisfied,
  because the checker's freshness datum says that $v_i$ occurs free in
  none of the old hypotheses (\textsf{iwshsatDT\_upd\_nfree\_hyps}).  The hypothesis
  $v_i \jnat$ becomes a numeral equation that the decider affirms
  (\textsf{numC\_dec}).  The hypothesis $c$ is the stage hypothesis.
  The three readings are joined at one fuel by monotonicity
  (\textsf{iwshsatDT\_iHIns2\_ext\_qf}).  The updated environment still
  satisfies $W_{\mathrm{qf}}$, since $\Num\cdot n$ is a $\forall$-free
  well-formed code.  The premise's motive gives the truth of the image
  of $\ssub{i}{\suc v_i}{c}$, and the law that a second update at $i$
  shadows the first (\textsf{updSubM\_shadow}) turns that truth into the next
  stage.
\item The stage at $\Num\cdot w$ is moved to the image of $a$, which
  has the same value $w$, by the swap engine
  (\textsf{cellSwap\_eqSubs}).
\end{enumerate}
$\forall$-freeness enters through the environment predicate
$W_{\mathrm{qf}}$ and through the fuel monotonicity that joins the
three readings.

\paragraph{The refuted quantifier nodes.}
At a $\forall$-introduction or a negated-universal introduction node,
the node's own conclusion contains an application of $\cA$; at a
universal elimination or a negated-universal elimination node, a cited
premise's conclusion does.  The node's scan passes, so the scan of
that conclusion is $\z$; pushing the scan down to the constant $\cA$ refutes
that equation (\textsf{qfScan\_codeA\_refute}).  The arm closes by eliminating the
proved negation.

\paragraph{Meta corollaries.}
Three HOL statements follow from the internal theorems and the
companion paper's HOL soundness.  These statements are sanity checks of the
internal result, not a new consistency argument.

\begin{thm}[\textsf{rga\_con\_qf\_meta}, \textsf{qf\_con\_meta\_num}, \textsf{rga\_consistent\_qf\_via\_chk}]
\label{thm:c-qf-meta}
\leavevmode
\begin{enumerate}
\item If $\emptyset \vdash p \jnat$ and $\emptyset \vdash \ChkQf(p,
  \bquo{\z = \tofn 1})$, then $\mathit{False}$.
\item For every natural $p$: if $\emptyset \vdash \ChkQf(\tofn{p},
  \bquo{\z = \tofn 1})$, then $\mathit{False}$.
\item If $P \proves \emptyset \vdash \z = \tofn 1$ and $P$ is
  $\forall$-free, then $\mathit{False}$.
\end{enumerate}
\end{thm}

\noindent
\emph{Proof.}  (1) Theorem~\ref{thm:c-qf-root} in the empty context
gives $\emptyset \vdash \Tru(\Subst(\z, \bquo{\z = \tofn 1}))$, and the
empty environment leaves the code unchanged
(\textsf{iSubM\_cEq01\_id}), so $\emptyset \vdash \Tru(\bquo{\z = \tofn
1})$.  The HOL refutation \textsf{iWtrue\_no\_false\_root} of
\S\ref{sec:pillars-con}, which goes through the companion paper's
soundness theorem, excludes this derivation.  (2) is (1) at a numeral.  (3) The
checker accepts the code of $P$ (Theorem~\ref{thm:con-qf-faithful}),
an acceptance that (2) excludes.  \qed

\subsection{The value record}
\label{app:fold-value}

The record folds of the remaining theorems use the value record of
\S\ref{sec:con-frag}.  The record's HOL mould is an inductive relation.

\begin{defi}[\textsf{orec}]
\label{def:c-orec}
$\mathrm{orec}(r, t)$ is the least relation between naturals and terms
closed under the following rules: $\mathrm{orec}(0, \z)$; from
$\mathrm{orec}(j, a)$, $\mathrm{orec}(j{+}1, \suc a)$; from
$\mathrm{orec}(m, a)$ and $\mathrm{orec}(m, b)$, $\mathrm{orec}(1, a =
b)$, and from $\mathrm{orec}(m, a)$, $\mathrm{orec}(m', b)$ and $m \ne
m'$, $\mathrm{orec}(0, a = b)$; the strong Kleene rules for $\neg$ and
$\lor$; from a head step $c \leadsto d$ and $\mathrm{orec}(k, d)$,
$\mathrm{orec}(k, c)$; the recursor rule, which records the value of
an applied recursor through the value of its counter; from
$\mathrm{orec}(1, u\cdot\tofn n)$ for every $n$, $\mathrm{orec}(1,
\cA\cdot u)$; and from $\mathrm{orec}(0, u\cdot\tofn n)$ for one $n$,
$\mathrm{orec}(0, \cA\cdot u)$.
\end{defi}

\noindent
The two quantifier rules read the numeral instances and ask for no
certificate.  The internal record $\Rcd_v(c)$ (\textsf{iRec}) is the
diagonal self-application of a record body with one row per rule of
$\mathrm{orec}$.  Each row is guarded by decided reads of the code
(operator head, head reduct, and the length and head of the spine,
where the \emph{spine} of an application is its head together with its
arguments), and each premise
is a bounded literal: the record of the sub-code at some fuel at most
the current fuel's predecessor.  The fuel-bounded core
$\Rcd^s_v(c)$ (\textsf{iRecAt}) is the disjunction of the rows at fuel
$s$, with no clause at fuel $0$.  The full record body adds an
exclusivity conjunct, which refutes the record at every other value.
The two quantifier rows read the numeral families
\[
  \mathrm{famT}(u) \ldef \forall n.\ \Rcd_1(\mathit{ap}(u, \Num\cdot n)) ,
  \qquad
  \mathrm{famF}(u) \ldef \exists n.\ \Rcd_0(\mathit{ap}(u, \Num\cdot n))
\]
outside the fuel binder (\textsf{famT}, \textsf{famF}), and
$\mathrm{famDec}(u)$ says that both families are decided.

\begin{thm}[\textsf{wosr\_orec}, \textsf{iRec\_of\_wosr}, \textsf{iRec\_faithful}, \textsf{orec\_iff\_wosr}]
\label{thm:c-orec}
For every term $t$, natural $r$ and fuel $s$:
\begin{enumerate}
\item if $t \wbev{s} r$ and $t$ is well-formed, then $\mathrm{orec}(r, t)$;
\item if $t \wbev{s} r$, and $t$ is well-formed and $\forall$-free, then
  $\emptyset \vdash \Rcd_{\tofn r}(\bquo{t})$;
\item if $\emptyset \vdash \Rcd_{\tofn r}(\bquo{t})$ and $t$ is
  $\forall$-free, then $\mathrm{orec}(r, t)$;
\item if $t$ is well-formed and $\forall$-free, then $\mathrm{orec}(r,
  t)$ holds exactly when $t \wbev{s'} r$ for some $s'$.
\end{enumerate}
\end{thm}

\noindent
Part~1 is a simultaneous induction over the witnessed evaluation and
head evaluation, so grounded truth is contained in the record at every
well-formed term.  Part~2 produces the record from the HOL relation,
row by row; part~3 reads the record back.  At a quantifier the record can hold
where grounded truth does not, since the record asks for every instance and no
certificate.

\begin{thm}[\textsf{vfalse\_root3}]
\label{thm:c-false-root}
If $\Gamma \vdash s \jnat$ and $\Gamma$ derives the record body at
value $1$, fuel $s$ and the code $\bquo{\z = \tofn 1}$, then $\Gamma$
derives every well-formed term.
\end{thm}

\noindent
\emph{Proof.}  The only row that fits an equation node is the equality
row.  At value $1$ the equality row needs a common value of the two sides.  The zero
side reads $0$ and the successor side reads a non-zero value, a
numeral disequality.  The two bounded literals are opened at fresh fuel
variables, which occur only in the pair slot; no machine and no fuel
arithmetic is involved.  \qed

\subsection{The record fold}
\label{app:fold-record}

\paragraph{Paths, environments and the guard.}
A proof list cites earlier entries, and one entry can be cited along
several paths from the root.  The record fold walks the tree that the
citation structure unfolds into, in post-order
(\texttt{RGA\_Record/Path\_Walk}).  For a proof code $q$ and a rank $r$,
$\mathit{Node}(q, r)$ is the code of the suffix list whose head is the
node at rank $r$, and $\mathit{Node}(q, r) \ne \z$ says that the rank is
occupied (a decided test).  A kid's rank is strictly below its
parent's, so the engine of Theorem~\ref{thm:c-cov} on ranks reaches
every premise of a node.  For an \emph{eigen tuple} $nb$, a coded list
of naturals indexed like the path, $\mathit{Env}(q, r, nb)$
(\textsf{iEnvAt}) is the path environment at $r$.  Walking down, a
$\forall$-introduction or negated-universal elimination node shifts the
environment and extends it with $\Num\cdot n$ for the tuple's entry
$n$; an induction node updates the environment with the stage, and a substitution node
with the substituted term's image; every other node leaves
the environment unchanged.  Each of the four case-dependent (rule, premise) pairs
(negated-universal elimination, induction, and the two branches of
disjunction elimination) adds one item to the path's guard list: the
added hypothesis and the environment where it is added.  For a fuel
tuple $sb$, $\mathit{Guard}(q, r, nb, sb)$ (\textsf{iGuardAt}) is the
bounded universal over the guard list of the fuel-bounded value-1
records $\Rcd^{s_j}_1$ of the added hypotheses' images, with fuels
$s_j$ from $sb$.

\begin{defi}[\textsf{fMRb}, \textsf{fMot}]
\label{def:c-recmot}
The motive of the record fold at the proof code $q$ and the rank $r$ is
\[
  \mathit{Node}(q, r) \ne \z \;\limp\; \forall nb.\ \forall sb.\
  \neg\mathit{Guard}(q, r, nb, sb) \lor
  \Rcd_1\bigl(\Subst(\mathit{Env}(q, r, nb),
     \mathit{concl}(\mathit{head}(\mathit{Node}(q, r))))\bigr) .
\]
\end{defi}

\noindent
The guard enters as a disjunct, not as the antecedent of an
implication, because the guard is not decided in general: \rga does not derive
the guard decided at every one of the codes $P_x$ of
Theorem~\ref{thm:con-guard-not} (\textsf{gx\_guard\_dec\_absurd}).  The
occupancy test in front is decided.

\paragraph{The arm shape and the generic fold.}
The fold is stated at an abstract \emph{residual}: a formula former
$R$, of which the arms may read $R(q)$ from the context.  An \emph{arm}
for row $t$ of the node table (\textsf{armT}) is a HOL statement: in
every well-formed context $\Delta$ that derives the root facts ($q
\jnat$, $\mathit{TreeValid}(q)$, an empty root context, $R(q)$), the
node facts at the rank $r$ (the node is a successor position $\suc k$,
accepted, on which row $t$ fires, with the row's kid positions), the
accumulator of the induction, $nb \jnat$ and $sb \jnat$, the context
$\Delta$ derives the motive's body at $r$, $nb$ and $sb$.
$\mathit{arms}(R)$ (\textsf{armsT}) is the conjunction of the 41 arms.

\begin{thm}[\textsf{fold\_root\_record}]
\label{thm:c-fold-root}
Let $R$ be a residual with $\mathit{arms}(R)$.  Let $\Gamma$ be a
context and $p$ a term, with $v_1$ free in neither.  If $\Gamma$
derives $p \jnat$, $\mathit{TreeValid}(p)$, $p \ne \z$,
$|\mathit{hyps}(\mathit{head}(p))| = \z$ and $R(p)$, then $\Gamma
\vdash \Rcd_1(\Subst(\z, \mathit{concl}(\mathit{head}(p))))$.
\end{thm}

\noindent
\emph{Proof.}  Theorem~\ref{thm:c-cov} at the root rank, with every
step by \textsf{fold\_step}: a split on the decided occupancy test,
the two universals at fresh variables, and the 41-way dispatch on the
node's first firing row, each branch closed by its arm.  At the root
the rank is occupied, the environment is empty and the guard list is
empty, so the guard holds and the disjunction yields the record.  \qed

\begin{thm}[\textsf{fold\_consistent\_ctx}, \textsf{fold\_sound\_ctx}]
\label{thm:c-fold-con}
Let $R$ be a residual with $\mathit{arms}(R)$, $\Gamma$ a context and
$p$ a term, with $v_1$ free in neither, and let $\Gamma \vdash p \jnat$
and $\Gamma \vdash R(p)$.  Then $\Gamma \vdash \neg\Chk(p, \bquo{\z =
\tofn 1})$.  If moreover $\Gamma$ is well-formed, then $\Gamma \vdash
\forall c.\ \Chk(\lft p, c) \limp \Rcd_1(\Subst(\z, c))$.
\end{thm}

\noindent
\emph{Proof.}  The decided-matrix split on $\Chk(p, \bquo{\z = \tofn
1})$ (Theorem~\ref{thm:in-chk-dec}).  In the accepting case,
Theorem~\ref{thm:c-fold-root} and transport along the conclusion
equation give $\Rcd_1(\Subst(\z, \bquo{\z = \tofn 1}))$; opening the
record's fuel at a fresh variable and applying
Theorem~\ref{thm:c-false-root} derives the goal.  The soundness form
splits on $\Chk(\lft p, c)$ at a fresh $c$ and introduces the universal.
\qed

\paragraph{The 41 arms.}
Table~\ref{tab:c-rec-arms} sorts the arms by what they need from the
residual.

\begin{table}[ht]
\centering
\footnotesize
\begin{tabular}{@{}>{\raggedright\arraybackslash}p{0.2\linewidth}>{\raggedright\arraybackslash}p{0.5\linewidth}>{\raggedright\arraybackslash}p{0.24\linewidth}@{}}
\hline
group & rules (the Isabelle arm is \textsf{arm\_}$\langle$Name$\rangle$ with the first letter capitalized: \textsf{trueEqI} $\mapsto$ \textsf{arm\_TrueEqI}) & theory \\
\hline
any residual (34) & the zero rule, weakening, cut; 24 citation rules (the equality,
  successor, negation and disjunction rules: \textsf{trueEqI}, \textsf{trueEqE}, \textsf{falseEqI},
  \textsf{falseEqE}, \textsf{negE}, \textsf{disjI1}, \textsf{disjI2}, \textsf{eqRefl}, \textsf{eqSym},
  \textsf{neSym}, \textsf{eqTI}, \textsf{eqTE}, \textsf{sucNe0}, \textsf{sucNeI}, \textsf{sucTI},
  \textsf{sucEqE}, \textsf{sucNeE}, \textsf{sucTE}, \textsf{botTE}, \textsf{dnegI}, \textsf{dnegE},
  \textsf{disjI3}, \textsf{disjE2}, \textsf{disjE3}); the context-changing substitution,
  disjunction-elimination and negated-universal-elimination rules; the four conversion rules &
  \texttt{ConRGA/Full\_Fold} \\
guard interface (2) & hypothesis, $\forall$-introduction & \texttt{ConRGA/Wfd\_Decision} \\
swap interface (5) & the recursor's successor step, universal elimination, negated-universal
  introduction, induction, equality substitution & \texttt{ConRGA/Wfd\_Static} \\
\hline
\end{tabular}
\caption{The 41 arms of the record fold.  The 34 arms of the first
  group hold at every residual (\textsf{armsT\_of\_fixed}).}
\label{tab:c-rec-arms}
\end{table}

\noindent
The seven arms of the last two groups consume the residual through two
interfaces, and a third, minor one.  The interfaces are HOL-level
properties of $R$, statements about derivations, not \rga formulas.

\begin{defi}[\textsf{gdIf}, \textsf{swIfG}, \textsf{xnIf}, \textsf{xgIf}, \textsf{xqIf}]
\label{def:c-interfaces}
Let $R$ be a residual.  In each clause, $\Delta$ ranges over
well-formed contexts that derive $R(q)$, $q \jnat$, $r \jnat$,
$nb \jnat$, $\mathit{TreeValid}(q)$ and $\mathit{Node}(q, r) \ne \z$.
\begin{itemize}
\item \emph{Guard interface} $\mathrm{gd}(R)$: for every such
  $\Delta$ and every $sb$ with $\Delta \vdash sb \jnat$, $\Delta \vdash
  \mathit{Guard}(q, r, nb, sb) \jbool$.
\item \emph{Swap interface} $\mathrm{sw}_G(R)$, for a guard $G$ on the
  swapped term: at every such $\Delta$ where the node has one of the
  five swap tags, with critical formula $C$ (\textsf{critAt}) and the
  arm's environment $Z$, and for all codes $w, X, Y, v, k, s$ that
  $\Delta$ proves natural, if $X$ and $Y$ are well-formed term codes,
  $\Delta \vdash \Rcd_v(X)$, $\Delta \vdash \Rcd_v(Y)$ and $X$ passes
  $G$, then $\Delta \vdash \Rcd^s_k(\Subst(Z[w{:=}X], C))$ implies
  $\Delta \vdash \Rcd_k(\Subst(Z[w{:=}Y], C))$.
\item $G$ passes every numeral code $\Num\cdot n$ (\textsf{xnIf}), is
  stable under weakening, and, wherever $R$ holds, passes the path
  image of the node's swapped term at the recursor's successor step,
  negated-universal introduction and equality substitution
  (\textsf{xgIf}).
\item \emph{Swapped-term interface} (\textsf{xqIf}): wherever $R$
  holds, the path image of the swapped term at those three rules is
  $\forall$-free.
\end{itemize}
\end{defi}

\begin{thm}[\textsf{fold\_consistent\_pathG}, \textsf{fold\_sound\_pathG}]
\label{thm:c-fold-pathG}
Let $R$ be a residual such that $v_1$ is not free in $R(\vz)$ and
$R(\vz)$ is well-formed, and let $G$ be a guard.  If $\mathrm{gd}(R)$
and $\mathrm{sw}_G(R)$ hold, and $G$ satisfies the conditions of the
third item of Definition~\ref{def:c-interfaces}, then
\[
  R(\vz),\ \vz\jnat \vdash \neg\Chk(\vz, \bquo{\z = \tofn 1}) .
\]
The soundness form holds under the hypotheses of
Theorem~\ref{thm:c-fold-con}, with these interfaces in place of
$\mathit{arms}(R)$.
\end{thm}

\noindent
\emph{Proof.}  The two guard arms (\textsf{arm\_Hyp\_path},
\textsf{arm\_AllI\_path}) and the five swap arms
(\textsf{arm\_RSuc\_pathG} and its siblings) hold at $R$, and
\textsf{armsT\_of\_fixed} completes $\mathit{arms}(R)$.
Theorem~\ref{thm:c-fold-con} finishes.  \qed

\noindent
Every consistency and soundness theorem of Section~\ref{sec:con}
other than the $\forall$-free ones is an instance of Theorem~\ref{thm:c-fold-pathG}.  The fragment
takes $R = \Frag$ (\S\ref{app:fold-frag}); the residual theorem takes
$R = \Res$ and the trivial guard (\S\ref{app:fold-res}); the level
theorems take a decided scan and a context hypothesis through the
shell \textsf{con\_shell} (Appendix~\ref{app:oga}).  The two exhibited
arms follow.

\paragraph{A guard arm: $\forall$-introduction.}
The rule concludes $\Gamma \vdash \cA\cdot f$ from $v_0\jnat,\ \lft\Gamma
\vdash \lft f\cdot v_0$.  Let $z$ be the path environment at the node
and $U = \Subst(z, f)$, so that the node's conclusion has image
$\cA\cdot U$.  The arm \textsf{arm\_AllI\_path} proceeds as follows.
\begin{enumerate}
\item By the guard interface, $\mathit{Guard}(q, r, nb, sb)$ is
  decided.  The arm splits on the guard before introducing any variable.  If
  the guard fails, the motive's first disjunct holds.
\item If the guard holds, the arm opens a fresh variable $j$ with $j
  \jnat$.  The kid's path environment is $z$ shifted and extended by
  $\Num\cdot j$ (\textsf{ir\_iEnvGd\_kid}).  The induction hypothesis
  at the kid, with the eigen tuple extended by $j$, gives
  $\neg\mathit{Guard}_{\mathrm{kid}} \lor \Rcd_1(\cdots)$.
  $\forall$-introduction is not case-dependent, so the kid's guard is
  the node's (\textsf{ir\_iGuardAt\_same\_negD}), and the first
  disjunct contradicts the case hypothesis.  The second disjunct is
  the record of the kid's conclusion, whose image is $\mathit{ap}(U,
  \Num\cdot j)$.
\item Universal introduction at $j$ gives $\mathrm{famT}(U)$, and the
  universal row of the record (\textsf{ir\_iRecP\_allTI}) gives
  $\Rcd_1(\mathit{ap}(\bquo{\cA}, U))$, which is the image of the
  node's conclusion.
\end{enumerate}
The split in step~1 must come before the universal introduction in
step~3.  Without a decided guard the arm could not commute the case
analysis with the eigenvariable, and that commutation is where a quantified case
hypothesis on the path stops the fold.  No formula that \rga derives
at the code of every valid proof $P_x$ of
Theorem~\ref{thm:con-guard-not} supplies the guard interface (\textsf{gdIf\_not\_of\_gx}).  Each $P_x$ is one
disjunction elimination whose case hypothesis is $\forall n.\
\mathit{ofam}(x)\cdot n$, ``program $x$ never halts'', and a guard
interface for all the proofs $P_x$ would make the complement of the halting
set recursively enumerable (\textsf{gx\_guard\_dec\_absurd},
\S\ref{app:fold-res}).

\paragraph{A swap arm: universal elimination.}
The rule concludes $\Gamma \vdash f\cdot a$ from $\Gamma \vdash
\cA\cdot f$ and $\Gamma \vdash a \jnat$.  With $z$, $U$ as above and
$A = \Subst(z, a)$, the cited premises give $\Rcd_1(\mathit{ap}(\bquo{\cA},
U))$ and $\Rcd_1(A = A)$.  The arm \textsf{arm\_AllE\_pathG}:
\begin{enumerate}
\item opens the equality row of the second record, which gives a
  value $m$ with $\Rcd_m(A)$;
\item opens the universal row of the first record at $m$, which gives
  $\Rcd_1(\mathit{ap}(U, \Num\cdot m))$ and $\Rcd_m(\Num\cdot m)$
  (\textsf{iRec\_allE\_inst\_num});
\item writes $\mathit{ap}(U, \Num\cdot m)$ as the image of the critical
  formula $C$, the code of $\lft f\cdot\vz$, under the arm's environment
  updated with $\Num\cdot m$, and opens the record's fuel to obtain the
  core $\Rcd^{s'}_1$ there;
\item applies the swap interface with $X = \Num\cdot m$ and $Y = A$:
  both are well-formed term codes with value $m$, and $X$ passes the
  guard because it is a numeral code.  The result is $\Rcd_1$ of the
  image of $C$ at $A$, which is $\mathit{ap}(U, A)$, the image of the
  node's conclusion.
\end{enumerate}
The universal row records $\cA\cdot U$ through numeral instances only,
and the node's conclusion is an instance at an arbitrary term.  The
swap bridges the two kinds of instance, and, informally (\S\ref{sec:con-break}; no
theorem locates this failure), the swap is the second place where the fold
stops: at a quantified critical formula the swap engine would need the
numeral families of its quantified subterms.
Stated as one formula, as the residual states it, the swap is not
derivable even at the $\forall$-free critical formula of the proof
$P_3$ of \S\ref{sec:con-frag}
(Theorem~\ref{thm:con-res-not-frag}).

\subsection{The fragment}
\label{app:fold-frag}

\begin{defi}[\textsf{iFrag}]
\label{def:c-frag}
$\Frag(p) \ldef \mathit{Static}(p) \wedge \mathit{FragV}(p)$, two decided
scans of the proof list coded by $p$ (\textsf{iStat},
\textsf{iFragV}).  Together the two scans require, at every position of the
list:
\begin{itemize}
\item at a substitution node, the substituted term that the checker's
  witness search returns is $\forall$-free;
\item at a negated-universal elimination, induction or disjunction
  elimination node, every hypothesis that the rule adds to a premise is
  $\forall$-free;
\item at a node of the five swap rules, the critical formula is
  $\forall$-free, well-formed and the code of a term;
\item at the recursor's successor step, negated-universal introduction
  and equality substitution, the swapped term is $\forall$-free.
\end{itemize}
\end{defi}

\noindent
$\Frag(p)$ is decided at every $p$ that a well-formed context proves
natural (\textsf{hq\_iFrag}).

\paragraph{The substitution clause.}
The checker does not read a substitution node's variable and term off
the proof.  The checker's detector for the substitution rule searches, below a
bound computed from the codes, for a variable and a term that make the
node's judgment the substitution instance of the cited premise, and it
returns the first hit in a fixed search order (\textsf{iSub1Ns}, built
on the first-hit fold \textsf{fhLtF}).  $\Frag$ tests that returned
term.  Several pairs of a variable and a term can produce the same
judgment, for instance when the variable does not occur in the
premise, so the returned term need not be the one the proof used.
The code of every accepted $\forall$-free proof passes the clause
nevertheless, through the minimality of the search, and \rga proves
this containment at an abstract code (\S\ref{app:fold-contain}).  At closed codes
$\Frag$ is thus a decided scan of the nodes: it contains the code of
every accepted $\forall$-free proof, and at each node it reads the
critical formula and the case-dependent hypotheses of the rule instance
that the checker's first-hit searches cite.  When two rule instances
justify the same node, membership can depend on which one the checker
finds first, and the paper gives no HOL characterization of
$\Frag(\bquo{P})$ that is independent of this choice.  HOL
evaluates the class at the codes of valid $\forall$-free proofs, at
$P_{\mathrm w}$ and $P_3$ (\S\ref{sec:con-frag}), at the
incomparability witnesses of Appendix~\ref{app:oga-frag}, and at the
proofs $P_x$ of Theorem~\ref{thm:con-guard-not}.

\paragraph{Discharging the interfaces.}
Lemma~Q (Theorem~\ref{thm:con-lemq}) is one internal induction on the
fuel (\textsf{lemQ\_level}).  The induction's step decides the core at fuel
$s{+}1$ from the decided bounded records of the kids at fuel $s$
(\textsf{core\_B\_step}).  At a $\forall$-free code no quantifier row
fires, and three laws keep the kids $\forall$-free: the head step, the
recursor spine's counter, and its dispatch at every value preserve
$\forall$-freeness (\textsf{qfScan\_gHR}, \textsf{qfScan\_gSpArg3},
\textsf{qfScan\_rDispS}).  The three interfaces then hold at $\Frag$
with no further hypothesis.
\begin{itemize}\raggedright
\item \textsf{gdIf\_F}: under $\Frag$ every path environment is
  $\forall$-free (\textsf{qfEnv\_iEnvAt}), so every guard item is the
  core of a $\forall$-free code, decided by Lemma~Q, and the bounded
  universal over the items is decided (\textsf{whole\_guard\_Bq}).
\item \textsf{swIfQ\_F}: the critical formula is a $\forall$-free,
  well-formed term code, and the two environments are related by the
  compiled value-congruence checker, so the swap engine applies
  directly (\textsf{swTopQ}).
\item \textsf{xqIf\_F}: the swapped-term clause of $\Frag$ and the
  $\forall$-free path environment.
\end{itemize}

\begin{thm}[\textsf{rec\_swap\_qf}]
\label{thm:c-swap}
Let $\Delta$ be a well-formed context and $X, Y, v, u, u', k, s$ codes
that $\Delta$ proves natural.  Suppose that $X$ and $Y$ are well-formed
term codes with $\Delta \vdash \Rcd_v(X)$ and $\Delta \vdash
\Rcd_v(Y)$; that $u$ is $\forall$-free and $u$ and $u'$ are well-formed
term codes; and that the compiled value-congruence checker relates $u$
and $u'$ at $X$ and $Y$, so that $u'$ is $u$ with some sites $X$
replaced by $Y$.  If $\Delta \vdash \Rcd^s_k(u)$, then $\Delta \vdash
\Rcd_k(u')$.
\end{thm}

\noindent
\emph{Proof sketch.}  One internal induction on the fuel over related
pairs (\textsf{swap\_level}).  Every constructor row is a direct step.
The application row goes through two simulation faces, each one
internal induction: a head step of the $\forall$-free source is matched
by a head step of the target or, at a recursor redex whose counter is a
literal numeral in the source only, by a recursor spine of the target;
and related recursor spines have related counters and dispatches.
Neither $X$ nor $Y$ needs to be $\forall$-free.  \qed

\noindent
Theorem~\ref{thm:c-fold-pathG} at $R = \Frag$ (through its form at
path points, \textsf{fold\_consistent\_swIfQP}) gives $\Frag(\vz),\ \vz\jnat
\vdash \neg\Chk(\vz, \bquo{\z = \tofn 1})$.  Discharging the decided
$\Frag(\vz)$ (\textsf{ir\_impI\_bool}, \textsf{hq\_iFrag}) and
introducing the universal gives Theorem~\ref{thm:con-frag}(1); the
soundness form gives Theorem~\ref{thm:con-frag}(2).  A third form
has both universals outermost:

\begin{thm}[\textsf{rga\_sound\_frag\_int\_prenex}]
\label{thm:c-frag-prenex}
$\emptyset \vdash \forall p.\ \forall c.\ \Frag(p) \limp \bigl(\Chk(p,
c) \limp \Rcd_1(\Subst(\z, c))\bigr)$.
\end{thm}

\paragraph{The search reading.}
Let $\mathit{conSearch}$ be the term that searches the numerals for a
code that $\Chk$ accepts as a proof of $\z = \tofn 1$, and
$\mathit{fragSearch}$ the same search restricted to $\Frag$.  If a
well-formed context $\Gamma$ derives $v \jnat$ and $\Chk(v, \bquo{\z =
\tofn 1})$, then $\Gamma \vdash \Rcd_1(\bquo{\mathit{conSearch}})$
(\textsf{con\_search\_found}).

\begin{thm}[\textsf{frag\_search\_diverges}, \textsf{rga\_con\_int\_of\_con\_search}, \textsf{con\_search\_unrefutable}]
\label{thm:c-search}
\leavevmode
\begin{enumerate}
\item For every natural $v$, $\emptyset \vdash
  \neg\Rcd_{\tofn v}(\bquo{\mathit{fragSearch}})$.
\item If $\emptyset \vdash \neg\Rcd_1(\bquo{\mathit{conSearch}})$, then
  $\emptyset \vdash \ConR$.
\item $\emptyset \nvdash \neg\Rcd_1(\bquo{\mathit{conSearch}})$.
\end{enumerate}
\end{thm}

\noindent
Part~1 runs an internal invariant of the search, and
Theorem~\ref{thm:con-frag} closes the invariant at the search's test.
Part~3 is part~2 with G\"odel's second theorem
(\textsf{rga\_goedel2}).

\paragraph{The Liar.}
Let $L$ be the companion paper's Liar, the fixed point of
$\lambda\neg\vz$.  For every well-formed context $\Gamma$ and every
$v$ with $\Gamma \vdash v \jnat$, $\Gamma \vdash \neg\Rcd_v(\bquo{L})$
(\textsf{liar\_iRec\_none}).  The proof refutes the core at every fuel
by an internal induction on the fuel and then refutes the record at
each value.

\subsection{The fragment contains the accepted \texorpdfstring{$\forall$}{A}-free codes}
\label{app:fold-contain}

This subsection proves Theorems~\ref{thm:con-qf-in-frag},
\ref{thm:con-qf-via-frag} and~\ref{thm:con-witness}, from the theories
\texttt{ConOGA/ConOGA\_Frag\_Witness}, \texttt{ConOGA/ConOGA\_Frag\_QF}
and \texttt{ConOGA/ConOGA\_Frag\_QF\_Int}.  The containment is proved
twice.  The internal proof introduces the scans of $\Frag$ inside \rga
at an abstract code.  The containment's numeral instances also have a HOL proof, by
evaluating the scans at a numeral, the method that also proves the
witness's membership.  For a natural $p$
and a position $m$, write $\mathit{drop}_m(p)$ for the code of the
suffix of the list coded by $p$ that starts at position $m$; the suffix's head
is the node at $m$.

\paragraph{The rule tag.}
The checker reads the node at the head of a list through its 41-row
node table.  The node's \emph{rule tag} is the first row whose detector
fires at the node's context and conclusion against the tail of the
list, or $0$ when no row fires.

\paragraph{The containment inside \rga.}
The internal proof is not an induction over the proof, and it does not
read the code back.  The core of that proof holds at an open code in every
well-formed context.

\begin{thm}[\textsf{qf\_tree\_in\_frag}]
\label{thm:c-qf-tree}
Let $\Gamma$ be a well-formed context and $p$ a term with $\Gamma
\vdash p\jnat$.  If $\Gamma \vdash \mathit{TreeValid}(p)$ and $\Gamma
\vdash \QfScan(p)$, then $\Gamma \vdash \Frag(p)$.
\end{thm}

\noindent
\emph{Proof.}  $\Frag(p)$ is the conjunction of $\mathit{Static}(p)$,
a bounded universal over the positions of the list coded by $p$, and
$\mathit{FragV}(p)$, a negated bounded search for a bad node.  \rga
introduces each conjunct from the node check at an arbitrary position below the
length, in every extension $\Delta$ of $\Gamma$ that makes the position
a natural below the length (\textsf{iAllLt\_allI}, by an internal
induction on the bound, and \textsf{nodeScan\_I}).  At such a position
$m$, $\mathit{drop}_m(p) = \suc k$ for a natural $k$, and the validity
test and the $\forall$-free scan both pass at $\suc k$
(\textsf{qf\_node\_at}).  At an accepted node \rga derives that some
row $i < 41$ fires and that the rule tag of the node is $\tofn i$
(\textsf{iTreeValid\_node\_row}, \textsf{nodeTbl\_tag}).  The proof
splits over the 41 values of $i$; at each value, the row's own detector fires
(\textsf{tag\_fire}).  The clauses of the two node checks leave
obligations only at the eight constrained rows.
\begin{itemize}
\item Universal elimination and negated-universal introduction and
  elimination: the row's detector is refuted at a $\forall$-free node
  (\textsf{allE\_node\_refuted}, \textsf{notAllI\_node\_refuted},
  \textsf{notAllE\_node\_refuted}).
\item The recursor's successor step: the node's $\forall$-free
  conclusion has the successor-step shape, so the critical formula,
  built from the conclusion's subterms and their lift, is $\forall$-free,
  well-formed and a code, and the swapped term is $\forall$-free
  (\textsf{rsuc\_c24k\_qf}).
\item Induction: the stage formula is a hypothesis of the
  $\forall$-free step premise, hence $\forall$-free
  (\textsf{natInd\_nc\_qf}).
\item The first disjunction elimination: both added hypotheses are
  disjuncts of a $\forall$-free cited conclusion
  (\textsf{disjE1\_slots\_qf}).
\item Equality substitution: the second cited conclusion is the image
  of the template under a substitution and is $\forall$-free, and the
  template is well-formed, so Theorem~\ref{thm:c-sub1-laws}(1) makes the
  template $\forall$-free (\textsf{eqsubs\_c40k\_qf}).
\item Substitution: the returned term is $\forall$-free
  (\textsf{sub1\_ns\_qf}, below).
\end{itemize}
Where $\mathit{Static}$ asks for a static term, a $\forall$-free term
passes (\textsf{iStatR\_of\_qf}).  \qed

\noindent
Two facts about the compiled substitution carry the last two rows.
Write $\mathit{Sub}_1(v, u, c)$ for the compiled image of
$\ssub{v}{u}{c}$ at codes (\textsf{iSub1S}).

\begin{thm}[\textsf{iSub1S\_qf\_reflect}, \textsf{iSub1S\_nfree\_indep}]
\label{thm:c-sub1-laws}
Let $\Gamma$ derive $v\jnat$, $u\jnat$, $u'\jnat$, $c\jnat$ and
$\Wfc(c) \ne \z$.
\begin{enumerate}
\item If $\Gamma$ derives that $\mathit{Sub}_1(v, u, c)$ is
  $\forall$-free, then $\Gamma$ derives that $c$ is $\forall$-free.
\item If $\Gamma$ derives that $v$ is not free in $c$, then $\Gamma
  \vdash \mathit{Sub}_1(v, u, c) = \mathit{Sub}_1(v, u', c)$.
\end{enumerate}
\end{thm}

\noindent
\emph{Proof.}  Each part is a course-of-values induction inside \rga
(\textsf{cov\_pair\_ind}) on a packed pair of the substitution's
configuration and the code, with a motive guarded by decided tests:
$c$ well-formed and its image $\forall$-free for~(1), $c$ well-formed
and $v$ not free in $c$ for~(2).  At a variable, (1) holds outright,
and in~(2) the variable is not $v$, so both images are the variable.
At an operator node with zero, one or two kids the image is the node
with the same operator over the kids' images
(\textsf{subm1\_opr\_read\_gen0} and its siblings).  For~(1) the
image's $\forall$-freeness gives the operator's and the kid images',
the induction hypothesis gives the kids, and the node is rebuilt.
For~(2) the kids meet the guard at the shifted variable, and the
induction hypothesis gives equal kid images.  Well-formedness refutes
three or more kids.  (2) needs no hypothesis that $c$ is the code of
a term: the proof of~(2) compares two images of one machine and never uses
$\mathit{Sub}_1(v, u, c) = c$.  \qed

\noindent
Part~2 lifts to the coded hypothesis sets of a judgment by an
induction on the fold that builds their images
(\textsf{iSub1Fset\_indep}).  At a substitution node
(\textsf{sub1\_ns\_qf}) the term that $\mathit{Static}$ reads is the
least witness $\mathit{ns}$ of the detector for the substituted
variable $x$, and \rga splits on whether $\mathit{ns}$ is
$\forall$-free.  If $\mathit{ns}$ is not, $x$ is free neither in the cited
premise's conclusion nor in any of the premise's hypotheses, since the images of that conclusion and those hypotheses
are the node's $\forall$-free judgment (\textsf{qfRevScan\_at0}).
Part~2 and its lift then make the images at the substituend $\z$ equal
to the node's judgment, so the witness check fires at $\z$
(\textsf{sub1Wit\_fire\_at0}), the least witness is $\z$
(\textsf{iFhLt\_hit0}), and the code $\z$ of $\vz$ is $\forall$-free:
a contradiction.

\noindent
\emph{Proof of Theorem~\ref{thm:con-qf-in-frag}}
(\textsf{rga\_qf\_in\_frag\_int}).  In the context of $p\jnat$, $c\jnat$
and the hypothesis $\ChkQf(p, c)$, whose conjuncts include
$\mathit{TreeValid}(p)$ and $\QfScan(p)$, Theorem~\ref{thm:c-qf-tree}
gives $\Frag(p)$.  The verdict $\ChkQf(p, c)$ is decided
(\textsf{hq\_iChkQf}), so the implication is introduced over the verdict
(\textsf{ir\_impI\_bool}), and two $\forall$-introductions give the
sentence.  The head-context conjunct and the conclusion $c$ are not
used.  \qed

\noindent
\emph{Proof of Theorem~\ref{thm:con-qf-via-frag}}
(\textsf{rga\_con\_qf\_via\_frag\_int}).  Split on the decided verdict
$\ChkQf(p, \bquo{\z = \tofn 1})$ at a fresh $p$ with $p\jnat$.  In the
accepting case, universal elimination of
Theorem~\ref{thm:con-qf-in-frag} at $p$ and $\bquo{\z = \tofn 1}$ and
modus ponens give $\Frag(p)$, Theorem~\ref{thm:con-frag} at $p$ gives
$\neg\Chk(p, \bquo{\z = \tofn 1})$, and \textsf{iChk\_of\_iChkQf}
gives $\Chk(p, \bquo{\z = \tofn 1})$.  The universal follows by
$\forall$-introduction.  \qed

\paragraph{The numeral instances, in HOL.}
For each natural $p$ and term $c$: if $\emptyset \vdash \ChkQf(\tofn p,
c)$, then $\emptyset \vdash \Frag(\tofn p)$
(\textsf{qf\_accepted\_in\_frag}).  The HOL proof of these instances
reads the accepted code back and uses soundness
(Appendix~\ref{app:machines-audit}).  The same instances follow
without soundness from Theorem~\ref{thm:c-qf-tree} in the empty context
(\textsf{qf\_accepted\_in\_frag\_tree}).  The HOL proof that reads the code back rests on
the rule tag, which HOL evaluates at closed codes.

\begin{thm}[\textsf{ruleTag\_code\_cases}]
\label{thm:c-ruletag}
Let $H$ be a context, $c$ a term and $Q$ a proof list, and let $t$ be
the HOL value of the rule tag at the code of the list $(H, c)$
followed by $Q$.  Either $t = 0$ and no detector fires at
$(\quo{Q}, \quo{H}, \quo{c})$, or $t < 41$, the detector of row $t$
fires there, and no detector of a row below $t$ fires.
\end{thm}

\noindent
\emph{Proof.}  The tag is a chain of conditionals over the 41 rows,
one test per row, whose value is $0$ when no test fires and otherwise
the value attached to the first test that fires
(\textsf{pr2fi\_ifNestIdx\_cases}, one list induction).  Row $i$'s
value carries the tag $i$, and its test is detector $i$ at the tail,
the head context and the head conclusion.  \qed

\noindent
The readers that locate a node's cited premises, its induction stage,
its template or its substitution witness are first-hit searches.  The
HOL value of such a search is the successor of the least hit below the search bound, or
$0$ with no hit (\textsf{pr2fi\_fhLtF\_least}).

\paragraph{Closing a scan at a closed code.}
The node scans $\mathit{FragV}$, $\mathit{StatLv}_N$ and
$\mathit{StatLv}$ are negated bounded searches for a bad node.  One lemma closes such a scan from
facts about each node.

\begin{thm}[\textsf{nodeScan\_closedI}]
\label{thm:c-nodescan}
Let $K$ be a node check, a compiled per-node test with the
substitution, well-formedness and naturality laws of the checker's
tests, and let $p$ be a natural number.  If $\emptyset \vdash K(\mathit{drop}_m(\tofn p)) \ne \z$
for every position $m$ below the length of the list coded by $p$, then
\rga derives, in the empty context, the node scan of $K$ at $\tofn p$.
\end{thm}

\noindent
\emph{Proof.}  The scan is decided at $\tofn p$
(\textsf{hq\_nodeScan}), so by the disjunction property \rga derives
the scan or derives the scan's negation.  In the second case double-negation
elimination gives the bounded search for a bad node, and the search's closed
witness gives a numeral $m$ at which \rga derives both conjuncts of a
bad node.  The first conjunct evaluates to $m$ being below the length,
and the second, $K(\mathit{drop}_m(\tofn p)) = \z$, contradicts the
hypothesis at $m$, against \rga's consistency.  \qed

\noindent
The compiled scans $\mathit{Static}$, $\mathit{HypLv}_N$,
$\mathit{HypLv}$ and $\mathit{Sub1Qf}$ are bounded universals over the
positions, and close at a closed code by a HOL induction on the bound
(\textsf{iAllLt\_closed\_ne0}).  All seven scans (those of $\Frag$,
$\LvC_N$ and $\LvC$) pass at a code each of
whose nodes has a tag outside the rows that the scan constrains
(\textsf{iStat\_closed\_off} and its siblings).  The constrained rows
are the recursor's successor step, universal elimination,
negated-universal introduction and elimination, substitution,
induction, the first disjunction elimination and equality
substitution.

\paragraph{Every accepted $\forall$-free code.}

\begin{thm}[\textsf{iStat\_afree}, \textsf{iFragV\_afree}]
\label{thm:c-afree-scans}
Let $P$ be a $\forall$-free proof list.
\begin{enumerate}
\item $\Gamma \vdash \mathit{Static}(\bquo{P})$ for every context $\Gamma$.
\item If $P$ is valid, then $\emptyset \vdash \mathit{FragV}(\bquo{P})$.
\end{enumerate}
\end{thm}

\noindent
\emph{Proof.}  Node by node.  At position $m$ the node code is the code
of $(H, c)$ followed by $Q$, where $(H, c)$ is the judgment of $P$ at
$m$ and $Q$ the rest of the list after that judgment.  By
Theorem~\ref{thm:c-ruletag} a nonzero tag fires its row, and the HOL
counterparts of the checker's detectors read the firing row in HOL.
\begin{itemize}
\item Universal elimination, negated-universal introduction and
  negated-universal elimination cannot be the tag
  (\textsf{ruleTag\_afree\_not\_quant}).  The rows of these three rules cite a premise
  $\cA\cdot f$, conclude $\neg(\cA\cdot f)$, or cite $\neg(\cA\cdot f)$,
  and every judgment of $P$ is $\forall$-free.
\item At the recursor's successor step the node's own conclusion is
  $\cR\cdot g\cdot h\cdot \suc m' = h\cdot m'\cdot(\cR\cdot g\cdot h\cdot m')$.
  The critical formula evaluates to the code of a term built from
  subterms of that $\forall$-free conclusion, well-formed because $P$
  is valid, and the swapped term to the code of $m'$.
\item At induction, the first-hit searches return the cited premises
  and a witness of the detector, so the first cited conclusion, the
  base, is the stage formula at $\z$.  The base is $\forall$-free, so
  the stage formula is; the witness's own checks
  make the stage formula well-formed and a code.  The stage formula is both the critical formula and
  the added hypothesis.
\item At the first disjunction elimination the first cited conclusion
  is $a \lor b$, so both added hypotheses $a$ and $b$ are
  $\forall$-free.
\item At equality substitution the template is the least hit of the
  template search, and the second cited conclusion is its instance,
  so the template is $\forall$-free; the swapped term is the left side
  of the cited equation.
\item At substitution (the $\mathit{Static}$ scan only) the returned
  term is the least witness for the variable $x$ that the detector
  found (\textsf{sub1Ns\_node}).  If $x$ is free in the premise's
  judgment, the node's $\forall$-free judgment (a hypothesis or the
  conclusion) contains the substituted term, which is therefore
  $\forall$-free.  Otherwise the code $\z$ of $\vz$ is
  also a witness, since substituting for $x$ then leaves the premise
  unchanged; so the least witness is $\z$, and $\vz$ is
  $\forall$-free.
\end{itemize}
For part~1 these facts are read through the compiled node check of
$\mathit{Static}$ in HOL; $\forall$-free terms are static.  For part~2
the internal node check of $\mathit{FragV}$ is evaluated at the closed
node code: the check's reads become numerals by numeral laws of the first-hit
fold, the suffix, the cited node and the lift, and the freeness,
well-formedness and code tests then evaluate at those numerals
(\textsf{fqOkT\_afree\_node}).  Theorem~\ref{thm:c-nodescan} closes
the scan.  \qed

\noindent
\emph{Proof of the numeral instances.}  From $\emptyset
\vdash \ChkQf(\tofn p, c)$, the conjunct $\mathit{TreeValid}(\tofn p)$
and the genuineness of the validity test give a valid $P$ with $p =
\quo{P}$, and the scan conjunct with Theorem~\ref{thm:c-qfscan} makes
$P$ $\forall$-free.  Theorem~\ref{thm:c-afree-scans} gives both
conjuncts of $\Frag(\tofn p)$.  The head-context conjunct of $\ChkQf$
is not used, and $c$ is arbitrary.  \qed

\noindent
The containment of the accepted $\forall$-free codes in the level
classes of Appendix~\ref{app:oga} is proved in the same two ways: in
the per-level classes at closed codes by reading the code back
(\textsf{qf\_accepted\_in\_lv}), and as one \rga sentence at an
abstract code, for $\LvC$ (\textsf{rga\_qf\_in\_lvi\_int},
Theorem~\ref{thm:con-lv-contain}(1)) and for every $\LvC_N$
(\textsf{rga\_qf\_in\_lv\_int}); Appendix~\ref{app:oga-frag} on the
fragment and the level classes gives the proofs.

\paragraph{The witness.}
Let $f = \lambda(\vz = \vz)$, which is $\lambda x.\ x = x$ in de
Bruijn form.
The proof list $P_{\mathrm w}$ (\textsf{fwP}) has three entries: $\cA\cdot
f$ in the empty context, $f\cdot\vz$ in the context $\vz\jnat$, and
$\vz\jnat$ in the context $\vz\jnat$.  The last entry is a hypothesis,
the middle entry a conversion of the last, and the first entry a
$\forall$-introduction over the middle entry.

\begin{thm}[\textsf{fwP\_accepted}, \textsf{fwP\_in\_frag}, \textsf{fwP\_in\_lvi}, \textsf{fwP\_in\_lv}, \textsf{fwP\_not\_qf}, \textsf{frag\_strictly\_beyond\_qf}]
\label{thm:c-witness}
\leavevmode
\begin{enumerate}
\item $\emptyset \vdash \Chk(\bquo{P_{\mathrm w}}, \bquo{\cA\cdot f})$.
\item $\emptyset \vdash \Frag(\bquo{P_{\mathrm w}})$, $\emptyset
  \vdash \LvC(\bquo{P_{\mathrm w}})$, and $\emptyset
  \vdash \LvC_N(\bquo{P_{\mathrm w}})$ for every $N$.
\item $P_{\mathrm w}$ is not $\forall$-free; $\emptyset \vdash
  \neg\QfScan(\bquo{P_{\mathrm w}})$; and for every term $c$, $\emptyset
  \nvdash \ChkQf(\bquo{P_{\mathrm w}}, c)$.
\item Hence there are naturals $p$ and $c$ with $\emptyset \vdash
  \Chk(\tofn p, \tofn c)$, $\emptyset \vdash \Frag(\tofn p)$ and
  $\emptyset \vdash \LvC_N(\tofn p)$ for every $N$, while $\emptyset
  \nvdash \QfScan(\tofn p)$.
\end{enumerate}
\end{thm}

\noindent
\emph{Proof.}  (1) $P_{\mathrm w}$ proves $\emptyset \vdash \cA\cdot f$,
and the checker is faithful (Theorem~\ref{thm:in-chk}(1)).  (2) At the
root, the $\forall$-introduction row fires, and the rows of the
recursor's successor step, universal elimination and negated-universal
introduction do not, since every entry of the tail has the non-empty
context $\vz\jnat$.  At the middle entry the conversion row fires.  At
the hypothesis entry, whose tail is empty, the hypothesis row fires,
and none of the rows from the recursor's successor step to induction
fires, since each cites a premise.  So every node's tag is at most a
firing row and lies outside the constrained rows (\textsf{fwP\_tag\_off}),
and the scans of $\Frag$, $\LvC$ and $\LvC_N$ pass.  (3) The scan's HOL value at
$\quo{P_{\mathrm w}}$ is $0$ (Theorem~\ref{thm:c-qfscan}), which \rga
derives and so refutes the scan; acceptance by $\ChkQf$ would derive
the scan, against \rga's consistency.  (4) is (1)--(3).  \qed

\noindent
The tags of the three nodes are bounded, not computed exactly:
$\Frag$ needs only that no constrained row is the tag.  The witness is
in the classes because no clause of the classes applies at the witness's nodes.

\paragraph{A universal elimination is read as a conversion.}
Extending $P_{\mathrm w}$ by a universal elimination does not give a
witness with a swap node at a quantified position.  Let $P_2$
(\textsf{fwP2}) add to $P_{\mathrm w}$ the entry $\z\jnat$ and the root
$f\cdot\z$, which follows from $\cA\cdot f$ and $\z\jnat$ by universal
elimination.  The node table reads that root as a conversion: its tag
is at most $20$, the conversion row (\textsf{fwP2\_root\_conv}).  The
conversion row comes before the universal-elimination row, and it
fires because the root's head reduct $\z = \z$ is the entry $\z\jnat$,
in the same context.  In general, a root whose head reduct is a later
judgment of the proof, in the same context, has tag at most~$20$
(\textsf{root\_conv\_of\_reduct}).  So the same holds at every root
that eliminates $\cA\cdot f$ at a term $a$, in any context, since
universal elimination cites $a\jnat$, which is $a = a$, the reduct of
$f\cdot a$ (\textsf{allE\_fwF\_root\_conv}).  A root that the table
reads as a universal elimination needs a function whose reduct at the
argument is not a judgment of the proof.  The proof $P_3$ of
\S\ref{sec:con-frag} is one: with $g = \lambda x.\ (\lambda y.\ y =
y)\cdot x$, $P_3$ eliminates $\forall x.\ g\cdot x$ at $\z$, the root's
reduct $(\lambda y.\ y = y)\cdot\z$ is not a judgment of the proof, and
the root's tag is exactly $25$ (\textsf{fwP3\_tag0}).  $P_3$ lies in
$\Frag$, in $\LvC$ and in every $\LvC_N$ through that swap node
(\textsf{fwP3\_in\_frag}, \textsf{fwP3\_in\_lvi},
\textsf{fwP3\_in\_lv}): at the root the critical formula $g\cdot\vz$
is $\forall$-free, well-formed and a term code, and its head reduct
$(\lambda y.\ y = y)\cdot\vz$ is static and $\forall$-free, so the critical formula is
certified at depth~$1$ and at every level.

\subsection{The residual}
\label{app:fold-res}

The residual theorem takes the two interfaces themselves as a
hypothesis about the proof code.

\begin{defi}[\textsf{swB}, \textsf{iRes}]
\label{def:c-res}
For a critical formula $C$, an environment $z$, a variable $w$, codes
$X, Y$, values $v, k$ and a fuel $s$, the \emph{unguarded swap body} is
\[
  \mathit{Wf}(z) \limp \Rcd_v(X) \limp \Rcd_v(Y) \limp
  \Rcd^s_k(\Subst(z[w{:=}X], C)) \limp \Rcd_k(\Subst(z[w{:=}Y], C)) ,
\]
where $\mathit{Wf}(z)$ says that every entry of $z$ is a well-formed
term code, a decided test (\textsf{sigWfGnChk}; the body is
\textsf{swB\_unguarded}).  The \emph{guarded swap body}
\textsf{swB} prefixes the four decided tests $\Wfc(X) \ne \z$,
$\IsCode(X) \ne \z$, $\Wfc(Y) \ne \z$, $\IsCode(Y) \ne \z$.  The
residual is $\Res(p) \ldef \mathit{ResGuard}(p) \wedge
\mathit{ResSwap}(p)$, where
\begin{itemize}
\item $\mathit{ResGuard}(p)$ says: if $\mathit{TreeValid}(p)$, then at
  every occupied rank $r$, for all $nb$ and $sb$, $\mathit{Guard}(p,
  r, nb, sb)$ is decided;
\item $\mathit{ResSwap}(p)$ says: if $\mathit{TreeValid}(p)$, then at
  every occupied rank $r$ and for all $nb$, if the node has one of the
  five swap tags, the guarded swap body holds at the node's critical
  formula and the arm's environment, for all $w, X, Y, v, k, s$.
\end{itemize}
\end{defi}

\noindent
Eliminating $\mathit{ResGuard}$ gives the guard interface
(\textsf{gdIf\_iRes}).  Eliminating $\mathit{ResSwap}$ and the four
guards gives the swap interface at the trivial guard
(\textsf{swIf\_iRes}).  Theorem~\ref{thm:c-fold-pathG} at $R = \Res$
(\textsf{fold\_consistent\_swIf}) gives Theorem~\ref{thm:con-res}.  In
the derivation of Theorem~\ref{thm:con-res} the residual is only
eliminated: it is not introduced there, and it is never the antecedent
of an internal implication.  (It is introduced at swap-free codes in
Theorem~\ref{thm:con-res-frag}(2,3), below.)

\begin{thm}[\textsf{rga\_sound\_under\_res}]
\label{thm:c-res-sound}
Let $\Gamma$ be a well-formed context and $p$ a term, with $v_1$ free in
neither.  If $\Gamma \vdash p \jnat$ and $\Gamma \vdash \Res(p)$, then
$\Gamma \vdash \forall c.\ \Chk(\lft p, c) \limp \Rcd_1(\Subst(\z, c))$.
\end{thm}

\noindent
At $p = \vz$ and the context of Theorem~\ref{thm:con-res} this theorem is
$\Res(\vz),\ \vz\jnat \vdash \forall c.\ \Chk(\lft\vz, c) \limp
\Rcd_1(\Subst(\z, c))$ (\textsf{rga\_sound\_under\_res\_v0}).

\noindent
The residual built from the unguarded body derives the guarded one
(\textsf{iRes\_unguarded\_iRes}), so the consistency theorem under the
unguarded residual follows by one cut
(\textsf{rga\_con\_under\_res\_unguarded}).  The guarded theorem is
therefore the stronger.  The unguarded body is refutable at a concrete
critical formula, and so is the unguarded residual at the code of an
accepted proof; the guards defeat that counterexample.
Appendix~\ref{app:traps-t01} on the unguarded swap body gives the
theorems.

\paragraph{The residual is not derivable.}
Corollary~\ref{thm:con-res-underivable} (\textsf{iRes\_underivable}) is
proved as in \S\ref{sec:con-res}: a cut of $\vz\jnat \vdash \Res(\vz)$
against Theorem~\ref{thm:con-res} and universal introduction would give
$\emptyset \vdash \ConR$, which G\"odel's second theorem
(\textsf{rga\_goedel2}) excludes.

\paragraph{The residual at rejected codes.}
Where \rga refutes acceptance, it derives the residual: if $\Gamma
\vdash \neg\mathit{TreeValid}(p)$, then $\Gamma \vdash \Res(p)$
(\textsf{iRes\_of\_rejected}), a derived rule.  Both halves of $\Res(p)$ are
implications from $\mathit{TreeValid}(p)$, that is, disjunctions with
first disjunct $\neg\mathit{TreeValid}(p)$, so one $\lor$-introduction
gives each half.  At the numeral of every code that the validity test
rejects, closed evaluation supplies the refutation, so $\emptyset
\vdash \Res(\tofn p)$ there (\textsf{iRes\_of\_rejected\_num}).

\paragraph{The residual and the fragment.}
\rga derives the residual on swap-free codes of the fragment
(Theorem~\ref{thm:con-res-frag}), but not from membership in the
fragment.

\begin{thm}[\textsf{fwP3\_iRes\_underivable}, \textsf{fwP3\_iResSwap\_}\allowbreak\textsf{underivable}, \textsf{iRes\_not\_from\_frag}, \textsf{iRes\_not\_from\_lvi}, \textsf{iRes\_not\_from\_lv}]
\label{thm:con-res-not-frag}
\leavevmode
\begin{enumerate}
\item $\emptyset \nvdash \Res(\bquo{P_3})$, and $\emptyset \nvdash
  \mathit{ResSwap}(\bquo{P_3})$.
\item $\emptyset \nvdash \forall p.\ \Frag(p) \limp \Res(p)$,
  $\emptyset \nvdash \forall p.\ \LvC(p) \limp \Res(p)$, and,
  for every $N$, $\emptyset \nvdash \forall p.\ \LvC_N(p) \limp
  \Res(p)$.
\item $\Frag(\vz),\ \vz\jnat \nvdash \Res(\vz)$, $\LvC(\vz),\
  \vz\jnat \nvdash \Res(\vz)$, and, for every $N$,
  $\LvC_N(\vz),\ \vz\jnat \nvdash \Res(\vz)$.
\end{enumerate}
\end{thm}

\noindent
\emph{Proof idea.}  The root of $P_3$ is a universal elimination with
critical formula $g\cdot\vz$.  The swap half at the root, at $X =
\bquo{\tofn 1}$ and at the code of a well-formed term $t$,
gives $\emptyset \vdash \neg\Rcd_1(\bquo{t}) \lor \Rcd_1(\bquo{g\cdot
t})$ (\textsf{fwP3\_res\_disj\_gen}).  Take for $t$ the search
$\mathit{conSearch}$ through all proof codes for an accepted proof
of $\z = \tofn 1$ (the unrestricted counterpart of the search of
\S\ref{sec:con-frag}; Theorem~\ref{thm:c-search} and
{}[\citen{goedel-paper}, \S\,``Consequences'']).  By the
disjunction property, \rga would derive one of the two disjuncts.  The
first disjunct is the refutation of the search's value-$\ctrue$ record, which
G\"odel's second theorem excludes (\textsf{con\_search\_}\allowbreak\textsf{unrefutable}).
The second disjunct would give $g\cdot\mathit{conSearch}$ a value, and the
search has no value at any fuel (\textsf{conSearch\_no\_value}, a HOL
fact that follows from \rga's consistency).  The guard half holds at
$P_3$ (\textsf{fwP3\_iResGuard}), so the swap half is what fails.
Parts~2 and~3 follow because $P_3$ lies in $\Frag$, in $\LvC$ and
in every $\LvC_N$: each sentence would give $\Res(\bquo{P_3})$ by
instantiation, and a derivation in context would give the sentence by
discharging the decided class and introducing the universal
(\textsf{iRes\_not\_from\_frag\_ctx}, \textsf{iRes\_not\_from\_lvi\_ctx},
\textsf{iRes\_not\_from\_lv\_ctx}).  \qed

\medskip
\noindent
Parts~2 and~3 are the implication and the entailment forms of the
same failure.  The class tests are decided, so by the deduction theorem
(Theorem~\ref{thm:in-dt}(1), \textsf{deduction\_theorem\_exact}) the two forms are equivalent, and
the residual fails to follow from the classes in either form.

\medskip
The instance involves no quantifier: $g\cdot\vz$, $\tofn 1$ and the
search are $\forall$-free.  What fails is the form in which the
residual states the swap.  The fold consumes its interfaces as rules
about derivations (Appendix~\ref{app:fold-record}), and for $\Frag$
Lemma~Q and the swap engine supply both interfaces as derived rules.  The residual
states the swap as one formula, whose antecedents include the records
of arbitrary terms, and at $P_3$ \rga does not derive that formula.
So within the fragment the swap positions are responsible: the guard
half holds on all of $\Frag$, the whole residual holds at the codes
without swap nodes, and the swap half fails at $P_3$, whose root is a
swap node.  The residual theorem and the fragment theorem are two
instances of one fold; the fragment supplies the two interfaces as
derived rules, and the residual assumes them as a formula.  Whether
the residual holds at the code of every genuine proof in a sense
weaker than derivability in \rga, for instance semantically through
the companion paper's soundness invariant, has a negative answer at
$P_3$.  \oga does not derive $\Res(\bquo{P_3})$, and $\Res(\bquo{P_3})$
is neither true in the $\omega$-semantics nor grounded-true
(Theorem~\ref{thm:con-res-lf}); at a closed sentence, the soundness
invariant asserts grounded truth.  Which genuine codes satisfy the
residual semantically is not characterized: the residual is grounded-true
wherever \rga derives it (\textsf{rga\_provable\_wosr}), as at
$P_{\mathrm w}$, and it fails at $P_3$.

\medskip
In detail, Theorems~\ref{thm:con-res-frag} and~\ref{thm:con-res-not-frag} are
proved in \texttt{ConOGA/ConOGA\_Residual\_Frag}.  The no-swap scan is
$\mathit{NoSwap}(p) \ldef \mathit{NodeScan}(K, p)$ with $K(x)$ the
conjunction of the five tests $\mathit{tag}(x) \ne t$ for $t \in \{24,
25, 26, 30, 40\}$, the tags of the five swap rules, built from the
compiled tag reader like the node checks of $\Frag$ (\textsf{iNoSwap});
the scan is decided (\textsf{hq\_iNoSwap}).
\begin{itemize}\raggedright
\item \emph{The guard half.}  If a residual $R$ supplies the guard
  interface, then $\Gamma \vdash R(p)$ and $\Gamma \vdash p\jnat$ give
  $\Gamma \vdash \mathit{ResGuard}(p)$ in every well-formed $\Gamma$
  (\textsf{iResGuard\_of\_gdIf}).  The two implications of
  $\mathit{ResGuard}$ have decided antecedents, the validity test and
  the occupancy test, so the implications are introduced by
  \textsf{ir\_impI\_bool}, and the universals by $\forall$-introduction;
  the guard interface closes the body.  At $R = \Frag$
  (\textsf{gdIf\_F}), discharging the decided $\Frag(\vz)$ gives
  Theorem~\ref{thm:con-res-frag}(1) (\textsf{iFrag\_iResGuard}).
\item \emph{The swap half at swap-free codes.}  If $\Gamma$ is
  well-formed, $\Gamma \vdash p\jnat$ and $\Gamma \vdash \mathit{NoSwap}(p)$, then $\Gamma \vdash
  \mathit{ResSwap}(p)$ (\textsf{iResSwap\_of\_noswap}).  At an occupied
  path point the node is the node of a list position
  (\textsf{poNode\_suffix}); there the scan's check holds, and the check
  refutes each of the five tag antecedents, so each implication
  $\mathit{tag} = t \limp \cdots$ is introduced over a refuted decided
  antecedent.  With the guard half, the swap half gives
  Theorem~\ref{thm:con-res-frag}(2); part~(3) is part~(2) at
  $P_{\mathrm w}$, whose three tags lie outside the five swap tags
  (\textsf{fwP\_tag\_off}, \textsf{fwP\_iNoSwap}).
\item \emph{The disjunction at $P_3$.}  From $\emptyset \vdash
  \Res(\bquo{P_3})$, for every well-formed $t$, $\emptyset
  \vdash \neg\Rcd_1(\bquo{t}) \lor \Rcd_1(\bquo{g\cdot t})$
  (\textsf{fwP3\_res\_disj\_gen}).  The root of $\bquo{P_3}$ is the
  post-order maximum, accepted, occupied and of tag $25$; its path
  environment is $\z$ (\textsf{path\_env\_root}).  The swap half there
  gives the guarded swap body at the critical formula
  $\mathit{ap}(\mathit{lift}(F), \bquo{\vz})$, where $F = \bquo{g}$
  and $\mathit{lift}(F)$ codes the lift of $g$ (\textsf{crit25\_kid\_node}).  At $X = \Num\cdot 1$ and $Y =
  \bquo{t}$ the two instances are $\bquo{g\cdot 1}$ and $\bquo{g\cdot
  t}$.  The $X$ side has a record by evaluation, hence a core at an
  eigen-fuel; eliminating the body up to the $Y$ record leaves
  $\neg\Rcd_1(Y) \lor (\mathit{core} \limp \Rcd_1(\bquo{g\cdot t}))$,
  and disjunction elimination gives the fuel-free disjunction.
\item \emph{The search has no value.}  For all $s$ and $r$,
  $\mathit{conSearch}$ does not evaluate to $r$ at fuel $s$
  (\textsf{conSearch\_no\_value}).  The search is $(W\,W)\cdot\z$ with
  $W = \lambda w.\lambda n.\ \Chk(n, \bquo{\z = \tofn 1}) \lor (w\,
  w)\cdot\suc n$.  By induction on the fuel, the stage $(W\,W)\cdot
  \tofn n$ evaluates only through its $\beta$-reduct, the disjunction
  $\Chk(\tofn n, \bquo{\z = \tofn 1}) \lor (W\,W)\cdot\suc\tofn n$.
  The disjunction's first disjunct never evaluates to $\ctrue$ (by soundness and
  consistency), and its second is the next stage, at a smaller fuel.
\item \emph{Assembly.}  The disjunction property (\textsf{ir\_disj\_prop})
  applied at $t = \mathit{conSearch}$ leaves two cases.  The first case
  contradicts \textsf{con\_search\_unrefutable}.  In the second case, the
  record of the $\forall$-free $g\cdot\mathit{conSearch}$ is faithful
  (\textsf{iRec\_faithful}, \textsf{orec\_iff\_wosr}), and two
  $\beta$-inversions and the inversion of an equation give a value of
  $\mathit{conSearch}$, against \textsf{conSearch\_no\_value}.  This argument is
  \textsf{fwP3\_iRes\_underivable}; with \textsf{fwP3\_iResGuard} the argument
  gives \textsf{fwP3\_iResSwap\_underivable}.  The sentences of
  Theorem~\ref{thm:con-res-not-frag}(2) would give
  $\Res(\bquo{P_3})$ by instantiation and modus ponens against
  \textsf{fwP3\_in\_frag} and \textsf{fwP3\_in\_lv}; a derivation in
  the context of part~(3) would give the sentence by
  \textsf{ir\_impI\_bool} over the decided class and
  \textsf{ir\_allI\_open}.  The same argument with
  \textsf{fwP3\_in\_lvi} gives the forms over $\LvC$ in
  Theorem~\ref{thm:con-res-not-frag}(2,~3): \rga derives
  neither the sentence $\forall p.\ \LvC(p) \limp \Res(p)$
  (\textsf{iRes\_not\_from\_lvi}) nor $\Res(\vz)$ under the
  hypotheses $\LvC(\vz)$ and $\vz\jnat$
  (\textsf{iRes\_not\_from\_lvi\_ctx}).
\end{itemize}
\rga also refutes the scan at $P_3$, deriving $\emptyset \vdash
\neg\mathit{NoSwap}(\bquo{P_3})$ (\textsf{fwP3\_not\_noswap}), so
the two witnesses $P_{\mathrm w}$ and $P_3$ are separated inside \rga.

\paragraph{The guard half beyond the fragment.}
The guard half holds on $\Frag$ (Theorem~\ref{thm:con-res-frag}(1)) but
not at every code.  The witnesses are proofs with a quantified case
hypothesis.  For a program $x$, let $a_x$ be the sentence $\forall n.\
\mathit{ofam}(x)\cdot n$, where the $\forall$-free halting family
$\mathit{ofam}(x)$ says at $n$ that program $x$ on input $x$ has not
halted by step $n$ (Appendix~\ref{app:oga}).  So $a_x$ says that the
program never halts.  The valid proof $P_x$ (\textsf{gxP}) derives
$(\lambda y.\ y = y)\cdot\z$ in the context $a_x \lor a_x$ by
disjunction elimination.  The case hypothesis of that elimination is $a_x$, and the branch
below that hypothesis is one conversion of $\z\jnat$.

\begin{thm}[\textsf{iResGuard\_not\_uniform}, \textsf{iResGuard\_not\_ctx}, \textsf{gdIf\_not\_of\_gx}, \textsf{gxP\_not\_frag}]
\label{thm:con-guard-not}
\leavevmode
\begin{enumerate}
\item $\emptyset \nvdash \forall p.\ \mathit{ResGuard}(p)$, and $\vz\jnat
  \nvdash \mathit{ResGuard}(\vz)$.
\item If $\emptyset \vdash R(\bquo{P_x})$ for every $x$, then the
  formula $R(p)$ does not supply the fold's guard interface: it is not
  the case that, in every well-formed context $\Delta$ that derives
  $R(q)$, $\mathit{TreeValid}(q)$ and that $r$ is a node of $q$,
  $\Delta$ derives that the guard item of $q$ at $r$ is decided
  (\textsf{gdIf}).
\item $\emptyset \nvdash \Frag(\bquo{P_x})$ for every $x$.
\end{enumerate}
\end{thm}

\noindent
\emph{Proof idea.}  The argument is Turing's, as for
Theorem~\ref{thm:con-not-lf}.  Suppose that \rga derived, for every
$x$, that the guard at the branch of $P_x$ is decided.  By the
disjunction property, \rga would derive the guard or its negation.  The
guard is the fuel-bounded record of $a_x$ at value $\ctrue$.  If \rga
derives the guard, \rga proves the family of $a_x$ true; if \rga refutes the guard,
the record's rows give an instance $n$ at which program $x$ has halted.
So the programs that never halt would be exactly those whose family
\rga proves, a recursively enumerable set, and the never-halting programs do not form such a set
(\textsf{gx\_guard\_dec\_absurd}).  Either statement of part~1 gives
the guard half at every $P_x$, and part~2 is the same argument from the
guard interface.  Part~3 holds because $a_x$ is not $\forall$-free.
\qed

\medskip
Theorem~\ref{thm:con-guard-not} exhibits the failure that \S\ref{sec:con-break}, on where the fold
stops, describes: a quantified
case hypothesis on the path stops the fold, at the guard.

\medskip
In detail, Theorem~\ref{thm:con-guard-not} is proved in
\texttt{ConOGA/}\allowbreak\texttt{ConOGA\_Guard\_NoGo}.  Write $a_x = \cA\cdot
\mathit{ofam}(x)$ and $H_x$ for the context $\{a_x \lor a_x\}$.  The
proof list $P_x$ (\textsf{gxP}) has four entries: $(\lambda y.\ y =
y)\cdot\z$ and $a_x \lor a_x$ in the context $H_x$, and $(\lambda y.\
y = y)\cdot\z$ and $\z\jnat$ in the context $H_x$ extended by $a_x$.
The second entry is a hypothesis, the fourth entry an axiom, the third
entry a conversion of the fourth, and the first entry a disjunction
elimination over the second and third entries.  The list is a valid proof for every $x$
(\textsf{gxP\_valid}), and the checker reads its root as a disjunction
elimination, since none of the checker's earlier rows fires there
(\textsf{gxP\_tag0}).
Inside \rga, the root's kid in slot~1 is occupied, and its added
hypothesis is $a_x$ (\textsf{gxP\_root\_kid}).
\begin{itemize}\raggedright
\item \emph{The guard at the kid.}  Disjunction elimination leaves the
  root's environment unchanged, so the guard list at the kid has one
  item: the record of $a_x$ at value $\ctrue$, bounded by fuel~$1$.
  If \rga derives the guard there, the universal row of the record
  gives the family $T_x \ldef \forall n.\
  \Rcd_1(\mathit{ap}(\bquo{\mathit{ofam}(x)}, \Num\cdot n))$
  (\textsf{gx\_rec1\_pos}).
\item \emph{A refuted guard means halting.}  If \rga refutes the guard,
  the record of $a_x$ at value $\ctrue$ is refuted at fuel~$1$, so its
  core or its exclusivity fails, and the disjunction property gives
  one of the two failures.  If the core fails, its universal row is refuted
  (its other conjuncts hold at the closed code), so $\neg T_x$, and
  the witness property gives an $n$ at which \rga refutes the
  value-$\ctrue$ record of the instance; then $x$ halts by step $n$.  If
  exclusivity fails, some value $w \ne 1$ has a record at $a_x$; at a
  quantifier node only the false row then remains, which gives an $n$
  with a value-$\cfalse$ record of the instance, and again $x$ halts by
  step $n$ (\textsf{gx\_rec1\_neg}).
\item \emph{Turing's argument.}  If the guard at the kid is decided for
  every $x$, the two cases give $x \notin K$ if and only if $\emptyset
  \vdash T_x$ (with \textsf{famT\_ofam\_sound}), so the complement of
  $K$ is recursively enumerable (\textsf{renp\_famT\_ofam}), which it
  is not (\textsf{gx\_guard\_dec\_absurd}), as in the proof of
  Theorem~\ref{thm:e-halting}.
\item \emph{Assembly.}  A derivation of $\forall p.\
  \mathit{ResGuard}(p)$ gives $\mathit{ResGuard}(\bquo{P_x})$ for every
  $x$, and eliminating that instance at the kid gives the decided guard
  (\textsf{iResGuard\_not\_uniform}).  A derivation at the free
  variable gives the sentence by $\forall$-introduction
  (\textsf{iResGuard\_not\_ctx}).  For part~2 the guard interface
  itself gives the decided guard (\textsf{gdIf\_not\_of\_gx}).  For
  part~3, the fragment's node check at the root asks the added
  hypothesis to be $\forall$-free, and $a_x$ is not
  (\textsf{gxP\_not\_frag}).
\end{itemize}

\paragraph{The residual under the level formulas.}
The hypotheses $\LF$ and $\LF_{\le N}$, which give the level theorem
and the per-level theorem (Theorems~\ref{thm:con-lv}
and~\ref{thm:con-lvN}) their swap interface, do not make the residual
derivable on $\LvC$ or on the per-level classes, and neither does any
set of hypotheses true in the $\omega$-semantics.

\begin{thm}[\textsf{iRes\_not\_}\allowbreak\textsf{from\_lv\_omega}, \textsf{iRes\_not\_}\allowbreak\textsf{from\_lv\_oga}, \textsf{iRes\_not\_from\_lv\_}\allowbreak\textsf{under\_lf\_ctx}, \textsf{fwP3\_iRes\_}\allowbreak\textsf{not\_oga}, \textsf{fwP3\_iRes\_}\allowbreak\textsf{not\_omega}, \textsf{fwP3\_iRes\_}\allowbreak\textsf{not\_grounded}]
\label{thm:con-res-lf}
\leavevmode
\begin{enumerate}
\item If $\Gamma$ is a well-formed context whose members are true in
  the $\omega$-semantics, in particular if every member is a theorem of
  \oga, then $\Gamma \nvdash \forall p.\ \LvC(p) \limp \Res(p)$
  (\textsf{iRes\_not\_from\_lvi\_omega}, \textsf{iRes\_not\_from\_lvi\_oga}),
  and $\Gamma \nvdash \forall p.\ \LvC_N(p) \limp \Res(p)$ for
  every $N$.
\item $\LvC(\vz),\ \vz\jnat,\ \LF \nvdash \Res(\vz)$
  (\textsf{iRes\_not\_from\_lvi\_under\_lf\_ctx}), and, for every $N$,
  $\LvC_N(\vz),\ \vz\jnat,\ \LF_{\le N} \nvdash \Res(\vz)$.
\item $\emptyset \nvdash_{\mathrm O} \Res(\bquo{P_3})$, and $\Res(\bquo{P_3})$ is
  neither true in the $\omega$-semantics nor grounded-true.
\end{enumerate}
\end{thm}

\noindent
\emph{Proof idea.}  A derivation of $\forall p.\ \LvC_N(p) \limp
\Res(p)$ from $\Gamma$ gives $\Res(\bquo{P_3})$ from $\Gamma$, since
$P_3$ is in $\LvC_N$ (and the same holds for $\LvC$, which also contains $P_3$), and so the disjunction $\neg\Rcd_1(\bquo{t}) \lor
\Rcd_1(\bquo{g\cdot t})$ from $\Gamma$ for every well-formed $t$
(\textsf{fwP3\_res\_disj\_ctx}).  \rga is sound for the
$\omega$-semantics, so the disjunction would be $\omega$-true.  Take
for $t$ the term $W\cdot W$, with $W = \lambda w.\ \forall x.\ w\,w$.
The evaluation of $W\cdot W$ passes through a universal and returns to itself, and
$W\cdot W$ has no $\omega$-value (\textsf{omv\_cycY\_unvalued}).  The record of $W\cdot W$
follows the same cycle, through the record's universal rows and back
to the code of $W\cdot W$, and is ungrounded: $\Rcd_1(\bquo{W\cdot W})$
has no $\omega$-value, and no $\Rcd_v(\bquo{W\cdot W})$ is $\omega$-true
(\textsf{iRec1\_cycY\_unvalued}, \textsf{iRec\_cycY\_not\_true}).
$g\cdot(W\cdot W)$ reduces to the equation $W\cdot W = W\cdot W$,
whose record needs a true record of $W\cdot W$, so the disjunction is
not $\omega$-true at $W\cdot W$ (\textsf{omv\_resDisj\_cycY},
below).  This argument gives part~1.  The case of part~1 for \oga
theorems holds because the theorems of \oga are $\omega$-true
(\textsf{oga\_omega\_true}), and $\LF_{\le N}$ is one of those theorems
(\textsf{o\_iLFs}, Theorem~\ref{thm:e-olf-N}), so $\LF_{\le N} \nvdash \forall p.\
\LvC_N(p) \limp \Res(p)$ (\textsf{iRes\_not\_from\_lv\_under\_lf}).  $\LF$ is one of
those theorems too (Theorem~\ref{thm:con-olf}), so $\LF \nvdash \forall p.\ \LvC(p)
\limp \Res(p)$ (\textsf{iRes\_not\_from\_lvi\_under\_lf}).
Part~2 is the same argument from a derivation at the free variable.
For part~3, an \oga derivation of $\Res(\bquo{P_3})$, or an
$\omega$-true $\Res(\bquo{P_3})$, would make the disjunction
$\omega$-true in the same way, and a grounded-true sentence is
$\omega$-true (\textsf{wosr\_omv}).  \qed

\medskip
Part~1 holds for every such $\Gamma$, so even the level formula of a
higher level, $\LF_{\le M}$ with $M > N$, or the level-free $\LF$,
together with any other theorems of \oga, does not make the residual
derivable on $\LvC_N$, nor on $\LvC$.
With $\Gamma$ empty, part~1 gives the level half of
Theorem~\ref{thm:con-res-not-frag}(2) again, by a different argument.
The proof of that theorem uses G\"odel's second theorem at the search, this proof the
$\omega$-semantics at $W\cdot W$.

\medskip
In detail, the theory \texttt{ConOGA/}\allowbreak\texttt{ConOGA\_Residual\_LF}
proves Theorem~\ref{thm:con-res-lf}.
\begin{itemize}\raggedright
\item \emph{The disjunction in any context.}  The disjunction at $P_3$
  holds in every well-formed context $\Gamma$: from $\Gamma \vdash
  \Res(\bquo{P_3})$, $\Gamma \vdash \neg\Rcd_1(\bquo{t}) \lor
  \Rcd_1(\bquo{g\cdot t})$ for every well-formed $t$
  (\textsf{fwP3\_res\_disj\_ctx}), since the other facts of that proof
  are closed and weaken.  With \textsf{fwP3\_in\_lv}, a derivation of
  $\forall p.\ \LvC_N(p) \limp \Res(p)$ under $\LF_{\le N}$ gives the
  disjunction under $\LF_{\le N}$ (\textsf{lf\_res\_disj}).
\item \emph{To the $\omega$-semantics.}  \oga proves $\LF_{\le N}$
  (\textsf{o\_iLFs}), so $\LF_{\le N}$ is $\omega$-true
  (\textsf{oga\_omega\_true}), and \rga derivations from $\omega$-true
  hypotheses have $\omega$-true conclusions (\textsf{rga\_omega\_sound}).
  So the disjunction is $\omega$-true (\textsf{lf\_res\_omega}; from a
  derivation at the free variable, \textsf{lf\_res\_omega\_ctx}).  The
  same holds if $\Res(\bquo{P_3})$ is $\omega$-true
  (\textsf{res\_fwP3\_omega\_disj}) or derivable in \oga
  (\textsf{oga\_res\_omega}).
\item \emph{The cycle.}  Let $W = \lambda w.\ \forall x.\ w\,w$
  (\textsf{cycY} is $W\cdot W$).  No $\omega$-derivation ends at
  $W\cdot W$, at $\forall x.\ W\cdot W$, or at $(\lambda x.\ W\cdot
  W)\cdot\tofn n$, because each such derivation needs one at the next term of this
  cycle; a rule induction over the $\omega$-semantics proves this claim, so
  $W\cdot W$ has no $\omega$-value (\textsf{omv\_cycY\_unvalued}).
\item \emph{The record of the cycle.}
  $\Rcd_1(\bquo{W\cdot W})$ unfolds through the application row to
  the record of $\bquo{\forall x.\ W\cdot W}$, whose universal rows
  ask for the records of the instances $(\lambda x.\ W\cdot W)\cdot
  \tofn n$, which unfold back to the record of $\bquo{W\cdot W}$.  At
  these three codes every other row is dead: the row's guard has a computed
  value that refutes the row.  The pairs of a term and a value that no
  $\omega$-derivation reaches are collected inductively
  (\textsf{omBlk}): the records at the three codes (from truth at every
  value; from falsity at value $\ctrue$ from fuel $2$, $1$, $3$
  respectively), the pairs whose term has another value, and their
  closure under the backward readings of the rules.  A negation swaps
  the polarity; a disjunction is blocked from truth when both disjuncts
  are, from falsity when one is; a $\beta$-redex is blocked when its
  contractum is; a universal is blocked from truth when one instance
  is, from falsity when all are.  One rule induction over the
  $\omega$-semantics, with one case per rule, shows that no blocked
  pair has an $\omega$-value (\textsf{omBlk\_never}).  The record's
  shapes do not enter that induction.  The induction's one specific case is the record's
  self-application, whose contractum at the next code of the cycle is
  blocked again; that blocking is proved forward through the record's rows, from
  truth by the live rows, and from falsity by the live row at value
  $\ctrue$ together with the exclusivity conjunct, whose instance at
  another value needs a true record at the same code.  The fuel
  thresholds close the false half around the cycle: from fuel~$2$ the
  bounded kid of the record at $\bquo{W\cdot W}$ is the record at
  $\bquo{\forall x.\ W\cdot W}$ at a fuel of at least~$1$; the
  universal's family is not fuel-bounded, so its instances are taken
  from fuel~$3$; and from there the cycle returns to $\bquo{W\cdot W}$
  from fuel~$2$.  So no record at
  $\bquo{W\cdot W}$ is $\omega$-true (\textsf{iRec\_cycY\_not\_true}),
  $\Rcd_1(\bquo{W\cdot W})$ is not $\omega$-false either, and it has no
  $\omega$-value (\textsf{iRec1\_cycY\_unvalued}).  The values other
  than $\ctrue$ and $\cfalse$ need no separate refutation; they are blocked
  from truth like every value.
\item \emph{The disjunction at the cycle.}  $g\cdot(W\cdot W)$
  head-reduces to $W\cdot W = W\cdot W$, whose record needs a true
  record of $W\cdot W$, so the record of $g\cdot(W\cdot W)$ is blocked
  from truth, and the disjunction is not $\omega$-true at $t = W\cdot
  W$ (\textsf{omv\_resDisj\_cycY}).
\item \emph{Assembly.}  From a derivation of $\forall p.\ \LvC_N(p)
  \limp \Res(p)$ in a well-formed context $\Gamma$ of $\omega$-true
  hypotheses, instantiation at $\bquo{P_3}$ and
  \textsf{fwP3\_res\_disj\_ctx} give the disjunction at $W\cdot W$ in
  $\Gamma$, which \textsf{rga\_omega\_sound} makes $\omega$-true,
  contradicting the previous step (\textsf{iRes\_not\_from\_lv\_omega});
  with \textsf{oga\_omega\_true}, the hypotheses may be any theorems of
  \oga (\textsf{iRes\_not\_from\_lv\_oga}).  The forms under $\LF_{\le
  N}$ and at the free variable follow from \textsf{lf\_res\_omega} and
  \textsf{lf\_res\_omega\_ctx} (\textsf{iRes\_not\_from\_lv\_under\_lf},
  \textsf{iRes\_not\_from\_lv\_under\_lf\_ctx}), part~3 from
  \textsf{oga\_res\_omega} and \textsf{res\_fwP3\_omega\_disj}
  (\textsf{fwP3\_iRes\_not\_oga}, \textsf{fwP3\_iRes\_not\_omega}) and
  \textsf{wosr\_omv} (\textsf{fwP3\_iRes\_not\_grounded}).
\item \emph{The forms over $\LvC$} (Theorem~\ref{thm:con-res-lf}(1,~2)).
  The same assembly with
  \textsf{fwP3\_in\_lvi} in place of \textsf{fwP3\_in\_lv} shows
  that no well-formed context of $\omega$-true hypotheses derives
  $\forall p.\ \LvC(p) \limp \Res(p)$
  (\textsf{iRes\_not\_from\_lvi\_omega}), nor one of theorems of
  \oga (\textsf{iRes\_not\_from\_lvi\_oga}).  Since \oga proves
  $\LF$ (\textsf{o\_iLFI}), the context $\{\LF\}$ is one such context:
  under the hypothesis $\LF$, \rga derives neither the sentence
  (\textsf{iRes\_not\_from\_lvi\_under\_lf}) nor $\Res(\vz)$
  under the further hypotheses $\LvC(\vz)$ and $\vz\jnat$
  (\textsf{iRes\_not\_from\_lvi\_under\_lf\_ctx}).  No
  level-indexed hypothesis is needed.
\end{itemize}

\paragraph{The guard half and the fold's motive in $\omega$-sound
systems.}
Theorem~\ref{thm:con-res-ext} is proved in
\texttt{ConOGA/}\allowbreak\texttt{ConOGA\_Motive\_NoGo}.  Both of that theorem's
semantic parts use only \rga's $\omega$-soundness and the ungrounded
term $Y = W\cdot W$ of the previous paragraph; no disjunction property
is used.  Write $K = \lambda.\,\z\jnat$ and $Z = K\cdot\z$.  The proof
$P_{\mathrm h}$ (\textsf{phP}) has nine judgments, all but the last two
in the empty context:
\begin{itemize}\raggedright
\item $Z \lor Y$, by disjunction elimination on $Y \lor Z$ with the
  branches $Y \vdash Z \lor Y$ and $Z \vdash Z \lor Y$;
\item $Y \lor Z$, by substitution of $\z$ for $\vz$ in $Y \lor
  K\cdot\vz$, which comes by disjunction introduction from $K\cdot\vz$,
  which is a conversion of $\z\jnat$, an axiom;
\item in the context $\{Z\}$: $Z \lor Y$ by disjunction introduction
  from the hypothesis $Z$;
\item in the context $\{Y\}$: $Z \lor Y$ (the \emph{branch node}) by
  disjunction introduction from the hypothesis $Y$ (the
  \emph{hypothesis node}).
\end{itemize}
The proof is valid (\textsf{phP\_valid}) and accepted by the checker
at the empty context (\textsf{phP\_accepted}).  The disjunction comes
from an open instance by substitution, not from $Z$ by disjunction
introduction, because of the checker's row order: with the judgment
$\emptyset \vdash Z$ in the list, the cut row would fire at the root
before disjunction elimination, and cut adds no guard.  The checker
reads the root as a disjunction elimination (\textsf{phP\_tag0}; the checker's
earlier rows are off by \textsf{lowRows\_orapp\_off}, a generic lemma
for a disjunction whose left disjunct is an application) and the
branch node as a disjunction introduction (\textsf{phQ2\_tag}).  The
node reader is evaluated in HOL at both nodes, and inside \rga
\textsf{phP\_walk} reads the walk: the two tags, the node codes at the
two ranks, the empty environment at both, so that the guard item and
the hypothesis node's conclusion are both exactly $\bquo{Y}$, and the
occupancy of the hypothesis node.
\begin{itemize}\raggedright
\item \emph{The guard half} (part~1).  In every context, the guard
  half at the branch node, with eigenvariable tuple $(\z)$ and fuel
  tuple $(2)$, yields $B \lor \neg B$ for the fuel-bounded record $B =
  \Rcd^2_1(\bquo{Y})$ (\textsf{ph\_guard\_dec\_disj},
  \textsf{phBK}).  $B$ has no value in the $\omega$-semantics
  (\textsf{phBK\_unvalued}): it is blocked from truth at every fuel and
  from falsity at fuel~$2$, by the blocked pairs of the previous
  paragraph.  At fuels $0$ and $1$ the record is decided false, so $2$
  is the least fuel that works.  So the guard half is not
  $\omega$-true at $\bquo{P_{\mathrm h}}$
  (\textsf{phP\_iResGuard\_not\_omega}), no well-formed context of
  $\omega$-true hypotheses derives the guard half there
  (\textsf{phP\_iResGuard\_underivable\_omega}), and no such context derives
  $R(\bquo{P_{\mathrm h}})$ for a residual $R$ that supplies the guard
  interface (\textsf{phP\_gdIf\_underivable\_omega}).  The open form
  $\vz\jnat \nvdash \mathit{ResGuard}(\vz)$
  (\textsf{iResGuard\_open\_underivable\_omega}) is a second proof,
  without the disjunction property, of the second statement of
  Theorem~\ref{thm:con-guard-not}(1).
\item \emph{The motive} (part~2).  In every context, the body of the
  motive (Definition~\ref{def:c-recmot}) at the hypothesis node, with
  eigenvariable tuple $(\z, \z)$ and fuel tuple $(2)$, yields $\neg B
  \lor \Rcd_1(\bquo{Y})$ (\textsf{ph\_motive\_disj}): the refuted
  guard item, or the value-$\ctrue$ record of the node's conclusion.
  Neither disjunct has an $\omega$-value
  (\textsf{phBK\_unvalued}, \textsf{iRec\_cycY\_not\_true},
  \textsf{iRec1\_cycY\_unvalued}), so the disjunction is not
  $\omega$-true (\textsf{ph\_E3\_disj\_not\_omega}), and by \rga's
  $\omega$-soundness neither is the motive at the hypothesis node
  (\textsf{phP\_fMot\_not\_omega}).  Nor is the induction hypothesis
  of the fold at $P_{\mathrm h}$, the motive at every rank
  (\textsf{phMotAll\_not\_omega}).  The fold derives that motive:
  Theorem~\ref{thm:c-cov} with every step by \textsf{fold\_step}, the
  run of Theorem~\ref{thm:c-fold-root}, gives the motive at the
  hypothesis node from $\Gamma \vdash R(\bquo{P_{\mathrm h}})$ for
  every residual $R$ with $\mathit{arms}(R)$, when $v_1$ is not free
  in $\Gamma$ (\textsf{fold\_motive\_at\_phP}).  So no such $\Gamma$ of
  $\omega$-true hypotheses derives $R(\bquo{P_{\mathrm h}})$
  (\textsf{fold\_residual\_underivable\_omega}), an
  $R(\bquo{P_{\mathrm h}})$ that is well-formed and has no free $v_1$
  is not $\omega$-true
  (\textsf{fold\_residual\_not\_omega}), and $\vz\jnat \nvdash
  R(\vz)$ when $R$ commutes with substitution
  (\textsf{fold\_residual\_open\_underivable}).
\item \emph{Any $\omega$-sound system} (part~3).  The locale
  \textsf{omega\_sound\_ext} fixes an entailment $\vdash_X$ over \rga's
  terms with one assumption: if $\Gamma \vdash_X c$, $\Gamma$ is
  well-formed and every member of $\Gamma$ is $\omega$-true, then $c$
  is $\omega$-true.  In the locale, parts~1 and~2 give
  \textsf{xder\_not\_iResGuard}, \textsf{xder\_not\_motive} and
  \textsf{xder\_not\_fold\_residual}.  \rga's entailment is an
  instance (\textsf{rga\_omega\_sound\_ext}, from
  \textsf{rga\_omega\_sound}).  \oga's \textsf{omega\_sound}
  (\texttt{OGA/OGA\_Sound}) has the locale's assumption as its
  statement; the instance for \oga is not stated separately.
\end{itemize}
The quantifier ``every residual'' ranges over the term formers $R$ for
which the 41 arms hold, $\mathit{arms}(R)$; that is the formal sense of
``a fold of this shape''.  At the hypothesis node the motive
internalizes the hypothetical judgment $Y \vdash Y$ as the
implication $\neg B \lor \Rcd_1(\bquo{Y})$, a conditional at an
ungrounded antecedent.  Curry's paradox needs one more ingredient,
contraction; the fold's induction reuses its hypothesis at every rank
(the accumulator premise of \textsf{fold\_step}), but the development
does not state the contraction half as a theorem of its own.

\endgroup

\section{The level theorems and OGA}
\label{app:oga}
\begingroup\sloppy

This appendix gives the definitions and proofs behind
\S\ref{sec:con-lv}: the certified depth and the class $\LvC$, the
family sentences $\LF$, $\LF(n)$ and $\LF_{\forall}$, the level
theorem in \rga (Theorem~\ref{thm:con-lv}), the proof that \rga does
not derive $\LF$ (Theorem~\ref{thm:con-not-lf}) and does not refute it
(Theorem~\ref{thm:con-lf-indep}), the proof of $\LF$ in \oga
(Theorem~\ref{thm:con-olf}) with the \oga form of the level theorem
(Theorem~\ref{thm:con-oga}), the per-level forms
(Theorem~\ref{thm:con-lvN}), the reach-guarded theorem
(Theorem~\ref{thm:con-e1}) with its \oga form, the containments in
the level classes (Theorem~\ref{thm:con-lv-contain}), and the
containment of the per-level classes in $\LvC$
(Theorem~\ref{thm:con-lv-closed}).
Statements are given at the granularity of the mechanized lemmas, with
their hypotheses.  The fold that all these theorems instantiate is the
record fold of Appendix~\ref{app:fold}; this appendix supplies only the
fold's interfaces.  Theory files are relative to \texttt{hol/}, and every name
resolves in the development (Section~\ref{sec:mech}, ``Artifact'').

\subsection{Theory files}
\label{app:oga-files}

\begingroup
\small
\begin{longtable}{@{}>{\raggedright\arraybackslash}p{0.33\linewidth}>{\raggedright\arraybackslash}p{0.62\linewidth}@{}}
\hline
theory & names used in this appendix \\
\hline
\endfirsthead
\hline
theory & names used in this appendix \\
\hline
\endhead
\hline
\endfoot
\texttt{ConRGA/Wfd\_Static} & \textsf{iStatC}, \textsf{iFltC}, \textsf{aLeafF} \\
\texttt{ConRGA/Wfd\_General} & \textsf{iNstC}, \textsf{ngGdT}, \textsf{iStatLv}, \textsf{nodeScan}, \textsf{con\_shell} \\
\texttt{ConRGA/Wfd\_Dynamic} & \textsf{iDynWFD}, \textsf{iSub1Qf}, \textsf{dyn\_core}, \textsf{gdIf\_D} \\
\texttt{ConRGA/Wfd\_Decision} & \textsf{lemQ\_int}, \textsf{famDec} \\
\texttt{ConOGA/ConOGA\_Seam} & \textsf{oent} ($\vdashO$), \textsf{o\_cutR},
  \textsf{posOk}, \textsf{o\_qfRec}, \textsf{o\_qfFam}, \textsf{core\_B\_nest},
  \textsf{o\_famDec\_nest}, \textsf{lv\_core}, \textsf{iLF}, \textsf{iLF\_E}, \textsf{o\_iLF},
  \textsf{iHypLv}, \textsf{iLvClass}, \textsf{hq\_iLvClass},
  \textsf{oga\_omega\_true}, \textsf{rga\_unrefutable\_of\_oga},
  \textsf{iLF\_unrefutable} \\
\texttt{ConOGA/ConOGA\_Fold} & \textsf{iLFs}, \textsf{o\_iLFs}, \textsf{gdIf\_L},
  \textsf{swIfQ\_L}, \textsf{swRelL\_holds}, \textsf{rga\_con\_lv\_under\_lf}, \textsf{rga\_sound\_lv\_under\_lf},
  \textsf{rga\_con\_lv\_oga}, \textsf{rga\_sound\_lv\_oga},
  \textsf{rga\_con\_stat\_under\_dwfd\_lf}, \textsf{rga\_sound\_stat\_under\_dwfd\_lf},
  \textsf{rga\_con\_stat\_oga\_under\_dwfd}, \textsf{rga\_sound\_stat\_oga\_under\_dwfd},
  \textsf{rga\_con\_lv\_under\_lf\_all}, \textsf{iLFs\_unrefutable}, \textsf{iLFs\_imp\_iff} \\
\texttt{ConOGA/ConOGA\_Otto} & \textsf{iLF\_decides\_halting}, \textsf{rga\_not\_iLF},
  \textsf{rga\_not\_iLFs}, \textsf{rga\_oga\_separation\_iLF}, \textsf{rga\_oga\_separation\_iLFs},
  \textsf{iLF\_independent}, \textsf{iLF\_undecided}, \textsf{iLFs\_independent},
  \textsf{iLFs\_undecided} \\
\texttt{ConOGA/ConOGA\_Containment} & \textsf{lyInv\_all}, \textsf{iHypLv\_iDynWFD},
  \textsf{rga\_con\_lv\_under\_lf\_via\_stat}, \textsf{o\_iDynWFD\_lv},
  \textsf{rga\_con\_lv\_oga\_via\_stat} \\
\texttt{OGA/OGA\_Proof}, \texttt{OGA/OGA\_Tier} & rule \textsf{pr\_allTI} (ATI),
  \textsf{ga\_allTI\_bridge}, \textsf{ga\_exTI\_bridge} \\
\texttt{OGA/OGA\_Sound}; \texttt{RGA/RGA\_Omega\_Sound} & \textsf{entails\_embed};
  \textsf{ati\_omega\_sound} \\
\texttt{ConOGA/ConOGA\_Frag\_Witness} & \textsf{fwP\_in\_lv} \\
\texttt{ConOGA/ConOGA\_Frag\_QF} & \textsf{qf\_accepted\_in\_lv}, \textsf{frag\_lv\_incomparable},
  \textsf{fwPR}, \textsf{fwPR\_in\_frag}, \textsf{fwPR\_not\_lv}, \textsf{frag\_not\_lv\_at\_root},
  \textsf{fwPF}, \textsf{fwPL} \\
\texttt{RGA/RGA\_Sound} & \textsf{closed\_val} \\
\texttt{ConOGA/ConOGA\_Frag\_Lv} & \textsf{closed\_code\_of},
  \textsf{frag\_tree\_in\_lv}, \textsf{frag\_accepted\_in\_lv} \\
\texttt{ConOGA/ConOGA\_Lv\_Int} & \textsf{lv\_tree\_of\_frag\_pos}, \textsf{iNstC\_of\_qf},
  \textsf{iStatC\_of\_qf}, \textsf{iSub1Qf\_of\_iStat}, \textsf{qf\_tree\_in\_lv},
  \textsf{rga\_qf\_in\_lv\_int}, \textsf{gHR\_wf\_gn}, \textsf{iHdStep}, \textsf{iHdStep\_of\_lv},
  \textsf{frag\_tree\_in\_lv\_open}, \textsf{rga\_frag\_in\_lv\_int}, \textsf{fwPR\_not\_hdstep} \\
\texttt{ConOGA/ConOGA\_Depth} & \textsf{dHIdx}, \textsf{dIdx}, \textsf{iDep}, \textsf{dep\_kid1\_one}, \textsf{dep\_kid\_two},
  \textsf{dep\_leaf}, \textsf{dep\_leaf\_nested}, \textsf{dep\_leaf\_flat}, \textsf{iDep\_qf},
  \textsf{OkI}, \textsf{OkU}, \textsf{OkI\_B}, \textsf{OkU\_B}, \textsf{OkI\_kid1\_one},
  \textsf{OkU\_kid1\_one}, \textsf{OkU\_nested}, \textsf{OkU0\_nested}, \textsf{okI\_okU},
  \textsf{okU\_okI}, \textsf{lv\_coreU}, \textsf{lv\_core\_iLFU}, \textsf{lv\_core\_iLFI},
  \textsf{iLFU}, \textsf{iLFA}, \textsf{iLFI}, \textsf{iLFI\_E}, \textsf{iLFA\_iLFI},
  \textsf{iLFI\_iLFA}, \textsf{o\_iLFU\_I}, \textsf{o\_iLFU\_0}, \textsf{o\_iLFU\_step},
  \textsf{o\_iLFU}, \textsf{o\_iLFA}, \textsf{o\_iLFI}, \textsf{aOfamC\_okI},
  \textsf{iLFI\_at\_ofam}, \textsf{iLFI\_decides\_halting}, \textsf{rga\_not\_iLFI},
  \textsf{iLFI\_unrefutable}, \textsf{iLFI\_independent}, \textsf{iLFI\_undecided},
  \textsf{rga\_oga\_separation\_iLFI}, \textsf{iLFI\_imp\_iff}, and the \textsf{iLFA}
  twins of the last six names \\
\texttt{ConOGA/ConOGA\_Level} & \textsf{okIIdx}, \textsf{ir\_okIIdx\_eq}, \textsf{lvChkIIdx},
  \textsf{lvOkIIdx}, \textsf{iHypLvI}, \textsf{hq\_iHypLvI}, \textsf{iStatLvI}, \textsf{hq\_iStatLvI},
  \textsf{iLvClassI}, \textsf{hq\_iLvClassI}, \textsf{lv\_slotI}, \textsf{whole\_guard\_BLI}, \textsf{gdIf\_LI},
  \textsf{sl\_levelU}, \textsf{OkU\_nested\_dep}, \textsf{swNestU}, \textsf{SwU},
  \textsf{swU\_all}, \textsf{slTopU}, \textsf{swIfQ\_VLI}, \textsf{xqIf\_VLI},
  \textsf{swIfQ\_LI}, \textsf{xqIf\_LI}, \textsf{rga\_con\_lvi\_under\_lf},
  \textsf{rga\_sound\_lvi\_under\_lf}, \textsf{rga\_con\_lvi\_under\_lf\_all},
  \textsf{rga\_con\_lvi\_under\_lfa\_all}, \textsf{nfree\_iLFI}, \textsf{rga\_con\_lvi\_oga},
  \textsf{rga\_sound\_lvi\_oga}, \textsf{lvI\_tree\_of\_frag\_pos}, \textsf{qf\_tree\_in\_lvi},
  \textsf{rga\_qf\_in\_lvi\_int}, \textsf{frag\_tree\_in\_lvi\_open}, \textsf{rga\_frag\_in\_lvi\_int}, \textsf{iHdStep\_of\_lvi}, \textsf{fwP\_in\_lvi},
  \textsf{ngR\_lam\_code}, \textsf{ngR\_genuine\_lt} \\
\texttt{ConOGA/ConOGA\_Level\_Ports} & \textsf{vOk1I\_25\_fails}, \textsf{not\_lvi\_at\_node},
  \textsf{fwPR\_not\_lvi}, \textsf{frag\_not\_lvi\_at\_root}, \textsf{frag\_accepted\_in\_lvi},
  \textsf{fwPL\_in\_lvi}, \textsf{fwPF\_not\_lvi}, \textsf{frag\_lvi\_incomparable} \\
\texttt{ConOGA/ConOGA\_Level\_Closed} & \textsf{dep\_asc\_two}, \textsf{dep\_asc\_one},
  \textsf{dep\_asc\_nested}, \textsf{lamN\_code}, \textsf{posOk\_closed\_dep},
  \textsf{posOk\_closed\_okU}, \textsf{iHypLv\_closed\_lvi}, \textsf{iStatLv\_closed\_lvi},
  \textsf{lv\_closed\_in\_lvi}, \textsf{lv\_closed\_imp\_lvi}, \textsf{nodeScanP},
  \textsf{iLvClassU}, \textsf{hq\_iLvClassU}, \textsf{lvNodeUT\_lvNodeIT}, \textsf{vOkTU\_vOkTI},
  \textsf{lvU\_in\_lvI}, \textsf{rga\_con\_lvu\_under\_lf\_all}, \textsf{rga\_con\_lvu\_oga},
  \textsf{lv\_closed\_in\_lvu} \\
\texttt{ConOGA/ConOGA\_Level\_Reach} & \textsf{lyInvI\_all}, \textsf{iHypLvI\_iDynWFD},
  \textsf{gdIf\_E1LI}, \textsf{swIfQ\_E1LI}, \textsf{xqIf\_E1LI},
  \textsf{rga\_con\_statI\_under\_dwfd\_lf}, \textsf{rga\_sound\_statI\_under\_dwfd\_lf},
  \textsf{rga\_con\_statI\_oga\_under\_dwfd}, \textsf{rga\_con\_lvi\_under\_lf\_via\_stat},
  \textsf{o\_iDynWFD\_lvi}, \textsf{rga\_con\_lvi\_oga\_via\_stat} \\
\end{longtable}
\endgroup

\paragraph{Notation.}
Besides the notation of Section~\ref{sec:con}, this appendix uses the
code constructor $\mathit{ap}(a,b)$ and the numeral codes $\Num\cdot n$
of [\citen{goedel-paper}, Appendix ``G\"odel II in full''].  For a code $x$ of an application,
$\mathit{arg}(x)$ is the code of its argument (Isabelle
\textsf{iKid2}).  $\mathit{ALeaf}(c)$ is the decided test ``$c$ codes an
application whose function part is the constant $\cA$''
(\textsf{aLeafF}); such a code is a \emph{quantifier application}
$\cA\cdot m$, and $m$ is its \emph{matrix}.  For an environment $z$
(a coded list of pairs of a variable and a code),
$\mathit{EnvWf}(z)$ says that every code in $z$ is the code of a
well-formed term (\textsf{sigWfGnChk}), and $\mathit{EnvQf}(z)$ that
every code in $z$ is $\forall$-free (\textsf{envQf}).  Both environment
tests are compiled and decided.  Finally, for a matrix code $u$,
\[
  \mathit{FamDec}(u) \;\ldef\;
  \bigl(\forall n.\ \Rcd_1(\mathit{ap}(u, \Num\cdot n))\bigr)\jbool
  \;\wedge\;
  \bigl(\exists n.\ \Rcd_0(\mathit{ap}(u, \Num\cdot n))\bigr)\jbool
\]
says that both numeral families of $u$ are decided
(\textsf{famDec}, with the two families \textsf{famT} and \textsf{famF}
of \texttt{RGA\_Record/Value\_Record.thy}).

\subsection{The certified depth and the class}
\label{app:oga-level}

\paragraph{The static scan.}
A code is \emph{static} (\textsf{iStatC}, value $0$ meaning ``pass'')
when it is $\forall$-free, or a successor, negation, disjunction or
equation node whose arguments are static, or a quantifier application
$\cA\cdot(\lambda b)$ whose body $b$ is static.  The four node kinds
$\suc$, $\neg$, $\lor$ and ${=}$ form the \emph{spine} of a code.

\paragraph{The depth machine.}
The depth machine refines the static scan at quantifier applications.
It is one compiled course-of-values machine over codes, with no level
parameter, built from an explicit index tree (\textsf{dHIdx}, the
head; \textsf{dIdx}, the tree).  At a quantifier application whose
matrix is not $\forall$-free it reads the depth of a code that is not
an argument, the head reduct of the matrix's instance, from its own
table of the values at smaller codes.

\begin{defi}[\textsf{iDep}, \textsf{OkI}, \textsf{OkU}]
\label{def:e-depth}
$\Dep(c) \ldef \cmp(\mathit{dIdx})\cdot\z\cdot c$, where at each code
$c$:
\begin{itemize}
\item if $c$ is $\forall$-free, $\Dep(c) = 1$;
\item if $c$ is a spine node, $\Dep(c)$ is the largest value of its
  arguments if none of them is $0$, and $0$ otherwise;
\item if $c$ is a quantifier application $\cA\cdot m$ with $m$
  $\forall$-free, $\Dep(c) = 1$;
\item if $c$ is a quantifier application $\cA\cdot m$ with $m$ not
  $\forall$-free, $\Dep(c) = 1 + \Dep(R)$ when $m$ passes the
  \emph{nested-leaf check}, $R < c$ and $\Dep(R) \ne \z$, and
  $\Dep(c) = 0$ otherwise.  The nested-leaf check asks that $m$ code
  an abstraction $\lambda b$; that its \emph{target}
  $\mathit{ap}(\mathit{lift}(m), \bquo{\vz})$, where
  $\mathit{lift}(m)$ codes the lift of the term coded by $m$, have a
  head step; and that the target's head reduct $R$ be static, code a
  well-formed term ($\Wfc(R)\ne\z$) and be the code of a term
  ($\IsCode(R)\ne\z$);
\item at every other code, $\Dep(c) = 0$.
\end{itemize}
A code is \emph{certified}, and \emph{certified at level $n$} for a
term $n$, when
\[
  \Ok(c) \;\ldef\; \Dep(c) \ne \z \;\wedge\;
  \Wfc(c) \ne \z \;\wedge\; \IsCode(c) \ne \z ,
  \qquad
  \Ok(n, c) \;\ldef\; \Dep(c) \le n + 1 \;\wedge\; \Ok(c) .
\]
\end{defi}

\noindent
Both certificates are decided at every $c$ with $\Gamma\vdash c\jnat$
(and, for $\Ok(n, c)$, every $n$ with $\Gamma\vdash n\jnat$)
(\textsf{OkI\_B}, \textsf{OkU\_B}).  The check computes $R$ with the
compiled head-step reader instead of projecting the body out of $m$,
so the machine needs no law about composing substitutions; at the
code of $\lambda b$ the reader and the projection agree (\textsf{ngR\_lam\_code}).  Quantifiers
in a code certified at level $n$ therefore nest at most $n+1$ deep,
and each nesting passes through one head step.  The laws hold at every
code $c$ with $\Gamma \vdash c\jnat$, abstract or closed:
\begin{itemize}
\item at a spine node of nonzero depth, the arguments have nonzero
  depth, at most $\Dep(c)$ (\textsf{dep\_kid1\_one},
  \textsf{dep\_kid\_two}); a code of nonzero depth that is neither
  $\forall$-free nor a spine node is a quantifier application
  (\textsf{dep\_leaf});
\item at a quantifier application of nonzero depth whose matrix is not
  $\forall$-free, the matrix passes the nested-leaf check, $R < c$,
  $\Dep(R) \ne \z$, and $\Dep(c) = 1 + \Dep(R)$ exactly
  (\textsf{dep\_leaf\_nested});
\item a $\forall$-free code, and a quantifier application with a
  $\forall$-free matrix, have depth~$1$ (\textsf{iDep\_qf},
  \textsf{dep\_leaf\_flat}).
\end{itemize}
The certificates descend in the same way (\textsf{OkI\_kid1\_one},
\textsf{OkU\_kid1\_one} and their companions).  At a nested leaf,
$\Ok(x + 1, c)$ gives $\Ok(x, R)$ (\textsf{OkU\_nested}), and a nested
leaf is never certified at level $0$ (\textsf{OkU0\_nested}).
$\Ok(n, c)$ gives $\Ok(c)$, and $\Ok(c)$ gives $\Ok(\Dep(c), c)$
(\textsf{okU\_okI}, \textsf{okI\_okU}).

The machine reads $\Dep(R)$ only behind the decided test $R < c$, so
the nested law reads $R$'s value from the table with no reachability
argument.  The test is no restriction at the codes of terms: at the
code of $\lambda b$ the head reduct of the target is the code of $b$
(\textsf{ngR\_lam\_code}), and $\bquo{b} < \bquo{\lambda b} <
\bquo{\cA\cdot(\lambda b)}$ (\textsf{ngR\_genuine\_lt}).

\paragraph{The class.}
Three rules of \rga carry a case-dependent hypothesis
(\S\ref{sec:con-frag}): elimination of a negated universal, induction,
and disjunction elimination.  The first two rules add one hypothesis to a
premise, and disjunction elimination adds one to each of its two
branches.  The five swap rules are the recursor's successor step,
universal elimination, introduction of a negated universal, induction
and equality substitution.  The \emph{static target} of a swap node is
its critical formula, or, at universal elimination and introduction of
a negated universal, the head reduct of the critical formula.

\begin{defi}[\textsf{iHypLvI}, \textsf{iSub1Qf}, \textsf{iStatLvI}, \textsf{iLvClassI}]
\label{def:e-class}
For a proof code $p$,
$\LvC(p) \ldef (\mathit{HypLv}(p) \wedge \mathit{Sub1Qf}(p)) \wedge
\mathit{StatLv}(p)$, where each conjunct is a bounded scan over the
nodes of the proof list coded by $p$:
\begin{itemize}
\item $\mathit{HypLv}(p)$: at every node of the three case-dependent
  rules, every added hypothesis is certified;
\item $\mathit{Sub1Qf}(p)$: at every substitution node, the substituted
  term that the checker reads is $\forall$-free;
\item $\mathit{StatLv}(p)$: at every node of the five swap rules,
  the static target is static and certified; at
  universal elimination and at introduction of a negated universal the
  critical formula has a head step; and at the recursor's successor
  step, introduction of a negated universal and equality substitution
  the swapped term is $\forall$-free.
\end{itemize}
\end{defi}

\begin{thm}[\textsf{hq\_iLvClassI}]
\label{thm:e-class-dec}
If $\Gamma$ is well-formed and $\Gamma \vdash p\jnat$, then $\Gamma
\vdash \LvC(p)\jbool$.
\end{thm}

\noindent
\emph{Proof.}  Each conjunct is decided (\textsf{hq\_iHypLvI},
\textsf{hq\_iSub1Qf}, \textsf{hq\_iStatLvI}), and a conjunction of
decided terms is decided.  \qed

\medskip
The class contains a proof with a quantifier: the witness
$P_{\mathrm w}$ of Theorem~\ref{thm:con-witness}, a three-line proof of
$\forall x.\ x = x$, satisfies $\emptyset \vdash
\LvC(\bquo{P_{\mathrm w}})$ (\textsf{fwP\_in\_lvi}).  No clause of the
class applies at the witness's nodes, so the three scans close at the
witness's code without any certificate being evaluated.

\medskip
$\mathit{HypLv}$ and $\mathit{Sub1Qf}$ are compiled bounded
universals over explicit index trees (\textsf{okIIdx},
\textsf{lvChkIIdx}, \textsf{lvOkIIdx}, \textsf{sub1QfChkIdx}); the
tree \textsf{okIIdx} reads to $\Ok(c)$ (\textsf{ir\_okIIdx\_eq}).
$\mathit{StatLv}$ is a node scan (\textsf{nodeScan}), decided in
well-formed contexts.  The class is a decided class of proof codes,
defined by where quantifiers may occur in a proof.  It is not a class
of sentences of a given arithmetical complexity.

\paragraph{The closed family sentences.}

\begin{defi}[\textsf{iLFI}, \textsf{iLFU}, \textsf{iLFA}]
\label{def:e-lf}
\[
  \LF \;\ldef\; \forall c.\ \Ok(c) \limp \mathit{ALeaf}(c)
  \limp \forall z.\ \mathit{EnvWf}(z) \limp \mathit{EnvQf}(z) \limp
  \mathit{FamDec}(\mathit{arg}(\Subst(z, c))) ,
\]
with implications associating to the right.  For a term $n$, $\LF(n)$
is the same sentence with $\Ok(n, c)$ in place of $\Ok(c)$, and
$\LF_{\forall} \ldef \forall n.\ \LF(n)$.
\end{defi}

\noindent
In words, $\LF$ says: at every certified quantifier application code,
under every environment of well-formed $\forall$-free codes, both
numeral families of the matrix of the instance are decided.  Every
antecedent is a decided test, and the sentence is closed; it mentions
no proof code and no level.  The elimination rule of $\LF$, \textsf{iLFI\_E}, takes
$\Gamma \vdash c\jnat$, $\Gamma \vdash z\jnat$ and the four
antecedents as derivations in $\Gamma$, and returns
$\Gamma \vdash \mathit{FamDec}(\mathit{arg}(\Subst(z, c)))$.  \rga
derives $\LF$ from $\LF_{\forall}$ by instantiating the level at the
code's own depth (\textsf{iLFA\_iLFI}, with \textsf{okI\_okU}), and
$\LF_{\forall}$ from $\LF$ (\textsf{iLFI\_iLFA}, with
\textsf{okU\_okI}); both derivations discharge only decided hypotheses.

\subsection{The level theorem in \rga}
\label{app:oga-lv}

The free-variable form of Theorem~\ref{thm:con-lv}
(\textsf{rga\_con\_lvi\_under\_lf}) is an instance of one shell around
the record fold.

\paragraph{The shell.}
The theorem \textsf{con\_shell} takes a decided scan $S$ of proof
codes and a context formula $H(p)$, both with the variable $v_1$ not
free.  Given the fold's three interfaces of
Appendix~\ref{app:fold-record} at the conjunction $S(p) \wedge H(p)$,
namely the guard interface, the swap interface and the swapped-term
interface, the shell concludes
\[
  H(\vz),\ \vz\jnat \;\vdash\; S(\vz) \limp \neg\Chk(\vz, \bquo{\z = \tofn 1}) .
\]
The implication is introduced over the decided $S(\vz)$, and $H(\vz)$
is only eliminated.  Theorem~\ref{thm:con-lv} is the shell at
$S = \LvC$ and $H = \LF$.

\paragraph{The guard interface.}
\textsf{gdIf\_LI} proves that the guard is decided at every point of
the path, for every accepted code in the class.  A guard item is the
fuel-bounded record of a case-dependent hypothesis under its path
environment.  On the class, $\mathit{HypLv}$ certifies every such
hypothesis (\textsf{lv\_slotI}), and $\mathit{Sub1Qf}$
makes every path environment $\forall$-free.  The \emph{level-certified
core} then decides each item (\textsf{lv\_coreU}): under a
$\forall$-free environment $z$, the fuel-bounded record
$\Rcd^s_k(\Subst(z, h))$ of a code $h$ certified at level $n$ is
decided at every fuel $s$ and value $k$, given the families at the
quantifier applications of $h$ that are themselves certified at level
$n$.  Under a well-formed $\forall$-free $z$ ($\mathit{EnvWf}(z)$ and
$\mathit{EnvQf}(z)$), $\LF(n)$ supplies those families
(\textsf{lv\_core\_iLFU}), and
$\LF$ supplies them at a certified $h$ at the level of its own depth
(\textsf{lv\_core\_iLFI}).  The proof of \textsf{lv\_coreU} is one
internal induction over the fuel, on a motive with the level $n$ as
a parameter, packing the fuel, the value and the position inside
$h$.  Every position
carries its own certificate, so the families are asked for only at
certified positions.  The guard walk itself is
\textsf{whole\_guard\_BLI}.

\paragraph{The swap interface.}
\textsf{swIfQ\_LI} rests on a relational swap layer at a term level
$n$.  At a static source code $C$ certified at level $n$ under a
$\forall$-free environment, a record of the instance transfers to any
related well-formed term code.  The proof is one internal induction on
the fuel bound (\textsf{sl\_levelU}).  The proof spends $\LF$ at two places:
the families at a certified quantifier application, and the decided
core at the item's own depth.  At a nested leaf the layer needs the
same transfer at the reduct $R$, one level down.  The fuel induction
takes it as a HOL premise, the \emph{oracle} \textsf{swNestU} at $n$:
the transfer at every static certified code $R$ with $\Dep(R) \le n$,
in every extension of the context.  At a nested leaf certified at
level $n$, $\Dep(c) = 1 + \Dep(R)$ and $\Dep(c) \le n + 1$ give
$\Dep(R) \le n$ (\textsf{OkU\_nested\_dep}), so the leaf needs no case
split on the level.

One internal induction on the level discharges the oracle
(\textsf{swU\_all}).  The motive of this induction, $\mathit{SwU}(x)$, says: for every
code $C$, environment $z$ and fuel bound $s$, if $z$ satisfies
$\mathit{EnvWf}$ and $\mathit{EnvQf}$, then every record of fuel at
most $s$ of the instance of a static position of $C$ certified at
level $x$ transfers to every related code (the position's test is
\textsf{slAU}).  The motive quantifies the triple
$\langle C, z, s\rangle$ with one universal, and its two antecedents
are decided.
\begin{itemize}
\item \emph{Level $0$.}  The oracle holds vacuously: $\Dep(R) \le 0$
  and $\Dep(R) \ne \z$ contradict each other.
\item \emph{Level $x+1$.}  In the context $x\jnat, \mathit{SwU}(x),
  \Gamma$, the oracle at level $x+1$ is an instance of the hypothesis
  $\mathit{SwU}(x)$: a certified $R$ with $\Dep(R) \le x+1$ is
  certified at level $x$.  The instance moves only atomic facts along
  the projections of the triple, so it needs no transport of formulas.
\end{itemize}
So \rga derives $\mathit{SwU}(l)$ at every level $l$ with $\Gamma
\vdash l\jnat$, under $\LF$
and the swap sites of the fixed proof.  At the top the level is the
code's own depth (\textsf{slTopU}, with \textsf{okI\_okU}).  The swap
and swapped-term interfaces use only $\mathit{StatLv}$,
$\mathit{Sub1Qf}$ and $\LF$ (\textsf{swIfQ\_VLI},
\textsf{xqIf\_VLI}); the forms at $\LvC$ are \textsf{swIfQ\_LI} and
\textsf{xqIf\_LI}.

\paragraph{The soundness companion.}

\begin{thm}[\textsf{rga\_sound\_lvi\_under\_lf}]
\label{thm:e-sound-lv}
Let $\Gamma$ be a well-formed context in which the variable $v_1$ is not
free, and let $p$ be a term in which $v_1$ is not free.  If
$\Gamma \vdash p\jnat$, $\Gamma \vdash \LvC(p)$ and
$\Gamma \vdash \LF$, then
\[
  \Gamma \vdash \forall c.\ \Chk(\lft p, c) \limp \Rcd_1(\Subst(\z, c)) .
\]
\end{thm}

\noindent
\emph{Proof.}  The soundness form of the record fold
(\textsf{fold\_sound\_swIfQP}) at the residual $\LvC(p) \wedge
\LF$, with the three interfaces above.  \qed

\paragraph{Closing the universal.}
The shell gives $\vz\jnat,\ \LF \vdash \LvC(\vz) \limp
\neg\Chk(\vz, \bquo{\z = \tofn 1})$ (\textsf{rga\_con\_lvi\_under\_lf}).
$\LF$ is closed (\textsf{nfree\_iLFI}), so one
$\forall$-introduction at $\vz$ gives the stated form of
Theorem~\ref{thm:con-lv}, $\LF \vdash \forall p.\ \LvC(p)
\limp \neg\Chk(p, \bquo{\z = \tofn 1})$
(\textsf{rga\_con\_lvi\_under\_lf\_all}), and a cut with
$\LF_{\forall} \vdash \LF$ gives the form under $\LF_{\forall}$
(\textsf{rga\_con\_lvi\_under\_lfa\_all}).  Each of these forms is an entailment, a
derivation under the hypothesis.  The implication form adds nothing: \rga neither
proves nor refutes $\LF$, and for every well-formed $S$, $\emptyset
\vdash \LF \limp S$ holds exactly when $\emptyset \vdash S$
(\textsf{iLFI\_imp\_iff}, Theorem~\ref{thm:con-lf-indep}).

\subsection{\rga does not derive the family sentence}
\label{app:oga-notlf}

Theorem~\ref{thm:con-not-lf} (\textsf{rga\_not\_iLFI}) is a
computability argument.  The argument uses no consistency sentence and not
G\"odel's second theorem.

Let $\mathit{ofam}(x)$ be the halting family of the companion
development (\texttt{RGA/RGA\_Omega}): a $\forall$-free matrix whose
instance at $n$ says that program $x$ on input $x$ has not halted by
step $n$.  Let $K$ be the set of programs $x$ that halt on input $x$
(Isabelle \textsf{diagK}).

\begin{thm}[\textsf{iLFI\_decides\_halting}]
\label{thm:e-halting}
If $\emptyset \vdash \LF$, then for every natural $x$: $x \notin K$
if and only if $\emptyset \vdash \forall n.\
\Rcd_1(\mathit{ap}(\bquo{\mathit{ofam}(x)}, \Num\cdot n))$.
\end{thm}

\noindent
\emph{Proof.}  The code $\bquo{\cA\cdot\mathit{ofam}(x)}$ is a
quantifier application with a $\forall$-free matrix, so it is
certified, at depth~$1$ (\textsf{dep\_leaf\_flat}, \textsf{aOfamC\_okI}),
and the empty environment satisfies $\mathit{EnvWf}$ and
$\mathit{EnvQf}$.  \rga derives all of this at the closed code, and
evaluates the $\forall$-free, well-formedness and code tests there
(\textsf{iLFI\_at\_ofam}).  So $\LF$ gives $\emptyset
\vdash \mathit{FamDec}(\bquo{\mathit{ofam}(x)})$, and in particular the
true family $T_x \ldef \forall n.\ \Rcd_1(\mathit{ap}(\bquo{\mathit{ofam}(x)},
\Num\cdot n))$ is decided.  By the disjunction property, \rga derives
$T_x$ or derives $\neg T_x$.  Two facts about the record at the halting
family, \textsf{famT\_ofam\_sound} and \textsf{famT\_ofam\_refuted},
settle the two cases.  If \rga derives $T_x$, then $x \notin K$.
Otherwise, at a step $n$ where the program has halted, $T_x$ gives the
value-$1$ record of the instance at $n$, while \rga refutes that record
(\textsf{ofam\_iRec\_false}), against consistency.  If \rga derives
$\neg T_x$, the witness property gives a numeral $n$ at which \rga
refutes the value-$1$ record of the instance, and then $x \in K$.
Otherwise the program has not halted by step $n$, and \rga derives that
record (\textsf{ofam\_iRec\_true}), again against consistency.  Both
facts turn a witnessed evaluation of the $\forall$-free instance into a
derivation about its record (\textsf{iRec\_of\_wosr}), and the second fact
uses the exclusivity of the record's values.  \qed

\medskip
\emph{Proof of $\emptyset\nvdash\LF$.}  The set of $x$ with $\emptyset
\vdash T_x$ is recursively enumerable, by an explicit proof-search
enumerator (\textsf{renp\_famT\_ofam}).  If \rga derived $\LF$,
Theorem~\ref{thm:e-halting} would make the complement of $K$
recursively enumerable, which it is not (\textsf{not\_rens\_compl\_diagK}).
\qed

\medskip
The obstruction is Turing's: a sound recursively enumerable theory
does not decide the halting problem.  The proof uses \rga's
consistency and its disjunction and witness properties, all theorems of
HOL about \rga, but no consistency sentence.  Since \rga derives $\LF$
from $\LF_{\forall}$, it does not derive $\LF_{\forall}$ either
(\textsf{rga\_not\_iLFA}).  The separation theorem
\textsf{rga\_oga\_separation\_iLFI} pairs the non-derivability with the
\oga derivation of \S\ref{app:oga-oga}.

\paragraph{Nor does \rga refute $\LF$.}
What \oga derives in the empty context, \rga does not refute
(\textsf{rga\_unrefutable\_of\_oga}): an \oga theorem is true in
\oga's $\omega$-semantics (\textsf{oga\_omega\_true}, from
\textsf{oga.omega\_sound}); a refutation in \rga would be true there
too (\textsf{rga\_omega\_sound}); and the $\omega$-value of a term is
unique (\textsf{omv\_not\_inv}, \textsf{omv\_unique}).  \oga
derives $\LF$ and $\LF_{\forall}$ (Theorem~\ref{thm:e-olf}), so \rga
refutes neither sentence (\textsf{iLFI\_unrefutable},
\textsf{iLFA\_unrefutable}).  With the non-derivability above, both
sentences are independent of \rga (\textsf{iLFI\_independent},
\textsf{iLFA\_independent}), and \rga does not prove them decided
(\textsf{iLFI\_undecided}, \textsf{iLFA\_undecided}): by the
disjunction property, a proof of $(\LF)\jbool$ would give a
proof of $\LF$ or of its negation.  \rga has no primitive implication:
$\LF \limp S$ abbreviates $\neg\LF \lor S$.  So, by the disjunction
property again, $\emptyset \vdash \LF \limp S$ gives $\emptyset \vdash
\neg\LF$, which is excluded, or $\emptyset \vdash S$, and the converse
is $\lor$-introduction (\textsf{iLFI\_imp\_iff},
\textsf{iLFA\_imp\_iff}).  Together these facts are
Theorem~\ref{thm:con-lf-indep}.  They are meta theorems, which use
soundness (the disjunction property
and $\omega$-soundness) by nature.

\subsection{\oga and the discharge}
\label{app:oga-oga}

\paragraph{\oga.}
\oga has \rga's 42 rules and one more, ATI (\textsf{pr\_allTI}):
\[
\infrule{\Gamma \vdash \cA\cdot(\lambda((\lft u\cdot\vz)\jbool))}
        {\Gamma \vdash (\cA\cdot u)\jbool}
\]
If every instance of $u$ is decided, the universal $\cA\cdot u$ is
decided.  The rule is sound for an $\omega$-semantics
(\textsf{ati\_omega\_sound}), and every \rga derivation is an \oga
derivation (\textsf{entails\_embed}).  Two derived forms carry every
use of ATI below.  For a body $b$, \textsf{ga\_allTI\_bridge} derives
$(\forall.\,b)\jbool$ from $\forall.\,(b\jbool)$, and
\textsf{ga\_exTI\_bridge} derives $(\exists.\,b)\jbool$, that is
$(\neg\cA\cdot(\lambda\neg b))\jbool$, from the same premise.  An
\oga fact enters an \rga argument only as a context hypothesis that
the \rga derivation eliminates: from $H, \Gamma \vdash c$ in \rga and
$\Gamma \vdashO H$, infer $\Gamma \vdashO c$ (\textsf{o\_cutR}).

\paragraph{One ATI step, exhibited.}
The basic step decides a value record from its fuel-bounded core.  Let
$c$ be a code and $k$ a value, and suppose that \rga decides the core
at a fresh fuel: $s\jnat, \Gamma \vdash \Rcd^s_k(c)\jbool$ for a fresh
variable $s$.
\begin{enumerate}
\item In \rga, $\forall$-introduction at $s$ gives
  $\Gamma \vdash \forall s.\ (\Rcd^s_k(c))\jbool$ (\textsf{rf\_prem}).
\item In \oga, \textsf{ga\_exTI\_bridge} applied to this premise gives
  $\Gamma \vdashO (\exists s.\ \Rcd^s_k(c))\jbool$, which is
  $\Gamma \vdashO \Rcd_k(c)\jbool$, since the value record is the
  existential over the fuel of its core (\textsf{o\_iRec\_B}).
\end{enumerate}
The lemma \textsf{o\_iRec\_of\_core} packages the two steps.  The same pattern,
one level up, decides a numeral family from its instance records: the
records of all instances $\mathit{ap}(u, \Num\cdot n)$ at a fresh $n$
are decided by the step above, run under the numeral's eigenvariable
(\textsf{o\_iRec\_all}); then \textsf{ga\_allTI\_bridge} decides the
true family $\forall n.\ \Rcd_1(\ldots)$, and \textsf{ga\_exTI\_bridge}
the false family $\exists n.\ \Rcd_0(\ldots)$
(\textsf{o\_famDec\_of\_prems}).  So deciding the families of a matrix
takes two nested ATI steps, one over the fuel and one over the
numeral, and needs from \rga only the cores of all numeral instances,
decided at every fuel and both values (\textsf{o\_famDec\_of\_cores}).

\paragraph{The $\forall$-free base.}

\begin{thm}[\textsf{o\_qfRec}, \textsf{o\_qfFam}]
\label{thm:e-qfbase}
\oga derives, in the empty context,
\[
  \forall c.\ \qfs(c) = \z \limp (\Rcd_1(c)\jbool \wedge \Rcd_0(c)\jbool)
  \qquad\text{and}\qquad
  \forall u.\ \qfs(u) = \z \limp \mathit{FamDec}(u) .
\]
\end{thm}

\noindent
\emph{Proof.}  First sentence: at a fresh $c$ with $\qfs(c) = \z$,
Lemma~Q (Theorem~\ref{thm:con-lemq}, \textsf{lemQ\_int}) decides the
core $\Rcd^s_k(c)$ at every fuel $s$ and both values, and the ATI step
over the fuel decides both records.  The implication is introduced
over the decided test $\qfs(c) = \z$, and likewise in the second
sentence.  Second sentence: under the first
sentence as an \rga hypothesis, the instance $\mathit{ap}(u, \Num\cdot
n)$ of a $\forall$-free matrix $u$ at a fresh numeral $n$ is
$\forall$-free (\textsf{qfScan\_mkAppS\_I}), so both of its records are
decided (\textsf{qfFam\_prem}).  ATI over the numeral decides the two
families, and a cut against the first sentence (\textsf{o\_cutR})
removes the hypothesis.  \qed

\medskip
The base needs no totality and no syntactic class of atoms.  The code
of the Liar and every unbounded $\forall$-free search are included.
Under a $\forall$-free environment $z$, the instance of a $\forall$-free
matrix is $\forall$-free, so the families of $\Subst(z, m)$ are decided
too (\textsf{o\_famDec\_qf\_env}).

\paragraph{The nested step.}
Let $m$ pass the nested-leaf check of Definition~\ref{def:e-depth},
with head reduct $R$, and let $z$ satisfy $\mathit{EnvWf}$.
The core of each numeral instance $\mathit{ap}(\Subst(z,m), \Num\cdot
n)$ is decided at every fuel and value as soon as the cores of $R$'s
images are, under the environment that extends $z$ by the numeral
code $\Num\cdot n$ for the bound variable (\textsf{core\_B\_nest}, in
\rga).  The only live row of the record at such an instance is the
application row, whose head reduct is that image; every other row is
excluded by the instance's shape.  The two ATI steps then decide both
families of $\Subst(z, m)$ in \oga (\textsf{o\_famDec\_nest}).

\begin{thm}[\textsf{o\_iLFU}, \textsf{o\_iLFA}, \textsf{o\_iLFI}]
\label{thm:e-olf}
If $\Gamma$ is well-formed and $\Gamma \vdash l\jnat$, then $\Gamma
\vdashO \LF(l)$.  Hence $\emptyset \vdashO \LF_{\forall}$ and
$\emptyset \vdashO \LF$.
\end{thm}

\noindent
\emph{Proof.}  By one internal induction on the level: \oga's
induction rule (\textsf{ir\_nat\_ind}) on the motive $\LF(v_i)$, with
$v_i$ fresh.  In each case \textsf{o\_iLFU\_I} introduces in \oga the
two universals of $\LF(n)$ at fresh variables $c$ and $z$, and the
four implications of $\LF(n)$ over their decided antecedents.  What remains is
$\Delta \vdashO \mathit{FamDec}(\mathit{arg}(\Subst(z, c)))$ for a
context $\Delta$ that derives $\Ok(n, c)$, $\mathit{ALeaf}(c)$,
$\mathit{EnvWf}(z)$ and $\mathit{EnvQf}(z)$.  The instance of a
quantifier application is the quantifier application of the instance
of its matrix, $\mathit{arg}(\Subst(z, c)) = \Subst(z, \mathit{arg}(c))$
(\textsf{aleaf\_img}), and \oga's derivations transport
$\mathit{FamDec}$ along this equation (\textsf{o\_famDec\_transp}).
Split on the decided test $\qfs(\mathit{arg}(c)) = \z$
(\textsf{o\_split0}).  If the matrix is $\forall$-free, the $\forall$-free
base under the $\forall$-free environment $z$ concludes
(\textsf{o\_fam\_flatleaf}).  Otherwise:
\begin{itemize}
\item \emph{Level $0$} (\textsf{o\_iLFU\_0}).  A nested leaf is not
  certified at level $0$ (\textsf{OkU0\_nested}), so this case is
  refuted in \rga.
\item \emph{Level $x+1$} (\textsf{o\_iLFU\_step}).  The induction
  hypothesis $\LF(x)$ sits in the context as an \rga hypothesis, which
  the body of \textsf{o\_iLFU\_I} receives through its \rga weakening
  map.  The matrix passes the nested-leaf check and its head reduct
  $R$ is certified at level $x$ (\textsf{OkU\_nested}).  The
  level-certified core at level $x$, fed with families by $\LF(x)$
  (\textsf{lv\_core\_iLFU}), decides the cores of $R$'s images under
  every extended environment; the extended environment is again
  well-formed and $\forall$-free (\textsf{sigWfGnChk\_updSubM},
  \textsf{qfEnv\_eigen}).  The nested step gives the families.
\end{itemize}
$\LF_{\forall}$ follows by one $\forall$-introduction in \oga
(\textsf{o\_rp\_allI}), and $\LF$ by a cut (\textsf{o\_cutR}) against
\rga's derivation of $\LF$ from $\LF_{\forall}$
(\textsf{iLFA\_iLFI}).  \qed

\medskip
Every implication in this proof is introduced over a decided
antecedent and every universal at an eigenvariable.  ATI is used only
through \textsf{ga\_allTI\_bridge} and \textsf{ga\_exTI\_bridge}, in the
$\forall$-free base and in the nested step.  The per-level proof of
\S\ref{app:oga-perlevel} needs an induction on $N$ in HOL and, at each
level, a cut against the induction hypothesis as an \oga theorem; here
the induction is \oga's own, and its hypothesis is an \rga hypothesis
of the context.  The derivations $\emptyset \vdashO \LF_{\forall}$ and
$\emptyset \vdashO \LF$ use no soundness induction: both derivations are roots of
the dependency check of \S\ref{sec:con-summary}.

\paragraph{The \oga forms of the level theorem.}
Theorem~\ref{thm:con-oga} (\textsf{rga\_con\_lvi\_oga}) is proved in
three steps.  The free-variable form of Theorem~\ref{thm:con-lv}
(\textsf{rga\_con\_lvi\_under\_lf}) embeds into \oga
(\textsf{entails\_embed}); a cut against $\emptyset \vdashO \LF$
(Theorem~\ref{thm:e-olf}, weakened to the context $\vz\jnat$)
removes the hypothesis; and $\forall$-introduction at $\vz$ closes the
universal.  The soundness form follows the same pattern from
Theorem~\ref{thm:e-sound-lv}, with the implication introduced over the
decided $\LvC(\vz)$ (Theorem~\ref{thm:e-class-dec}).

\begin{thm}[\textsf{rga\_sound\_lvi\_oga}]
\label{thm:e-sound-oga}
\[
  \emptyset \vdashO \forall p.\ \LvC(p) \limp
    \forall c.\ \Chk(p, c) \limp \Rcd_1(\Subst(\z, c)) .
\]
\end{thm}

\noindent
Neither theorem is \oga's consistency or \rga's: both are \oga
derivations about \rga's tree checker, restricted to the class
$\LvC$.  Whether \oga proves $\ConR$, the consistency sentence for
all of \rga's proofs, is an open question.

\subsection{The per-level forms}
\label{app:oga-perlevel}

The per-level forms certify codes with one compiled course-of-values
machine per meta level $N$, built from an explicit index tree defined
by primitive recursion on $N$ (\textsf{lvG}); quantifiers in a code
that passes level $N$ nest at most $N+1$ deep.  With that scan
$\Ok_N$ in place of the certificate, the per-level forms define the per-level class
$\LvC_N$ and the closed per-level sentence $\LF_N$, and set $\LF_{\le
N} \ldef \LF_N \wedge \cdots \wedge \LF_0$
(Definitions~\ref{def:e-scan} and~\ref{def:e-class-N}).  The
reach-guarded theorem (\S\ref{app:oga-reach}), the fragment comparison
(\S\ref{app:oga-frag}) and the residual results of
\S\ref{sec:con-res} hold both over $\LvC$ and over each $\LvC_N$.

\begin{defi}[\textsf{iNstC}, \textsf{ngGdT}, \textsf{posOk}]
\label{def:e-scan}
$\mathit{Nst}_N(c) \ldef \cmp(\mathit{lvG}_N)\cdot\z\cdot c$, with
the value $0$ meaning that $c$ passes level $N$.  A code passes level
$N$ when it is $\forall$-free; or it is a spine node and every argument
passes level $N$; or it is a quantifier application $\cA\cdot m$ whose
matrix $m$ passes the \emph{leaf test} of level $N$.
\begin{itemize}
\item At level $0$ the leaf test asks that $m$ be $\forall$-free.  Level
  $0$ is the flat scan \textsf{iFltC}.
\item At level $N+1$ the leaf test asks that $m$ be $\forall$-free, or
  that $m$ pass the nested-leaf check of Definition~\ref{def:e-depth}
  with, in addition, the head reduct $R$ passing level $N$.
\end{itemize}
A code is \emph{certified at level $N$} when
\[
  \Ok_N(c) \;\ldef\; \mathit{Nst}_N(c) = \z \;\wedge\;
  \Wfc(c) \ne \z \;\wedge\; \IsCode(c) \ne \z .
\]
\end{defi}

\noindent
The test $\Ok_N(c)$ is decided at every $c$ with $\Gamma\vdash c\jnat$
(\textsf{posOk\_B}), and descends like $\Ok$ (\textsf{posOk\_kid1\_one},
\textsf{posOk\_kid\_two}, \textsf{posOk\_leaf}).  Quantifiers in a
code certified at level $N$ nest at most $N+1$ deep.  The two
certificates rest on different machines.  At a closed code the two certificates agree
in one direction: if $\emptyset \vdash c\jnat$ and $\emptyset \vdash
\Ok_N(c)$, then $\emptyset \vdash \Ok(\tofn N, c)$
(\textsf{posOk\_closed\_okU}), a rule proved in HOL by induction on the
term that $c$ codes (\S\ref{app:oga-closed}).  At an abstract code the
development relates $\Ok_N(c)$ to $\Ok(\tofn N, c)$ by no theorem.

\begin{defi}[\textsf{iHypLv}, \textsf{iStatLv}, \textsf{iLvClass}, \textsf{iLF}, \textsf{iLFs}]
\label{def:e-class-N}
The per-level class $\LvC_N(p) \ldef (\mathit{HypLv}_N(p) \wedge
\mathit{Sub1Qf}(p)) \wedge \mathit{StatLv}_N(p)$ is
Definition~\ref{def:e-class} with $\Ok_N$ in place of $\Ok$.  The
per-level sentence $\LF_N$ is $\LF$ of Definition~\ref{def:e-lf} with
$\Ok_N$ in place of $\Ok$, and $\LF_{\le N}$ is defined by primitive
recursion on $N$: $\LF_{\le 0} \ldef \LF_0$ and $\LF_{\le N+1} \ldef
\LF_{N+1} \wedge \LF_{\le N}$.
\end{defi}

\begin{thm}[\textsf{rga\_con\_lv\_under\_lf\_all}, \textsf{rga\_con\_lv\_under\_lf}, \textsf{rga\_con\_lv\_oga}, \textsf{iLFs\_independent}]
\label{thm:con-lvN}
For every $N$:
\begin{enumerate}
\item $\LF_{\le N} \vdash \forall p.\ \LvC_N(p) \limp \neg\Chk(p,
  \bquo{\z = \tofn 1})$;
\item $\emptyset \vdashO \forall p.\ \LvC_N(p) \limp \neg\Chk(p,
  \bquo{\z = \tofn 1})$;
\item $\emptyset \nvdash \LF_{\le N}$ and $\emptyset \nvdash
  \neg\LF_{\le N}$, hence $\emptyset \nvdash (\LF_{\le N})\jbool$
  (\textsf{iLFs\_undecided}); and for every well-formed $S$,
  $\emptyset \vdash \LF_{\le N} \limp S$ if and only if $\emptyset
  \vdash S$ (\textsf{iLFs\_imp\_iff}).
\end{enumerate}
\end{thm}

\noindent
Each part is a schema over the meta-level number $N$, one result for
each $N$.  Theorems~\ref{thm:con-lv} to~\ref{thm:con-oga} say the same
for one class and one closed hypothesis, with the level internal.  The
two kinds of class rest on different machines: $\LvC_N$ certifies
with the scan $\Ok_N$, $\LvC$ with the depth machine.

\noindent
$\LvC_N$ is decided (\textsf{hq\_iLvClass}), and its scans are
compiled like those of $\LvC$ (\textsf{posOkIdx}, \textsf{lvChkIdx},
\textsf{lvOkIdx}); $\mathit{StatLv}_N$ is a negated bounded search for
a bad node, like $\mathit{FragV}$ in Appendix~\ref{app:fold}.
$P_{\mathrm w}$ lies in every $\LvC_N$ (\textsf{fwP\_in\_lv}).
Theorem~\ref{thm:con-lvN} is proved as the level theorem, with the
level fixed in HOL.  The theorem's first part is the shell at $S = \LvC_N$ and
$H = \LF_{\le N}$ (\textsf{rga\_con\_lv\_under\_lf}, closed by
$\forall$-introduction as \textsf{rga\_con\_lv\_under\_lf\_all}).
The guard interface \textsf{gdIf\_L} uses the core \textsf{lv\_core}
at level $N$, fed by $\LF_N$.  The swap interface \textsf{swIfQ\_L}
rests on the layer \textsf{swRelL\_holds}, proved by induction on $N$
in HOL: the nested instance at level $N+1$ comes from the layer at
level $N$.  The depth machine replaces that
induction in HOL by the internal induction \textsf{swU\_all} of
\S\ref{app:oga-lv}.

\begin{thm}[\textsf{o\_iLF}, \textsf{o\_iLFs}, \textsf{rga\_not\_iLFs}, \textsf{rga\_oga\_separation\_iLFs}]
\label{thm:e-olf-N}
For every meta level $N$: $\emptyset \vdashO \LF_N$ and
$\emptyset \vdashO \LF_{\le N}$, while $\emptyset \nvdash \LF_{\le N}$.
\end{thm}

\noindent
\emph{Proof.}  By induction on $N$ in HOL, with the cases of
Theorem~\ref{thm:e-olf} at the certificate $\Ok_N$
(\textsf{o\_iLF\_I}).  At level $N+1$ the nested case runs in the
context extended by the hypothesis $\LF_N$, and a cut against the
induction hypothesis $\emptyset \vdashO \LF_N$ removes that hypothesis.
$\LF_{\le N}$ follows by conjunction introduction in \oga
(\textsf{oga\_conj\_provI}).  $\emptyset \nvdash \LF_N$ is the
halting-problem argument of \S\ref{app:oga-notlf} at $\Ok_N$
(\textsf{iLF\_decides\_halting}, \textsf{rga\_not\_iLF}), and
$\LF_{\le N}$ derives $\LF_N$ (\textsf{iLFs\_top}).  \qed

\medskip
As for $\LF$, \rga refutes neither sentence (\textsf{iLF\_unrefutable},
\textsf{iLFs\_unrefutable}), so both sentences are independent of \rga and not
provably decided (\textsf{iLF\_independent}, \textsf{iLF\_undecided},
and Theorem~\ref{thm:con-lvN}(3) for $\LF_{\le N}$).  As for $\LF$, the
internal implication with the antecedent $\LF_{\le N}$ adds nothing
(Theorem~\ref{thm:con-lvN}(3)), a meta theorem proved by the
disjunction property.

\begin{thm}[\textsf{rga\_sound\_lv\_under\_lf}, \textsf{rga\_sound\_lv\_oga}]
\label{thm:e-sound-N}
For every meta level $N$: Theorem~\ref{thm:e-sound-lv} holds with
$\LvC_N(p)$ and $\LF_{\le N}$ in place of $\LvC(p)$ and $\LF$, and
\[
  \emptyset \vdashO \forall p.\ \LvC_N(p) \limp
    \forall c.\ \Chk(p, c) \limp \Rcd_1(\Subst(\z, c)) .
\]
\end{thm}

\noindent
The \oga form of Theorem~\ref{thm:con-lvN}(2)
(\textsf{rga\_con\_lv\_oga}) is proved as Theorem~\ref{thm:con-oga},
with the cut against $\emptyset \vdashO \LF_{\le N}$
(Theorem~\ref{thm:e-olf-N}).

\paragraph{The per-level classes inside $\LvC$.}
\label{app:oga-closed}
Between the per-level classes and $\LvC$ sits a third class.  For a
term $n$, write $\LvC(n, p)$ for the class $\LvC$ with ``certified at
level $n$'' ($\Ok(n, c)$) in place of ``certified''
(\textsf{iLvClassU}): the class at a level that is an object term.

\begin{thm}[\textsf{lv\_closed\_in\_lvi}, \textsf{lv\_closed\_imp\_lvi}, \textsf{lvU\_in\_lvI}, \textsf{rga\_con\_lvu\_under\_lf\_all}, \textsf{rga\_con\_lvu\_oga}]
\label{thm:con-lv-closed}
\leavevmode
\begin{enumerate}
\item For every $N$ and every number $p$: if $\emptyset \vdash
  \LvC_N(\tofn p)$, then $\emptyset \vdash \LvC(\tofn p)$; and
  $\emptyset \vdash \LvC_N(\tofn p) \limp \LvC(\tofn p)$.
\item If $\Gamma$ is well-formed, $\Gamma \vdash n\jnat$, $\Gamma
  \vdash p\jnat$ and $\Gamma \vdash \LvC(n, p)$, then $\Gamma \vdash
  \LvC(p)$.
\item $\LF \vdash \forall n.\ \forall p.\ \LvC(n, p) \limp \neg\Chk(p,
  \bquo{\z = \tofn 1})$, and $\emptyset \vdashO \forall n.\ \forall p.\
  \LvC(n, p) \limp \neg\Chk(p, \bquo{\z = \tofn 1})$.
\end{enumerate}
\end{thm}

\noindent
Part~1 is a rule at closed codes, and a schema of internal
implications, one for each $N$ and each numeral.  Part~2 is a derived
rule at abstract codes and levels: a code certified at level $n$ is
certified, scan by scan.  Part~3 combines part~2 with
Theorems~\ref{thm:con-lv} and~\ref{thm:con-oga}; it is one sentence
over every level and every code, and its implication is introduced
over the decided class $\LvC(n, p)$.  Both kinds of class contain the
accepted $\forall$-free codes (Theorem~\ref{thm:con-lv-contain}) and
$P_{\mathrm w}$ (\textsf{fwP\_in\_lvi}, \textsf{fwP\_in\_lv}).

Theorem~\ref{thm:con-lv-closed} is proved in
\texttt{ConOGA/}\allowbreak\texttt{ConOGA\_Level\_Closed}.  At a closed
code the depth machine can be read with numeral values: the ascent laws
give the depth of a node from the depths of its kids, as the larger of
the two at a node with two kids, the kid's depth at a node with one, and
$d+1$ at a nested quantifier application whose head reduct has depth
$d$ (\textsf{dep\_asc\_two}, \textsf{dep\_asc\_one},
\textsf{dep\_asc\_nested}).  By induction on the term $t$ in HOL: if
$\emptyset \vdash \Ok_N(\bquo t)$, then $\emptyset \vdash
\Dep(\bquo t) = \tofn d$ for some $d$ with $0 < d \le N+1$
(\textsf{posOk\_closed\_dep}).  At a nested leaf the matrix is the code
of a $\lambda$-term (\textsf{lamN\_code}), so the head reduct is the
code of its body, and the machine's code-order test holds there
(\textsf{ngR\_genuine\_lt}).  Hence $\emptyset \vdash \Ok(\tofn N,
c)$ at every closed code $c$ that $\Ok_N$ passes
(\textsf{posOk\_closed\_okU}), and the scans of $\LvC_N(\tofn p)$ give
those of $\LvC(\tofn p)$ node by node (\textsf{iHypLv\_closed\_lvi},
\textsf{iStatLv\_closed\_lvi}), which is part~1's rule
(\textsf{lv\_closed\_in\_lvi}).  Both classes are decided at
numerals, so the disjunction property splits on $\LvC_N(\tofn p)$, and
$\lor$-introduction gives the implication $\LvC_N(\tofn p) \limp
\LvC(\tofn p)$ in either case (\textsf{lv\_closed\_imp\_lvi}).  The
same reading gives $\emptyset \vdash \LvC(\tofn N, \tofn p)$
(\textsf{lv\_closed\_in\_lvu}).

The class $\LvC(n, p)$ (\textsf{iLvClassU}) has the scans of $\LvC$
with $\Ok(n, c)$ in place of $\Ok(c)$, built with a node-scan former
that takes the level term $n$ as a parameter (\textsf{nodeScanP}); it
is decided (\textsf{hq\_iLvClassU}).  Since $\Ok(n, c)$ includes
$\Ok(c)$ as a conjunct, each node check of $\LvC(n, \cdot)$ gives the
node check of $\LvC$ at the same node, in any context
(\textsf{lvNodeUT\_lvNodeIT}, \textsf{vOkTU\_vOkTI}), and the two
node scans give part~2 (\textsf{lvU\_in\_lvI}), a derived rule at an
abstract level and code.  For part~3, in the context of $\LF$,
$v_1\jnat$, $\vz\jnat$ and $\LvC(v_1, \vz)$, part~2 and
Theorem~\ref{thm:con-lv} (\textsf{rga\_con\_lvi\_under\_lf}) give
$\neg\Chk(\vz, \bquo{\z = \tofn 1})$.  The implication is introduced
over the decided $\LvC(v_1, \vz)$ (\textsf{ir\_impI\_bool}), and two
$\forall$-introductions give the sentence
(\textsf{rga\_con\_lvu\_under\_lf\_all}); a cut against $\emptyset
\vdashO \LF$ gives the sentence's \oga form (\textsf{rga\_con\_lvu\_oga}).
Parts 2 and~3 use no soundness theorem (\S\ref{sec:con-summary}).  The
single sentence $\forall p.\ \LvC_N(p) \limp \LvC(p)$ for each $N$
remains open (\S\ref{sec:con-lv}).  At an abstract code that sentence needs an
internal induction on codes in place of the HOL induction on terms, and
the code-order test of a nested leaf inside \rga at an abstract matrix.

\subsection{The reach-guarded theorem and its \oga form}
\label{app:oga-reach}

\paragraph{The reach-indexed family formula.}
The per-level class asks for families at every certified quantifier
application of a case hypothesis.  The formula $\DynFam(p)$
(\textsf{iDynWFD}) asks for the families only where evaluation reaches.  A
\emph{walk} from a code $c$ is a coded list of moves (first argument,
second argument, head reduct, recursor counter, or dispatch at a
value); it is legal when every move spends one unit of fuel and
satisfies the guard of its row of the record.  $\DynFam(p)$ says: if
$p$ is accepted, then at every occupied point of the path tree of $p$,
every ancestor on its path and every case-dependent slot, for every
instantiation of the eigenvariables and every fuel $s$, both numeral
families are decided at every quantifier application that a legal walk
from the fuel-$s$ core of the added hypothesis's instance reaches with
fuel left.  Every antecedent is decided, but the formula as a whole is
not shown decided, and it mentions the proof code $p$.  The
introduction and elimination rules of $\DynFam$ are \textsf{iDynWFD\_I} and
\textsf{iDynWFD\_E}.

\paragraph{The guard from the reach-indexed formula.}
The engine is a core theorem with no shape condition
(\textsf{dyn\_core}): the fuel-bounded core of any code, at every fuel
up to $s$ and every value, is decided, given both numeral families at
every quantifier application that the fuel-$s$ evaluation of the code
reaches with one unit of fuel left.  The proof is one internal
induction on the fuel bound over walk positions
(\textsf{dyz\_level}), whose step is the one-step core decision of
\texttt{ConRGA/Wfd\_Decision.thy}.  With each guard item decided this
way, the guard walk gives the whole guard (\textsf{whole\_guard\_BD})
and the guard interface \textsf{gdIf\_D}.

\paragraph{The reach-guarded theorem.}
The level classes ask for families at every certified position ($\LF$),
or at every position certified at level $N$ ($\LF_{\le N}$); the
reach-guarded theorem asks for the families through $\DynFam$, only where the
evaluation of a case hypothesis's instance reaches a quantifier
application.  Write $\Stat(p)$ for the static half of $\LvC$, the
swap-node and substitution clauses of Definition~\ref{def:con-lv}, and
$\Stat_N(p)$ for the same clauses with the per-level scan $\Ok_N$ in
place of the certificate.

\begin{thm}[\textsf{rga\_con\_statI\_}\allowbreak\textsf{under\_dwfd\_lf}, \textsf{rga\_con\_stat\_}\allowbreak\textsf{under\_dwfd\_lf}]
\label{thm:con-e1}
\leavevmode
\begin{enumerate}
\item $\vz\jnat,\ \DynFam(\vz) \wedge \LF \vdash \Stat(\vz)
  \limp \neg\Chk(\vz, \bquo{\z = \tofn 1})$.
\item For every $N$:
  $\vz\jnat,\ \DynFam(\vz) \wedge \LF_{\le N} \vdash \Stat_N(\vz)
  \limp \neg\Chk(\vz, \bquo{\z = \tofn 1})$.
\end{enumerate}
\end{thm}

\noindent
The guard interface comes from $\DynFam$, a hypothesis about the proof
code, which no theorem shows decided; the swap interface comes from
the closed $\LF$, or $\LF_{\le N}$.  Part~2 contains
Theorem~\ref{thm:con-lvN}(1), and part~1 contains
Theorem~\ref{thm:con-lv}, each containment by one cut (below).  In \oga the swap
half is discharged and $\DynFam(\vz)$ remains the only hypothesis
(Theorem~\ref{thm:e-e1-oga} and the level-free form below).

Theorem~\ref{thm:con-e1}(2) (\textsf{rga\_con\_stat\_under\_dwfd\_lf}) is
the shell at the decided class $S = \Stat_N \ldef \mathit{StatLv}_N
\wedge \mathit{Sub1Qf}$ and the context formula $H(p) = \DynFam(p)
\wedge \LF_{\le N}$.  The guard interface of this instance is \textsf{gdIf\_D}, weakened
to $H$ (\textsf{gdIf\_E1L}); its swap and swapped-term interfaces are
the per-level ones of \S\ref{app:oga-perlevel}, which need only
$\Stat_N$ and $\LF_{\le N}$ (\textsf{swIfQ\_E1L}, \textsf{xqIf\_E1L}).
The soundness companion \textsf{rga\_sound\_stat\_under\_dwfd\_lf} is
the derived rule of Theorem~\ref{thm:e-sound-lv} with $\Stat_N(p)$ and
$\DynFam(p) \wedge \LF_{\le N}$ in place of $\LvC(p)$ and $\LF$.

\paragraph{Containment of the per-level theorem.}

\begin{thm}[\textsf{iHypLv\_iDynWFD}]
\label{thm:e-contain}
Let $\Gamma$ be well-formed.  If $\Gamma$ derives $p\jnat$,
$\mathit{HypLv}_N(p)$, $\mathit{Sub1Qf}(p)$ and $\LF_N$, then
$\Gamma \vdash \DynFam(p)$.
\end{thm}

\noindent
\emph{Proof.}  By the introduction rule \textsf{iDynWFD\_I}.  The core
fact is a correspondence between walks and positions
(\textsf{lyInv\_all}): every legal walk from the instance of a code
certified at level $N$, under a $\forall$-free environment, ends at a
$\forall$-free code or at the instance of a certified node of that code
at the same position.  The proof of \textsf{lyInv\_all} is one internal induction
(\textsf{ir\_nat\_ind}) over the suffixes of the move list.  So every
quantifier application that a walk reaches is the instance of a
certified quantifier application of the hypothesis, where $\LF_N$
supplies the families.  $\mathit{HypLv}_N(p)$ certifies the
hypotheses, and $\mathit{Sub1Qf}(p)$ keeps the path environments
$\forall$-free.  \qed

\medskip
With Theorem~\ref{thm:e-contain}, one cut re-derives
the free-variable form of Theorem~\ref{thm:con-lvN}(1) from Theorem~\ref{thm:con-e1}(2)
(\textsf{rga\_con\_lv\_under\_lf\_via\_stat}, same statement as
\textsf{rga\_con\_lv\_under\_lf}): on the per-level class,
$\LF_{\le N}$ gives $\DynFam(\vz)$, and $\LvC_N(\vz)$ gives
$\Stat_N(\vz)$.  The reach-guarded theorem therefore contains the
per-level theorem.

\paragraph{The \oga form of the reach-guarded theorem.}

\begin{thm}[\textsf{rga\_con\_stat\_oga\_under\_dwfd}]
\label{thm:e-e1-oga}
For every meta level $N$:
$\DynFam(\vz),\ \vz\jnat \vdashO \Stat_N(\vz) \limp \neg\Chk(\vz,
\bquo{\z = \tofn 1})$.
\end{thm}

\noindent
\emph{Proof.}  In \rga, weaken Theorem~\ref{thm:con-e1}(2) to the context
$\LF_{\le N}, \DynFam(\vz), \vz\jnat$ and cut away the theorem's conjunctive
hypothesis, which that context derives.  Then \textsf{o\_cutR} against
$\vz\jnat, \DynFam(\vz) \vdashO \LF_{\le N}$ (Theorem~\ref{thm:e-olf-N},
weakened) removes $\LF_{\le N}$.  \qed

\medskip
The soundness form \textsf{rga\_sound\_stat\_oga\_under\_dwfd} is the
corresponding \oga derived rule: in a well-formed context $\Gamma$ in
which $v_1$ is not free, for a $p$ in which $v_1$ is not free, from
$\Gamma \vdash p\jnat$, $\Gamma \vdash \Stat_N(p)$ and $\Gamma \vdash
\DynFam(p)$ it gives $\Gamma \vdashO \forall c.\ \Chk(\lft p, c) \limp
\Rcd_1(\Subst(\z, c))$.  On the per-level class \oga also discharges
the reach-indexed formula:

\begin{thm}[\textsf{o\_iDynWFD\_lv}]
\label{thm:e-odyn}
Let $\Gamma$ be well-formed.  If $\Gamma$ derives $p\jnat$,
$\mathit{HypLv}_N(p)$ and $\mathit{Sub1Qf}(p)$, then $\Gamma \vdashO
\DynFam(p)$.
\end{thm}

\noindent
\emph{Proof.}  Theorem~\ref{thm:e-contain} with $\LF_N$ as an \rga
hypothesis, and \textsf{o\_cutR} against $\emptyset \vdashO \LF_N$
(Theorem~\ref{thm:e-olf-N}).  \qed

\medskip
With Theorem~\ref{thm:e-odyn}, Theorem~\ref{thm:e-e1-oga} re-derives
Theorem~\ref{thm:con-lvN}(2) (\textsf{rga\_con\_lv\_oga\_via\_stat}, same
statement as \textsf{rga\_con\_lv\_oga}).

\paragraph{The level-free form.}
Theorem~\ref{thm:con-e1}(1) (\textsf{rga\_con\_statI\_under\_dwfd\_lf})
is the same shell at the decided class $\Stat(p) \ldef
\mathit{StatLv}(p) \wedge \mathit{Sub1Qf}(p)$, the static half of
$\LvC$, and the context formula $\DynFam(p) \wedge \LF$.  The guard
interface of this instance is again \textsf{gdIf\_D}, weakened (\textsf{gdIf\_E1LI});
its swap and swapped-term interfaces are those of
the level theorem in \S\ref{app:oga-lv}, which need only $\Stat$ and $\LF$
(\textsf{swIfQ\_E1LI}, \textsf{xqIf\_E1LI}).  The soundness companion
is \textsf{rga\_sound\_statI\_under\_dwfd\_lf}.  Theorem~\ref{thm:e-contain}
holds with the certificate $\Ok$ in place of $\Ok_N$: if a well-formed
$\Gamma$ derives $p\jnat$, $\mathit{HypLv}(p)$, $\mathit{Sub1Qf}(p)$
and $\LF$, then $\Gamma \vdash \DynFam(p)$
(\textsf{iHypLvI\_iDynWFD}).  The proof of \textsf{iHypLvI\_iDynWFD} is the walk correspondence
of \textsf{lyInv\_all} with the descent laws of $\Ok$
(\textsf{lyInvI\_all}).  One cut then re-derives the free-variable
form of Theorem~\ref{thm:con-lv} from part~1
(\textsf{rga\_con\_lvi\_under\_lf\_via\_stat}).  In \oga, a cut
against $\emptyset \vdashO \LF$ (Theorem~\ref{thm:e-olf}) leaves
$\DynFam(\vz)$ as the only hypothesis:
\[
  \DynFam(\vz),\ \vz\jnat \vdashO \Stat(\vz) \limp \neg\Chk(\vz,
  \bquo{\z = \tofn 1})
  \qquad (\textsf{rga\_con\_statI\_oga\_under\_dwfd}) ,
\]
\oga discharges $\DynFam(p)$ on $\LvC$ (\textsf{o\_iDynWFD\_lvi}, as
Theorem~\ref{thm:e-odyn}), and the two derivations give
Theorem~\ref{thm:con-oga} a second time
(\textsf{rga\_con\_lvi\_oga\_via\_stat}).  Like the per-level
theorems of this subsection, these derivations use no soundness theorem
(\S\ref{sec:con-summary}).

\subsection{The fragment and the level classes}
\label{app:oga-frag}

The fragment $\Frag$ of Appendix~\ref{app:fold-frag} and the level
classes, $\LvC$ and every $\LvC_N$, are incomparable, and on accepted
codes in $\Frag$ the level classes add one clause: the head step that $\mathit{StatLv}_N$ and
$\mathit{StatLv}$ ask of the critical formula at universal elimination
and at introduction of a negated universal (Definitions~\ref{def:e-class}
and~\ref{def:e-class-N}).  Theorem~\ref{thm:e-frag-lv} states the
comparison at closed codes for the per-level classes; the paragraph
after its proof gives the same three parts for $\LvC$.  The internal
containments are proved for $\LvC$ and for every $\LvC_N$.

\begin{thm}[\textsf{frag\_accepted\_in\_lv}, \textsf{frag\_not\_lv\_at\_root}, \textsf{frag\_lv\_incomparable}]
\label{thm:e-frag-lv}
\leavevmode
\begin{enumerate}
\item Let $p$ and $c$ be naturals with $\emptyset \vdash \Chk(\tofn p,
  \tofn c)$ and $\emptyset \vdash \Frag(\tofn p)$.  Suppose that at
  every node of $p$ that the checker reads as a universal elimination
  or a negated-universal introduction, \rga proves that the compiled
  head-step reader finds a step of the node's critical formula.  Then
  $\emptyset \vdash \LvC_N(\tofn p)$ for every $N$.
\item There are naturals $p$ and $c$ with $\emptyset \vdash \Chk(\tofn
  p, \tofn c)$ and $\emptyset \vdash \Frag(\tofn p)$, such that the
  checker reads the root of $p$ as a universal elimination and, for
  every $N$, $\emptyset \vdash \neg\LvC_N(\tofn p)$, because \rga does
  not derive the level-$N$ swap check at the root.
\item There is an accepted code that \rga proves in $\Frag$ and
  refutes in every $\LvC_N$, and an accepted code that \rga proves in
  every $\LvC_N$ and refutes in $\Frag$.
\end{enumerate}
\end{thm}

\noindent
\emph{Proof sketch.}  (1) (\textsf{frag\_tree\_in\_lv}, at a code that
the validity test accepts.)  The proof works node by node at the closed
code.  The fragment's node check gives internal facts about the
node: the critical formula and the added hypotheses are
$\forall$-free, and the critical formula is well-formed and a term
code.  Tree acceptance gives the well-formedness and the term-code
property of the added hypotheses.  Each closed term of these checks has
a numeral value (\textsf{closed\_val}, by soundness and completeness),
and a value that \rga proves to be a well-formed, $\forall$-free term
code is the code of a well-formed $\forall$-free term
(\textsf{closed\_code\_of}).  The level checks are then built from
these terms as in the containment of the accepted $\forall$-free codes
(\textsf{qf\_accepted\_in\_lv}): a $\forall$-free term is certified at
every level.  At universal elimination and negated-universal
introduction the hypothesis supplies the head step, and the reduct of a
well-formed $\forall$-free term is again well-formed and
$\forall$-free, hence an admissible static target.  The substitution
scan follows from the static scan of $\Frag$, which has the same
clause.  (2) The witness $P_{\mathrm R}$ (\textsf{fwPR}) proves
$f\cdot\z$ for the recursor matrix $f = \cR\cdot g\cdot h$, with $g =
(\suc\z)\jnat$ and $h = \lambda\lambda\vz$, by universal elimination at
$\z$ of $\forall x.\ f\cdot x$.  The universal comes from
$\forall$-introduction over an induction on the eigenvariable, and the
sub-proof below it is $\forall$-free.  The critical formula $f\cdot\vz$
is $\forall$-free, so $P_{\mathrm R}$ is in $\Frag$
(\textsf{fwPR\_in\_frag}).  The critical formula has no head step, since $\cR$ steps only
at $\z$ and at syntactic successor numerals, so the level-$N$ swap check
fails at the root (\textsf{fwPR\_not\_lv}).  (3) is
\textsf{frag\_lv\_incomparable}, with two further witnesses.  \qed

\paragraph{The comparison for $\LvC$.}
The three parts hold for $\LvC$ in place of every $\LvC_N$, with the
same witnesses (\texttt{ConOGA/}\allowbreak\texttt{ConOGA\_Level\_Ports}).
Part~(1): under the same hypotheses, $\emptyset \vdash
\LvC(\tofn p)$ (\textsf{frag\_accepted\_in\_lvi}); the head-step
hypotheses give $\emptyset \vdash \mathit{HS}(\tofn p)$ node by
node, and the core of Theorem~\ref{thm:con-lv-contain}(2) does the
rest (\textsf{frag\_tree\_in\_lvi\_open}).
Part~(2): at the root of $P_{\mathrm R}$, whose tag is $25$, the head-step
reader finds no step of the critical formula, so the node check of
$\mathit{StatLv}$ fails there (\textsf{vOk1I\_25\_fails}), and \rga
refutes the class at the closed code: $\emptyset \vdash
\neg\LvC(\bquo{P_{\mathrm R}})$ and $\emptyset \nvdash
\LvC(\bquo{P_{\mathrm R}})$ (\textsf{not\_lvi\_at\_node},
\textsf{fwPR\_not\_lvi}, \textsf{frag\_not\_lvi\_at\_root}).
Part~(3) (\textsf{frag\_lvi\_incomparable}): the witness
$P_{\mathrm F}$ (\textsf{fwPF}) is accepted and in $\Frag$, and \rga
refutes it in $\LvC$ (\textsf{fwPF\_not\_lvi}): the universal
elimination of $P_{\mathrm F}$, in the context $\{\cA\cdot\z\}$, has the $\forall$-free
matrix $\z$, and $\z\cdot\vz$ has no head step.  The witness
$P_{\mathrm L}$ (\textsf{fwPL}) is accepted and in $\LvC$
(\textsf{fwPL\_in\_lvi}), and \rga refutes it in $\Frag$: the
universal elimination of $P_{\mathrm L}$ has the matrix $\lambda x.\ \cA\cdot(\lambda y.\
y = y)$, which is not $\forall$-free, while the head reduct
$\cA\cdot(\lambda y.\ y = y)$ is static and certified at depth~$1$.
The same three witnesses separate $\Frag$ from every $\LvC_N$.

\medskip
The head-step clause is the only clause of $\LvC_N$ that $\Frag$ does
not imply at accepted codes, by part~(1), and it is needed, by
part~(2).  Parts~(1) and~(2) are schemata over closed codes.  Part~(1)
also holds at an abstract code, as one \rga sentence, with the
head-step hypothesis written as the decided scan $\mathit{HS}$, and
so does the containment of the accepted $\forall$-free codes
(\textsf{qf\_accepted\_in\_lv} at closed codes), both for $\LvC$
(Theorem~\ref{thm:con-lv-contain}) and for every $\LvC_N$.

\paragraph{The internal containments.}
Theorem~\ref{thm:con-lv-contain} is proved in
\texttt{ConOGA/ConOGA\_Level}, and its per-level form in
\texttt{ConOGA/ConOGA\_Lv\_Int}, without reading the code back and
without an induction over the proof.  The core
(\textsf{lvI\_tree\_of\_frag\_pos}; per level,
\textsf{lv\_tree\_of\_frag\_pos}) derives the class at an
abstract code $p$ that the validity test accepts and that lies in
$\Frag$, given the class's swap check at the nodes read as universal
elimination or negated-universal introduction as a premise.  The core
introduces the three scans of the class node by node, with the
introduction rules of the internal fragment containment
(Appendix~\ref{app:fold-contain}).  At each node, the fragment's node
check makes the critical formula and the added hypotheses
$\forall$-free, and tree acceptance makes them well-formed term codes.
A $\forall$-free code is static (\textsf{iStatC\_of\_qf}), certified
at depth~$1$ (\textsf{iDep\_qf}) and certified at every level $N$
(\textsf{iNstC\_of\_qf}), and the substitution clause is a conjunct of
the fragment's static check (\textsf{iSub1Qf\_of\_iStat}).  The proof
for $\LvC$ is the per-level proof with these two certificate steps
changed.
\begin{itemize}\raggedright
\item \emph{Part~1.}  At an accepted $\forall$-free code the premise
  holds vacuously: \rga refutes that a $\forall$-free node is read as a
  universal elimination or a negated-universal introduction
  (\textsf{qf\_tree\_in\_lvi}; per level, \textsf{qf\_tree\_in\_lv}).
  The sentence follows as in
  Theorem~\ref{thm:con-qf-in-frag} (\textsf{rga\_qf\_in\_lvi\_int};
  per level, \textsf{rga\_qf\_in\_lv\_int}).
\item \emph{Part~2.}  At those nodes $\mathit{HS}$ gives the head step
  of the critical formula, and the head reduct is $\forall$-free
  (\textsf{qfScan\_gHR}) and a well-formed term code
  (\textsf{gHR\_wf\_gn}), hence static and certified
  (\textsf{frag\_tree\_in\_lvi\_open}, \textsf{rga\_frag\_in\_lvi\_int}; per level,
  \textsf{frag\_tree\_in\_lv\_open}, \textsf{rga\_frag\_in\_lv\_int}).
  The law \textsf{gHR\_wf\_gn} holds at every well-formed term code with
  a head step.  Its proof is a course-of-values induction inside \rga
  along the compiled head-step reader, with one case for each kind of
  head step: $\beta$, the recursor at $\z$, the recursor at a
  successor, and a step in the head of an application.
\item \emph{Part~3.}  $\mathit{HS}(p)$ is a node scan whose check, at
  the two rules, is the head-step clause of $\mathit{StatLv}$, so the
  static scan of $\LvC$ gives $\mathit{HS}$ at every code
  (\textsf{iHdStep\_of\_lvi}; per level, \textsf{iHdStep\_of\_lv}).
\item \emph{Part~4.}  $\mathit{HS}$ is decided at $\bquo{P_{\mathrm
  R}}$.  The positive case of this decision would give $\LvC_0(\bquo{P_{\mathrm R}})$ by
  the per-level part~2, which \rga refutes (\textsf{fwPR\_not\_lv}), so
  \rga refutes $\mathit{HS}(\bquo{P_{\mathrm R}})$
  (\textsf{fwPR\_not\_hdstep}).
\item \emph{Part~5.}  The mechanized proof reads the class at the
  closed code, as in part~(2) of the comparison for $\LvC$ above
  (\textsf{fwPR\_not\_lvi}).  Parts~3 and~4 give the same conclusion
  by the head-step argument of \S\ref{sec:con-break} on where the fold
  stops.
\end{itemize}

\paragraph{Limits and open questions.}
\rga does not derive the residual $\Res$
(Definition~\ref{def:con-res}) from membership in $\Frag$, in $\LvC$ or
in any per-level class (Theorem~\ref{thm:con-res-not-frag}).  \rga does
not derive $\Res$ from $\LvC$ under the hypothesis $\LF$, or from
$\LvC_N$ under $\LF_{\le N}$.  Nor does \rga derive $\Res$ from these
classes under any hypotheses true in the
$\omega$-semantics (Theorem~\ref{thm:con-res-lf}, proved in
Appendix~\ref{app:fold-res}).  These non-derivability results are
theorems, not gaps.  The
reach-guarded theorem (\S\ref{app:oga-reach}) and the comparison with
$\Frag$ (\S\ref{app:oga-frag}) hold both for $\LvC$ and for every
$\LvC_N$.  The containment shown in \S\ref{app:oga-reach} relates
theorems: the reach-guarded theorem contains the level theorem and each
per-level theorem.  Between the classes: a closed code that \rga
proves in $\LvC_N$ it also proves in $\LvC$, \rga proves every numeral
instance of $\LvC_N(p) \limp \LvC(p)$, and \rga derives $\LvC(p)$ from
$\LvC(n, p)$ at an abstract level and code (Theorem~\ref{thm:con-lv-closed},
\S\ref{app:oga-closed}).  Whether \rga derives, for each $N$, the
single sentence $\forall p.\ \LvC_N(p) \limp \LvC(p)$ over an
abstract code is an open question.  So is the proof-theoretic strength
of the level classes, that is, which arithmetic proofs have codes in
$\LvC$.

\endgroup

\section{The traps catalog}
\label{app:traps}

A \emph{trap} is a machine-checked negative result: a theorem that
refutes a step or a design that looks natural, often because it is
valid classically, but fails in \rga.  Some traps are theorems about
\rga itself, such as the failure of the unconditional deduction
theorem.  The other traps
are counterexamples to a candidate motive, invariant, residual or
machine for the internal proofs of Sections~\ref{sec:inside}
to~\ref{sec:con}, each exhibited at a concrete proof, term or code.
The development keeps these counterexamples because each trap explains why the final
proof takes the path it does.  A reader who tries the natural step can
look up the counterexample instead of rediscovering it, and a proposed
alternative can be checked against the precise shape of the
obstruction.

\paragraph{Where the traps live.}
The two logical traps of \S\ref{app:traps-logic} are theorems of the
theory \textsf{RGA/RGA\_Paradox}.  Apart from the later T01 theorems
of \S\ref{app:traps-t01}, all the other traps are in eight trap theories:
\textsf{RGA/}\allowbreak\textsf{RGA\_Traps}, \textsf{RGA\_Kit/}\allowbreak\textsf{RGA\_Kit\_Traps},
\textsf{RGA\_Decider/}\allowbreak\textsf{RGA\_Decider\_Traps},
\textsf{RGA\_Checker/}\allowbreak\textsf{RGA\_Checker\_Traps},
\textsf{RGA\_Record/}\allowbreak\textsf{RGA\_Record\_Traps},
\textsf{Shadows/}\allowbreak\textsf{Shadows\_Traps}, \textsf{ConRGA/}\allowbreak\textsf{ConRGA\_Traps} and
\textsf{ConOGA/}\allowbreak\textsf{ConOGA\_Traps} (paths relative to \texttt{hol/}, with the
extension \texttt{.thy} omitted).  For each trap, the theory text states
the tempting route, the theorems that refute it, and the design choice
of the final proof that they forced.  Labels such as T01 are the
development's own, taken from the theory texts, and are not
consecutive; each trap can be found by searching for its label.  Names resolve in the development
(Section~\ref{sec:mech}, ``Artifact'').  The theorems of \S\ref{app:traps-t01} live
outside the trap theories, in
\texttt{ConOGA/}\allowbreak\texttt{ConOGA\_Frag\_QF} and
\texttt{ConOGA/}\allowbreak\texttt{ConOGA\_Residual\_Frag}; the summary
of trap T01 in \texttt{ConRGA\_Traps} points to these theorems.

\paragraph{Models in HOL.}
Several traps are stated on a model in HOL of an internal construction
(the development calls such a model a \emph{shadow}): a HOL definition
that mirrors an internal formula of \rga, such as a motive for the
consistency fold, with HOL's own logic in place of \rga's.  A
refutation at the model is cheap.  Such a refutation is not a theorem about the
internal formula, but it locates the failure before the internal
formula is built, and the final internal proofs avoid each refuted
shape.

\paragraph{Not repeated here.}
The theorems of [\citen{goedel-paper}, Section ``G\"odel's proof in four forms''], which reads G\"odel's proof
in four forms, and
Theorem~\ref{thm:pl-no-uniform} are negative results too, but they
concern the classical ingredients of G\"odel's proof rather than design
choices, and the main text states them in full.

\subsection{Discharge, excluded middle, and the satisfaction wall}
\label{app:traps-logic}

Three theorems explain the main restrictions on reasoning in \rga.  The
first two theorems are about the proof system; the third theorem is the obstruction that
shaped the motives of the internal consistency proofs of
Section~\ref{sec:con}.  All three theorems turn on terms without a grounded
value, and the Liar is the standard example.  Write $L$ for the Liar, the term-level fixed point whose
evaluation unfolds in two head steps to $\neg L$.  $L$ has no grounded
value (\textsf{liar\_no\_value}): a value $r$ of $L$ would transfer
along the two steps to $\neg L$, so $L$ would have both the value $r$
and its negation, against the uniqueness of values.

\begin{thm}[\textsf{rga\_deduction\_theorem\_fails}]
\label{thm:tr-dt}
$L \vdash L$, but $\emptyset \nvdash \neg L \lor L$, that is,
$\emptyset \nvdash L \limp L$.
\end{thm}

\noindent
\emph{Proof idea.}  The first judgment is an instance of the hypothesis
rule.  If $\emptyset \vdash \neg L \lor L$, the companion paper's
soundness theorem (\textsf{rga\_provable\_wosr}) would give the
disjunction the grounded value $\ctrue$.  The disjunction's clause needs a
true disjunct or two false ones, and either way a value of $L$, which
$L$ does not have.

\emph{The trap.}  Classically, from $\Gamma, a \vdash c$ one infers
$\Gamma \vdash a \limp c$.  In \rga hypothesizing is free and
discharging is not.  An implication is a disjunction, and at the empty
context \rga proves a disjunction only by proving one of its disjuncts
(Theorem~\ref{thm:br-mech}(1)), so a theorem $a \limp c$ commits \rga to
refuting $a$ or proving $c$.  \rga's rule for discharging a
hypothesis needs the hypothesis to be provably decided (\textsf{ir\_impI\_bool} in
\textsf{RGA/}\allowbreak\textsf{RGA\_Sound}).  Otherwise $a \limp c$
is derived by refuting $a$ or by proving $c$ outright
(\textsf{ir\_impI1}, \textsf{ir\_impI2} in
\textsf{RGA/}\allowbreak\textsf{RGA\_Proof}), and at the empty context
these are the only ways.
So the deduction theorem holds exactly at decided hypotheses: for a
well-formed $a$, every entailment $a, \Gamma \vdash c$ gives $\Gamma
\vdash a \limp c$ if and only if $\Gamma \vdash a \lor \neg a$
(\textsf{deduction\_theorem\_exact}, \S\ref{sec:inside-codes}), and at a
hypothesis $\exists s.\, D(s)$ with $D$ decided an entailment is
equivalent to one sentence (\textsf{sigma\_deduction\_theorem}).
This restriction on discharge is why the consistency results of Section~\ref{sec:con} are stated
as entailments, derivations under hypotheses, and why the internal checker keeps
hypothesis contexts as coded data rather than as an object-level
implication.  The checker form of the trap is T18 in
Table~\ref{tab:tr-logic}: the one-judgment proof of $L \vdash L$ is
valid, but its encoded leaf $L \limp L$ is unprovable.

\begin{thm}[\textsf{rga\_liar\_lem\_explosive}, \textsf{rga\_lem\_scheme\_explosive}]
\label{thm:tr-lem}
For every context $\Gamma$ and every well-formed $c$,
$L \lor \neg L,\ \Gamma \vdash c$.  Consequently, if every context
derived $a \lor \neg a$ for every well-formed $a$, every context would
derive $\z = \tofn 1$.
\end{thm}

\noindent
\emph{Proof idea.}  Eliminate the hypothesis $L \lor \neg L$ by cases.
In the branch $L$, two conversion steps turn $L$ into $\neg L$, and the
pair $L$, $\neg L$ derives every well-formed $c$.  In the branch $\neg
L$, the same two steps run backwards and give $L$.  The second
statement is a cut of the first statement against the assumed instance at $L$.

\emph{The trap.}  Classical case analysis on an arbitrary sentence is
unavailable in \rga, and adding excluded middle as an axiom scheme
would make \rga inconsistent.  A case split on a single sentence $a$
needs $a \lor \neg a$, which is the statement $a\jbool$ that $a$ is
decided, so \rga splits on $a$ only once it has derived $a\jbool$.  The
case split of [\citen{goedel-paper}, theorem~\textsf{goedelG\_of\_expl\_d3g\_ctx}] is an example:
it splits on the checker's verdict $\Chk(v, S_G)$ at a variable $v$
with $v\jnat$, which is decided.

The third theorem is the general form of
Theorem~\ref{thm:tr-dt}.  [\citen{goedel-paper}, Section ``One mechanism''], on the one mechanism behind
the failures of the implications, cites this general form.

\begin{thm}[\textsf{sat\_wall}, \textsf{sat\_wall\_imp}, \textsf{sat\_wall\_closure}]
\label{thm:tr-wall}
Let $G$ and $V$ be terms.
\begin{enumerate}
\item If neither $G$ nor $V$ has a grounded value, then
  $\emptyset \nvdash \neg G \lor V$, that is, $\emptyset \nvdash G \limp V$.
\item Let $S$ be a set of well-formed terms, none of which has a
  grounded value.  Every term built from members of $S$ with $\neg$
  and $\lor$ is neither provable nor refutable: $\emptyset \nvdash t$
  and $\emptyset \nvdash \neg t$.
\end{enumerate}
\end{thm}

\noindent
\emph{Proof idea.}  Part~1 is the argument of
Theorem~\ref{thm:tr-dt} with $G$ and $V$ in place of $L$.  For part~2,
an induction over the construction shows that the strong-Kleene clauses
of $\neg$ and $\lor$ cannot produce a value from arguments without one.
A well-formed term with no value is independent: a proof of it, or of
its negation, would give it a value by soundness
(\textsf{no\_value\_independent}).  Since $\land$, $\limp$ and $\neq$
are abbreviations over $\neg$ and $\lor$, part~2 covers every formula
built from such terms with \rga's connectives.

\begin{cor}[\textsf{sat\_wall\_liar}, \textsf{sat\_wall\_eigen}, \textsf{owsjt\_conj\_false\_at\_liar}]
\label{cor:tr-wall}
$\emptyset \nvdash \neg L \lor L$, and $\emptyset \nvdash L\jnat \limp
L\jnat$.  The reading of the judgment $L \vdash L$ as ``under every
environment, the hypotheses and the conclusion are both true'' is
false.
\end{cor}

\noindent
\emph{The trap.}  The companion paper's soundness proof shows, for
every judgment $H \vdash c$ of a valid proof, that under every
environment, if all hypotheses in $H$ are true then so is $c$.  The
natural internal motive copies that shape as a formula: a guard ``the
hypotheses hold'' implying the conclusion's truth.  Hypothetical
positions of valid proofs can have instances where neither the guard
nor the consequent has a value.  At $L \vdash L$ that instance is the
identity environment.  At the eigen-hypothesis position $\vz\jnat
\vdash \vz\jnat$ that every $\forall$-introduction adds, the instance is the
environment that maps $\vz$ to $L$, since free variables range over
terms without a value as well as numerals.  The corollary gives both instances, and
Theorem~\ref{thm:tr-wall} shows that no propositional rearrangement of
the guard and the consequent escapes.  The final motive of
Section~\ref{sec:con-frag} is therefore a guarded disjunction over a
value-indexed record, with guards only on the hypotheses that a proof
rule introduces case by case.  The guard half of the residual of
\S\ref{sec:con-res} assumes the decidedness of those guards.

\subsection{The unguarded swap body (trap T01)}
\label{app:traps-t01}

Section~\ref{sec:con-res} states consistency of full \rga under the
residual $\Res$, whose swap half asks for the \emph{guarded swap body}
at every occupied path point of the five value-congruence rules
(Definition~\ref{def:con-res}).  The body has four decided tests in
front: $X$ and $Y$ are well-formed, and they are codes of terms.  The
\emph{unguarded} body omits the four tests.  For a critical formula
$C$ and an environment code $z$ it says: for every variable $w$, codes
$X$ and $Y$, values $v$ and $k$, and fuel $s$, if the environment is
well-formed and $\Rcd_v(X)$ and $\Rcd_v(Y)$ hold, then the fuel-bounded
record $\Rcd^s_k$ of $C$'s instance at $z$ updated with $w \mapsto X$
gives the record $\Rcd_k$ of $C$'s instance at $z$ updated with
$w \mapsto Y$.  This unguarded body is what the fold's swap arms consume, and the body looks
harmless, because $X$ and $Y$ carry records of the same value.  The
four theorems below refute the unguarded body at a concrete critical formula.

Write $J$ for a fixed ill-formed term, and take
\[
  X = (\z = \z), \qquad Y = (\z = \z) \lor J, \qquad
  C = (\lambda(\z = \z))\cdot v_i .
\]
$X$ has the record value $\ctrue$, and so has $Y$, since the disjunction
clause reads one true disjunct.  The instance of
$C$ at $X$ is $(\lambda(\z = \z))\cdot(\z = \z)$, which has the record
value $\ctrue$ after one head step.  The instance at $Y$ has no record at
any value and fuel (\textsf{swapCex\_no\_record\_CY}): the head step
refuses a $\beta$-redex whose argument is ill-formed, and every other
clause of the record is excluded by the shape of the term.

\begin{thm}[\textsf{swap\_body\_unguarded\_refutable}]
\label{thm:tr-t01-body}
Let $\Gamma$ be well-formed, and let $\Gamma \vdash z\jnat$ and $\Gamma
\vdash z = \z$, so that $z$ is the empty environment.  If $\Gamma$
derives the unguarded swap body at the critical formula
$\bquo{(\lambda(\z = \z))\cdot v_i}$ and the environment $z$, then
$\Gamma$ derives every well-formed formula.
\end{thm}

\noindent
\emph{Proof idea.}  Instantiate the body at $X$, $Y$, $v = k = \ctrue$ and
a fuel $s$ at which the instance at $X$ has its record.  The
well-formedness of the empty environment and the three records are
derivable, so the body yields $\Rcd_1$ at the instance at $Y$, whose
record \textsf{swapCex\_no\_record\_CY} refutes.  The instance codes
are computed by the closed-code law of the substitution machine for
codes of bounded arity, which covers the ill-formed $Y$
(\textsf{iSubM\_closed\_ar}).

The critical formula $C$ occurs in a genuine proof, as HOL derives that proof:

\begin{thm}[\textsf{swapCexProof\_proves}]
\label{thm:tr-t01-proof}
The five-judgment list
\[
  \emptyset \vdash (\lambda(\z = \z))\cdot\z,\quad
  \emptyset \vdash (\lambda\vz)\cdot\z = \z,\quad
  \emptyset \vdash (\lambda(\z = \z))\cdot((\lambda\vz)\cdot\z),\quad
  \emptyset \vdash \z = \z,\quad
  \emptyset \vdash \z\jnat
\]
is a valid proof of $\emptyset \vdash (\lambda(\z = \z))\cdot\z$.
\end{thm}

\noindent
In HOL the root step is an equality substitution from the second and
third judgments, at the template $(\lambda(\z = \z))\cdot\vz$, so the
HOL derivation's critical formula is $C$ with $i = 0$.  The checker
reads this root differently.  The checker's node reader reports the first firing
row of the node table (\textsf{ruleTag\_code\_cases},
Appendix~\ref{app:fold-contain}), and the conversion row (row~20)
fires before the equality-substitution row (row~40): the root's
$\beta$-reduct $\z = \z$ is itself a judgment of the proof, in the
same context.  So the root's rule tag is at most~20
(\textsf{swapCexProof\_root\_conv}).

\begin{thm}[\textsf{residual\_unguarded\_refutable}]
\label{thm:tr-t01-res}
Suppose $\Gamma$ is well-formed, $\Gamma \vdash p\jnat$, $\Gamma
\vdash \mathit{TreeValid}(p)$ (the tree checker's validity test), and
$\Gamma$ derives the unguarded residual of $p$.  Suppose further that
$\Gamma$ derives three facts about the root node of $p$: it is
occupied, it carries the rule tag of equality substitution, and its
critical formula is $\bquo{(\lambda(\z = \z))\cdot v_i}$.  Then
$\Gamma$ derives every well-formed formula.
\end{thm}

\noindent
\emph{Proof idea.}  At the root the path environment is the empty one
(\textsf{path\_env\_root}).  The swap half of the residual, read at the
root, gives the unguarded body there, which
Theorem~\ref{thm:tr-t01-body} refutes after a transport along the
equation for the critical formula.

Theorem~\ref{thm:tr-t01-res} takes the three facts about the root node
as hypotheses.  At the code of the proof of
Theorem~\ref{thm:tr-t01-proof} the root is read as a conversion: its
rule tag is at most~$20$ (\textsf{swapCexProof\_root\_conv}), so no
consistent $\Gamma$ derives that the root carries the tag of equality
substitution, $40$.  A longer proof has a root that carries that tag.  Let $P_{\mathrm{T01}}$
(\textsf{t01P}) be the nine-judgment proof that derives
$(\lambda\vz)\cdot\z = \z$ and $(\lambda(\z = \z))\cdot((\lambda\vz)
\cdot\z)$ in the context $\forall x.\ x = x$, brings both judgments to the empty
context by cuts against $\emptyset \vdash \forall x.\ x = x$ (the three
judgments of $P_{\mathrm w}$, Theorem~\ref{thm:con-witness}), and
concludes $(\lambda(\z = \z))\cdot\z$ from them by equality
substitution.  The root's reduct $\z = \z$ now occurs only in the
context $\forall x.\ x = x$, so the conversion row does not fire at the
root.

\begin{thm}[\textsf{t01P\_accepted}, \textsf{t01P\_tag0}, \textsf{t01P\_crit}, \textsf{t01P\_residual\_}\allowbreak\textsf{unguarded\_refutable}]
\label{thm:tr-t01-genuine}
\leavevmode
\begin{enumerate}
\item $\emptyset \vdash \Chk(\bquo{P_{\mathrm{T01}}},
  \bquo{(\lambda(\z = \z))\cdot\z})$.
\item The checker reads the root of $P_{\mathrm{T01}}$ as equality
  substitution (rule tag exactly~40), and for some $i$, \rga proves
  that the root's critical formula is $\bquo{(\lambda(\z = \z))\cdot
  v_i}$.
\item If $\Gamma$ is well-formed and derives the unguarded residual of
  $\bquo{P_{\mathrm{T01}}}$, then $\Gamma$ derives every well-formed
  formula; in particular \rga does not derive that residual
  (\textsf{t01P\_residual\_}\allowbreak\textsf{unguarded\_underivable}).
\end{enumerate}
\end{thm}

\noindent
\emph{Proof idea.}  (1) The proof is valid, and the checker is faithful.
(2) No row before 40 fires at the root: rows 0--24 find no reduct or
redex of the root among the judgments of the empty context, and rows
25--39 are refuted one by one (no judgment of the empty context is a
universal with the root's function part, the root is no negation or
universal, the empty context has no negation, no $a\jnat$ and no
weakening, and no earlier judgment concludes the root).  Row 40 fires.  The checker's equality-substitution
reader returns a template $s$ with root $= s[v_i{:=}\z]$ and a judgment
$s[v_i{:=}(\lambda\vz)\cdot\z]$ of the empty context; the only template
with both properties is $(\lambda(\z = \z))\cdot v_i$
(\textsf{t01\_template}).  (3) is Theorem~\ref{thm:tr-t01-res}, with acceptance from~(1) and
the three facts about the root from~(2).

Theorem~\ref{thm:tr-t01-body} does not depend on
Theorem~\ref{thm:tr-t01-genuine}: it refutes the unguarded body at $C$
directly.

\begin{thm}[\textsf{swap\_body\_guarded\_vacuous}]
\label{thm:tr-t01-guarded}
For the ill-formed $Y = (\z = \z) \lor J$ above, every $X$ with
$\Gamma \vdash X\jnat$, and any well-formed critical formula,
environment, variable, values and fuel, $\Gamma$ derives the guarded
swap body.
\end{thm}

\noindent
\emph{Proof idea.}  \rga computes the well-formedness test at the
closed code of $Y$ and derives that the test's value is $0$.  The four tests are
decided, so the body's implications can be introduced one at a time
(\textsf{ir\_impI\_bool} in \textsf{RGA/}\allowbreak\textsf{RGA\_Sound}); the third test's hypothesis contradicts the
computed value, and the remaining implication follows by
$\neg$-elimination.

\emph{The trap and what it forced.}  An endpoint whose hypothesis
derives every formula at some proofs says nothing about those proofs,
and the unguarded body derives every formula at the critical formula
$C$.  The guarded body is therefore what $\Res$ asks for.  Two facts make the
guarded endpoint the stronger theorem: the unguarded residual implies
the guarded one (\textsf{iRes\_unguarded\_iRes}), and consistency under
the unguarded residual follows from Theorem~\ref{thm:con-res} by one
cut (\textsf{rga\_con\_under\_res\_unguarded}).
Theorem~\ref{thm:tr-t01-guarded} shows that this counterexample
does not refute the guarded body.  \rga derives the guarded residual
at the code of $P_{\mathrm w}$, but not at the code of every genuine
proof: at the proof $P_3$ of \S\ref{sec:con-frag} the residual's swap half is not
derivable (Theorems~\ref{thm:con-res-frag}
and~\ref{thm:con-res-not-frag}).  Trap T02 in Table~\ref{tab:tr-fold} is the
same defect in a HOL model of the residual.

\subsection{The catalog}
\label{app:traps-table}

Tables~\ref{tab:tr-logic} to~\ref{tab:tr-fold} list every trap of the
development, one row per trap, grouped by theme.  The first column
gives the trap's label where the theory text uses one.  The second column
gives the Isabelle names of the trap's main theorems and the theory
that holds them.  The third column states the tempting step, and the
fourth column what the theorems say, in the theory's own terms.  Rows marked
\emph{model} are stated on a model in HOL.  Notation: $L$ is the Liar,
$\Rcd$ the value record of Definition~\ref{def:con-rcd}, $\Subst$ the
compiled substitution, $\Dec$ the compiled decider, and $\Chk$ the
tree checker.

\begingroup
\footnotesize
\renewcommand{\arraystretch}{1.15}

\begin{longtable}{@{}>{\raggedright\arraybackslash}p{0.05\linewidth}>{\raggedright\arraybackslash}p{0.24\linewidth}>{\raggedright\arraybackslash}p{0.24\linewidth}>{\raggedright\arraybackslash}p{0.38\linewidth}@{}}
\caption{Traps in the logic and the proof system.}
\label{tab:tr-logic}\\
\hline
label & theorems (theory) & tempting step & what the theorems say \\
\hline
\endfirsthead
\hline
label & theorems (theory) & tempting step & what the theorems say \\
\hline
\endhead
\hline
\endfoot
--- & \textsf{rga\_deduction\_}\allowbreak\textsf{theorem\_fails} (\textsf{RGA/}\allowbreak\textsf{RGA\_Paradox})
  & From $\Gamma, a \vdash c$, infer $\Gamma \vdash a \limp c$.
  & $L \vdash L$, but $\emptyset \nvdash L \limp L$
    (Theorem~\ref{thm:tr-dt}); the step is valid exactly at decided
    $a$. \\
--- & \textsf{rga\_liar\_lem\_}\allowbreak\textsf{explosive},
  \textsf{rga\_lem\_scheme\_}\allowbreak\textsf{explosive} (\textsf{RGA/}\allowbreak\textsf{RGA\_Paradox})
  & Reason by cases on any sentence; add excluded middle.
  & The hypothesis $L \lor \neg L$ derives every well-formed formula;
    excluded middle for all well-formed terms derives $\z = \tofn 1$
    (Theorem~\ref{thm:tr-lem}). \\
T06 & \textsf{sat\_wall}, \textsf{sat\_wall\_imp}, \textsf{sat\_wall\_}\allowbreak\textsf{closure},
  \textsf{sat\_wall\_liar}, \textsf{sat\_wall\_eigen},
  \textsf{owsjt\_conj\_}\allowbreak\textsf{false\_at\_liar} (\textsf{Shadows/}\allowbreak\textsf{Shadows\_Traps})
  & State the internal soundness motive as ``the hypotheses hold implies
    the conclusion holds''.
  & An implication between terms without a grounded value is not
    provable, and any $\neg$/$\lor$ combination of well-formed terms
    without a grounded value is neither provable nor refutable; instances at the
    Liar hypothesis and at the eigen-hypothesis $L\jnat$
    (Theorem~\ref{thm:tr-wall}, Corollary~\ref{cor:tr-wall}). \\
T18 & \textsf{enc\_leaf\_refuted}, \textsf{enc\_nc\_refuted},
  \textsf{enc\_imp\_no\_refl}, \textsf{enc\_witness\_valid}
  (\textsf{RGA\_Kit/}\allowbreak\textsf{RGA\_Kit\_Traps})
  & Encode each judgment $H \vdash c$ of a checked proof as the object
    formula $\bigwedge H \limp c$.
  & The one-judgment proof of $L \vdash L$ is valid, but
    $\emptyset \nvdash L \limp L$ and $\emptyset \nvdash \neg(L \land
    \neg L)$; so implication is not reflexive on all well-formed terms.
    The checker keeps contexts as coded data. \\
--- & \textsf{cnry\_valid\_not\_canon}, \textsf{cnry\_afree\_normal}
  (\textsf{RGA/}\allowbreak\textsf{RGA\_Traps})
  & Read the normal-form theorems as properties of every valid proof,
    or normalize a proof by inspecting its own judgments.
  & A five-judgment valid proof of $(\lambda(\z\jnat))\cdot\z$, which
    introduces $\cA\cdot(\lambda(\z\jnat))$ with a genuine certificate and
    eliminates it again, is not canonical, nor in the $\forall$-free
    normal form, although a proof of that form with the same conclusion
    exists.  The normal-form theorems are existence statements. \\
--- & \textsf{beta\_erases\_to\_A},
  \textsf{conv\_reaches\_}\allowbreak\textsf{nonA\_conclusion} (\textsf{RGA/}\allowbreak\textsf{RGA\_Traps})
  & Define truth by dispatching on the syntactic head: the quantifier
    clause exactly at formulas of the shape $\cA\cdot u$.
  & For well-formed $u$, $(\lambda\vz)\cdot(\cA\cdot u)$ head-steps to
    $\cA\cdot u$, and every derivation of $\cA\cdot u$ gives a derivation of
    $(\lambda\vz)\cdot(\cA\cdot u)$ by conversion.  The decider reduces
    first and dispatches afterwards. \\
\end{longtable}

\begin{longtable}{@{}>{\raggedright\arraybackslash}p{0.05\linewidth}>{\raggedright\arraybackslash}p{0.24\linewidth}>{\raggedright\arraybackslash}p{0.24\linewidth}>{\raggedright\arraybackslash}p{0.38\linewidth}@{}}
\caption{Traps in truth predicates, records and families.}
\label{tab:tr-sem}\\
\hline
label & theorems (theory) & tempting step & what the theorems say \\
\hline
\endfirsthead
\hline
label & theorems (theory) & tempting step & what the theorems say \\
\hline
\endhead
\hline
\endfoot
T09 & \textsf{qf\_record\_undecided}, \textsf{var\_record\_}\allowbreak\textsf{undecided},
  \textsf{neg\_num\_undecided}, \textsf{fob\_strictly\_}\allowbreak\textsf{inside\_qf}
  (\textsf{Shadows/}\allowbreak\textsf{Shadows\_Traps}; model)
  & Every $\forall$-free term has a grounded value, so consistency of
    the $\forall$-free proofs follows cheaply.
  & The Liar is $\forall$-free and well-formed but has no grounded value;
    nor have a free variable and $\neg\tofn 2$.  A decided syntactic
    class excludes the Liar.  The final proofs carry decidedness as a
    decided scan ($\Frag$) and a fuel-bounded record. \\
T08 & \textsf{twin\_gap}, \textsf{twin\_permits},
  \textsf{swap\_core\_}\allowbreak\textsf{not\_enough} (\textsf{Shadows/}\allowbreak\textsf{Shadows\_Traps}; model)
  & Make the value record single-valued by refuting only the boolean
    twin of the value.
  & At a true term the twin conjunct refutes only the value 0; if a
    term had the record values 1 and 2 at once, which the conjunct does
    not exclude, the false clause of $a = \tofn 1$ would hold.  A swap
    delivered at the record's core gives only the core at the target.
    The record carries ``no other value'' in full. \\
T10 & \textsf{valued\_guard\_}\allowbreak\textsf{liar\_no\_value},
  \textsf{vinst\_liar\_undecided}, \textsf{valued\_family\_}\allowbreak\textsf{not\_affirmable},
  \textsf{valued\_family\_}\allowbreak\textsf{undecided} (\textsf{Shadows/}\allowbreak\textsf{Shadows\_Traps}; model)
  & State the family of a universal over valued codes: every instance
    $u\cdot a$ whose argument $a$ has a value is true.
  & The guard ``$L$ has a value'', the term $L\jnat$, has no value;
    for $u = \lambda x.\,x = \z$ the guarded instance at $L$ has no
    value and the affirming clause fails; for $u = \lambda x.\,x = x$
    both clauses fail.  The final families range over numerals. \\
T11 & \textsf{soundEB\_false}, \textsf{cert\_escapes\_}\allowbreak\textsf{code\_level},
  \textsf{guard\_fuel\_escapes} (\textsf{Shadows/}\allowbreak\textsf{Shadows\_Traps}; model)
  & Induct on a bound for proof codes, with guards and consequents
    that are decider verdicts at fuels bounded by the level.
  & If the fuels are bounded by $F(y)$ with $F(y) \ge 2$, the motive is
    false at every level $y$ that contains the code of a two-judgment
    proof of $\vz\jnat \vdash (\suc\vz)\jnat$.  If they are free, then
    for every bound $B$ an environment makes the hypothesis
    $\cA\cdot\vz$ affirmed while every proof of its open instance has a
    code above $B$, and some such environment needs a guard fuel above
    $B+1$.  Fuel is quantified existentially in the final
    proofs. \\
T13 & \textsf{afree\_rfree\_}\allowbreak\textsf{not\_decided},
  \textsf{oga\_all\_liar\_}\allowbreak\textsf{undecided},
  \textsf{oga\_guarded\_}\allowbreak\textsf{family\_undecided} (\textsf{ConOGA/}\allowbreak\textsf{ConOGA\_Traps})
  & Let \oga discharge ``the witness families are decided'', with a
    syntactic freeness condition as the base.
  & Not every closed, well-formed term free of $\cA$ and of the recursor
    is decided in \oga (the Liar); $\forall x.\,L$ is not decided in
    \oga; a guarded family $\cA\cdot(\lambda(g \limp c))$ is not decided
    in \oga when $g$ and $c$ have no value in \oga's $\omega$-model at
    some instance and no instance is false. \\
T05 & \textsf{pair\_affirm\_A\_}\allowbreak\textsf{eq0\_never},
  \textsf{pair\_affirming\_}\allowbreak\textsf{half\_incomplete},
  \textsf{refuting\_half\_}\allowbreak\textsf{affirms\_A\_eq\_n0},
  \textsf{pmRdVA\_not\_}\allowbreak\textsf{monotone}, \textsf{pmRdVB\_not\_}\allowbreak\textsf{monotone},
  \textsf{pair\_both\_halves\_}\allowbreak\textsf{incomplete},
  \textsf{pair\_blind\_to\_}\allowbreak\textsf{wosr\_truth} (\textsf{RGA\_Record/}\allowbreak\textsf{RGA\_Record\_Traps})
  & Use as internal truth motive the verdict of a pair machine: an
    affirming half that searches only for certificates and a refuting
    half that searches only for refuted instances.
  & The affirming half never affirms $(\cA\cdot u) = \z$, yet \rga proves
    $(\cA\cdot(\lambda\z)) = \z$, which the refuting half affirms;
    reading the other half as a fallback is not monotone, in either
    priority order; a provable well-formed formula is affirmed by
    neither half at any fuel, and so is a formula that is true in the
    grounded semantics.  The final
    record has no verdicts. \\
\end{longtable}

\begin{longtable}{@{}>{\raggedright\arraybackslash}p{0.05\linewidth}>{\raggedright\arraybackslash}p{0.24\linewidth}>{\raggedright\arraybackslash}p{0.24\linewidth}>{\raggedright\arraybackslash}p{0.38\linewidth}@{}}
\caption{Traps in the coding, the decider and the checker.}
\label{tab:tr-mach}\\
\hline
label & theorems (theory) & tempting step & what the theorems say \\
\hline
\endfirsthead
\hline
label & theorems (theory) & tempting step & what the theorems say \\
\hline
\endhead
\hline
\endfoot
R08 & \textsf{iSubM\_zero\_}\allowbreak\textsf{not\_identity}
  (\textsf{RGA\_Kit/}\allowbreak\textsf{RGA\_Kit\_Traps})
  & $\Subst$ at the empty environment is the identity on every total
    numeral.
  & Not every $X$ with $\emptyset \vdash X\jnat$ has $\emptyset \vdash
    \Subst(\z, X) = X$: at $\tofn 3$, which codes no term, $\Subst$
    returns $\tofn 1$.  The identity law holds at codes of well-formed
    terms. \\
--- & \textsf{liftX\_truncates\_}\allowbreak\textsf{arity3},
  \textsf{prDecTabTr\_}\allowbreak\textsf{redexAr3\_affirms},
  \textsf{prDecTab\_}\allowbreak\textsf{redexAr3\_refuses},
  \textsf{iPrDecV\_}\allowbreak\textsf{redexAr3\_refuses},
  \textsf{redexAr3\_no\_wosr},
  \textsf{truncating\_decider\_}\allowbreak\textsf{contradicts\_converse},
  \textsf{iPrDecV\_not\_}\allowbreak\textsf{truncating}
  (\textsf{RGA\_Decider/}\allowbreak\textsf{RGA\_Decider\_Traps})
  & Build the decider table on substitution coders that are correct on
    well-formed codes only.
  & On an ill-formed equation node $e$ with three operands, the two-slot
    lift coder truncates the code to that of $\z = \z$; the table built
    on that coder affirms $(\lambda\vz)\cdot e$ at fuel 3, while the final table
    and $\Dec$ refuse that term and the term has no grounded value.  A compiled decider
    faithful to the truncating table, an introduction rule of the shape
    of Theorem~\ref{thm:in-wtrue-intro}, and a converse without side
    conditions are jointly contradictory. \\
--- & \textsf{sEq01\_motive\_}\allowbreak\textsf{not\_inductive}
  (\textsf{RGA\_Decider/}\allowbreak\textsf{RGA\_Decider\_Traps})
  & Prove that the decider never affirms $\z = \tofn 1$ by a fuel
    induction with that statement as the motive.
  & A table whose value cells all read ``value 0'' up to the code of
    $\z = \tofn 1$ does not affirm $\z = \tofn 1$, but one step of the
    table machine applied to that table does.  The proof uses two slot invariants
    instead. \\
T17 & \textsf{hctx\_bound\_}\allowbreak\textsf{bridge\_refuted},
  \textsf{hctx\_list\_}\allowbreak\textsf{bridge\_refuted} (\textsf{RGA\_Checker/}\allowbreak\textsf{RGA\_Checker\_Traps})
  & Read a hypothesis-context code as a counted pair, with the length in
    the first Cantor projection and the list in the second.
  & At the code of the set $\{1, 2\}$ the first projection is not the
    list's length, and the second projection is not the list's code.  A
    set code is the code of its sorted list. \\
R15 & \textsf{first\_proj\_zero\_}\allowbreak\textsf{not\_empty\_list}
  (\textsf{RGA\_Checker/}\allowbreak\textsf{RGA\_Checker\_Traps})
  & Test that a context is empty by testing that its first projection
    is 0.
  & The code of the list $[1]$ has first projection 0 and length 1.  The
    checker tests the coded length. \\
T15 & \textsf{p0\_frame\_refuted}, \textsf{p0PC\_not\_genuine},
  \textsf{p0PC\_treeAccept}, \textsf{p0\_guard\_}\allowbreak\textsf{internally\_zero},
  \textsf{p0PC\_old\_accept\_}\allowbreak\textsf{new\_reject} (\textsf{RGA\_Checker/}\allowbreak\textsf{RGA\_Checker\_Traps})
  & Assume that every code accepted by the checker spine, with empty root
    context, is the code of a proof.
  & A spine without a node-code test accepts a code listing two
    judgments that conclude $\z\jnat$, and \rga derives the acceptance;
    yet the code is not the code of any proof, and \rga derives that the
    in-range test returns 0.  The final spine rejects the code. \\
T16 & \textsf{env\_pin\_key\_}\allowbreak\textsf{refuted},
  \textsf{subm\_shift\_no\_}\allowbreak\textsf{finite\_witness},
  \textsf{env\_code\_}\allowbreak\textsf{independent\_refuted},
  \textsf{segDom\_upd\_out\_}\allowbreak\textsf{refuted}; repair
  \textsf{subm\_shift\_cancel} (\textsf{RGA\_Checker/}\allowbreak\textsf{RGA\_Checker\_Traps})
  & Cancel a lift under a binder by extending the coded environment by a
    fixed entry, or by a list computed from the environment alone.
  & For any two finite lists $\zeta, \zeta'$ some term $t$ has
    $\sof{\zeta}{\lft t} \neq \sof{\zeta'}{t}$; an internal slot of that
    shape derives $\z = \tofn 1$; an update outside an initial segment
    of keys leaves no segment.  Shifting the keys and prepending the new
    entry cancels the lift on covered terms. \\
\end{longtable}

\begin{longtable}{@{}>{\raggedright\arraybackslash}p{0.05\linewidth}>{\raggedright\arraybackslash}p{0.24\linewidth}>{\raggedright\arraybackslash}p{0.24\linewidth}>{\raggedright\arraybackslash}p{0.38\linewidth}@{}}
\caption{Traps in the consistency fold and its residual.}
\label{tab:tr-fold}\\
\hline
label & theorems (theory) & tempting step & what the theorems say \\
\hline
\endfirsthead
\hline
label & theorems (theory) & tempting step & what the theorems say \\
\hline
\endhead
\hline
\endfoot
T07 & \textsf{branch\_unguarded\_}\allowbreak\textsf{false},
  \textsf{branch\_filter\_wall}, \textsf{left\_sibling\_}\allowbreak\textsf{not\_in\_suffix}
  (\textsf{Shadows/}\allowbreak\textsf{Shadows\_Traps}; model)
  & A path-indexed motive without guards; or a guard that is the
    hypothesis's own truth; and an induction over suffixes of the proof
    list.
  & A valid proof contains the position $\z = \tofn 1 \vdash \z = \tofn 1$
    inside a branch of a disjunction elimination, where the conclusion
    is false at the identity environment; at a Liar branch the filtered
    motive $L \limp L$ is unprovable; in a valid proof ending in a cut,
    the suffix at the hypothesis node does not contain the cut's first
    premise. \\
T21 & \textsf{P3\_no\_position\_}\allowbreak\textsf{order} (\textsf{RGA\_Record/}\allowbreak\textsf{RGA\_Record\_Traps})
  & Index the fold by list position.
  & In the valid three-judgment proof of $\z\jnat$ by a cut, the cut's
    first premise can only be position 2, its second premise only
    position 1, and position 1's only premise is again position 2.  So
    no ranking of positions puts the first premise below every node of
    the second premise's derivation, which contains the first premise.  The fold indexes
    the derivation tree by post-order rank. \\
T04 & \textsf{allI\_in\_branch\_}\allowbreak\textsf{guarded} (\textsf{Shadows/}\allowbreak\textsf{Shadows\_Traps}; model)
  & The decidedness of the guard is needed at one site per rule only,
    where the case hypothesis is introduced.
  & In a derivation that splits $\emptyset \vdash \z\jnat \lor \z\jnat$
    by cases and concludes $\forall x.\,x = x$ by
    $\forall$-introduction in each branch, the introduction node in the
    first branch has a non-empty guard.  The residual decides the whole
    guard at every occupied path point. \\
T12 & \textsf{dyWit\_not\_static}, \textsf{dyWit\_reach},
  \textsf{dyn\_strictly\_}\allowbreak\textsf{larger}, \textsf{dyX\_selfreach},
  \textsf{dyX\_no\_level} (\textsf{Shadows/}\allowbreak\textsf{Shadows\_Traps}; model)
  & State the decidedness of the witness families over a static
    syntactic shape; or state the swap half over the reach by a level of
    nesting.
  & $(\lambda\vz)\cdot(\cA\cdot u)$ is not static for any leaf predicate,
    yet its evaluation reaches only the quantifier node $\cA\cdot u$; a
    recursion whose step builds a quantifier re-reaches its own
    quantifier, in the model's reach relation, and has no level.
    Theorem~\ref{thm:con-e1} uses the reach for the guard half, and the
    static class $\Stat$ with $\LF$, or $\Stat_N$ with the per-level formula $\LF_{\le N}$, for the swap half. \\
T02 & \textsf{swapAt\_unguarded\_}\allowbreak\textsf{fails},
  \textsf{sw5\_unguarded\_unsat}, \textsf{shRes\_unguarded\_}\allowbreak\textsf{unsat}
  (\textsf{Shadows/}\allowbreak\textsf{Shadows\_Traps}; model)
  & Swap two codes of the same record value without asking that they be
    well-formed.
  & An ill-formed code can carry a record: $\cR\cdot\z\cdot J\cdot\z$ has
    value 0 for the ill-formed $J$.  The swap from $\z$ to that code fails
    in the context $(\lambda(\z = \z))\cdot{-}$, for every predicate that
    satisfies the record laws, and a valid three-judgment proof exhibits
    the context, so the unguarded residual is unsatisfiable there. \\
T01 & \textsf{swap\_body\_}\allowbreak\textsf{unguarded\_refutable},
  \textsf{residual\_}\allowbreak\textsf{unguarded\_refutable},
  \textsf{swapCexProof\_proves}, \textsf{swap\_body\_}\allowbreak\textsf{guarded\_vacuous}
  (\textsf{ConRGA/}\allowbreak\textsf{ConRGA\_Traps});
  \textsf{t01P\_residual\_}\allowbreak\textsf{unguarded\_refutable}
  (\textsf{ConOGA/}\allowbreak\textsf{ConOGA\_Residual\_Frag})
  & State the swap body of the residual without the four tests on the
    swapped codes.
  & The unguarded body at the empty environment and the critical formula
    $(\lambda(\z = \z))\cdot v_i$ derives every formula; so does the
    unguarded residual at an accepted code whose root is read as an
    equality substitution with that critical formula, and the accepted
    proof $P_{\mathrm{T01}}$ is such a code; the guarded body holds
    vacuously at the counterexample (\S\ref{app:traps-t01}). \\
T03 & \textsf{junk\_position}, \textsf{junk\_hyp\_}\allowbreak\textsf{position},
  \textsf{wvP3\_overasks} (\textsf{ConOGA/}\allowbreak\textsf{ConOGA\_Traps})
  & State the decidedness of the witness families at path points, with
    formers that quantify over every quantifier position of a case
    hypothesis or swap target.
  & The position reader descends into the code of a variable.  For every
    $w$, a $\forall$-free equation $v_V = \z$ with a suitable index $V$
    passes the static scan and has the code of $\cA\cdot w$ at a position,
    so the leaf former asks for the decidedness of the families of an
    arbitrary code.  The level theorems assume the closed level
    formula $\LF$ (per level, $\LF_{\le N}$) instead. \\
\end{longtable}
\endgroup

\section{The endpoint table}
\label{app:endpoints}

This appendix maps the paper's results to the Isabelle development.
Theory files are given relative to \texttt{hol/}, with the extension
\texttt{.thy} omitted.  Every name below has its defining line (a
theorem, lemma, corollary, definition, abbreviation, recursive function
or inductive definition) in exactly one theory file of the development
(Section~\ref{sec:mech}, ``Artifact''), except \textsf{ir\_disj\_prop},
which \texttt{PGA\_Con} also defines, for the propositional system
\pga; the table gives the definition for \rga.
The last subsection, \S\ref{app:endpoints-statements}, shows the
statements of the headline results as Isabelle prints them.

\paragraph{Results of the companion formalization.}
Fourteen names belong to the companion formalization, the Isabelle
formalization of \rga that the companion paper~\cite{rga-paper}
describes: these names are defined in the companion formalization's theories (the base theories and
those under \texttt{RGA/}) and were there before this paper's
development began.  This paper states these results where later proofs need
their exact form: the fact about chosen indices and the compiler's
correctness and totality (Section~\ref{sec:inside} and
Appendix~\ref{app:machines}), the semantics behind the soundness
invariant and the invariant itself (Appendix~\ref{app:pillars}), the
failures of the unconditional deduction theorem and of excluded middle
(Section~\ref{sec:inside} and Appendix~\ref{app:traps}), and the
disjunction and witness properties behind the failing implications
([\citen{goedel-paper}, Section ``G\"odel's proof in four forms'']).  The companion paper itself states the
compiler's correctness and totality, the soundness invariant, the
disjunction and witness properties and the two failures; the other names,
the compiler's value equation at numerals, the fact about chosen
indices and the semantic notions behind the invariant, are facts and
definitions of the formalization.  Table~\ref{tab:g-paper} marks all
fourteen ``companion''.

\paragraph{Names added after the campaign.}
The statement fingerprints of the endpoints of
\S\ref{app:endpoints-list} were recorded when the proof campaign and
its cleanup closed.  Two hundred and forty-seven names that the paper's
numbered results cite were added after that point, and the
fingerprints show that no endpoint statement has changed since
(Appendix~\ref{app:machines-audit}, ``Fingerprints'').

\subsection{The paper's theorems}
\label{app:endpoints-paper}

Table~\ref{tab:g-paper} lists every theorem, corollary and definition
of the paper, 212 in all (177 theorems, 5 corollaries
and 30 definitions), in the order in which they occur, grouped by
section and appendix.  For each result the table gives the Isabelle names in the result's
heading and, after them, the further names that the result's statement cites or
that the text attaches to the result (the main text gives at most two names in
a heading).
The table has 515 rows for 479 distinct names, since
36 names occur in two results, for instance a theorem of
{}[\citen{goedel-paper}, Section ``G\"odel's theorems for \rga''] and its full form in
in [\citen{goedel-paper}].  For each name the table gives the theory
file, a status, and the endpoint group of Table~\ref{tab:g-endpoints}
when the name is one of the development's endpoints.  The status is
``campaign'' for a name present when the campaign closed (229 rows),
``added'' for a name added after the campaign (254 rows), ``trunk''
for a name of this development, of either kind, that has since moved
into the companion formalization's theories (17 rows), and
``companion'' for a name of the companion formalization (15 rows).

\begingroup
\footnotesize

\begin{longtable}{@{}>{\raggedright\arraybackslash}p{0.15\linewidth}>{\raggedright\arraybackslash}p{0.33\linewidth}>{\raggedright\arraybackslash}p{0.25\linewidth}>{\raggedright\arraybackslash}p{0.10\linewidth}>{\raggedright\arraybackslash}p{0.06\linewidth}@{}}
\caption{The paper's theorems, corollaries and definitions, with their
  Isabelle names.}
\label{tab:g-paper}\\
\hline
result & Isabelle name & theory & status & group \\
\hline
\endfirsthead
\hline
result & Isabelle name & theory & status & group \\
\hline
\endhead
\hline
\endfoot
\multicolumn{5}{@{}l}{\textbf{Section~\ref{sec:inside}}} \\*
Theorem~\ref{thm:in-dt} & \textsf{deduction\_}\allowbreak\textsf{theorem\_exact} & \textsf{RGA/\allowbreak RGA\_Sound} & trunk &  \\
 & \textsf{sigma\_deduction\_}\allowbreak\textsf{theorem} & \textsf{RGA/\allowbreak RGA\_Sound} & trunk &  \\
 & \textsf{rga\_deduction\_}\allowbreak\textsf{theorem\_fails} & \textsf{RGA/\allowbreak RGA\_Paradox} & companion &  \\
 & \textsf{prc\_sentence\_iff} & \textsf{RGA\_Goedel/\allowbreak Goedel\_Cut} & added &  \\
Corollary~\ref{cor:in-contrapos} & \textsf{dec\_contrapos} & \textsf{RGA/\allowbreak RGA\_Sound} & trunk &  \\
 & \textsf{sigma\_contrapos} & \textsf{RGA/\allowbreak RGA\_Sound} & trunk &  \\
 & \textsf{sigma\_reductio} & \textsf{RGA/\allowbreak RGA\_Sound} & trunk &  \\
\end{longtable}
\endgroup

\subsection{The endpoints}
\label{app:endpoints-list}

The development designates 131 theorems as its \emph{endpoints}: the
results that a paper about it would state.  Tier~A holds 28 headline
results (groups A1 to A10) and Tier~B 103 further results (groups B1 to
B11), among them five design checks on models in HOL (group B11).  The
development fingerprints each endpoint's statement.  Since the
endpoints were fixed, every fingerprint is unchanged except one, which
records a strengthening:
\textsf{rga\_sound\_sch} (Theorem~\ref{thm:pl-sound}) lost its premise
that the conclusion $c$ is well-formed, which follows from the proof
premise.

Table~\ref{tab:g-endpoints} lists the endpoints by group, with each
name, the name's theory, and the numbered results of the paper that state the endpoint.
41 endpoints are stated in this way (21 of Tier~A and
20 of Tier~B), in 32 results.  The seven Tier~A endpoints
without a reference are cited in the text.  Two of these seven are soundness forms that
accompany a stated consistency theorem (\textsf{rga\_sound\_stat\_}\allowbreak\textsf{under\_dwfd\_lf}
and \textsf{rga\_sound\_stat\_}\allowbreak\textsf{oga\_under\_dwfd}, Appendix~\ref{app:oga}).
Four are forms of stated theorems that the paper names beside those theorems: the
packed form of the $\forall$-free soundness theorem
(\textsf{rga\_sound\_}\allowbreak\textsf{qf\_int}, Appendix~\ref{app:fold}), completeness and
truth reflection for the opaque provability predicate $\Prv$ (\textsf{rga\_complete\_}\allowbreak\textsf{sch} and
\textsf{rga\_truth\_}\allowbreak\textsf{reflect\_sch}, Section~\ref{sec:pillars}), and the
packed form of internal completeness (\textsf{rga\_complete\_}\allowbreak\textsf{int},
Appendix~\ref{app:pillars}).  The seventh is the
proof of G\"odel's second theorem through a search witness
(\textsf{rga\_goedel2\_}\allowbreak\textsf{via\_search}, [\citen{goedel-paper}, Section ``G\"odel's theorems for \rga'']).  Of the 103
Tier~B endpoints, 33 are cited somewhere in the paper; most of
the others belong to results that this paper does not discuss.

\begingroup
\footnotesize

\begin{longtable}{@{}>{\raggedright\arraybackslash}p{0.40\linewidth}>{\raggedright\arraybackslash}p{0.38\linewidth}>{\raggedright\arraybackslash}p{0.14\linewidth}@{}}
\caption{The 131 endpoints, by group.}
\label{tab:g-endpoints}\\
\hline
name & theory & stated in \\
\hline
\endfirsthead
\hline
name & theory & stated in \\
\hline
\endhead
\hline
\endfoot
\multicolumn{3}{@{}l}{\textbf{A1.}\ Consistency and soundness of full \rga under the residual} \\*
\textsf{rga\_con\_under\_res} & \textsf{ConRGA/\allowbreak Guarded\_Residual} & \ref{thm:con-res} \\
\textsf{rga\_sound\_under\_res} & \textsf{ConRGA/\allowbreak Guarded\_Residual} & \ref{thm:c-res-sound} \\
\multicolumn{3}{@{}l}{\textbf{A2.}\ The level-$N$ class under the closed level formula} \\*
\textsf{rga\_con\_lv\_under\_lf} & \textsf{ConOGA/\allowbreak ConOGA\_Fold} & \ref{thm:con-lvN} \\
\textsf{rga\_sound\_}\allowbreak\textsf{lv\_under\_lf} & \textsf{ConOGA/\allowbreak ConOGA\_Fold} & \ref{thm:e-sound-N} \\
\multicolumn{3}{@{}l}{\textbf{A3.}\ The static class under the reach-indexed family, in \rga and in \oga} \\*
\textsf{rga\_con\_stat\_}\allowbreak\textsf{under\_dwfd\_lf} & \textsf{ConOGA/\allowbreak ConOGA\_Fold} & \ref{thm:con-e1} \\
\textsf{rga\_sound\_stat\_}\allowbreak\textsf{under\_dwfd\_lf} & \textsf{ConOGA/\allowbreak ConOGA\_Fold} &  \\
\textsf{rga\_con\_stat\_}\allowbreak\textsf{oga\_under\_dwfd} & \textsf{ConOGA/\allowbreak ConOGA\_Fold} & \ref{thm:e-e1-oga} \\
\textsf{rga\_sound\_stat\_}\allowbreak\textsf{oga\_under\_dwfd} & \textsf{ConOGA/\allowbreak ConOGA\_Fold} &  \\
\multicolumn{3}{@{}l}{\textbf{A4.}\ The level-$N$ class in \oga} \\*
\textsf{rga\_con\_lv\_oga} & \textsf{ConOGA/\allowbreak ConOGA\_Fold} & \ref{thm:con-lvN} \\
\textsf{rga\_sound\_lv\_oga} & \textsf{ConOGA/\allowbreak ConOGA\_Fold} & \ref{thm:e-sound-N} \\
\multicolumn{3}{@{}l}{\textbf{A5.}\ The fragment theorem} \\*
\textsf{rga\_con\_frag\_int} & \textsf{ConRGA/\allowbreak Wfd\_General} & \ref{thm:con-frag} \\
\textsf{rga\_sound\_frag\_}\allowbreak\textsf{int\_prenex} & \textsf{ConRGA/\allowbreak Wfd\_General} & \ref{thm:c-frag-prenex} \\
\textsf{rga\_sound\_frag\_int} & \textsf{ConRGA/\allowbreak Wfd\_General} & \ref{thm:con-frag} \\
\multicolumn{3}{@{}l}{\textbf{A6.}\ The $\forall$-free fragment and its checker} \\*
\textsf{rga\_con\_qf\_int} & \textsf{RGA\_QF/\allowbreak QF\_Last\_Arms} & \ref{thm:con-qf} \\
\textsf{rga\_sound\_qf\_int} & \textsf{RGA\_QF/\allowbreak QF\_Last\_Arms} &  \\
\textsf{rga\_con\_qf\_meta} & \textsf{RGA\_QF/\allowbreak QF\_Last\_Arms} & \ref{thm:c-qf-meta} \\
\textsf{iChkQf\_faithful} & \textsf{RGA\_Checker/\allowbreak Chk\_QF\_Faithful} & \ref{thm:con-qf-faithful} \\
\multicolumn{3}{@{}l}{\textbf{A7.}\ Schematic soundness and completeness} \\*
\textsf{rga\_sound\_sch} & \textsf{RGA\_Decider/\allowbreak Decider\_Truth} & \ref{thm:pl-sound} \\
\textsf{rga\_consistent\_}\allowbreak\textsf{via\_wtrue} & \textsf{RGA\_Decider/\allowbreak Decider\_Truth} & \ref{thm:pl-con} \\
\textsf{rga\_complete\_sch} & \textsf{RGA\_Decider/\allowbreak Decider\_Truth} &  \\
\textsf{rga\_truth\_}\allowbreak\textsf{reflect\_sch} & \textsf{RGA\_Decider/\allowbreak Decider\_Truth} &  \\
\multicolumn{3}{@{}l}{\textbf{A8.}\ Internal completeness} \\*
\textsf{rga\_complete\_int} & \textsf{RGA\_Checker/\allowbreak Chk\_Completeness} &  \\
\multicolumn{3}{@{}l}{\textbf{A9.}\ Faithfulness of the tree checker} \\*
\textsf{iChk\_faithful} & \textsf{RGA\_Checker/\allowbreak Chk\_Faithful} & \ref{thm:in-chk} \\
\textsf{iChk\_genuine} & \textsf{RGA\_Checker/\allowbreak Chk\_Faithful} & \ref{thm:in-chk} \\
\multicolumn{3}{@{}l}{\textbf{A10.}\ G\"odel's theorems} \\*
\textsf{rga\_goedel2} & \textsf{RGA\_Goedel/\allowbreak Goedel\_II\_PR} & thm.\ in [\citen{goedel-paper}] \\
\textsf{rga\_goedel2\_}\allowbreak\textsf{via\_search} & \textsf{RGA\_Goedel/\allowbreak Goedel\_II} &  \\
\textsf{rga\_goedel1} & \textsf{RGA\_Goedel/\allowbreak Goedel\_II} & thm.\ in [\citen{goedel-paper}] \\
\textsf{rga\_goedel1\_}\allowbreak\textsf{unrefutable} & \textsf{RGA\_Goedel/\allowbreak Goedel\_II} & thm.\ in [\citen{goedel-paper}] \\
\multicolumn{3}{@{}l}{\textbf{B1.}\ The Gentzen shape} \\*
\textsf{rga\_con\_int\_}\allowbreak\textsf{of\_gentzen} & \textsf{ConRGA/\allowbreak Gentzen\_Shape} &  \\
\textsf{rga\_con\_int\_of\_}\allowbreak\textsf{gentzen\_gen} & \textsf{ConRGA/\allowbreak Gentzen\_Shape} &  \\
\textsf{rga\_con\_int\_}\allowbreak\textsf{of\_gentzen\_ex} & \textsf{ConRGA/\allowbreak Gentzen\_Shape} &  \\
\textsf{no\_gentzen\_pair} & \textsf{ConRGA/\allowbreak Gentzen\_Shape} &  \\
\textsf{gentzen\_pair\_recipe} & \textsf{ConRGA/\allowbreak Gentzen\_Shape} &  \\
\textsf{inf\_of\_gentzen} & \textsf{ConRGA/\allowbreak Gentzen\_Shape} &  \\
\textsf{rga\_con\_int\_of\_INFF} & \textsf{ConRGA/\allowbreak Gentzen\_Shape} &  \\
\textsf{rga\_not\_INF} & \textsf{ConRGA/\allowbreak Gentzen\_Shape} &  \\
\multicolumn{3}{@{}l}{\textbf{B2.}\ Separations between \rga and \oga} \\*
\textsf{rga\_not\_iLF0} & \textsf{ConOGA/\allowbreak ConOGA\_Otto} &  \\
\textsf{rga\_not\_iLF} & \textsf{ConOGA/\allowbreak ConOGA\_Otto} &  \\
\textsf{rga\_not\_iLFs} & \textsf{ConOGA/\allowbreak ConOGA\_Otto} & \ref{thm:e-olf-N} \\
\textsf{rga\_oga\_}\allowbreak\textsf{separation\_iLF0} & \textsf{ConOGA/\allowbreak ConOGA\_Otto} &  \\
\textsf{rga\_oga\_}\allowbreak\textsf{separation\_iLF} & \textsf{ConOGA/\allowbreak ConOGA\_Otto} &  \\
\textsf{rga\_oga\_}\allowbreak\textsf{separation\_iLFs} & \textsf{ConOGA/\allowbreak ConOGA\_Otto} & \ref{thm:e-olf-N} \\
\textsf{iLF\_diag\_recipe} & \textsf{ConOGA/\allowbreak ConOGA\_Otto} &  \\
\textsf{rga\_oga\_separation\_}\allowbreak\textsf{goedelW} & \textsf{ConOGA/\allowbreak ConOGA\_Otto} &  \\
\textsf{goedelW\_decided\_oga} & \textsf{ConOGA/\allowbreak ConOGA\_Otto} &  \\
\textsf{goedelW\_}\allowbreak\textsf{unrefutable\_oga} & \textsf{ConOGA/\allowbreak ConOGA\_Otto} &  \\
\textsf{goedelW\_true\_}\allowbreak\textsf{oga\_model} & \textsf{ConOGA/\allowbreak ConOGA\_Otto} &  \\
\multicolumn{3}{@{}l}{\textbf{B3.}\ Value records and halting status} \\*
\textsf{iRec\_decided\_}\allowbreak\textsf{iff\_orec} & \textsf{RGA\_Liar\_Div/\allowbreak Record\_Halting} &  \\
\textsf{halting\_status\_}\allowbreak\textsf{iff\_orec} & \textsf{RGA\_Liar\_Div/\allowbreak Record\_Halting} &  \\
\textsf{iRec\_faithful} & \textsf{RGA\_Liar\_Div/\allowbreak Record\_Halting} & \ref{thm:c-orec} \\
\textsf{iRec\_of\_wosr} & \textsf{RGA\_Liar\_Div/\allowbreak Record\_Halting} & \ref{thm:c-orec} \\
\textsf{iRec\_exact\_of\_wosr} & \textsf{RGA\_Liar\_Div/\allowbreak Record\_Halting} &  \\
\textsf{iRec\_neither\_of\_list} & \textsf{RGA\_Liar\_Div/\allowbreak Record\_Halting} &  \\
\textsf{iRec\_decided\_iff} & \textsf{RGA\_Liar\_Div/\allowbreak Record\_Halting} &  \\
\textsf{iRec\_true\_}\allowbreak\textsf{decided\_iff} & \textsf{RGA\_Liar\_Div/\allowbreak Record\_Halting} &  \\
\textsf{halting\_status\_iff} & \textsf{RGA\_Liar\_Div/\allowbreak Record\_Halting} &  \\
\textsf{orec\_iff\_wosr} & \textsf{RGA\_Liar\_Div/\allowbreak Record\_Halting} & \ref{thm:c-orec} \\
\multicolumn{3}{@{}l}{\textbf{B4.}\ The Liars and their records} \\*
\textsf{liar\_record\_decided} & \textsf{RGA\_Liar\_Div/\allowbreak Liar\_Record} &  \\
\textsf{liar\_iRec\_none} & \textsf{RGA\_Liar\_Div/\allowbreak Liar\_Record} &  \\
\textsf{omega\_record\_decided} & \textsf{RGA\_Liar\_Div/\allowbreak Liar\_Record} &  \\
\textsf{goedelW\_no\_value} & \textsf{RGA\_Liar\_Div/\allowbreak Liar\_Record} & \ref{thm:pl-no-value} \\
\textsf{goedelW\_unprovable} & \textsf{RGA\_Liar\_Div/\allowbreak Liar\_Record} & \ref{thm:pl-liar-indep} \\
\textsf{goedelW\_unrefutable} & \textsf{RGA\_Liar\_Div/\allowbreak Liar\_Record} & \ref{thm:pl-liar-indep} \\
\textsf{goedelW\_undecided} & \textsf{RGA\_Liar\_Div/\allowbreak Liar\_Record} & \ref{thm:pl-liar-indep} \\
\textsf{goedelW\_instances\_}\allowbreak\textsf{refuted} & \textsf{RGA\_Liar\_Div/\allowbreak Liar\_Record} &  \\
\textsf{goedelW\_omega\_}\allowbreak\textsf{incomplete} & \textsf{RGA\_Liar\_Div/\allowbreak Liar\_Record} & \ref{thm:br-omega} \\
\textsf{goedelW\_true\_orec} & \textsf{RGA\_Liar\_Div/\allowbreak Liar\_Record} &  \\
\textsf{liarRec\_prov\_iff} & \textsf{RGA\_Liar\_Div/\allowbreak Liar\_Record} &  \\
\textsf{liarRec\_refute\_iff} & \textsf{RGA\_Liar\_Div/\allowbreak Liar\_Record} &  \\
\textsf{liarRec\_rec\_}\allowbreak\textsf{decided\_iff} & \textsf{RGA\_Liar\_Div/\allowbreak Liar\_Record} &  \\
\multicolumn{3}{@{}l}{\textbf{B5.}\ Normal forms} \\*
\textsf{rga\_normal\_form} & \textsf{RGA/\allowbreak RGA\_Normal} &  \\
\textsf{canj\_unchecked\_}\allowbreak\textsf{leaf\_0eq1} & \textsf{RGA/\allowbreak RGA\_Normal} &  \\
\textsf{rga\_afree\_normal} & \textsf{RGA/\allowbreak RGA\_Normal\_QF} &  \\
\textsf{rga\_afree\_}\allowbreak\textsf{normal\_strong} & \textsf{RGA/\allowbreak RGA\_Normal\_QF} &  \\
\textsf{rga\_afree\_normal\_chk} & \textsf{RGA\_Checker/\allowbreak Chk\_QF\_Faithful} & \ref{thm:con-normal} \\
\multicolumn{3}{@{}l}{\textbf{B6.}\ Fuel and size calibration} \\*
\textsf{fuel\_escape\_recipe} & \textsf{RGA/\allowbreak RGA\_FuelBound} &  \\
\textsf{Gmin\_not\_ack\_bounded} & \textsf{RGA/\allowbreak RGA\_FuelBound} &  \\
\textsf{fuel\_escape\_}\allowbreak\textsf{recipe\_gen} & \textsf{RGA/\allowbreak RGA\_FuelBound} &  \\
\textsf{Gmin\_escapes} & \textsf{RGA/\allowbreak RGA\_FuelBound} &  \\
\textsf{qf\_size\_lower} & \textsf{RGA/\allowbreak RGA\_NormalSize} &  \\
\textsf{nf\_size\_}\allowbreak\textsf{escape\_recipe} & \textsf{RGA/\allowbreak RGA\_NormalSize} &  \\
\textsf{nf\_size\_escape\_}\allowbreak\textsf{recipe\_ack} & \textsf{RGA/\allowbreak RGA\_NormalSize} &  \\
\textsf{nf\_size\_escapes} & \textsf{RGA/\allowbreak RGA\_NormalSize} &  \\
\textsf{nf\_size\_not\_}\allowbreak\textsf{ack\_bounded} & \textsf{RGA/\allowbreak RGA\_NormalSize} &  \\
\multicolumn{3}{@{}l}{\textbf{B7.}\ Divergence and invariants} \\*
\textsf{div\_inv} & \textsf{RGA\_Liar\_Div/\allowbreak Divergence} &  \\
\textsf{div\_inv\_of\_list} & \textsf{RGA\_Liar\_Div/\allowbreak Divergence} &  \\
\textsf{grow\_no\_list} & \textsf{RGA\_Liar\_Div/\allowbreak Divergence} &  \\
\textsf{grow\_diverges} & \textsf{RGA\_Liar\_Div/\allowbreak Divergence} &  \\
\textsf{div\_inv\_trivial\_not} & \textsf{RGA\_Liar\_Div/\allowbreak Divergence\_Rosser} &  \\
\textsf{div\_inv\_not\_iff} & \textsf{RGA\_Liar\_Div/\allowbreak Divergence\_Rosser} &  \\
\textsf{div\_inv\_}\allowbreak\textsf{commit\_eq\_mem} & \textsf{RGA\_Liar\_Div/\allowbreak Divergence\_Rosser} &  \\
\textsf{div\_inv\_commit\_eq} & \textsf{RGA\_Liar\_Div/\allowbreak Divergence\_Rosser} &  \\
\textsf{div\_inv\_pure} & \textsf{RGA\_Liar\_Div/\allowbreak Divergence\_Rosser} &  \\
\textsf{div\_inv\_ghost} & \textsf{RGA\_Liar\_Div/\allowbreak Divergence\_Rosser} &  \\
\textsf{div\_inv\_ghost\_ex} & \textsf{RGA\_Liar\_Div/\allowbreak Divergence\_Rosser} &  \\
\textsf{div\_inv\_via\_ghost} & \textsf{RGA\_Liar\_Div/\allowbreak Divergence\_Rosser} &  \\
\textsf{grow\_diverges\_}\allowbreak\textsf{via\_ghost} & \textsf{RGA\_Liar\_Div/\allowbreak Divergence\_Rosser} &  \\
\textsf{rga\_con\_int\_}\allowbreak\textsf{of\_con\_search} & \textsf{RGA\_Liar\_Div/\allowbreak Divergence} & \ref{thm:c-search} \\
\textsf{con\_search\_}\allowbreak\textsf{unrefutable} & \textsf{RGA\_Liar\_Div/\allowbreak Divergence} & \ref{thm:c-search} \\
\textsf{con\_search\_}\allowbreak\textsf{no\_invariant} & \textsf{RGA\_Liar\_Div/\allowbreak Divergence} &  \\
\textsf{search\_diverges} & \textsf{RGA\_Liar\_Div/\allowbreak Divergence} &  \\
\textsf{con\_search\_}\allowbreak\textsf{closure\_iff} & \textsf{RGA\_Liar\_Div/\allowbreak Divergence} &  \\
\textsf{frag\_search\_diverges} & \textsf{RGA\_Liar\_Div/\allowbreak Divergence} & \ref{thm:c-search} \\
\textsf{rosser\_search\_split} & \textsf{RGA\_Liar\_Div/\allowbreak Divergence\_Rosser} &  \\
\textsf{rosser\_diverges} & \textsf{RGA\_Liar\_Div/\allowbreak Divergence\_Rosser} &  \\
\textsf{rosser\_a\_unprovable} & \textsf{RGA\_Liar\_Div/\allowbreak Divergence\_Rosser} &  \\
\textsf{rosser\_b\_unprovable} & \textsf{RGA\_Liar\_Div/\allowbreak Divergence\_Rosser} &  \\
\textsf{rosser\_no\_invariant} & \textsf{RGA\_Liar\_Div/\allowbreak Divergence\_Rosser} &  \\
\textsf{rosser\_no\_}\allowbreak\textsf{invariant\_pure} & \textsf{RGA\_Liar\_Div/\allowbreak Divergence\_Rosser} &  \\
\textsf{rosser\_}\allowbreak\textsf{counterexample} & \textsf{RGA\_Liar\_Div/\allowbreak Divergence\_Rosser} &  \\
\textsf{infinite\_collapse} & \textsf{RGA\_Liar\_Div/\allowbreak Divergence\_Rosser} &  \\
\textsf{infinite\_collapse\_}\allowbreak\textsf{diverges} & \textsf{RGA\_Liar\_Div/\allowbreak Divergence\_Rosser} &  \\
\textsf{width\_clause\_}\allowbreak\textsf{diverges} & \textsf{RGA\_Liar\_Div/\allowbreak Divergence\_Rosser} &  \\
\textsf{width\_subclause\_}\allowbreak\textsf{unprovable} & \textsf{RGA\_Liar\_Div/\allowbreak Divergence\_Rosser} &  \\
\textsf{width\_counterexample} & \textsf{RGA\_Liar\_Div/\allowbreak Divergence\_Rosser} &  \\
\multicolumn{3}{@{}l}{\textbf{B8.}\ The HOL model of the residual, unconditional} \\*
\textsf{shadow\_consistent} & \textsf{Shadows/\allowbreak Shadow\_Integrated} &  \\
\textsf{shadow\_sound} & \textsf{Shadows/\allowbreak Shadow\_Integrated} &  \\
\textsf{shadow\_res\_holds} & \textsf{Shadows/\allowbreak Shadow\_Integrated} &  \\
\multicolumn{3}{@{}l}{\textbf{B9.}\ Companions of G\"odel's second theorem} \\*
\textsf{goedelSigma\_is\_pr} & \textsf{RGA\_Goedel/\allowbreak Goedel\_II\_PR} & thm.\ in [\citen{goedel-paper}] \\
\textsf{D3G\_goedelSigma} & \textsf{RGA\_Goedel/\allowbreak Goedel\_II\_PR} & thm.\ in [\citen{goedel-paper}] \\
\textsf{pr\_cert\_builder} & \textsf{RGA\_Goedel/\allowbreak Goedel\_II\_PR} & thm.\ in [\citen{goedel-paper}] \\
\textsf{goedelE\_chk} & \textsf{RGA\_Goedel/\allowbreak Goedel\_II} & thm.\ in [\citen{goedel-paper}] \\
\textsf{D3Gdec\_holds} & \textsf{RGA\_Goedel/\allowbreak Goedel\_II} &  \\
\textsf{goedelG\_of\_con\_}\allowbreak\textsf{via\_search} & \textsf{RGA\_Goedel/\allowbreak Goedel\_II} &  \\
\multicolumn{3}{@{}l}{\textbf{B10.}\ Machine results stated in the paper} \\*
\textsf{prDecTab\_correct} & \textsf{RGA\_Decider/\allowbreak RGA\_Decide\_PR} & \ref{thm:in-tab} \\
\textsf{lemQ\_int} & \textsf{ConRGA/\allowbreak Wfd\_Decision} & \ref{thm:con-lemq} \\
\multicolumn{3}{@{}l}{\textbf{B11.}\ Design checks on HOL models} \\*
\textsf{shadow\_consistent\_qf} & \textsf{Shadows/\allowbreak Shadow\_QF} &  \\
\textsf{shadow\_consistent\_}\allowbreak\textsf{of\_tot} & \textsf{Shadows/\allowbreak Shadow\_Findings} &  \\
\textsf{shadow\_consistent\_}\allowbreak\textsf{of\_wfd} & \textsf{Shadows/\allowbreak Shadow\_Findings} &  \\
\textsf{arB\_decided} & \textsf{ConOGA/\allowbreak ConOGA\_Shadow} &  \\
\textsf{oga\_liar\_undecided} & \textsf{ConOGA/\allowbreak ConOGA\_Shadow} &  \\
\end{longtable}
\endgroup

\subsection{The headline statements, as Isabelle prints them}
\label{app:endpoints-statements}

This subsection shows twenty-five theorems and five definitions behind
the paper's headline results as Isabelle prints them.  These statements are
typeset from an Isabelle build of the development: a
document-only theory prints each fact with Isabelle's \texttt{thm}
command, and the cartouches tool\footnote{\url{https://github.com/astahfrom/cartouches}.} cuts the generated \LaTeX\ into the
displays below, so no statement here is transcribed by hand.  Pollack
separates two questions about a machine-checked proof, whether the proof is a
derivation and whether its statement means what is
claimed~\cite{pollack98believe}; generated statements let a reader
answer the second question against the formula that the prover checked.  Generated statements
do not settle the second question alone, since a printer can print two different terms
alike, and Wiedijk shows that Isabelle's printer can~\cite{wiedijk12pollack};
the fingerprints of the trusted definitions (Section~\ref{sec:mech})
cover what a printout cannot.  The notation is Isabelle's
(Appendix~\ref{app:notation} maps the main text's symbols to Isabelle
constants).  Variables are de~Bruijn indices
(\isa{\isactrlbold v\isactrlsub {\isadigit{0}}} is bound by the
innermost binder, and \isa{{\isasymupharpoonleft}} lifts the indices
under a new binder), and \rga's operators are bold, to set them apart
from HOL's.  \isa{\isavar{A}\ {\isasymtriangleright}\ \isavar{{\isasymGamma}}}
adds the hypothesis $A$ to the context $\Gamma$; postfix
\isa{\isactrlbold N} and \isa{\isactrlbold B} are the judgments
$\jnat$ and $\jbool$; \isa{\isactrlbold n{\isasymlbrace}\isavar{t}{\isasymrbrace}}
is the main text's code numeral $\bquo{t}$, and
\isa{\isactrlbold n\isaconst{Suc}\ {\isadigit{0}}} the numeral
$\tofn 1$; \isa{\isavar{P}\ {\isasymP}{\isasymlbrakk}{\isasymGamma}\ {\isasymturnstile}\ \isavar{c}{\isasymrbrakk}}
says that the proof $P$ proves $\Gamma \vdash c$, written $P \proves \Gamma \vdash c$ in
the main text.  A statement holds
for every value of its free variables.  Abbreviations print expanded:
the universal $\forall$ prints as \textbf{A} applied to a
\textbf{$\lambda$}, the existential as $\neg$\textbf{A} applied to a
\textbf{$\lambda$} of a negation, the 0/1 tests through
\textsf{gapSub}, and \oga's entailment as
\textsf{prf\_derives.entails oga.derives}.  Where a numbered result
has several Isabelle facts, those shown here state the result's
headline form; Table~\ref{tab:g-paper} lists the others.

\begingroup
\newcommand{\GenRef}[1]{{\normalfont\normalsize #1}\enspace}
\newcommand{\GenFact}[2]{%
  \par\noindent\begin{minipage}{\linewidth}
  \begin{isabelle}\raggedright
  #1\GetSnip{ConRGA-Paper-Statements:thm:#2-0}%
  \hangindent=1.5em\hangafter=1\relax
  \GetSnip{ConRGA-Paper-Statements:thm:#2-1}%
  \end{isabelle}
  \end{minipage}\par\medskip}
\newcommand{\GenSec}[1]{\par\bigskip\noindent\textbf{#1}\par\nopagebreak\smallskip}

\GenSec{Section~\ref{sec:inside}.}
\GenFact{\GenRef{Theorem~\ref{thm:in-dt} (also [\citen{goedel-paper}, theorem~\textsf{sigma\_intro}]).}}{deduction-theorem-exact}
\GenFact{}{sigma-deduction-theorem}
\GenFact{\GenRef{Definition~\ref{def:in-truth}.}}{iPrDecTrueAt-def}
\GenFact{}{iWtrue-def}
\GenFact{\GenRef{Definition~\ref{def:in-chk}.}}{iChk-def}
\GenFact{\GenRef{The consistency sentence $\ConR$ (\S\ref{sec:inside-chk}, [\citen{goedel-paper}, definition~\textsf{goedelD}]).}}{ConRGA-def}
\GenFact{\GenRef{The code $\bquo{\z = \tofn 1}$ (\S\ref{sec:inside-chk}).}}{cEq01-def}

\GenSec{Section~\ref{sec:pillars}.}
\GenFact{\GenRef{Theorem~\ref{thm:pl-con}.}}{rga-consistent-via-wtrue}
\GenFact{\GenRef{Theorem~\ref{thm:pl-nec}.}}{establishment-transparency}
\GenFact{\GenRef{Theorem~\ref{thm:pl-int}.}}{rga-complete-int2}
\GenFact{\GenRef{Theorem~\ref{thm:pl-int-ent}.}}{complete-ent}
\GenFact{\GenRef{Theorem~\ref{thm:pl-tb-generic}.}}{tbicond-liar-explosive-diag}
\GenFact{\GenRef{Theorem~\ref{thm:pl-liar-neither}.}}{ordinary-liar-neither}
\GenFact{\GenRef{Theorem~\ref{thm:pl-refuted}.}}{refuted-unestablished-trichotomy}
\GenFact{\GenRef{Theorem~\ref{thm:pl-tracking}.}}{con-tschema-decided}

\GenSec{[\citen{goedel-paper}, Section ``G\"odel's proof in four forms''].}
\GenFact{\GenRef{[\citen{goedel-paper}, theorem~\textsf{goedelG\_under\_ConRGA}] (also [\citen{goedel-paper}, theorem~\textsf{goedelG\_of\_expl\_d3g\_ctx}]).}}{goedelG-under-ConRGA}

\GenSec{[\citen{goedel-paper}, Section ``G\"odel's theorems for \rga''].}
\GenFact{\GenRef{[\citen{goedel-paper}, theorem~\textsf{rga\_goedel1}] (also [\citen{goedel-paper}, theorem~\textsf{rga\_goedel1}]).}}{rga-goedel1}
\GenFact{}{rga-goedel1-unrefutable}
\GenFact{\GenRef{[\citen{goedel-paper}, theorem~\textsf{rga\_goedel2}] (also [\citen{goedel-paper}, theorem~\textsf{rga\_goedel2\_of\_expl\_d3g}]).}}{rga-goedel2}
\GenFact{\GenRef{[\citen{goedel-paper}, theorem~\textsf{rga\_goedel2\_internal\_contra}].}}{rga-goedel2-internal-contra}
\GenFact{}{rga-goedel2-internal}
\GenFact{\GenRef{[\citen{goedel-paper}, theorem~\textsf{rga\_goedel2\_sentence}].}}{rga-goedel2-sentence}
\GenFact{}{goedel2T-sentence}
\GenFact{}{rga-goedel2-sentence-witness}

\GenSec{Section~\ref{sec:con}.}
\GenFact{\GenRef{Theorem~\ref{thm:con-qf}.}}{rga-con-qf-int}
\GenFact{\GenRef{Theorem~\ref{thm:con-frag}.}}{rga-con-frag-int}
\GenFact{\GenRef{Theorem~\ref{thm:con-lv}.}}{rga-con-lvi-under-lf-all}
\GenFact{\GenRef{Theorem~\ref{thm:con-oga}.}}{rga-con-lvi-oga}
\GenFact{\GenRef{Theorem~\ref{thm:con-res}.}}{rga-con-under-res}

\GenSec{[\citen{goedel-paper}, Appendix ``G\"odel II in full''].}
\GenFact{\GenRef{Beklemishev and Shamkanov's conditions ([\citen{goedel-paper}, \S\,``The entailments and the cut transformers'']).}}{bs-goedel2}
\endgroup

\bibliographystyle{alpha}
\bibliography{logic}

\end{document}